\documentclass[a4paper,11pt]{article}
\usepackage[a4paper,left=23mm,right=23mm,top=22mm,bottom=23mm,footskip=11mm]{geometry}
\usepackage[T1]{fontenc}
\usepackage{lmodern,amsmath,amssymb,amsthm,mathtools,array,booktabs}
\usepackage[unicode,hidelinks,bookmarksopen=false]{hyperref}
\usepackage{orcidlink}
\hypersetup{pdfauthor={Philipp Bogdan},pdftitle={On fourteen problems in the theory of composites},pdfdisplaydoctitle=true}
\newtheorem{theorem}{Theorem}
\newtheorem{lemma}[theorem]{Lemma}
\newtheorem{proposition}[theorem]{Proposition}
\newtheorem{corollary}[theorem]{Corollary}
\newlength{\tableusable}
\newcommand{\Milton}[1]{\par\smallskip\begingroup\leftskip=1em\rightskip=.6em\itshape
\hypertarget{milton-problem-#1}{}\textup{\textbf{Problem #1.}}\enspace\csname question#1\endcsname\par\endgroup\smallskip}
\expandafter\def\csname question1\endcsname{Is rank convexity equal to quasiconvexity for the two well problem in two spatial dimensions?}
\expandafter\def\csname question2\endcsname{For two-phase composites in two spatial dimensions, such that phase 1 occupies a volume fraction $f$, is the $G_f$-closure equal to its lamination closure when the fields on the right of the constitutive law have $n$ components, each being the sum of a real 2 component vector and the gradient of a scalar periodic potential, while the fields on the left of the constitutive law also have $n$ components, each having zero divergence, in which $n$ is an arbitrary positive integer?}
\expandafter\def\csname question3\endcsname{For three-dimensional isotropic composites each having an effective conductivity $\sigma_* I$ and being built from two isotropic materials having conductivities $\sigma I$ and $I$, can one characterize all possible conductivity functions $\sigma_*(\sigma)$?}
\expandafter\def\csname question4\endcsname{For three-dimensional isotropic composites of two isotropic phases, are all possible conductivity functions $\sigma_*(\sigma)$ achievable by multiple rank laminate microstructures and, if so, does it suffice to consider laminate microstructures where one laminates only in mutually orthogonal directions?}
\expandafter\def\csname question5\endcsname{For two-dimensional anisotropic composites of two anisotropic phases, are all possible conductivity functions $\boldsymbol\sigma_*(\boldsymbol\sigma_1,\boldsymbol\sigma_2)$ achievable by multiple rank laminate microstructures?}
\expandafter\def\csname question6\endcsname{Are there other geometries with isotropic effective permittivities that lie on the arc?}
\expandafter\def\csname question7\endcsname{Can any of these discovered geometries, having maximum absorption, can be replaced by simpler ones?}
\expandafter\def\csname question8\endcsname{Can one complete the characterization of the $G$-closure for a void (or rigid) phase mixed with an isotropic elastic phase?}
\expandafter\def\csname question9\endcsname{The bounds (4.1) imply $\mu_* -c\kappa_*\le\mu$ for all $c>0$. Can this inequality be improved, in 2 and/or 3 dimensions, for a range of $c>0$? Alternatively, can one construct composites of an elastic phase with void with $(\kappa_*,\mu_*)$ arbitrarily close to $(0,\mu)$?}
\expandafter\def\csname question10\endcsname{Identify, for given $c>0$, in 2 and/or 3 dimensions, isotropic microstructures of an elastic phase with void that have the largest possible value of $\mu_*-c\kappa_*$ (or a sequence of isotropic microstructures with moduli such that $\mu_*-c\kappa_*$ converges to its largest possible value).}
\expandafter\def\csname question11\endcsname{Is $4\mu/5$ the largest possible value of $\mu_*$ for a two-dimensional elastic material with voids, given that $\kappa_*=1/\kappa=0$?}
\expandafter\def\csname question12\endcsname{Can one obtain bounds that correlate the possible compressive and shear deformations of composites when these deformations are not infinitesimal?}
\expandafter\def\csname question13\endcsname{What are the possible $(\kappa_*,\mu_*)$-pairs for three-dimensional isotropic elastic polycrystals (composites built from a single crystal in various orientations)? The bounds of Hill [47] are optimal for $\kappa_*$ [7], but improved bounds for $\mu_*$ or $(\kappa_*,\mu_*)$ pairs are lacking. Hashin and Shtrikman obtained improved bounds on $\mu_*$ [46], but only under additional assumptions about crystal orientations, that are not generally valid.}
\expandafter\def\csname question14\endcsname{Can one obtain improved bounds on the elastic moduli pairs of isotropic composites of two isotropic elastic phases, and ultimately find the optimal ones?}

\usepackage{longtable,etoolbox,calc,graphicx}
\usepackage{microtype}
\makeatletter
\newcommand{\dashint}{\mathop{\mathchoice{\vcenter{\hbox{$\displaystyle\int$}}\kern-.72em-\kern.28em}{\vcenter{\hbox{$\textstyle\int$}}\kern-.53em-\kern.10em}{\int}{\int}}\nolimits}
\makeatother
\newcommand{\tightlist}{\setlength{\itemsep}{0pt}\setlength{\parskip}{0pt}}

\providecommand{\mathscr}{\mathcal}
\providecommand{\fint}{\dashint}

\newcommand{\SSection}[2]{\section{#2}\label{sup:#1}}
\begin{document}
{\raggedright\fontsize{23}{27}\selectfont On fourteen problems in the theory of composites\par}
\smallskip
{\large Philipp Bogdan\,\orcidlink{0009-0009-7405-3555}\par}
{\itshape Imperial College London}\quad
\href{mailto:philipp.bogdan25@imperial.ac.uk}{\texttt{philipp.bogdan25@imperial.ac.uk}}\par
{\small Revised 21 November 2026}\par
\smallskip
\textbf{Abstract.}
We study fourteen questions of Milton on the effective properties of composites.
For one scalar potential every effective tensor of two fixed planar phases
has a laminate representative. Two coupled potentials suffice for failure.
In space, dyadic resistor networks characterize the isotropic conductivity
closure, and some laws in it have no finite-perimeter periodic representative.
We construct intervals on the single-pole absorption arc and hole-free
planar inclusions realizing optimal double-coating laws. For a solid with
void or perfectly rigid inclusions we characterize the two- and
three-dimensional tensor closures. Porous solids can approach full shear
at vanishing bulk in every fixed dimension $d\ge2$. Prescribing a finite
planar shortening gives a sharp loss of initial shear. We give the complete
shear interval at each bulk for a polycrystal. For two positive elastic
phases, changing only the bulk moduli can reverse which dimension admits
a shear common to every attainable bulk.

\section*{Results and open questions}
\pdfbookmark[1]{Results and open questions}{results-guide}
Here $\kappa_*,\mu_*$ are the effective bulk and shear moduli, $\mu_s$ is the shear modulus of the solid, and $f$ denotes a phase fraction.

\begingroup\small
\newcommand{\guide}[3]{\hyperlink{milton-problem-#1}{\textbf{#1}} & #3 & \pageref{#2}\\[2pt]}
\begin{tabular}{@{}p{1.8em}@{\hspace{.55em}}>{\raggedright\arraybackslash}p{\dimexpr\linewidth-4.9em\relax}@{\hspace{.55em}}r@{}}
\guide{1}{problem-1}{Two positive quadratic wells on $\mathbb R^{2\times4}$ have distinct rank-one convex and quasiconvex envelopes.}
\guide{2}{problem-1}{One potential always admits a laminate at fixed fraction. Two coupled potentials suffice for failure at $11/25$.}
\guide{3}{problem-3}{A graph limit characterizes the scalar conductivity closure. Some laws have no finite-perimeter cell.}
\guide{4}{problem-4}{Three directions per tree do not suffice in space. The unrestricted laminate question is open.}
\guide{5}{problem-5}{Finite real tests separate a four-rectangle response even from convex averages of laminate responses.}
\guide{6}{problem-6}{Laminates attain intervals on the single-pole arc and miss three open gaps at quarter fraction.}
\guide{7}{problem-7}{Connected inclusions replace double coatings. Hole-free planar maximizers exist over a fraction interval and on a second branch.}
\guide{8}{problem-8}{A finite moment test characterizes the two- and three-dimensional void and rigid tensor closures.}
\guide{9}{problem-9}{Exactly isotropic porous solids approach zero bulk and full constituent shear for each fixed $d\ge2$.}
\guide{10}{problem-9}{The same sequence reaches the limiting supremum $\mu_s$ of $\mu_*-c\kappa_*$ for every $c>0$.}
\guide{11}{problem-11}{An incompressible planar solid has full limiting shear at zero bulk. A finite stroke costs the sharp factor $\lambda^2/(2-\lambda^2)$.}
\guide{12}{problem-12}{Mixed loading depends on the phase laws. Two specified rubbers have sharp shear supremum $0.531637\ldots$ at the stated compression.}
\guide{13}{problem-13}{Four variational families give the sharp shear interval at each bulk. Algebraic data admit two-sided computation.}
\guide{14}{problem-14}{Changing only the phase bulks reverses which dimension admits a common shear across the entire bulk range.}
\end{tabular}
\endgroup

\clearpage
The questions are quoted by number from Milton \cite[Problems 1 to 14, pp.~2 to 16 of the preprint]{milton}. Reference and equation numbers inside the quotations are his. His bounds (4.1) are $0\le\kappa_*\le\kappa$ and $0\le\mu_*\le\mu$. For an energy $W$, $QW$ and $RW$ denote its quasiconvex and rank-one convex envelopes. We normalize cell volume to one and write $\langle\cdot\rangle$ for the mean. At fixed fraction $f$, each cell of an approximating sequence has that fraction. A response function uses the same geometry at every contrast. In a hierarchical laminate each finer ply is homogenized before the coarser one. We use $e(u)=(\nabla u+\nabla u^T)/2$, $E_d=E-(\operatorname{tr}E)I/d$ and $A:B=\operatorname{tr}(A^TB)$. An isotropic stiffness acts by $CE=\kappa\operatorname{tr}(E) I+2\mu E_d$. Energies use $E:CE$ unless a factor $1/2$ is specified. Void carries zero stress, and a perfectly rigid phase admits only zero strain.

\section{5. A checkerboard outside the lamination closure}\label{problem-5}
\Milton{5}

The answer is no, even when all the measurements use real, positive-definite conductivities. Cut the unit square at $x_1=2/5$ and $x_2=1/5$. Take phase 1 where exactly one coordinate is below its cut. Use this same checkerboard in every conductivity test.
\begin{theorem}\label{p5:obstruction}
For real symmetric positive-definite phase matrices, the checkerboard response is outside even the closed convex hull of the responses of planar finite laminate trees of arbitrary fractions, in the topology of pointwise convergence. Finitely many energy measurements separate it strictly from that hull.
\end{theorem}
\begin{proof}\phantomsection\label{proof:p5:obstruction}
\emph{What the measurements determine.}

Write \(\chi\in\{0,1\}\) for the phase indicator, \(f=\mathbb E\chi\), and \(D=\mathbb {RP}^1\)
for the unoriented directions. In a periodic cell, let \(P_j\) project onto the nonconstant
Fourier modes in direction \(n_j\). On a finite laminate, give leaves their volume probabilities
and use one \(P_j\) per split with normal \(n_j\). It takes the child average minus the parent
average on that node and is zero elsewhere. In both cases the projections are mutually orthogonal
and sum to the projection removing the mean. Set
\[
T_a=\sum_j a(n_j)P_j,\qquad
\mu=\sum_j\|P_j\chi\|_2^2\delta_{n_j},\qquad
K(a,b)=\mathbb E[\chi(T_a\chi)(T_b\chi)] .
\]

Thus \(\mathbb E[\chi T_a\chi]=\mu(a)\), \(\mathbb E[(T_a\chi)(T_b\chi)]=\mu(ab)\), and
\[
\mu(D)=f(1-f)\le\tfrac14,\qquad
|K(a,b)|\le\|a\|_{L^2(\mu)}\|b\|_{L^2(\mu)}. \tag{5.1}
\]

We extend \(K\) complex bilinearly. Then (5.1) holds for complex symbols as well, with
\(\|a\|_{L^2(\mu)}^2=\mu(|a|^2)\).

Use \(\sigma_2=A_s=\left(\begin{smallmatrix}1&s\\s&1\end{smallmatrix}\right)\), \(0<s<1/2\), and
\(\sigma_1=A_s+t vv^T\), \(t>0\), for \(v=e_1,e_2,e_1+e_2\). Prescribe the mean field
\(e=v/|v|^2\). For the local field \(E\), put \(z=v^TE\). On each nonconstant mode the field is
parallel to its frequency and the flux is perpendicular to it. On each ideal, scale-separated
laminate split, continuity of tangential field and normal flux gives exactly \(P_j E\parallel
n_j\) and \(n_j^T P_j(A_s E+t\chi v z)=0\) for the child averages. Solving this scalar equation in
each direction gives
\[
z=1-tT_{m_s}(\chi z),\qquad
m_s(n)=\frac{(v^Tn)^2}{1+2s n_1n_2},\quad |n|=1 .
\]

The measured response is \(e^T(\sigma_*-A_s)e=t\mathbb E[\chi z]\). Since \(\|T_{m_s}\|\le4\), the
equation for \(z\) has a unique square-integrable solution when \(|t|\le1/8\), namely its Neumann
series. Moving one self-adjoint \(T_{m_s}\) across the inner product in the cubic term gives
\[
t f-t^2\mu(m_s)+t^3K(m_s,m_s)+O(t^4). \tag{5.2}
\]

The remainder is at most \(128|t|^4\) for \(|t|\le1/8\), at every depth.

Suppose that convex combinations of tree responses approximate the cell pointwise. In each combination average the coefficients $f,\mu,K$ with the same weights. Inequality (5.1) for the resulting measure and form follows by Cauchy-Schwarz in the weights, and the remainder in (5.2) stays uniform. Choose combinations approximating increasing finite lists
of these three probes at rational \(s\in(0,1/2)\) and \(t=2^{-k}\le1/8\), with errors tending to
zero. On a subsequence these fractions tend to some \(f'\), these measures tend weakly to
some \(\lambda\), and these forms tend to some \(H\) on all continuous pairs. Indeed \(D\) is
compact, and (5.1) permits diagonal extraction on a countable dense subset of \(C(D)\), followed
by extension in the uniform norm. The limiting bound (5.1), for real and complex symbols, extends
\(H\) to \(L^2(\lambda)\). Passing to this subsequence in (5.2) at fixed \(t\), then successively
comparing coefficients as \(t\downarrow0\), identifies \(f'=f_\mathrm{cell}\) and gives
\[
\lambda(m_s)=\mu_\mathrm{cell}(m_s),\qquad H(m_s,m_s)=K_\mathrm{cell}(m_s,m_s).
\]

Both sides are holomorphic in \(s\) on \(|\Re s|<1\). Indeed \(|1+2sn_1n_2|\ge1-|\Re s|\) there,
so \(s\mapsto m_s\) is holomorphic into \(C(D)\) with \(\|m_s\|_\infty\le |v|^2/(1-|\Re s|)\), and
(5.1) bounds the forms. The identities extend from the real tests to \(s=iT\). Let \(g,h\)
indicate the individual directions \([e_1],[e_2]\). As \(T\to\infty\), \(m_{iT}\) tends boundedly
to \(v_1^2g+v_2^2h\), hence in \(L^2(\lambda+\mu_\mathrm{cell})\). The three probes therefore give
the quadratic forms at \(g,h,g+h\). Polarization of the symmetric forms \(H,K_\mathrm{cell}\)
yields
\[
\lambda(g)=\mu_\mathrm{cell}(g),\qquad \lambda(h)=\mu_\mathrm{cell}(h),\qquad
H(g,h)=K_\mathrm{cell}(g,h). \tag{5.3}
\]
\emph{A separating inequality.}

For any two disjoint direction sets with indicators \(a,b\), put \(Y=25T_a\chi,\ Z=25T_b\chi\) and
\[
r_0=Y-3Z+3,\qquad r_1=Y+3Z-6,\qquad
R=\mathbb E[(1-\chi)r_0^2+\chi r_1^2].
\]

The two transforms have zero means and zero mutual inner product. Using
\(r_1^2-r_0^2=12YZ-18Y-18Z+27\) gives
\[
R=9+27f+175\mu(a)+5175\mu(b)+7500K(a,b). \tag{5.4}
\]

Checkerboard. On the unit torus take \(u=\mathbf1_{\{x_1<2/5\}},\ w=\mathbf1_{\{x_2<1/5\}},\
\chi=u+w-2uw\). Projecting onto the two axes gives the centered column and row averages, so
\(Y=15u-6,\ Z=5w-1\). Thus \(r_0=15(u-w)\) vanishes where \(\chi=0\), and \(r_1=15(u+w-1)\)
vanishes where \(\chi=1\). In particular \(R=0\). The four entries \((f,\mu(g),\mu(h),K(g,h))\)
are \((11/25,54/625,4/625,-144/15625)\).

Laminate. At each node let \(F\) be its conditional phase fraction. Starting from
\((X,Y,Z)=(25f-11,0,0)\), add \(25(F_\mathrm{child}-F_\mathrm{parent})\) in coordinate \(Y\),
\(Z\) or \(X\) according as the normal belongs to the first set, the second set or their
complement. Each parent is the barycenter of its children on one coordinate line. At a leaf
\(Y,Z\) are the transforms in (5.4) and \(\ell=X+Y+Z=25\chi-11\).

The polynomial
\[
\begin{aligned}
P&=r_0r_1-2r_0(\ell-14)+r_1(\ell+11)\\
 &=-XY+9XZ+8YZ-12X+24Y-36Z
\end{aligned}
\]
is affine along each coordinate. To bound its leaf average, put \(p=P/10\),
\((d_X,d_Y,d_Z)=(3/2,3/4,1/2)\), and \(\rho^2=(X/d_X)^2+(Y/d_Y)^2+(Z/d_Z)^2\). Use \(2|p|\) as a
comparison function, raising it to \(\max(2|p|,1+p^2-\rho^2)\) inside \(\{|p|<1,\rho<1\}\). This
function is separately convex. Indeed the three derivatives of \(P\) are
\[
-12-Y+9Z,\qquad 24-X+8Z,\qquad -36+9X+8Y .
\]

For \(\rho<1\), Cauchy-Schwarz bounds their variable parts in absolute value by \(5,5,15\),
respectively. Each derivative is constant along its own coordinate. Hence \(d_i|\partial_iP|>10\)
on every coordinate line that meets \(\{\rho<1\}\), making \(1+p^2-\rho^2\) convex along such a
line. On the region boundary \(1+p^2-\rho^2\le2|p|\). On any coordinate segment the chord between
the comparison values at its two ends dominates \(2|p|\) by convexity. The part of the segment
inside the indicated region is an interval. An end of this interval is either an end of the
segment lying in the region, where the chord equals the comparison value, or a point of the region
boundary, where the chord is at least \(2|p|\). In both cases the chord is at least
\(1+p^2-\rho^2\) at that end, and hence throughout the interval.

At the root \(|p|=6|X|/5\ge\rho\), so the comparison function is at least \(1\). At a leaf
\(11\le|\ell|\le(7/4)\rho\), so its value is \(|P|/5\). Iterating Jensen's inequality gives
\(\mathbb E|P|\ge5\). On the two leaf phases \(P=r_0(r_1+50)\) and \(P=r_1(r_0+25)\),
respectively. Cauchy-Schwarz now gives
\[
25\le(\mathbb E|P|)^2
\le R\,\mathbb E[(r_1+50)^2+(r_0+25)^2]
=R(44^2+28^2+2\mathbb E Y^2+18\mathbb E Z^2).
\]

Since \(\mathbb E[Y^2+Z^2]\le625/4\), the last factor is at most \(11065/2<75^2\). Consequently
\(R\ge1/225\) for every finite tree, every phase fraction and every disjoint pair of direction
sets.

\emph{Passage to the closure.}

Take continuous \(0\le a,b\le1\) with \(ab=0\). Apply the tree bound to the disjoint sets
\(\{a>\theta\},\{b>\eta\}\) and integrate over \((\theta,\eta)\in[0,1]^2\). By linearity of
\(\mu\) and bilinearity of \(K\), the right side of (5.4) with these continuous \(a,b\) is still
at least \(1/225\). The inequality is affine in the data, so it passes first to convex combinations and then to the limiting data \(f',\lambda,H\). Now
take disjoint continuous tent functions equal to one at the respective axes, with supports
shrinking to those axes. They converge to \(g,h\) in \(L^2(\lambda)\). By (5.1) and (5.3), every
periodic cell in the closure must therefore satisfy \(R(g,h)\ge1/225\) on its two axes. For the
four rectangles this would require \(0\ge1/225\).
The argument used only the scalar probe energies in (5.2). A basic neighborhood in their product topology involves a finite list. In that finite-dimensional space the checkerboard point is outside the closed convex hull, and strict separation by a linear functional follows by projection onto the hull.
\phantomsection\label{proofend:p5:obstruction}\end{proof}

\subsection*{Separation at finite precision}
\phantomsection\label{item:p5-precision}
There is a quantitative distinction between strict separation and a
test with specified error. Normalize both phase matrices so that their
eigenvalues lie in $[0.1,1]$. For every unit $e$, measure $e^TA_*e$.
There is one 255-split laminate whose reading differs from that of the
checkerboard by less than $0.000360$ throughout this class. Its fraction
is approximately $0.4400052$. A second tree, with the fraction corrected
to $11/25$, gives the outward bound $0.000434$.

In either case the coordinatewise midpoint of the predictions lies
within $0.001$ of both for every reading. Thus every finite panel admits
overlapping error boxes of half-width $0.001$. The bounds follow from
primal and complementary checkerboard trials, the exact laminate tree
and a uniform polynomial enclosure on the set of phase matrices.
Supplement~\ref{sup:planar} specifies the two certificates. This
estimate leaves the sharp noise threshold open.

\section{1 and 2. Two wells and coupled fields in the plane}\label{problem-1}
\Milton{1}
\Milton{2}

Cao first answered both questions in the negative. His Theorem~1.1 uses 16 potentials at every prescribed $0<f<1$, and his Theorem~6.1 uses two homogeneous quadratic wells on $\mathbb R^{2\times32}$ \cite{cao}. The checkerboard in Problem~5 gives an explicit two-potential example and a separate construction which allows all competing phase fractions.

\subsection*{The first field count at which lamination fails}

For a fixed pair of phases, the number of scalar potentials matters. In
Problem~5 one scalar potential is measured for many choices of the two
conductivities. Here the two constitutive tensors are fixed once and the
entire effective tensor is observed.

\begin{theorem}[Sharp field count]\label{new:threshold}
With one scalar potential, each effective tensor of two real symmetric
positive-definite planar phases is obtained by a laminate of the same
fraction and rank at most three. With two coupled potentials there is a
fixed phase pair for which the four-rectangle cell of Problem~5 at fraction
$11/25$ is outside the closed laminate set at that fraction. Moreover two
positive-definite homogeneous quadratic forms on $\mathbb R^{2\times4}$
have a minimum $W$ with $QW(\bar F)<RW(\bar F)$ at some $\bar F$.
\end{theorem}

\begin{proof}\phantomsection\label{proof:new:threshold}
Write a two-potential gradient in column order and let $U$ select its first
and fourth coordinates. Put $M=10^{40}$, $t=10^{-9}$ and
\[
 L_0=\operatorname{diag}(1,M,M,1),\qquad
 L_1=L_0+tU\begin{pmatrix}0&1\\1&0\end{pmatrix}U^T .
\]
The smallest eigenvalue of $L_1$ is $1-t$. For an effective tensor $C$
at fraction $s$, put $H=U^T(C-L_0)U$ and measure
\[
 \Phi(C,s)=9+27s-\frac{175H_{22}+5175H_{11}}{t^2}
                   +7500\,\frac{H_{12}-ts}{t^3}.
\]
Then $\Phi(C_{\square},11/25)<10^{-5}$, whereas
$\Phi(C_{\rm lam},s)>1/250$ for every finite laminate with its own
fraction $s$. At $s=11/25$ the gap persists in the closure.
To see how two fields encode this test, set
\[
 p(n)=\frac{n_x^2}{n_x^2+M n_y^2},\qquad
 r(n)=\frac{n_y^2}{M n_x^2+n_y^2}.
\]
The quadratic terms of $H$ measure the directional moments
$\mu(p),\mu(r)$ from the checkerboard proof. The off-diagonal cubic term
measures their mixed moment $K(p,r)$. The remainder changes $\Phi$ by
less than $8\cdot10^{-6}$. Since $p+r\le1$, randomly assigning each
tree split to the first class with probability $p$, to the second with
probability $r$, and otherwise to neither transfers the tree residual
bound $1/225$ to this statistic. On the checkerboard, as $M$ grows, the
two multipliers tend to the horizontal and vertical projections. Their
error in the statistic is below $8\cdot10^{-7}$ at the displayed $M$.
The expansion and all constants are derived in Supplement~\ref{sup:planar}.

For the scalar assertion let $A,B$ be the phase matrices, $f$ the
fraction of $B$, $g=1-f$, and $\Delta=B-A$. A stratified cell is already
a simple laminate. Otherwise eliminating the phase means in the cell
energy gives a positive-definite $2\times2$ matrix $Y$ and
\[
 C(Y)=gA+fB-fg\Delta(fA+gB+Y)^{-1}\Delta,\qquad
 \min(\det A,\det B)\le\det Y\le\max(\det A,\det B).
\]
The determinant bounds follow by applying planar duality to the constrained
energy defining $Y$. Each determinant endpoint is realized by a rank-two
sequential laminate. On a line $Y+svv^T$ meeting both determinant levels,
the inverse formula shows that the endpoint effective tensors differ by a
symmetric rank-one matrix. Layering them with normal in its kernel gives
every intermediate arithmetic mean, including $C(Y)$, at the same fraction.
The limiting argument for a singular $\Delta$ and the compactness of the
rank-three class are given in Supplement~\ref{sup:planar}.

For the wells, write the indefinite linear part of $\Phi$ as a difference
of two quadratic measurements. Planar duality converts the negative
measurement into a complementary energy on the same geometry. The fraction
term becomes a constant phase cost $c_i$. Both combined quadratic phase
laws have a common eigendirection $e$ of eigenvalue $M$. Adding
$\varepsilon^2c_i(e\cdot Z)^2$ to each homogeneous form and translating
$Z$ by $e/\varepsilon$ recovers the phase costs up to a relative
$O(\varepsilon)$ energy error, after subtracting a common affine function.
The strict gap survives. This uses four potential columns; minimality
of that number for the two-well assertion remains open.
\phantomsection\label{proofend:new:threshold}\end{proof}

For one scalar potential, the constructions for well-ordered phases and
simultaneously diagonal phases precede this argument \cite{grab,nsi}.
A change of spatial coordinates simultaneously diagonalizes two positive
matrices by congruence. The proof above gives the rank bound directly,
including unordered phases. At two potentials the fraction term in $\Phi$
matters. For a block-diagonal system with an unspecified number of
potentials the following statement gives exclusion even if a competing
laminate changes its fraction.

\begin{theorem}\label{p1:coupled}
There exist a finite $n$ and two block-diagonal positive-definite phase laws $L_0,L_1$ on $\mathbb R^{2\times n}$ with the following property. In the unit square put $L_1$ where exactly one of $x_1<2/5$ and $x_2<1/5$ holds, and $L_0$ elsewhere. This periodic cell has fraction $11/25$ of $L_1$, and its effective tensor lies outside the closure of finite laminates of those phases over all fractions. This block-diagonal form of the all-fractions assertion does not specify the component count. For some finite $m$ there are two positive-definite homogeneous quadratic forms on $\mathbb R^{2\times m}$ whose minimum $W$ satisfies $QW(F)<RW(F)$ at a mean gradient $F$.
\end{theorem}
\begin{proof}\phantomsection\label{proof:p1:coupled}
By Theorem~\ref{p5:obstruction} there are finitely many scalar tests, each with its own pair of positive phase matrices and mean field $e_i$, and real numbers $c_i\ne0$, such that
\[
  \sum_i c_i e_i\cdot B_i e_i + \delta
    \le \sum_i c_i e_i\cdot A_i e_i\qquad (\delta>0).
\]
Here $B_i$ are the effective matrices of the checkerboard and $A_i$ those of any one laminate tree in the same tests, including trees of different fractions. Use the two block-diagonal phases formed from these tests. The independent scalar potentials share the same geometry, so the diagonal blocks of their effective law are the measured matrices. The strict gap proves the first assertion.

For the wells, retain each positive-weight test and replace each negative-weight test by its planar dual. If $D$ is a phase matrix, its dual is $JD^{-1}J^T$, where $J$ is a quarter-turn. Rotating flux and gradient gives a compatible gradient and an equilibrated flux for the dual cell. At a layer interface the rotation interchanges the continuous normal flux and tangential field conditions. Thus this replacement sends every effective matrix $A$ to $JA^{-1}J^T$, for the same geometry. In a positive block prescribe $e_i$ and in a negative block prescribe $JB_i e_i$. Weight its energy by $|c_i|$. Let $q_0,q_1$ be the sums of the resulting phase energies, $F$ the matrix of these mean columns, and $W=\min(q_0,q_1)$. For any positive $A$ and $p=B_i e_i$,
\[
 p\cdot A^{-1}p\ge 2 e_i\cdot p-e_i\cdot A e_i ,
\]
with equality at $A=B_i$. Thus the total laminate energy at $F$ exceeds the corresponding checkerboard energy by at least $\delta$.

A single leaf is itself a tree. Let $V(G)$ be the infimum of the total effective quadratic energy at $G$ over all finite trees. Then $V\le W$. For $G-H=n\otimes a$, laminate two near-minimizing trees normal to $n$. The two mean gradients $G,H$ are admissible ply fields at the mixed mean, showing rank-one convexity of $V$ on letting the errors vanish. Thus $RW\ge V$. The independent periodic correctors on the checkerboard give a trial for $QW(F)$ with energy at most its total cell energy, since $W$ is bounded above by the energy of the selected phase pointwise. Periodic mollification gives the same bound with smooth trials by quadratic growth. The two inequalities give the asserted strict gap.
\phantomsection\label{proofend:p1:coupled}\end{proof}

\section{3. The entire three-dimensional isotropic conductivity function}\label{problem-3}

\Milton{3}
The vertices of a cubic resistor network can be colored by two phases. The weak-contrast curvature of the network measures the density of edges across a change of color. The criterion below characterizes the full fixed-fraction closure. Beyond the two explicit derivative identities, the constraints on a proposed function take the form of a limit of finite network tests.

Fix \(0\le f\le1\), let \(\Omega=\mathbb C\setminus(-\infty,0]\), and put \(J_0=[1,2]\). On the unit
three-torus, a binary indicator \(\chi\) of mean \(f\) has coefficient \(a_\chi(z)=1+(z-1)\chi\).
Its effective tensor \(A_\chi(z)\) maps an imposed mean electric field to the mean current. Write
\(\mathcal C_f\) for the functions \(s\) for which one sequence of such cells of exact fraction
\(f\) satisfies \(A_{\chi_n}\to sI_3\) locally uniformly on \(\Omega\). The class is the same if
every approximating tensor is required to be scalar. Exact attainment by a single cell would
require the stronger identity \(A_\chi(z)=s(z)I_3\) for every \(z\in\Omega\).

For \(m=2^k\), \(k\ge1\), put \(T_m=(\mathbb Z/m)^3\) and \(N=m^3\). Use every coloring
\(P:T_m\to\{0,1\}\) invariant under the 48 signed coordinate permutations. With one labeled edge
\((\ell,i)\) for each positive coordinate step, set
\[
\theta_{\ell i}=\frac{P_\ell+P_{\ell+e_i}}2,\qquad
 c_{\ell i}(z)=1+(z-1)\theta_{\ell i},\qquad b_{\ell i}=\delta_{i1}.
\]

The graph response for real \(x>0\) is
\[
Y_P(x)=\frac1N\min_v\sum_{\ell,i}c_{\ell i}(x)
           (b_{\ell i}+v_{\ell+e_i}-v_\ell)^2 .
\]

Parallel edges keep their labels when \(m=2\). The minimum uses periodic vertex potentials. Pinning one vertex gives a rational
computation at each sample \(x_j=1+j/m\). Use one coloring in the maximum
\[
d_m(s)=\min_{P\text{ cubic on }T_m}\ \max_{0\le j\le m}|s(x_j)-Y_P(x_j)|.
\]

\begin{theorem}\label{p3:network}
For $0\le f\le1$ and $s$ holomorphic on $\Omega$ and real on its positive axis,
\[
s\in\mathcal C_f\quad\Longleftrightarrow\quad
 s'(1)=f,\qquad s''(1)=-\frac23 f(1-f),\qquad
 \lim_{k\to\infty}d_{2^k}(s)=0.
\]

\end{theorem}
The test involves an unbounded sequence of finite graphs. Whether Milton's explicit function constraints \cite[Eqs.~(3.8), (3.9), (3.11), (3.12)]{milton} already suffice remains open.
\begin{proof}\phantomsection\label{proof:p3:network}
The sample at \(x=1\) gives \(s(1)=1\). All fractions use the same finite graph libraries. The
endpoints give just \(s=1\) for \(f=0\) and \(s(z)=z\) for \(f=1\).

\emph{Compactness.}

Both the graph and continuum cell equations extend holomorphically to \(\Omega\), with bounds
uniform in the geometry and period on every compact \(K\subset\Omega\). Indeed, write \(F=b+q\)
with \(q\) a periodic gradient and require \(cF\) to be orthogonal to the gradients, using
\(a_\chi\) in the continuum. For \(0\le\theta\le1\),
\[
\Re\!\left(e^{-i\arg(z)/2}[1+(z-1)\theta]\right)
 \ge \min(1,|z|)\cos(\arg(z)/2).
\]

The right side has a positive minimum on \(K\). Testing with \(q\) gives a uniform \(L^2\) bound
on the field and hence on the mean current. Coercivity on the gradient space and the analytic
inverse of its affine operator give holomorphy. Covariance under coordinate sign changes kills the
off-diagonal tensor entries, and coordinate permutations equalize the diagonals. This gives a
scalar tensor for every cubic graph and cell, including at complex contrasts.

On \(J_0\), the real graph responses are nondecreasing and concave, as minima of affine increasing
energies. The constant field gives \(Y_P(x)\le1+p(x-1)\) where \(p=N^{-1}\sum_\ell P_\ell\), and
\(Y_P(1)=1\). Thus they are \(1\)-Lipschitz on \(J_0\). The same argument gives Lipschitz constant
\(f\) for physical scalar responses. Let \(\mathcal L\) be the locally uniform closure of the
union of the dyadic graph libraries. Periodically repeating \(P\) repeats its Kirchhoff solution,
so the libraries are nested. The uniform bounds give compactness by Montel's theorem, both for
\(\mathcal L\) and for the scalar part of the closed fixed-fraction continuum tensor family.
Replacing the sampled maximum by \(\|s-Y_P\|_{J_0}\) costs at most \(\omega_s(1/(2m))+1/(2m)\), where
\(\omega_s\) is the modulus of continuity on \(J_0\). Consequently
\[
\lim_{k\to\infty} d_{2^k}(s)=\min_{Y\in\mathcal L}\|s-Y\|_{J_0} .
\]

If this vanishes, any sequence of minimizing colorings, or colorings with discrepancy
\(d_m(s)+o(1)\), converges in response to \(s\) locally uniformly on \(\Omega\). Every normal
subsequential limit agrees with \(s\) on \(J_0\), and the identity theorem identifies it throughout
\(\Omega\). Cauchy's formula now gives convergence of the derivatives at one.

\emph{Counting interfaces by curvature.}

Give edge fields the inner product \(N^{-1}\sum_e\overline{u_e}v_e\). Let \(\Pi\) project onto
periodic graph gradients and let \(\Theta=\operatorname{diag}(\theta_e)\). For unit loading
\(b^{(r)}_{\ell i}=\delta_{ir}\), the Kirchhoff field obeys \(F=b^{(r)}-h\Pi\Theta F\) at
\(z=1+h\). Expanding the mean flux gives
\[
Y'_P(1)=p,\qquad -\tfrac12Y_P''(1)=\|\Pi\Theta b^{(r)}\|_2^2 .
\]

The latter is the same for all three loads. Define
\[
\Delta_P=p(1-p)+\tfrac32 Y_P''(1),\qquad
 \beta_P=N^{-1}\sum_{\ell,i}|P_{\ell+e_i}-P_\ell|^2.
\]

Thus \(N\beta_P\) counts unlike bonds. At a nonzero Fourier mode \(\xi\in T_m\), put
\(t_i=\sin^2(\pi\xi_i/m)\) and \(T=\sum_i t_i\). The gradient direction has components \(e^{2\pi
i\xi_i/m}-1\), so its projection has diagonal \(t_i/T\). Averaging the two ends of an \(i\)-edge
contributes the squared multiplier \(1-t_i\). With the Fourier transform normalized by \(1/N\),
Parseval gives \(\sum_{\xi\ne0}|\widehat P(\xi)|^2=p(1-p)\). Summing the three projected norms and
also transforming the three differences gives
\[
\Delta_P=\sum_{\xi\ne0}|\widehat P(\xi)|^2\frac{\sum_i t_i^2}{T},
 \qquad
 \beta_P=4\sum_{\xi\ne0}|\widehat P(\xi)|^2T,\qquad
 \ \beta_P/12\le\Delta_P\le\beta_P/4\  .
\]

The bounds follow from \(T^2/3\le\sum_i t_i^2\le T^2\). If the discrepancies tend to zero and the two derivative tests hold, the minimizing graphs therefore satisfy
\(p\to f\) and \(\beta_P\to0\). 

\emph{A continuum cell at the prescribed fraction.}

Center the periodic voxel \(V_\ell=(\ell+[-1/2,1/2)^3)/m\) at site \(\ell/m\) and give it color
\(P_\ell\). If \(p<f\), fill color-zero voxels in whole symmetry orbits until the remaining volume
fits in one orbit. If the last orbit has \(o\) voxels and volume \(v\) remains, insert a centered
cube of side \((v/o)^{1/3}\) in each. If \(p>f\), remove phase one by the same rule. If \(p=f\),
keep the voxel cell. At most \(N|p-f|+48\) voxels are touched. The resulting indicator
\(\chi_{P,f}\) is binary and invariant under signed coordinate permutations, has mean exactly
\(f\) even for irrational \(f\), and is chosen independently of \(z\). Denote its scalar response
by \(g\).

Fix \(x\in J_0\), use \(a=1+(x-1)\chi_{P,f}\), and let \(F_{\ell i}=b_{\ell i}+v_{\ell+e_i}-v_\ell\)
and \(J_{\ell i}=c_{\ell i}F_{\ell i}\) at the graph minimum. For edge \(e=(\ell,i)\) its stencil
consists of
\[
\{V_{\ell-\sigma}:\sigma\in\{0,1\}^3,\ \sigma_i=0\}\ \cup\ \{V_{\ell+e_i}\}.
\]

Mark \(e\) if a stencil voxel was touched or if the original colors in its stencil differ. Write
\(B\) for the marked edges. Each stencil is connected by nearest-neighbor steps and has at most
five voxels. A mixed stencil thus contains an unlike bond. Each voxel occurs in at most five
stencils per orientation, so a particular voxel or bond occurs in at most fifteen stencils. It
follows that
\[
|B|/N\le 15\big(\beta_P+|p-f|+48/N\big).
\]

Voltage trial. Place \(v_\ell/m\) at the lower corner of its site-centered voxel \(V_\ell\) and
interpolate trilinearly. This placement gives the five-voxel stencil above. In each voxel the \(i\)th field
component is a convex combination of the four parallel edge slopes, with transverse bilinear
weights. Integrate Jensen's inequality against \(a\). Edge \(e\) then has weight \(c_e^+/N\),
where \(c_e^+\) is \(N\) times the sum of its four integrals of \(a\) times the respective
bilinear weight. Those weights integrate to \(1/N\) in total before multiplication by \(a\). Hence
\(1\le c_e^+\le x\), with \(c_e^+=c_e\) on unmarked edges.

Current trial. In \(V_\ell\), with local coordinates \(0\le y_i\le1\), use
\[
j_i=(1-y_i)J_{\ell-e_i,i}+y_iJ_{\ell i}.
\]

Normal components match across faces and the divergence is \(m\sum_i(J_{\ell
i}-J_{\ell-e_i,i})=0\) by Kirchhoff's equations. Integration preserves the graph mean current. Its
loading-direction component is \(Y_P\) because \(J\) pairs to zero with the minimizing potential
difference, and the transverse components vanish by symmetry. Jensen's inequality for its
complementary energy assigns \(J_e^2\) a weight \(r_e^-/N\) from \(V_\ell\) and \(V_{\ell+e_i}\).
The two linear weights integrate to \(1/N\) together, so \(1/x\le r_e^-\le1\), with
\(r_e^-=1/c_e\) off \(B\).

Put \(Y=Y_P(x)\) and \(\delta_P(x)=(x-1)N^{-1}\sum_{e\in B}c_e F_e^2\). The voltage minimum at
mean \(e_1\) is \(g(x)\). The minimum of \(\int a^{-1}|j|^2\) over divergence-free currents of
mean \(Ye_1\) is \(Y^2/g(x)\), attained by the physical current for load \((Y/g(x))e_1\). Its
electric field pairs to zero with each solenoidal variation of mean zero. The ratios \(c_e^+/c_e\)
and \(r_e^- c_e\) are at most \(x\) on \(B\). Since \(N^{-1}\sum_e c_e F_e^2=Y\), the trials give
\[
\frac{Y^2}{Y+\delta_P(x)}\le g(x)\le Y+\delta_P(x),
 \qquad |Y-g(x)|\le\delta_P(x).
\]

To bound the energy on \(B\) uniformly in the coloring, use the edge norm
\(\|U\|_q=(N^{-1}\sum_e|U_e|^q)^{1/q}\). Write \(c_e=\bar c(1+\kappa\nu_e)\) with \(\bar
c=(1+x)/2\), \(\kappa=(x-1)/(x+1)\le1/3\) and \(\nu_e=2\theta_e-1\). Since \(|\nu_e|\le1\), the
only operator estimate needed is for the unweighted lattice. Indeed, \(\Pi(cF)=0\) and \(\Pi
F=F-b\) give
\[
F=b-\kappa\Pi(\nu F).
\]

For dyadic periods the periodic discrete Poisson solution operator is the orthogonal gradient projector $\Pi$. Caputo and Ioffe give its period-independent weak $(1,1)$ bound and the interpolation argument \cite[Eq.~(4.4), Appendix A]{ci}. Interpolating weak $(1,1)$ with the $L^2$ contraction gives a uniform strong $L^{4/3}$ bound. Self-adjointness gives the dual $L^4$ bound. The edge norm and the sitewise vector norm are equivalent at exponent four, uniformly in the period. Interpolation with the norm-one estimate at two now gives a fixed $2<q\le4$ with $\|\Pi\|_{q\to q}\le2$ for every dyadic period.

Now \(\|b\|_q=1\), so \(\|F\|_q\le(1-2\kappa)^{-1}\le3\). H\"{o}lder's inequality with
\(\alpha=1-2/q\) gives \(\delta_P(x)\le 2\|F\|_q^2(|B|/N)^\alpha\le18(|B|/N)^\alpha\). In
particular,
\[
\|Y_P-g_{\chi_{P,f}}\|_{J_0}
 \le 18\big[15(\beta_P+|p-f|+48/N)\big]^\alpha
 \le C'\big(\Delta_P+|p-f|+48/N\big)^\alpha .
\]

For the minimizing graphs above the right side tends to zero. The recovered cells thus converge to \(s\) on
\(J_0\) and, by normality and analytic uniqueness, locally uniformly on \(\Omega\). This proves sufficiency.

\emph{Necessity.}

First, the continuum sum rule applies even to periodic cells with anisotropic tensor. In the
expansion above replace \(\Pi\) by the periodic continuum gradient projection \(\Gamma\) and
\(\Theta\) by multiplication by \(\chi\). Then \(A'_\chi(1)=fI_3\). The nonzero-mode symbol of
\(\Gamma\) is \(kk^T/|k|^2\), of trace one, whence
\[
-\tfrac12\operatorname{tr} A''_\chi(1)
 =\sum_{i=1}^3\|\Gamma(\chi e_i)\|_2^2
 =\sum_{k\ne0}|\widehat\chi(k)|^2=f(1-f).
\]

Passing to any scalar limit gives the two prescribed derivatives.

\emph{Cubic density.} Fix a finite list \(x_j>0\) and \(s\in\mathcal C_f\). Choose one
exact-fraction periodic cell \(\chi\) with all \(A_j=A_\chi(x_j)\) within \(\varepsilon\) of
\(s(x_j)I_3\) in operator norm. In the chamber \(D=\{0<y_3<y_2<y_1<1/2\}\) of the centered unit
cube, insert \(\chi(Ly)\) in all whole period cubes of side \(1/L\) contained in \(D\). The
leftover set lies within \(\sqrt{3}/L\) of the boundary. Fill proportion \(f\) of this set with
phase one, for example by slicing it with a plane. Every whole period cube already has fraction
\(f\). Reflect and permute this construction to the other 47 chambers and periodize. The resulting
cell has exact fraction \(f\) and scalar response \(g_L\).

For each \(j\) these cells satisfy
\[
\lambda_{\min}(A_j)\le\liminf_L g_L(x_j)
 \le\limsup_L g_L(x_j)\le\lambda_{\max}(A_j).
\]

To see the upper bound, on chamber \(RD\) use a smoothed periodic voltage corrector of \(\chi\)
for load \(R^T e_1\), rotated back and repeated at period \(1/L\). For the lower bound use a
smoothed complementary minimizing current of that same mean. Its mean-zero part is the curl of a
periodic vector potential, by the Fourier representation of solenoidal fields. Smooth
approximation in \(L^2\) makes the two trial energies as close as desired to the respective
quadratic forms of \(A_j\) and \(A_j^{-1}\). Multiply each repeated potential by
\(L^{-1}\eta_\tau(y)\), where the cutoff is one outside a chamber boundary layer of width
\(\tau\), zero near its boundary and \(|\nabla\eta_\tau|\le C/\tau\). Take the gradient or curl
after the cutoff. Under an orthogonal map \(R\), a curl potential transforms by \((\det R)R\).
Thus the fluctuations paste to periodic gradient and curl fields of mean zero. For the fixed
smooth trials and \(L\tau\ge\sqrt3\) the fields are bounded, giving boundary energy \(O(\tau)\). On
each chamber interior outside the strip, periodic averaging gives its volume times the
corresponding period-averaged trial energy. Taking \(L\to\infty\), then \(\tau\downarrow0\) and
the smoothing error to zero bounds \(g_L\) by \(\lambda_{\max}(A_j)\) from above and \(1/g_L\) by
\(\lambda_{\max}(A_j^{-1})\) from above, in the stated limit sense. One large \(L\) works for the
entire finite list. Letting \(\varepsilon\downarrow0\) along longer prefixes of a list of distinct
nodes accumulating at one, normality gives cubic density for the full function.

\emph{Dyadic approximation.} For any fixed measurable cubic cell \(\chi\) of
mean \(f\), color each \(V_\ell\) by its majority phase, using the same color at every tie. Its
voxel indicator \(b_m\) is cubic and converges to \(\chi\) in \(L^1\). Indeed the voxel averages
converge in \(L^1\) by approximation with continuous functions and the vanishing voxel diameter,
and majority rounding costs at most twice the averaging error against a binary target. In
particular \(p_m\to f\). Translation continuity in \(L^1\) gives
\[
\beta_{P_m}=\sum_i\|b_m(\cdot+e_i/m)-b_m\|_1
 \le 6\|b_m-\chi\|_1+\sum_i\|\chi(\cdot+e_i/m)-\chi\|_1\ \longrightarrow\ 0 .
\]

The orbit adjustment changes volume \(|p_m-f|\), so
the exact-fraction indicators still converge to \(\chi\) in \(L^1\). At each fixed \(x>0\), test
the voltage and complementary minima of these cells with smooth trials approaching the minimizers
of \(\chi\). The bounded trial integrands can be integrated against the converging conductivities
and resistivities. This bounds the limsup of both \(g_m(x)\) and \(1/g_m(x)\) by \(g_\chi(x)\) and
\(1/g_\chi(x)\) respectively, giving convergence. The common Lipschitz bound makes it uniform on
\(J_0\). The displayed graph-to-cell estimate now shows \(Y_{P_m}\to g_\chi\) uniformly on \(J_0\)
along all dyadic levels. Cubic density and a triangle inequality give \(d_m(s)\to0\) for every
physical \(s\). Together the two directions identify \(\mathcal C_f\) with the subset of
\(\mathcal L\) satisfying the two derivative identities.
\phantomsection\label{proofend:p3:network}\end{proof}

\subsection*{One cell at all contrasts}

The graph criterion tests a limit of responses. A separate question asks
for one phase indicator whose conductivity is the prescribed function at
every contrast. The distinction already appears at every nontrivial
volume fraction.

\begin{theorem}[A closure law with no finite-perimeter cell]\label{new:perimeter}
Fix $0<f<1$, put $q=1-f$, and for real $z>0$ write
\[
 A=q+fz,\qquad H=(q+f/z)^{-1},\qquad
 L_f(z)=\frac{\sqrt{H^2+8HA}-H}{2}.
\]
The law $L_f$ belongs to $\mathcal C_f$. No single periodic cell of
finite phase perimeter has tensor $L_f(z)I$ for every positive $z$.
\end{theorem}
\begin{proof}\phantomsection\label{proof:new:perimeter}
Let $C_a^p(k)=a+3ap(k-a)/(3a+(1-p)(k-a))$ be the dipole
conductivity of a sphere with shell $a$ and core $k$ of relative volume
$p$. Retain share $1-\delta$ of a ball of phase fraction $f$, add a
$z$ shell of share $f\delta$ and a unit shell of share $q\delta$.
One step acts by
$M=C_1^{\,1-q\delta}\circ C_z^{(1-\delta)/(1-q\delta)}$.
Its positive fixed point attracts iteration. On the interval between
$1$ and $z$,
\[
 M(k)-k=\tfrac{\delta}{3}(2A-k-k^2/H)+O_z(\delta^2).
\]
The fixed points therefore tend to $L_f$. Coated-ball packings and a
diagonal choice on finite contrast sets give one exact-fraction
approximating sequence for the whole law.

The identity $L_f(1/z)=L_f(z)/(AH)$ gives
\[
 s(z)s(1/z)+\frac{s(z)+z\,s(1/z)}{1+z}=2 \qquad(s=L_f).
\]
For an exact finite-perimeter cell with this reciprocal equality there is
one indecomposable component of one phase, of volume $m>0$, with
conductivity $a(z)=1$ or $a(z)=z$, such that
\[
 m\bigl(2+s(z)/a(z)\bigr)^2\le9s(z)/a(z) \qquad (z>0).
\]
The component can be chosen before the contrast. There are only countably many positive-volume components. Their mean field traces are analytic in $z$. If each component identity failed identically, one contrast near one would avoid the discrete zero sets of them all. The periodic gluing argument in Supplement~\ref{sup:cells} rules this out. Supplement~\ref{sup:cells} derives the component estimate from the trace
of the harmonic-coordinate fields. It forces $s/a$ to stay bounded above
and away from zero. But $L_f(z)\sim\sqrt{2f/q}\,\sqrt z$ as
$z\to\infty$, contradicting either choice of phase.
\phantomsection\label{proofend:new:perimeter}\end{proof}

A different law excludes a class with arbitrarily rough stratification.
Set $\tau=q/2$ and
$F_\tau(1+h)=1+fh-fq h^2/[3(1+\tau h)]$. This is the isotropized
coated-cylinder law. It has no exact cell whose phase, off a closed set
of zero $W^{1,2}$ capacity, is locally a measurable function of a smooth
coordinate with nonvanishing gradient. In each regular chart the normal
flux and tangential field have a common Sobolev decomposition. The
resulting sign condition is incompatible with the cylinder pole. See
Supplement~\ref{sup:cells}. Exact scalar realization by a completely
arbitrary measurable periodic cell remains open.

An earlier response-only hierarchy for this same exact-fraction closure
uses Fourier templates, binary rounding and positive matrices
\cite[Theorem 1.1]{eve}. Theorem~\ref{p3:network} instead tests cubic
resistor networks and uses their weak-contrast curvature to control the
density of unlike neighbouring colors.

\subsection*{Three readings from one geometry}
\begin{proposition}\label{new:three}
At $f=1/2$, consider the following proposed readings. Each pair is
jointly attainable in closure. The least possible maximum absolute error
on all three readings is exactly $1/200000$.
\[
\begin{array}{c|ccc}
 z&2&3&4\\ \hline
 y&120800419/83800000&362399/200000&30400071/14200000\\
 s_*(z)&604/419&453/250&152/71
\end{array}
\]
\end{proposition}
\begin{proof}\phantomsection\label{proof:new:three}
Pack identical, rescaled four-layer balls with successive material
conductivities $1,z,1,z$ and relative volumes $1/20,1/4,9/20,1/4$.
Their dipole law gives exactly the last row. The three errors have
alternating signs. The conductivity representation gives a fractional
linear transform of the readings which is concave in a common contrast
variable. This particular ball gives an affine transformed law. An error
strictly below $1/200000$ would move the two outer transformed values up
and the middle value down, contradicting concavity. Supplement~\ref{sup:cells}
gives the transform, the two-reading constructions and the general
two-sided approximation procedure for rational input.
\phantomsection\label{proofend:new:three}\end{proof}

\section{4. Three lamination directions in space}\label{problem-4}

\Milton{4}
Three orthogonal directions do not suffice. The question of completeness with unrestricted directions remains open.
\begin{theorem}\label{p4:directions}
Let $X_9$ be the nine-split laminate of $zI$ and $I$ in the table below. Three further layering steps make its response $s_X(z)I$, where $s_X=\operatorname{tr}(K_{X_9})/3$. No sequence of trees using at most three directions per tree can converge to this tensor at every real $z>0$, even when the directions and fractions vary with the tree.
\end{theorem}
\begin{proof}\phantomsection\label{proof:p4:directions}
For a normalized trace \(s(z)=\operatorname{tr}K(z)/3\) define
\[
f=s'(1),\qquad m=(1-f)/3,\qquad
 v=\frac{s'''(1)}{6f}-m^2,\qquad a(z)=\frac{s(z)-1}{f(z-1)}.
\]

If all normals lie on three fixed orthogonal axes and $0<f<1$, then
\begin{equation}\label{p4:covariance}
\left(\frac{1-a(z)}{z-1}-ma(z)\right)^2
       \le \frac{2}{3}v\frac{s'(z)}{f}\quad(z>0,\ z\ne1).                    
\end{equation}

Choose the orthogonal directions as coordinate
axes. On a tree leaf the three loadings form a diagonal field with positive entries \(L_i(z)\).
Across a split the tangential components are transmitted unchanged, and the normal component is
multiplied by the positive ratio of the parent and child normal conductivities. At \(z=1\) write
\(W_i=-L_i'(1)\). The first derivative of every subtree tensor is its phase fraction times \(I\).
At a split the weak field therefore increases in its normal coordinate by
\(f_\mathrm{child}-f_\mathrm{parent}\). Telescoping gives \(\sum W_i=1-f=3m\) on each phase-\(z\)
leaf.

Let \(\mathbb E_1\) take the phase-conditional volume mean and the uniform mean over the three
coordinates. Flux averaging gives \(\mathbb E_1 L_i(z)=a(z)\). Pairing the flux at \(z\) with the
gradient at \(w\) and interchanging the two contrasts gives
\[
\mathbb E_1 L_i(z)L_i(w)=\frac{s(z)-s(w)}{f(z-w)} .
\]

The identity holds because all mean-zero field fluctuations are orthogonal to the flux
fluctuations. The telescoping sum gives \(\mathbb E_1 W_i=m\). Differentiation at one and the diagonal limit then give \(\mathbb
E_1(W_i-m)^2=v\), \(\mathbb E_1 W_iL_i(z)=(1-a(z))/(z-1)\) and \(\mathbb E_1L_i(z)^2=s'(z)/f\).
On every leaf we may subtract the
coordinate mean \(\bar L\) from \(L_i\) in its pairing with \(W_i-m\). Positivity implies
\[
\sum_i(L_i-\bar L)^2=\sum_iL_i^2-\tfrac13(\sum_i L_i)^2
                   \le\tfrac23\sum_i L_i^2 .
\]

Cauchy-Schwarz proves \eqref{p4:covariance}. The conductivity functions form a normal family on the complex plane
slit along the nonpositive real axis, by the rotated coercivity estimate in the proof of Theorem~\ref{p3:network}, applied to the finite-dimensional space of compatible tree fields.
Pointwise convergence on all positive contrasts thus gives local uniform convergence and
convergence of the derivatives in \eqref{p4:covariance}. This proves the obstruction for closure. 

For a stack of tensors \(U,V\), write \(L(U,V;p,n)\) for the tensor with first share \(p\) and
layer normal \(n\). Continuity of tangential field and normal flux gives the rational recursion
\[
L=pU+(1-p)V-p(1-p)\frac{[(U-V)n][(U-V)n]^T}{n^T[(1-p)U+pV]n}.
\]

In the following construction \(z\) and \(1\) stand for the corresponding isotropic tensors.

\par\addvspace{4pt}\begingroup\small\renewcommand{\arraystretch}{1.20}\setlength{\tabcolsep}{4pt}
\setlength{\tableusable}{\dimexpr\linewidth-32pt\relax}
\noindent\begin{tabular}{@{}>{\raggedright\arraybackslash}p{\dimexpr 0.18\tableusable\relax}>{\raggedright\arraybackslash}p{\dimexpr 0.2\tableusable\relax}>{\raggedright\arraybackslash}p{\dimexpr 0.2\tableusable\relax}>{\raggedright\arraybackslash}p{\dimexpr 0.19\tableusable\relax}>{\raggedright\arraybackslash}p{\dimexpr 0.23\tableusable\relax}@{}}\toprule
\bfseries Material & \bfseries Left child & \bfseries Right child & \bfseries Left share & \bfseries Layer normal \\\midrule
\(P_1\) & \(z\) & \(1\) & \(87/88\) & \((0,1,0)\) \\
\(P_2\) & \(z\) & \(1\) & \(134/135\) & \((-3,0,4)\) \\
\(Q\) & \(P_1\) & \(P_2\) & \(9/20\) & \((4,0,3)\) \\
\(A\) & \(1\) & \(Q\) & \(1/42\) & \((1,0,25)\) \\
\(P_3\) & \(z\) & \(1\) & \(88/89\) & \((0,1,0)\) \\
\(P_4\) & \(1\) & \(z\) & \(1/147\) & \((10,0,23)\) \\
\(B\) & \(P_3\) & \(P_4\) & \(21/32\) & \((23,0,-10)\) \\
\(M\) & \(A\) & \(B\) & \(3/7\) & \((-1,0,5)\) \\
\(X_9\) & \(M\) & \(1\) & \(79/80\) & \((1,0,0)\) \\
\bottomrule\end{tabular}\par\endgroup\addvspace{4pt}

Put \(X=X_9\). Reflection in the \(x\)-\(z\) plane preserves every layer direction of \(X\) up to
sign. Thus \(e_y\) is an eigenvector of \(K_X(z)\) at every contrast. Mix \(X\) equally with its
quarter-turn about \(e_y\), in layers normal to \(e_y\). Both tensors give the same flux for this
normal loading, so their effective tensor equals their arithmetic average. The quarter-turn
averages the transverse block to half its trace times the transverse identity. In coordinate order
\(y,x,z\) the result has diagonal \((a,b,b)\).

Mix this result in share \(1/3\) with its exchange of the first two coordinates, normal to the
third coordinate. The diagonal becomes \(((a+2b)/3,(2a+b)/3,b)\). Mix it equally with its exchange
of the last two coordinates, normal to the first. Again the identical normal columns give the
arithmetic mean. The resulting tensor is
\[
K_\mathrm{iso}(z)=\frac{\operatorname{tr}K_X(z)}{3} I
\]

Denote by $R(z)$ the left side of \eqref{p4:covariance} divided by $v s'(z)/f$ when this denominator is positive. Insert the nine rows in the rational recursion, take one third of the trace, and differentiate this rational function at $1$ and at $10^6$. The exact arithmetic required here is
\[
 f=1123548936049/1160514432000,\quad
 v s'(10^6)/f>0,\qquad 0.69196432<R(10^6)<0.69196434 .
\]
The last comparison violates \eqref{p4:covariance}. All three are asserted in \path{certificates/baseline/p04.py}.

\emph{An arbitrary triple of directions.}
Give the leaves their volume probabilities. For each internal node $j$ let $H_j$ project onto the normalized scalar function supported on that node, constant on each child, with mean zero on the node. These Haar functions form an orthonormal basis of all mean-zero leaf functions. A compatible gradient fluctuation can jump between the child means only parallel to the node normal $n_j$. Its orthogonal projector on leaf vector fields is thus
$\Pi=\sum_j H_j\otimes n_j n_j^T$ (unit normals). Let $P$ be multiplication by the phase-one indicator. For a constant loading $e$ flux orthogonality on all compatible jumps gives $F=e-(z-1)\Pi PF$. Consequently
\[
 e^TK(z)e=|e|^2+(z-1)\langle Pe,[I+(z-1)P\Pi P]^{-1}Pe\rangle .
\]
Expanding at one gives
\[
 -K''(1)/2=\sum_j V_j p_j(1-p_j)(f_{j,L}-f_{j,R})^2 n_j n_j^T .
\]
Indeed the coefficient of $n_j n_j^T$ is $\|H_j\chi\|^2$, where $V_j,p_j$ are node volume and left share, and $f_{j,L},f_{j,R}$ are the child fractions. The weights sum to $f(1-f)$ by Parseval. In an isotropic limit with $0<f<1$, normalize by this sum. Compactness gives limits $n_i$ and $w_i$ for the three directions and their total weights. The limit matrix is $I/3$, by its trace and isotropy. Thus the square matrix with columns $\sqrt{3w_i} n_i$ has product with its transpose $I$, so the three limiting directions are orthogonal.

Replace the three normals of each tree by the limiting directions. Since the $H_j$ have mutually orthogonal ranges, $\|\Pi-\widetilde\Pi\|$ is bounded by the largest norm of the three differences of normal projectors, tending to zero regardless of tree size. The positive contractions $P\Pi P$ and $P\widetilde\Pi P$ have resolvent norms at most $1/\min(1,z)$ for $z>0$. The resolvent identity in the displayed energy formula thus bounds the change on a unit loading by $f(z-1)^2\|\Pi-\widetilde\Pi\|/\min(1,z)^2$. Any such approximation would give the same limit with three orthogonal directions, contradicting \eqref{p4:covariance}.
\phantomsection\label{proofend:p4:directions}\end{proof}

\subsection*{A test for unrestricted laminate normals}
\phantomsection\label{item:p4-score}
For a half-volume isotropic specimen with conductivity $s(z)$ define
\[
 \mathcal S(s)=\frac{841}{16900}\frac{s(25)}{25}
                +\frac{10201}{16900}s(1/25).
\]
The two readings use the same specimen. All limits of hierarchical
laminates satisfy a bound strictly below $0.2659042923076924$.
Consequently a physical specimen reaching this value would resolve the
unrestricted question negatively. An explicit half-volume laminate
gives $0.26548649548<\mathcal S<0.26548649549$.

The origins of the two constants differ. In the complementary energy
principle use two divergence-free matrix fluxes, one per contrast.
After fixing the fluxes, the optimal phase assignment fills the half
of the cell on which the difference of their weighted squared norms is
largest. This gives a trial from below for the physical supremum.
A variance bound for laminate spectral measures gives the smaller
laminate ceiling. The explicit trial leaves a gap of about
$4.18\cdot10^{-4}$ to the separating value. The precise inequality and
the flux normalization are given in Supplement~\ref{sup:cells}.

\section{6. Points and intervals on the single-pole arc}\label{problem-6}
\Milton{6}

Yes. A six-stage laminate gives a further point on the arc. With $z=\varepsilon_1/\varepsilon_2$, $s=z-1$, the normalized single-pole arc is $1+fs-f(1-f)s^2/[3(1+us)]$ \cite[Eq.~(7)]{p7:kmm}.
\subsection*{Intervals and exclusions}
Write $h=z-1$ and
\[
 K_f(u,h)=1+fh-\frac{f(1-f)h^2}{3(1+uh)} .
\]
To attain $u$ means that a single sequence of two-phase geometries with
phase-one fraction $f$ has isotropic response tending to $K_f(u,h)I$
for all positive contrasts. The scale-separated realizations below also
give locally uniform convergence off the nonpositive real axis.

\begin{theorem}[Bands of attainable poles]\label{new:arc}
Let $0<f<1$, $a=(2-f)/3$, $b=1-f$ and
\[
 R(f)=\frac{3-2f-\sqrt{1-4f(1-f)/3}}4 .
\]
Every pole in $[R(f),\min(a,b)]$ and in
$[\max(a,b),1-R(1-f)]$ is attainable. At $f=1/4$ the upper interval
extends to $1-r_7=0.8522624015603\ldots$, where $r_7$ is the unique
root of $672r^3-520r^2+96r-5$ in $(147/1000,148/1000)$.
No limit of hierarchical laminates with arbitrary normals and unbounded
depth attains a pole in any of
\[
 (1/4,3/8),\qquad(7/12,3/4),\qquad(7/8,11/12)
 \quad\text{at }f=1/4 .
\]
\end{theorem}
\begin{proof}\phantomsection\label{proof:new:arc}
In a coordinate stack, a diagonal conductivity is averaged harmonically
in the normal direction and arithmetically in the two tangential
directions. Retaining a common denominator $1+uh$ in the three
directional responses gives matching conditions for their numerators.
Two plies can thus agree in a chosen direction at every contrast. Three
further stacks carry that matching condition to a new pole with the
quantity $f+3u-1$ fixed. Iteration and phase interchange give the
two intervals in the statement. At small $f+3u-1$ the construction takes
a limit of successively deeper stacks. Its truncations are filled at the
correct fraction. Uniform ellipticity at any fixed positive contrast
controls the change caused by the final filling.

For the additional upper segment, first work at the interchanged
fraction $3/4$. The seven-stage contact starts at $r_7$. Two matched
recipes join it through $5/32$ to the interval ending at $1/4$.
Phase interchange sends $r$ to $1-r$. Supplement~\ref{sup:arc} gives
the matching rules, the algebraic endpoints and the verification of the
joins; both phase shares and the response identities are tested.

For necessity, equality on the arc concentrates the positive spectral
measure at one pole. Along a laminate tree the conditional phase fields
are matrix martingales. Their increments across a split share its layer
normal. At contact the limiting field in the phase with trace between
zero and one has rank at most two. Its squared trace and its
second moment consequently give
\[
       (u-(1-f))\bigl(u-(2-f)/3\bigr)\ge0 .
\]
The limiting-rank argument and the exterior support constraints are
given in Supplement~\ref{sup:arc}. At $f=1/4$ these yield the three
displayed gaps. The argument for the gaps uses the laminate structure.
\phantomsection\label{proofend:new:arc}\end{proof}

\begin{figure}[htbp]\centering
\includegraphics[width=.89\linewidth]{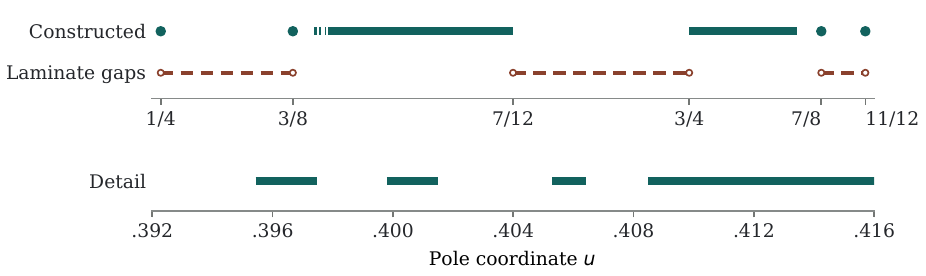}
\caption{Pole locations at $f=1/4$. Solid segments and filled points
have explicit constructions. Dashed segments with open ends exclude
limits of hierarchical laminates of arbitrary depth and with arbitrary normals. General two-phase geometries remain unclassified at the three stated probes.}\label{fig:arc}
\label{figure:arc}\end{figure}

Figure~\ref{fig:arc} shows the sets at quarter fraction. The detailed
view resolves three short constructed intervals next to the first long band.

\begin{center}\small
\begin{tabular}{@{}p{.46\linewidth}p{.47\linewidth}@{}}\toprule
Quarter-fraction poles & Realizations\\ \midrule
$1/4,\ 11/12$ and $3/8,\ 7/8$ & Coated spheres and isotropized coated cylinders, respectively.\\
$[(5-\sqrt3)/8,\,7/12]$ and $[3/4,\,1-r_7]$ & Two complete bands. Their combined length is about $0.27710$ in the pole coordinate.\\
$[0.395432,0.39748]$, $[0.3998,(39-\sqrt{177})/64]$, $[0.405301,0.406411]$ &
Additional certified subintervals. Decimal endpoints here are rounded into the attained sets.\\ \bottomrule
\end{tabular}
\end{center}

At $f=5/12$ one more point lies below the general band. A tree with
rational layer shares has normalized trace
\[
 \frac{1+3h/4+25h^2/432}{1+h/3}=K_{5/12}(1/3,h).
\]
Isotropizing its diagonal response gives the pole $1/3$ at this fraction.
It lies strictly below $R(5/12)=(39-\sqrt{219})/72$.
Supplement~\ref{sup:arc} gives the construction and exact check.

At quarter fraction, existence of a general physical geometry for the poles $1/3$, $2/3$
or $9/10$ remains open. A fixed smooth periodic cell is a different
admissible class from a laminate limit. At half fraction the interior
pole $u=1/2$ gives $(z^2+4z+1)/(3(z+1))$. Torquato and Kim
\cite[Eq.~(30)]{tk24} identify this expression with the Schwarz P and D
cells. Supplement~\ref{sup:arc} sets out the apparent conflict with a
smooth-interface obstruction and its hypotheses. The band and
laminate-exclusion statements above are independent of that issue.

\subsection*{A six-stage example}

\begin{proposition}\label{p6:arc}
At fraction $f=1/4$ a six-stage hierarchical laminate has isotropic response $(z^2+10z+5)/[8(z+1)]$ on $\mathbb C\setminus(-\infty,0]$. Its parameter on the arc is $u=1/2$.
\end{proposition}
\begin{proof}\phantomsection\label{proof:p6:arc}
Start with $P=zI$ and $Q=I$ at real $z>0$. Write $[U,V]_{i,p}$ for plies of shares $p,1-p$ with normal $e_i$. For diagonal matrices the normal entry is the harmonic mean and the two tangential entries are arithmetic means, by continuity of normal flux and tangential field. Form
\[
 A=[P,Q]_{3,3/4},\quad B=[Q,A]_{1,2/3},\quad
 C=[B,A]_{2,1/2},\quad D=[C,Q]_{3,1/2}.
\]
In that order their phase fractions are $3/4,1/4,1/2,1/4$. The entries of $A$ are $((3z+1)/4,(3z+1)/4,4z/(z+3))$, and those of $B$ are $((3z+1)/(2z+2),(z+3)/4,2(z+1)/(z+3))$. Thus $C=\mathrm{diag}(p,p,q)$, where $p=(z+3)(3z+1)/[8(z+1)]$ and $q=(3z+1)/(z+3)$. Write $D=\mathrm{diag}(T,T,N)$ with $T=(p+1)/2$, $N=2q/(1+q)$.

Rotate two copies of $D$ so that $N$ lies along axes 1 and 2. Stack these copies normal to axis 3 in shares $1/3,2/3$. With $h=N-T$ the entries are $(T+h/3,T+2h/3,T)$. An equal stack of this material and its quarter-turn about axis 1, normal to that axis, gives $(2T+N)I/3$. Substitution yields the stated function, equal to $1+s/4-s^2/[16(1+s/2)]$. The rational identities extend throughout the slit plane by analyticity. In the catalogue of \cite[Eq.~(9) and the following paragraph, Eqs.~(24)-(28)]{p7:kmm} the admissible parameters at this fraction are $1/4,3/8,7/16,7/12,3/4,7/8,11/12$. The displayed construction supplies a distinct response.
\phantomsection\label{proofend:p6:arc}\end{proof}

\section{7. One inclusion per period cell}\label{problem-7}
\Milton{7}

A double coating can be replaced by one connected inclusion per cell if the inclusion may have cavities. In the plane we also give a replacement with no cavities.

For a core of coefficient $k$ and relative volume $t$, surrounded by
material $q$ in dimension $n$, write
\[
 C_n(k,q,t)=q+\frac{ntq(k-q)}{nq+(1-t)(k-q)} .
\]
Indeed, an affine loading in direction $e$ has potential $A e\cdot x$ in
the core and $(B+D|x|^{-n})e\cdot x$ in the coating. At outer radius one
and inner radius $t^{1/n}$, continuity gives $A=B+D/t$ and
$kA=q(B+(1-n)D/t)$. Thus the outer normal flux divided by the normal
derivative of the affine potential with the same trace is
$q(B+(1-n)D)/(B+D)=C_n(k,q,t)$. Such a ball is neutral in a medium of
this coefficient. Successive coatings give the double-coating formula below.

\begin{lemma}[Periodic obstacle]\label{p7:obstacle}
Let $a\ge0$ and $c>0$. Suppose a bounded lattice-periodic obstacle $\psi$ has $\psi(x)+a|x|^2/2$ convex on $\mathbb R^n$. Among periodic $H^1$ functions $v\ge\psi$, the functional $\int (|\nabla v|^2/2+cv)$ has a unique minimizer. It is the least periodic $H^1$ majorant with $\Delta v\le c$, belongs to $C^{1,1}$, satisfies $D^2v\ge-aI$, and has $\Delta v=c$ on its open noncontact set.
\end{lemma}
\begin{proof}\phantomsection\label{proof:p7:obstacle}
Poincaré's inequality,
the lower bound on the mean and the energy bound on a minimizing sequence
give existence by weak compactness. Strict convexity modulo constants and the
positive linear term give uniqueness. The variational inequality gives
$\Delta v\le c$, with equality on $v>\psi$ once continuity is established.
Moreover $v$ is the smallest admissible majorant with $\Delta v\le c$.
Indeed, testing the inequality with $\min(v,V)$ for another such $V$ and
testing $\Delta V\le c$ with $(v-V)_+$ gives
$\int|\nabla(v-V)_+|^2\le0$. A positive constant difference would decrease
the energy. Taking $V=(v(x+h)+v(x-h)+a|h|^2)/2$ now gives $D^2v\ge-aI_n$
distributionally. The positive matrix measure $D^2v+aI_n$ has trace bounded
by $(c+na)\,dx$. Thus $v\in C^{1,1}$.
\phantomsection\label{proofend:p7:obstacle}\end{proof}

\begin{theorem}\label{p7:shell}
Let $n=2$ or $3$ and fix a Bravais lattice of covolume one,
$0<f<b<1$, and $\rho=1-f/b$. There is one bounded, path-connected open
inclusion of volume $f$ per cell, separated from its translates. It is a convex open body $O$ with countably infinitely many disjoint closed convex cavities removed. The cavities lie inside $O$ and accumulate exactly at its outer boundary. The total interface measure is finite. For coefficient $p$ in the inclusion and $q$ in the cavities and exterior, its effective tensor equals
\[
 C_n\big(C_n(q,p,\rho),q,b\big)I_n
\]
for every pair of positive real coefficients and every pair with positive
imaginary parts. In particular, whenever such a double coating maximizes
absorption, so does this same cell.
\end{theorem}

\begin{proof}\phantomsection\label{proof:p7:shell}
All potentials below are periodic on the lattice. Write $D^2$ for the Hessian.
We construct a filled body $O$, cavities $K\subset O$, and two
$C^{1,1}$ potentials $u,w$ satisfying
\begin{align}
 \Delta u&=\mathbf1_O-b,& D^2u|_O&=(1-b)I_n/n,\label{p7:outer}\\
 \Delta w&=\rho\mathbf1_O-\mathbf1_K,&
 D^2w|_K&=-aI_n,\qquad a=(1-\rho)/n.\label{p7:inner}
\end{align}
Indicators denote the periodic extensions.

The filled outer body is the periodic E-inclusion of Liu, James and Leo
\cite[Sec.~II, Eqs.~(4), (7)-(9)]{p7:ljl}.
In the plane Vigdergauz constructed periodic extremal inclusions
\cite[p.~390, abstract]{p7:vig}, whose boundary Grabovsky and Kohn
parametrized explicitly \cite[Sec.~3, Eq.~(3.18)]{p7:gk}.
Our second obstacle places the cavities inside this body to reproduce the double coating.

For the outer body apply Lemma~\ref{p7:obstacle} with $\psi(x)=-\min_{\ell}|x-\ell|^2/2$, $a=1$ and
$c=nb/(1-b)$, where $\ell$ runs over lattice points throughout. Write the
minimizer as $-U$. The central contact
$\{U=|x|^2/2\}$ is convex, since $|x|^2/2-U$ is nonnegative and convex.
It is compactly inside the Voronoi cell, since at a tie between two closest lattice
points $U$ would have to touch two upper quadratics with different gradients.
Its translates are all the contacts. On their interiors $D^2U=I_n$, and
their boundaries have zero measure by convexity (including the
empty-interior case). Elsewhere $\Delta U=-c$.
Integrating fixes the contact fraction at $b$, so the central interior
$O$ is nonempty, and $u=(1-b)U/n$ satisfies \eqref{p7:outer}.
Real analyticity of the boundary follows from the local obstacle equation. Here $|x|^2/2-U$ is the nonnegative classical
obstacle solution with $\Delta(|x|^2/2-U)=(n+c)\mathbf1_{\{U<|x|^2/2\}}$.
At every boundary point the convex contact contains the cone to a fixed
interior ball. The positive relative thickness at every small radius is
the regular-point hypothesis in \cite[corollary to Theorem~III]{p7:cks}.

Write $d=\operatorname{dist}(x,\partial O)$ inside $O$. We first choose a
periodic $C^{1,1}$ function $v_0$, positive in $O$ and zero outside, such that
\[
 D^2v_0\ge-aI_n/2,\qquad
 0<\mu=\rho-\Delta v_0\text{ a.e. in }O,\qquad
 \mu=O(d^4)\quad(d\downarrow0).
\]
In inward normal coordinates $(y,d)$ use $P=\sum_{k=2}^{6}A_k(y)d^k$.
Set $A_2=\rho/2$. In $\Delta P$ the coefficient of $d^j$ involves $A_{j+2}$
only through $(j+2)(j+1)A_{j+2}$, since the Laplacian has leading normal
part $\partial_d^2$. Choose successively the coefficients of $d,d^2,d^3,d^4$
in $\Delta P$ to be $0,0,0,-1$. Then $s=(2P/\rho)^{1/2}=d(1+O(d))$ is
smooth up to the inner side of the boundary, and $D^2P=\rho\nabla d\otimes
\nabla d+O(d)$. Set $v_0=P$ up to $s=\tau$ and then continue $v_0$ beyond that level as
a function of $s$ with derivative $\max(\rho\tau-a(s-\tau)/16,0)$,
continued by its final constant into the interior. For small $\tau$
the joins are in a collar where $1/2\le|\nabla s|\le2$ and $D^2s$
is bounded. On the transition the Hessian is
$-(a/16)\nabla s\otimes\nabla s+O(\tau)$. The values and first derivatives
match at both joins and at extension by zero. This proves the stated bounds.

Take the maximal dyadic cubes of side $l_j\le L(c_j)$ centered at $c_j\in O$,
where $L(x)=\eta\min(d(x),d_*)^2$ with $d_*>0$ fixed inside the collar
and $L=0$ outside $O$. Their closures cover $O$, since fine cubes
containing an interior point qualify. Their interiors are disjoint. The
parent fails. The estimate $l_j\asymp L(x)$ on each cube follows from
$L(c_{\rm parent})<2l_j$ and the Lipschitz bound on $L$, made arbitrarily
small by $\eta$. Also $L/d$ can be made arbitrarily small. At each center
use the tangent quadratic
\[
 Q_j=v_0(c_j)+\nabla v_0(c_j)\cdot(x-c_j)-a|x-c_j|^2/2 .
\]
Use all translates of these caps. The global Hessian bounds on $v_0$ give
\[
 a|x-c_j|^2/4\le v_0(x)-Q_j(x)\le C|x-c_j|^2 .
\]
Minimize with obstacle $\max(0,\sup_j Q_j)$ and $c=\rho$, writing $w$ for
the minimizer. The cap supremum is locally Lipschitz because $v_0,\nabla v_0$ are
bounded and caps contributing to a positive maximum have uniformly bounded
distance to their centers. The obstacle is semiconvex with constant $a$.
Comparison with $v_0$ gives $0\le H=v_0-w\le C l_j^2$ on cube $j$.
Since $v_0\asymp d^2$ near the outer boundary and is positive inside, reducing
$\eta$ ensures $w>0$ throughout $O$ and $w=0$ outside.

The contact $F_j=\{w=Q_j\}$ is convex by semiconvexity, contains $c_j$,
and lies in $B(c_j,C_1l_j)\Subset O$. Indeed $w(c_j)=v_0(c_j)$, and at
any contact $x$ one has $a|x-c_j|^2/4\le H(x)\le C L(x)^2$.
Since $L$ has arbitrarily small Lipschitz constant this yields the stated
ball, with a constant independent of $\eta$ (use $l_j\asymp L(c_j)$).
Off $O$ equality is impossible by the strict cap bound.
Thus contacts are locally finite in $O$. The gradients of the affine
functions $Q_j+a|x|^2/2$ are distinct because they are the gradients at the distinct
centers of the strongly convex $v_0+a|x|^2/2$. Differentiability of $w$ then
makes the $F_j$ disjoint. At a positive obstacle value the supremum is a
maximum, since a sequence of centers tending to $\partial O$ has
$v_0,\nabla v_0\to0$. Remove the contacts of positive volume, with their
boundaries, and let $K$ be the union of their interiors. The other contacts
and all contact boundaries have zero volume. The free equation, the cap
Hessian in $K$, and $C^{1,1}$ regularity prove \eqref{p7:inner}. Integration
gives $|K|=\rho b$. In fact every $F_j$ has positive volume. Otherwise take a ball about $c_j$ meeting no other contact. There $\Delta H=-\mu<0$ almost everywhere, whereas $H\ge0$ takes its minimum zero at $c_j$, contradicting the strong minimum principle. The cube centers approach every point of $\partial O$ because their cubes cover $O$ and their side lengths are $O(d^2)$ near its boundary. Local finiteness gives the asserted accumulation set.
The remainder in $O$ is open. To join two of its points, take a compact
neighborhood of the segment between them inside the body. It meets finitely
many removed bodies. Enlarge the bodies hit by the segment slightly, keeping
them disjoint inside $O$, away from its endpoints and all other bodies.
Replace each chord by a path on the enlarged convex body's connected boundary.

We verify finiteness of interface measure. Put
$M(r)=|K\cap\{r<d<2r\}|$, $l=r^2$. In this band the cube sides are
comparable to $l$, so there are $O(r/l^n)$ such cubes by the smooth collar
volume estimate. Balls $B_j$ about their centers of a sufficiently large fixed
multiple of $l_j$ have bounded overlap, stay at comparable depth, and their
half balls contain both the cubes and contacts. Write $m_j=|B_j\cap K|$.
Since $\Delta H=\mathbf1_K-\mu$ in $O$ and $H(c_j)=0$, the harmonic
extension of $H|_{\partial B_j}$ has center value at most the integral of
the centered Dirichlet Green kernel over $K\cap B_j$. Its integral is
maximal on a concentric ball of volume $m_j$. In dimension three
the $C/|x-c_j|$ kernel gives $C l^2(m_j/l^n)^{2/3}$. In dimension two
the logarithmic kernel gives $C l^2 u(1+|\log u|)$, $u=m_j/|B_j|$,
also bounded by $C l^2 u^{2/3}$. Harnack's inequality for the nonnegative
harmonic extension and the Green representation of $H$ therefore give
\[
 \sup_{B_j/2}H\le C l^2\big[(m_j/l^n)^{2/3}+r^4\big] .
\]
If $M(r)=O(r^{1+\beta})$, bounded overlap and concavity show that the
integral of $H$ over the band is bounded by
$Cr l^2(r^{2\beta/3}+r^4)$. Here adjoining dyadic bands contribute only
a fixed factor. Test $\Delta H=\mathbf1_K-\mu$ with a smooth cutoff
equal to one on the band, supported at comparable depths with Laplacian
$O(r^{-2})$. This improves $\beta$ to $\min(4,2+2\beta/3)$.
Starting at $\beta=0$ gives $\beta=10/3$ after two steps.
On $F_j$, $a|x-c_j|^2/4\le H$, so the same estimate and concavity give
\[
 \sum_{j\text{ in band}}\operatorname{diam}(F_j)^{n-1}
 \le C\frac r l
 \left(r^{(n-1)\beta/3}+r^{2(n-1)}\right).
\]
Both exponents exceed one. Summation over dyadic $r$ converges since
$l=r^2$. Each convex body's surface measure is bounded by a dimensional
constant times its diameter to power $n-1$ (metric projection of an enclosing
sphere is a surjective contraction onto its boundary). Including the smooth
outer surface proves the interface assertion.

It remains to verify the tensor. Set $k=C_n(q,p,\rho)$,
$\lambda=(p-q)/(p-a(p-q))$, so $k=p(1-\rho\lambda)$.
With $k$ throughout $O$, use the matrix of electric fields
$E_0=I+sD^2u$, where $s=-(k-q)/(q+(k-q)(1-b)/n)$. Its interior value is
$AI$ with $A=1+s(1-b)/n$. Indeed, its flux $J_0=q E_0+(k-q)A\mathbf1_O I$
has divergence $(qs+(k-q)A)\nabla\mathbf1_O=0$ and mean
$[q+(k-q)Ab]I=C_n(k,q,b)I$. For the actual phases take
$E=E_0-A\lambda D^2w$. Its columns are gradients with the prescribed means.
Using \eqref{p7:inner} separately in $K$, $O\setminus K$ and the exterior,
and $(p-q)(1+a\lambda)=p\lambda$, one obtains for the new flux
\[
 J-J_0=Ap\lambda\big[(\Delta w)I-D^2w\big].
\]
This has zero mean and divergence-free columns in distributions, including
at every cavity accumulation point. These formulas first hold for positive
coefficients. They also hold for strictly passive ones, since the inner
coated ball has strictly passive $k$ by its energy identity, and all
denominators are positive convex combinations of the relevant coefficients
up to positive factors. Uniqueness of the cell fields gives the result.
\phantomsection\label{proofend:p7:shell}\end{proof}

\subsection*{A simply connected planar maximizer}

\begin{theorem}\label{p7:star}
There is one simply connected inclusion of area $f=1/4$ per unit square,
with a real-analytic fourfold-symmetric boundary lying strictly inside the
square. Its effective tensor $Z=zI$ satisfies
\begin{equation}\label{p7:disk}
 \frac{z-q}{z+q}=\frac{ft}{1-rt^2},\qquad
 t=\frac{p-q}{p+q},\qquad r=15/32 ,
\end{equation}
for all positive real and all strictly passive phase pairs $(p,q)$, with $p$
in the inclusion. These are the disk coatings with cumulative area fractions
$15/68$ and $8/17$. For $p=-15+i$, $q=1+i$ it gives $Z=4iI$. Among all periodic two-phase cells with area fraction $1/4$ of $p$ at these data, and all limits of their effective tensors, the largest least directional loss $\lambda_{\min}((Z-Z^*)/(2i))$ is $4$. The competitors may be anisotropic.
\end{theorem}

\begin{proof}\phantomsection\label{proof:p7:star}
Write $P$ for projection onto constants and $\Gamma$ for projection onto
mean-zero gradients in the normalized periodic $L^2$ space of vector fields.
If $\chi$ is the inclusion indicator, put $B e=\chi e$ and
$T=\chi(I-P-2\Gamma)\chi$. This is a self-adjoint contraction. For normalized
flux $J/q$ and electric field $E$ put $F=(J/q+E)/2$, $v=(J/q-E)/2$.
The phase law gives $v=t\chi F$ and the field constraints give
$F=\langle F\rangle+(I-P-2\Gamma)v$. Taking means after eliminating $v$ gives
\begin{equation}\label{p7:cayley}
 \Phi=(Z-qI)(Z+qI)^{-1}=t B^*(I-tT)^{-1}B .
\end{equation}
The inverse exists for nonreal $t$ and real $|t|<1$, covering these phase
pairs. The two identities $T^2B=rB$, $B^*TB=0$ therefore suffice
for \eqref{p7:disk}, since $B^*B=fI$.

We first verify optimality. Rotating vector components by $90^\circ$
interchanges gradients and divergence-free fields of zero mean and hence
changes $T$ to $-T$. The trace spectral measure
$(2f)^{-1}\operatorname{Tr}(B^*dE_T B)$ is thus a symmetric probability
measure on $[-1,1]$. With its image $\nu$ under $\lambda\mapsto\lambda^2$,
write $s=t^2$, $u=\operatorname{Tr}\Phi/2$ and
$Q=u/(ft)=\int(1-sx)^{-1}d\nu$. The identity
$\operatorname{Im}(s/(1-sx))=(\operatorname{Im}s)|1-sx|^{-2}$
and Cauchy-Schwarz give
$|Q|^2\le\operatorname{Im}(sQ)/\operatorname{Im}s$ when
$\operatorname{Im}s\ne0$. At the stated phases $t=(28+4i)/25$, whence
\[
 7|u|^2\le\operatorname{Re}u+7\operatorname{Im}u,\qquad
 |z_0-7i/3|\le5/3,\qquad z_0=q(1+u)/(1-u).
\]
To pass to anisotropic cells, put $C=(Z+qI)^{-1}$, $c=\operatorname{Tr}C/2$.
Passivity makes $-\operatorname{Im}c>0$, and $z_0=1/c-q$. With
$m=\lambda_{\min}(\operatorname{Im}Z)$ one has
$-\operatorname{Im}c=\operatorname{Tr}[C^*(\operatorname{Im}Z+I)C]/2
\ge(m+1)|c|^2$.
Thus $m\le\operatorname{Im}z_0\le4$. Continuity gives the bound in the
closure. Formula \eqref{p7:disk} gives $u=(7+4i)/13$ and $z=4i$ at
these data. The coating formula identifies \eqref{p7:disk} by taking
$b=f/(1-r)$ and core fraction $rb$.

\smallskip\noindent\emph{From the cell equation to a boundary equation.}
In length units $1/\delta$ times the physical lengths, represent a real
two-vector by $u_x-iu_y$. In the inclusion the action of $T$ becomes
the real-antilinear operator
\[
 \mathcal Kp(z)=\frac1\pi\operatorname{p.v.}\int_D\wp_\delta(z-y)
           \overline{p(y)}\,dA(y),\qquad
 \wp_\delta(z)=\delta^2\wp(\delta z)
\]
for the square-lattice Weierstrass function. Indeed $\Gamma$ has the Hessian
of the periodic Green function as kernel, with singular part
$\nabla^2\log|z|/(2\pi)$. Its distributional diagonal term cancels $I$ in
$I-P-2\Gamma$. Off the poles the matrix kernel is symmetric and traceless
(the constant Laplacian cancels $-P$). Its complex entry is holomorphic
by the Cauchy-Riemann equations and has principal part $1/(\pi z^2)$.
Its difference from $\wp_\delta/\pi$ is constant by double periodicity,
and a quarter-turn changes its sign. This proves the representation with the
correspondingly scaled cell mean.

For holomorphic $p$, Green's formula gives
\[
 \mathcal Kp(z)=\frac{1}{2\pi i}\int_{\partial D}
  L_\delta(y-z)\overline{p(y)}\,d\bar y,\qquad
 L_\delta(u)=\delta\zeta(\delta u),\quad L_\delta'=-\wp_\delta .
\]
In excising a small circle, the leading singular boundary term integrates
to zero by its angular frequency $-2$. Thus the two successive fields are
holomorphic. Parameterize the boundary by $\gamma(w)=w(R(w)+h(w))$,
$w=e^{i\theta}$, and set $N=w\partial_w$, $G=N\gamma$,
$p^\#(w)=p(1/w)$ (conjugation on real Laurent coefficients).
If $\mathcal P$ denotes the interior trace of the Cauchy integral with
kernel $L_\delta(y-z)\,dy/(2\pi i)$, the two boundary traces are
\[
 M=-G^\#/G=d\bar\gamma/d\gamma,\qquad
 g=\mathcal PM,\qquad W=\mathcal P(Mg^\#).
\]
A quarter-turn changes the sign of $g$, so its inclusion mean is zero.
If $\operatorname{Re}W=r$ on the boundary, the maximum principle and
reflection symmetry give $W=r$ throughout. Antilinearity then proves
$B^*TB=0$ and $T^2B=rB$ for both coordinate loads.

For the remaining boundary equation use an even fourfold radial graph.
Its trial profile is a cosine polynomial with 76 integer coefficients.
The area constraint fixes the constant mode. The Cauchy boundary operator
acts on exponentially weighted Fourier series. After enclosing the low
modes of its derivative, an explicit inverse for the high modes controls
the entire tail. The residual and inverse defect give a contracting ball
in residual coordinates. Positivity and strict separation from translated
copies persist throughout that ball. The numerical bounds, the construction
of the approximate inverse and the complete proof of this step are in
Supplement~\ref{sup:inclusions}. The resulting exact profile satisfies
both operator identities above and hence gives \eqref{p7:disk}.
\phantomsection\label{proofend:p7:star}\end{proof}

The integers $n_j$ in the finite witness above, in order, are
\begin{center}\scriptsize\ttfamily
\begin{tabular}{rrrr}
46787862359913840640 & 16907746505615800320 & 5056678727937039360 & 756857360158470272 \\
-412653040809127936 & -497030898027754816 & -326830609237935296 & -168709141971578944 \\
-72005095969870176 & -23620206017806516 & -3163627280048868 & 3590534490281608 \\
4527381467139087 & 3495841473788878 & 2163448507710242 & 1118369296549592 \\
463327111178553 & 120499617642592 & -24052764701161 & -63640320311843 \\
-58244254975671 & -39934354067588 & -22660756515290 & -10607001436259 \\
-3648631278142 & -314312384934 & 889942053706 & 1049308312740 \\
813874366881 & 510127631528 & 266717459357 & 110557508953 \\
26872143920 & -9184765911 & -19090844684 & -17317942680 \\
-12005127677 & -6911408918 & -3272789482 & -1122496307 \\
-69925172 & 316620580 & 366713420 & 285727129 \\
180598640 & 95214129 & 39647755 & 9444865 \\
-3745237 & -7407386 & -6725959 & -4690159 \\
-2720853 & -1298718 & -447496 & -24603 \\
131723 & 151540 & 117720 & 75893 \\
43115 & 22152 & 7763 & -927 \\
-5033 & -6761 & -6449 & -5707 \\
-2754 & 1843 & 2483 & 2272 \\
652 & -1427 & -2714 & 751 \\
\end{tabular}
\end{center}

\begin{corollary}\label{p7:spatial-max}
At $p=-15+i$, $q=1+i$ and $f=1/4$ in three dimensions, one connected
perforated inclusion with finite interface measure on each prescribed
Bravais lattice attains least directional loss $(51+10\sqrt{17})/19$. This is the largest value over all periodic two-phase cells with that fraction of $p$, including anisotropic ones and limits of their effective tensors.
\end{corollary}
\begin{proof}\phantomsection\label{proof:p7:spatial-max}
For $z=\operatorname{Tr}Z/3$, set $m=fp+(1-f)q=-3+i$,
$n=(1-f)p+fq=-11+i$, and $z=m-48/(y+n)$.
To realize $zI$ at these phase data, choose a unit $n$ where $\det[n,(\operatorname{Re} Z)n,(\operatorname{Im}Z)n]=0$, possible since this determinant is odd on the sphere. Choose a real unit $v$ perpendicular to all three columns. Reflection $S$ in $v^\perp$ fixes $n$ and $Zn$. Equal plies of $Z$ and $SZS^t$ normal to $n$ have identical normal columns, so their tensor is the arithmetic mean. Reflection invariance makes $v$ a real principal axis of that mean. The three trace-preserving layering steps of Theorem~\ref{p4:directions} give $zI$. The fraction is preserved throughout.

The Kern-Miller-Milton bound puts $y$ in the disk of center
$c=-14+15i$ and radius $5\sqrt{17}$, since $2p,2q,-q$ lie on its boundary
and $-p$ is exterior \cite[p.~5 and Eqs.~(13), (14), (69), (77)-(80)]{p7:kmm}.
Since $c+n=-25+16i$ and $|c+n|^2-425=456$, inversion gives
\[
 \left|z-\frac{-7+51i}{19}\right|\le \frac{10\sqrt{17}}{19}.
\]
The least directional loss is at most $\operatorname{Im}z$. In
Theorem~\ref{p7:shell} choose $b=832/(778+285\sqrt{17})$, which lies in
$(1/4,1)$. Substitution in $C_3(C_3(q,p,1-f/b),q,b)$ gives the top of
this disk. The exact substitution in $\mathbb Q(\sqrt{17},i)$ is asserted in \path{certificates/baseline/p07/optical.py}.
\phantomsection\label{proofend:p7:spatial-max}\end{proof}

\begin{proposition}\label{p7:nohole}
At the three-dimensional phase data and fraction of
Corollary~\ref{p7:spatial-max}, a periodic cell attaining the maximal
least directional loss cannot have the phase $q=1+i$ equal, up to null
sets, to a connected open subset of the period torus. In particular a
single smooth inclusion with connected exterior and no cavities cannot
attain this loss, in either phase assignment.
\end{proposition}
\begin{proof}\phantomsection\label{proof:p7:nohole}
Use $Z$ for the effective tensor and $z_*=(-7+i(51+10\sqrt{17}))/19$.
Equality at the top of the trace disk in Corollary~\ref{p7:spatial-max}
forces $\operatorname{Tr}Z/3=z_*$, and the three eigenvalues of
$\operatorname{Im}Z$ must all equal $\operatorname{Im}z_*$.
By reciprocity $Z=z_*I+D$ for a real symmetric traceless $D$.
Take matrix fields with mean electric field $iI$, and put
$X=-\operatorname{Re}J$, $Y=\operatorname{Re}E$. In the phase
$\varepsilon_j=p_j+i$, $p_1=-15$, $p_2=1$, one has
$\operatorname{Im}E=X+p_jY$, so the local energy is
$L_j=|X+p_jY|^2+|Y|^2$. Its mean is $\operatorname{Tr}\operatorname{Im}Z$.

Write subscripts $b,d,a$ for the scalar, symmetric traceless and skew
parts of real matrices, respectively. With $t=(\sqrt{17}-3)/40$ set
\[
 T=t[-|X_b|^2+2|X_d|^2-28X:Y+4|Y_b|^2-2|Y_d|^2+2|Y_a|^2].
\]
Then $\langle T\rangle\ge T(\langle X\rangle,\langle Y\rangle)$.
Indeed, rotate any nonzero Fourier direction to the first axis. The
first row of $X$ vanishes by equilibrium, and only the first row of $Y$
can be nonzero by compatibility. The $Y$ and cross terms vanish,
and the rest divided by $t$ is
$|X_{21}|^2+|X_{31}|^2+|X_{22}-X_{33}|^2+|X_{23}+X_{32}|^2$.
For $S_j=L_j-T$, its $2\times2$ blocks in the three sectors have common
off-diagonal $p_j+14t$ and respective diagonal pairs
\[
 (1+t,p_j^2+1-4t),\quad
 (1-2t,p_j^2+1+2t),\quad (1,p_j^2+1-2t).
\]
The second determinant is zero in both phases since $200t^2+30t=1$, and
the first is $t(p_j-1)(p_j-27)$, equal to $3t(p_j^2-1)$ at the two specified phase values. The last block is two-thirds the first
plus one-third the second, positive definite even in phase 2 because
the two singular blocks there have distinct null lines. Thus $S_j\ge0$.

Put $F_j=(\langle X\rangle_j,\langle Y\rangle_j)$ for the conditional
phase means, $f_1=1/4$, $f_2=3/4$. The difference
\[
 \mathcal E=\operatorname{Tr}\operatorname{Im}Z
 -\sum f_j S_j(F_j)-T(\langle X\rangle,\langle Y\rangle)
\]
is nonnegative, and dominates $\sum f_j\langle S_j((X,Y)-F_j)\rangle_j$.
For $D=0$ the mean electric fields are $i\beta_j I$, where
$\beta_1=(z_*-q)/(f_1(p-q))$, $\beta_2=(p-z_*)/(f_2(p-q))$,
$p=-15+i$. Write $x=-\operatorname{Re}(ip\beta_1-iq\beta_2)$ and
$y=\operatorname{Re}(i\beta_1-i\beta_2)$. Then
\[
 \mathcal E/3=\operatorname{Im}z_*-\sum f_j |\beta_j|^2
             +f_1 f_2 t(-x^2-28xy+4y^2)=0 .
\]
For verification $57x=302+110\sqrt{17}$, $57y=32+10\sqrt{17}$, and
the first difference is $(2995+1060\sqrt{17})/2166$.
For general $D$ the phase means of $E$ acquire the purely imaginary
additions $iD/[f_1(p-q)]$ and $-iD/[f_2(p-q)]$ from the two prescribed
means of $E,J$. The $Y$ means are unchanged, and the $X$ means acquire
the corresponding real additions. Thus the sole additional term in
$\mathcal E$ is $-(1-2t)|D|^2/[f_1f_2(p-q)^2]$, forcing $D=0$.

Vanishing of both skew-block energies implies that $E$ is symmetric
everywhere. The positive scalar block in phase 1 forces its trace
constant there. Divide all fields by $i$ to restore mean $I$.
In the interior of phase 2 the trace is constant on each connected
component, since $\partial_k\operatorname{tr}E=\sum_j\partial_j E_{kj}
=\sum_j\partial_j E_{jk}=0$ by symmetry, compatibility and equilibrium.
If phase 2 agrees with one connected open set, the trace therefore takes
only two values. For a symmetric compatible field, every nonconstant
Fourier mode is determined by its trace, being a multiple of
$\xi\otimes\xi/|\xi|^2$. It follows that $E=I+cH$,
$H=D^2\Delta^{-1}(\chi-f_1)$ real, with
$c=3(\beta_1-\beta_2)=-[50+i(32+10\sqrt{17})]/19$.
Let $\Gamma$ denote mean-zero gradient projection, columnwise, so
$H=\Gamma(\chi I)$. Applying it to $(q+(p-q)\chi)E$ gives
$\Gamma(\chi H)=(-q/(p-q)-1/c)H$. The left side is real and $H\ne0$.
The scalar on the right has imaginary part
$1/16-19s/(2500+s^2)<0$, $s=32+10\sqrt{17}\in(72,82)$, a contradiction. The radical identities and strict sign in this equality case are also asserted in \path{certificates/baseline/p07/optical.py}.
\phantomsection\label{proofend:p7:nohole}\end{proof}

\subsection*{Hole-free replacements on a range of phase data}
For a square-periodic inclusion of area $F$, put
$t=(p-q)/(p+q)$, with $p$ inside and $q$ outside the inclusion. The
target disk law at core ratio $r$ is
\[
       \frac{z_*-q}{z_*+q}=\frac{Ft}{1-rt^2},
       \qquad 0\le r<1-F .
\]
Matching this identity for every contrast is stronger than matching the
loss at a single frequency. The fraction $F$ always denotes the
inclusion area here, including when the phase labels are interchanged.

\begin{theorem}[A continuum of simply connected periodic replacements]\label{new:firstbranch}
One square-symmetric, simply connected real-analytic inclusion per
unit square realizes the disk law for every parameter pair in each row
below. The shape may vary with the parameters.
\[
\begin{array}{c|c}
 F & r \\ \hline
 {[4/25,\,2/5]} & 65/96-5F/6\\
 1/4 & [1/50,\,101/200]\\
 {[1/10,\,1/5]} & [1/50,\,1/5]\\
 {[1/5,\,3/10]} & [1/50,\,13/100]
\end{array}
\]
The first row attains the maximum least directional loss for
$p=-15+i$, $q=1+i$ at every indicated phase-one fraction.
\end{theorem}
The proof uses 88,880 node records on 69 pieces. One node per piece
was regenerated. For the original 22,576 records the pre-submission
inputs were recovered; the additional 66,304 contain both input hashes.
Section~\ref{detail:node-example} works through one whole interval.
Section~\ref{detail:square-extension} gives the new rectangle bounds.

\begin{proof}\phantomsection\label{proof:new:firstbranch}
Use the boundary equation in Theorem~\ref{p7:star}. On each parameter
piece approximate the boundary by a polynomial in the parameter with
Fourier coefficients in its angle. A polynomial family of approximate
right inverses acts on the residual. Ball arithmetic bounds the error
and derivative at the Chebyshev nodes. A bound on complex parameter
ellipses controls the tails of the two interpolation series. It gives a
uniform bound between nodes, including both endpoints of each piece.

In residual coordinates let $\epsilon$ bound the error, $d$ the
linear inverse defect, $K$ the inverse norm and $H$ the second
derivative. On a residual ball of radius $\beta$ it suffices that
$c=d+HK^2\beta<1$ and $\epsilon\le(1-c)\beta$. Bounds on the radial
profile separate the corrected curve from its translates. There are
30 adjoining pieces in the first row, 17 in the second, and 22
closed rectangles whose unions give the last two rows. Supplement~\ref{sup:inclusions} gives the
interpolation argument and the files needed to reproduce the bounds.
At the specified phase data the upper supporting point of the planar
disk arcs has $r=65/96-5F/6$ on the first branch, as follows by
differentiating the imaginary part of the disk law on both arcs.
\phantomsection\label{proofend:new:firstbranch}\end{proof}

\begin{proposition}[A second planar branch]\label{new:secondbranch}
On the triangular lattice of unit cell area, there is one star-shaped
real-analytic inclusion of phase $2$, of area $F=1/100$, whose full
law is the doubly coated disk law with
\[
 r_*=\frac{35\sqrt{113}-100}{352}
             +\frac{25+5\sqrt{113}}{88}\,F .
\]
For phase $1$ of permittivity $-15+i$ in the matrix and phase $2$ of
permittivity $1+i$ in the inclusion, it realizes the second maximizing
disk branch at phase-one fraction $99/100$.
\end{proposition}
\begin{proof}\phantomsection\label{proof:new:secondbranch}
The reduction and constants are given in Supplement~\ref{sup:inclusions}.
Threefold symmetry removes the mean two-direction field. The real part
of the remaining boundary field is prescribed to equal $r_*$ and the
area to equal $F$. In an analytic normal parametrization, write the
unknown displacement as $Ay$ for an approximate right inverse $A$.
With $\|A\|<48.50$, the residual is below $3.391\cdot10^{-18}$, the
linear defect is below $0.0082$ and the second-derivative bound is below
$1.04\cdot10^{10}$ on the required ball. At $\|y\|\le6\cdot10^{-18}$
the contraction factor is below $0.0085$ and the ball maps strictly
inside itself. The displacement is at most $2.91\cdot10^{-16}$ in the
weighted profile norm. The corrected angular derivative is positive and
the boundary lies inside radius $0.2075$, below half the shortest period.
The optical comparison for the two disk arcs is included there as well.
\phantomsection\label{proofend:new:secondbranch}\end{proof}

\begin{figure}[htbp]\centering
\includegraphics[width=.79\linewidth]{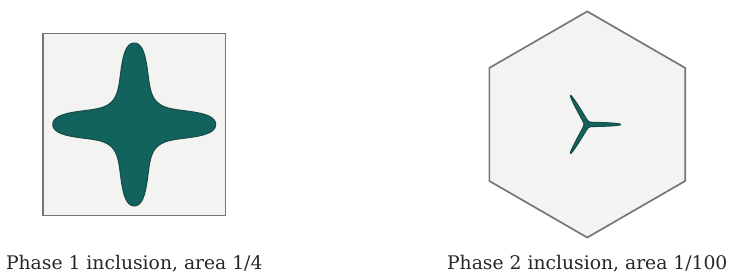}
\caption{Reference boundaries for Theorem~\ref{p7:star} and
Proposition~\ref{new:secondbranch}, plotted from the certificate
coefficients on the same length scale in unit-area period cells.
The left solid is fourfold symmetric on a square lattice. The right
solid is threefold symmetric on a triangular lattice. The
certificates bound the correction from each trial curve to an exact
double-coating law.}\label{fig:boundaries}
\label{figure:boundaries}\end{figure}

\subsection*{Finite cavities in space}
\begin{proposition}[An obstruction from the period lattice]\label{new:flat}
At the spatial quarter-fraction data of Corollary~\ref{p7:spatial-max},
let $\lambda_1,\lambda_2$ be the first two successive minima of the
Bravais lattice. There is a finite $M$ such that, if
$\lambda_2/\lambda_1>M$, the maximum cannot be attained by one connected
phase-one shell with finitely many separated Lipschitz cavities. Its
filled body is assumed bounded with separated translates and connected
exterior, with disjoint closed cavities in its interior.
\end{proposition}
\begin{proof}\phantomsection\label{proof:new:flat}
Normalize the shortest period to one. In a hypothetical sequence with
$\lambda_2/\lambda_1\to\infty$, the exterior obstacle field on a
central segment converges, after continuation inward from the boundary,
to the flat-layer field. Its Hessian tends uniformly to $e\otimes e$,
with $e$ along the shortest period. At the first cavity the shell
identities instead require a rank-one contact matrix, while the
continued matrix there is positive definite. If no cavity intervenes,
the two endpoint gradients on the segment agree, contradicting the
strict positive change of their $e$ component. The estimates and
continuation up to first contact are given in
Supplement~\ref{sup:inclusions}.
\phantomsection\label{proofend:new:flat}\end{proof}

The obstruction leaves existence on some Bravais lattice open. On the
simple cubic lattice at least 45 cavities are necessary and only 45 and
49 survive the count tests below 50. A countably perforated shell
continues to exist on every lattice by Theorem~\ref{p7:shell}.

\section{8. A finite test for the full elastic tensor}\label{problem-8}
\Milton{8}

Yes. In two and three dimensions a finite test characterizes the full
tensor at every phase fraction, with void or with perfectly rigid material.
For planar void the test reduces to two scalar inequalities in closed form.

Fix the dimension \(d=2,3\), the isotropic solid tensor
\(CE=2\mu E+\lambda\operatorname{tr}(E)I\), with \(\mu>0\) and
\(d\lambda+2\mu>0\), and its volume fraction \(0<f<1\).
On \(V=\mathrm{Sym}(d)\) use the Frobenius product and denote the identity
operator by \(I_V\). With void let
\(D=C\) and take stiffness as the response. With a perfectly rigid phase
let \(D=C^{-1}\) and take compliance. A zero eigenvalue of the latter
means an infinite stiffness. An attainable response belongs to the
closure of cubic-periodic cell responses at the prescribed fraction,
allowing measurable phases. Lemma~\ref{p8:substitution} permits any
full-rank period lattice. All averages are per unit cell volume.

For a unit normal \(n\) put \(J_n=\{a\odot n:a\in\mathbb R^d\}\),
where \(a\odot n=(a n^t+n a^t)/2\), and \(N_n=J_n^\perp=\{H:Hn=0\}\).
Write \(L_n=J_n\) for void and \(L_n=N_n\) for rigid. Let \(M_n\)
be the orthogonal projector onto \((D^{1/2}L_n)^\perp\), and set
\(\mathcal M=\mathrm{conv}\{M_n:|n|=1\}\).

\begin{theorem}\label{p8:main}
In either model the entire response closure is the compact convex set of
symmetric \(T\succeq0\) for which there exists \(M\in\mathcal M\) with
\begin{equation}\label{p8:lmi}
\begin{pmatrix} I-D^{-1/2}TD^{-1/2} & \sqrt{1-f}\, I\\
\sqrt{1-f}\,I &(1-f)I+fM\end{pmatrix}\succeq0.
\end{equation}
At the pure fractions
\(0,1\) the only responses are respectively \(0,D\).
\end{theorem}

The test includes all couplings and all ranks. It is finite even with
the elementary parametrization of
\(\mathcal M\) by at most \(1+\dim\mathrm{Sym}(V)\) normals and
nonnegative weights summing to one (Carath\'eodory's theorem).
Corollaries~\ref{p8:planar} and~\ref{p8:spatial} give smaller
semidefinite descriptions. The distance of a candidate to the closure
is the minimum Hilbert-Schmidt distance to a tensor satisfying
\eqref{p8:lmi}.

Avellaneda's directional bound and its attainment by sequential laminates
are classical \cite[Sec.~23.3, Eqs.~(23.20)-(23.25)]{book}.
Mixed energy and complementary-energy minima were obtained for several
families of loads in the compliant and rigid limits
\cite[Thms.~10.1 to 10.2, pp.~86 to 87]{p8void}, \cite[Thms.~6.1 to 6.2, pp.~110 to 111]{p8rigid}.
The additional step in Lemmas~\ref{p8:short} and~\ref{p8:rigshort}
removes the response to an arbitrary loading \(e\), a strain or a stress,
while preserving the dual energies for loads in \(e^\perp\).

\subsection{Cell fields and substitution at exact fraction}
Let \(S\) be the elastic part of a period cell. Write \(z\) for a
compatible strain \(F+e(v)\) in the void model and an equilibrated
symmetric stress of mean \(F\) in the rigid model. Here \(e(v)=
(\nabla v+\nabla v^t)/2\), with \(v\) periodic in \(H^1\).
The other type of field is denoted by \(y\). Equilibrium is rowwise
distributional divergence zero and all stresses are in \(L^2\).
For these two types of fields \(\langle z:y\rangle=\langle z\rangle:
\langle y\rangle\) because at nonzero Fourier frequency \(k\) their spaces
are \(L_{k/|k|}\) and its orthogonal complement. Conversely this
Fourier characterization gives every field of each type, by inverting
the symmetric-gradient symbol on its range.
Write \(|u|_D^2=u:Du\).

\begin{lemma}\label{p8:projection}
Define the cell response \(A\) by \(F:AF=\inf_z\langle\mathbf1_S|z|_D^2\rangle\).
There is a dual field \(y_F\) vanishing off \(S\), of mean \(AF\) and
cost \(\langle |y_F|_{D^{-1}}^2\rangle=F:AF\). There are primal trials
of mean \(F\) with \(\langle\mathbf1_S|z-D^{-1}y_F|_D^2\rangle\)
arbitrarily small. The dual minimum for prescribed mean \(H\) is
\(H:A^\dagger H\) on \(\mathrm{ran}\,A\), with \(A^\dagger\) the inverse on
that range extended by zero on the kernel. The minimum is infinite
outside the range.
\end{lemma}
\begin{proof}\phantomsection\label{proof:p8:projection}
Project \(F\), in the \(D\)-metric on \(L^2(S;V)\), onto the
orthogonal complement of the restricted zero-mean primal fluctuations.
Call its projection \(p_F\). Then \(y_F=\mathbf1_S Dp_F\) has dual
type by orthogonality on the full cell. Polarizing
\(F:AF=\|p_F\|_D^2\) gives \(\langle y_F\rangle=AF\).
Primal trials approach \(p_F\) on \(S\) by the definition of projection.
Every supported dual field of mean \(H\) obeys
\(|F:H|^2\le (F:AF)\langle |y|_{D^{-1}}^2\rangle\) by pairing with
these trials. Thus \(H\perp\ker A\), and for \(H=AF\) the field
\(y_F\) realizes equality. In particular \(0\preceq A\preceq fD\)
by the constant trial.
\phantomsection\label{proofend:p8:projection}\end{proof}

\begin{lemma}[Changing the fraction]\label{p8:volume}
Suppose a sequence of cells has elastic fraction \(s\ge f\) with
\(s-f\to0\). Suppose finitely many primal trials and supported dual
trials are given, the latter uniformly bounded in \(L^2\).
There are polyhedral cells at exactly \(f\) and trials of the same
means whose primal energies increase by \(o(1)\) at most and whose
dual fields change by \(o(1)\) in \(L^2\). For a fixed measurable
cell at \(s=f\) the errors can be made arbitrarily small.
\end{lemma}
\begin{proof}\phantomsection\label{proof:p8:volume}
On a small cube \(Q\) of side \(h\) a dual field can be cleared from a
concentric cube of side \((1-2a)h\), with change supported in \(Q\),
zero change of mean and \(L^2\) error at most \(c_a\|y\|_{L^2(Q)}\)
for any \(0<a\le1/4\). In the strain case cut its displacement down to
a rigid motion on the central cube, using Korn's inequality
\(h^{-1}\|u-r\|_2+\|\nabla(u-r)\|_2\le c\|e(u)\|_2\).
In the stress case take \((1-\zeta)y\), with \(\zeta=1\) on the
central cube, zero near \(\partial Q\), and correct its divergence in
the frame by a symmetric stress with zero boundary traction.
The required force \(g=y\nabla\zeta\) annihilates rigid motions
by testing equilibrium against \(\zeta r\). Solve
\(\int e(v):e(w)=-\int g\cdot w\) on the connected frame modulo
rigid motions. Korn's inequality and scaling give
\(\|e(v)\|_2\le c_a h\|g\|_2\). Extending the correction by zero
preserves equilibrium. The total change, a compactly supported
equilibrium stress, has zero mean by testing against affine functions.

In a partition of the cell into \(N\) equal cubes, select the
\(\lfloor(1-s)N\rfloor\) with the largest inactive proportions.
As the dyadic partition is refined their union tends in measure to
the inactive set. Indeed the cube averages of its indicator tend to
the indicator in \(L^1\). Thresholding at \(1/2\) first gives the
same conclusion, and changes the selected count by \(o(N)\).
For fixed \(a\), the joint squared norm of the supported fields in
the selected cubes tends to zero. Clear their central cubes of side
\((1-2a)h\) and mark them as holes. The additional volume \(\Delta\) of holes
needed to reach \(1-f\) is at most \(s-f+2da+N^{-1}\).
Among the other cubes take the \(\lceil2^dN\Delta\rceil\) with least
joint squared norm. For small \(a,s-f\) and large \(N\) there are
enough, and their total squared norm is at most
\((2^d\Delta+N^{-1})\sum_j\|y_j\|_2^2/s\).
Clear their central half-side cubes and mark them as holes too.
Make the complement of all holes elastic. Refill up to \(2^{-d}/N\)
of hole volume by a slab to get exactly \(f\). The whole change of
the dual fields tends to zero, first refining at fixed \(a\), then
taking \(a\downarrow0\) and \(s-f\to0\).
The newly elastic part outside the original elastic phase has volume
at most the misclassified volume, \(2da\), and the refilled excess.
For the given primal trials choose \(a\) and the refinement separately
at each index so that their original squared norms on this added
volume are arbitrarily small (absolute continuity). They themselves
are unchanged. A fixed affine change of the cube lattice gives the
same statement on a parallelepiped, with constants depending on its
shape.
\phantomsection\label{proofend:p8:volume}\end{proof}

\begin{lemma}[Substitution]\label{p8:substitution}
In either model write \(G_g\) for its periodic response closure at
elastic fraction \(g\), including \(G_0=\{0\}\) and \(G_1=\{D\}\).
Partition a unit-volume periodic cell into finitely many regions
\(R_i\) with null boundaries. For each \(i\) take \(0\le f_i\le1\)
and \(A_i\in G_{f_i}\) with \(A_i\succ0\) or \(A_i=0\).
Then the cell response with coefficient \(A_i\) on \(R_i\), defined
by minimizing \(\sum_i\int_{R_i}z:A_i z\) over primal fields at each
prescribed mean, belongs to \(G_{\sum_i |R_i|f_i}\).
For a simple laminate of normal \(n\) this also holds for two
arbitrary semidefinite children in their respective closures when
their common kernel meets \(L_n\) only at zero. For a single region
the conclusion holds for any actual periodic cell response, including
a nonzero singular response. Both the inserted cells and the outer
cell may use arbitrary full-rank period lattices. In particular, these
lattices give the same closure as the cubic torus at each fraction.
\end{lemma}
\begin{proof}\phantomsection\label{proof:p8:substitution}
Here are the field estimates for scale separation. For a fixed macro
loading use a nearly minimizing primal trial and the optimal supported
dual field of Lemma~\ref{p8:projection}, applied with the positive
coefficient on each active region. The same projection proof applies
using the piecewise positive energy inner product. On a zero region
the dual field vanishes. Tile almost all the interiors of the
regions by small disjoint equal lattice cubes. Leave a frame of
relative width \(a>0\) in each cube. Each macro field can be made
constant, equal to its cube average, throughout the concentric core,
with joint \(L^2\) error tending to zero at fixed \(a\). For a strain
subtract from its displacement the affine function with the average
gradient, choosing the constant term for mean-zero difference, and
cut off the difference to zero in the core. Poincar\'e bounds the
change by \(c_a\) times the gradient deviation on the cube. For
a stress subtract its average and use exactly the frame correction
in Lemma~\ref{p8:volume}. Its bound is \(c_a\) times the stress
deviation. Summed squared deviations tend to zero by convergence
of cube averages in \(L^2\). These changes preserve means.

In each core of constant macro fields \(Z,Y\), insert many whole
periods of a child with primal trial of mean \(Z\) and supported dual
trial of mean \(Y\). In a positive definite child the two costs can
be taken arbitrarily close to \(Z:A Z\) and \(Y:A^{-1}Y\) by
Lemma~\ref{p8:projection}, also when \(A\) is obtained as a positive
definite limit. In a zero child only the primal trial is needed because
the dual mean and field there are zero. At any fraction strictly
between zero and one an exact zero child exists. For void isolate the
solid in a concentric cube with a void collar. The displacement can
be made constant on that cube at zero cost. For rigid isolate the
elastic cube inside a connected rigid matrix. Any admissible displacement
is one rigid motion on this matrix in the periodized space, forcing
its symmetric mean to vanish. The assertion then follows from
Lemma~\ref{p8:projection}.

For the pasting, keep a second thin elastic frame in the constant-field
region. Cut the scaled displacement fluctuation to zero across it.
At a fixed cutoff the term from its derivative tends to zero in \(L^2\), and
the repeated gradient costs at most the frame volume times a
fixed child-field constant in the limit. For a repeated zero-mean stress
fluctuation \(s(kx)\), cut it off by \(\zeta\) and correct its
divergence in this second frame by the zero-traction construction
above. The force \(s(kx)\nabla\zeta\) annihilates rigid motions. It
tends weakly to zero in \(L^2\) and hence strongly in the dual of
\(H^1\) on the fixed frame, by compactness. The correction tends
strongly to zero. Both pasted fields retain their prescribed means.
Use full child periods in the part where a phase is extreme and
elastic material in the frames and uncovered sets. Align the inner
frames with child cells, or leave a further all-elastic strip of width
one child period. This gives an elastic surplus tending to zero.
First fix the child cells and
linear trial families on a load basis, then choose the frame fractions
so small that all their costs are small, then fix a fine macro grid
and finally take sufficiently many repeats. The first frame and
uncovered volume have arbitrarily small original macro energy by
absolute continuity (the first frame occupies at most \(2da\) of
the cell). The flattening estimates apply at that fixed \(a\).
Lemma~\ref{p8:volume} restores the fraction. Pairing the primal and
dual energy estimates against each macro load gives convergence to
the asserted tensor on a finite polarization basis.

For a simple layer with fractions \(1-\theta,\theta\) and child
responses \(A,B\), convexity and averaging along the layers give
the exact quadratic
\begin{equation}\label{p8:layer}
\min_{j\in L_n}(1-\theta)|F+\theta j|_A^2+
                  \theta|F-(1-\theta)j|_B^2 .
\end{equation}
The minimizer gives constant primal fields whose local dual means
are their images under the respective children. Stationarity makes
the dual jump perpendicular to \(L_n\), so they are admissible.
For fixed actual children even if singular, use precisely these
fields in the pasting just proved. Their local dual means belong
to the ranges by construction. For limiting children the quadratic
depends continuously on the matrices whenever its jump coefficient
is positive definite on \(L_n\), as in the hypothesis. Arbitrary
normals are limits of rational normals. For a single actual child
of response \(A\), take the constant macro fields \(F,AF\). Even at
singular \(A\) the required child trials exist by Lemma~\ref{p8:projection}.
Use fine whole period parallelepipeds in the core and leave every
cell meeting the boundary of the core in the elastic frame. At fixed
frame width these cut cells occupy vanishing volume. Periodic averaging
and the frame corrections above apply to any fixed lattice, including
a rotated or rescaled one. Using a cube as outer cell and restoring
the fraction gives cubic-periodic approximants. Interchanging the
two period lattices gives the reverse approximation.
\phantomsection\label{proofend:p8:substitution}\end{proof}
\subsection{A common directional moment}
\begin{lemma}\label{p8:upper}
Put \(q=(1-f)/f\) and
\[
 B(M)=D^{1/2}\bigl[I-q(M+qI)^{-1}\bigr]D^{1/2}.
\]
Every cell is dominated by one \(B(M)\), and every such
\(B(M)\) is attainable at \(f\).
\end{lemma}
\begin{proof}\phantomsection\label{proof:p8:upper}
Use normalized mean \(Z=D^{1/2}F\) and write \(g\) for a zero-mean
normalized primal fluctuation. If \(\eta\) indicates the extreme
phase and \(\phi=1-f\), then, for any constant \(p\), use
\(-|Z+g|^2\le 2p:(Z+g)+|p|^2\) in that phase. Minimizing the
resulting full-cell quadratic over \(g\) bounds the normalized
response at \(Z\) by
\[
 |Z|^2+2\phi p:Z+\phi|p|^2-\|\Pi(\eta p)\|_2^2 ,
\]
where \(\Pi\) is the projection onto normalized fluctuations.
Parseval gives \(\|\Pi(\eta p)\|_2^2=f\phi\,p:(I-M)p\)
with \(M=\sum_{k\ne0}|\widehat\eta(k)|^2 M_{k/|k|}/(f\phi)\),
a convex average since \(\eta^2=\eta\).
Minimizing over \(p\) gives the bound. Compactness passes it to the closure.

For recovery put \(M=\sum w_iM_{n_i}\). Laminate the previous
normalized response \(X\) in proportion \(t\) with the solid \(I\)
in proportion \(1-t\). Write \(P=I-M_n\) and \(H=I-X\succ0\).
Eliminating the normalized primal jump in \(\mathrm{ran}\,P\) from
\eqref{p8:layer} gives
\[
 I-X'=t\big[H+(1-t)HP
   (I-(1-t)PHP)^{-1}PH\big],
 \qquad t(I-X')^{-1}=H^{-1}-(1-t)P;
\]
the second equality follows by multiplication. Start with \(X=0\),
put \(v_0=1,\ v_i=1-f\sum_{j\le i}w_j\), and take
\(t_i=v_i/v_{i-1}\). Multiplying the inverse update by \(v_{i-1}\)
and summing gives \((1-f)(I-X_{\rm last})^{-1}=(1-f)I+fM\).
This is the stated \(B(M)\), and the final extreme fraction is \(1-f\).
Lemma~\ref{p8:substitution} realizes every step. The Schur complement
shows that \(T\preceq B(M)\) is precisely \eqref{p8:lmi}.
\phantomsection\label{proofend:p8:upper}\end{proof}

\subsection{Removing one response direction}
For \(A\succ0\) and a unit loading \(e\in V\) write
\[
 A_e=A-\frac{(Ae)\otimes(Ae)}{e:Ae}.
\]
It is the largest \(U\preceq A\) that is positive semidefinite and
annihilates \(e\), as one sees by minimizing \((F+t e):A(F+t e)\) over \(t\).
Its inverse on \(R=e^\perp\) is the compression of \(A^{-1}\).

\begin{lemma}[Void]\label{p8:short}
For each positive definite attainable stiffness \(A\) and each
unit strain \(e\), \(A_e\) is attainable at the same fraction.
\end{lemma}
\begin{proof}\phantomsection\label{proof:p8:short}
Use the porous cells of Lemma~\ref{p9:fields} in principal axes of
\(e\), rotating the whole cell back afterwards. They admit
displacements of mean strain \(e\) with
\(\int_G|e(u)|^2=o(1)\) on the solid \(G\), and for each \(H:e=0\)
equilibrium stresses of mean \(H\), supported in \(G\), with
\(\|\tau_H-H\|_2=o(1)|H|\). Both estimates are geometric.
Put material \(A\) in \(G\) and the zero child at fraction \(f\)
in its complement. The boundaries have zero measure.
Lemma~\ref{p8:substitution} puts each response in the original
two-phase closure at exactly \(f\). A uniform strain trial bounds
these responses above by \(A\), and the displacement trial forces every
limit \(U\) to annihilate \(e\). Conversely, for every \(H\in R\)
the stress costs \(\int_G\tau_H:A^{-1}\tau_H\to H:A^{-1}H\).
Thus \(F:UF\ge\sup_{H\in R}(2F:H-H:A^{-1}H)=F:A_e F\).
The opposite order follows from \(Ue=0,\ U\preceq A\).
\phantomsection\label{proofend:p8:short}\end{proof}
\begin{lemma}[Perfectly rigid phase]\label{p8:rigshort}
For each positive definite attainable compliance \(A\) and each
unit stress \(e\), \(A_e\) is attainable at the same fraction.
\end{lemma}
\begin{proof}\phantomsection\label{proof:p8:rigshort}
We construct a thin rigid framework forcing every supported mean
strain to belong to \(R=e^\perp\). For every \(E\in R\) it admits
\(u=Ex+w\) rigid on the framework with
\(\|\nabla w\|_{L^2}\to0\). Put the zero-compliance child on this
framework and \(A\) outside, both at fraction \(f\). Lemma~\ref{p8:substitution}
applies. The constant stress trial gives \(U\preceq A\) for any limiting
compliance. Every approximating compliance annihilates \(e\) by
Lemma~\ref{p8:projection}. The displayed strain trials cost
\(E:A^{-1}E+o(1)\) and imply \(U\succeq A_e\) by the same conjugation
as in Lemma~\ref{p8:short}.

\emph{Plane.}
Write \(J(x,y)=(-y,x)\) and \(h=1/N\), with \(N\ge 16\) even.
At each square-lattice vertex put a rigid square hub of side \(2h\)
with four rectangular arms of half-width \(h/16\) reaching distance
\(4h\) along the coordinate axes. On each positive coordinate link
\(t\) put \(n=Jt\). In local coordinates based at the hub place
squares at \(c_{ij}=h(it+jn)\) for \(4\le i\le N-4,\ j=0,1\).
Their four vertices are \(c_{ij}\pm r_\sigma,c_{ij}\pm q_\sigma\)
where \(\sigma=(-1)^{i+j}\) and
\[
 r_\sigma=h(t/2+\sigma p n),\qquad q_\sigma=h(n/2-\sigma p t).
\]
Take \(p_x=1/4\) on horizontal links and \(p_y=\pm1/4\) on vertical
links. Each square touches its box of side \(h\) only at the four
specified corners, and adjoining squares touch there. The two end
squares of row zero overlap the opposing arms in open sets. The
other squares are disjoint from the arms. The different strips meet
only at the hubs. The area of their \(h\)-neighborhood is \(O(h)\).

A full-cell \(H^1\) displacement on touching squares gives their
rigid motions equal velocities at the corner. Indeed a nonzero
mismatch forces, on almost every circle of radius \(s\) about the
corner, a change bounded away from zero between two occupied arcs.
Angular Cauchy-Schwarz then costs at least \(c\int_0^\epsilon ds/s\).
On a cycle of four squares the successive contacts, in units of \(h\)
relative to the first center, are
\((1/2,p),(1+p,1/2),(1/2,1-p),(-p,1/2)\), with \(p\) negated for odd
initial parity. These are a nondegenerate parallelogram.
If the consecutive sides are \(d_i\), summing velocity differences gives
\(\sum_i\omega_i Jd_i=0\), so opposite squares have equal rotations.
The adjacent cycles propagate this equality through both rows.
Thus the rotations have form \(m+\sigma k\) on each strip. The
difference of center velocities on a neighboring step \(d\) is
\((mJ-2pkI)d\), by substituting the shared offsets \(r_\sigma\)
or \(q_\sigma\). Conversely, these differences and spins satisfy every
contact. The two welded end spins equal the periodic hub spin
\(\omega=m+k\). The end centers are \(\ell t\) apart, where
\(\ell=1-8h\). Subtracting the rigid hub offsets at the ends gives
\[
 (E-\omega J)t=\ell k(-2pI-J)t .
\]
The two values of \(k\) on the perpendicular links may differ.
These equations reduce exactly to
\(2E_{12}=E_{11}/(2p_x)-E_{22}/(2p_y)\), since they prescribe
\(E_{12}-\omega=E_{11}/(2p_x)\) and
\(E_{12}+\omega=-E_{22}/(2p_y)\).
Conversely these equations determine \(\omega,k,m\). Assigning the
first center its hub velocity and propagating the differences gives
the required weld at the far end. The annihilating stresses are
therefore proportional to the matrices with
\((e_{11},e_{22},e_{12})=(-2,2,1)\) and \((-2,-2,1)\).
These have opposite determinant signs.

For an allowed \(E\) the spins are bounded by \(c|E|\) independent
of \(N\). With \(H=mJ-2pkI\) on a strip, the link equation gives
\((H-E)t=(1-\ell)(H-\omega J)t=O(h)|E|\).
All centers are within \(h\) of its axis. Subtracting \(Ex\) from
the prescribed rigid displacements thus gives a periodic function
of amplitude \(O(h)|E|\), Lipschitz on the pieces with common
constant \(c|E|\). At a shared corner this bound across pieces
follows from equality there and the positive intervening angles.
Welded pieces have the same motion. Every noncontacting pair not
separated by at least \(h\) belongs to finitely many scaled local
configurations of positive separation. Prescribe this correction
to be zero outside the \(h\)-neighborhood of the pieces.
The same Lipschitz bound holds up to a fixed factor. Extend each
component by the infimum of prescribed values plus their Lipschitz
constant times the torus distance. The gradient bound and
\(O(h)\) support give the required \(L^2\) estimate.

\emph{Space.}
Place hubs at \(\mathbb Z^3\) and its translate by \(c=(t,t,t)\).
Connect \(c\) to \(q_0=0,q_j=e_j\) for \(1\le j\le3\), and repeat.
Put \(a_i=q_i-c\) and use \(t=1/4\) or \(1/2\).
The four open edges in a cube meet only at \(c\), and their interiors
lie in the open cube. The first and last quarters of each edge are
arms of the respective hubs, the middle half is a separate bar.
Thicken these segments with radius \(r\). For an arm from a vertex
\(v\) to contact \(p\) take the convex hull of the balls at \(v\)
and \(p-rn\), both of radius \(r\), where \(n\) points toward \(p\).
For a central bar with contacts \(p,p'\) take the convex hull of
the balls centered \(r\) inward from its endpoints. Adjacent arms
of a hub overlap in its vertex ball. For small \(r\) the only
contacts between distinct bodies are the intended tangencies.

At a contact the two spherical caps bound a gap of axial width
\(g(s)=2(r-\sqrt{r^2-s^2})\asymp s^2/r\) at transverse radius
\(s\). A mismatch in contact velocity would cost at least
\(c\int_0^\epsilon s/g(s)\,ds=\infty\) by axial Cauchy-Schwarz.
Conversely two rigid motions agreeing at the contact differ by
\(O(s)\) in this gap. Interpolate them linearly in the axial
coordinate between the two cap surfaces. The axial derivative is
\(O(1+r/s)\), and the transverse derivatives are bounded since the
cap slopes are \(O(s/r)\). Their squares are integrable against
\(s g(s)\,ds\). The resulting field is bounded and extends in
\(H^1\) across the point (a cutoff on a shrinking ball costs
\(O(\epsilon)\) in squared gradient norm). Use disjoint contact
neighborhoods and smooth cutoffs away from them, where different
bodies have positive separation, to obtain a bounded full-cell
extension.

The rotational velocity at an arm endpoint is perpendicular to the
edge axis. Taking axial projections
of the contact equations on each bar gives
\(a_i\cdot(Ea_i+b)=0\), with \(b\) the difference between the
two hub translation correctors. To see that they suffice, take zero hub spins
and fit a bar motion to the two endpoint velocities, whose difference
is perpendicular to the bar. The sole linear relation among the
four \(a_i\) has coefficients \((1-3t,t,t,t)\). Hence solvability
in \(b\) amounts to \(e:E=0\) with
\(e=\operatorname{diag}(c)-cc^t\). Its eigenvalues are
\(t,t,t(1-3t)\), giving definite and indefinite cases.

To obtain vanishing energy cost, fix a small reference radius
\(r_0\). The above extension gives for an allowed strain a bounded
periodic corrector \(w\in H^1\) rigidifying the reference bodies.
Every smaller body lies in its reference body because each moving end
ball stays internally tangent at the same contact. Let \(s\) be
distance to the fixed union of axes and multiply \(w\) by the cutoff
that is one for \(s\le r\), zero for \(s\ge\sqrt r\) and linear
in \(\log s\) between. Its squared gradient integral is
\(O(1/|\log r|)\), since cylindrical shells of each segment and its
two caps have volume density at most \(c(s+s^2)\).
Thus the correction gradient tends to zero in \(L^2\), by boundedness
of \(w\) and absolute continuity of \(\int|\nabla w|^2\) on the
shrinking support.

In each dimension the change \(x'=Lx,\ u'(x')=L^{-t}u(x)\) preserves
rigidity, takes strain to \(L^{-t}EL^{-1}\) and annihilating stress
to \(LeL^t\). The estimates persist for fixed nonsingular \(L\).
By inertia the two templates, up to sign, cover every invertible
stress. Approximate singular stresses by invertible ones and use
continuity of \(A_e\) and closedness.
\phantomsection\label{proofend:p8:rigshort}\end{proof}

\begin{proof}[Completion of Theorem~\ref{p8:main}]\phantomsection\label{proof:p8:main}
From \(A\succ0\) one can now reach
\((A^{-1}+k e\otimes e)^{-1}\) for every \(k\ge0\). Choose a
normal for which \(e\notin L_n\). Such a normal exists since the
intersection of all \(L_n\) is zero. Indeed, in the strain case
\(x^tEx\) would have to be divisible by every \(n\cdot x\), and
in the stress case \(En=0\) for every \(n\). Laminate \(A,A_e\)
at varying proportions. The result varies continuously from \(A\)
to \(A_e\), including the endpoints since both matrices are positive
on \(L_n\). It belongs to the order interval \([A_e,A]\) by
monotonicity of \eqref{p8:layer}. This interval is the line segment
between its ends because their difference has rank one. Each displayed
inverse lies on that segment by the rank-one inverse formula.

Now if \(0\prec T\preceq B(M)\), start with \(B(M)\) and
add separately the rank-one summands of \(T^{-1}-B(M)^{-1}\) to its
inverse. For a singular \(T\) or \(B(M)\), average the moment with
a fraction \(\epsilon\) of the uniform spherical moment \(M_0\).
It is positive definite since the intersection of all normalized
jump spaces is zero. The inverse is operator convex on positive
matrices, as follows from the Schur complement, so \(B\) is
concave and the new bound dominates
\((1-\epsilon)T+\epsilon B(M_0)\succ0\). Take \(\epsilon\downarrow0\).
This proves sufficiency in every rank and completes the proof.
\phantomsection\label{proofend:p8:main}\end{proof}
\subsection{Explicit moment variables}
\begin{corollary}[Plane]\label{p8:planar}
Let \(k=\lambda+\mu\), \(h=I_2/\sqrt2\), \(P_h=h\otimes h\) and
\(P_d=I_V-P_h\). The moment variable in \eqref{p8:lmi} is
\(M=W\) for void, \(M=I_V-W\) for rigid, where
\[
 W\succeq0,\qquad h:Wh=\frac{\mu}{k+\mu},\qquad
 \operatorname{Tr}(P_dW)=\frac{k}{k+\mu}.
\]
For a void stiffness \(T\succeq0\) require \(T\prec C\), put
\(P=T+T(C-T)^{-1}T\), and let \(Y=4k\mu/(k+\mu)\).
The two remaining scalar inequalities are
\[
 h:Ph\le\frac{fY}{2(1-f)},\qquad
 \operatorname{Tr}(P_dP)\le\frac{fY}{2(1-f)} .
\]
In particular for a nonzero such candidate, writing
\(m=\max\{h:Ph,\operatorname{Tr}(P_dP)\}\), its least solid fraction
among \(0<f<1\) is \(2m/(Y+2m)\).
\end{corollary}
\begin{proof}\phantomsection\label{proof:p8:planar}
For \(t\perp n\) of unit length, \(C^{-1/2}N_n\) is spanned by
\(u_t=\sqrt{Y}C^{-1/2}tt^t\). In a hydrostatic and orthonormal
deviatoric basis this has coordinates
\((a,b\cos 2\alpha,b\sin 2\alpha)\), where \(\alpha\) is the angle
of \(t\), \(a^2=\mu/(k+\mu)\), \(b^2=1-a^2\).
Thus the two trace conditions hold on the convex hull of
\(u_t\otimes u_t\). Conversely factor a passing \(W=ZZ^t\) and
choose an orthonormal column basis making
\(Z^t(b^2P_h-a^2P_d)Z\) have zero diagonal. Since its trace is zero,
there is a unit zero of its quadratic form, and one repeats in
the orthogonal complement. Each nonzero column normalized and signed
is then some \(u_t\), and their squared norms sum to one.
The Schur complement of \eqref{p8:lmi} for void is
\(q X(I-X)^{-1}\preceq W\) for \(X=C^{-1/2}TC^{-1/2}\) and
\(q=(1-f)/f\). The two trace upper bounds suffice to extend this
positive matrix to a moment by adding the missing hydrostatic mass along
\(P_h\) and splitting the missing deviatoric mass equally along \(P_d\).
Since \(X(I-X)^{-1}=C^{-1/2}PC^{-1/2}\), rescaling gives the
stated formulas.
\phantomsection\label{proofend:p8:planar}\end{proof}

\begin{corollary}[Space]\label{p8:spatial}
Let \(W\) be a fully symmetric fourth-order tensor, positive
semidefinite on \(V\), with \(\sum_{ij}W_{iijj}=1\).
Put \(N_{ij}=\sum_kW_{ijkk}\) and define \(Q(W)\) by
\[
 G:Q(W)G=\frac{\operatorname{tr}(Z^2N)}{\mu}
 -\left(\frac1\mu-\frac1{\lambda+2\mu}\right) Z:WZ,
 \qquad Z=C^{1/2}G.
\]
Then \eqref{p8:lmi} with \(M=I-Q(W)\) for void or \(M=Q(W)\)
for rigidity, and the stated constraints on \(W\), is a semidefinite
membership test with fifteen scalar moment variables.
\end{corollary}
\begin{proof}\phantomsection\label{proof:p8:spatial}
For \(W=n^{\otimes4}\), set \(L a=a\odot n\). Then
\(L^t Z=Zn\), \(L^t C L=\mu I+(\lambda+\mu)nn^t\) and
\(C^{1/2}L(L^tCL)^{-1}L^tC^{1/2}\) is the orthogonal projector
onto \(C^{1/2}J_n\). Substituting its inverse gives the displayed
quadratic form. It remains to show that the stated constraints
are exactly the hull of the \(n^{\otimes4}\).

The dual cone to the fourth powers consists of the nonnegative ternary
quartics. By Hilbert's ternary quartic theorem each is a sum of squares
of real quadratic forms \cite[abstract, p.~1]{p8quartic}.
Positivity of \(W\) on all quadratic squares
is therefore precisely membership in the closed cone of fourth
powers, by finite-dimensional separation. The cone is closed since
the convex hull of unit fourth powers is compact with double trace
one. Double trace fixes the mass.
\phantomsection\label{proofend:p8:spatial}\end{proof}

\subsection{Limits of finite contrasts}

\begin{corollary}\label{p8:contrast}
The stiffness closure with true void equals the set of limits obtained
by replacing void with a positive definite phase whose stiffness tends
to zero in operator norm, allowing the geometry to vary at fixed \(f\).
The analogous statement holds for the compliance closure with perfect
rigidity when the other phase's compliance tends to zero.
\end{corollary}
\begin{proof}\phantomsection\label{proof:p8:contrast}
Normalize by \(D\) as in Lemma~\ref{p8:upper}, so the regularized
phase costs at most \(\delta\) times the squared normalized loading,
with \(\delta\to0\). For each geometry use the same trial
\(g=-\Pi(\eta p)\) as in that proof at the minimizing
\(p=-[(1-f)I+fM]^{-1}Z\). Since
\(\|g\|_2\le |p|\le |Z|/(1-f)\), its extra energy is at most
\(\delta\langle|Z+g|^2\rangle\le\delta[1+(1-f)^{-2}]|Z|^2\).
Compactness of the moments and Theorem~\ref{p8:main} put every
regularized limit in the extreme-phase closure. Conversely, on a
fixed cell its regularized response dominates its extreme response.
Each fixed near-minimizing primal trial has finite full-cell squared
norm, so its extra energy tends to zero. A polarization basis gives
convergence of cell tensors. Diagonal selection over cells proves
the reverse inclusion. The pure fractions follow directly.
\phantomsection\label{proofend:p8:contrast}\end{proof}

\subsection{Exact void and perfect rigidity}

\begin{theorem}\label{p14s:rect}
Suppose an isotropic elastic solid of bulk and shear moduli \(k,m>0\)
occupies a prescribed fraction \(0<f<1\). Define
\[
 (b_2,s_2)=\left(m,\frac{km}{k+2m}\right),\qquad
 (b_3,s_3)=\left(\frac{4m}3,\frac{m(9k+8m)}{6(k+2m)}\right).
\]
In dimension \(d=2,3\) the closed set of isotropic effective pairs
\((k_*,m_*)\) for this solid and true void is
\[
 \left[0,\frac{f k b_d}{b_d+(1-f)k}\right]\times
 \left[0,\frac{f m s_d}{s_d+(1-f)m}\right].
\]
With a perfectly rigid second phase it is
\[
 \left[\frac{k+(1-f)b_d}{f},\infty\right]\times
 \left[\frac{m+(1-f)s_d}{f},\infty\right],
\]
where each interval includes its infinite endpoint in compliance
coordinates. Attainment in both statements is in the fixed-fraction
tensor closure.
\end{theorem}
\begin{proof}\phantomsection\label{proof:p14s:rect}
The response closure is rotation invariant because \(D\) is isotropic
and rotating a moment rotates its normals. Since \eqref{p8:lmi} is
affine in its two variables, any passing isotropic \(T\) also passes
with the rotation average of its moment, the moment of a uniform
normal. By Lemma~\ref{p8:upper} this is equivalent to \(T\preceq B(M)\)
for that moment. Theorem~\ref{p8:main} gives attainment of every
positive semidefinite \(T\) satisfying this inequality.

For this uniform normal the mean projector onto \(C^{1/2}J_n\)
has hydrostatic eigenvalue \(p=k/(\lambda+2m)\). For
\(h=I_d/\sqrt d\), take \(Z=C^{1/2}h=\sqrt k I_d\) in the projection
formula of Corollary~\ref{p8:spatial}. The formula applies in both
dimensions because the acoustic matrix is \(mI+(\lambda+m)nn^t\). Each projector has trace \(d\), so
by isotropy its other eigenvalue is \(r=(d-p)/(d(d+1)/2-1)\).
Substitution of \(\lambda=k-2m/d\) gives
\(b_d=k(1/p-1)\), \(s_d=m(1/r-1)\).
For rigid material \(M\) has the eigenvalues \(p,r\), and for void they
are \(1-p,1-r\). Applying \(t\mapsto t/(t+(1-f)/f)\) to those
eigenvalues and rescaling by \(D\) gives the two rectangles.
The isotropic stiffness eigenvalues are \(d k_*,2m_*\).
The compliance eigenvalues are their reciprocals.
\phantomsection\label{proofend:p14s:rect}\end{proof}

\begin{corollary}[Exactly incompressible solid and void]\label{p14s:incompressible}
Let \(S\) be a solid of shear modulus \(m>0\) and fraction
\(0<f<1\) in a unit-volume periodic cell, with true void elsewhere.
Its exactly incompressible cell stiffness has quadratic form
\[
 Q_\infty(E)=\inf\left\{2m\int_S|e(u)|^2:
   u-Ex\in H^1_{\rm per},\ \operatorname{div}u=0\ \text{a.e. on }S\right\}.
\]
As \(S\) varies at fixed \(f\), the isotropic section of the closure
of these tensors, in bulk and shear moduli, is
\[
 \begin{cases}
 [0,fm/(1-f)]\times[0,fm/(2-f)],& d=2,\\[2pt]
 [0,4fm/(3(1-f))]\times[0,3fm/(5-2f)],&d=3 .
 \end{cases}
\]
\end{corollary}
\begin{proof}\phantomsection\label{proof:p14s:incompressible}
Write \(Q_K\) for the cell energy on the same geometry when the
solid has bulk modulus \(K\). For any affine-periodic displacement
\(u\) let \(a=\langle\mathbf1_S\operatorname{div}u\rangle\)
and solve on the torus
\[
 \Delta\phi=\mathbf1_S\operatorname{div}u-
        \mathbf1_{S^c}a/(1-f)=g.
\]
The right side has mean zero and
\(\|\nabla^2\phi\|_2^2=\|g\|_2^2
\le (1-f)^{-1}\int_S|\operatorname{div}u|^2\) by Fourier series
and Cauchy-Schwarz. The displacement \(u-\nabla\phi\) has the
same mean strain and zero divergence on \(S\). Applying this to
\(u=Ex\) first gives \(Q_K(E)\le Q_\infty(E)\le c_0(f,m,d)|E|^2\)
uniformly in the geometry. Next take a trial of finite-\(K\) energy
at most \(Q_K(E)+\epsilon\). Its divergence on \(S\) has squared
norm at most \((Q_K(E)+\epsilon)/K\). By the triangle inequality
for the corrected deviatoric strain,
\[
 Q_\infty(E)^{1/2}\le
 (Q_K(E)+\epsilon)^{1/2}
 \left(1+\sqrt{\frac{2m}{(1-f)K}}\right).
\]
It follows that \(Q_\infty\) is a finite quadratic form, being the
pointwise limit of the \(Q_K\). Their tensors differ in operator norm by
at most \(c_1(f,m,d)/\sqrt K\) for \(K\ge1\), independently of the
cell, even when the solid has rough boundary.

If incompressible cell tensors converge to an isotropic limit,
replace them on the same geometries by their finite-\(K\) tensors
and average those tensors over rotations. These averages belong
to the finite-\(K\) closure by its convexity in Theorem~\ref{p8:main}
and its rotation invariance proved in Theorem~\ref{p14s:rect}. Averaging contracts operator norm, so the
isotropic limit lies within \(c_1(f,m,d)/\sqrt K\) of the rectangle
in Theorem~\ref{p14s:rect}, in stiffness coordinates. Conversely,
approximate any point in the limiting rectangle by points of
the finite-\(K\) rectangles and then by cells of exact fraction
\(f\). Taking the hard constraint on these very cells changes their
tensors by a quantity tending uniformly to zero. It remains to
take \(k=K\to\infty\) in the two endpoints of
Theorem~\ref{p14s:rect}, where \(s_2\to m\) and \(s_3\to3m/2\).
\phantomsection\label{proofend:p14s:incompressible}\end{proof}

\paragraph{Example at one-half solid.}
\phantomsection\label{item:p8-example}
For unit three-dimensional stiffness on symmetric strains, the
constituent has $k=1/3$ and $m=1/2$. The shear shift in
Theorem~\ref{p14s:rect} is $s_3=7/16$. At $f=1/2$ the rectangle
therefore reaches zero bulk with
\[
 \mu_*=\frac{f m s_3}{s_3+(1-f)m}=\frac{7}{44}.
\]
All five trace-free stiffness eigenvalues are $2\mu_*=7/22$ there.
This endpoint belongs to the fixed-fraction closure.

\section{9 and 10. Full limiting shear at vanishing bulk}\label{problem-9}
\Milton{9}\Milton{10}

\begin{theorem}\label{p9:sharp}
Fix $d\ge2$ and finite positive moduli $\kappa_s,\mu_s$ of an isotropic
solid. Allow any period lattice and any measurable solid set in its
period torus, with void elsewhere, and require the effective tensor
to be isotropic. Over all such cells and all solid fractions,
\[
 \sup(\mu_*-c\kappa_*)=\mu_s\qquad(c>0).
\]
No individual cell attains the supremum. For all sufficiently large
integers $M$ there is a cell, independent of $c$ and isotropic at that
index, with $0\le\kappa_*\le C/M$ and $0\le\mu_s-\mu_*\le C/M$,
where $C$ depends only on $d,\kappa_s,\mu_s$. Its solid fraction tends to one.
\end{theorem}

Positivity and affine strain trials give Milton's separate bounds (4.1), $0\le\kappa_*\le\kappa_s$ and
$0\le\mu_*\le\mu_s$.
The sequence answers Problem 10 and gives the second alternative in Problem 9
for every $c>0$. The planar construction prescribes the six grain orientations.
In higher dimensions the last orientation is a root of a continuous scalar
equation whose existence we prove below.

We use the whole-cell displacement formulation on a period torus, with normalized
averages and $L^2$ norms. Write $e(u)=(\nabla u+\nabla u^T)/2$,
$W_s(E)=E:C_sE=2\mu_s|\operatorname{dev} E|^2+\kappa_s(\operatorname{tr} E)^2$,
and $S_s=C_s^{-1}$. For solid $G$ of indicator $\chi$ the doubled energy is
\begin{equation}\label{p9:energy}
 Q_G(E)=\inf_{w\in H^1_{\rm per}}\langle\chi W_s(E+e(w))\rangle=E:C_GE.
\end{equation}
A stress is admissible if it is symmetric, periodic, in $L^2$, zero in the pores
and distributionally divergence free. The exact dual formula is
\begin{equation}\label{p9:dual}
 Q_G(E)=\max_{\tau\ {\rm admissible}}
 \{2 E:\langle\tau\rangle-\langle\tau:S_s\tau\rangle\}.
\end{equation}
The Hilbert-space projection argument of Lemma~\ref{p8:projection}, valid in every dimension, gives the maximizing stress $\tau_E$ of mean $C_G E$, at which the minimum complementary energy at that mean equals $E:C_G E$. Completing the square bounds the right side of \eqref{p9:dual} by $Q_G(E)$ for every admissible stress. We use the displacement and traction estimates of Lemma~\ref{basic:korn} on the solid bodies.

\subsection*{A block joined to its neighbours by bending rods}

The geometric tests used in Problem 8 follow from rod bending and equilibrium on the bodies of a single cell.
\begin{lemma}\label{p9:fields}
Let $F=\operatorname{diag}(s_1,\ldots,s_d)$ with $|s_i|\le1$ and $d\ge2$.
For every integer $M\ge1$ there is a polyhedral solid $G$ on the unit
cube torus, of fraction $1-O_d(M^{-1})$, with both of the following properties.
\begin{enumerate}
\item There is a corrector $w\in H^1_{\rm per}$ on the whole torus,
including the pores, such that $u=Fx+w$ satisfies
$\langle\chi|e(u)|^2\rangle\le C_d/M$.
\item For every $T=T^T$ with $T:F=0$ there is an admissible stress, linear in $T$,
of mean $T$ with $\|\tau-T\|_2^2\le C_d|T|^2/M$.
\end{enumerate}
\end{lemma}

\begin{proof}\phantomsection\label{proof:p9:fields}
The stressed rods carry axial loads through the gaps between blocks.
The small change of direction across each collar leaves a small residual
force that its connections equilibrate. The same rods can bend as the
collars translate, accommodating the prescribed block increments $Fe_i$.
All constructions use the block $B=(0,L)^d$ with
$h=(M+c)^{-1}$, $L=Mh$, $3H=1-L=ch$ and $c>0$ fixed small.
A face normal is $n=e_i$ and $z=x_i-L$ is height across its gap, $0<z<3H$.
Tile the positive face by cubes of side $h$. On each tile we will specify
small tangential cubes called patches. A patch at $z=0$ sweeps in direction
$m$ to height $2H$ and from there in direction $l$ to height $3H$,
where $m\cdot n=l\cdot n=1$. In the first and last thirds retain only the
central half-side cube of each member of an $N\times\cdots\times N$
equal subdivision of the patch, transported along the sweep. These prisms are
the rods. $N$ is even and will be chosen last. The middle third lies in a
solid body called a collar. Between the block and collar, and between the
collar and next block, only the rod end faces are welded. Their total
retained area is $2^{1-d}$ times the patch area, independent of $N$.
A bundle carrying scalar force $q$ in direction $b$ has stress
$(q/A)b\otimes b$, with $A$ its total retained area. This gives zero traction
on the sides and force $\pm qb$ at the ends. The same stress may continue
through a full collar along the straight extensions of the retained rods.

Two elementary estimates will be used throughout. If the prescribed end
translations of a rod differ by $D\perp b$, take $X_0$ at its lower-end
center and write $X=X_0+\xi+zb$ with $\xi$ in the end face and $0\le z\le H'$
(the rod's height). For $r=b\cdot(X-X_0)/(H'|b|^2)$ choose a smooth $\psi$
equal to zero for $r\le1/4$ and one for $r\ge3/4$. Set
\[
 u(X)=U_0+\psi(r)D-\psi'(r)\frac{D\cdot(X-X_0)}{H'|b|^2} b .
\]
The terms containing $\psi'$ in the gradient are skew together. Since
$D\cdot(X-X_0)=D\cdot\xi$, the strain is bounded by
$C|D|\sup|\xi|/(H')^2$ when $|b|\ge1$. At fixed heights and translations
we can choose $N$ large enough to make the strain at most one and the
field constant near both entire ends. With constant translations on bodies
this defines a full-cell $H^1$ displacement because each rod interior and each body
has a smooth local extension, the constant extensions agree near welds, and
all nonincident pieces are separated at fixed parameters. Subtract $Fx$
and glue local extensions by a partition of unity on the torus.
Only rods strain, on a volume at most $1-L^d=O_d(h)$.

For a connected bounded Lipschitz body $P$, traction $g\in L^2(\partial P)$
of zero total force and moment has a symmetric equilibrium stress with
\begin{equation}\label{p9:neumann}
 \|\sigma\|_{L^2(P)}\le C_P\|g\|_{H^{-1/2}(\partial P)}.
\end{equation}
Indeed $v\mapsto\int_{\partial P}g\cdot v$ annihilates rigid motions,
and the trace and Korn inequalities bound it by the right side times
$\|e(v)\|_2$. Its Riesz representative for the symmetric-gradient inner
product modulo rigid motions gives $\sigma=e(v_g)$.

\smallskip\noindent\emph{Planar gaps.}
Here take $c=1/1000$, let $t=e_{3-i}$ on face $i$ and write $y=x\cdot t$.
In each length-$h$ tile put two patches of width $p=h/24$ centered at the
quarter and three-quarter marks, with signs $\epsilon=+1,-1$ respectively.
The directions are
\[
 m=n+\epsilon t,\qquad l=m+\eta s_i t,\qquad 0<\eta\le1 .
\]
Number patches along each face from $1$ to $2M$ and denote their centers by $y_r$.
Use the full rectangles $H\le z\le2H$, $|y-y_r|\le2p$ as collars.
Between the facing ends of successive collars run connector rods in direction
$t$ at heights from $3H/2-H/8$ to $3H/2+H/8$, retaining the central half
of each of $N$ equal height intervals. From the last collar on each face
extend this same height interval as a solid arm up to $y=L+3H/2+H/8$.
The two arms overlap in a square and make the last two collars one corner body.
Include all lattice translates and fill the rest with void. The rod sweeps
stay within $p/2+6H<2p$ of their center in $y$ and distinct centers are
$h/2$ apart. All collars lie strictly between $0$ and $L$ in $y$. The
connector ends lie on the collar sides away from its corners, and the arms
meet only in the terminal square beyond both faces. Thus the prescribed welds
are the only contacts. The ordinary collars are rectangles and the corner
is a connected polygon with Lipschitz boundary, independent of $N,\eta$.

Give the block in cell $R$ the translation $FR$. Give its collar of sign
$\epsilon$ the additional translation $(t-\epsilon n)/\eta$. At both last
collars this is $(e_1+e_2)/\eta$. On the first rods
$(t-\epsilon n)\cdot m=0$, on the last
$[s_in-(t-\epsilon n)/\eta]\cdot l=0$, and consecutive collar translations
differ in the $n$ direction. The bending test above therefore works on
all rods, including connectors.

Write $p_i=T_{ii}$, $t_i=t\cdot Tn$, and choose
$q_{i\epsilon}=(p_i+\epsilon(t_i-\eta s_ip_i))/2$.
The last bundle takes load $q_{i\epsilon}/M$, so the two last bundles of one
tile deliver $Tn/M$. Start with the same scalar load on each first bundle,
and continue its stress straight along $m$ across the collar. This reference
stress matches the first end. On its upper end the additional force needed
is $b_{ir}t$ with $b_{ir}=\eta s_i q_{i\epsilon}/M$. Give the connector
from $r$ to $r+1$ the scalar force $-\sum_{j\le r} b_{ij}$. At every
ordinary collar this cancels the residual force. On the corner the two
faces leave the net force
$\eta(s_1p_1 e_2+s_2p_2 e_1)$, parallel to $e_1-e_2$ because $F:T=0$.
Cancel it by changing the first-bundle load at the last collar of face one,
whose direction is $e_1-e_2$. All remaining moments, computed relative to
the reference stresses, are $O_M(\eta)|T|$ because the residual forces have that
size and act at bounded positions. On a first bundle add opposite equal scalar
loads to its two halves. Their retained centroids differ by $(p/2)t$, giving
a pure moment of either sign since $t$ and $m$ are independent.
Use one such bundle per ordinary collar and one on the corner to cancel its
moment. The added densities are $O_M(\eta)|T|$ independently of even $N$,
because retained half-areas and centroids are independent of $N$.
The residual traction now includes the first-bundle adjustments on the lower
collar ends. It balances on each body. Its $L^2$ norm is $O_M(\eta)|T|$ because
the upper-end mismatch has density $O(\eta)|T|$, and the connector and
first-bundle areas at fixed $M$ are independent of $N$.
Apply \eqref{p9:neumann} and add the reference stress. These body corrections
and all extra rod densities can thus be made arbitrarily small at fixed $M$
by decreasing $\eta$ before choosing $N$.

\smallskip\noindent\emph{Gaps in dimensions $d\ge3$.}
Two connected frames have translation directions
$w_1=(1,\ldots,1)$ and $w_2=(1,2,\ldots,d)$. On face $i$ project onto
$n^\perp$ by $P_i$. For $a=1,2$ let
\[
 v_{ia}=P_iw_a/|P_iw_a|^2,\qquad V_{ia}=n^\perp\cap w_a^\perp,\qquad
 m_{ia0}=n-(n\cdot w_a)v_{ia}.
\]
For each member $t_{iar}$ of a fixed basis of $V_{ia}$ put
$m_{iar}=m_{ia0}+t_{iar}$. The resulting $d-1$ directions for each $a$
span $w_a^\perp$ and have unit normal component. Use
$l_{iar}=m_{iar}+\eta s_i v_{ia}$ as their last directions.
Since $P_iw_1,P_iw_2$ are independent (at least two unequal coordinates
remain), $V_{i1}+V_{i2}=n^\perp$.

Put $K=2(d-1)$ patches of side $b=h/(100K)$ in each tile, one per direction.
In tangential coordinates scaled by $h$ their centers are
$(j/(K+1),1/2,\ldots,1/2)$ for $1\le j\le K$. A patch's protective prism
is the gap over the cube of tangential sup-radius $4b$ about its center.
These prisms are disjoint and strictly inside their respective tiles, at
distance greater than $b$ from the tile edges. Choose $c$ so small that
$c(1+\max_{iar}(|m_{iar}|+|v_{ia}|))<1/(10^4 K)$.
The sweeps then stay within $b/2+b/100$ of their center in each tangential
coordinate. Use the full swept middle third as each collar.

To connect all collars of label $a$, use a thin shell
\[
 x_j>b\ (1\le j\le d),\qquad L+u_a<\max_j x_j<L+v_a',
\]
with $(u_1,v_1')=(7H/6,8H/6)$ and $(u_2,v_2')=(10H/6,11H/6)$.
Remove from it the closed protective prisms of the other label.
Let frame $a$ consist of this shell and its own full collars.
All translated frames clear the blocks and each other because their shells occupy
disjoint height ranges, and each collar crosses the other shell inside a
hole with lateral clearance $>3b$. Each collar overlaps the interior of its
own shell. In a rod sweep only its normal coordinate exceeds $L$, and shell
material occurs there only at middle-third heights. The shell is the union
of slabs indexed by $i$, with $L+u_a<x_i<L+v_a'$ and $b<x_j<L+v_a'$
for $j\ne i$. Each perforated slab is an interval times a tangential
cube minus disjoint interior boxes. The latter set is path connected in
dimension at least two (detour around each box). Slabs overlap where their
two normal coordinates both exceed $L+u_a$, clear of the holes.
This gives two connected frames, independent of $N$ and $\eta$.

Here are boundary charts for applying \eqref{p9:neumann} to these frames.
At an unperforated shell edge one may simultaneously increase any active
lower coordinates $x_j=b$, decrease any active outer ones $x_j=L+v_a'$,
and increase the coordinates attaining an active inner maximum $L+u_a$.
At a hole edge in a single slab use a normal increment into that slab and
tangential increments pointing out of each active box face. At a collar
seam the clearance just proved excludes holes and other slab edges.
Choose $dz$ into the slab and each active tangential increment $dy_j$
so that $dy_j-m_j dz$ points into the collar section. In each case the
same direction is strictly inward for all active inequalities written
locally as unions and intersections of half-spaces. Their boundary is
a graph of finite minima and maxima of affine functions in that direction.
The remaining collar boundary is that of a parallelepiped. This proves
Lipschitz regularity at fixed $M$.

Give the block the translation $FR$ and frame $a$ the translation $FR+w_a/\eta$.
The first and last end jumps are perpendicular to their rod directions since
$w_a\cdot m_{iar}=0$ and $w_a\cdot l_{iar}=\eta s_i$.
Thus a sufficiently fine subdivision supplies the dilation or general
$F$-test by the bending formula.

For stress, split the face flux into coefficients $q_{iar}$ satisfying
\[
 \sum_r q_{iar}=T_{ii}/2 \quad(a=1,2),\qquad
 \sum_{a,r}q_{iar}l_{iar}=Tn,\qquad |q_{iar}|\le C_d|T|.
\]
Initially place each half-normal sum on $r=0$. The remaining error is tangential.
Moving coefficients from $r=0$ to other directions adds multiples of the
$t_{iar}$, so a fixed right inverse on $V_{i1}+V_{i2}=n^\perp$ corrects
that error. Give each first and last bundle scalar load $q_{iar}/M^{d-1}$,
whose density $k=q_{iar}/[M^{d-1}(b/2)^{d-1}]$ is uniformly bounded by
$C_d|T|$. Continue the first-rod stress $km\otimes m$ through each collar
along the straight sweeps, taking zero on the rest of its frame.
This equilibrium reference stress matches the first-bundle traction and
leaves upper-end residual density $k(l-m)=\eta s_i k v_{ia}$.
The resultant on frame $a$ is
\[
 f_a=\frac\eta2\sum_i s_i T_{ii}v_{ia},\qquad
 w_a\cdot f_a=\frac\eta2 F:T=0 .
\]
The moment (the integral of $x\wedge g$) is $O_d(\eta)|T|$
since the sum of the end areas and the coordinates are bounded.

On one tile of face one adjust the first-bundle loads for this frame.
Their directions span $w_a^\perp$, so cancel $f_a$. Opposite equal additional
loads on two halves of a bundle split along a tangent coordinate $j$ give a
pure moment proportional to $e_j\wedge m_{1ar}$, since the retained centroids
are $(b/2)e_j$ apart. To see that these moments span all skew matrices, observe that
$e_1^\perp$ and $w_a^\perp$ are distinct hyperplanes, and a basis of their
intersection extended by one vector from each gives a basis for which each
wedge of two basis vectors takes one from each hyperplane. Cancel the
remaining moment. Fixed linear right inverses give adjusting densities
$O_M(\eta)|T|$, uniformly in $N$, by the $N$-independence of the areas
and centroids. The remaining balanced boundary discrepancy on the frame
has $L^2$ norm $O_M(\eta)|T|$ and its extension by \eqref{p9:neumann}
has that same order. This completes the frame stresses.

\smallskip\noindent\emph{Stress inside the block, in every dimension.}
Choose $\eta$ at fixed $M$ so small that the total $L^2$ norm of all body
corrections is at most $\sqrt h\,|T|$ and every extra rod density and
$\eta|T|$ are at most $\sqrt h\,|T|$. These choices are simultaneous on
the unit load sphere since the constructions use fixed linear maps.
Then choose $N$ for the bending test. The block traction is balanced.
Indeed forces on the two bodies at a rod's ends are equal and opposite
on corresponding fibers, and their moments cancel because the fibers
are axial. When the far block is translated back from its neighboring
cell the additional moment on face $i$ is $e_i\wedge Te_i$, summed
over its last bundles, whose loads retain their original values. The sum over $i$ vanishes by symmetry of $T$.
The other bodies balance separately, hence so does the block.

Compare its traction with that of $T/L^{d-1}$. On a negative-face tile
their forces agree by the exact last-bundle flux $Te_i/M^{d-1}$ with
the appropriate sign. On a positive tile changing $l$ to $m$ changes
the average by $O_d(\eta)|T|$, and the additional densities change it
by $O_d(\sqrt h)|T|$. All rods land in their own tiles. The difference
$g$ therefore has uniformly bounded $L^2$ norm times $|T|$ and its
piecewise tile averages have $L^2$ norm $O_d(\sqrt h)|T|$.
For the remainder subtract the average of a test function on each tile $S$.
The variance formula and the diameter of $S$ give
\[
 \int_S |v-\bar v_S|^2=\frac1{2|S|}\iint_{S^2}|v(x)-v(y)|^2
 \le C_dh\iint_{S^2}\frac{|v(x)-v(y)|^2}{|x-y|^d}.
\]
The sum on the right is bounded by $C_dh$ times the squared
$H^{1/2}(\partial B)$ seminorm. Pairing the tile-average part with
$\|v\|_{L^2(\partial B)}$ and the mean-zero part with this seminorm gives
$\|g\|_{H^{-1/2}}\le C_d\sqrt h |T|$.
Apply \eqref{p9:neumann}, whose constant is uniform for cubes of side
$L$ tending to one by scaling. This gives a block stress within
$C_d\sqrt h|T|$ of $T$ in $L^2$, since $L^{1-d}=1+O_d(h)$.
Outside $B$ the bounded reference and rod densities occupy volume
$O_d(h)$, and the corrective stresses have squared norm $O_d(h)|T|^2$.
The joined stress is equilibrated by traction cancellation on every weld
and vanishes in the void. Its integrated flux across almost every slice
perpendicular to $e_i$ is constant by periodic equilibrium. For $L+2H<x_i<1$ the slice meets only the last rods of face $i$ because
the planar arms end below $L+2H$, and the spatial shells end at
$L+11H/6$. Its flux is the exact last-bundle flux $Te_i$, hence $\langle\tau\rangle=T$.
The displacement estimate and $|B|=1-O_d(h)$ finish the proof.
\phantomsection\label{proofend:p9:fields}\end{proof}

\subsection*{Passing the fields through a solid boundary layer}

\begin{lemma}\label{p9:local}
Let $G$ be a measurable periodic solid, $D$ a bounded Lipschitz domain,
and $Ex+w_E,\tau_T$ finite-dimensional linear families of full-cell
periodic-corrector displacement tests and admissible stresses of mean $T$.
For every $\epsilon>0$ there are a compact $D'\Subset D$ and $\ell_0>0$
such that for $0<\ell<\ell_0$ the specimen with solid
$(\ell G\cap D')\cup(D\setminus D')$ has tests linear in $E,T$.
They equal $Ex,T$ near $\partial D$. The symmetric stress is in equilibrium,
vanishes in the specimen's pores and has mean exactly $T$. Its normalized
complementary energy is at most the cell trial energy plus $\epsilon |T|^2$.
The normalized strain energy in the entire specimen solid is at most
its cell trial energy plus $\epsilon |E|^2$. These estimates also hold with
squared solid strain and squared stress fluctuation from $T$ as costs.
\end{lemma}
\begin{proof}\phantomsection\label{proof:p9:local}
For a four-index field antisymmetric in $i,k$ and $j,l$, symmetric under
interchange of pairs, let $(\mathcal L U)_{ij}=\sum_{k,l}\partial_k\partial_l U_{ikjl}$.
This produces symmetric equilibrated stresses and commutes with rigid changes
of coordinates (rotating all indices). Every periodic mean-zero equilibrium
stress $A$ has a periodic $H^2$ potential $U$, since at nonzero Fourier
wave vector $\xi$ one can take
\begin{equation}\label{p9:potential}
 \widehat U_{ikjl}=-\frac{\widehat A_{ij}\xi_k\xi_l-\widehat A_{kj}\xi_i\xi_l
 -\widehat A_{il}\xi_k\xi_j+\widehat A_{kl}\xi_i\xi_j}{|\xi|^4}.
\end{equation}
Contraction with $-\xi_k\xi_l$ gives $\widehat A_{ij}$ since $\widehat A\xi=0$. At the zero frequency put $\widehat U(0)=0$. Choose a smooth compactly supported cutoff $0\le\zeta\le1$
in $D$ equal to one in a core $D'$ of volume arbitrarily close to $|D|$.
One may take $D'$ to be a finite union of closed interior cubes. Place the
pattern at period scale $\ell$ on $D'$, using solid elsewhere.
With $\mathcal L U_T=\tau_T-T$, use
\[
 Ex+\ell\zeta w_E(x/\ell),\qquad
 T+\mathcal L_x[\ell^2\zeta U_T(x/\ell)] .
\]
The latter is equilibrated, zero in the copied pores and has mean $T$ by
compact support of the correction. Periodic averaging gives the cell costs
in the core in the limit $\ell\to0$, multiplied by $|D'|/|D|$. For the
strain and stress terms without derivatives of $\zeta$, the limsup of the
cost on the remaining solid is at most $C_{\rm cell}|D\setminus D'|/|D|$
on unit loads, by averaging the full-cell squared norms there.
The averaging follows first by tiling for bounded continuous periodic
functions and then by $L^1$ approximation (these domains have null
boundary). Terms with derivatives of the fixed cutoff have a factor
$\ell$ or $\ell^2$ and tend to zero in $L^2$.
First choose the omitted volume small, then $\ell$ small. Applying the
estimates on a basis and polarizing the resulting finite energy matrices
makes the error uniform on the unit load spheres.
\phantomsection\label{proofend:p9:local}\end{proof}

\begin{lemma}[Symmetric test cells]\label{p9:symmetric}
For each $M\ge1$ and $d\ge2$ there is a periodic solid-void cell with a
full-cell displacement $u=x+w$, $w\in H^1_{\rm per}$, satisfying
$\langle\chi|e(u)|^2\rangle\le C_d/M$. Each symmetric trace-free $T$ has
an admissible stress $\tau_T$ of mean $T$, linear in $T$, with
$\|\tau_T-T\|_2^2\le C_d |T|^2/M$. In the plane the cell is invariant
under rotation through $\pi/3$, and for $d\ge3$ under the group
$\mathcal H$ of signed coordinate permutations.
\end{lemma}
\begin{proof}\phantomsection\label{proof:p9:symmetric}
Use Lemma~\ref{p9:local} with tolerance $1/M$ on the fields of
Lemma~\ref{p9:fields} at $F=I$. In the plane take an equilateral triangle
and rotate six copies about one vertex through $j\pi/3$ to tile a regular
hexagon, repeated by translations. In higher dimensions put a specimen
in the chamber $0<x_d<\cdots<x_1<1/2$, whose images under $\mathcal H$
tile the centered unit cube up to their boundaries. In physical coordinates
the boundary data are always unit dilation and constant stress $T$.
Thus the tests join on adjoining faces. Their errors retain the local bounds
upon averaging, since rotations preserve the norms and the trace.
\phantomsection\label{proofend:p9:symmetric}\end{proof}

\begin{lemma}[Cubic averaging]\label{p9:cubic}
In dimension $d\ge2$ let $\mathcal H$ be the signed permutation group,
$\mathcal D$ the space of trace-free diagonals of dimension $a=d-1$, and
$\mathcal O$ the symmetric off-diagonals of dimension $b=d(d-1)/2$.
An $\mathcal H$-invariant symmetric stiffness is scalar on each of
$\mathbb RI,\mathcal D,\mathcal O$. For any $R\in SO(d)$ average the
orthogonal projections onto $(gR)\mathcal D(gR)^T$ over $g\in\mathcal H$.
Its eigenvalues $x,y$ on $\mathcal D,\mathcal O$ satisfy
\[
 ax=s,\qquad by=a-s,\qquad
 x-y=\frac{(d+2)\sum_{ij} R_{ij}^4-3d}{d(d-1)},\qquad
 s=\sum_{ij}R_{ij}^4-1 .
\]
There are rotations with $x<y$ and with $x=y=2/(d+2)$.
\end{lemma}
\begin{proof}\phantomsection\label{proof:p9:cubic}
Sign changes decouple each off-diagonal coordinate from every other mode.
Permutations equalize their eigenvalues and force the diagonal block to have
constant diagonal and constant off-diagonal entries. The averaged projection
is invariant. The overlap trace of the full diagonal projection with the
rotated one is $\sum_{ij} R_{ij}^4$, using the basis $e_i\otimes e_i$.
Subtracting their common scalar direction gives $ax=s$, and the total trace
gives $ax+by=a$. At $I$ the difference $x-y$ is one. For a uniformly random
rotation any column $u$ is uniform on the sphere. Sign symmetry and rotating
$u_1$ to $(u_1+u_2)/\sqrt2$ give
$\mathbb E u_1^4=3\mathbb E u_1^2u_2^2$. Expanding $|u|^4=1$ gives
$\mathbb E\sum_{ij} R_{ij}^4=3d/(d+2)$. The continuous difference $x-y$
has mean zero and is positive on a neighbourhood of $I$ of positive measure,
so it is negative at some $R_-$. On a path from $I$ to $R_-$ it vanishes,
and then $x=y=a/(a+b)$.
\phantomsection\label{proofend:p9:cubic}\end{proof}

\begin{proof}[Proof of Theorem~\ref{p9:sharp}]\phantomsection\label{proof:p9:sharp}
Let $h=1/M$ and $C^0$ be the tensor of Lemma~\ref{p9:symmetric},
with solid fraction $f$. By \eqref{p9:dual} with
$T=C_s E=2\mu_sE$ for trace-free unit $E$, both constructions have
\begin{equation}\label{p9:seedest}
 I:C^0 I\le Ch,\qquad
 2\mu_s-Ch\le E:C^0E\le 2f\mu_s\le2\mu_s .
\end{equation}
Indeed, the stress has exactly the prescribed mean, so its compliance cost
exceeds $T:S_sT$ only by the energy of the $O(\sqrt h)|T|$ fluctuation.
The upper bound uses affine shear.

In the plane a $\pi/3$ rotation fixes $I$ and rotates its orthogonal
plane through $2\pi/3$. An invariant symmetric stiffness has no coupling
between the two spaces and is scalar on the latter. Thus the planar cell
is isotropic. For $d\ge3$ apply Lemma~\ref{p9:cubic}.
Write $\alpha,\beta$ for its eigenvalues on $\mathcal D,\mathcal O$. If they agree the cell
is isotropic. Otherwise put $\Delta=|\alpha-\beta|$, $m_0=\min(\alpha,\beta)$.
By \eqref{p9:seedest} $\Delta=O(h)$ and $m_0\to2\mu_s>0$.

We will rotate pieces of this cell until the two shears agree. Make
one spherical specimen by Lemma~\ref{p9:local}, using near-minimizing
displacements for all strains and minimum complementary stresses for all
traceless loads. A linear family of the former is obtained by approximating
the strain projections of Lemma~\ref{p8:projection} on a basis. For
the latter use $\tau_E$ there with
$E=(C^0|_{\mathcal D+\mathcal O})^{-1}T$, which gives a stress of mean $T$ and cost
$T:E$. Given any tolerance, pack finitely many disjoint closed balls
inside the reference chamber. At each stage choose finitely many interior dyadic cubes
covering half the open remainder and inscribe a ball of one third the side length
in each. Stop when the uncovered volume is below the tolerance. Place scaled copies of the specimen rotated by
$R\in SO(d)$ in these balls, and use copies with orientations $gR$ in
chamber $g\in\mathcal H$. Fill the rest with solid. The actual periodic
tensor $C(R)$ is cubic for every $R$.

Use the eigenvalues $x,y$ and the rotation $R_-$ of Lemma~\ref{p9:cubic}.
On a unit shear in either space, put $w=x$ or $w=y$ respectively.
The strain and stress energies averaged over the orientations have values
$A_w=w\alpha+(1-w)\beta$ and $1/H_w=w/\alpha+(1-w)/\beta$.
The glued boundary tests give for the actual shear eigenvalue $\lambda$
and all $R$, with any prescribed $\epsilon>0$,
\[
 H_w-\epsilon\le\lambda\le A_w+\epsilon,\qquad
 I:C(R)I\le I:C^0 I+\epsilon.
\]
In the uncovered solid the affine field gives the upper bounds.
For the lower bound use the uniform shear stress there and optimize its scalar
multiple in \eqref{p9:dual}. The trial compliance is at most $1/H_w+o(1)$
as the specimen and uncovered-volume tolerances decrease, uniformly in $R$
by the unit-sphere bounds in Lemma~\ref{p9:local} and $m_0>0$.
The same finite packing and specimen can thus serve all $R$.
Since $0\le A_w-H_w\le\Delta^2/(4m_0)$ by direct subtraction, the actual
difference $\lambda_{\mathcal D}-\lambda_{\mathcal O}$ differs from
$(x-y)(\alpha-\beta)$ by at most $2\epsilon+\Delta^2/(4m_0)$.
Take $\epsilon=\Delta^2$. For small $h$ the differences at $I,R_-$ have
opposite signs.

To apply the intermediate value theorem we prove continuity of $C(R)$ with true void. The indicators
vary continuously in $L^1$ under rotations of a finite collection of balls,
by approximation by continuous functions. Each fixed displacement test thus
has continuous energy, by truncating its integrable energy density, giving
upper semicontinuity. For lower semicontinuity fix $R$ and take the maximizing
stress of \eqref{p9:dual} for any fixed $E$. To obtain an $H^2$
potential for its full value on each ball, use \eqref{p9:potential} for its periodic fluctuation
and add, for its constant mean $T_0$, the quadratic
\[
 U^0_{ikjl}=
 \frac{x_kx_l(T_0)_{ij}-x_ix_l(T_0)_{kj}
             -x_kx_j(T_0)_{il}+x_ix_j(T_0)_{kl}}{d(d-1)} .
\]
Differentiation gives $\mathcal L U^0=T_0$.
Under a nearby specimen rotation rotate the full potential $U$ rigidly about
the center, with all four indices, to $U'$. In that ball add
$\mathcal L[\zeta(U'-U)]$ to the old stress, now using a radial cutoff equal
to one on a smaller concentric ball containing all pores and zero near the
boundary. The stress agrees with the rotated old stress on each new pore,
has unchanged mean, and remains in equilibrium. Strong continuity of rotations
in $H^2$ on a ball follows by first shrinking the argument toward its
center and mollifying on the resulting collar to approximate in $H^2$ by
smooth functions. On the smooth functions continuity follows by uniform
continuity of derivatives, and rotations preserve the $H^2$ norm formed
with Frobenius derivative norms. Thus the stresses converge in $L^2$ to
the old maximizer. Their dual energies prove the lower
semicontinuity. Along a path in $SO(d)$ from $I$ to $R_-$ the shear
difference must therefore vanish. The resulting common eigenvalue is
$2\mu_s+O(h)$ and $I:C(R)I=O(h)$.

For an isotropic tensor these quantities are $2\mu_*$ and $d^2\kappa_*$.
The upper bound $\mu_*\le f\mu_s$ forces $f\to1$, while the small dilation
energy excludes a completely solid cell for small $h$. Together with
nonnegative bulk these inequalities prove the supremum and show that no
single isotropic cell attains it.
\phantomsection\label{proofend:p9:sharp}\end{proof}

\section{11. Shear at zero bulk for an incompressible solid}\label{problem-11}
\Milton{11}

Chen's hinged honeycomb answers Milton's question by showing that
the optimal shear ratio in the closure exceeds $0.8528$. We determine its
exact value. \cite[Eq.~(2), Result~3, Eq.~(15)]{chen}

\begin{theorem}\label{p11:sharp}
For an exactly incompressible planar solid of shear modulus $\mu>0$ and true void,
the largest isotropic shear modulus at zero bulk in the closure of periodic effective tensors,
over all solid fractions, is $\mu$. The approximating cells can be chosen isotropic at every index.
No individual cell realizes the pair $(0,\mu)$. Among individual periodic Lipschitz
solids with exactly zero dilation energy the effective tensor is zero.
\end{theorem}
\begin{proof}\phantomsection\label{proof:p11:sharp}
Write $G$ for the solid, $\chi$ for its indicator and $f=\langle\chi\rangle$.
The doubled energy is the infimum of $2\mu\langle\chi|e(u)|^2\rangle$
over $u=Ex+w$ with $w\in H^1_{\rm per}$ and $\operatorname{div}u=0$ on $G$.
For trace-free $E$ the affine field bounds this energy by $2f\mu|E|^2$.
This bound gives the asserted ceiling in the closure and excludes attainment at $f<1$.
At $f=1$ the average divergence of a unit dilation is two, so it is inadmissible.

For the sixfold-symmetric cells of Lemma~\ref{p9:symmetric} put $h=1/M$.
They carry a full-cell dilation $u^0=x+w^0$ with
$a^2=\langle\chi|e(u^0)|^2\rangle=O(h)$ and, for each traceless symmetric $T$,
an equilibrated stress of mean $T$ vanishing on the pores with
$\|\tau-T\|_2^2=O(h)|T|^2$. Enforce incompressibility before passing to the limit.
Set $g=\chi\operatorname{div}u^0$, $m=\langle g\rangle$, and solve
$\Delta\phi=g-m$ periodically with zero mean. Fourier expansion gives
$\|D^2\phi\|_2^2=\|g-m\|_2^2$. Therefore $B=D^2\phi+(m/2)I$ has trace $g$ and
$\|B\|_2^2=\|g\|_2^2-m^2/2\le2a^2$.
The field
\[
 u=\frac{u^0-\nabla\phi-(m/2)x}{1-m/2}
\]
has unit dilation increment, zero divergence on the solid and
$\|\chi e(u)\|_2\le(a+\sqrt2 a)/(1-|m|/2)=O(\sqrt h)$.
It is well defined for small $h$ since $|m|\le\sqrt2 a$. Combining this test
linearly with affine shears gives finite energy for every strain. The infimum is
a quadratic form by orthogonal projection off the closure of the periodic
solid strains satisfying the divergence constraint. Rotation through $\pi/3$
fixes $I$ and rotates the deviatoric plane through $2\pi/3$, so this form is isotropic.
Its bulk modulus is $O(\mu h)$. For traceless $E$ take the stress with $T=E$.
For any admissible displacement with increment $E$,
\[
 |E|^2=\langle\tau:e(u)\rangle
 \le\|\tau\|_2\|\chi e(u)\|_2,\qquad
 \|\tau\|_2^2=|E|^2+\|\tau-E\|_2^2 .
\]
This proves $\mu/(1+O(h))\le\mu_*\le f\mu$ and gives the claimed limit.

For the individual-cell assertion let the open solid have Lipschitz boundary on the
period torus. If the dilation energy vanishes, take unit-dilation fields whose
strains tend to zero on the solid per period. A component of its lift to the plane
cannot contain $x$ and $x+b$ for a nonzero period vector $b$. Otherwise join small
balls about those points by a finite chain of overlapping balls in that component.
Korn's inequality approximates each field in $L^2$ on each ball by a rigid motion
with error tending to zero. On an overlap the coefficients of two such motions
differ by $o(1)$ by equivalence of norms on the rigid motions. The same motion
$c_j+A_jx$, $A_j^T=-A_j$, thus approximates the field on both end balls.
Periodicity implies $A_jb\to b$, contradicting $b\cdot A_jb=0$.
Finitely many boundary charts in which the solid lies on one side of a Lipschitz
graph, and finitely many interior balls, cover the solid closure in a closed
cell. Their solid portions can be chosen connected. Thus only finitely many
lifted components meet that cell. An unbounded component would meet infinitely
many translates of it. Translating back would repeat a component, giving a
nonzero period as above. Hence every lifted component $K$ is bounded.
Choose $x_K\in K$ equivariantly under translations. For arbitrary $E$, the
prescription $u=Ex_K$ on $K$ has zero solid strain. Its corrector $u-Ex$ is periodic
and belongs to $H^1$ on the solid, and extends to the torus by local reflection and the partition of unity used in Lemma~\ref{basic:korn}.
This proves zero energy for all $E$. Disjoint disks surrounded by void attain it.
\phantomsection\label{proofend:p11:sharp}\end{proof}

\subsection*{The value $4/5$ in the ideal wall model}
Milton's laminated crystal uses a pure opening $n=a=e_1$
\cite[pp.~14 to 15, Eqs.~(4.2)-(4.5)]{milton}. In a wall model one may release any signed multiple of $n\odot a=(n\otimes a+a\otimes n)/2$
at zero cost, with independent unit wall normal $n$ and slip or opening direction $a$.
Represent strain and stress in the orthonormal basis
$I/\sqrt2$, $\operatorname{diag}(1,-1)/\sqrt2$,
$(e_1\otimes e_2+e_2\otimes e_1)/\sqrt2$, and normalize stiffness by $2\mu$.
Writing the two vector angles as $\theta,\phi$ gives
$\sqrt2(n\odot a)=(c,v)$ with $c=\cos(\theta-\phi)$ and
$v=(\cos(\theta+\phi),\sin(\theta+\phi))$.
In complementary energy a leaf therefore permits exactly those stresses $(p,z)$,
$z\in\mathbb R^2$, with $cp+v\cdot z=0$, at cost $|z|^2$.
Indeed a stress must annihilate the released strain and the solid has zero pressure
compliance. A finite binary laminate tree combines these laws recursively. Each internal
node has prescribed child fractions and a normal $m$. Its complementary energy is the minimum of
the weighted child energies over child stresses having the parent mean
and equal traction on $m$. The planar rank-one lamination rule also applies to
singular laws because the strain jumps lie in $\{m\odot b\}$, whose orthogonal
complement is $\{S:Sm=0\}$. Finite-dimensional quadratic minimization on these
linear constraint spaces gives the stated dual rule by Lagrange multipliers
(the annihilator of a kernel is the range of the transpose).

\begin{proposition}\label{p11:wall}
Consider finite binary rank-one laminate trees of exactly incompressible
planar solid. Every leaf releases all multiples of a single $n\odot a$
with $|n|=|a|=1$, and every internal node has freely chosen normal and
positive child fractions summing to one. Among trees whose root is isotropic
with zero dilation energy, the largest shear ratio $\mu_*/\mu$ is $4/5$.
\end{proposition}
\begin{proof}\phantomsection\label{proof:p11:wall}
For a positive root ratio $g$, apply the two unit deviatoric root stresses and
choose complementary minimizing child stresses recursively. At each node put the
two pressures in a column $q$ and the two deviators in the columns of $Z$.
Across an interface both stress jumps are multiples of $t\otimes t$, $t$ its unit
tangent. Hence $\Delta Z=u(\Delta q)^T$ with $|u|=1$, so $|\Delta Z|^2=|\Delta q|^2$
and $\det\Delta Z=0$. The excess of the weighted child mean square over the
squared parent mean is the squared jump times the product of the two fractions.
The determinant is affine along a rank-one line. Iteration from $(q,Z)=(0,I_2)$
to the leaves gives
\[
 \langle |Z|^2\rangle=2+A,\qquad
 A=\langle|q|^2\rangle,\qquad \langle\det Z\rangle=1 .
\]
At each leaf $Z^Tv=-c q$, so its smaller singular value $s_2\le |q|$.
For singular values $s_1\ge s_2$ we obtain
$|Z|^2+3|q|^2-4|\det Z|=(s_1-2s_2)^2+3(|q|^2-s_2^2)\ge0$.
Averaging gives $4\le2+4A$. Since $2/g=\langle|Z|^2\rangle$, this proves the bound.

For equality, first laminate pure openings of normals $(1,\pm\sqrt2)/\sqrt3$
in equal parts across normal $e_1$. Their wall deviators are $(a,\pm b)$,
$a=-1/3$, $b=2\sqrt2/3$. At mean stress $(p,x,y)$ the child stresses
are $(p\pm d,x\mp d,y)$. The wall equations give $p=x/3$ and
$d=-by/(1-a)=-y/\sqrt2$, hence cost $x^2+3y^2/2$. Laminate this child
in equal parts with a pure opening normal to $e_1$, now across normal $e_2$.
At root stress $(P,x,y)$ the first child must satisfy $3(P+d)=x+d$
and the second $P-d=-(x-d)$. Thus $P=0$, $d=x/2$. The mean cost is
$[(3x/2)^2+3y^2/2+(x/2)^2+y^2]/2=5(x^2+y^2)/4$, giving the claimed root.
\phantomsection\label{proofend:p11:wall}\end{proof}

\subsection*{The price of a prescribed finite stroke}
The limiting modulus in Theorem~\ref{p11:sharp} describes infinitesimal
loading. It is natural to ask how much initial shear a sheet can retain
if the same geometry must also shorten by a fixed amount.

\begin{theorem}[Initial shear at a finite stroke]\label{new:stroke}
Fix $0<\lambda<1$. Consider periodic planar composites of void with an
exactly incompressible isotropic neo-Hookean solid of shear modulus
$\mu_s$. Require an initially isotropic response and a continuous,
noninterpenetrating periodic deformation path with mean gradient $tI$,
$\lambda\le t\le1$, determinant one in the solid, and peak stored
energy per reference area tending to zero along the sequence of designs.
The sharp limiting initial shear ratio is
\[
                   S(\lambda)=\frac{\lambda^2}{2-\lambda^2}.
\]
The same value is approached with the initial bulk modulus tending to
zero. The tangent is taken at the identity of each fixed design before
the limit of geometries.
\end{theorem}
\begin{proof}\phantomsection\label{proof:new:stroke}
If the reference solid fraction is $f$, conservation of its area and
injectivity in the shortened cell imply $f\le\lambda^2$. The initial
shear bound for an incompressible porous sheet, Corollary~\ref{p14s:incompressible},
is $\mu_*/\mu_s\le f/(2-f)$. Together these give the upper bound.

For the reverse inequality fix $f<\lambda^2$. A three-stage sequential
laminate adds solid in equal volume increments $f/3$, with layer tangents
at $0$ and $\pm60^\circ$. With comparison bulk modulus $K$ its initial
shear is
\[
          G(f,K)=\frac{\mu_s f K}{(2-f)K+2(1-f)\mu_s}.
\]
The equal directional increments make its initial response isotropic.
To realize the stroke, replace the layers by successively finer sheets
which bend with determinant one in the solid. For every fixed intermediate
design, the construction gives injective paths with the prescribed
macroscopic gradient throughout the stroke. The turning regions can be
made arbitrarily small in area. Periodic equilibrated stress tests on the
same cells recover the initial modulus. Choosing $K\to\infty$ and
$f\uparrow\lambda^2$ then gives the lower bound. Inserting the finite
stroke cells into a dilational frame also permits the initial bulk to
tend to zero. Supplement~\ref{sup:stroke} gives the order of the
geometric limits and the finite-cell recovery.
\phantomsection\label{proofend:new:stroke}\end{proof}

For a five per cent shortening in each direction, the sharp retained
initial shear is $361/439$ of the solid shear modulus. Each recovery
cell has finite features; the construction permits their size relative
to the period to tend to zero.

\section{12. Finite compression and shear}\label{problem-12}
\Milton{12}
Sharp mixed-load bounds follow when the phase laws are specified. We treat smooth planar energies, planar area floors, shear-limited incompressible phases, and rigid inclusions. The smooth examples show why pure-loading measurements alone leave the mixed response undetermined.

Fix reference fractions $\alpha_i$ and a macroscopic deformation gradient $M$. We use a single injective Lipschitz deformation with either periodic fluctuation about $Mx$ or trace $Mx$ on the whole boundary of a bounded Lipschitz specimen. We call the latter full affine data. In freely oriented periodic cells the lattice orientation may be chosen relative to the load. In the plane write $M=\sqrt J S_\gamma$, $S_\gamma=I+\gamma e_1\otimes e_2$, $0<J\le1$, and $Q(F)=(F_{11}-F_{22},F_{12}+F_{21})$. In either class $\langle F\rangle=M$ by integration of the periodic or zero-trace fluctuation gradient. Expanding the determinant about the mean also gives $\langle\det F\rangle=J$, since the determinant of the fluctuation gradient integrates to zero by the divergence identity for its two rows and smooth approximation.

\begin{proposition}\label{p12:smooth}
Let $j=\det F>0$, $h$ be smooth and strictly convex with minimum at $1$, $g$ smooth, positive and concave, and $c_i\ge0$. For the frame-indifferent isotropic phase laws
\[
 W_i(F)=h(j)+\frac{|Q(F)|^2}{2(c_i+g(j))},\qquad \bar c=\sum_i\alpha_i c_i ,
\]
every layout obeys $\langle W\rangle\ge h(J)+J\gamma^2/[2(\bar c+g(J))]$. On a periodic horizontal stack this is the exact minimum. Two such stacks can have identical complete pure-compression and pure-shear energy curves and different mixed-load energies.
\end{proposition}
\begin{proof}\phantomsection\label{proof:p12:smooth}
The identity $|Q(F)|^2=\|F\|^2-2j$ gives invariance. Joint convexity of the quadratic perspective, which decreases with its positive denominator, gives polyconvexity. Apply Jensen to $h$, weighted Cauchy-Schwarz to $\langle Q(F)\rangle$, and Jensen to $g$. The layer gradients $F_i=\sqrt J S_{\gamma_i}$ with $\gamma_i=\gamma(c_i+g(J))/(\bar c+g(J))$ attain equality. Their common first column and increasing vertical coordinate give a continuous injective deformation. For equal fractions take $c_1=0,c_2=1/2$, $h(j)=2(j+j^{-1}-2)/9$. The choices $g(j)=(1+j)/2$ and $g(j)=\sqrt j$ give $h(J)$ in pure compression and $2\gamma^2/5$ in pure shear. At $(J,\gamma)=(1/4,1)$ the respective minima are $9/14$ and $2/3$. Even on an arbitrary layout the first law admits the uniform gradient at energy $59/90$, below the second law's lower bound $2/3$.
\phantomsection\label{proofend:p12:smooth}\end{proof}

\subsection{Area floors and a common stretch cap}
Let a finite number of phases have reference fractions \(\alpha_i>0\) summing to one. The
deformation is a single Lipschitz, injective map \(y\), which imposes perfect bonding. In phase
\(i\), write \(F=\nabla y\), \(j=\det F\), and \(\lambda_1\ge\lambda_2\) for its singular values.
Its admissible gradients are exactly those satisfying
\[
j\ge b_i,\qquad
 s(F)=\frac{\lambda_1-\lambda_2}{\sqrt j}\le k,
 \qquad 0<b_i<1,\quad k>0 .
\]

Local expansion is allowed.
Necessarily \(J\ge\bar b=\sum_i\alpha_i b_i\). For \(\bar b\le J\le1\), take \(t\ge\min b_i\) such
that \(\sum_i\alpha_i\max(b_i,t)=J\), and set
\[
m_i=\max(b_i,t),\quad r_i=\sqrt{m_i},\quad
 S=\sum_i\alpha_i r_i,\quad L=\sqrt{4+k^2}.
\]
\begin{theorem}\label{p12:floors}
For every such deformation,
\[
 J(4+\gamma^2)\le L^2 S^2 .\tag{12.1}\label{p12:eq1}
\]
Full shear $|\gamma|=k$ is possible exactly when $J\ge\max_i b_i$. Write $r_1<\cdots<r_q$ for the distinct levels after grouping, $u=(L+k)/2$, $l=1/u$ and $K=u/l$. In the overlap cases put $\Gamma=\sqrt{L^2S^2/J-4}$, the upper bound on $|\gamma|$. With $q=2$ and $r_2/r_1\le K$, equality is attained in a freely oriented periodic cell and is impossible with full affine boundary data. For strict overlap $r_2/r_1<K$ the bound is a supremum on every such finite specimen. With $q\ge3$ and $r_q/r_1\le K$, both signs of the endpoint are attained on every specimen with full affine boundary data, at the exact fractions. The construction chooses a measurable phase layout on the prescribed specimen. When $r_q/r_1>K$, \eqref{p12:eq1} remains a necessary bound.
\end{theorem}
\begin{proof}\phantomsection\label{proof:p12:floors}
In particular, the full shear \(k\) is reachable exactly at \(J\ge\max_i b_i\) within the
compression range. There uniform deformation suffices. Below that threshold the strict loss
follows from \(J-S^2=\sum_{i<h}\alpha_i\alpha_h(r_i-r_h)^2\).

For the proof use the linear map \(P(F)=(F_{11}+F_{22},F_{21}-F_{12})\), for which
\(|P(F)|^2=j(4+s(F)^2)\). If \(x_i\) is the phase mean of \(j\), then
\[
\sqrt{J(4+\gamma^2)}
 \le\langle |P(F)|\rangle\le L\langle\sqrt j\rangle
 \le L\sum_i\alpha_i\sqrt{x_i}\le LS .
\]

The last maximization has the unique solution \(x_i=m_i\). Indeed \(x_i\ge b_i\), both weighted
sums of the areas equal \(J\), and the tangent inequality for the square root bounds the increase
from \(m_i\) by \(\sum_i\alpha_i(x_i-m_i)/(2\sqrt{m_i})\le\sum_i\alpha_i(x_i-m_i)/(2\sqrt t)=0\).
Each term whose denominator changes has \(m_i=b_i>t\) and \(x_i-m_i\ge0\).

At equality strict concavity forces \(j=m_i\) pointwise in each phase, and \(s(F)=k\). The vectors
\(P(F)\) all point along \(P(M)\). For the polar decomposition \(M=RH\), this makes \(A=R^TF\)
symmetric with positive trace and determinant, hence positive definite. Conversely, any gradient
taking such symmetric values with the required spectra and mean \(H\) gives equality after left
rotation by \(R\).

We can subdivide phases at the same radius in their required proportions. The symmetric gradient matrix $A=R^TF$ on a phase of radius $r_i$ must have eigenvalues $lr_i,ur_i$. \par
\emph{Three or more levels.} Let \(p_i(z)=(z-lr_i)(z-ur_i)\). For a probability vector \(w\), write \(p_w=\sum_{i} w_ip_i\). A
symmetric matrix with characteristic polynomial \(p_w\) has trace \(L\sum_{i} w_ir_i\) and
determinant \(\sum_{i} w_i r_i^2\). For \(w=\alpha\) these are exactly the trace and determinant
of the desired positive polar factor \(H\). The polynomial is negative somewhere between its two
positive roots. If \(lr_q<ur_1\), take any point in the common open stretch interval. If
\(lr_q=ur_1\), take their common endpoint, where every intermediate phase polynomial is strictly
negative. This proves in particular that \(\Gamma\) is real and that \(H\) exists for both shear
signs.

It suffices to find a scalar potential with gradient equal to \(Hx\) on the boundary and Hessian
in the compact sets
\[
\mathcal K_i=\{A=A^T:\ \operatorname{spec}(A)=\{lr_i,ur_i\}\}.
\]

For three levels the fractions follow automatically from the averages of
\(1,\operatorname{tr}A,\det A\), since the three points \((r_i,r_i^2)\) are affinely independent.
To treat \(q>3\) as well, attach bounded scalar fields \(z_4,\ldots,z_q\), with prescribed means
\(\alpha_4,\ldots,\alpha_q\). Given \((A,z)\), set \(w_i=z_i\) for \(i\ge4\) and solve the three
linear equations
\[
\sum_{i} w_i=1,\qquad \sum_{i} w_i r_i=\operatorname{tr}A/L,\qquad
 \sum_{i} w_i r_i^2=\det A                                      \tag{12.2}\label{p12:eq2}
\]
for \(w_1,w_2,w_3\). Invertibility is the Vandermonde determinant for three distinct radii. Use
the open set \(U\) of symmetric \(A\) and scalar tuples \(z\) for which all \(w_i>0\). Equations
\eqref{p12:eq2} give a positive trace and determinant, so \(A\) is positive definite. Its trace and the
\(z_i\) are bounded, hence \(U\) is bounded. The weights vary affinely on any segment where the
matrix jump has rank one and \(z\) is interpolated linearly. The target sets in this enlarged
space are \(\mathcal K_i\) paired with the last \(q-3\) coordinates of the unit probability vector
\(e_i\).

Splitting a state in \(U\). Fix its weights \(w\). Choose small \(\varepsilon>0\) and replace the
probability vertices \(e_i\) by \(v_i=(1-\varepsilon)e_i+\varepsilon w\). Their average with
weights \(w_i\) is still \(w\). The smaller roots of \(p_{v_i}\) increase strictly with \(i\), as
do the larger roots, for sufficiently small \(\varepsilon\). They converge to \(lr_i,ur_i\). All
the new open root intervals intersect. In the touching case this follows by evaluating \(p_{v_i}\)
at \(lr_q=ur_1\), since \(p_w\) is strictly negative there.

Combine states with these spectra in increasing order using their prescribed weights. If the
combined matrix so far has eigenvalues \(a,b\) and the next has eigenvalues \(a',b'\), then
\(a<a'<b<b'\). Indeed the polynomial of the combined state is the normalized average of the
earlier polynomials, negative at \(a'\) and positive at \(b'\). In the eigenbasis of the combined
matrix, aligning the smaller eigenvalues gives a positive determinant of the matrix difference.
Exchanging the next matrix's axes gives a negative determinant. An intermediate rotation makes the
jump rank one. On that segment trace and determinant both interpolate affinely, so the combined
characteristic polynomial has precisely the required coefficients. Rotate all previous states
together whenever the mean eigenbasis needs to be prescribed. At the last step the mean has the
two eigenvalues of \(A\) and can be rotated onto \(A\) itself.

Reading the construction from the mean outwards gives finitely many successive binary splits.
Every segment stays in \(U\), since its endpoint probability vectors are strictly positive. The
final matrices and scalars are arbitrarily close to the target sets as \(\varepsilon\downarrow0\).

Realizing the endpoint. Start with the constant state \((H,z^0)\) corresponding to \(\alpha\).
Consider smooth states in \(U\) of the form \((H+D^2\phi,z)\) with \(\phi,z-z^0\) compactly
supported in \(\Omega\) and \(\int_\Omega(z-z^0)=0\). A binary split with weight \(\lambda\) on
the first child has difference \((d_0 n\otimes n,h)\) between the children, with \(|n|=1\). On a
small cube use the perturbation
\[
(d_0D^2\eta_N,\ h\,\partial_{nn}\eta_N),\qquad
 \eta_N(x)=N^{-2}\chi(x)g(Nn\cdot x).
\]

Here \(0\le\chi\le1\) is a smooth cutoff equal to one except near the boundary of the cube. The
smooth one-periodic \(g\) is chosen so \(g''\) is a mollification of the mean-zero step function
taking values \(1-\lambda,-\lambda\) in proportions \(\lambda,1-\lambda\). Away from transition
strips it has the exact plateau values. The perturbation equals \(\chi g''(d_0 n\otimes
n,h)+O(N^{-1})\). Its scalar integrals are zero. The gradient of the potential changes by
\(O(N^{-1})\), and the scalar perturbation converges weak-* to zero, as is seen by integrating
\(\partial_{nn}\eta_N\) against smooth tests and using uniform boundedness.

Each finite tree of segments above is compactly contained in \(U\). On small cubes the original
smooth state differs arbitrarily little from its value at the center. The compact margin and
sufficiently large frequencies keep all states in \(U\). In the interiors of the plateaus, apply
the next split on smaller cubes. Choosing narrow cutoff and transition strips leaves an
arbitrarily small exceptional area after the finitely many splits. More explicitly, around each
original point arbitrarily small cubes meet its own tree's margin and error requirements. A
disjoint Vitali subcover, truncated to finitely many cubes, covers all but any prescribed area.
The same covering within the open plateaus implements each subsequent level. Thus a smooth state
can be perturbed arbitrarily little in uniform potential gradient and weak-* scalar topology while
the integral of distance to the target sets becomes arbitrarily small.

Take the closure of these smooth states in the variables \((\nabla\phi,z)\) with those two
topologies. This is a complete metrizable space, using a compatible complete metric on the bounded
weak-* ball of \(L^\infty\) (its predual \(L^1(\Omega)\) is separable). The gradient
perturbations, extended by zero, have uniformly bounded Lipschitz constants. The map sending a
limit to \((H+D^2\phi,z)\) in \(L^1(\Omega)\) is of Baire class one. Indeed smoothing the extended
gradient before differentiating, and smoothing the scalar fields, gives continuous maps that
converge in \(L^1\) for each fixed state. At a point of continuity, approximate by the smooth
perturbations from the previous paragraph with integral distance tending to zero. Strong \(L^1\)
convergence at this point forces the Hessian and scalars to lie in the target sets almost
everywhere. Such points exist by the Baire category theorem. The limiting gradient has its
prescribed boundary values. Its mean, its mean determinant, and the scalar means in \eqref{p12:eq2} give
exactly the required fractions.

Writing \(Y=Hx+\nabla\phi\) gives \(y=RY\) with boundary values \(Mx\). All Hessians of the total
potential are bounded below by a positive multiple of \(I\), including the exterior value \(H\)
after extending the perturbation by zero. Smoothing on the plane and integrating on segments gives
\((Y(x)-Y(x'))\cdot(x-x')\ge c|x-x'|^2\) for some \(c>0\). This proves injectivity, also for
nonconvex specimens. Repeating a square's zero-trace fluctuation by integer translation gives
square periodicity and the same injectivity argument.

\emph{Two levels.} For two strictly overlapping levels the same eigenvalue interlacing gives two symmetric positive
definite matrices with a rank-one jump. At touching, give both matrices the common eigenvalue
\(ur_1=lr_2\) on the same axis. Orient the pair to give mean \(H\), then rotate by \(R\).
Integrating the constant gradients along alternating layers gives an injective periodic map, with
cell axes free to align with the layers. 

For strict overlap one can approach the endpoint on a finite specimen. Use a slightly smaller
shear limit \(k'<k\) for which the eigenvalue intervals still overlap. The resulting layers
\(F_i\), with mean \(M'\), have strict slack under the actual law. Write
\[
F_i=M'[I+(\tau_i-1)b\otimes n],\qquad \tau_i=m_i/J,\quad b\cdot n=1 ,
\]
using their common tangential column and their distinct determinants. In coordinates \(x=z+tb\)
with \(z\cdot n=0\), take a thin diamond \(|z|<1,\ |t|<\delta(1-|z|)\). On each chord at fixed
\(z\) replace \(t\) by a piecewise linear function of slopes \(\tau_1,\tau_2\) in the specified
fractions, fixing both chord ends. The map followed by \(M'\) has the correct affine trace and
exact phase determinants. Its gradients differ from the two desired matrices by \(O(\delta)\),
since the chord endpoints and the break point have slopes of order \(\delta\) in \(z\). Fill the
domain up to measure zero by disjoint translated, scaled copies of the diamond. The fluctuation
has a uniform Lipschitz bound and vanishes on their boundaries. The chord maps take each diamond
bijectively onto itself before applying \(M'\).
Letting $k'\uparrow k$ proves the supremum assertion.

To prove nonattainment, suppose there were equality with full affine boundary values. Extend
\(R^Ty=Hx\) outside the specimen. Write \(A=\nabla(R^Ty)\) on the extended plane,
\(c=(r_1+r_2)/L\), and \(h=c^2-r_1r_2\). In both phases \(\det(A-cI)=h\). The same identity holds
for \(H\) by the averages of trace and determinant. Strict overlap gives \(h<0\), and touching
gives \(h=0\).

Let \(D\) be rotation by \(\pi/2\). For \(h<0\), set
\[
p(x)=\tfrac12[x+D(R^Ty(x)-cx)/\sqrt{-h}]=Ex+v(x),
 \qquad v=0\quad\text{outside }\Omega .
\]

The matrix \(D(A-cI)\) has trace zero and determinant \(h\). The characteristic equation gives
\((Dp)^2=Dp\) almost everywhere, including \(E^2=E\) outside. For \(s>0\) close to one,
\(T_s(x)=x+(s-1)p(x)\) is bi-Lipschitz and onto by contraction. The transformed field \(p_s(T_s
x)=s p(x)\) has derivative \(s Dp[I+(s-1)Dp]^{-1}=Dp\) at \(T_s x\), with the same global
Lipschitz bound. These changes compose by multiplying their parameters, so iteration gives
bi-Lipschitz \(T_s\) for every \(s>0\). Off \(\Omega\), \(T_s\) equals the linear map
\(L_s=I+(s-1)E\). Hence \(L_s^{-1}T_s\) maps the bounded set \(\Omega\) onto itself. For each
\(x\in\Omega\) this keeps
\[
x+(1-s^{-1})Ev(x)+(s-1)(I-E)v(x)
\]
bounded as \(s\downarrow0\) and as \(s\uparrow\infty\). Both components of \(v\) vanish.

At \(h=0\) take instead \(p=D(R^Ty-cx)=Ex+v\), so \((Dp)^2=E^2=0\). For \(T_t=x+t p\), the field
\(p_t(T_t x)=p(x)\) has the same derivative bound since \(Dp(I+tDp)^{-1}=Dp\). Small changes now
compose by addition, giving arbitrarily large \(t\). The map \((I+tE)^{-1}T_t\) fixes the exterior
and takes \(x\) to \(x+t v(x)-t^2Ev(x)\). Boundedness again forces \(v=0\). In both cases the
resulting affine deformation contradicts the two distinct determinants on positive-area sets.
\phantomsection\label{proofend:p12:floors}\end{proof}

\subsection{An incompressible phase}
\begin{proposition}\label{p12:incomp}
Let the fraction $\theta$ have $j=1$ and the fraction $f=1-\theta$ have $j\ge j_0>0$, with $0<\theta<1$. In the plane let both phases obey the preceding shear cap $k$. Then $m=(J-\theta)/f\ge j_0$ and
\[
 \gamma^2+\frac{(4+k^2)\theta f}{J}(1-\sqrt m)^2\le k^2 .
\]
For $(\theta,j_0,k)=(1/3,1/5,1)$, the entire freely periodic interval at $J=1/2$ is $[-2/3,2/3]$.
In space impose $M=J^{1/3}S_\gamma$, with $S_\gamma=I_3+\gamma e_1\otimes e_2$, and impose on both phases the cap
$\max(\sum_i\nu_i,\sum_i\nu_i^{-1})\le L_3=1+\sqrt{4+k^2}$ on the stretches $\nu_i=\lambda_i/j^{1/3}$. Then
\[
 1+\sqrt{4+\gamma^2}\le L_3\min\left(\frac{\theta+fm^{1/3}}{J^{1/3}},\frac{\theta+fm^{2/3}}{J^{2/3}}\right).
 \tag{12.3}\label{p12:eq3}
\]
For the same phase data and $127/256\le J\le1$ this gives the exact largest $|\gamma|$ for freely oriented periodic cells.
\end{proposition}
\begin{proof}\phantomsection\label{proof:p12:incomp}
The two planar phase means of $j$ are $1,m$. Averaging $P(F)$ as in Theorem~\ref{p12:floors} gives the bound after squaring. At half area the layers
$2^{-1/2}\bigl(\begin{smallmatrix}1&t\\0&2\end{smallmatrix}\bigr)$ and
$2^{-1/2}\bigl(\begin{smallmatrix}1&t/2\\0&1/2\end{smallmatrix}\bigr)$, $|t|\le1$, have the required areas $1,1/4$, squared shear $(1+t^2)/2$ and mean shear $2t/3$.

In space the entries of $\operatorname{cof}F$ and $\det F$ are minors. Their integrals equal the minors of $M$, by integration by parts for smooth maps followed by Lipschitz approximation (equivalently, differentiate each minor along $M+t\nabla v$ and integrate its divergence). The nuclear norms of $F$ and $\operatorname{cof}F$ are bounded by $L_3 j^{1/3}$ and $L_3 j^{2/3}$, respectively. Convexity of the two norms and concavity of the powers of $j$ give \eqref{p12:eq3}, since the stretches of $S_\gamma$ are $b,1,b^{-1}$ with $b-b^{-1}=|\gamma|$.

For sharpness put $a=J^{1/3}$, $q=m^{1/3}$, $A=\theta+fq$, $B=\theta+fq^2$, $w=B/(aA)$, $\rho=\min(w,1/w)$, and $L_3=1+\sqrt5$. The polynomial $c(x)=x^3-L_3x^2+\rho L_3x-1$ has positive roots $x_1\ge x_2\ge x_3$ with $x_{i+1}\le qx_i$ on the interval of the statement. Here are the sign checks. We have $5/8\le q\le1$ and
$w^3=3(1+2q^2)^3/[(1+2q^3)(1+2q)^3]$.
Differentiation gives the sign of $(\log w)'$ as that of $2q^2-1$ for $q<1$. Hence the minimum of $\rho^3$ is at $5/8,1/\sqrt2$ or $1$, and the first two values are $6858/6859$ and $48/(24+17\sqrt2)$. Thus $9997/10000<\rho\le1$. Substitution gives $c(0)<0$, $c(31/50)>0$, $c(124/125)>0$, $c(1)\le0$, $c((1+\sqrt5)/2)\le0$. The coefficient of $\rho$ is positive for $x>0$, so the positive signs use its lower bound and the two nonpositive signs use $\rho=1$. The roots consequently satisfy $x_3<31/50$, $124/125<x_2\le1$ and $x_1\ge(1+\sqrt5)/2$, which give the desired interlacing. The derivative and all endpoint sign comparisons above are asserted in \path{certificates/baseline/p12.py}.

For $w\ge1$ use the roots, and for $w<1$ their reciprocals, in decreasing order $d_i$. Give the two phases the spectra $d_i$ and $e_i=q d_i$. For $q<1$ the symmetric positive matrices $H_1=\operatorname{diag}(d_i)$ and $H_2=H_1-vv^T$ have these spectra if
$v_i^2=\prod_j(d_i-e_j)/\prod_{j\ne i}(d_i-d_j)$. Interlacing makes these numbers nonnegative, and evaluation of the characteristic polynomial at $d_i$ verifies the formula. Both phases meet the cap. All minors interpolate on this rank-one segment. The mean $H$ has determinant $J$, trace $A\sum d_i$, and cofactor trace $B\sum d_i^{-1}=a\operatorname{tr}H$. Thus $H/a$ has determinant one and equal trace and cofactor trace, giving one unit eigenvalue and two reciprocal eigenvalues. Proper rotations on the left and right give the specified simple shear at equality in \eqref{p12:eq3}. The common tangential columns give compatible layers. Their map is injective because the positive determinants make the normal coordinate of the images strictly increasing. At $q=1$ uniform deformation suffices.
\phantomsection\label{proofend:p12:incomp}\end{proof}

\subsection{Rigid inclusions}
Here a reference fraction \(\theta\in(0,1)\) must have \(F\in SO(2)\). The matrix phase admits
\(j\ge j_0>0\) and stretch ratio \(\lambda_1/\lambda_2\le K\), where \(K>1\). Let \(b\ge a>0\) be
the singular values of \(M\), so \(ab=J\), and write
\[
p=\frac{2K}{K-1},\qquad
 \nu=\frac{a/b-1/K}{1-1/K},\qquad
 b=\sqrt J\,\frac{\sqrt{4+\gamma^2}+|\gamma|}{2}.
\]
\begin{proposition}\label{p12:rigid}
Every periodic or full affine specimen satisfies
\[
\qquad \theta\le b^p\nu .\qquad                       \tag{12.4}\label{p12:eq4}
\]

In the sector \(b\le1,\ j_0\le b^2/K\), \eqref{p12:eq4} is necessary and sufficient over freely oriented
periodic cells, including equality. It gives \(J\ge\theta^{1-1/K}\) in pure compression. The
pure-compression endpoint is also attained on any full affine specimen provided $j_0\le b^2/K$.
\end{proposition}
\begin{proof}\phantomsection\label{proof:p12:rigid}
For necessity define the positive-sign Burkholder function
\[
\mathcal B_p(F)=\big[\tfrac p2\det F-(\tfrac p2-1)\|F\|_\mathrm{op}^2\big]\|F\|_\mathrm{op}^{p-2}.
\]

It equals one in the rigid phase and is nonnegative on the matrix distortion cone
\(\|F\|_\mathrm{op}^2\le K\det F\). For periodic data write $y=Mx+w(x)$. The maps $Mx+w(jx)/j$ are $K$-quasiregular on the plane and converge locally uniformly to the affine map $Mx$. The limit is also $K$-quasiregular, since $\|M\|_{\rm op}^2\le\langle\|F\|_{\rm op}^2\rangle\le K\det M$, and is locally in $W^{1,p}$. Their gradients have the same distribution on a period parallelogram for every positive integer $j$. The lower-semicontinuity inequality for $-\mathcal B_p$ therefore gives $\langle\mathcal B_p(F)\rangle\le\mathcal B_p(M)$ \cite[Thm.~1.5, Eq.~(1.15)]{burkholder}. For full affine data the same bound on $\|M\|_{\rm op}$ follows from the mean and determinant identities. Enclose the specimen inside a period square, extend by $Mx$ on its complement and periodize the fluctuation. Matching traces give a Lipschitz map, and the distortion inequality holds on both pieces. Apply the periodic inequality and cancel the identical exterior term. Since $\mathcal B_p(M)=b^p\nu$, the result follows.

For sufficiency, use the singular bases of \(M\). Condition \eqref{p12:eq4} gives \(0<\nu\le1\). Split
\(\operatorname{diag}(b,a)\) between \(bI\) in fraction \(\nu\) and \(\operatorname{diag}(b,b/K)\)
in fraction \(1-\nu\), with common first column. The second material region consists entirely of
matrix. In the first region fill the open layer up to measure zero by disjoint disks, repeating
the pattern with the cell. In a unit-radius copy put \(\rho=b^{K/(K-1)}\). Map its core of radius
\(\rho\) isometrically and its annulus by preserving angle and sending radius \(r\) to
\[
f(r)=\rho(r/\rho)^{1/K},\qquad \rho\le r\le1 .
\]

This matches scale \(b\) on the boundary of the disk. The annular stretches are \(f(r)/r\) and
\(f(r)/(Kr)\), so the stretch ratio is \(K\) and the smallest determinant is \(b^2/K\). Scaled,
translated copies match the background affine layer on every circle. The maps are uniformly
Lipschitz and injective since each disk maps onto its affine disk image. The cores occupy total
fraction \(\nu\rho^2=\nu b^p\). Assign any required fraction \(\theta\) of those cores to the
rigid phase. The rest also obey the matrix law. Proper rotations give the prescribed \(M\). For
\(\nu=1\) the same disk filling applies directly to the entire finite domain. For \(b=1\) the
conformal region itself is isometric.
\phantomsection\label{proofend:p12:rigid}\end{proof}

\subsection*{A sharp three-dimensional energy-budget example}
Pure-load tests describe the phase laws only along the tested paths.
When both laws and the loading path are fixed, upper and lower bounds can
instead meet. Here $\gamma$ denotes engineering shear and energy is
measured per unit reference volume.

\begin{theorem}[Shear after axial compression]\label{new:rubber}
Mix equal volumes of incompressible neo-Hookean phases with shear
moduli $\mu$ and $10\mu$. First compress to
$D_\lambda=\operatorname{diag}(\lambda^{-1/2},\lambda^{-1/2},\lambda)$
with $\lambda=9/10$, then shear to
$(I+\gamma e_1\otimes e_3)D_\lambda$. Impose mean energy at most
$\mu/4$ throughout a continuous, globally injective path on a fixed
periodic layout. Over all such layouts the sharp shear supremum obeys
\[
           0.53163717318 < \Gamma < 0.53163717320 .
\]
Every $0\le\gamma<\Gamma$ is feasible. No piecewise-smooth layout
attains $\Gamma$; attainment in the arbitrary measurable class is open.
\end{theorem}

\begin{proof}\phantomsection\label{proof:new:rubber}
Set $\mu=1$, $b=\lambda^{-1/2}$. The local energy in a phase of
modulus $k$ is $k(|F|^2-3)/2$ on $\det F=1$. If $x\ge b\ge z$ are
the singular values of the imposed gradient, then
$d=x-z=\sqrt{(b-\lambda)^2+\lambda^2\gamma^2}$. With the two
outer singular directions chosen consistently, the improper orthogonal
matrix $R$ with signed singular basis $(1,1,-1)$ gives
\[
 |F|^2-2R:\operatorname{cof}F
 =(\operatorname{tr}R^TF)^2+
   \sum_{i<j}\big((R^TF)_{ij}-(R^TF)_{ji}\big)^2 .
\]
Both the gradient and cofactor have prescribed means. Using their two
linear tests and minimizing the hard-phase excess over its middle stretch
leads to the following scalar bound. For $r>1$ put
\[\begin{gathered}
 u=11r-10r^{-2},\quad \tau=2r-r^{-2},\quad D=(u+\tau)/2-b,
\\
 L(r)=b^{-1}-bD-\tfrac32+\tfrac{u^2+\tau^2}{4}
                     +\tfrac94(2r^2+r^{-4}-3).
\end{gathered}\]
At $d=D(r)$ the mean energy is at least $L(r)$. Its derivative there
satisfies $L'(r)=(u-b)D'(r)$. The equation $L(r_*)=1/4$ has its
root in $(1.03278619490856271,1.03278619490856273)$; keeping its
test $u(r_*)$ fixed excludes all larger shears. Thus
\[
       \Gamma=\lambda^{-1}\sqrt{D(r_*)^2-(b-\lambda)^2}.
\]
For the matching sequence the hard phase has stretches
$(r,r,r^{-2})$ and the soft phase satisfies the trace equality in the
test. Two rank-one mixing steps realize the endpoint energy as a limit.
For $r<r_*$ sufficiently close to the root its strict energy slack pays
for the transition and boundary regions. Switching successive rows of
one fixed periodic layout keeps the entire compression-then-shear path
inside the budget and preserves injectivity. The detailed path and the
endpoint nonattainment argument are in Supplement~\ref{sup:rubber}.
\phantomsection\label{proofend:new:rubber}\end{proof}

\section{13. The shear interval of an isotropic polycrystal}\label{problem-13}
\Milton{13}

Let $C$ be any positive-definite elastic stiffness on $V=\mathrm{Sym}(3)$, $S=C^{-1}$, and $\lambda I\preceq C\preceq\Lambda I$. A grain has stiffness $C_R$, obtained by rotating $C$ by $R\in SO(3)$. Write $G$ for the closure of effective tensors of periodic measurable orientation fields. Put $e=I_3/\sqrt3$, $J=e\otimes e$, $D=I_V-J$ and use eigenvalue coordinates $a=3\kappa$, $b=2\mu$. The Hill bulk endpoints in these units are $a_R=(e:Se)^{-1}$ and $a_V=e:Ce$ \cite[\S22.3, p.~442]{book}. Avellaneda and Milton proved both attainable by multilayer laminates \cite[p.~824]{avellaneda}. Here the shear is prescribed along with the bulk.

Four one-parameter families of weighted cell minima determine the region. For a sign $+$ take a compatible strain, for $-$ an equilibrated stress. In a unit periodic cube take six such fields, $X_0$ of mean $e$ and type $\epsilon$, and $X_j$ of means $E_j$, $1\le j\le5$, of type $\eta$, where the $E_j$ form an orthonormal basis of $DV$. Define
\begin{equation}\label{p13:support}
 H_{\epsilon\eta}(s,t)=\inf\int\min_{R\in SO(3)}
 \left(sX_0:C_R^\epsilon X_0+\frac t5\sum_{j=1}^5 X_j:C_R^\eta X_j\right),\qquad s,t\ge0 .
\end{equation}
Exponents stand for $+1,-1$. The same rotation serves all six loads at a point. Fourier polynomials suffice in the infimum, with each nonconstant strain amplitude in $\{k\odot v:v\in\mathbb R^3\}$ and each stress amplitude annihilating $k$ at frequency $k$, separately for its real and imaginary parts. Indeed truncation preserves those constraints and gives strong $L^2$ convergence. The integrand is a minimum of uniformly bounded quadratic forms, so the energies converge.

\begin{lemma}[Local homogenization]\label{p13:locality}
In dimension two or three let $G$ be the periodic response closure of the rotations of a
positive-definite crystal. Any measurable periodic tensor field with
values in $G$ has its effective tensor in $G$. The same assertion holds
for the closure $G_f$ of two fixed finite positive-definite phases at fraction $f$. More generally, let a unit period cell be partitioned up to null sets into finitely many regions $R_i$ with boundaries of measure zero. For $0\le f_i\le1$ take $A_i\in G_{f_i}$, allowing pure phases at the endpoints. Then the tensor of the medium taking value $A_i$ on $R_i$ belongs to $G_{\sum_i |R_i|f_i}$.
\end{lemma}
\begin{proof}\phantomsection\label{proof:p13:locality}
All coefficients and their inverses have common positive lower and
finite upper bounds. First suppose the field is constant on each of
finitely many dyadic cubes. For a finite basis of strain and stress
loads approximate the minimizing macroscopic fields by smooth admissible fields, extending the approximations linearly in the load. Subdivide until each smooth field varies by at
most $\xi$ in a subcube $B$, where we use its center value as the
local load. Choose a child cell approximating the tensor on $B$ and
trigonometric admissible correctors for these local loads with trial
energies within $\xi$ of the respective child minima. Use whole fine
repeats in $B$. Dyadic subcubes and cubic child cells permit an exact
tiling at common integer frequencies. On a cube assigned fraction $f_i$, every repeat has exactly that fraction.

For either field type, multiply its mean-zero corrector $v(kx)$ by
$1_B(x)$ and project onto the mean-zero admissible space. This changes
it by $o(1)$ in $L^2$ at fixed partition as $k\to\infty$. Indeed a
trigonometric polynomial multiplier moves each admissible mode $k\nu$
by finitely many offsets $j$, and the wave projector satisfies
$\Pi(k\nu+j)\to\Pi(\nu)$. Approximate $1_B$ in $L^2$ using boundedness
of the fixed $v$. Add the projected fluctuations to the original
smooth macroscopic fields. Periodic averaging inside each $B$ gives
the selected child energies at the center loads. Replacing the center
values by the smooth macroscopic fields costs $O(\xi)$ by the quadratic
energy bound and Cauchy-Schwarz. The squared norms of the correctors
are bounded in terms of the bounded center loads by uniform ellipticity.
First take enough repeats for this fixed partition, then let the
approximation errors tend to zero. The minimizing fields in each cell depend linearly on their loads. An admissible error above a minimizer costs exactly its squared energy norm by stationarity. Thus choosing the approximations on a basis controls the excess on all unit loads. For the resulting sequence of tensors $A_k$ and the target effective tensor $A_*$, the two tests give, for every $E$,
\[
 \limsup_k E:A_k E\le E:A_*E,\qquad
 \limsup_k E:A_k^{-1}E\le E:A_*^{-1}E .
\]
Uniform ellipticity and inversion give $A_k\to A_*$. For the finite partition $R_i$, exhaust each interior by whole dyadic subcubes. On the omitted part of each $R_i$ fill fraction $f_i$ with phase one. The total fraction is then exactly $\sum_i|R_i|f_i$. The smooth macro fields contribute vanishing energy on the omitted parts as their measure tends to zero, and the same localized-corrector estimates apply on the retained cubes. If the chosen child sequences already converge for a finite list of phase laws on the same geometries, use the union of their strain and stress trial families before choosing the repeat count. This gives the joint substitution on that list.

For a measurable $G$-valued field select nearest members of the compact
set $G$ to its dyadic cell averages. These piecewise constant choices
converge in $L^2$ by conditional averaging and hence along a subsequence
almost everywhere. Testing both energies against the original minimizing
fields gives convergence of the effective tensors by dominated convergence
with the fixed squared fields as weights. For the crystal orbit the
closure is rotation invariant. A rational rotation of a cubic cell
has a cubic supercell, and the general case follows by density of
rational rotations and continuity of the rotation action on tensors.
\phantomsection\label{proofend:p13:locality}\end{proof}

\begin{lemma}[Four support inequalities]\label{p13:four}
Let $\mathcal R\subset(0,\infty)^2$ be nonempty and compact. For
$\epsilon,\eta\in\{1,-1\}$ write
$P_{\epsilon\eta}(a,b)=(a^\epsilon,b^\eta)$. Assume that every
convex combination of two points of $P_{\epsilon\eta}(\mathcal R)$
dominates a point of the same set, coordinate by coordinate, for
each sign pair. Also assume that any two points of $\mathcal R$
can be joined within it by a continuous path whose coordinates
stay between their respective endpoint values.
Then $(a,b)\in(0,\infty)^2$ belongs to $\mathcal R$ precisely when,
for all four sign pairs and all $s,t\ge0$,
\[
 s a^\epsilon+t b^\eta\ \ge \min_{(\alpha,\beta)\in\mathcal R}(s\alpha^\epsilon+t\beta^\eta).
\]
The projection onto the first coordinate is a closed interval, and
over each of its points, including its ends, the section in the
second coordinate is a closed interval.
\end{lemma}
\begin{proof}\phantomsection\label{proof:p13:four}
Each upper set $P_{\epsilon\eta}(\mathcal R)+[0,\infty)^2$ is closed
by compactness and convex by the first hypothesis. A separating
linear functional bounded below on this set has nonnegative
coefficients. Thus its membership test is exactly the stated
inequality for that sign pair.

For a passing pair, fix $\eta$. The two choices of $\epsilon$ give
points of $\mathcal R$ on opposite sides of the proposed $a$, both
with second observable at most $b^\eta$. The path between them reaches
$a$ with the same bound. Applying this once for each $\eta$ gives
two attainable second coordinates bracketing $b$ over exactly $a$.
Their path has first coordinate fixed and reaches $b$. Conversely
every member satisfies the inequalities by definition. The path
property fills the first projection and every vertical section,
whose ends exist by compactness.
\phantomsection\label{proofend:p13:four}\end{proof}

\begin{theorem}\label{p13:interval}
Every bulk eigenvalue $a_R\le a\le a_V$ occurs in the isotropic closure and these are all the possible bulks. For such $a$ put
\[
 l_\eta(a)=\sup_{\epsilon=\pm1,\ r\ge0}
       [H_{\epsilon\eta}(r,1)-r a^\epsilon].
\]
The closed shear interval consists of precisely the $\mu$ with
$l_+(a)/2\le\mu\le1/[2l_-(a)]$. Here $l_-(a)\ge H_{+-}(0,1)\ge1/\Lambda$. Every point has a realizing sequence of periodic orientations of the given crystal.
\end{theorem}
\begin{proof}\phantomsection\label{proof:p13:interval}
The cell strain and complementary stress principles give
$\lambda I\preceq A\preceq\Lambda I$ for $A\in G$, hence compactness.
Constant trials bound the stiffness and compliance by their respective
arithmetic means.

\emph{Isotropization.} Write $\Pi B$ for the rotation average of $B$. Laminate copies of $A_0=A$ in one direction with their orientations distributed by normalized Haar measure. One can realize this distribution by a measurable function of the layer coordinate, so Lemma~\ref{p13:locality} applies. Repeat the operation, calling the successive tensors $A_j$. With $U_j=\Pi A_j$, $L_j=(\Pi A_j^{-1})^{-1}$, the two constant-field trials give $L_j\preceq A_{j+1}\preceq U_j$. Thus $U_j$ decreases and $L_j$ increases. Along any jointly convergent subsequence $A_j\to B$, $A_{j+1}\to B'$, the traces of the stiffness and compliance have consecutive differences tending to zero. Equality in trace in the two matrix bounds gives $B'=\Pi B=(\Pi B^{-1})^{-1}$. Equality of the matrix arithmetic and harmonic means forces all matrices to agree almost surely. Indeed, with mean $P$ the average of $(B_R-P)B_R^{-1}(B_R-P)$ is $P\langle B_R^{-1}\rangle P-P$. The nonnegative integrand must vanish almost surely at equality, and continuity in $R$ gives $B=\Pi B$. The limiting brackets thus coincide at an isotropic tensor in $G$. Its hydrostatic and averaged deviatoric energies, in either convention, are at most the respective energies of the input. A common hydrostatic eigenvalue of all copies is preserved by the two brackets.

\emph{Bulk endpoints and interpolation.} The constant trials on $e$ give the two Hill bounds. To realize them start from $C$ on the strain side or $S$ on the stress side, and diagonalize its response to $e$ with eigenvalues $x\le y\le z$ and mean $d_0$. Mix with a rotated copy exchanging the first and last axes in fraction $(d_0-x)/(z-x)$ when $x<z$, then exchange axes 2 and 3 in equal parts. The hydrostatic loading is preserved because the response jump is proportional to $n_i n_i^T-n_j n_j^T$, which annihilates the unused axis and is also $(n_i+n_j)\odot(n_i-n_j)$. Use the unused axis as normal for a stress jump, and $(n_i+n_j)/\sqrt2$ for a strain jump. The eigenvalues become $(d_0,y,2d_0-y)$ and then $(d_0,d_0,d_0)$. Isotropizing gives the required endpoints.

Between two isotropic materials there is a continuous path with each modulus between its two original values, fixing any shared modulus. Here is an explicit construction. Let $B$ be the host and $A$ the core. Adding host layers of normal $n$ to a tensor $B+L$ with retained share $u$ gives an increment $L'$ with
$u(L')^{-1}=L^{-1}+(1-u)\Gamma_{B,n}$, where $U_n v=n\odot v$ and $\Gamma_{B,n}=U_n(U_n^t B U_n)^{-1}U_n^t$. This follows by minimizing over the layer jump $U_n v$ and multiplying out the inverse. Continuity of the energy minimum extends the resulting formula to singular increments. To retain core fraction $p$, start at $v_0=1$ and set $v_i=v_{i-1}-(1-p)c_i$, $u_i=v_i/v_{i-1}$, for a convex list of normals with weights $c_i$. Telescoping gives
\[
 B+p[(A-B)^{-1}+(1-p)\textstyle\sum_i c_i\Gamma_{B,n_i}]^{-1}.
\]
Choose $c_i$ so that the sum equals the spherical average (the compact convex hull contains the average). For isotropic $A,B$ this tensor is isotropic and continuous in $p$, including singular contrasts by the direct layer formula. The constant trials put both coordinates between their two constituent values. This proves the assertion about the path, fills the bulk interval and connects any two attainable shears at a common bulk.

\emph{Separation.} In \eqref{p13:support} the first rotation meeting a given tolerance in a countable dense list is a measurable near-minimizer. Taking these two infima in either order shows that $H$ minimizes $s\,e:A^\epsilon e+(t/5)\operatorname{Tr}(D A^\eta)$ on $G$. Isotropization can only reduce both terms, so $H$ is the minimum of $s a^\epsilon+t b^\eta$ on the compact isotropic region $\mathcal R$. For each sign pair, laminate two isotropic tensors in the prescribed
proportions and isotropize. The constant strain and stress trials give
a point whose signed observables are at most their respective arithmetic
means. The preceding interpolation supplies the path required by
Lemma~\ref{p13:four}. Its four inequalities therefore characterize
$\mathcal R$. The tests with $t=0$ impose the Hill interval, and
taking $r=s/t$ for $t>0$ gives the two stated endpoint formulas.
\phantomsection\label{proofend:p13:interval}\end{proof}
\subsection*{Which shapes can the region have?}
The four support families determine a vertical interval at each bulk.
There are only four possible shapes for their union. Let $e=I_3/\sqrt3$,
$J=e\otimes e$, $D=I_V-J$ on $V=\mathrm{Sym}(3)$ and $T=D-2J$.

\begin{theorem}[Shape of the polycrystal region]\label{new:shapes}
Apply the following tests in order to a positive-definite single-crystal
stiffness $C$. Each interval in the second and third rows has positive
length.
\begin{center}\small
\begin{tabular}{@{}p{.53\linewidth}p{.40\linewidth}@{}}\toprule
Crystal condition & Isotropic $(\kappa_*,\mu_*)$ region\\ \midrule
$C=aJ+bD$ & One point $(a/3,b/2)$.\\
$C$ anisotropic and $Ce=ae$ & A vertical segment at $\kappa_*=a/3$.\\
$Ce$ nonhydrostatic and $C=bT+H\otimes H$ for $b>0$, $H\in V$
& A horizontal segment from $\kappa_R$ to $\kappa_V$ at $\mu_*=b/2$.\\
All other crystals & A closed topological disk with positive-length vertical sides at both Hill bulks.\\ \bottomrule
\end{tabular}
\end{center}
\end{theorem}
\begin{proof}\phantomsection\label{proof:new:shapes}
An isotropic crystal retains its moduli under every rotation. The
arithmetic and harmonic bounds on the hydrostatic load coincide exactly
when that load is an eigenvector of $C$; a nonisotropic crystal with
this property admits different isotropic shears. For the rank-one
identity in the third row, elimination of the normal traction in a
layer preserves the fixed isotropic shear. A laminate of rotated grains
at either Hill bulk gives the two ends and interpolation fills the
segment.

In the remaining case the hydrostatically uniform constructions on each
Hill side give two distinct shears. Mix either segment with a material
at the opposite bulk. The core fraction at prescribed intermediate bulk
is independent of the core shear when the host is fixed; the final
shear is strictly increasing in the core shear. Thus every vertical
section has positive length, uniformly over the closed bulk interval.
The endpoint graphs are continuous by compactness and bulk repair.
Linear interpolation between the graphs is a homeomorphism from a
rectangle to the region. Supplement~\ref{sup:crystal} gives the
construction and algebraic exceptional-case tests underlying the side
segments.
\phantomsection\label{proofend:new:shapes}\end{proof}

\subsection{Two-sided bounds on each support}
Use the normalized periodic cube \(\mathbb T^3\). A constant field \(X\) consists of six symmetric
\(3\times3\) tensors. In each tensor slot use the standard axis-aligned Kelvin basis \(e_i e_i^T\)
and \((e_i e_j^T+e_j e_i^T)/\sqrt2\) for \(i<j\), with dimension \(d=36\) in all. The prescribed
hydrostatic and deviatoric mean loads are specified separately below. Fix the signs
\(\sigma=(\epsilon,\eta,\ldots,\eta)\). For nonzero frequency \(k\) write \(n=k/|k|\). The
admissible orthogonal projection \(\Pi(k)\) acts separately on the six loads. Its stress block and
strain block are respectively
\[
Q_n Y=(I-nn^T)Y(I-nn^T),\qquad I_V-Q_n .
\]

The projection is zero-homogeneous in frequency and continuous on the unit sphere. Its ranges
consist of the stresses with \(Yn=0\) and the strains \(n\odot v\). Write \(\mathcal F\) for the
real, mean-zero \(L^2\) fields with every nonzero Fourier coefficient in the prescribed range.
Haar measure on the torus has mass one. Define
\[
f(X)=\min_{R\in SO(3)}\sum_{\alpha=0}^5 X_\alpha:C_R^{\sigma_\alpha}X_\alpha,\qquad
 W(X)=\inf_{F\in\mathcal F}\int f(X+F).
 \tag{13.1}\label{p13:eq1}
\]

Take rational \(0<c\le L\) bounding all eigenvalues of both \(C\) and \(C^{-1}\) between them.
Then \(c|X|^2\le W(X)\le f(X)\le L|X|^2\). With \(X_0=\sqrt{s}\,e\) and \(X_i=\sqrt{t/5}\,E_i\)
for the five orthonormal deviatoric loads, \(W(X)=H_{\epsilon\eta}(s,t)\). Keeping the weights in
\(X\) allows zero weights without changing \(c\).

\begin{lemma}[Regularity of the relaxed energy]\label{p13:regularity}
For the integrand and wave constraints of \eqref{p13:eq1}, at every mean $X$,
\[
\int W(X+F)\ \ge W(X)\qquad(F\in\mathcal F).                 \tag{13.2}\label{p13:eq2}
\]
\[
W\in C^{1,1}(\mathbb R^d),\qquad
 \operatorname{Lip}(DW)\le M_0=2(d-1)L.                       \tag{13.3}\label{p13:eq3}
\]
\end{lemma}
\begin{proof}\phantomsection\label{proof:p13:regularity}
For each fixed orientation field the minimum over its six correctors is a quadratic form of the
constant \(X\), bounded above by \(L|X|^2\) using constant trials. Minimizing these forms over
orientation fields gives \eqref{p13:eq1}. One can select a measurable pointwise near-minimizing rotation
from a countable dense list, and then exchange the two infima. Consequently \(W-L|\cdot|^2\) is
concave, \(W\) is homogeneous of degree two, and
\[
|W(X)-W(Y)|,\ |f(X)-f(Y)|\ \le L(|X|+|Y|)|X-Y|.
\]
To prove \eqref{p13:eq2}, approximate \(F\)
in \(L^2\) by an admissible trigonometric polynomial. Partition the cube into small cubes \(B_i\)
on which \(X+F\) differs from a constant \(X_i\) by at most \(\xi\). For each \(i\), take a
mean-zero admissible trigonometric polynomial \(V_i\) with \(\int f(X_i+V_i)<W(X_i)+\xi\). The
bound \(\|V_i\|_2^2\le (L/c-1)|X_i|^2+\xi/c\) holds uniformly in the partition because the given
trigonometric field is bounded. Project \(\sum_i 1_{B_i}(x)V_i(hx)\) onto \(\mathcal F\), for
positive integers \(h\). By the projection estimate in Lemma~\ref{p13:locality}
the change tends to zero in \(L^2\) at fixed partition. Periodic averaging gives
\(\int_{B_i} f(X_i+V_i(hx))\to |B_i|\int f(X_i+V_i)\), and the analogous identity for \(|V_i|^2\).
Add the projected corrector to the original \(F\) in \eqref{p13:eq1}. The quadratic continuity bound costs
\(O(\xi)\) after sending \(h\to\infty\). Sending \(\xi\to0\), then using \(L^2\) density for the
original field, proves \eqref{p13:eq2}.

Every one of the \(d\) individual Kelvin/load basis directions is a wave at an integer axial
frequency. For a diagonal strain take its own axis, and for an off-diagonal strain take either
index axis. For each stress take an axis outside its indices. Two-valued, one-directional periodic
fluctuations in \eqref{p13:eq2} prove convexity on each coordinate line. Thus \(\Delta W\ge0\) in
distributions. Convolve \(W\) with a smooth mollifier in field space. Its Hessian has eigenvalues
at most \(2L\) by semiconcavity, and its trace is nonnegative. Each eigenvalue is therefore at
least \(-2(d-1)L\). The smooth gradients are equicontinuous and locally bounded, and their limits are
the distributional derivatives of \(W\). Thus \eqref{p13:eq3} follows, and two-homogeneity gives \(DW(0)=0\).

\phantomsection\label{proofend:p13:regularity}\end{proof}

\begin{lemma}[Finite-frequency reduction]\label{p13:band}
Fix any of the six-load projection spaces $\mathcal F$ defined above.
For each integer \(m\ge1\) put
\[
N_m=2^{8m^3+4m+1},\qquad
 E_m=\frac1m+\frac1{2m^2}+3\,2^{-m}+2\,2^{-2m}+26\,2^{-4m}.
 \tag{13.4}\label{p13:eq4}
\]

Suppose \(g:\mathbb R^d\to\mathbb R\) has globally \(M\)-Lipschitz gradient and satisfies the following
inequality at every constant \(A\) and every real admissible mean-zero polynomial \(P\) with
frequencies \(0<|k|\le N_m\),
\[
\int_{\mathbb T^3} g(A+P)\ge g(A).
\]

Then every \(F\in\mathcal F\) satisfies
\[
\int g(X+F)\ge g(X)-M E_m\|F\|_2^2.    \tag{13.5}\label{p13:eq5}
\]

\end{lemma}
\begin{proof}\phantomsection\label{proof:p13:band}Smooth the admissible field with a periodic approximate identity at a
scale chosen by an averaging argument. The projection onto the allowed
Fourier fields commutes with convolution. Bounds on the discarded high
frequencies and on the quadratic growth control the energy error at the
specified cutoff. The detailed estimates are in Supplement~\ref{sup:crystal}.\phantomsection\label{proofend:p13:band}\end{proof}

Enumerate polynomials \(p\) in \(d\) variables with rational coefficients. Put
\[
h(0)=0,\qquad h(X)=|X|^2p(X/|X|)\quad (X\ne0).
\]

The certificate also gives rationals \(\tau,M>0\) and integers \(m,K\ge1\). Let
\(g=h-\tau|\cdot|^2\), \(N=N_m\), and \(q=16 M N(2N+1)^3/K\). Check the following finite
conditions.

\begin{enumerate}\item On the field-space unit sphere \(h(X)\le\sum_\alpha X_\alpha:C_R^{\sigma_\alpha}X_\alpha\) for
every \(R\in SO(3)\).
\item On that sphere \(-MI\preceq D^2g(X)\preceq MI\). Require \(M E_m\le\tau\) and \(K>2N\).
\item For all constant \(A\) and all real admissible mean-zero trigonometric polynomials \(P\)
supported in \(0<|k|\le N\), check\end{enumerate}
\[
\frac1{K^3}\sum_{z\in\{0,\ldots,K-1\}^3/K}
           [g(A+P(z))-g(A)]\ \ge q\|P\|_2^2 .                \tag{13.8}\label{p13:eq8}
\]

\begin{lemma}[Finite lower certificates]\label{p13:certificate}
For real algebraic crystal data with isolating information the three
conditions above are decidable. Every $h$ satisfying them obeys
\[
h(X)\le W(X)\quad\text{for every }X.   \tag{13.9}\label{p13:eq9}
\]
\end{lemma}
\begin{proof}\phantomsection\label{proof:p13:certificate}The displayed test on trigonometric inputs is a finite polynomial test
on their real and imaginary coefficients. Truncation and the Hessian
bound carry its error to arbitrary admissible fields. The estimates,
including the norm in which the error is controlled, are given in
Supplement~\ref{sup:crystal}.\phantomsection\label{proofend:p13:certificate}\end{proof}

\begin{lemma}[Exhaustion of the certificates]\label{p13:exhaustion}
For the envelope \eqref{p13:eq1} and the polynomials accepted by
Lemma~\ref{p13:certificate}, at every $X$,
\[
W(X)=\sup_{\text{accepted }h}h(X).           \tag{13.12}\label{p13:eq12}
\]
\end{lemma}
\begin{proof}\phantomsection\label{proof:p13:exhaustion}Regularize the relaxed cell energy and approximate its restriction to the
unit sphere together with its first two derivatives by a polynomial.
Positive homogeneity recovers the trial minorant on the full field space.
The slack in the finite-frequency test absorbs the smoothing and rational
approximation errors. Supplement~\ref{sup:crystal} gives the quantified
approximation and the order of the choices.\phantomsection\label{proofend:p13:exhaustion}\end{proof}

\begin{proposition}\label{p13:computeH}
For positive-definite $C$ specified by real algebraic entries with isolating data, every support \eqref{p13:support} at nonnegative algebraic weights can be bracketed to any prescribed positive rational error by finite rational or real-algebraic tests.
\end{proposition}
\begin{proof}\phantomsection\label{proof:p13:computeH}
For fixed trigonometric trials, dyadic piecewise rotations from a dense
algebraic list approximate a measurable orientation in measure and
give converging energies. Fourier truncation preserves admissibility,
and rational coefficients in a rational Cartesian basis of each wave
space are dense. Enumerating these trials and orientations gives upper
bounds at the prescribed algebraic mean. Conservative Riemann sums on
the dyadic pieces compute each trial energy to any tolerance using
derivative bounds for its fixed trigonometric polynomials.
Lemmas~\ref{p13:certificate} and~\ref{p13:exhaustion} give an
exhaustive enumeration of lower bounds. Run the two enumerations
until the gap is below the prescribed tolerance, using one-sided
rounding at each evaluation.
\phantomsection\label{proofend:p13:computeH}\end{proof}

Already at $m=4$ the cutoff is $N_4=2^{529}$. The quantifier
elimination test includes Fourier coefficients for every lattice point
in this frequency ball, far beyond practical computation.

\subsection*{Inward estimates at the two Hill bulks}

Let $e=I/\sqrt3$, let $D$ be orthogonal projection onto the five-dimensional
space $\mathcal D$ of symmetric traceless three-by-three matrices, and take
$L=C$ at the Voigt bulk and $L=C^{-1}$ at the Reuss bulk. If $L_R$ denotes
the rotated tensor, set
\[
 \alpha=e:Le,\quad M=|DLe|,\quad g_R=D L_R e/M,\quad
 K=\operatorname{co}\{g_R:R\in SO(3)\}.
\]
We first suppose $M>0$. We call $g$ the hydrostatic deviator of a local
material when its response to the constant loading $e$ is $\alpha e+Mg$
in the chosen stiffness or compliance convention. All norms on
$\mathcal D$ are Frobenius norms.

\begin{lemma}\label{p13:steer}
Suppose $M>0$, and write $A_L=A$ in the stiffness convention and $A_L=A^{-1}$ in the compliance convention. For each $g\in K$ there is an $A$ in the crystal closure $G$ with $A_L e=\alpha e+Mg$. There is a universal effectively selectable integer $m$ such that, for each $0<u\le1/4$, rotation $R$ and $z\in\mathcal D$ with $|z-g_R|\le u$, one can choose $A\in G$ satisfying $A_L e=\alpha e+M(1-u)^{4m}z$ and, with $\varphi=(1-u)^{5m}$,
\[
 A\preceq\varphi C_R+(1-\varphi)\lambda_{\max}(C)I,\qquad
 A^{-1}\preceq\varphi C_R^{-1}+(1-\varphi)\lambda_{\max}(C^{-1})I .
\]
\end{lemma}

\begin{proof}\phantomsection\label{proof:p13:steer}
For a half-turn $U_n=2nn^t-I$, let $M_n h=U_n hU_n^t$.
The jump $M_n h-h$ has only mixed entries between $n$ and $n^\perp$.
It is traceless, of rank at most two. For a strain jump it is a symmetric
rank-one product with normal $n$. For a stress jump, choose any normal in its
kernel. Thus the constant hydrostatic loading in a laminate of a material
with its half-turned copy solves the transmission conditions because the response
jump is equilibrated in the stiffness convention and compatible in the
compliance convention. The laminate updates
$h\mapsto h+s(M_n h-h)$ with switched fraction $s$.

Let $a\ge b\ge c$ be the eigenvalues of $g_I$. Every $h\in K$ has
eigenvalues $x\ge y\ge z$ in $[c,a]$ and zero trace, by the Rayleigh
inequalities. Conversely every such diagonal triple is obtainable in two
layers. A half-turn about $(e_i+e_j)/\sqrt2$ swaps the $i,j$ entries.
For $x\ge b$ average $a,b$ with their interchange to get one entry $x$. The other two are $a+b-x,c$, between which $y,z$ lie since $z\ge c$.
For $x<b$ average $b,c$ to get $x$. The remaining entries $a,b+c-x$
enclose $y,z$ since $y\le a$. One more swap and a rotation give the
target. In particular $K$ contains a ball about zero. This proof works
for the deviator of any starting material, including when the target is zero.

For unit $h\in\mathcal D$, averaging $M_n h$ over the unit sphere gives
$h/5$ by expanding $U_nhU_n^t$ using $\langle n_i n_j\rangle=\delta_{ij}/3$
and $\langle n_i n_j n_k n_l\rangle=
(\delta_{ij}\delta_{kl}+\delta_{ik}\delta_{jl}+\delta_{il}\delta_{jk})/15$.
There exists a uniform $c_0>0$ such that the convex hull of these
$M_n h$ contains $h/5+c_0\mathbb B_{\mathcal D}$. Otherwise compactness
of the two unit spheres in $(h,A)$ and separation give
$A:M_nh=A:h/5$ for all $n$, since a continuous function with mean equal to
its maximum is constant. On expanding and homogenizing this identity in
$n$, the real polynomial $|n|^2$ would divide
$(n^tAn)(n^thn)$. The irreducible quadratic $|n|^2$ is prime over the
real polynomial ring, and neither nonzero traceless quadratic factor is
its multiple. This proves the assertion. A rational $c_0>0$ is found
by decreasing dyadic search, since the requirement
$\max_{|n|=1}A:(M_nh-h/5)\ge c_0$ for all unit $h,A\in\mathcal D$
is a first-order real algebraic test.

Choose $m$ with $5m c_0\ge2$. Follow
$p(t)=g_R+t(z-g_R)$, $0\le t\le1$. Put
$\lambda=-5m\log(1-u)$ and $h=p/|p|$. Then
$|p'|/|p|\le u/(1-u)\le\lambda c_0$.
At each $t$, choose a convex combination of at most six $M_n h$
equal to $h/5+p'/(\lambda|p|)$. Its generators $M_n-I$, multiplied
by $\lambda$, therefore carry $p$ to $p'-(4/5)\lambda p$.
On a time step $\Delta t$, apply the chosen half-turn layers successively,
with switched fractions $\lambda\Delta t$ times their weights at the
left endpoint. The product of their linear updates differs from the
Euler increment by $O((\lambda\Delta t)^2)$ in norm, uniformly in the
choices, since $\|M_n-I\|\le2$. Each update is a contraction. The
smooth path $\exp(-4\lambda t/5)p(t)$ has precisely this first-order
increment, so summing the local errors gives convergence to it on uniform
meshes. Let $\varphi_N$ be the product of the unswitched fractions on a mesh of $N$ time steps. The logarithm expansion gives $\varphi_N\to e^{-\lambda}=\varphi$. At least this weight $\varphi_N$ of the leaves remains in orientation $R$. Constant strain and stress trials give the two matrix bounds in the statement with $\varphi_N$ for each finite laminate. In a convergent subsequence of their tensors these bounds pass to the limit, as does the computed hydrostatic response. This yields the asserted element of $G$.
\phantomsection\label{proofend:p13:steer}\end{proof}

\begin{theorem}\label{p13:inward}
For algebraic positive definite crystal data, at each exact Hill bulk
the least attainable isotropic shear has computably enumerable upper
bounds converging down to it, and the greatest shear has computably
enumerable lower bounds converging up to it.
\end{theorem}

\begin{proof}\phantomsection\label{proof:p13:inward}
We describe a finite witness and prove its density below. On the unit
three-torus choose grain orientations $R(x)$ on a rational cubical
partition. Supply a traceless tensor field $z(x)$ whose nonconstant
Fourier modes are transverse stresses at Voigt, compatible strains at
Reuss, with unrestricted mean. Choose rational $0<u\le1/4$,
$\tau=(1-u)^{4m}$ and require $\tau z\in K$ everywhere. On a selected
union $G_{\rm good}$ of subcubes require $|z-g_R|\le u$. Use
$\ell=1-(1-u)^{5m}$ on $G_{\rm good}$ and $\ell=1$ otherwise.
For $T=C$ or $T=C^{-1}$ take five trial fields $X_i$ of orthonormal
deviatoric means, compatible strains or equilibrated stresses respectively.
For $\Delta_T=\lambda_{\max}(T)-\lambda_{\min}(T)$ define
\[
 Q_T=\frac15\sum_{i=1}^5\int\big[X_i:T_R X_i+\Delta_T\ell |X_i|^2\big].
\]
Then $Q_C/2$ bounds the minimum shear above and $1/(2Q_{C^{-1}})$
bounds the maximum below.

Indeed, at each point apply Lemma~\ref{p13:steer} to attain the response
$\alpha e+M\tau z$, using its two matrix bounds on $G_{\rm good}$ and
the whole convex hull elsewhere. Both local stiffness and compliance
are bounded above by $T_R+\ell\Delta_T I$ in their respective
conventions by constant strain and constant stress trials. For a measurable choice, observe that for each of the two flags the tensors satisfying the
response and energy inequalities form a nonempty compact fiber in the
compact crystal closure with closed graph. To select one, intersect the
fiber successively with the first closed ball in a finite cover of mesh
$2^{-j}$ for which the intersection with all earlier balls is nonempty.
Nonemptiness is a closed parameter condition by compactness. The nested
intersections give a Borel limit. The effective tensor of this measurable
field belongs to the same closure by local homogenization.

The constant applied loading solves the hydrostatic equation exactly
because its response field meets the required wave constraint. Two
whole-material swaps as in the proof of Lemma~\ref{p13:steer}, with target
zero for its mean deviator, leave $e$ a hydrostatic eigenvector of
eigenvalue $\alpha$. In both energy conventions the traces on the five
shear loads can only decrease in these swaps, and in isotropizing the
resulting tensor with rotated copies, by the constant trial fields.
Rotations preserve $e$ and the shear-space trace. The isotropic result
has the exact required bulk, and its averaged shear-space eigenvalue
in stiffness (respectively compliance) is at most $Q_C$ (respectively
$Q_{C^{-1}}$). This proves both inward estimates.

\smallskip\noindent\emph{Restoring the hydrostatic differential constraint.}
Take periodic cells approaching an isotropic tensor on the specified side.
Orientations may be taken piecewise constant on a rational grid, with
algebraic rotations, since such approximations converge in measure and testing
the two energies on fixed fields gives convergence of the cell tensor.
Write $g=g_{R(x)}$, $A_L$ for the effective operator in the chosen convention,
and $h=\alpha-e:A_Le\to0$. If $P_0$ is the orthogonal projection onto
the zero-mean fluctuations allowed for the applied loading, the variational
principle gives
$M^2\|P_0g\|_2^2\le\|L\|h$ by inserting the trial fluctuation
$-M P_0g/\|L\|$ into the constant-field energy.
On a nonzero frequency with unit direction $n$ put $U=nn^t$, $P=I-U$.
The projection $Q$ onto admissible traceless responses, keeping the mean,
acts at this frequency by
\[
 Q_V Z=PZP-\tfrac12\operatorname{tr}(PZ)P,\qquad
 Q_R Z=UZP+PZU .
\]
In coordinates $n=e_3$, $P_0$ retains the mixed and longitudinal entries
at Voigt, or the planar block at Reuss. Since $Z_{33}=-\operatorname{tr}
(PZP)$, $\|g-Qg\|_2^2\le3\|P_0g\|_2^2\to0$.

The following estimate permits contraction into $K$ even though the
projected field $Qg$ need only be close in $L^2$. For sufficiently small
$t>0$ with $\|g-Qg\|_2\le t^{36}$ there is a traceless admissible response
$z$ for which
\begin{equation}\label{p13:truncate}
 \operatorname{dist}(z,K)\le Ct\quad\text{a.e.},\qquad
 |z-g|\le t^{28}\quad\text{outside a set of measure }Ct^8 .
\end{equation}
All constants here are independent of the cell, and only their finiteness is needed. Set $w=(-\Delta)^{-2}(g-\bar g)$ on the torus.
There is a fourth-order homogeneous differential operator $\mathcal B$
with $Qg=\bar g+\mathcal B w$. Its symbol is $|\xi|^4 Q(\xi/|\xi|)$,
a polynomial. Thus $\mathcal B$ of any symmetric traceless potential
has the required wave constraint.

Each component of $D^4w$ is a composition of four periodic Riesz
transforms applied to $g-\bar g$, so $\|D^4w\|_8\le C$ since $|g|=1$.
Indeed the symbol $-i\xi_j/|\xi|$, set to zero at the origin, satisfies
the three-torus multiplier criterion \cite[Eq.~(2.10)]{p13:riesz}, since
its first and second derivatives away from zero are $O(|\xi|^{-1})$
and $O(|\xi|^{-2})$, giving the required discrete differences by
integration on unit segments. This also gives the $L^4$ bound used below.

Let $\mathcal M$ denote the periodic centered cube maximal operator.
Declare $O$ to be the open set where either
$\mathcal M|Qg-g|>t^{28}$ or $\mathcal M|D^4w|^8>t^{-8}$.
The weak $(1,1)$ inequality gives $|O|\le Ct^8$, using
$\|Qg-g\|_1\le t^{36}$. Indeed disjoint cubes selected from a level
cover cover the level set upon fivefold enlargement, which increases their
total measure by at most $5^3$. Both tests bound the corresponding
averages on every cube centered in $O^c$.
On $O$ let $d$ be a smooth regularized distance to $O^c$, comparable
to the true distance $d_0$ in the supremum norm, with
$|D^i d|\le C_i d_0^{1-i}$, $i\le4$. To construct it, take maximal
dyadic cubes whose concentric sixfold enlargements are inside the
periodic lift of $O$. On their $3/2$ enlargements their sides are
comparable to $d_0$ by failure of the parent, and packing gives bounded overlap. Average the side lengths using a
smooth partition obtained by normalizing rescaled nonnegative bump
functions equal to one on each original cube and supported on its
enlargement. The quotient rule gives the derivative bounds. For small
$t$ the cubes are smaller than the period, so the construction is
periodic.

Fix a smooth nonnegative centrally symmetric probability density $\chi$
supported in $[-1,1]^3$, and set
$\phi=(t^{-3}\chi(\,\cdot\,/t)-t^2\chi)/(1-t^2)$.
Its moments of orders one through three vanish, by symmetry and cancellation
of the second moment. Its negative mass is $O(t^2)$ and its supremum
norm is $O(t^{-3})$. With $\epsilon=t^8$, set $r=\epsilon d$
(scaling $d$ down by a fixed constant so $r<d_0/2$). Define
\[
 \widetilde w(x)=\int\phi(y)w(x-r(x)y)\,dy\quad(x\in O),
 \qquad \widetilde w=w\quad(x\notin O).
\]
The cube of radius $r$ at $x$ lies in one of radius $2d_0$ centered
at a nearest good point. Consequently, on the smaller cube the means
of $|Qg-g|$ and $|D^4w|^8$ are at most $C\epsilon^{-3}t^{28}$
and $C\epsilon^{-3}t^{-8}$, respectively. If the radius is held
fixed when applying $\mathcal B$, the resulting convolution of $Qg$
is within $Ct^{-3}\epsilon^{-3}t^{28}$ of the same convolution of $g$.
The latter lies within $Ct^2$ of $K$ by convexity and the negative-mass
bound.

Fix $x\in O$, write $C_r=x+[-r,r]^3$ for $r=r(x)$, and put
$A_r=(|C_r|^{-1}\int_{C_r}|D^4w|^8)^{1/8}\le C\epsilon^{-3/8}t^{-1}$.
Subtract a cubic polynomial from $w$ so that the remainder $R$
has mean-zero derivatives of orders zero through three on $C_r$,
choosing coefficients in descending order. Write
$\|F\|_{8,r}=(|C_r|^{-1}\int_{C_r}|F|^8)^{1/8}$. Iterated Poincaré and the
scaled cube embedding $W^{1,8}\hookrightarrow C^0$ give
\begin{align*}
 \|D^jR\|_{8,r}&\le C r^{4-j}A_r &&(0\le j\le4),\\
 \|D^jR\|_{L^\infty(C_r)}
 &\le C\big(\|D^jR\|_{8,r}+r\|D^{j+1}R\|_{8,r}\big)
 \le C r^{4-j}A_r &&(0\le j\le3).
\end{align*}
The polynomial is held fixed when differentiating at $x$. Its variable-radius
average equals the polynomial at every point by the three moment cancellations.

For $|y|_\infty\le1$ set $F_y(x)=x-r(x)y$. Then
$DF_y=I-y\otimes Dr$ and $D^iF_y=-y\otimes D^ir$ for $2\le i\le4$.
The distance estimates give $|D^ir|\le C_i\epsilon r^{1-i}$,
$1\le i\le4$. For $1\le k\le4$ the chain rule, evaluated on $k$ unit
vectors, has one term per partition $\pi$ of $\{1,\ldots,k\}$,
\[
 D^k(R\circ F_y)=\sum_{\pi}D^{|\pi|}R(F_y)
           \big[\,D^{|B|}F_y:B\in\pi\,\big],
\]
where each block uses the vectors indexed by $B$. The singleton partition
has coefficient $(DF_y)^{\otimes k}$, within $C\epsilon$ of $I^{\otimes k}$.
For any other partition, writing $j=|\pi|$ and $m$ for the number of
blocks of size at least two gives
\[
 \prod_{B\in\pi}\|D^{|B|}F_y\|
 \le C\epsilon^m r^{\sum_{B}(1-|B|)}
 = C\epsilon^m r^{j-k},\qquad m\ge1 .
\]
At $k=4$ remove the term with four identity factors, which gives the
frozen-radius convolution of $D^4w$. Using $\|\phi\|_\infty\le Ct^{-3}$
and changing variables $v=x-r y$ in each remaining term yields
\[
 \left|D^4\widetilde w(x)-\int\phi(y)D^4w(x-r y)\,dy\right|
 \le Ct^{-3}\epsilon\sum_{j=1}^4 r^{j-4}\|D^jR\|_{8,r}
 \le Ct^{-3}\epsilon^{1-3/8}t^{-1}=Ct .
\]
This calculation extends from smooth $w$ by local $W^{4,8}$ approximation.
On every compact subregion of $O$, $r$ is bounded below, so the same
estimates give convergence of the averaged derivatives there.
The frozen term is bounded by $Ct^{-3}A_r$, hence
$|D^4\widetilde w|\le C_t$ in $O$ for fixed $t$. Lower derivatives
of $b=\widetilde w-w$ obey, for $0\le k\le3$,
\[
 |D^k b(x)|\le C_t\sum_{j=0}^k r^{j-k}\|D^jR\|_{L^\infty(C_r)}
       \le C_t d_0(x)^{4-k},
\]
by the same chain rule and the supremum bound above.

To glue across $\partial O$, choose smooth $\eta$ equal to zero on
$(-\infty,1]$ and one on $[2,\infty)$. Extend
$b_h=\eta(d/h)b$ by zero off $O$. It belongs to $W^{4,8}$ on the
torus since its support stays away from $\partial O$. The derivatives
of the cutoff of order $i\ge1$ are bounded by $C_i h^{-i}$ and
supported in $E_h=\{h<d<2h\}\subset O$, where $d_0\asymp h$.
For $1\le i\le k\le4$, the product rule terms there satisfy
\[
 |D^i(\eta(d/h))|\,|D^{k-i}b|\le C_t h^{-i} h^{4-k+i}
 =C_t h^{4-k}.
\]
Since $|E_h|\to0$, these terms tend to zero in $L^8$.
On $O$, $D^4b$ belongs to $L^8$ because $D^4w\in L^8$ and
$D^4\widetilde w$ is bounded. Dominated convergence for the
undifferentiated cutoff now gives $D^k b_h\to\mathbf1_O D^k b$
in $L^8$ for every $0\le k\le4$, with the right side extended by
zero. These are precisely the weak derivatives of the extension
of $b$, including at the boundary. Hence $z=\bar g+\mathcal B\widetilde w$
is globally admissible and equals $Qg$ almost everywhere on $O^c$.
There $|Qg-g|\le t^{28}$ by differentiation of integrals. In $O$
the frozen and variable-radius bounds prove \eqref{p13:truncate}.

Choose $t\to0$ along the approaching cells. Since $K$ contains a
fixed centered ball, $K+Ct\mathbb B\subset(1+C't)K$.
Choose rational $u\to0$ with $t^{28}<u/3$ and
$(1-u)^{4m}(1+C't)<1$, leaving strict margin in the range test.
Mollify $z$ with a positive approximate identity, keeping its range in
$K+Ct\mathbb B$ and its waves admissible. For each fixed cell choose
the smoothing scale so that $|z-g|<2u/3$ for the smoothed field outside
a set of measure tending to zero along the sequence. Approximate
$\widetilde w$ after smoothing, as well as its independent constant mean
$\bar g$, to obtain a piecewise-polynomial field satisfying the same
strict range test and $|z-g|<u$ except on a set of vanishing measure.

For precision, smooth periodic potentials may first be approximated in
$C^4$ by finite Fourier sums. Each such sum is a sum of tensor products
of smooth univariate periodic functions. Approximate the fourth derivative
of each factor uniformly by a periodic rational piecewise-linear function,
subtract its rational mean and integrate periodically four times,
subtracting the mean after each integration. Restore a rational
approximation of the original mean of the factor after the last step.
The mean-zero primitive is bounded in the uniform norm, giving convergence
through four derivatives by rational splines. Project the potentials onto
the rational symmetric traceless subspace. The operator $\mathcal B$ has
rational coefficients, and a rational traceless mean approximates $\bar g$.
For the shear trials Fourier truncation approximates each minimizing field
in $L^2$. Strain fluctuations at a nonzero frequency are symmetric gradients.
Stress fluctuations are double curls of symmetric potentials. The double-curl symbol
is, up to sign, $V\mapsto [\xi]_\times V[\xi]_\times^t$, whose image
is the transverse symmetric space since $[\xi]_\times$ is invertible on
$\xi^\perp$. The same spline construction gives admissible trial fields
on a common rational mesh by differentiation and addition of their known
constant means.

It remains to control their energy on the exceptional sets of varying
cells. All the minimizing shear fields are bounded uniformly in $L^p$
for some $p>2$. To see this, if $T(x)$ ranges in $[\lambda I,\Lambda I]$,
$X=E+P_0[(I-2T/(\Lambda+\lambda))X]$ for the minimizing field of mean
$E$. The fluctuation projector on strains has symbol $Z\mapsto UZ+ZU-UZU$, and
on stresses it has symbol $Z\mapsto PZP$. Both have $L^2$ norm one
and finite $L^4$ norm by the Riesz estimate above. Interpolation gives
$\|P_0\|_{p\to p}\to1$ as an upper bound as $p\downarrow2$. The
multiplication factor has norm at most $(\Lambda-\lambda)/(\Lambda+\lambda)<1$.
Choose $p>2$ so their product is less than one. The Neumann series then
gives the asserted uniform bound. Thus the $L^2$ energy on sets of
vanishing measure tends to zero by H\"older. Select closed good subcubes
covering all but a vanishing fraction of $\{|z-g_R|<u\}$. This is possible because inside a
constant grain this set is open and compact subregions admit such a
cover on a fine grid. The extra term $\int\ell |X_i|^2$ in either
energy tends to zero. Approximating each minimizer in $L^2$ as above
changes the corresponding bounded-coefficient energy by arbitrarily
little. This proves density of the inward estimates.

All finite witnesses can be enumerated. Use rational grid coordinates,
algebraic rotations, algebraic load means and the rational potentials
just described. Range membership is the pair of extreme eigenvalue
inequalities proved in Lemma~\ref{p13:steer}, so good-cube and range tests
for piecewise polynomials on closed boxes are first-order real algebraic
tests. The energies are exact algebraic integrals of polynomials on
rational boxes. Keeping running minima or maxima gives the stated
monotone estimates. Equality in $(e:Ce)(e:C^{-1}e)\ge1$ holds exactly when $e$ is an eigenvector of $C$. Thus $M=0$ in either convention exactly when $a_R=a_V$. In that case every grain has $L_Re=\alpha e$.
Every cell already has that hydrostatic eigenvalue. The ordinary shear
trials and isotropization give the same density and bounds with
$\ell=0$, without any marker tests.
\phantomsection\label{proofend:p13:inward}\end{proof}

\begin{theorem}\label{p13:algorithm}
For algebraic $C$ as in Proposition~\ref{p13:computeH} and any $a\in[a_R,a_V]$ specified to arbitrary accuracy, there are terminating two-sided approximations to both sharp shears at $a$. Periodic crystal orientations realizing any named attainable pair to a prescribed tensor tolerance can also be found.
\end{theorem}
\begin{proof}\phantomsection\label{proof:p13:algorithm}
For each sign $\eta$ bracket the minimum $l_\eta(a)$ of $b^\eta$ at the exact bulk. For lower bounds enumerate rational $r\ge0$ in Theorem~\ref{p13:interval}, use the certified lower bounds on $H(r,1)$, and evaluate the subtraction with one-sided rounding. Continuity of the compact-set support in $r$ makes these bounds exhaustive even when the supremum is reached only at infinity.

For an upper bound, seek attainable points below a rational shear-observable ceiling $v$, one on either side of $a$. A strict side test uses a periodic orientation and strain or stress trigonometric trials. Require its hydrostatic trial energy to be strictly below $a^\epsilon$ and its mean deviatoric trial energy to be strictly below $v$. The two arithmetic bounds after isotropization give a point on that side. Alternatively use an exact Hill-endpoint witness of Theorem~\ref{p13:inward} with shear observable below $v$, since the prescribed bulk interval puts that endpoint on the requisite side. Two successful side tests certify $l_\eta(a)<v$ by interpolation. Conversely, if $v>l_\eta(a)$ and $a$ is interior, mix an exact minimizer with each Hill endpoint in arbitrarily small proportions on the path constructed above. The resulting bulks lie strictly to either side because the bulk increment in the laminate formula is strictly monotone in the core share whenever the constituent bulks differ. The shear observables stay below $v$. Density of paired trials verifies both sides. At an endpoint, use the exact-side witness on that side and the same strict-side argument for the other. If $a_R=a_V$, the single support $H_{+\eta}(0,1)$ gives the answer directly. One can decide equality of the two Hill bounds from the specified crystal. Dovetail all tests. This searches over the weights themselves and needs no unverified upper cutoff.

For realization at $B=aJ+bD$, enumerate dyadic orientation fields from rational quaternion lists and six pairs of trigonometric strain and stress trials with respective means $F_i,BF_i$, for an orthonormal basis of $V$. Retain cells with small total residual $r=\sum_i\langle(Y_i-C_RX_i):C_R^{-1}(Y_i-C_RX_i)\rangle$. For fixed cell with tensor $A$, the infimum over its fields equals $\operatorname{Tr}((A-B)A^{-1}(A-B))$. This follows by expanding the square, using $\langle X_i:Y_i\rangle=F_i:BF_i$, and minimizing the two energies separately. Consequently $\|A-B\|_{\rm HS}^2\le\Lambda r$. Conversely if $B\in G$, this residual can be arbitrarily small, by periodic approximation followed by Fourier truncation of its fields and approximation of its orientations in measure with the bounded trials fixed. The target means may be computable reals, since strict residual comparisons are semidecidable by enclosing the displayed integral. The resulting sequence realizes $B$ including at a Hill endpoint.
\phantomsection\label{proofend:p13:algorithm}\end{proof}

\subsection*{An obstruction for data given by approximations}

\begin{theorem}\label{p13o:impossible}
At the fixed physical bulk \(1/3\), neither endpoint of the physical
shear interval is uniformly computable to absolute error \(10^{-6}\)
from an oracle returning the crystal entries with arbitrarily prescribed
errors. This holds even with the promise that the crystal is rational
and that the bulk is attainable.
\end{theorem}
\begin{proof}\phantomsection\label{proof:p13o:impossible}
Here $1/3$ is the physical bulk $\kappa$, so $a=1$. The constants $L,H$ below bound the physical shear $\mu=b/2$. Keep the notation \(e=I/\sqrt3\), \(D\) for deviatoric projection.
Let \(t=1/50\) and \(U=3nn^t-I\), where \(n\) is the crystal axis.
For \(0\le s\le s_0=10^{-4}\) use the crystal
\[
 C_s E=E+s[(\operatorname{tr}E)U+(U:E)I]+tD[U(DE)+(DE)U].
\]
With \(n\) the third axis its Kelvin entries are rational for rational
\(s\). On \(\operatorname{span}(e,U/\sqrt6)\) its matrix has diagonal
\((1,1+2t)\) and off-diagonal \(3\sqrt2s\). Its remaining eigenvalues
are \(1+t\) on the mixed axial shears and \(1-2t\) on the in-plane
traceless matrices. Consequently \(lI_V\preceq C_s\preceq uI_V\),
where \(l=24/25\), \(u=1041/1000\), since the row sums of the off-diagonal
part of the two-dimensional block are below \(1/1000\).
In eigenvalue coordinates the Hill bulks are
\(a_V=1\), \(a_R=1-18s^2/(1+2t)\).

Fix \(s>0\) and a sequence of orientation fields whose effective
tensors tend to an isotropic tensor on one of its Hill bulks.
Write \(U(x)\) for the rotated marker and \(C(x)\) for its stiffness.
At the Voigt bulk the constant-strain trial \(e\) has energy tending
to the minimum. Its minimizing corrector thus tends to zero in
\(L^2\) by coercivity and stationarity. The optimal stress is within
\(o(1)\) in \(L^2\) of \(C(x)e=e+\sqrt3s U(x)\), and its mean tends
to \(e\). Thus \(\langle U\rangle\to0\), and at nonzero Fourier frequencies
the summed squared projections of \(\widehat U(k)\) onto
\(J_n=\{a\odot n\}\), \(n=k/|k|\), tend to zero by equilibrium.
At Reuss the constant-\emph{stress} trial \(e\) has vanishing
complementary energy deficit. The optimal strain is now within
\(o(1)\) of
\[
 C(x)^{-1}e=a_R^{-1}e-\frac{\sqrt3s}{(1+2t)a_R}U(x),
\]
and its mean tends to \(a_R^{-1}e\). Again \(\langle U\rangle\to0\).
The vanishing projections are now onto \(J_n^\perp\), by compatibility.
In each case \(s\) is held fixed while taking these cell limits.

It follows that, with total squared error \(o(1)\), the nonzero modes
of \(U\) in the Voigt case belong to the two-dimensional space of
traceless tensors annihilating \(n\), and in the Reuss case to the
two-dimensional space \(\{a\odot n:a\perp n\}\).
Indeed on the traceless matrices the relevant projection bounds the
distance to its kernel up to a constant, uniform by rotation from
one fixed normal. The squared norms of the retained modes sum to
\(6+o(1)\), since \(|U|^2=6\) and the mean tends to zero.

For a unit deviator \(E\), let \(P\) project onto mean-zero
compatible strains on the torus. The deficit of the constant strain trial
is \(\sup_{z=Pz}(2\langle z:CE\rangle-\langle z:Cz\rangle)\).
It is therefore bounded below by \(\|P(CE)\|_2^2/u\) and above by
\(\|P(CE)\|_2^2/l\). Sum these comparisons over an orthonormal
deviatoric basis and divide by five. The average uniform energy is one because
the deviatoric block of every grain has trace
\(1+2t+2(1+t)+2(1-2t)=5\). The limit of the minimized average is
the shear eigenvalue \(2\mu_*\).

For clarity the following table computes the summed squared projection
on a unit marker mode \(W\). Take the third axis as normal. With
\(S_{ij}=(e_i e_j^t+e_j e_i^t)/\sqrt2\) the five load columns
are \(D_1=\mathrm{diag}(1,-1,0)/\sqrt2\),
\(D_2=\mathrm{diag}(1,1,-2)/\sqrt6\), \(S_{12},S_{13},S_{23}\).
In each entry form \(H=s(W:E)I+tD(WE+EW)\) and sum
\(H_{33}^2+2H_{13}^2+2H_{23}^2\).
\[
\begin{array}{c|ccccc}
 W&D_1&D_2&S_{12}&S_{13}&S_{23}\\ \hline
 S_{13}&t^2/2&t^2/6&t^2/2&(3s+t)^2/9&0\\
 D_1&(3s-2t)^2/9&0&0&t^2/2&t^2/2
\end{array}
\]
The sum is invariant under rotations about the normal, transitive on the
unit circle in each of the two real spaces. For complex Fourier
coefficients add the contributions of the real and imaginary parts.
By Parseval and the preceding concentration, the limits of the
five-load squared projection averages at Reuss and Voigt are respectively
\[
 y_1(s)=\tfrac65[(3s+t)^2/9+7t^2/6],\qquad
 y_2(s)=\tfrac65[(3s-2t)^2/9+t^2].
\]
The discarded mode error is \(o(1)\) by boundedness of the displayed
linear maps. Since \(y_1\) increases and \(y_2\) decreases on this interval,
every Voigt-side shear is at most \(L=(1-y_2(s_0)/u)/2\), and every
Reuss-side shear is at least \(H=(1-y_1(s_0)/l)/2\). Their separation is
\[
 H-L=\frac{1865119}{166560000000}>\frac1{90000}.
\]

The table, the crystal block and this rational difference are checked
in \path{certificates/baseline/p13obstruction/verify.py}.

Both sides are nonempty by Theorem~\ref{p13:interval}. On replacing
\(C_s\) by \(C_0\) on any texture, the cell stiffness changes in operator
norm by \(O(s)\) uniformly, since the local tensors differ by \(O(s)\)
in relative order and the strain principle preserves that order.
As \(s\downarrow0\), subsequences of the two sets of realizing tensors
therefore give isotropic tensors in the closure of \(C_0\) at its
only bulk \(1/3\), one of shear at most \(L\), the other at least \(H\).
At the same bulk, \(C_s\) has upper shear at most \(L\) for \(s>0\), and
the rescaled crystal \(C_s/a_R\) has lower shear at least
\(H/a_R\ge H\). Both families tend to the same rational crystal.

For a machine that first halts at step \(N\ge0\), put
\(s=2^{-15-N}\), and for a machine that never halts put \(s=0\).
These crystals all have rational entries. There is a uniform oracle
for either family. To approximate \(s\) to error \(2^{-j}\), simulate
\(j\) steps and return zero if the machine has not halted. The error
is at most \(2^{-15-j}\). The displayed rational matrix formulas on
this interval have bounded derivatives and \(a_R\) bounded away from
zero, giving entrywise names to any precision by rational interval
arithmetic. For the first family the upper endpoint, and for the
rescaled family the lower endpoint, distinguish the two outcomes
at error \(10^{-6}\) by comparison with \((L+H)/2\).
This would decide halting.
\phantomsection\label{proofend:p13o:impossible}\end{proof}

\subsection*{An exact bulk slice for zinc}
\phantomsection\label{item:p13-zinc}
Use the following nominal hexagonal single-crystal stiffnesses in GPa.
\[
(C_{11},C_{12},C_{13},C_{33},C_{44})=(165,31.1,50,61.8,39.6).
\]
At effective bulk modulus $66$ GPa, the finite certificates give
\[
\begin{aligned}
 36.80014385&\ \le\mu_{\min}\le 37.62336469 \quad{\rm GPa},\\
 43.35521987&\ \le\mu_{\max}\le 44.07248692 \quad{\rm GPa}.
\end{aligned}
\]
Every shear between $37.62336469$ and $43.35521987$ GPa is attained
in the closure at this bulk. To obtain the inner bounds, construct
anisotropic textures with certified elastic tensors, isotropize them
and bracket the target bulk with a pair of trees. Mix the pair at the
core fraction fixed by the bulk equation. The outer ends follow from
joint strain and stress translations valid for arbitrary orientations.
The two endpoint brackets enclose different extrema. They assume the
stated stiffnesses as exact inputs. Supplement~\ref{sup:crystal} gives
the arithmetic and the finite witnesses.

\section{14. Joint bounds on bulk and shear}\label{problem-14}
\Milton{14}
The singular isotropic rectangles are given after the full-tensor
criterion in Problem~8. For two positive elastic phases we first state
an implicit variational characterization. A change of the two phase bulk
moduli can reverse which dimension admits a shear common to the entire
bulk range. A spatial cubic obstruction also separates the physical
upper boundary from every finite quadratic translation at particular
interior bulks. The explicit full boundary for general phase data
remains open.

\subsection{Sharp variational bounds for positive phases}
\begin{corollary}\label{p14:support}
Fix two positive-definite isotropic elastic phases $C_1,C_2$ in dimension $d=2$ or $3$, and a phase-one fraction $0<f<1$. Put $a=d\kappa_*$, $b=2\mu_*$, and $m=d(d+1)/2-1$. The following four tests specify the joint isotropic region in the periodic closure. Take a unit hydrostatic mean and an orthonormal basis of $m$ deviatoric means, with compatible or equilibrated periodic fields according to signs $\epsilon,\eta$ as in \eqref{p13:support}. For weights $s,t\ge0$ put
\[
 W_i(x)=s X_0:C_i^\epsilon X_0+\frac tm\sum_{j=1}^m X_j:C_i^\eta X_j.
\]
Let $H^f_{\epsilon\eta}(s,t)$ be the infimum, over the fields, of
\[
 \int W_2+\int_0^f (W_1-W_2)^\uparrow(u)\,du,
\]
where $g^\uparrow$ is the nondecreasing quantile of the values of $g$ in the unit-volume cell. Then $(a,b)$ occurs exactly when $s a^\epsilon+t b^\eta\ge H^f_{\epsilon\eta}(s,t)$ for all signs and nonnegative weights. More explicitly, order the phases by shear in the definitions of $K_0,K_1$ in Lemma~\ref{p14:mix}, keeping their assigned fractions. The admitted bulks are $dK_0\le a\le dK_1$, and the full shear interval is
\[
 l_+(a)\le 2\mu_*\le 1/l_-(a),\qquad
 l_\eta(a)=\sup_{\epsilon=\pm1,\;r\ge0}
      [H^f_{\epsilon\eta}(r,1)-r a^\epsilon].
\]
\end{corollary}
\begin{proof}\phantomsection\label{proof:p14:support}
For any fields, assign phase one where $W_1-W_2$ is smallest, splitting a level set to fill the exact fraction. Moving mass between two ordered levels proves the quantile formula. Truncating the fields in Fourier space gives a decreasing sequence of infima on finite mode spaces with the same limit, since the change of the objective is bounded by the sum of the two $L^1$ energy changes. Minimizing first in the fields expresses $H^f_{\epsilon\eta}$ as the infimum of the corresponding weighted hydrostatic and deviatoric cell energies. Uniform ellipticity gives compactness at the prescribed fraction. Lemma~\ref{p13:locality} gives closure under further homogenization at that exact fraction. Rotating a cell rotates its tensor because both underlying phases are isotropic.

Let $\mathcal R_f$ denote the isotropic region. Rotation averaging preserves the hydrostatic contraction and the deviatoric trace in either energy convention. Iterating the rotated-layer construction of Problem~13 gives an isotropic member with both chosen energies at most those of the original tensor. In two dimensions rotation through $\theta$ acts by $2\theta$ on the deviatoric plane, so the rotation average is again scalar on it. Consequently $H^f_{\epsilon\eta}$ is the minimum of $s a^\epsilon+t b^\eta$ over $\mathcal R_f$. For two points of $\mathcal R_f$, first laminate them and then isotropize. Constant strain and stress trials bound both selected energies by their prescribed arithmetic means. To join the points continuously, add successive layers of one isotropic material around the other using the spherical mean of the operators $\Gamma_{B,n}$ for host $B$ (the circular mean in the plane), as in the laminate formula in the proof of Theorem~\ref{p13:interval}. The eigenvalues stay between the two input values by the uniform strain and stress bounds and are continuous in the total core fraction. The two hypotheses of Lemma~\ref{p13:four} therefore hold for $\mathcal R_f$. Its conclusion gives the four inequalities and the formula for both shear endpoints.
\phantomsection\label{proofend:p14:support}\end{proof}

\subsection*{Complete bulk sections for positive phases}

\begin{lemma}[Hydrostatic field]\label{p14:hydro}
In dimension $d$ fix $m>0$, let $b=2(d-1)m/d$ and let $k(x)>0$ be a
two-valued periodic bulk modulus. If the shear is everywhere $m$, then
$I/d$ is a hydrostatic eigenstrain with mean stress $KI$, where
\[
 K=\langle(k+b)^{-1}\rangle^{-1}-b .
\]
Write $\chi$ for the core indicator, $f=\langle\chi\rangle$ and
$Q=(D^2\Delta^{-1}-I/d)(\chi-f)$, with the inverse Laplacian taken
on mean-zero periodic scalars. If $\chi Q=0$, changing the core shear to
any positive value preserves this hydrostatic solution. If
$\|\chi Q\|_2\to0$ at fixed phases and fraction, the hydrostatic
response tends to the same value.
\end{lemma}
\begin{proof}\phantomsection\label{proof:p14:hydro}
Let $t=(K+b)/(k+b)$, of mean one, and set
$e=I/d+D^2\Delta^{-1}(t-1)$. Its trace is $t$. Its stress is
$2me+(k-2m/d)tI$ and its divergence is $\nabla[(k+b)t]=0$.
Its mean is $\langle kt\rangle I=KI$. This gives Hill's common-shear
identity \cite[\S5.3, (5.13)-(5.19)]{book}.
Writing $t_c,t_m$ for the two phase traces, the strain deviator is
$(t_c-t_m)Q$. A shear change in the core adds stress proportional to
$\chi Q$. When this tends to zero in $L^2$, test the difference of the
two weak equilibrium equations against the displacement difference.
Uniform ellipticity bounds its strain norm by a constant times this
residual norm, and hence the mean stresses have the same limit.
\phantomsection\label{proofend:p14:hydro}\end{proof}

\begin{lemma}[Bulk endpoints and interpolation]\label{p14:mix}
In dimension $d=2$ or $3$ let positive isotropic phases $(k_i,m_i)$ occupy fractions $f_i>0$ summing to one, $i=0,1$, with $m_0\le m_1$. Define
\[
 K_j=\left[\sum_{i=0}^1\frac{f_i}{k_i+2(d-1)m_j/d}\right]^{-1}
       -\frac{2(d-1)m_j}{d}.
\]
The isotropic bulk projection is exactly $[K_0,K_1]$. Between two isotropic
responses of the same fraction and the same shear, every intermediate bulk
is attained in the closure at that shear. The analogous interpolation holds
at common bulk.
\end{lemma}
\begin{proof}\phantomsection\label{proof:p14:mix}
Replace both shears by their smaller, then their larger value.
The cell energy is monotone in stiffness and Lemma~\ref{p14:hydro}
gives the two bounding dilation energies. These are the Hashin-Shtrikman bulk bounds in three dimensions and their planar analogues \cite[\S23.5, Eqs.~(23.42), (23.45), pp.~470 to 473]{book}. The argument applies in either ordering of
$k_0,k_1$. 

With matrix $(k,m)$, core $(l,n)$, core fraction $p$, a sequential matrix laminate
with distribution $\nu$ of normals has stiffness
\begin{equation}\label{p14:matlam}
\begin{aligned}
 C&=C_m+pL\big[I+(1-p)\bar\Gamma L\big]^{-1},\qquad L=C_c-C_m,\quad \bar\Gamma=\int\Gamma_N\,d\nu(N),\\
 \Gamma_N T&=N\odot\left[\frac{I-NN^T}{m}
                       +\frac{NN^T}{k+2(d-1)m/d}\right]TN .
\end{aligned}
\end{equation}
In this elastic coating formula, successive matrix additions give arbitrary
weights to the normal directions. \cite[p.~167, (9.22) and its converse, pp.~174 to 175, (9.55), (9.59)]{book} Uniformly distributed normals can be approximated by
finite sums. 
At uniform normal distribution the result is isotropic. For $e=I/\sqrt d$ the displayed operator satisfies
$e:\Gamma_N e=1/[d(k+2(d-1)m/d)]$. Thus the bulk modulus is
\[
 k+\frac{p(l-k)}{1+(1-p)(l-k)/(k+2(d-1)m/d)}.
\]
Taking phase $j$ as matrix gives $K_j$. Formula \eqref{p14:matlam} includes a common bulk or shear. Indeed $C_m^{1/2}\Gamma_N C_m^{1/2}$ is the orthogonal projector onto the normalized jump space, so $0\preceq\bar\Gamma\preceq C_m^{-1}$. Choose $0<\delta\le1$ with $C_c\succeq\delta C_m$. The symmetric matrix $\bar\Gamma^{1/2}L\bar\Gamma^{1/2}$ is bounded below by $-(1-\delta)I$ and has the same nonzero eigenvalues as $\bar\Gamma L$. Thus the displayed inverse exists even when $L$ is singular. For invertible $L$ the expression is the usual $(L^{-1}+(1-p)\bar\Gamma)^{-1}$ correction multiplied by $p$. Continuity gives the general formula. Between these two endpoints the continuous isotropic path in the proof of Theorem~\ref{p13:interval} preserves the original fractions by Lemma~\ref{p13:locality}.

At a common shear all affine deviatoric strains have uniform
stress, and hence that shear is retained under further mixing.
Lemma~\ref{p14:hydro} traverses the intervening bulks continuously as
the mixing fraction varies. At common bulk affine dilation has uniform
stress. A matrix laminate with isotropically distributed normals is
isotropic and depends continuously on its core fraction from zero to one
by the laminate formula \cite[Eqs.~(9.55), (9.59), pp.~174 to 175]{book}. These operations
preserve the original fraction. Lemma~\ref{p13:locality} allows substitution on finer scales at the specified fractions.
\phantomsection\label{proofend:p14:mix}\end{proof}

\begin{proposition}\label{p14:sections}
At half of each of the two indicated positive phases, the entire attainable
bulk interval at each shear in the middle column is the last column, in the
isotropic closure.
\[
\begin{array}{c|c|c|c}
d &(k_0,m_0),\ (k_1,m_1)&\mu_*&\kappa_*\\\hline
2 &(7,1),\ (3,2)&[57/41,46/33]&[13/3,31/7]\\
3 &(3,1),\ (1/2,3/2)&[71/58,333/272]&[46/37,4/3]\\
3 &(8,1),\ (2,2)&[1.421570,1.423160]&[68/19,88/23]
\end{array}
\]
\end{proposition}

\begin{proof}\phantomsection\label{proof:p14:sections}
By Lemma~\ref{p14:mix} it suffices to attain the two extreme bulks at
the indicated shears. Lemma~\ref{p13:locality} allows the grains in each layer and chamber below to have their respective fractions. We first give two particularly short constructions.
Write $(k,m)$ for the matrix moduli, $(l,n)$ for the core moduli and $p$ for its fraction. Apply \eqref{p14:matlam}.
For equal weights on coordinate axes the hydrostatic
eigenvalue gives the matrix bulk $K_j$ by substitution. On the trace-free
diagonal and off-diagonal spaces the two shear moduli are, respectively,
$g_1=g(\gamma_1)$, $g_2=g(\gamma_2)$, where
\[
 g(z)=m+\frac{p(n-m)}{1+(1-p)(n-m)z},\quad
 \gamma_1=\frac{2}{d[k+2(d-1)m/d]},\quad \gamma_2=\frac2{dm}.
\]
Indeed averaging $\Gamma_N$ on the axes multiplies these two spaces by
$\gamma_1/2,\gamma_2/2$. For uniform normals, averaging over rotations
before taking the inverse gives shear $M_j=g(w\gamma_1+(1-w)\gamma_2)$
with $w=2/(d+2)$, the relative dimension of the trace-free diagonal space.
The bulk is again $K_j$.

We record how to isotropize the axis grain or any grain having $I$ as
eigenstrain. Its rotated copies have the same stress under uniform dilation.
Arrange signed-permuted copies in the chambers $0<x_d<\cdots<x_1<1/2$
and their images to get a cubic tensor. For a cubic tensor with
shears $a,z$, use copies in orientations $gR$, one per chamber $g$,
with $R$ chosen so $x(R)=y(R)=w$ by Lemma~\ref{p9:cubic}.
On every trace-free unit load the uniform strain and stress tests put
both new shears between $H=(w/a+(1-w)/z)^{-1}$ and $A=wa+(1-w)z$.
This interval has width at most $\max(w,1-w)|a-z|$ since $H\ge\min(a,z)$.
Iteration therefore converges to an isotropic tensor at unchanged bulk.
For any original grain $C$, with $q=d(d+1)/2-1$ and $P$ the deviatoric
projection, the resulting shear obeys
\begin{equation}\label{p14:VR}
 \frac{q}{2\operatorname{tr}(PC^{-1})}\ \le\ M \le\ \frac{\operatorname{tr}(PC)}{2q}.
\end{equation}
Indeed the traces are unchanged under rotation and the Voigt and Reuss
bounds, summed on an orthonormal deviatoric basis, hold through every
hierarchy. In particular for the axis grain these bounds are the
harmonic and arithmetic means of $g_1,g_2$ with weights $w,1-w$.

For row one the arithmetic mean with matrix 1, the isotropic matrix
shears $M_0,M_1$, and the harmonic mean with matrix 0 are respectively
$25/18,57/41,46/33,200/143$, in increasing order. In row two the corresponding
numbers are $1371/1120,71/58,333/272,1139/930$, also increasing.
Thus the interval $[M_0,M_1]$ is attained on each bulk line.
All entries follow by substitution in $g$, at $p=1/2$.

\smallskip\noindent\emph{A finite construction for row three.}
Write $(k,m)$ for the selected matrix and $(l,n)$ for the core, and $x$
for a core fraction. Choose a mesh $0=b_0<b_1<\cdots$ and at each step write
$x=b_j$, $b=b_{j-1}$ and $w=(x-b)/(1-b)$. The material at fraction
$b$ will be used inside the one at $x$. We describe its layers starting
from the coarsest, and compute their stiffnesses in the reverse order.

A convenient way to keep the hydrostatic field of Lemma~\ref{p14:hydro}
is to track a matrix $D$ down the tree. Start at $D=0$, fraction $x$.
In a split normal to $N$ from fraction $y$ to child fraction $y'$ add
$(y'-y)NN^T$ to $D$. At the leaves the common-shear strain under mean
$I/d$ is exactly $I/d+(t_c-t_m)D$, with $t_c,t_m$ from Lemma~\ref{p14:hydro}
at fraction $x$. Indeed the strain increments are compatible layer jumps,
the child mean equals the parent and $\operatorname{tr}D=y-x$.
The pressure $kt$ has conditional mean $y l t_c+(1-y) k t_m$ and changes
by $(y'-y)(l t_c-k t_m)$. Hence the mean traction jump is zero by
\[
 l t_c-k t_m+2m(1-1/d)(t_c-t_m)=0 .
\]
This verifies the layer equations at every scale. A terminal core with
$D=(1-x)I/d$ thus has zero deviator, and changing its shear keeps the
bulk exactly at the matrix endpoint.

For cardinal layers put $R_i=x+dD_{ii}$, initially all $x$.
Their average is the current fraction $y$. A split normal to $e_i$ changes
only $R_i$, by $d(y'-y)$. Split the coordinates in order, each from $x$
to $1$ with weight $w$ (stop) or to $b$ with weight $1-w$ (continue).
If all choices are $b$ the resulting fraction is $b$ and $D$ is scalar.
Insert the already constructed isotropic material at $b$ there, with freely
rotated axes, since the inherited common-shear strain is hydrostatic.
The base case $b=0$ is pure matrix. In two dimensions a stopped state has
coordinates $1,r$, $r=b$ or $x$. Fill it by a simple laminate of
core fraction $h=(1+r)/2$, normal to the second of these coordinates
(the one with value $r$). Its final core has both coordinates $1$.
Isotropize the complete step as above. This gives a planar texture at each mesh
fraction with hydrostatic eigenstrain and zero common-shear deviator on its core.
These field properties persist under rotation aggregation because all children
have the same fraction. Thus the same constant macroscopic dilation strain
and its same hydrostatic stress are compatible and equilibrated throughout.
Lemma~\ref{p13:locality} gives substitution at the required fractions for these trial fields
at both values of the core shear.

The same planar texture at fraction $h$ has scalar conductivity
\begin{equation}\label{p14:scalar}
 G=m+\frac{h(n-m)}{1+(1-h)(n-m)/(2m)}
\end{equation}
when the two scalar conductivities are $m,n$. In fact for macroscopic
electric field $v$ give a node field $v-Da/m$ and flux $m v-Da+y a$.
Across a layer the field jump is normal and the flux jump has zero normal
component. At a matrix leaf the conductivity law holds. At a core leaf
$D=(1-h)I/2$, so it holds on choosing
$a=(n-m)v/[1+(1-h)(n-m)/(2m)]$. The root flux is $mv+ha$.
This works for every $v$ and is preserved by the rotations at each isotropic
state, which has scalar $D$ on entering and the same mean flux in every
orientation.

In the spatial construction a stopped state has one coordinate $1$ (the
polar axis) and two transverse coordinates, each $b$ or $x$. If the two
are unequal laminate with pure matrix across the lower transverse coordinate.
Use a core-bearing child of fraction $h=(1+2r)/3$ where $r$ is the greater
of the two transverse coordinates, with weight $y/h$.
The increase $3(h-y)=x-b$ makes both transverse coordinates $r$ in this
child. For already equal transverse coordinates use that child directly.
Fill it by the extrusion in the polar direction of the above planar texture
at fraction $h$. On its final core the additional matrix $D$ in the
transverse plane is $(1-h)I_2/2$. Each transverse $R_i$ therefore
increases by $3(1-h)/2=1-r$ and becomes $1$, while the polar coordinate stays
$1$. This proves the required hydrostatic property for the whole step.
For a surviving node filled from fraction $b$, the increments raising all coordinates from $b$ to $1$ are isotropic on
each final core, even upon rotating the subtree. Hence every spatial step
has exact bulk
\[
 K(x)=\big[(1-x)/(k+4m/3)+x/(l+4m/3)\big]^{-1}-4m/3 .
\]

We give the stiffness calculation for the extruded child so the numerical
test specifies a finite sequence of layer operations. Index matrix and core
by $i=0,1$, put $a_i=k_i+m_i/3$, $b_i=k_i-2m_i/3$,
$c_i=k_i+4m_i/3$, let a bar mean their arithmetic average at core fraction $h$,
and let $\Delta v=v_1-v_0$. Evaluate the planar texture with phase
data $(a_i,m_i)$ and denote its shear by $P$. For axial strain $e$ and planar
trace $u$ the effective normal energy, apart from the planar deviator,
is $A u^2+2Bue+C e^2$, where
\[
 t=\frac{h(1-h)}{(1-h)a_1+h a_0+m},\qquad
 (A,B,C)=(\bar a,\bar b,\bar c)
       -t\big((\Delta a)^2,\Delta a\Delta b,(\Delta b)^2\big).
\]
To see this, average fluctuations along the extrusion axis by convexity.
The phase normal energy is $a_i u^2+2b_i ue+c_i e^2$.
Since $\Delta a\ne0$, we can write $b_i=\lambda a_i+\nu$ and shift the planar
trace to $u+\lambda e$. The remaining minimization uses the planar bulk,
$A=\bar a-t(\Delta a)^2$, known from the construction. Expanding back
gives $B=\bar b+\lambda(A-\bar a)$ and $C=\bar c+\lambda^2(A-\bar a)$.
The two shears involving the polar axis reduce to scalar conductivity
and equal \eqref{p14:scalar}, since the axial displacement contributes
the energy of a scalar gradient with coefficient $m_i$.
Thus, in transverse coordinates 1 and 2, the normal stiffness matrix has
$N_{11}=N_{22}=A+P$, $N_{12}=A-P$, $N_{13}=N_{23}=B$, $N_{33}=C$.
The three ordinary shear moduli are $P,G,G$.

For clarity, stacking two such orthotropic tensors across a cardinal
axis $i$ takes the arithmetic means of
\[
 1/N_{ii},\qquad N_{ji}/N_{ii},\qquad
 N_{jk}-N_{ji}N_{ik}/N_{ii}\quad(j,k\ne i).
\]
Recover the new normal block by solving those same three relations.
This follows by eliminating the normal strain using the common normal
traction and the common tangential strains. Shears involving $i$ average
harmonically, the others arithmetically, by the same two interface conditions.

Take 256 spatial steps of length $1/512$ and planar steps of length
$1/768$ through fraction $2/3$, so all stopping fractions occur on the
planar mesh. At each complete planar or spatial step use the lower bound
in \eqref{p14:VR} when $n>m$ and the upper when $n<m$. Round the shear
outwards to a multiple of $10^{-12}$ and retain the exact bulk. These
replacements preserve the one-sided inequalities by monotonicity under
lamination and by \eqref{p14:VR}. The only unknown modulus in the normal
block of an extruded stopping child is its planar shear $P$. Increasing
$P$ increases the tensor by a nonnegative operator on the transverse
deviators. Thus the same comparison applies across both dimensions.
The rational check \path{certificates/baseline/p14sections/sections.py} checks the
normal row sums $dK(x)$ at every completed step and certifies the final
one-sided shear bounds in the table.
\[
 \begin{array}{c|c|c}
 (k,m) &K(1/2)&\text{constructed shear}\\\hline
 (8,1)&68/19&\ge1423161638739/10^{12}\\
 (2,2)&88/23&\le1421569435537/10^{12}
 \end{array}
\]
The isotropic-normal matrix laminates at these bulks have shear $24/17$
and $87/61$, respectively. Interpolation on each bulk line, then at
the prescribed shear, proves row three.
\phantomsection\label{proofend:p14:sections}\end{proof}

\subsection*{Common shears in two and three dimensions}
Write
$H_f(x,y;s)=[f/(x+s)+(1-f)/(y+s)]^{-1}-s$.
For fixed phases and fraction define a common shear to be a shear value
attainable at every bulk in the bulk projection. The geometry is allowed
to depend on the chosen bulk. At an individual common shear the entire
horizontal section is therefore known.

\begin{theorem}[Changing bulk moduli reverses dimensional separation]\label{new:reversal}
Take fractions $1/4,3/4$ in the order shown. Keeping the two constituent
shears fixed at $1,2$ gives the following exact common-shear sets.
\[
\begin{array}{c|c|c}
(\kappa_1,\mu_1),\;(\kappa_2,\mu_2)&d=2&d=3\\\hline
(1/3,1),\;(310/967,2)&[21/13,\,6294/3895]&\varnothing\\
(18,1),\;(1/4,2)&\varnothing&[215/128,\,919/547]
\end{array}
\]
More generally, for any $0<f<1$ and any finite shear ratio $m/n>1$,
positive phase bulks can be chosen with a common interval of positive
length in the plane and an empty common-shear set in space. The opposite
choice is possible at every such fraction for $1<m/n\le9/4$.
\end{theorem}

\begin{proof}\phantomsection\label{proof:new:reversal}
Let the phase of fraction $f$ have moduli $(l,n)$, and its complement
$(k,m)$ with $m>n$. The sharp bulk edges in each dimension are
the two shifted harmonic means from Lemma~\ref{p14:mix}.
Suppose a low-bulk material and a high-bulk material with the same
shear have been constructed. Coating at their common shear varies the
bulk continuously between the two extreme values. Fixed-fraction
substitution keeps the prescribed phase volumes. Thus the main task
is to compare the two edge shear intervals.

For the first family let $P=\max(f,1-f)$, $R=3P^2(1-P)$,
$r=3R/(2+R)$, $\delta=1/n-1/m$, $c=(1-r)/2$, and set
\[
        l=n/3,\qquad 1/k=3/n+c\delta .
\]
The two candidate common shears are $H_f(n,m;s_2(l,n))$ and
$H_f(n,m;s_2(k,m))$, where $s_2(K,M)=KM/(K+2M)$. On a planar bulk
edge a staircase of layers has fourth angular statistic
$w_g=(1-g)[-2\log(1-g)-g]/g$ at core fraction $g$. The edge comparisons
show that this construction passes the opposite matrix shear at both
ends. Interpolation at each fixed bulk fills the edge intervals.
Their matrix endpoints also exclude every shear outside the displayed
range. In three dimensions the necessary edge bounds are disjoint.
Supplement~\ref{sup:common} gives the two inequalities, including their
uniform dependence on $f$ and $\delta$.

In space the corresponding shift is
$s_3(K,M)=M(9K+8M)/[6(K+2M)]$. For $1<m/n<9/4$ its two ranges,
$(2n/3,3n/2)$ and $(2m/3,3m/2)$, overlap. Choose the bulks first
to give the same shift, then perturb them slightly. A cubic axial
laminate at each bulk edge brackets the opposite matrix shear with
strict slack. The planar necessary condition
$1/k-1/l\le(1/n-1/m)/2$ for a common shear fails. Limits followed by
finite choices give the endpoint ratio $9/4$ as well. Substitution of
the two sets of rational data into the shifts gives the four endpoints
in the table.
\phantomsection\label{proofend:new:reversal}\end{proof}

\subsection*{A constraint beyond quadratic translations}
For a positive phase pair the variational tests in
Corollary~\ref{p14:support} use the full cell energies. Restricting a
bound to quadratic strain and stress translations loses information
even at a strictly interior bulk.
\begin{theorem}[A cubic obstruction at an interior spatial bulk]\label{new:cubic}
Take isotropic spatial phases $A=(1/10,1)$ and $B=(1,n)$ with
$n=1/40$ or $1/50$. At every $B$ fraction $0<f<1$ put
$K=H_f(1,1/10;17/28)$. The largest shear in the physical closure at
this bulk is strictly below the upper shear admitted by all finite
quadratic strain and stress translations. At $f=1/2$, $K=3/8$.
\end{theorem}
\begin{proof}\phantomsection\label{proof:new:cubic}
At the computed quadratic face, equality in the translation inequalities
forces the phase-field fluctuations toward two explicit linear spaces.
The compatible strains and equilibrated stresses give a wave constraint
for each nonzero Fourier direction. From these fields one recovers
gradient rows for which the mean of every cubic determinant is the
determinant of the three means. A rational linear combination of these
determinants and lower-order identities vanishes on both phase spaces.
At their required mean it equals $f(1-f)$. The corrected field is an
order-zero Fourier multiplier of the phase indicator. Its symbol is
smooth on the sphere, so $\|Y\|_4\le C\|\chi\|_4+C$ uniformly in
the geometry \cite[Thm.~2.1]{p13:riesz}. A cubic term with an erroneous
factor $e$ is bounded in $L^1$ by $\|e\|_2$ times the $L^4$ norms
of the other two factors. The vanishing $L^2$ projection error thus
passes the cubic identity to the limit. The
nonzero mean supplies the contradiction. The certificate verifies the
polynomial identities on the phase spaces and on the wave constraints
for the two stated contrasts.
\phantomsection\label{proofend:new:cubic}\end{proof}

At equal fractions the point $(3/8,93/1171)$ for $n=1/50$ is in the
physical region in both dimensions. The status in space of
$(3/8,83/941)$ for $n=1/40$ is open. A bound at the quadratic face alone
does not give a distance large enough to exclude this interior point.
These examples leave the general finite-phase boundary open.

\subsection{Coupled bounds for two spatial solids}

\begin{theorem}\label{p14t:bound}
In three dimensions take positive isotropic phases
\((\kappa_1,\mu_1)=(2,2)\) and \((\kappa_2,\mu_2)=(1,1)\), with
phase-one fraction \(0<f<1\). Define the increasing shifted harmonic
mean \(H(s)=[f/(2+s)+(1-f)/(1+s)]^{-1}-s\).
Every isotropic tensor in the fixed-fraction closure has
\(\kappa_*=H(X),\ \mu_*=H(Y)\), where \(4/3\le X\le8/3\) and
\begin{equation}\label{p14t:curve}
 \frac{459X^2-1896X-1960}{81X^2+786X-5080}
 \ \le Y\le \frac{6(1043X^2-158X-256)}{1821X^2+3004X+448}.
\end{equation}
\end{theorem}
These are outer bounds. Determining the attainable shear endpoints at a general prescribed $X$ remains a realization problem.

\begin{proof}\phantomsection\label{proof:p14t:bound}
Lemma~\ref{p14:mix} gives the lower and upper bulks $H(4/3)$ and $H(8/3)$ from the common-shear hydrostatic solution.
We will also use strict arithmetic and harmonic gaps for shear.
Here the local stiffness and compliance are scalar multiples
\(a(x)D\) of fixed operators, with \(a\) taking two distinct positive
values. At mean strain \(E\) the stiffness energy is
\(\min_w\langle(E+e(w)):C(x)(E+e(w))\rangle\) over periodic
\(H^1\) displacements. At mean stress \(F\) the complementary energy
is \(\min_\sigma\langle\sigma:C(x)^{-1}\sigma\rangle\) over
equilibrated periodic symmetric stresses with \(\langle\sigma\rangle=F\).
In each principle multiply the fields by \(D^{1/2}\). Write \(Z\) for
the resulting mean and \(\Pi\) for projection onto the space of
normalized mean-zero fluctuations. The trial
\(-\Pi(aZ)/\max a\) improves the constant-field energy by at least
\(\|\Pi(aZ)\|_2^2/\max a\), by expansion of the square. Summing over
an orthonormal deviatoric basis for \(Z\), Parseval gives the sum of
these squared norms as \(\mathrm{var}(a)\) times the deviatoric
trace of the wave projector. This trace is independent of direction
and positive because the projector has rank three and its hydrostatic
diagonal is at most one. Applied to stiffness and compliance,
this gives strict gaps bounded uniformly at fixed fraction for
an isotropic limit. Thus \(Y\) is finite and positive. The inverse
of the shifted harmonic mean for general scalar phase values
\(z_1,z_2\) is
\[
 s=\frac{h(z_1-z_2)^2}{\bar z-z_*}-\widetilde z,
 \quad h=f(1-f),\quad \bar z=fz_1+(1-f)z_2,\quad
 \widetilde z=(1-f)z_1+fz_2 .
\]
Replacing all three moduli by reciprocals replaces \(s\) by \(1/s\),
by direct substitution.

Use six simultaneous periodic fields \(F=(A,M)\), where \(A\) is
symmetric and \(M\) maps the traceless symmetric space linearly into
symmetric matrices. Let \(D_j\) be an orthonormal basis of that
five-dimensional space and \(M_j=M(D_j)\). In the upper test take
\(A\) a strain and the \(M_j\) stresses with phase energy
\(L_i[F]=A:C_iA+\sum_j M_j:C_i^{-1}M_j\). In the lower test
interchange strains and stresses and interchange \(C_i,C_i^{-1}\).
The six minimizations at fixed means are independent. Their sum
uses \(L_*\) formed in the same way from the one effective tensor.

Suppose the translation \(T\) is a constant quadratic form with
\(L_i-T\succeq0\) and its restriction to every nonzero Fourier
wave space positive semidefinite. Parseval gives
\(\langle T[F]\rangle\ge T[\langle F\rangle]\).
For a minimizer let \(F_0\) be its overall mean and \(d\) its
phase-one mean minus its phase-two mean. Set
\(\Delta=L_1-L_2\) and use overbars and tildes for the same
two weighted means as above. Averaging its six constitutive
relations gives \((L_*-\bar L)F_0=h\Delta d\) because the
mean conjugate field is \(L_*F_0\) by stationarity and polarization.
Apply phasewise Jensen to \(L_i-T\), using local means
\(F_0+(1-f)d,F_0-fd\), and add the mean bound for \(T\).
Expanding the two squares and using the averaged relation gives
\(-F_0:\Delta d\ge d:(\widetilde L-T)d\). These inequalities for
\(d=(h\Delta)^{-1}(L_*-\bar L)F_0\) pass to limits in \(L_*\).
At an isotropic limit take \(d=(pI/3,D\mapsto qD/5)\) with arbitrary
\(p,q\). The strict gaps and invariance of this plane under all the
operators permit solving the averaged relation for \(F_0\) there.
Consequently on this plane
\[
 h\Delta(\bar L-L_*)^{-1}\Delta-\widetilde L+T\succeq0 .
\]
Before adding \(T\), its quadratic in \((p,q)\) is
\(Xp^2+q^2/(10Y)\) in the upper test and
\(p^2/(9X)+2Yq^2/5\) in the lower, by the reciprocal identity
and the eigenvalues \(3\kappa,2\mu\) of spatial stiffness.

Each translation subtracts from a reference energy a quadratic penalty
on one deviatoric combination of the six fields and a quadratic form in
two scalar combinations. To give rational coefficients
throughout we use independent field coordinates \((A',M)\), meaning
the physical first field is \(A=\sqrt{g}\,A'\), with \(g>0\).
Set \(p'=\operatorname{tr} A'\), \(q=\sum_jD_j:M_j\), and put
\[
 U=\operatorname{dev}A'+\sum_j\left[
  (\operatorname{tr}M_j)D_j+
  t\operatorname{dev}\frac{D_j\operatorname{dev}M_j+
                        (\operatorname{dev}M_j)D_j}{2}\right].
\]
With \(v=(\alpha,\beta)^t\) define
\[
 B=\operatorname{diag}(g b_1,b_2)-\frac{r vv^t}{6(1-\ell)},\qquad
 T=L_{\rm ref}-r|U|^2-(p',q) B(p',q)^t .
\]
Here \(g=r/a\), where \(a\) is the positive rational function in the table.
The reference phase is the first in the upper test and second in the
lower.
\[
\begin{array}{c|cc}
 &\text{upper}&\text{lower}\\ \hline
 t &(121X+16)/(20(X+1)) &8(73-9X)/(316-87X)\\
 r &21/[2(47t^2-36t+405)] &168/(107t^2-96t+312)\\
 a &2[1+(18+7t^2/3)r] &[1+(1+7t^2/24)r]/4\\
 \ell &16(3-t)^2r/17 &(6+t)^2r/108\\
 (b_1,b_2)&(14/3,7/68)&(7/36,7/9)\\
 (\alpha,\beta)&(2,-(18+11t)/17)&(-1,(12+11t)/18)
\end{array}
\]

We detail the finite polynomial certificate. The parametrizations
increase and have ranges \([19/5,254/55]\), \([61/25,14/3]\).
On each respective interval \(r,a,1-\ell\) are positive.
\(B\succeq0\) and hence \(L_{\rm ref}-T\succeq0\).
The full \(36\)-by-\(36\) matrix \(L_{\rm other}-T\) and the
\(18\)-by-\(18\) restriction of \(T\) to the wave space normal
to axis three are positive semidefinite. To replicate these
calculations use tensor coordinates \((11,22,33,12,13,23)\).
An off-diagonal basis matrix has two entries equal to one.
Use loads \(\operatorname{diag}(1,-1,0),\operatorname{diag}(1,1,-2)\)
and the three off-diagonal basis matrices. Divide each summand
indexed by a load by that load's squared norm. Compute the
matrices by polarization of the displayed forms (in \(L_i\)
multiply the first diagonal block by \(g\)). For a wave retain
coordinates \((33,13,23)\) in each strain and \((11,22,12)\) in
each stress. Successively remove the first variable. For a nonzero
diagonal pivot \(p\), replace the remaining block \(H\) by
\(H-vv^t/p\), with \(v\) the first column below the pivot.
All these pivots are strictly positive on the range of \(t\), and
a zero pivot has identically zero row. Each rational sign is tested
by positive numerator and denominator at the left endpoint after
choosing their common sign, and zero real roots on the closed
interval by Sturm's rule. These steps prove positivity by congruence.
The sums defining the form contract a basis and its dual, so the form
is invariant under \(A'\mapsto RA'R^t,\ M(D)\mapsto RM(R^tDR)R^t\).
The one wave test therefore proves all Fourier wave inequalities.

On the scalar plane in the new coordinates \(U=0\)
because \(\sum D_j^2\) is rotation invariant. The two final matrices
to test are respectively
\[
 \operatorname{diag}\bigl(g(X+2),(1/Y+1/2)/10\bigr)-B,
 \qquad
 \operatorname{diag}\bigl(g(1/X+1)/9,2(Y+1)/5\bigr)-B .
\]
In each, the first diagonal is strictly positive throughout
\([4/3,8/3]\). Substitution of the claimed bound into its respective
determinant gives identically zero after clearing denominators,
and the determinant increases strictly with \(1/Y\) for the upper
test, with \(Y\) for the lower. This gives \eqref{p14t:curve}
including the endpoints. The matrix construction, pivot signs,
scalar identities and denominator signs are asserted with rational
arithmetic in \path{certificates/baseline/p14translation/certificate.py}.
\phantomsection\label{proofend:p14t:bound}\end{proof}

\begin{corollary}\label{p14t:comparison}
For the spatial phases \((2,2)\) and \((1,1)\) and any \(0<f<1\),
the Berryman-Milton comparison region $\mathcal B$ consists of all shifted pairs \((X,Y)\) for which
some \(\zeta\in[0,1]\), with \(u=1+\zeta\) and \(v=1-\zeta/2\),
satisfies
\[
 \frac4{3v}\le X\le\frac{4u}3,\qquad
 \frac{15+104v}{46v+80v^2}\le Y\le\frac{208u+15u^2}{160+46u}.
\]
At every \(4/3<X<8/3\), the lower bound in \eqref{p14t:curve}
strictly exceeds \(\min\{Y:(X,Y)\in\mathcal B\}\), and the upper bound
is strictly less than \(\max\{Y:(X,Y)\in\mathcal B\}\).
\end{corollary}
\begin{proof}\phantomsection\label{proof:p14t:comparison}
This is the three-point bulk-shear comparison of Berryman and Milton for the two phases in the statement \cite[authors' manuscript, \S\S4 to 5, Eqs.~(18)-(19), (34)-(42), (61)-(62)]{bm88}. In their variables $5\zeta/21\le\eta\le(16+5\zeta)/21$. With phase 1 the stiffer phase, the upper shear uses $\alpha=1+\eta=(32+5u)/21$ and the lower uses $\beta=1-\eta/2=(16+5v)/21$. Substitution gives
\[
 \frac{39\alpha u-5u^2}{6(u+5\alpha)}=\frac{208u+15u^2}{160+46u},
 \qquad
 \frac{227v+45\beta}{6(25v^2+23\beta v)}=\frac{15+104v}{46v+80v^2}.
\]
The normalized shear parameter gives the same elimination in \cite[Eq.~(1.8), Eqs.~(2.24), (2.33)-(2.35), pp.~1589, 1593 to 1594]{GT95}.
Both shear endpoints increase with \(\zeta\), by differentiation.
At fixed \(X\) the greatest ceiling uses \(u=3-8/(3X)\)
and the least floor uses \(v=3/2-3X/8\).
Subtract the new upper function from this ceiling, and this floor
from the new lower function. After dividing each difference by
\((X-4/3)(8/3-X)\), the resulting rational functions are strictly
positive on the closed interval. The same Sturm test in the
certificate checks their numerators and denominators. For example
at \(X=16/9\) the ceiling decreases from \(1383/916\) to \(1397/974\),
and at \(X=2\) the floor increases from \(62/53\) to \(979/796\). Since \(H\) increases strictly, the
improvement holds for the physical shear at every fraction.
\phantomsection\label{proofend:p14t:comparison}\end{proof}

Two further complete finite-phase strips, including the proofs of their
Fourier seed enclosures, are given in Supplement~\ref{sup:common}.
\section*{Displacement and trace estimates}
\begin{lemma}\label{basic:korn}
On a bounded connected Lipschitz domain $U\subset\mathbb R^d$, $d\ge2$, every $v\in H^1(U;\mathbb R^d)$ satisfies
\[
 \inf_{r(x)=a+Sx,\ S^t=-S}\|v-r\|_{H^1(U)}
      \le C_U\|e(v)\|_{L^2(U)},\qquad
 \|v|_{\partial U}\|_{H^{1/2}(\partial U)}\le C_U\|v\|_{H^1(U)} .
\]
\end{lemma}
\begin{proof}\phantomsection\label{proof:basic:korn}
The rigid-motion estimate is Korn's inequality on a fixed Lipschitz domain
\cite[Thm.~A.1(ii), p.~460]{lew-muller}. For the trace apply
Gagliardo's inequality componentwise with $p=2$. His boundary norm on
p.~289 includes the double integral of squared differences divided by
$|x-y|^d$ \cite[Thm.~1.I, Eq.~(1.4), p.~289]{gagliardo-trace}.
For the $H^1$ extension in Theorem~\ref{p11:sharp}, flatten the boundary in Lipschitz graph
charts, reflect each local field across the coordinate plane and paste the
extensions with a partition of unity. The charts and reflections are
bi-Lipschitz and hence preserve $H^1$.
\phantomsection\label{proofend:basic:korn}\end{proof}

\par\smallskip
\textbf{Use of generative AI.} Generative AI tools were used to develop and check mathematical arguments, write computational checks, and draft the manuscript.

\clearpage
\phantomsection\label{supp:begin}
\setcounter{page}{1}
\renewcommand{\thepage}{S\arabic{page}}
\setcounter{section}{0}
\setcounter{secnumdepth}{2}
\setcounter{tocdepth}{2}\hypersetup{bookmarksdepth=2}
\renewcommand{\thesection}{S\arabic{section}}
\renewcommand{\theHsection}{supp.\arabic{section}}
\setcounter{theorem}{0}
\renewcommand{\thetheorem}{S\arabic{theorem}}
\renewcommand{\theHtheorem}{supp.\arabic{theorem}}
\pdfbookmark[0]{Mathematical supplement}{supplement-start}
{\raggedright\fontsize{23}{28}\selectfont Mathematical supplement\par}
\smallskip
{\large On fourteen problems in the theory of composites\par}
\smallskip
Philipp Bogdan \quad Imperial College London \quad 21 November 2026

The supplement follows the proof dependencies of the article. Equations
identified by an author-supplied numerical tag within a construction are
local to that construction. Page numbers beginning with S refer to this
supplement. The accompanying files contain the profile coefficients,
laminate trees and numerical bounds used in the computer-assisted results.
In the standalone supplement the numbered article pages, including the
shorter proofs and references, follow the S pages. In the combined
edition those article pages appear first. The PDF bookmarks and the
companion claim register locate each proof.

\begin{center}\small
\begin{tabular}{@{}llr@{}}\toprule
 & Derivation & Page\\\midrule
\ref{sup:planar} & Four rectangles and the field threshold & \pageref{sup:planar}\\
\ref{sup:cells} & Exact conductivity cells and finite readings & \pageref{sup:cells}\\
\ref{sup:arc} & Arc constructions and laminate exclusions & \pageref{sup:arc}\\
\ref{sup:inclusions} & Planar boundary certificates and a spatial obstruction & \pageref{sup:inclusions}\\
\ref{sup:stroke} & Finite-stroke recovery and initial tangent & \pageref{sup:stroke}\\
\ref{sup:rubber} & Mixed-load supremum and endpoint paths & \pageref{sup:rubber}\\
\ref{sup:crystal} & Polycrystal shapes, truncation and a zinc slice & \pageref{sup:crystal}\\
\ref{sup:common} & Complete bulk sections and cubic constraints & \pageref{sup:common}\\
\bottomrule\end{tabular}
\end{center}

\medskip
\textit{Reading the certificates.}
The archive records three distinct kinds of calculation. Rational
identities establish the constants in the field-count and elastic
examples. Interval arithmetic bounds entire functions or entire
parameter sets in the optical constructions. Floating-point
transmission solves test the reference boundaries at individual
contrasts. Only the first two enter the proofs.

\clearpage
\SSection{planar}{Four rectangles and the field threshold}
Use a single four-rectangle cell throughout this section. The two
potentials in the main example share the phase geometry. The scalar
rank-three argument permits a different laminate for each given tensor.
\phantomsection\label{file:two-potentials}
\paragraph*{A fixed pair with two potentials}

Write a gradient in column order as \((\partial_x u_1,\partial_y u_1,\partial_x u_2,\partial_y u_2)\). Let \(U\) have columns the first and fourth standard basis vectors of \(\mathbb R^4\). Set

\[
 M=10^{40},\qquad t=10^{-9},\qquad
 S=\begin{pmatrix}0&1\\1&0\end{pmatrix},\qquad
 L_0=\operatorname{diag}(1,M,M,1),\qquad L_1=L_0+tUSU^T .
\]

The eigenvalues of \(L_1\) are \(1-t,1+t,M,M\). On the unit torus, put \(u=\mathbf1_{[0,2/5)}(x)\), \(w=\mathbf1_{[0,1/5)}(y)\) and \(\chi=u+w-2uw\). Use \(L_1\) where \(\chi=1\) and \(L_0\) elsewhere. The phase-one fraction is \(f=11/25\).

For any effective tensor \(C\), put \(H=U^T(C-L_0)U\). The separating functional on a tensor and its phase-one fraction \(s\) is

\[
 \Phi(C,s)=9+27s-\frac{175H_{22}+5175H_{11}}{t^2}
                    +7500\,\frac{H_{12}-ts}{t^3}.                 \tag{1}
\]

Writing \(C_\square\) for the checkerboard tensor, we will prove

\[
 \Phi(C_\square,11/25)<10^{-5},\qquad
 \Phi(C_{\rm lam},s)>1/250                                      \tag{2}
\]

for every finite laminate, with its own fraction \(s\). At \(s=11/25\), this puts the operator-norm distance from \(C_\square\) to the closed laminate set at least \((1/250-10^{-5})/(5350/t^2+7500/t^3)>0\). For every \(n>2\) the same example works by adding uncoupled potentials with identical coefficients in both phases.

\paragraph*{Directional moments of a cell or a laminate}

Use normalized cell average, or the probability on tree leaves given by their volume fractions, for \(\langle\cdot\rangle\). For a bounded real even symbol \(b\) on unit directions, \(T_b\) on the torus is the Fourier multiplier \(b(k/|k|)\) on \(k\ne0\), zero on constants. On a finite binary laminate tree let \(D_I v\) equal the child mean of \(v\) minus its parent mean on the two children of node \(I\), and zero off \(I\). With \(n_I\) its layer normal, set \(T_b=\sum_I b(n_I)D_I\). The \(D_I\) are mutually orthogonal projections and sum to the mean-zero projection. The tensor of this tree is obtained by applying the layer energy minimization recursively from the leaves to the root.

Define \(\mu(b)=\langle\chi T_b\chi\rangle\) and \(K(b,c)=\langle\chi(T_b\chi)(T_c\chi)\rangle\). Thus \(\mu\) is the measure giving mass \(|\widehat\chi(k)|^2\) at each nonzero Fourier direction, or mass \(\langle(D_I\chi)^2\rangle\) at each tree node's direction, with coincident masses added. In either setting

\[
 \mu(1)=s(1-s)\le\tfrac14,\qquad
 \langle(T_b\chi)(T_c\chi)\rangle=\mu(bc),\qquad
 |K(b,c)|\le\|b\|_{L^2(\mu)}\|c\|_{L^2(\mu)}.                  \tag{3}
\]

The two background conductivity blocks give the symbols

\[
 p(n)=\frac{n_x^2}{n_x^2+Mn_y^2},\qquad
 r(n)=\frac{n_y^2}{Mn_x^2+n_y^2}.                               \tag{4}
\]

They obey \(p,r\ge0\) and \(p+r\le1\), since \(M\ge1\) and \(p\le n_x^2,\ r\le n_y^2\) for a unit \(n\). The statistic we need is

\[
 \Psi=9+27s+175\mu(p)+5175\mu(r)+7500K(p,r).                    \tag{5}
\]

\paragraph*{Reading the moments from one tensor}

Prescribe mean gradient \(Ua\), \(a\in\mathbb R^2\), and let \(v\) be the first and fourth coordinates of the local gradient. Then

\[
 v=a-t\mathcal G(\chi S v),\qquad
 \mathcal G=\operatorname{diag}(T_p,T_r),\qquad Ha=tS\langle\chi v\rangle .
\]

Indeed, at each nonzero Fourier mode the gradient of each potential is parallel to \(n\) and its flux perpendicular to \(n\). For background block \(A\) the gradient correction from a polarization \(h\) is consequently \(-nn^Th/(n^TAn)\). Restricting to the coordinate \(x\) in the first block and \(y\) in the second gives (4). On a tree the same calculation applies to \(D_I\) of the fields, since jumps of child-average gradients are normal and jumps of child-average fluxes are tangential. Summing the projections gives the displayed equation. Averaging the constitutive law gives the formula for \(H\).

The operators \(\mathcal G\), multiplication by \(\chi\), and \(S\) have norm at most one. Their Neumann expansion gives, in operator norm,

\[
 H=tsS-t^2
     \begin{pmatrix}\mu(r)&0\\0&\mu(p)\end{pmatrix}
     +t^3K(p,r)S+E,\qquad \|E\|\le\frac{t^4}{1-t}.              \tag{6}
\]

For the cubic coefficient, move the left Fourier multiplier or tree projection across the inner product by self-adjointness. Between the two outer factors \(S\), the \(ij\) entry is \(S_{ij}K(b_i,b_j)\) with \((b_1,b_2)=(p,r)\). This gives \(K(p,r)S\). In particular

\[
 |\Phi-\Psi|\le \frac{5350t^2+7500t}{1-t}<8\cdot10^{-6}         \tag{7}
\]

for the cell and for every tree. The estimate uses the single value of \(t\) in the given phase pair.

\paragraph*{A bound for every laminate tree}

First assign every node to one of three classes \(X,Y,Z\) in any manner. Write \(Y=25\sum_{I\text{ in }Y}D_I\chi,\ Z=25\sum_{I\text{ in }Z}D_I\chi\) at a leaf. Define

\[
 r_0=Y-3Z+3,\qquad r_1=Y+3Z-6,\qquad
 D=\langle(1-\chi)r_0^2+\chi r_1^2\rangle .
\]

We prove \(D\ge 1/225\) at every rank. Start \((X,Y,Z)\) at \((25s-11,0,0)\) at the root. At a split add 25 times the increment in conditional phase fraction to the assigned coordinate. Thus each split varies one coordinate, preserves its conditional mean, and the leaf sum \(\ell=X+Y+Z\) equals \(25\chi-11\). Use the polynomial

\[
 P=r_0r_1-2r_0(\ell-14)+r_1(\ell+11)
   =-XY+9XZ+8YZ-12X+24Y-36Z.                                   \tag{8}
\]

Put \(b=P/10,\ v=(X/(3/2),Y/(3/4),Z/(1/2))\), and \(\Omega=\{|b|<1,\ |v|<1\}\). The conditional mean of the leaf cost \(2|b|\) dominates \(2|b|\) everywhere and \(G=1+b^2-|v|^2\) in \(\Omega\). To see this, all leaves lie outside \(\Omega\), since \(11\le|\ell|\le(7/4)|v|\). Along a split \(b\) is affine, and \(G\) is convex within \(\Omega\). Indeed the three slopes of \(P\) are \(-Y+9Z-12,\ -X+8Z+24,\ 9X+8Y-36\). When \(|v|<1\), their absolute values exceed \(7,19,21\) by Cauchy-Schwarz. After division by 10 these exceed the corresponding norms of the slopes of \(v\).

The region inside \(\Omega\) on each split is an interval. Interpolating the two child costs dominates \(2|b|\) by convexity. It also dominates \(G\) on that interval, by convexity of \(G\) and \(G\le 2|b|\) on the boundary. Induction proves the bounds. On the root axis \(|b|=6|X|/5\ge |v|\). A root inside \(\Omega\) has \(G\ge1\), and a root outside has \(|b|\ge1\). Hence \(\langle |P|\rangle\ge5\).

On phase 0, (8) equals \(r_0(r_1+50)\), and on phase 1 it equals \(r_1(r_0+25)\). The two transforms are centered, with \(\langle Y^2+Z^2\rangle\le625/4\) by orthogonality. This gives

\[
 \langle |P|\rangle
 \le \sqrt D\,\big(\|r_1+50\|_2^2+\|r_0+25\|_2^2\big)^{1/2}
 \le \sqrt D\,\big(18\cdot625/4+44^2+28^2\big)^{1/2}
 <75\sqrt D
\]

whenever \(D>0\), and the first inequality rules out \(D=0\). This proves the tree bound. On a single pure leaf the costs \(D\) are 9 and 36 in phases 0 and 1 respectively.

Expanding the squares uses \(\langle YZ\rangle=0\) and \(r_1^2-r_0^2=12YZ-18Y-18Z+27\). With node weights \(m_I=\langle(D_I\chi)^2\rangle\), the result is

\[
 D=9+27s+175\sum_{I\text{ in }Y}m_I+5175\sum_{I\text{ in }Z}m_I
       +7500\left\langle\chi
          \left(\sum_{I\text{ in }Y}D_I\chi\right)
          \left(\sum_{J\text{ in }Z}D_J\chi\right)\right\rangle .
\]

Now assign each node independently to \(Y\) with probability \(p(n_I)\), to \(Z\) with probability \(r(n_I)\), and to \(X\) otherwise. Distinct nodes give products of probabilities. A single node cannot be in both classes, so the expected mixed term before multiplication by 7500 is

\[
 K(p,r)-\sum_I p(n_I)r(n_I)\langle\chi(D_I\chi)^2\rangle\le K(p,r).
\]

Thus \(\Psi\ge\mathbb E_{\rm assignments}D\ge1/225\), uniformly in the tree depth and in the normal at every node.

\paragraph*{The four-rectangle value}

Let \(g,h\) indicate the exact horizontal and vertical Fourier directions. Then \(T_g\chi\) is the column average minus \(f\), and \(T_h\chi\) is the row average minus \(f\). Hence \(Y=25T_g\chi=15(u-2/5)\) and \(Z=25T_h\chi=5(w-1/5)\). Where \(u=w\) we have \(r_0=0\), and where \(u\ne w\) we have \(r_1=0\). In particular

\[\begin{gathered}
 \mu(g)=54/625,\quad \mu(h)=4/625,\quad K(g,h)=-144/15625,\\
 9+27f+175\mu(g)+5175\mu(h)+7500K(g,h)=0 .
\end{gathered}\]

The functions \(p,r\) coincide with \(g,h\) on the axes. Off the axes, \(|\widehat\chi(k,l)|\le2/(\pi^2|kl|)\), by integrating the product \(-2uw\). Consequently

\[
 \mu(p-g)\le\frac{16}{\pi^4}\sum_{l\ge1}\frac1{l^2}
                 \sum_{k\ge1}\frac1{k^2+Ml^2}
 \le\frac{8}{\pi^3\sqrt M}\sum_{l\ge1}l^{-3}<M^{-1/2}.
\]

The same estimate holds for \(\mu(r-h)\). Both differences take values in \([0,1]\), so their \(L^2(\mu)\) norms are at most \(M^{-1/4}\). By (3), \(|K(p,r)-K(g,h)|\le M^{-1/4}\). Therefore the checkerboard satisfies

\[
 |\Psi|\le5350M^{-1/2}+7500M^{-1/4}<8\cdot10^{-7}.
\]

Together with (7) and the tree bound, this proves (2).

\paragraph*{Scalar realization at every fixed fraction}

Let the two scalar conductivity matrices be \(A,B\in\operatorname{Sym}^+(2)\). The indicator \(\chi\) can be measurable and periodic on any planar lattice. Write \(f=\langle\chi\rangle\) for the fraction of \(B\), and \(g=1-f\). The cases \(fg=0\) are immediate. With \(\Delta=B-A\), put

\[
 \sigma=A+\chi\Delta,\qquad V=gA+fB,\qquad S_0=fA+gB,\qquad
 q=\frac{\chi-f}{\sqrt{fg}} .
\]

We use an auxiliary matrix \(Y\) to retain the phase means in the scalar variational problem. A determinant bound on \(Y\) and the sequential lamination formula will suffice. The following derivation covers singular and indefinite phase differences.

\paragraph*{A determinant bound from planar duality}

Let \(\mathcal E\) be the closed subspace of mean-zero periodic gradient fields in \(L^2\). If all nonzero frequencies of \(\chi\) are collinear, the medium is stratified and the simple laminate at fraction \(f\) with that normal has its tensor. Otherwise \(K_0:\mathcal E\to\mathbb R^2,\ e\mapsto\langle q e\rangle\), is onto. Indeed, if a vector \(z\) annihilates its range, \(qz\) is orthogonal to gradients. Thus \(k\cdot z=0\) for every \(k\ne0\) in the support of \(\widehat q\), forcing \(z=0\).

Define \(Y=Y(\sigma)\) by

\[
 z^TYz=\min_{e\in\mathcal E,\,K_0 e=z}
            \left\langle(e-qz)\cdot\sigma(e-qz)\right\rangle .    \tag{8a}
\]

The minimizer exists, is unique and depends linearly on \(z\) by coercivity on the closed constraint space. Its remainder \(w=e-qz\) has zero average in each phase, since \(\langle w\rangle=\langle q w\rangle=0\). The matrix \(Y\) is positive definite. A zero cost would imply \(qz\in\mathcal E\), which requires \(z\) to be parallel to every frequency in that same support. For noncollinear support this forces \(z=0\). With the geometry fixed, (8a) also makes \(Y\) monotone in \(\sigma\) in the positive-semidefinite order and positively homogeneous of degree one.

Let \(R_\perp\) rotate planar vectors by \(90^\circ\) and write \(\sigma^{\rm d}=R_\perp\sigma^{-1}R_\perp^T\). We claim

\[
 Y(\sigma^{\rm d})=R_\perp Y(\sigma)^{-1}R_\perp^T .             \tag{9}
\]

For the minimizing \(e=qz+w\), the field \(j=\sigma w-qYz\) is mean-zero and solenoidal. In fact, \(\sigma w\) has zero mean in each phase, and the Euler equation and polarization in (8a) give \(\langle\sigma w\cdot h\rangle=(Yz)\cdot K_0 h\) for every \(h\in\mathcal E\). Now \(R_\perp j\) is an admissible mean-zero gradient at \(z'=-R_\perp Yz\) with remainder \(w'=R_\perp\sigma w\). The dual flux \(\sigma^{\rm d}w'-q(-R_\perp z)=R_\perp e\) is solenoidal. This is the same Euler equation for the dual minimum, and surjectivity of \(K_0\) gives \(Y(\sigma^{\rm d})z'=-R_\perp z\), proving (9).

For a symmetric positive-definite \(2\times2\) matrix \(D\), \(R_\perp D^{-1}R_\perp^T=D/\det D\). Set \(a=\min(\det A,\det B)\) and \(b=\max(\det A,\det B)\). Pointwise \(\sigma/b\preceq\sigma^{\rm d}\preceq\sigma/a\). Monotonicity, homogeneity and (9) give \(Y/b\preceq Y/\det Y\preceq Y/a\). Consequently

\[
                         a\le\det Y\le b.                         \tag{10}
\]

For a mean field \(d\), minimize the cell energy first over \(w\) at fixed \(z\). The terms involving \(w\) against phasewise constants vanish. The remaining energy is \(d^TVd+2\sqrt{fg}\,d^T\Delta z+z^T(S_0+Y)z\). Minimizing in \(z\) gives the effective tensor

\[
               C(Y)=V-fg\Delta(S_0+Y)^{-1}\Delta .               \tag{11}
\]

\paragraph*{Realizing the determinant interval}

Suppose first that \(Y>0\) and \(\det Y=\det A\). Then \(T=(A+Y)^{-1}\) satisfies \(\operatorname{tr}(AT)=1\). Indeed the eigenvalues of \(A^{-1/2}Y A^{-1/2}\) have product one. Diagonalizing \(A^{1/2}TA^{1/2}\) now gives weights \(\omega_i\ge0\) of sum one and at most two unit directions \(n_i\) such that

\[
 T=\sum_i\omega_i\frac{n_i n_i^T}{n_i^T A n_i}.
\]

Here one takes \(n_i\) proportional to \(A^{-1/2}\) times the corresponding normalized eigenvector.

Start with \(B\) and laminate successively with \(A\) along these normals. For invertible \(\Delta\), the resulting tensor \(C\) at fraction \(f\) can be chosen to satisfy

\[
                 f(C-A)^{-1}=\Delta^{-1}+gT.                    \tag{12}
\]

To verify this, retaining fraction \(\theta\) of an old tensor \(D\) and adding \(A\) with layer normal \(n\) gives \(\theta(D_{\rm new}-A)^{-1}=(D-A)^{-1}+(1-\theta)nn^T/(n^TAn)\). This follows by solving for the common normal flux and tangential field in the two layers, or by a rank-one inverse update in the layer formula. If the old \(B\) fraction is \(f_{\rm old}\), the contribution to the right of (12) is \((f_{\rm old}-\theta f_{\rm old})nn^T/(n^TAn)\). Choose successive decreases \(g\omega_i\) starting from fraction one.

Equation (11), with \(S_0=A+g\Delta\), gives exactly \(f(C(Y)-A)^{-1}=\Delta^{-1}+g(A+Y)^{-1}\). Thus this laminate gives \(C(Y)\). For singular \(\Delta\), keep \(A,Y\), the normals and the proportions fixed and approximate \(B\) by positive-definite matrices with invertible difference. Both the finite lamination formula and (11) are continuous, giving the same conclusion for \(B\). Interchanging the phases takes \(q,z,\Delta\) to \(-q,-z,-\Delta\) and preserves both (8a) and (11). Using \(B\) as the coating thus realizes every \(C(Y)\) with \(\det Y=\det B\), also at the required fraction.

It remains to realize \(C(Y)\) for any \(Y>0\) satisfying (10). The endpoints already suffice if \(a=b\). Otherwise fix a unit vector \(u\) and vary \(Y_s=Y+suu^T\) over its positive-definite interval. The determinant \(\det Y_s=\det Y(1+s u^TY^{-1}u)\) ranges from zero to infinity. The two levels \(a,b\) are reached on opposite sides of \(s=0\), allowing equality at zero. By the rank-one inverse formula,

\[
 C(Y_s)=C(Y)+fg\,\frac{s}{1+s u^T(S_0+Y)^{-1}u}\,zz^T,\qquad
 z=\Delta(S_0+Y)^{-1}u .
\]

The denominator equals \(\det(S_0+Y_s)/\det(S_0+Y)>0\) throughout the interval. The coefficient increases with \(s\), so \(C(Y)\) is an arithmetic mixture of the two endpoint tensors. If \(z\ne0\), choose a layer normal perpendicular to \(z\). Every constant gradient then gives continuous normal flux across these layers, so their effective tensor is precisely their arithmetic mixture. If \(z=0\) either endpoint suffices. Both endpoints already have the same phase fraction \(f\). This last layering has rank at most three and proves the scalar assertion for the full tensor.

The rank bound also passes to limits of periodic responses. A tree of depth at most three has finitely many phase labelings, and its normals and proportions range over compact sets, allowing zero weights. Both the effective tensor and total fraction depend continuously on these data. Thus the rank-three-or-less tensors at a specified fraction form a closed set.

\paragraph*{\texorpdfstring{A homogeneous two-well energy on \(\mathbb R^{2\times4}\)}{A homogeneous two-well energy on \textbackslash mathbb R\^{}\{2\textbackslash times4\}}}

The fraction term in (1) matters when taking a two-well envelope, since competitors can choose their own fraction. We first incorporate that term into the phase costs, then absorb the costs into homogeneous quadratic forms in the same dimension.

In (1), the symmetric coefficient of the active block \(H\) has diagonal entries \(-5175/t^2,-175/t^2\) and off-diagonal entry \(3750/t^3\). Its determinant is negative. Its spectral decomposition therefore gives loadings \(u,v\in\operatorname{ran} U\) and constants \(b,d\) such that

\[
            \Phi(C,s)=u^TCu-v^TCv+b+ds,\qquad d=27-7500/t^2.    \tag{13}
\]

Choose a common shift \(c\) so that \(c_i=b+di+c\ge1\) for \(i=0,1\). Let \(\mathcal R\) rotate the spatial vector in each potential by \(90^\circ\), and let \(L_i^{\rm d}=\mathcal R L_i^{-1}\mathcal R^T\). For mean-gradient arguments \(F,G\in\mathbb R^{2\times2}\), use the two wells

\[
 W(F,G)=\min_{i=0,1}\big(F\cdot L_i F+G\cdot L_i^{\rm d}G+c_i\big).
\]

The effective tensor for the dual coefficients on a given geometry is \(\mathcal R C^{-1}\mathcal R^T\). Indeed rotation exchanges a solenoidal flux and a compatible gradient, taking the solution pair \((e,j)\) to \((\mathcal Rj,\mathcal Re)\). The same argument works at every node of a laminate tree. This is the primal and complementary energy construction used in Cao's two-well argument.

Take \(p=C_\square v\) and \(\bar Z=(u,\mathcal R p)\). On any finite tree with assigned phases, effective tensor \(C\) and fraction \(s\), the minimum of the two block energies plus the average phase cost at this mean is

\[
 u^TCu+p^TC^{-1}p+b+ds+c
 \ \ge\ 2v^Tp+c+\Phi(C,s),
\]

by completing the square in \(p-Cv\). For the checkerboard equality holds, and its minimizing gradients in the two blocks give an admissible periodic trial field for \(W\). Thus \(QW(\bar Z)\le2v^Tp+c+\Phi(C_\square,11/25)\). Periodic \(L^2\) gradients suffice by smooth approximation and quadratic growth.

For the other direction, \(RW\) equals the infimum of the average leaf costs over finite rank-one splitting trees rooted at the prescribed mean. Indeed this infimum is below \(W\), is rank-one convex by joining two trees at a new root, and bounds every rank-one convex minorant from above by Jensen's inequality. Choose a minimizing phase at each leaf. Each rank-one split supplies a common layer normal for both gradient blocks. The specified leaf fields are therefore admissible in the corresponding laminate energy minimization, and their cost is at least its minimum. By (2), \(RW(\bar Z)\ge2v^Tp+c+1/250\). This proves the strict gap for \(W\).

Finally write \(Z=(F,G)\), \(D_i=\operatorname{diag}(L_i,L_i^{\rm d})\). There is a common unit eigenvector \(e\) with \(D_i e=Me\), namely the \(\partial_y u_1\) direction of the first block. For \(\varepsilon>0\) define \textbf{homogeneous} positive-definite quadratic forms

\[
 Q_i^\varepsilon(Z)=Z^TD_i Z+\varepsilon^2 c_i(e\cdot Z)^2 .
\]

Translate their minimum by evaluating at \(Z+e/\varepsilon\) and subtract the common affine function \(2M(e\cdot Z)/\varepsilon+M/\varepsilon^2\). The result is \(\min_i\{ Z^TD_iZ+c_i(1+\varepsilon e\cdot Z)^2\}\). Its difference from \(W(Z)\) is bounded by \(C\varepsilon(1+|Z|^2)\) for \(0<\varepsilon\le1\). Since \(W(Z)\ge 1+\min_i\lambda_{\min}(D_i)|Z|^2\), the translated function lies between \((1-C'\varepsilon)W\) and \((1+C'\varepsilon)W\). Both envelopes respect order, positive scaling, translations of the argument, and addition of affine functions. Their strict gap at \(\bar Z\) therefore persists, for all sufficiently small positive \(\varepsilon\), for \(\min_i Q_i^\varepsilon\) at \(\bar Z+e/\varepsilon\). The dimension four here is a sufficient dimension for two homogeneous wells.
\phantomsection\label{endfile:two-potentials}

\subsection{Precision of a finite panel}\label{route:precision-of-a-finite-panel}
For the checkerboard and a fixed laminate $T$, write
$D(A_1,A_2)=C_\square(A_1,A_2)-T(A_1,A_2)$. On the set
$0.1I\preceq A_1,A_2\preceq I$, bounds for both $\lambda_{\max}(D)$
and $\lambda_{\max}(-D)$ control \emph{all} unit load directions.
The certificate encloses the checkerboard by primal and complementary
trial matrices, computes the laminate by nodewise series and encloses
the difference on the full normalized phase cone. Degrees up to 17
and an angular grid of 256 enter the polynomial part. Taylor remainders
and coefficient errors are added before the final comparison. These
are uniform bounds over the phase cone, not maxima over sampled probes.
\[
\begin{array}{l|cc}
 & \text{fitted tree} & \text{after correction to }11/25\\\hline
 \text{first tree} & <0.000360 & <0.000466\\
 \text{second tree}& <0.000420 & <0.000434
\end{array}
\]
The second value in the last column is bounded in the recorded
calculation by $0.00043340419$. Each reading may therefore be displaced
to the midpoint by less than half the corresponding bound. The
certificates, computed with the same checkerboard enclosure and
different trees, accompany the article.

\subsection{Prescribing the fraction in the functional counterexample}\label{route:prescribing-the-fraction-in-the-functional-counterexample}
Let $\chi_0$ be the checkerboard with cuts $2/5,1/5$ and put $m=(1,5)$.
For $0<\theta\le1$ let
\[
 B_\theta=\mathbf1_{[0,\theta)}(\{5m\cdot x\}),\qquad
 \eta_\theta=1-B_\theta+B_\theta\chi_0 .
\]
Each rectangle spans an integer number of band periods horizontally.
Consequently $\langle\eta_\theta\rangle=1-14\theta/25$.
Taking $\theta=(25/14)\min(f,1-f)$ and interchanging the phases
when $f<1/2$ gives every $0<f<1$, including $\theta=25/28$ at half
volume. The direction introduced by the bands supplies the third
directional statistic needed at equal areas.

Here is the tree inequality used in the extension. Partition the
lamination normals into three selected classes and the remainder.
The martingale increments of the fraction, multiplied by 25, produce
coordinates $(X,Y,Z,M)$ starting at $(0,0,0,25f-11)$, with leaf sum
$25\chi-11$. For $r_0=Y-3Z+3$ and $r_1=Y+3Z-6+3M/7$ set
\[
 R=\mathbb{E}\big[(1-\chi)(r_0^2+M^2)+\chi r_1^2\big].
\]
A separately affine polynomial comparison gives
$1-f\le400\sqrt R+R/100$ on every tree. The horizontal, vertical
and band-normal Fourier projections of $\eta_\theta$ make the
corresponding residual zero. The response near a scalar background
determines these limiting moments. This excludes its entire function
from the pointwise laminate closure at every prescribed fraction.
The directional projections, leaf costs and passage to the limit are derived below.
For the original checkerboard we give an explicit separating
list of $6\cdot10^{10}$ positive phase pairs. A different checkerboard
admits a list of $3\cdot10^{15}$ rational pairs using only the horizontal
current. The proven separation can be extremely small.
\subsection{The directional moments and the separating finite panels}\label{detail:five-band-proof}\label{file:five-band-proof}
\paragraph*{1. Two coefficients of the scalar response}

Write \(D=\mathbb {RP}^1\) for the circle of unoriented directions. For a bounded function \(g\) on \(D\), define \(T_g\) on scalar periodic fields as the Fourier multiplier with symbol \(g([k])\) at \(k\ne0\) and symbol zero at the constant mode. Write \(\mathbb E\) for the cell average. Define an angular measure and a real symmetric bilinear form by

\[
\mu=\sum_{k\ne0}|\widehat\chi(k)|^2\delta_{[k]},
\qquad
K(g,h)=\mathbb E\big[\chi (T_g\chi)(T_h\chi)\big].
\tag{1}
\]

We use the complex bilinear extension of \(K\) when testing with complex functions. Parseval's identity gives

\[
\begin{split}
\mathbb E[\chi T_g\chi]&=\int g\,d\mu,\qquad
\mathbb E[(T_g\chi)(T_h\chi)]=\int gh\,d\mu,\\
|K(g,h)|&\le \|g\|_{L^2(\mu)}\|h\|_{L^2(\mu)},\qquad
\mu(D)=f(1-f).
\end{split}
\tag{2}
\]

There is an exact version on a finite laminate tree. Give each leaf its volume probability and let \(\chi\) be its phase indicator. At each internal node \(I\) with normal \(n_I\), let \(D_I b\) equal the child mean of a leaf field \(b\) minus its node mean, on each child, and zero outside \(I\). These are mutually orthogonal projections whose sum removes the constant part of a field. Use

\[
T_g=\sum_I g([n_I])D_I,\qquad
\mu=\sum_I \mathbb E[(D_I\chi)^2]\delta_{[n_I]}.
\tag{3}
\]

With leaf expectation in (1), equations (2) hold on the tree as well.

Take a real symmetric positive-definite background \(\sigma_2=A\) and a phase \(\sigma_1=A+t vv^T\), with \(v\in\mathbb R^2\) and \(t\) near zero. For a unit representative \(n\) of a direction, set

\[
m(n)=\frac{(v^T n)^2}{n^T A n}.
\]

In either setting, the effective tensor has the expansion

\[
\sigma_* = A+vv^T
\left(t f-t^2\int m\,d\mu+t^3 K(m,m)+O(t^4)\right).
\tag{4}
\]

To verify it, prescribe mean electric field \(e\), let \(E\) be the local electric field and put \(q=v^T E\). At each nonzero Fourier frequency, the field correction is parallel to \(k\) and the flux is perpendicular to \(k\). The two conditions give

\[
q=v^Te-tT_m(\chi q).
\tag{5}
\]

For a tree, the corresponding conditions are \(D_I E\) parallel to \(n_I\) and \(n_I^TD_I( A E+t\chi v q)=0\) at every node. They follow recursively from continuity of tangential field and normal flux at each layering step, applied to the child averages. Solving for \(D_I E\), taking its scalar product with \(v\), and summing gives (5) with (3). This argument also permits a node with several layers sharing one normal.

The average flux is \(Ae+tv\mathbb E[\chi q]\). Expanding (5) by a Neumann series and using self-adjointness of the real multiplier \(T_m\) gives (4). In both settings \(\|T_m\|\le \|v\|^2/\lambda_{\min}(A)\). Thus the analytic expansion is bounded on a common disk about zero for fixed \(A,v\), independently of the cell or the tree. Convergence of responses for real \(t\) near zero implies convergence of the displayed coefficients. For instance, every subsequential analytic limit on a smaller disk is determined by those real values.

\paragraph*{Passing tree inequalities to a limit}

If a response belongs to the pointwise lamination closure, take any countable list of phase pairs and choose trees matching increasing finite sublists to increasing accuracy. Along a subsequence their fractions tend to \(f\), their measures converge weakly to \(\lambda\), and their forms converge on all continuous pairs to a bilinear form \(\mathcal K\). Compactness of the direction circle, separability of its continuous functions, and (2) give this extraction. The bound in (2) extends \(\mathcal K\) to \(L^2(\lambda)\). In particular \(\lambda(D)=f(1-f)\) and \(\mathcal K(1,g)=(1-f)\lambda(g)\), since \(T_1\chi=\chi-f\) on each tree.

A closed inequality in the fraction, the masses of finitely many direction sets, and their pair moments passes to \(f,\lambda,\mathcal K\) whenever it holds for every tree and every disjoint choice of sets. First take finite unions of arcs with \(\lambda\)-null boundaries. Their masses converge by weak convergence. Approximate each arc indicator in \(L^2(\lambda)\) by a continuous function. For each fixed approximation its squared error has discontinuities only at the endpoints, so its integral against the tree measure converges to its integral against \(\lambda\). Equation (2) then gives convergence of the pair moments. Finally, disjoint Borel sets can be approximated in \(\lambda\)-measure by disjoint finite unions of these arcs. Indeed, choose disjoint compact subsets and open supersets with arbitrarily small mass between each pair. Cover the compact sets by disjoint neighborhoods consisting of finite unions of arcs with null endpoints inside the respective open sets. Continuity in \(L^2(\lambda)\) proves the assertion.

\paragraph*{A single real shear parameter}

For \(n=(n_1,n_2)\), write \(x=n_1^2,\ y=n_2^2,\ q=2n_1n_2\). Thus \(x+y=1\) and \(4xy=q^2\). Use the backgrounds

\[
A_s=\begin{pmatrix}1&s\\s&1\end{pmatrix},\qquad -1<s<1,
\] and define four functions on the direction circle by \[
Q_s=(1+s q)^{-1},\qquad H_s=(x-y)Q_s,\qquad a_s=xQ_s,\qquad b_s=yQ_s.
\]

The first contrast uses \(v_+=(1,1)^T/\sqrt2\). The second uses any \(v_-\) with \[
v_-v_-^T=A_s+\sqrt{1-s^2}\,\mathrm{diag}(1,-1).
\] The matrix on the right has trace 2, determinant zero, and nonnegative diagonal entries, so this is a real rank-one contrast. Their symbols in (4) are \[
m_+=(1+q)Q_s/2,\qquad m_-=1+\sqrt{1-s^2}\,H_s .
\] In particular \(Q_s=(1-2s m_+)/(1-s)\) and \(H_s=(m_-1)/\sqrt{1-s^2}\).

Write \(\mu(g)=\int g\,d\mu\). Equation (4) gives \(f,\mu(m_\pm),K(m_\pm,m_\pm)\). Since \(T_1\chi=\chi-f\), we also know \(K(1,g)=(1-f)\mu(g)\). The affine formulas for \(Q_s,H_s\) therefore give their first moments and their individual quadratic moments. Thus the two probes determine \[
U(s)=\mu(Q_s),\qquad
W(s)=9+27f+175\mu(a_s)+5175\mu(b_s)+7500K(a_s,b_s).
\] Indeed \(a_s=(Q_s+H_s)/2,\ b_s=(Q_s-H_s)/2\), so \(4K(a_s,b_s)=K(Q_s,Q_s)-K(H_s,H_s)\).

These formulas also define analytic functions \(U,W\) on \(|\Re s|<1\), with bounds on every compact subset independent of the geometry. This follows directly from \(|1+s q|\ge1-|\Re s|\) and (2). Values at imaginary \(s\) below are analytic evaluations of functions determined by the real probes.

If the target belongs to the pointwise lamination closure on any neighborhood of the identity, choose successively accurate trees on a countable set of tests with real \(s>0\) dense in a small interval and \(t=1/n\) sufficiently small. Both phase tensors are positive-definite and anisotropic on these tests. For each selected \(s\) and contrast, the common Neumann bound in (5) and uniqueness at the accumulating \(t\)-values give convergence of the first three coefficients of (4). In particular \(f\) and the real values of \(U,W\) converge to those for the target. The geometry-independent analytic bounds give locally convergent subsequences on the connected strip. Every subsequential limit agrees with the target on the real interval and hence throughout the strip. Thus \(U,W\) converge at every fixed imaginary argument as well.

\paragraph*{Recovery from a positive scalar increment}

Take \(\sigma_1=A_s+tI\) and \(\sigma_2=A_s\), with \(0<s<1\). Write \(a=a_s,\ b=b_s,\ c=n_1n_2 Q_s\). The gradient and flux conditions used in (5) give \(E=e-t\,\mathcal G_s(\chi E)\) in both the cell and the tree, where \[
\mathcal G_s=\begin{pmatrix}T_a&T_c\\T_c&T_b\end{pmatrix}.
\] Its symbol at a direction \(n\) is \(nn^T/(n^T A_s n)\), so \(\|\mathcal G_s\|\le(1-s)^{-1}\) independently of the geometry. Expanding \(A_s e+t\mathbb E[\chi E]\) gives \[
\sigma_*=A_s+t fI-t^2 L+t^3 B+O(t^4).\]\[
L=\begin{pmatrix}\mu(a)&\mu(c)\\ \mu(c)&\mu(b)\end{pmatrix},
\quad
B=\begin{pmatrix}
K(a,a)+K(c,c)&K(a+b,c)\\
K(a+b,c)&K(b,b)+K(c,c)
\end{pmatrix}.
\] For the third-order term, self-adjointness gives \(\mathbb E[\chi T_g(\chi T_h\chi)]=K(g,h)\) for each pair of real entries of the symbol.

Since \(a+b=1-2sc\), solving first for \(d=K(c,c)\) from the off-diagonal entry and then for the mixed moment from the trace gives \[
\begin{aligned}
d&=\frac{(1-f)L_{12}-B_{12}}{2s},\\
2K(a,b)&=f(1-f)^2-4s(1-f)L_{12}+(4s^2+2)d-\operatorname{tr} B .
\end{aligned}
\tag{5a}
\] Indeed \(2K(a,b)=K(a+b,a+b)-\operatorname{tr}B+2d\), and \(K(1,1)=f(1-f)^2\). Thus \(U(s)=\operatorname{tr}L\) and all the terms of \(W(s)\) follow from the first three coefficients of this single matrix response. The bound on \(\mathcal G_s\) gives convergence of those coefficients under pointwise approximation on the same real \(t\)-tests as above. Analytic continuation of \(U,W\) therefore applies to this family with fixed principal axes.

\paragraph*{One current measurement with proportional phases}

Fix a nonzero mean field \(e=(p,r)^T\) and use \(\sigma_2=A,\ \sigma_1=(1+t)A\). Write \(d=e^TAe\) and \(N_A(n)=nn^T/(n^TAn)\). The measured scalar has the expansion \[
e^T\sigma_* e=d(1+ft)-L_e t^2+B_e t^3+O(t^4).
\] Let \(J_A\) be the symmetric matrix field obtained by applying \(T_{(N_A)_{ij}}\) to \(\chi\) in each entry, and put \(w_e(n)=(p n_2-r n_1)^2\). The quadratic and cubic coefficients determine the weighted angular mass and the mean determinant through \[
\begin{aligned}
\mu\!\left(\frac{w_e(n)}{n^TAn}\right)
  &=\frac{d f(1-f)-L_e}{\det A},\\
\Delta(A):=\mathbb E[\chi\det J_A]
  &=\frac{(1-f)L_e-B_e}{d\det A}.
\end{aligned}
\tag{5b}
\] Indeed the field equation is \(E=e-t\mathcal G_A(\chi A E)\), where \(\mathcal G_A\) has matrix symbol \(N_A\). The Neumann expansion gives \(L_e=\mathbb E[\chi (Ae)^T J_A(Ae)]\) and, by self-adjointness, \(B_e=\mathbb E[\chi (Ae)^T J_A A J_A(Ae)]\). The Gram determinant identity \(d(n^TAn)-(e^TAn)^2=(\det A)w_e(n)\) gives the first formula. For the second, \(\operatorname{tr}(AJ_A)=T_1\chi=\chi-f\). The two-dimensional Cayley-Hamilton identity gives \[
J_A A J_A=\operatorname{tr}(AJ_A)J_A-(\det J_A)(\det A)A^{-1}.
\] Multiply on the left by \(\chi e^TA\), on the right by \(Ae\), and average to obtain (5b). The bound \(\|\mathcal G_A\chi A\|\le\|A\|/\lambda_{\min}(A)\) gives a common disk of convergence for the cell and every tree.

Now use \(A=A_s\). The two quantities in (5b) are \(\mu(w_e Q_s)\) and \(\Delta(A_s)=K(a_s,b_s)-K(c_s,c_s)\), with \(c_s=n_1n_2 Q_s\). Both are analytic on \(|\Re s|<1\). If trees approximate the readings, choose a dense countable set of small positive \(s\) and a sequence of positive \(t\) tending to zero at each \(s\). Coefficient convergence gives the cell's \(f\), and the limiting data \(\lambda,\mathcal K\) of the direction-moment argument agree with the cell on the two analytic quantities. The identity theorem gives agreement at \(s=iT\). As \(T\to\infty\), \(a_{iT},b_{iT}\) converge boundedly to the horizontal and vertical singleton indicators \(g,h\), while \(|c_{iT}|\le1/(2T)\). Consequently \[
\lambda(r^2g+p^2h)=\mu(r^2g+p^2h),\qquad
\mathcal K(g,h)=K(g,h).
\tag{5c}
\] Take \(e=(\sqrt{207},\sqrt7)^T\). The two mass terms in (7) are \(25\mu(7g+207h)\), so (5c) identifies the entire residual for the limiting data. the direction-moment argument passes the tree bound (11) to the two singleton indicators. Substitution of the four-rectangle values from the four-rectangle calculation would then give \[
\frac1{225}\ \le\ 9+27\frac{11}{25}+25\frac{1206}{625}-7500\frac{144}{15625}=0 .
\] This proves exclusion using one longitudinal reading per phase pair with both phases proportional and their principal axes fixed.

\paragraph*{Measuring only the horizontal current}

Fix the mean field \(e_1=(1,0)^T\) and measure \(F_B(t)=e_1^T\sigma_*(A_s+tB,A_s)e_1\). Put \(p_\pm=(1,\pm1)^T/\sqrt2\), \(\Pi=p_+p_+^T\), and \(\alpha=1+s,\ \beta=1-s\). In the \(p_\pm\) basis the symbol \(nn^T/(n^TA_sn)\) has diagonal entries \(g_\pm=(p_\pm^Tn)^2Q_s\) and off-diagonal entry \(k=(p_+^Tn)(p_-^Tn)Q_s\). They satisfy \[
\alpha g_++\beta g_-=1,\qquad Q_s=g_++g_-,\qquad H_s=2k .
\] Equation (4) for \(B=\Pi\) gives \(f,\mu(g_+),K(g_+,g_+)\), since \((e_1^Tp_+)^2=1/2\). Both \(g_-\) and \(Q_s\) are affine in \(g_+\), so their first and quadratic moments follow using \(K(1,g)=(1-f)\mu(g)\).

For the second contrast take \(B=A_s\), so that \(\sigma_1=(1+t)A_s\). Write \(F_{A_s}(t)=1+ft-Mt^2+Nt^3+O(t^4)\). Then \[
\begin{aligned}
2M&=\alpha^2\mu(g_+)+2\alpha\beta\mu(k)+\beta^2\mu(g_-),\\
2N&=\alpha^3K(g_+,g_+)+\beta^3K(g_-,g_-)
       +2\alpha\beta\big((1-f)\mu(k)+K(k,k)\big).
\end{aligned}
\] The quadratic formula of the direction-moment argument gives \(M\). For the cubic term put \(h_+=\alpha g_++\beta k\) and \(h_-=\alpha k+\beta g_-\). The first field correction has components \(-T_{h_\pm}\chi/\sqrt2\) because \(e_1=(p_++p_-)/\sqrt2\). Moving one self-adjoint multiplier across the inner product gives \(2N=\alpha K(h_+,h_+)+\beta K(h_-,h_-)\). Its cross terms combine to \(2\alpha\beta K(1,k)\), and \(\alpha+\beta=2\) gives the displayed coefficient of \(K(k,k)\).

Since \(\alpha\beta>0\), these two coefficients give first \(\mu(k)\) and then \(K(k,k)\). Together with the moments of \(Q_s\) they determine \(U,W\) on \(0<s<1\). The same common Neumann bounds and analytic continuation therefore apply to this single measured tensor entry, with both phases sharing the fixed axes \(p_\pm\).

\paragraph*{Sampling at rational phase pairs}

Let \(J\) be any infinite set of positive integers. Suppose the checkerboard response belongs to the pointwise lamination closure at the pairs \[
\sigma_{2,j}=A_{2^{-j}},\qquad \sigma_{1,j}=A_{2^{-j}}+3^{-j}I,\qquad j\in J .
\] Choose a sequence of trees converging at these countably many tests by approximating progressively larger finite subsets. On the complex bidisk \(|s|,|t|<1/4\), the symbol \(nn^T/(1+s q)\) is analytic in \(s\) with norm at most \(4/3\). The Neumann series for the field is jointly analytic there. The mean flux gives the uniform bound \(\|\sigma_*\|\le1+|s|+|t|/(1-4|t|/3)\le13/8\). Cauchy's coefficient estimates and diagonal extraction give a subsequence of tree responses converging locally uniformly on the bidisk.

An entry of the difference between this limit and the checkerboard is analytic, say \(H(s,t)=\sum_{m,n\ge0} h_{mn}s^m t^n\). It vanishes at \((2^{-j},3^{-j})\) for all \(j\in J\) with \(j\ge3\). If some coefficient were nonzero, take the least integer \(d=2^m3^n\) having \(h_{mn}\ne0\). Unique prime factorization gives a single pair \((m_0,n_0)\) with this product. Then \[
\lim_{j\to\infty}d^j H(2^{-j},3^{-j})=h_{m_0n_0}.
\] Every other nonzero term has a factor \((d/(2^m3^n))^j\) tending to zero. Absolute convergence at \((2^{-3},3^{-3})\) supplies a summable bound after multiplication by \(d^3\) for all \(j\ge3\). Taking the limit along \(J\) contradicts \(h_{m_0n_0}\ne0\). Thus the same subsequence of trees would converge to the checkerboard on the whole bidisk. the scalar-gap calculation excludes this on its positive real parameters, which proves exclusion on \(J\) as well. By the definition of pointwise closure, some finite subset of each such \(J\) already has a positive lower bound on the maximum tensor error over that subset, uniform over all trees.

The same uniqueness argument applies to the single reading in (5b) with \(\sigma_{1,j}=(1+3^{-j})A_{2^{-j}}\) and \(\sigma_{2,j}=A_{2^{-j}}\). On the same bidisk the field correction operator has norm at most \(|t|(1+|s|)/(1-|s|)\le5/12\). The scalar responses are therefore uniformly bounded and jointly analytic. Agreement on any infinite subsequence of these rational pairs would give agreement throughout the bidisk, contradicting (5c), (6), and (11) for \(e=(\sqrt{207},\sqrt7)^T\).

For the horizontal current alone on the original checkerboard, use the rational phase pairs \[
\sigma_{2,j}=A_{2^{-j}},\qquad
\sigma_{1,j}=(1+5^{-j})A_{2^{-j}}+3^{-j}\Pi,\qquad j\in J .
\] The scalar response for \(\sigma_2=A_s,\ \sigma_1=(1+w)A_s+z\Pi\) is uniformly bounded and jointly analytic in \(|s|,|z|,|w|<1/8\). Indeed the Neumann operator has norm at most \((8/7)(|z|+|w|(1+|s|))\le17/56\). The preceding argument applies with \(d=2^m3^n5^\ell\), using absolute convergence at \((2^{-4},3^{-4},5^{-4})\). Convergence of the measured entries on \(J\) would give a common subsequential limit equal to the checkerboard throughout this polydisk. Restricting to \(w=0\) and \(z=0\) gives exactly the two scalar contrasts above, so the scalar-gap calculation again gives a contradiction. The phase difference at each sampled pair is positive-definite.

\paragraph*{3. The four rectangles}

Separate \(D\) into the frequency direction \([e_1]\), the frequency direction \([e_2]\), and their complement. Let \(Y,Z,X\), respectively, be 25 times the three corresponding transforms \(T_g\chi\), using their indicator functions for \(g\). With \(p_1=2/5\) and \(p_2=1/5\), direct Fourier separation gives

\[
Y=25(1-2p_2)(u-p_1),\qquad Z=25(1-2p_1)(w-p_2),\qquad
X=-50(u-p_1)(w-p_2).
\]

Indeed the product of the two centered one-variable functions has both frequency indices nonzero. The values on the four rectangles are

{\def\LTcaptype{none} 
\begin{longtable}[]{@{}rrrrrr@{}}
\toprule\noalign{}
\(u\) & \(w\) & \(\chi\) & \(X\) & \(Y\) & \(Z\) \\
\midrule\noalign{}
\endhead
\bottomrule\noalign{}
\endlastfoot
0 & 0 & 0 & $-4$ & $-6$ & $-1$ \\
0 & 1 & 1 & 16 & $-6$ & 4 \\
1 & 0 & 1 & 6 & 9 & $-1$ \\
1 & 1 & 0 & $-24$ & 9 & 4 \\
\end{longtable}
}

In particular, with \(\ell=X+Y+Z\), put

\[
r_0=Y-3Z+3,\qquad r_1=Y+3Z-6.
\]

The table gives \(\ell=25\chi-11\) and

\[
\mathbb E\big[(1-\chi)r_0^2+\chi r_1^2\big]=0.
\tag{6}
\]

The mixed moment in this identity is \(\mathbb E[\chi YZ]=-144/25\). The two rectangles with \(\chi=1\) have areas \(3/25\) and \(8/25\) and carry values \(YZ=-24\) and \(-9\), respectively.

For any cell or tree, take two disjoint direction sets with indicators \(g,h\). Put \(Y=25T_g\chi\), \(Z=25T_h\chi\), and assign the complementary directions to \(X\). Denote the left side of (6) by \(R(g,h)\) for these variables. Expanding the squares and using (2) gives

\[
R(g,h)=9+27f+175\mu(g)+5175\mu(h)+7500K(g,h).
\tag{7}
\]

Indeed \(\mathbb E Y=\mathbb E Z=\mathbb E[YZ]=0\) and \(r_1^2-r_0^2=12YZ-18Y-18Z+27\). For the checkerboard and the two exact axis indicators, \(\mu(g)=54/625\), \(\mu(h)=4/625\), and \(K(g,h)=-144/15625\), giving zero in (7). On a tree the same expansion applies when each node is assigned directly to one of the three classes. In that case \(T_g,T_h\) sum the projections in their respective classes, even if nodes with equal normals receive different assignments.

\paragraph*{A rank-independent obstruction for laminate trees}

Follow a random root-to-leaf path with the volume probabilities. At each node let \(F\) be the conditional phase-1 fraction. Its increments have conditional mean zero. Start the three coordinates at \((25f-11,0,0)\) and add \(25\,\Delta F\) to the coordinate assigned to each split. The leaf values of \(Y,Z\) are the two transforms in (7). The offset in \(X\) gives the same two possible leaf sums at every phase fraction, since

\[
\ell=X+Y+Z=(25f-11)+25(\chi-f)=25\chi-11 .
\tag{8}
\]

In \(r_0r_1+c\,r_0(\ell-14)+d\,r_1(\ell+11)\), the coefficients of \(Y^2\) and \(Z^2\) are \(1+c+d\) and \(-9-3c+3d\). Cancelling both pure squares gives the polynomial

\[
\begin{split}
P(X,Y,Z)
&=r_0r_1-2r_0(\ell-14)+(\ell+11)r_1\\
&=-XY+9XZ+8YZ-12X+24Y-36Z .
\end{split}
\tag{9}
\]

We use a local comparison on each split. Let \(p\) be a vector or scalar affine in each coordinate separately, \(v\) an affine vector function, and \(H\) an affine function taking values in \([0,1]\) on the allowed convex domain. Put \[
\Omega=\{\|p\|<1,\ \|v\|<1\},\qquad G=H+\|p\|^2-\|v\|^2 .
\] Suppose \(\|\partial_i p\|\ge\|\partial_i v\|\) throughout \(\Omega\) in each coordinate \(i\). If a leaf cost dominates \(2\|p\|\) everywhere and \(G\) inside \(\Omega\), its conditional mean satisfies the same two bounds at every node. Indeed, interpolate the two child costs linearly along the split. The interpolation dominates \(2\|p\|\) by convexity of the norm. Its segment inside \(\Omega\) is an interval, and \(G\) is convex there by the slope bound. At an endpoint on the boundary of \(\Omega\), either \(\|p\|=1\) or \(\|v\|=1\), in both cases giving \(G\le2\|p\|\). At an endpoint that is a child in \(\Omega\), use the child bound. The interpolation consequently dominates \(G\) throughout the interval. Induction proves the claim. A split with several children can be grouped into binary splits on the same coordinate.

In (9), \(P\) is affine in each coordinate separately. Use \(p=P/10,\ H=1\), and \(v=(X/r_X,Y/r_Y,Z/r_Z)\) for \((r_X,r_Y,r_Z)=(3/2,3/4,1/2)\). The three slopes of \(P\) are \[
-Y+9Z-12,\qquad -X+8Z+24,\qquad 9X+8Y-36 .
\] For \(\|v\|<1\) their absolute values exceed \(7,19,21\), respectively. Indeed Cauchy-Schwarz bounds the variable parts by \(\sqrt{r_Y^2+81r_Z^2}<5,\ \sqrt{r_X^2+64r_Z^2}<5,\ \sqrt{81r_X^2+64r_Y^2}<15\). Thus \(r_i|\partial_iP|>10\) as required.

Take the leaf cost \(|P|/5=2|p|\). All leaves are outside \(\Omega\), since (8) gives \(11\le|\ell|\le(7/4)\|v\|\). On the root axis, \(|p|=6|X|/5\ge\|v\|=2|X|/3\). A root inside \(\Omega\) therefore has \(G\ge1\), and one outside has \(2|p|\ge2\). The comparison yields, at every rank and fraction, \[
\mathbb E|P_{\rm leaf}|\ \ge\ 5 .
\tag{10}
\]

At the leaf sums (8), the first expression in (9) gives \(r_0(r_1+50)\) where \(\chi=0\) and \(r_1(r_0+25)\) where \(\chi=1\). By Cauchy-Schwarz,

\[
\mathbb E|P_{\rm leaf}|\le
\sqrt{R(g,h)}\,
\big(\|r_1+50\|_2^2+\|r_0+25\|_2^2\big)^{1/2}.
\]

Here \(Y,Z\) are centered and \(\mathbb E[Y^2+Z^2]\le625/4\) by (2). Thus the sum of squared norms is \(2\mathbb E[Y^2]+18\mathbb E[Z^2]+44^2+28^2\le18\cdot625/4+44^2+28^2<75^2\). Combining with (10) gives, for every rank, every phase fraction, and every disjoint pair of direction sets,

\[
R(g,h)\ \ge\ \left(\frac5{75}\right)^2=\frac1{225}.
\tag{11}
\]

The argument also applies to direct assignments of nodes to classes.

\paragraph*{A scalar statistic with a strict gap}

Fix \(T\ge1\) and put \[
\alpha_T=\Re Q_{iT}=|Q_{iT}|^2=(1+T^2 q^2)^{-1},\qquad
D_T=\Re\{U(iT)-U(2iT)\}=\mu(\alpha_T-\alpha_{2T}).
\] The statistic will be \[
S_T=\Re W(iT)+5000 D_T .
\tag{12}
\]

\textbf{Laminate bound.} Independently at each internal node with normal \(n\), assign its class to \(Y\) with probability \(p(n)=x\alpha_T\), to \(Z\) with probability \(r(n)=y\alpha_T\), and to \(X\) with the remaining probability. These are valid probabilities since \(p+r=\alpha_T\le1\). Write \(\mathbb E_c\) for the average over class assignments and \(R_c\) for the square residual (7) on each assigned tree.

The averaged first moments in (7) are \(\mu(p),\mu(r)\). For the mixed term set \(B_{IJ}=\mathbb E[\chi(D_I\chi)(D_J\chi)]\). Independence at distinct nodes and mutual exclusivity at a single node give the averaged mixed term exactly as \[
K(p,r)-\sum_I p(n_I)r(n_I)B_{II}\ \le\ K(p,r).
\] Here \(B_{II}\ge0\). Since \(p=\Re a_{iT}\) and \(r=\Re b_{iT}\), complex bilinearity gives \[
K(p,r)=\Re K(a_{iT},b_{iT})+K(\Im a_{iT},\Im b_{iT}).
\] By (2) and \(2uv\le u^2+v^2\), the absolute value of the last term is at most \[
\tfrac12\mu\big((x^2+y^2)(\Im Q_{iT})^2\big)
\ \le\ \tfrac12\mu\big(\alpha_T(1-\alpha_T)\big)
\ \le\ \tfrac23 D_T .
\] The last inequality follows pointwise by subtracting \(1/(1+4T^2q^2)\) from \(1/(1+T^2q^2)\). Equation (7) gives \(\mathbb E_c R_c\le S_T\), and (11) holds for each assignment. Therefore \[
S_T\ \ge\ \mathbb E_c R_c\ \ge\ \frac1{225} .
\tag{13}
\] Section 2 passes (13) to pointwise limits at each fixed \(T\). The lower bound is uniform in the phase fractions of the approximants.

\textbf{A test using line averages.} For a unit-square periodic cell, write \(m_1(x_1)=\int_0^1\chi(x_1,y)\,dy\) and \(m_2(x_2)=\int_0^1\chi(x,x_2)\,dx\). Let \(g_0,h_0\) indicate the exact horizontal and vertical frequency directions. Projection onto those modes gives \(T_{g_0}\chi=m_1-f\) and \(T_{h_0}\chi=m_2-f\). Set \(Y=25(m_1-f)\) and \(Z=25(m_2-f)\).

If the cell's response is in the lamination closure, it satisfies (13) at each fixed \(T\). Take \(T\to\infty\) on the cell. The functions \(a_{iT},b_{iT}\) converge boundedly to \(g_0,h_0\), hence in \(L^2(\mu)\), and \(D_T\to0\). By (2), (7), and (12), \(S_T\to R(g_0,h_0)\). This gives the necessary line-average condition \[
\left\|Y+3(2\chi-1)Z+3-9\chi\right\|_{L^2}\ \ge\ \frac1{15}.
\] For the checkerboard the field inside the norm equals \(15(u-w)\) on \(\chi=0\) and \(15(u+w-1)\) on \(\chi=1\), by the four-rectangle calculation. Both vanish pointwise, and \(f=11/25\). This violates the necessary condition and proves exclusion even for limits of unbounded rank.

\textbf{Explicit checkerboard bound.} Set \(E_T=\int_{\text{off-axis}}\alpha_T\,d\mu\) and abbreviate \(a=a_{iT},\, b=b_{iT}\). These symbols equal \(g_0,h_0\) on both axes, and \[
\|a-g_0\|_{L^2(\mu)}^2+\|b-h_0\|_{L^2(\mu)}^2
=\int_{\text{off-axis}}(x^2+y^2)\alpha_T\,d\mu\ \le\ E_T .
\]

Let \(B\) be the expression on the left of (6) with \(Y=25T_a\chi,\ Z=25T_b\chi\), using the ordinary complex square. The corresponding real residual with \(g_0,h_0\) vanishes pointwise. Subtract it inside each square. By (2), the triangle inequality, and \((u+3v)^2\le10(u^2+v^2)\) for nonnegative \(u,v\), this gives \(|B|\le6250 E_T\). More explicitly, the difference inside the appropriate square at each point is \(25(T_{a-g_0}\chi+3(2\chi-1)T_{b-h_0}\chi)\).

Expanding the two squares before taking real parts gives \[
W(iT)-B=625\mu\big(a+9b-(a-3b)^2\big).
\] The real part of its integrand inside \(\mu\) is \[
2(x-3y)^2\alpha_T(1-\alpha_T)+4q^2\alpha_T,
\] using \(x+9y-(x-3y)^2=4q^2\). It vanishes on the axes and is at most \(22\alpha_T\) elsewhere. Also \(D_T\le E_T\). Therefore \(S_T\le(6250+625\cdot22+5000)E_T=25000 E_T\).

On off-axis frequencies \((k,l)\), the checkerboard has \(\widehat\chi(k,l)=-2\widehat u(k)\widehat w(l)\), with \(|\widehat u(k)|\le1/(\pi|k|)\) and the corresponding bound for \(w\). For \(k\ge l\ge1\), \(|q|=2kl/(k^2+l^2)\ge l/k\). Accounting for signs and both orderings of the magnitudes gives \[
E_T\le \frac{32}{\pi^4}\sum_{l\ge1}\frac1{l^2}
                     \sum_{k\ge1}\frac1{k^2+T^2 l^2}
\le \frac{16}{\pi^3 T}\sum_{l\ge1} l^{-3}
\le \frac{24}{\pi^3 T}<\frac1T .
\] The middle estimate bounds the inner decreasing sum by the integral from zero to infinity. Taking \(T=10^7\) gives \(S_T<1/400\) for the checkerboard. By (13), every laminate limit exceeds the checkerboard on this same statistic by more than \(1/225-1/400=7/3600\), for any phase fraction.

\paragraph*{An explicit gap for the horizontal current}

Take \(p_1=1/4,\ p_2=(1-1/\sqrt2)/2\), so \(\alpha=1/2,\ \beta=1/\sqrt2\). The coefficient \(\beta^{-2}-2\) of \(\mu(g)\) vanishes. On every tree the residual for these cuts satisfies \[
R=\frac{3-2\sqrt2}{16}+\frac{\sqrt2}{4}f+2\mu(h)+8\sqrt2 K(g,h)
   \ \ge\ 2^{-16}.
\tag{13a}
\] The polynomial in the proof above becomes \[
P=-2XY+4XZ+2YZ+\tfrac{\sqrt2}{4}(-X+2Y-3Z).
\] Use the comparison of the tree argument with \(p=64P,\ H=1,\ v=(8X,16Y,32Z)\). For \(\|v\|<1\), Cauchy-Schwarz gives lower bounds on \(|\partial_XP|/8,\ |\partial_YP|/16,\ |\partial_ZP|/32\) of \(\sqrt2/64,\ (8\sqrt2-\sqrt{17})/256,\ (6\sqrt2-\sqrt{17})/256\), respectively. All exceed \(1/64\). At the root \(|p|=16\sqrt2|X|\ge\|v\|\), and both possible leaf sums have magnitude greater than \(\sqrt{21}/32\), whereas \(|\ell|\le(\sqrt{21}/32)\|v\|\). The leaf cost \(2|p|\) therefore gives \(\mathbb E|P|\ge1/128\). The two leaf factors satisfy \[
\|r_1-c\|_2^2+\|r_0+d\|_2^2
 =4\mathbb E Y^2+8\mathbb E Z^2+11/8\le27/8<4 .
\] Thus \(\mathbb E|P|\le2\sqrt R\), proving (13a) for every rank and every assignment of splits to classes.

The two direction data in (13a) are determined by (5c) for \(e=(1,0)^T\). On this checkerboard \(f=1/2-\sqrt2/8,\ \mu(h)=1/32,\ K(g,h)=-3\sqrt2/256\), giving \(R=0\). the direction-moment argument passes (13a) to the limiting singleton data. The bidisk argument therefore excludes approximation of the \(11\) entry alone at \(\sigma_{2,j}=A_{2^{-j}},\ \sigma_{1,j}=(1+3^{-j})A_{2^{-j}}\) on every infinite set of positive \(j\).

\paragraph*{Five diagonal bands at any prescribed phase fraction}

Keep the original checkerboard \(\chi_0\) of fraction \(f_0=11/25\) at its original size. For \(m=(1,5)\) and \(0<\theta\le1\), set \[
B_\theta(x)=\mathbf1_{[0,\theta)}(\{5m\cdot x\}),\qquad
\eta_\theta(x)=1-B_\theta(x)+B_\theta(x)\chi_0(x).
\] Here braces denote the fractional part. The five bands retain the checkerboard and their complement contains pure phase 1. Each horizontal side of a checkerboard rectangle has length \(2/5\) or \(3/5\), an integer number of periods of \(B_\theta\) in \(x_1\). Hence the bands cover exactly the fraction \(\theta\) of each rectangle, and the phase-1 fraction is \[
f_\theta=1-(1-f_0)\theta=1-14\theta/25 .
\] Every \(\eta_\theta\) has its response outside the pointwise lamination closure. For a prescribed \(0<f<1\), take \(\theta=(25/14)\min(f,1-f)\) and interchange the two phases when \(f<1/2\). A phase interchange preserves exclusion because swapping the two leaf labels in every laminate exchanges its two tensor arguments. In particular, \(\theta=25/28\) gives a single cell with exactly equal areas.

For an equal-area four-rectangle checkerboard alone, one of the cuts is at \(1/2\). Translation by half a period in that coordinate interchanges the phases and changes the sign of every \(T_g\chi\). Averaging the two translates gives \(K(g,h)=\mu(gh)/2\), the same cubic moment as the laminate limit in the direction-moment argument. The bands supply a third direction for the test below.

\textbf{A tree bound for three direction classes.} Assign each split to \(Y,Z,M\), or the remaining class \(X\). Start at \((X,Y,Z,M)=(0,0,0,25f-11)\) and add \(25\,\Delta F\) in the assigned coordinate. At every node \(\ell=X+Y+Z+M=25F-11\) lies in \([-11,14]\). Put \[
k=3/7,\qquad r_0=Y-3Z+3,\qquad r_1=Y+3Z-6+kM,\qquad
R=\mathbb E[(1-\chi)(r_0^2+M^2)+\chi r_1^2].
\] The extra coordinate records the difference between the mean fractions in the two types of strip. On any tree we have \[
1-f\ \le\ 400\sqrt R+R/100 .
\tag{14}
\] To prove it, define \[
P=r_0r_1-2r_0(\ell-14)+(\ell+11)r_1-kM(\ell-14),
\qquad Q=M\big(r_1-k(\ell-14)\big).
\] With \(P_0(X,Y,Z)\) denoting the polynomial in (9), these expand to \[
P=P_0+M\big((k-1)Y+(9-3k)Z\big),\qquad
Q=M\big(-kX+(1-k)Y+(3-k)Z\big).
\] Both are affine in each coordinate separately. In the comparison of the tree argument take \(p=(P,Q)/2,\ v=4(X,Y,Z),\ H=(14-\ell)/25\) on the slab \(-11\le\ell\le14\). For \(\|v\|<1\), the norm of \(\partial_i(P,Q)\) exceeds \(8\) in each coordinate \(i=X,Y,Z\). This follows from \[
\partial_X P=-Y+9Z-12,\qquad
\partial_Y(P+Q)=-X+8Z+24,\qquad
\partial_Z(P-3Q)=9X+8Y-36
\] by dividing the respective absolute values by \(1,\sqrt2,\sqrt{10}\). Along \(M\) we have \(\partial_M v=0\), so the slope bound holds in all four coordinates.

Set \(L=2\|p\|\) and use the leaf cost \(L+(1-\chi)M^2/100\). At a phase-1 leaf inside \(\Omega\), \(H=0\) gives \(G\le\|p\|^2\le L\). At a phase-2 leaf inside \(\Omega\), \(\|v\|<1\) and \(\ell=-11\) imply \(|M|>10\), giving \(G\le L+1\le L+M^2/100\). The root has \(p=v=0,\ G=1-f\). The comparison therefore gives \(1-f\le\mathbb E L+R/100\). The pair \((P,Q)\) at phase-2 leaves is \((r_0(r_1+50)+25kM,\,M(r_1+25k))\), and at phase-1 leaves it is \((r_1(r_0+25),\,Mr_1)\). Using \(L\le|P|+|Q|\) and Cauchy-Schwarz gives \[
\mathbb E L\le\sqrt R\,\bigl(\|r_1+50\|_2+25k+\|r_1+25k\|_2+\|r_0+25\|_2+\|M\|_2\bigr).
\] The parenthesis is less than \(400\), since (2) gives \(\|Y\|_2,\|Z\|_2\le25/2\) and \(\|M\|_2\le14+25/2\). This proves (14).

\textbf{Angular data from one fixed current direction.} For every periodic cell, knowing \(e^T\sigma_*((1+t)A,A)e\) for any fixed \(e\ne0\), all \(A\) near \(I\), and small \(t>0\) determines its entire \(\mu,K\). Compare the cell with another periodic cell having those readings, or with the limiting \(\lambda,\mathcal K\) of the direction-moment argument if trees approximate them. In the latter case take a countable dense set of real symmetric positive-definite backgrounds in the neighborhood and \(t\) tending to zero at each background. The common field expansion gives agreement in (5b), and continuity extends it to the neighborhood. Write \(\lambda,\mathcal K\) for the comparison data in either case.

For unit \(n\), put \(X=n_1^2-n_2^2,\ Y=2n_1n_2\). Use backgrounds with \(n^TAn=1+uX+vY\) for small real \(u,v\). The power series of the first observable in (5b) gives \(\lambda(w_e X^jY^k)=\mu(w_e X^jY^k)\) for all \(j,k\ge0\). The pair \((X,Y)\) parametrizes the direction circle, and polynomials in these coordinates are dense in its continuous functions. Hence \(w_e\lambda=w_e\mu\). The weight is positive off the single direction \([e]\), and equality of the total masses \(f(1-f)\) then gives \(\lambda=\mu\).

The determinant observable is \(K((N_A)_{11},(N_A)_{22})-K((N_A)_{12},(N_A)_{12})\). It is real analytic in \(A\) on the positive-definite cone, since the denominator \(n^TAn\) is uniformly positive on the circle at each background. Agreement near \(I\) therefore gives agreement throughout that cone. For distinct directions \(i,j\) supporting Fourier modes, let \(g_i,g_j\) be their singleton indicators. Choose an invertible real matrix \(C\) sending their lines to \([e_1],[e_2]\), respectively, and put \(\pi(n)=Cn/|Cn|\). The identity \[
N_{C^TA_s C}(n)=C^{-1}N_{A_s}(\pi(n))C^{-T}
\] gives \[
(\det C)^2\Delta(C^TA_s C)
 =K(a_s\circ\pi,b_s\circ\pi)-K(c_s\circ\pi,c_s\circ\pi).
\] This identity holds also for the limiting form by bilinearity. As in (5c), continue in \(s\) to \(iT\) and let \(T\to\infty\). This gives \(\mathcal K(g_i,g_j)=K(g_i,g_j)\). For \(i=j\), use \(K(1,g_i)=(1-f)\mu(g_i)\) and sum the already identified terms with \(j\ne i\). Indeed, (1) puts all the measure on the countable Fourier directions, so their singleton indicators sum to \(1\) in \(L^2(\mu)\). Their span is dense there. Thus \(\mathcal K=K\) on the whole space.

Rank-one tests also determine these data. The real symbols \(\epsilon(n_i^Tn)^2/[\epsilon(n_i^Tn)^2+((n_i^\perp)^Tn)^2]\) in (4) converge boundedly to \(g_i\) as \(\epsilon\downarrow0\), giving each mass and diagonal moment. The two rank-one probes of Section 2 after congruence by \(C\) give the mixed moments by the same imaginary limit.

In particular (14) must hold for a periodic cell approximable on the proportional scalar tests, for any three disjoint Borel direction sets \(E_Y,E_Z,E_M\), with \[
Y=25T_{\mathbf1_{E_Y}}\chi,\qquad Z=25T_{\mathbf1_{E_Z}}\chi,\qquad
M=25T_{\mathbf1_{E_M}}\chi+25f-11.
\] Expanding \(R\) uses only the fraction, the three masses, and the six pair moments, so the limit principle in the direction-moment argument applies.

\textbf{The three projections of the banded cell.} For each rational direction write \((a,b)\) for either primitive integer representative. Take \(E_M=\{[m]\}\), let \(E_Y\) consist of the directions with \(5\mid b\) other than \([m]\), and let \(E_Z\) consist of those with \(5\mid a\). These are disjoint countable Borel sets. Decompose \(\chi_0-f_0=\phi_Y+\phi_Z+\phi_X\) as in the four-rectangle calculation, with \[
\phi_Y=\tfrac35(u-\tfrac25),\qquad
\phi_Z=\tfrac15(w-\tfrac15),\qquad
\phi_X=-2(u-\tfrac25)(w-\tfrac15).
\] The nonzero Fourier coefficients of the centered one-variable indicators occur only at indices indivisible by 5. Thus the frequencies of \(\phi_Y,\phi_Z,\phi_X\) have residues modulo 5 of the forms \((j,0),(0,l),(j,l)\), respectively, with each displayed letter nonzero modulo 5. Multiplication by \(B_\theta\) shifts a frequency by \(5h m\), \(h\in\mathbb Z\), preserving these residues. At least one shifted component is indivisible by 5, so passing to its primitive representative preserves which components are divisible by 5. Also none of the original frequencies lies on \([m]\), since a nonzero frequency \((j,5j)\) has a nonzero second index divisible by 5. A shift parallel to \(m\) preserves this exclusion. These support statements for the products follow by approximating \(B_\theta\) in \(L^2\) with trigonometric polynomials in \(5m\cdot x\), since the three \(\phi_i\) are bounded.

Write \(\bar\eta=1-(1-f_0)B_\theta\). Its nonconstant frequencies all lie on \([m]\), and \(\eta_\theta=\bar\eta+B_\theta(\phi_Y+\phi_Z+\phi_X)\). The three projections therefore give exactly \[
Y=25 B_\theta\phi_Y,\qquad Z=25 B_\theta\phi_Z,\qquad
M=25(\bar\eta-f_\theta)+25f_\theta-11=25\bar\eta-11 .
\] Inside the bands, \(M=0\) and (6) gives zero residual pointwise. Outside, \(\eta_\theta=1,\ Y=Z=0,\ M=14\), hence \(r_1=-6+(3/7)14=0\). Thus \(R=0\), whereas (14) requires \(14\theta/25\le400\sqrt R+R/100\). This proves exclusion for every \(\theta>0\) in the stated range. Interchanging phases sends the data to \(1-f,\mu,\ (g,h)\mapsto\mu(gh)-K(g,h)\) by (2). Consequently the same single-current exclusion holds for the relabelled examples as well. Swapping the approximating trees' leaf labels would otherwise give limiting data for the original cell in violation of (14).

\textbf{One sequence for all fractions.} Fix \(e=(1,0)^T\) and let \[
A(u,v)=\begin{pmatrix}1+u&v\\v&1-u\end{pmatrix},\qquad
\sigma_{2,j}=A(2^{-j},3^{-j}),\qquad
\sigma_{1,j}=(1+5^{-j})\sigma_{2,j}.
\] Measuring just the \(11\) entry on any infinite subsequence distinguishes every banded example above from laminate limits. Indeed, for \(|u|,|v|,|t|<1/8\), the field correction operator for phases \((1+t)A,A\) has norm at most \(|t|(1+|u|+|v|)/(1-|u|-|v|)\le5/24\). The scalar responses are uniformly bounded and jointly analytic there. The uniqueness argument of Section 2 with \(2^m3^n5^\ell\) would turn approximation on the subsequence into convergence along a further subsequence throughout the tridisk. Multiplying both phase tensors by the same positive scalar multiplies the response by that scalar. Thus these trace-two backgrounds give the readings for all backgrounds near \(I\), contradicting the directional obstruction just proved.

\paragraph*{A finite list of ordinary real tensor values}

For the checkerboard with cuts \(2/5,1/5\), the following finite list separates its response from every laminate using only positive contrasts. Set \(T=10^7\), \(N=1000T\), and choose any \(0<\beta\le1/(100T)\). Use \[
z_j=j\beta/N,\qquad s_j=\frac{2z_j}{1+z_j^2}\quad(1\le j\le N),\qquad
h=10^{-30}\left(\frac{\beta}{3N}\right)^N .
\] At each \(s_j\) take both contrasts of Section 2 and \(t=h,2h,3h\). These are \(6N\) real positive-definite phase pairs. Both phases have positive off-diagonal entry and hence distinct eigenvalues. Decreasing \(\beta\) places the entire list in any prescribed neighborhood of the identity.

To verify sufficiency, write \(s(z)=2z/(1+z^2)\). The functions \(U(s(z)),W(s(z))\) are analytic for \(|z|<1\). Indeed \(Q_{s(z)}=(1+z^2)/(1+2qz+z^2)\), whose denominator factors into \(1+\xi z\) and \(1+\bar\xi z\) with \(|\xi|=1\). On \(|z|=\rho=1-1/(10T)\) the symbols \(Q,H,a,b\) have magnitude at most \(2/(1-\rho)^2=200T^2\). By (2) both analytic functions are bounded there by \(10^9T^4\) for every geometry.

For \(L=T,2T\), write \(w=iL/(1+\sqrt{1+L^2})\), so \(s(w)=iL\) and \(|w|\le1-1/(3T)\). For either analytic function the error at \(w\) of polynomial interpolation at the \(N\) nodes \(z_j\) is bounded by \[
\frac{\rho\,10^9T^4}{\rho-|w|}
       \left(\frac{|w|+\beta}{\rho-\beta}\right)^N
 \le 5T\,10^9T^4 e^{-N/(5T)}<10^{-40}.
\] This is Cauchy's interpolation remainder on \(|z|=\rho\), using the product \(\prod_j(z-z_j)\). The ratio in parentheses is at most \(1-1/(5T)\). The sum of the absolute Lagrange basis values at \(w\) is at most \(N(2N/\beta)^{N-1}<(3N/\beta)^N\), since the real nodes have spacing \(\beta/N\).

At each real node the response gives the scalar \[
G(t)=\frac{v^T(\sigma_*-A_s)v}{\|v\|^4}.
\] Here \(0\le m_\pm\le2\), for the minus sign by \((1+s q)^2-(1-s^2)(x-y)^2=(s+q)^2\). The series (5) bounds the remainder in \(G(t)/t=f-\mu(m)t+K(m,m)t^2+\cdots\) by \(16t^3\) for \(0<t\le1/4\). Quadratic interpolation at \(h,2h,3h\) estimates each of the three coefficients to within \(4000h\). Indeed the absolute row sums of the coefficient formulas are \(7,8/h,2/h^2\), and the largest remainder is \(432h^3\).

Suppose a tree matched every checkerboard tensor on the list to within \(h^4\) in operator norm. Since \(\|v_\pm\|^2\ge1\), the sampled values of \(G(t)/t\) would differ by at most \(h^3\). The coefficient estimates would differ by at most \(2h\) by the same row-sum bounds. Thus the differences in the true \(f,\mu(m_\pm),K(m_\pm,m_\pm)\) at each node would be bounded by \(d=9000h\).

The affine coefficients expressing \(Q_s,H_s\) in terms of \(1,m_\pm\) have absolute value at most 2. Since \(K(1,m)=(1-f)\mu(m)\), the resulting discrepancies are at most \(4d\) for each first moment of \(Q_s,H_s\) and \(24d\) for each quadratic moment. Here \(f(1-f)\) and \(f(1-f)^2\) are 1-Lipschitz on \([0,1]\), and \(0\le\mu(m_\pm)\le1/2\). Substitution into \[
W=9+27f+2675\mu(Q_s)-2500\mu(H_s)
                +1875\big(K(Q_s,Q_s)-K(H_s,H_s)\big)
\] bounds the discrepancy of either \(U\) or \(W\) by \(10^{10}h\) at the real nodes. At either auxiliary argument the interpolation and Lagrange bounds then give a discrepancy below \(2\cdot10^{-20}\). Consequently the discrepancy in \(S_T\) would be less than \(10^{-14}\). Such a tree would have \(S_T<1/400+10^{-14}<1/225\), contradicting (13).

Every finite laminate therefore has operator-norm error greater than \(h^4\) on at least one of the \(6\cdot10^{10}\) specified pairs. The maximum over that fixed finite list is continuous for pointwise convergence, so the error bound \(h^4>0\) persists in the closure.

\textbf{Three scalar increments per background.} At the same \(s_j\), take \(\sigma_2=A_{s_j}\) and \(\sigma_1=A_{s_j}+tI\) at \(t=k,2k,3k\), where \(k=\beta h/N\). Suppose a tree matched the checkerboard to within \(k^4\) on this list. Since \(\|\mathcal G_{s_j}\|\le2\), the \(16t^3\) remainder estimate applies to \((\sigma_*-A_{s_j})/t\) in operator norm. The same quadratic interpolation then bounds the discrepancies in \(f,L,B\) by \(d_0=9000k\).

At each node, \(\beta/N\le s_j\le1/2\), and (2) gives \(|L_{12}|\le1/4\) since \(|c|\le1\). Hence the discrepancy in \((1-f)L_{12}\) is at most \(2d_0\). Equation (5a) bounds the discrepancies in \(K(c,c)\) and \(K(a,b)\) by \(3d_0/(2s_j)\) and \(4d_0/s_j\), respectively. In the second bound, use \(|\Delta\operatorname{tr}B|\le2d_0\) and the 1-Lipschitz bound for \(f(1-f)^2\). The discrepancies in \(U,W\) at the real nodes are therefore at most \(4\cdot10^4d_0/s_j<10^{10}h\). Interpolation at \(iT,2iT\) gives the same contradiction to (13). Thus the maximum tensor error on these \(3\cdot10^{10}\) phase pairs is at least \(k^4>0\) throughout the lamination closure.

\paragraph*{A finite list using only horizontal current}

For the checkerboard with cuts \(1/4,(1-1/\sqrt2)/2\), set \[
T=10^{12},\quad N=1000T,\quad b=\frac1{100T},\quad
z_j=\frac{jb}{N},\quad s_j=\frac{2z_j}{1+z_j^2},\quad
h=10^{-30}\left(\frac{b}{3N}\right)^N .
\] At \(1\le j\le N\), take \(\sigma_2=A_{s_j}\) and \(\sigma_1=(1+t)A_{s_j}\) for \(t=h,2h,3h\). The maximum error in the \(11\) entry on these \(3N\) rational phase pairs is at least \(h^4>0\) throughout the lamination closure.

Write \(V(s)=\mu(b_s)\), \(J(s)=\Delta(A_s)\) from (5b), and \(D_y=\Re(V(iT)-V(2iT))=\mu(y(\alpha_T-\alpha_{2T}))\). At a fixed \(T\), define \[
S^y_T=\frac{3-2\sqrt2}{16}+\frac{\sqrt2}{4}f+2\Re V(iT)
 +8\sqrt2\left(\Re J(iT)+\sqrt{D_y/12}+\frac1{16T^2}\right).
\] Every tree has \(S^y_T\ge2^{-16}\). To see this, use the independent assignment in the scalar-gap calculation in (13a). The mean vertical mass is \(\Re V(iT)\), and the mean mixed term is at most \(K(\Re a_{iT},\Re b_{iT})\). The latter equals \(\Re J(iT)+\Re K(c_{iT},c_{iT})+K(\Im a_{iT},\Im b_{iT})\). Since \(\mu(D)\le1/4\), the \(L^2(\mu)\) norms satisfy \[
\|c_{iT}\|_2^2\le\frac1{16T^2},\qquad
\|\Im a_{iT}\|_2\le\tfrac14,\qquad
\|\Im b_{iT}\|_2^2\le\tfrac43D_y .
\] The last bound uses \(y^2\le y\) and \((\Im Q_{iT})^2\le(4/3)(\alpha_T-\alpha_{2T})\). Equation (2) now bounds the mean residual above by \(S^y_T\).

On the checkerboard, the scalar-gap calculation's estimate \(E_T<1/T\) uses the interval Fourier bound independently of the cut positions. For the two axis indicators \(g_0,h_0\), one has \(0\le\Re V(iT)-\mu(h_0)\le E_T\) and \(0\le D_y\le E_T\). Also \[
|K(a_{iT},b_{iT})-K(g_0,h_0)|
\le\tfrac12(\|a_{iT}-g_0\|_2+\|b_{iT}-h_0\|_2)
\le\sqrt{E_T/2}.
\] Subtracting the cell's exact residual \(R=0\) from \(S^y_T\) gives \[
S^y_T\le\frac{2}{T}+\frac{8+8/\sqrt6}{\sqrt T}+\frac{\sqrt2}{T^2}
 <13\cdot10^{-6}<2^{-16}.
\]

Suppose a tree matched all \(3N\) readings to within \(h^4\). Let \(F_s(t)\) be its horizontal reading on the proportional pair. At a real node, conjugating \(\mathcal G_{A_s}\) by \(A_s^{1/2}\) gives an orthogonal projection with matrix symbol \((A_s^{1/2}n)(A_s^{1/2}n)^T/(n^TA_sn)\). The transformed mean field \(A_s^{1/2}(1,0)^T\) has norm one. Hence every Neumann coefficient of \((F_s(t)-1)/t\) has magnitude at most one, for trees and cells alike. The three-point coefficient estimate above gives discrepancies in \(f,L_e,B_e\) of at most \(9000h\). Since \(\det A_{s_j}\ge3/4\) and \(|L_e|\le1\), (5b) bounds discrepancies in \(V,J\) at the real nodes by \(36000h\).

On the interpolation circle used above, both \(V(s(z))\) and \(J(s(z))\) have magnitude at most \(10^9T^4\) by (2). With the present \(N,T\) the interpolation remainder at either auxiliary point is at most \(5\cdot10^9T^5 e^{-200}<7\cdot10^{-18}\). The Lagrange bound \((3N/b)^N\) therefore gives discrepancies below \(2\cdot10^{-17}\) in \(V(iT),V(2iT),J(iT)\). The inequality \(|\sqrt u-\sqrt v|\le\sqrt{|u-v|}\) bounds the resulting discrepancy in \(S^y_T\) by \(10^{-7}\). This contradicts \(13\cdot10^{-6}+10^{-7}<2^{-16}\). Continuity of the maximum over the fixed finite list gives the asserted error bound in the pointwise closure.

The question is posed as Problem 5 in the original preprint.
\phantomsection\label{endfile:five-band-proof}

\subsection{Uniform enclosure on the conductivity cone}\label{detail:finite-panel-errors}\label{file:finite-panel-errors}
For an affine loading $e$, a conforming periodic voltage fluctuation
gives an upper quadratic form $U$ for the cell response $C$.
Applying the same argument to the rotated reciprocal phases and using
planar duality gives a lower form $L$. Restrict each variational
problem to the span of the trial potentials in the certificate.
Assembling its Gram matrices and eliminating their potential block by
a Schur complement gives the two finite energy pencils. The assembly
is checked on the stored dyadic nodal data by integer dot products on
each triangle. The reference lower eigenvalue and the accumulated
coefficient radii bound the perturbation of the Schur complement.

For the other material, compute the quadratic energy of the specified
tree of depth eight. At a split of share $t$ and normal $n$, minimize
\[
 t(e+(1-t)c n)^T A(e+(1-t)c n)
 +(1-t)(e-t c n)^T B(e-t c n)
\]
over $c\in\mathbb R$. Working from the leaves to the root gives the
laminate response $T$. The rational shares and integer normals in the
input also determine the phase fraction exactly by the same tree
with arithmetic averaging of the leaf fractions.

It suffices to bound the largest eigenvalue of each signed pencil
difference $U-T$ and $T-L$. The normalization used in the calculation
is $\rho=533/1024$ and $c=(1-\rho)/(1+\rho)$. Expand each finite
pencil through degree 17 in the scalar and the two conjugate planar
anisotropy coordinates of each phase. Inserting the coefficients into
the finite linear equations bounds the coefficient errors using the
recorded positive reference eigenvalue. The Neumann tail, the model
error and the rounding of matrix operations enter separately. The
integer checker quantizes the signed polynomial on the four radial
sign choices. On each page it elevates the radial monomials to the
Bernstein basis. Its nonnegative basis weights sum to one, so it
suffices to bound every matrix coefficient as a function of the two
angles. A grid of 256 points per angle, an integer Fourier transform
with enclosed roots of unity, and second-derivative bounds on each
trigonometric mode control the intervals between grid points.
For a Hermitian matrix $B$, perturbing it by $E$ changes
$\lambda_{\max}(B)$ by at most $\|E\|$. The sum of the radial,
Fourier, interpolation, tail and model errors therefore bounds the
entire phase cone. The two degree-17 certificates record those
summands for each signed difference.

Finally laminate the tree with the small share of pure phase needed
to attain the prescribed fraction. The arithmetic and harmonic
bounds for a positive tensor between $0.1I$ and $I$ bound the resulting
operator error by nine times that share. Adding this exact-fraction
correction gives the second column of the preceding table.
\phantomsection\label{endfile:finite-panel-errors}

\SSection{cells}{Exact conductivity cells and finite readings}
All conductivities below are measured with the phase values $z$ and
$1$. The fraction $f$ and the binary phase indicator are fixed across
all contrasts in a single-cell assertion.

\subsection{Reciprocal equality and a finite-perimeter component}\label{route:reciprocal-equality-and-a-finite-perimeter-component}
We first derive the equalities for fields in an exact cell. We then
identify one component of positive volume which must obey the same
identity at every contrast. This is the step that makes the asymptotic
growth of the law in Theorem~\ref{new:perimeter} decisive.
\phantomsection\label{file:perimeter-fields}
\paragraph*{What reciprocal equality forces in a single cell}

\textbf{Finite-component obstruction.} Suppose a measurable periodic cell has effective tensor \(s(z)I\) for all \(z>0\), and \(\Phi_z(s)=0\) throughout the positive axis. If either phase is a finite union of connected open sets up to a null set, \(s\) is rational.

Let \(E\) be the matrix of the three electric fields as columns, with mean \(I\), and \(J=a_z E\) the currents, with mean \(sI\). Here \(a_z=1+(z-1)\chi\). The columns of \(E\) are constant vectors plus periodic gradients. The columns of \(J\) are divergence free. The fields are holomorphic as \(L^2\)-valued functions on the slit plane, by the cell equation and rotated coercivity of \(a_z\). Thus isotropy and reciprocal equality extend to the upper unit semicircle.

At \(z=e^{i\theta}\), \(0<\theta<\pi\), use the Hermitian quadratic form

\[
 B(X)=|\operatorname{tr}X|^2-\sum_{ij}X_{ij}\overline{X_{ji}}.
\]

Pointwise \(B(E)=B(J)\) since \(|a_z|=1\). A nonzero Fourier mode of \(E\) has the form \(n v^T\), with \(n\) a real unit wavevector, so \(B=0\) on it. For a nonzero mode \(X\) of \(J\), use orthonormal coordinates with \(n=e_3\). Its last row vanishes, and

\[
 \|X\|_F^2-B(X)=|X_{11}-X_{22}|^2+|X_{12}+X_{21}|^2
                       +|X_{13}|^2+|X_{23}|^2.
 \tag{2}
\]

The sum of these nonnegative defects over nonzero current modes is exactly \(3\Phi_z(s)\). Indeed Parseval gives \(\int B(J)=\int B(E)=6\). The cell energy identity gives \(\int a_z\|E\|_F^2=3s\). Taking real parts after rotation by \(e^{-i\theta/2}\) gives \(\int\|J\|_F^2=\int\|E\|_F^2=3(s+z\overline{s})/(1+z)\). The mean current contributes \(3|s|^2-6|s|^2\) to the full defect, and \(s(1/z)=\overline{s(z)}\) on this arc. This is the equality case of the archived Fourier energy inequality.

Every term in (2) therefore vanishes. In particular the axial vector of the skew part of each current mode is parallel to \(n\). The axial vector of the skew part of the corresponding electric mode is perpendicular to \(n\), since the mode has the form \(nv^T\). Both skew means vanish. Writing \(v_E\) for the axial vector field of the skew part of \(E\), the corresponding field for \(J\) is \(a_z v_E\). Orthogonality gives \(\int a_z|v_E|^2=0\). Rotated coercivity gives \(v_E=0\), so both matrices are symmetric.

Consequently the nonconstant modes of \(E\) are multiples of \(nn^T\), and those of \(J\) are multiples of \(I-nn^T\) by (2). In particular both \(E\) and

\[
 K=\tfrac12(\operatorname{tr}J)I-J
\]

are scalar constant matrices plus Hessians of periodic potentials. These are closed complex-linear constraints on the fields. Analytic continuation gives them at every positive contrast.

Write \(p=\operatorname{tr}E\). Let \(D\) act on scalar functions by the Fourier multiplier \(nn^T-I/3\) at nonzero modes and zero at the mean. Since \(E\) and \(K\) have the Hessian forms just proved, their deviatoric parts are \(Dp\) and \(D(a_zp)/2\), respectively. The definition of \(K\) also gives \(\operatorname{dev}K=-a_z Dp\). With \(h=z-1\), this yields

\[
 (3+2h\chi)\,Dp=-h\,D(\chi p).
 \tag{3}
\]

Now suppose the \(z\)-phase agrees up to a null set with an open set with finitely many connected components \(U_j\). In each of them, \(\operatorname{div}E=0\) column by column since \(a_z\) is constant there. For a Hessian plus a scalar mean, \(\operatorname{div}E=\nabla p\). Thus \(p\) is constant on each \(U_j\) and \(\chi p\) belongs, for every contrast, to the span of their indicators. \textbf{This is the only step that uses the hypothesis on the phase sets.} Equation (3) puts \(Dp\) in a fixed finite-dimensional space, since \((3+2h\chi)^{-1}=(1-\chi)/3+\chi/(2z+1)\). Also \(D^*D=\tfrac23(I-P_{\rm mean})\), so \(p\) and \(E=(p/3)I+Dp\) lie in fixed finite-dimensional spaces.

The corrector gradients for all contrasts and loadings now span a finite-dimensional subspace of periodic gradients. Restrict the cell equation to that subspace in a fixed basis. Its matrix and load are affine in \(z\), and its matrix is positive definite for \(z>0\). Solving this system and taking the mean current gives a rational tensor. If the unit phase satisfies the component hypothesis, interchange the phases and use \(z\,s(1/z)\), which also has zero reciprocal defect. This proves the obstruction and applies it to (1).

\textbf{A bound on the number of spectral atoms.} There is also a direct geometric formula for the parameter measure \(\pi\) in the coating representation. Regard \(P=\chi\) as the orthogonal projection by multiplication, on scalar or matrix fields as appropriate. On \(P L^2\) put

\[
 S=\tfrac32 P D^*P D P,\qquad
 w=\frac{9z}{(z+2)(2z+1)} .
\]

This is a positive contraction. Put \(u=Pp\). Applying \(D^*\) to (3) after division by \(3+2h\chi\), using \(p=3+\tfrac32D^*Dp\), and projecting by \(P\) gives

\[
 \big[I-(1-w)S\big]u=\frac{9+h\langle u\rangle}{z+2}\,\chi
                       =\frac{3(s+2)}{z+2}\,\chi .
\]

Here \(h\langle u\rangle=3(s-1)\) follows by tracing the mean current. For \(z\ne1\), taking the mean of the resolved equation gives

\[
 \left\langle\chi,\,[I-(1-w)S]^{-1}\chi\right\rangle_{L^2}
   =\frac{(s-1)(z+2)}{(s+2)(z-1)}
   =\int_{[0,1]}\frac{\lambda\,d\pi(\lambda)}
                     {1-(1-w)(1-\lambda)} .
\]

The last identity follows from (5). Expanding at \(w=1\) identifies the spectral measure of \(S\) on \(\chi\) as the image of \(\lambda\,d\pi(\lambda)\) under \(\lambda\mapsto1-\lambda\). If there are \(N\) open components in the \(z\)-phase up to null sets, all the displayed resolvent vectors lie in the span of their \(N\) indicators. Their Taylor coefficients \(S^j\chi\) lie in the same span. Thus \(\pi\) has at most \(N\) atoms in \((0,1]\) and no continuous part. The formula for \(S\) itself holds without the component hypothesis.

At finite perimeter this atom count also holds with measure-theoretic indecomposable components, without any open representative. Even with infinitely many such components, one component must satisfy (1) of the finite-perimeter report. Its volume bounds a pure outer shell share of the law's relaxed radial profile.

If \textbf{both} phases have connected open representatives, the only possibilities within the reciprocal-equality class are the two singly coated sphere laws. Indeed the one-component bound leaves \(\pi\) with just one positive atom and possibly an atom at zero. With no atom at zero, its mean fixes \(\pi=\delta_f\), giving \(K=w/(f+qw)\). With positive mass at zero, the interchanged function \(w/K\) vanishes at zero. The same argument applied to that phase gives \(w/K=w/(q+fw)\), hence \(K=q+fw\). The two coating rules identify these functions.

Both possibilities have cells meeting this connectedness restriction. On the cubic lattice take the periodic E-inclusion of Liu, James and Leo with shape matrix \(I/3\) and volume fraction \(f\). Their obstacle construction gives one convex inclusion inside the cube with analytic boundary, so its interior and the exterior in the torus are connected. If \(\Delta v=\chi-f\), its defining property is \(\nabla^2v=(q/3)I\) in the inclusion. The unit-load field matrix \(I-h\nabla^2v/(1+hq/3)\) then solves the cell equations and gives \(s(z)=1+fh/(1+hq/3)\). Interchanging the phases in the construction at fraction \(q\) gives the other endpoint.

\paragraph*{Finite perimeter without open phase representatives}

Let the indicator have finite perimeter on the torus. A finite-perimeter set is indecomposable if every partition into two sets of positive volume increases its total perimeter. The components here are counted on the torus. Their decomposition is at most countable and their perimeters add. The Euclidean variational argument uses the relative isoperimetric inequality, which on the flat torus bounds \(\min(|G|,1-|G|)^{2/3}\) by a constant times the perimeter. No open representative is required.

If a cell attains full reciprocal equality and the \(z\)-phase has \(N<\infty\) indecomposable components, its coating measure has at most \(N\) atoms in \((0,1]\) and no continuous part. The conductivity is rational. Finiteness for the other phase gives rationality by interchange. Thus any finite-perimeter cell for \(L_f\) must have infinitely many indecomposable components in \textbf{each} phase.

\paragraph*{A perimeter estimate for the field trace}

Write \(a=1+h\chi\), \(z=1+h\), \(E\) for the harmonic coordinate matrix of mean \(I\), and \(p=\operatorname{tr}E\). At any sufficiently small real \(|h|\) for which \(E\) is symmetric, \(p>0\) a.e. and

\[
           \log p\in BV(\mathbb T^3),\qquad
           |D\log p|\le C|h|\,|D\chi|.
 \tag{5}
\]

Here \(D\) in a variation measure denotes the distributional derivative, distinct from the deviatoric Fourier operator in (1)-(4). In the application to full equality \(E\) is symmetric throughout the contrast interval. For (5) symmetry at the contrast under consideration suffices.

For the proof, at the fixed contrast the ball means \(E_Q\) are all symmetric and nonsingular. Along concentric balls they vary continuously with the radius up to a ball covering the torus, where the mean is \(I\). Hence \(E_Q\) is positive definite, \(E\ge0\) by differentiation, and \(E\) is invertible almost everywhere by the small-contrast normalized inverse estimate. Those estimates also give \(\|E E_Q^{-1}-I\|_{2,Q}\le C|h|\) and \(\|E E_Q^{-1}\|_{3,Q}\le C\), in normalized norms. For a ball, use the cubic estimates on a containing cube of comparable volume, then transfer the mean normalization to the ball by the first estimate and inversion near \(I\). Consequently

\[
   \|E\|_{3,Q}\le C p_Q,\qquad p_{Q/2}\ge c p_Q,
   \qquad p_Q=\dashint_Q p .
\]

The second inequality follows by averaging the normalized field over the smaller ball and using \(\|E_Q\|\le p_Q\). The whole-cell estimates give \(p\) and \(1/p\) in \(L^2\).

Take a nonnegative smooth mollifier \(\rho_r\) supported on \(B_r\) and bounded below by \(c r^{-3}\) on \(B_{r/2}\), and write \(p_r=p*\rho_r\). Symmetry and the gradient columns give \(\operatorname{div}E=\nabla p\). Solenoidality of \(aE\) therefore gives \(a_0\nabla p=\operatorname{div}[(a_0-a)E]\) for any constant \(a_0\). At each center \(x\) use the majority-phase conductivity as that constant on \(Q=B_r(x)\) in estimating the convolution. If \(\delta_Q\) is the smaller of the two phase shares in \(Q\), Hölder and the relative isoperimetric inequality on balls give

\[
 \frac{|\nabla p_r(x)|}{p_r(x)}
 \le \frac{C|h|}{r}\frac{\delta_Q^{2/3}\|E\|_{3,Q}}{p_Q}
 \le \frac{C'|h|}{r^3}|D\chi|(B_r(x)).
\]

Here \(p_r(x)\ge c' p_Q\). Integrating against continuous nonnegative tests gives domination by \(C|h|\,|D\chi|\) in the limit, since the normalized ball kernels approximate the identity. Also \(\log p_r\to\log p\) in \(L^1\) by differentiation and uniform integrability, using \(1/p_r\le(1/p)*\rho_r\). This proves (5).

For any indecomposable component \(G\) of a phase, its density-one part \(G^1\) is contained in the density-one part of that phase. Thus \(|D\log p|(G^1)=0\). By BV coarea almost every level set \(A=\{\log p>t\}\) has \(|D1_A|(G^1)=0\). Splitting \(G\) by \(A\) adds twice the perimeter of \(A\) on \(G^1\) to the sum of the two resulting perimeters relative to the original one. Indecomposability therefore makes \(\log p\), and hence \(p\), constant almost everywhere on \(G\). This is the usual constancy criterion on an indecomposable set.

Being constant on each member of the fixed component partition is a closed linear condition in \(L^2\). The \(L^2\)-valued field is holomorphic in contrast, so at full equality the component conclusion extends from the small real interval to the entire positive axis. If there are \(N\) occupied-phase components then \(\chi p\) lies in their \(N\)-dimensional span. The resolvent formula for the coating measure now gives the claimed atom bound exactly as for open components. If both phases are indecomposable, only the two Hashin-Shtrikman laws occur at full equality.

\paragraph*{Upper variance equality on an indecomposable phase}

There is already a conclusion from the fourth-order equality in (2), without assuming full reciprocal equality. Write \(\operatorname{dev}D^2\Delta^{-1}g=\chi\,\operatorname{dev}B\) and \(\langle g\rangle=0\). The right side is bounded by the upper-variance Hessian estimate. Taking divergence gives

\[
             \nabla g=\tfrac32\operatorname{div}(\chi\,\operatorname{dev}B),
             \qquad |Dg|\le C |D\chi|.
\]

Indeed \(\operatorname{div}(\operatorname{dev}B)=\tfrac23D\chi\). Replace \(\chi\) in its product with the bounded tensor \(\operatorname{dev}B\) by smooth convolutions. The product rule gives one term dominated by \(C|D\chi|\) and a second by \(C(|D\chi|*\rho_r)\,dx\). In the distributional limit both bounds yield the asserted measure domination, without multiplying two unspecified boundary representatives. Thus \(g\) is constant on each indecomposable component by the same coarea argument. If the occupied phase is indecomposable this gives \(S\chi=q\theta\chi\) in the notation of (1) and (4). The rank-one cyclic-space proof above yields full equality and the stated single-positive-atom formula. If both phases are indecomposable only the two endpoints are possible even at this fourth-order equality.
\phantomsection\label{endfile:perimeter-fields}

\phantomsection\label{file:perimeter-component}
\phantomsection\label{detail:common-component}
We use the decomposition of a set of finite perimeter into positive-volume
indecomposable components \cite[Sec.~4, Thm.~1]{acmm}. The two phases
have finite perimeter in the periodic cell. All cuts below are regional
BV cuts, the faces of the period cube are identified, and the relative
isoperimetric inequality on the flat torus gives the same component
decomposition there. On any such component, a BV scalar with derivative
zero in its measure-theoretic interior is constant almost everywhere
by coarea and indecomposability \cite[Sec.~4, Prop.~4 and its consequence]{acmm}.
\paragraph*{Fields on individual components}

Write \(U\) for the vector of coordinate potentials with gradient matrix \(E\) and \(V\) for the vector with gradient matrix \(K=(\operatorname{tr}J)I/2-J\), where \(J=aE\) is the current. Reciprocity makes both matrices symmetric Hessians plus scalar means, with \(\langle E\rangle=I\), \(\langle K\rangle=sI/2\). Hence on the periodic lift

\[
 U(x+k)=U(x)+k,\qquad V(x+k)=V(x)+(s/2)k,\qquad k\in\mathbb Z^3 .
\]

The finite-perimeter trace estimate shows that \(p=\operatorname{tr}E\) is constant on each measure-theoretic indecomposable component at every contrast. Denote the value on component \(i\) by \(p_i(z)\) and put \(\alpha_i(z)=a_i(z)p_i(z)/2\). The functions \(p_i\) are analytic since they are component averages of the analytic fields. On that component \(D(V+a_iU)=\alpha_i I\). We use a small positive \(h=z-1\) below. Then \(E\ge0\) and the fields \(E,K\) are locally in \(L^p\) for some \(p>3\), so \(U,V\) have continuous representatives.

Assume that no component satisfies the first identity in (1) identically. Since there are at most countably many components, choose a single \(z\) arbitrarily close to and greater than one at which

\[
                     b_i=\alpha_i-a_i-s/2\ne0\quad\text{for every }i .
 \tag{2}
\]

All the following fields and constants refer to this contrast. The aim is to contradict (2).

\paragraph*{A partition by bounded lifted pieces}

For component \(i\) the continuous vector map

\[
             W_i=\big[V-(s/2)x+a_i(U-x)\big]/b_i
\]

is periodic and \(DW_i=I\) almost everywhere on that component. The torus-valued map \(W_i(x)-x\bmod\mathbb Z^3\) is constant there. Indeed the real and imaginary parts of each coordinate exponential have zero Sobolev gradient on the component, and the indecomposable-set constancy criterion applies. By continuity the equality modulo the lattice holds on its essential closed support. Choose a constant \(c_i\) representing the difference. The bounded map \(W_i-c_i\) now selects one lift of every point of that support. In a sufficiently small torus chart the continuous difference \(W_i-c_i-x_{\rm chart}\) has oscillation less than one. On the support it takes lattice values, hence just one value on the chart. Above this chart the chosen lift is a translate of the whole component portion on just one sheet and is absent on all other sheets. Summing the variations over sheets and using a finite partition of unity gives equality of its perimeter and the torus component's perimeter. If the lift split into two pieces with additive perimeters, projection would give disjoint pieces of equal respective volumes downstairs, each with perimeter no larger than its upstairs perimeter. This follows by differentiating the sum of translated indicators over sheets and taking total variations. The two projected perimeters would add without excess, contradicting indecomposability downstairs.

Use the selected pieces and all lattice translates to partition \(\mathbb R^3\) in measure. Each has positive volume and is bounded and indecomposable. Their perimeters are locally summable because in each chart the sum over the lifts is the periodic copy of the original sum of component perimeters. We call two pieces adjacent when their reduced boundaries have common positive area. At perimeter-almost every point of each piece there is exactly one neighbor on the other side. To see this, \(D1_G=-\sum_{G'\ne G}D1_{G'}\) locally as measures by local perimeter summability. The reduced half-space limits and disjoint interiors in measure allow only one other piece at such a point. Same-phase pieces share no perimeter since their original component perimeters add without cancellation. These arguments also apply when component supports accumulate in a compact region. The finite-perimeter decomposition and constancy properties used here hold on the flat torus with its relative isoperimetric inequality.

The adjacency graph of the lifted pieces is connected. In a disconnected graph the union of all pieces in one graph component would have zero distributional perimeter in \(\mathbb R^3\), since the opposing normals cancel on every shared interface in that union and all interfaces are accounted for. Its volume and that of its complement are both positive, a contradiction. Also, a finite chain of successively adjacent pieces has indecomposable union. Indeed a perimeter-additive splitting of the union introduces no interface area in its density-one part. It therefore cannot split any one of its indecomposable pieces, and assigning two adjacent pieces to opposite sides would introduce their positive shared area into that part.

On a bounded lifted piece \(G_i\) there is a single constant vector \(d_i\) with \(V+a_iU=\alpha_i x+d_i\) there. This follows by the same Sobolev constancy criterion, now with ordinary Euclidean coordinates. Continuity gives the identity on its entire closed support. On any edge \(ij\), subtracting the two identities gives single scalar affine traces \(U=t_{ij}x+u_{ij}\) and \(V=v_{ij}x+w_{ij}\) on the shared support. Here a scalar affine map has linear part a scalar times \(I\). Its coefficients are unique on a positive-area set.

\paragraph*{Why the affine traces agree around a cycle}

Fix a lifted piece \(G_i\) and one coordinate potential \(U_k\). With distributional curls in \(\mathbb R^3\),

\[
 \operatorname{curl}(1_{G_i}\nabla U_k)=\mu_i\times e_k,\qquad
 \mu_i=t\, D1_{G_i},\qquad t=t_{ij}\ \text{on interface }ij .
 \tag{3}
\]

Here are details allowing arbitrary finite-perimeter interfaces. On \(G_i\), \(|E|\le p_i\) since \(E\ge0\). On a small ball \(Q=B_r(x)\) one also has \(\|E\|_{3,Q}\le C p_Q\) and \(\|E E_Q^{-1}-I\|_{2,Q}\le C|h|\). If \(G_i\) occupies at least half of \(Q\), the second estimate and \(p=p_i\) on \(G_i\) give \(p_Q\le C p_i\). To bound the mollified left side of (3), use \(\operatorname{curl}\nabla U_k=0\) to subtract the majority label from \(1_{G_i}\) in the convolution. Use \(|E|\le p_i\) when \(G_i\) is the minority and the ball \(L^3\) estimate when its complement is the minority. By the relative isoperimetric inequality the resulting bound is \(C p_i r^{-3}|D1_{G_i}|(B_r(x))\). Passing to the distributional limit proves that the curl is a measure dominated by \(C p_i|D1_{G_i}|\). For \(V\) use also \(|K|\le C_z |E|\).

At a typical reduced-boundary point shared with \(G_j\), the two blown-up indicators converge locally in measure to opposite half-spaces. On each fixed rescaled ball the potentials, after subtracting their values at the center and dividing by the scale, have bounded \(W^{1,p}\) seminorm for some \(p>3\). Use the local normalization estimate with \(p=4\) and \(|h|\) below a fixed absolute threshold. The cube estimate transfers to balls as for \(p=3\). The limiting \(G_i\) volume share on centered balls is \(1/2\), so \(p_Q\) is bounded there. Subsequential limits are continuous, with local uniform convergence. The two component identities pass to their respective half-spaces, hence \(U\) in the limit has planar trace \(t_{ij}x\) after centering. The curl of its gradient cut off to the half-space is the normal boundary measure crossed with \(t_{ij}e_k\). The rescaled products inside the curl converge weakly by the strong convergence of the indicators. Since the original curl is already dominated by perimeter, differentiation of measures on the reduced boundary identifies its density there with this tangent expression. This proves (3), including its sign with the convention that \(D1_{G_i}\) is the distributional gradient. The three curl bounds give \(|t|\le C p_i\) perimeter-almost everywhere because \(\sum_k|n\times e_k|^2=2\) for a unit normal. In particular \(\mu_i\) is finite.

Taking divergences in (3) for the three coordinates gives \(\operatorname{curl}\mu_i=0\). On \(\mathbb R^3\) take a scalar primitive \(g_i\in BV_{\rm loc}\) with \(Dg_i=\mu_i\), by the distributional Poincaré lemma. It has a constant value \(g_{\rm in}\) on \(G_i\). On a finite connected chain \(H\) of adjacent other pieces it likewise has a single value \(g_{\rm out}\) because \(H\) is indecomposable and \(|Dg_i|\) vanishes on its density-one part. By (3) the columns of \(1_{G_i}E-g_i I\) are curl free, giving a vector potential \(Z_i\in W^{1,1}_{\rm loc}\). Component constancy gives \(Z_i=U-g_{\rm in}x+c\) on \(G_i\) and \(Z_i=-g_{\rm out}x+c'\) on \(H\). A Sobolev \(W^{1,1}\) function is approximately continuous at \(\mathcal H^2\)-almost every point. At a typical shared boundary point both pieces have positive density, so these two continuous expressions give its same approximate limit. Thus the common trace of \(U\) on \textbf{all} such interfaces is \((g_{\rm in}-g_{\rm out})x+c'-c\). Using \(K\) in place of \(E\) gives the same conclusion about the full scalar affine expression for \(V\).

In a simple graph cycle, delete a vertex \(i\) and take the remaining chain as \(H\). Both incident cycle edges at \(i\) therefore have identical affine data for each potential. Moving around the cycle makes its two affine maps common to all edges. The same data are consequently shared by any edges related by a finite chain of simple cycles with successive cycles sharing an edge.

\paragraph*{Periodicity}

Take the equivalence classes of edges generated by the preceding cycle relation, with each bridge a singleton. The incidence graph between the original vertices and the edge classes is a tree. For clarity, each class is connected, and two distinct classes cannot share two vertices, since paths between them inside the two classes would contain a simple cycle using edges from both. Any simple cycle in the incidence graph would likewise connect distinct edge classes into one. For a shortest such cycle, paths joining the consecutive incidence vertices within their respective classes have disjoint interiors from nonconsecutive paths, or one obtains a shorter incidence cycle. Neighboring classes share only their incidence vertex. These paths would give a simple cycle in the original graph across different classes. Connectedness is inherited from the original graph. Only finite cycles and paths enter this argument.

The period lattice acts on this tree preserving node types, and fixes no original vertex by a nonzero translation since all lifted pieces are bounded. Some edge class must be fixed by a nonzero period \(k\). Otherwise the action would be free on tree vertices without edge inversion. A nontrivial translation in such an action has an invariant bi-infinite axis and acts on it by a nonzero shift (take a vertex minimizing its displacement and concatenate the translates of the joining path). Every commuting translation preserves the same axis and acts by a shift, since a reversal cannot commute with a nonzero shift. The resulting homomorphism \(\mathbb Z^3\to\mathbb Z\) has a nontrivial kernel, fixing vertices of the axis, a contradiction.

Write \(t x+d\) and \(v x+e\) for the common \(U,V\) traces on the fixed edge class. Translating an edge by \(k\) keeps it in the class while adding \(k\) to \(U\) and \((s/2)k\) to \(V\). Crucially both offsets \(d,e\) stay the same, hence \((1-t)k=0\) and \((s/2-v)k=0\). The fixed class cannot be a singleton bridge since the translation fixes no incident piece and preserves the two phase labels. On a positive-area interface of any incident component the relation \(V+a_i U=\alpha_i x+d_i\) now forces \(\alpha_i=s/2+a_i\). This contradicts the choice (2). At least one torus component therefore has this analytic identity at every positive contrast.

Finally \(\int a |E|^2=3s\) and \(|E|^2\ge |p|^2/3\) give (1) on that component. Neither a smoothness bound on the individual surfaces nor a finite number of components was used. Finite perimeter was used to obtain the component traces, the locally summable interface partition, and the interface measures (3).

\paragraph*{The volume of the component and the coating measure}

One can also locate a definite pure outer shell in the unique relaxed radial profile of this scalar law. If the components satisfying (1) have material \(z\) and total volume \(m_*>0\), their shell share obeys

\[
   \pi(\{0\})=0,\qquad
   \left(\int_{(0,1]}\frac{d\pi(\lambda)}{\lambda}\right)^{-1}
        =\lim_{w\downarrow0}\frac{w}{K(w)}
        \ \ge\ m_* .
 \tag{4}
\]

For unit-phase components the corresponding bound is \(K(0+)=\pi(\{0\})\ge m_*\). The two types cannot both occur by (1). The share is a statement about the radial transfer representation of the law and does not impose spherical geometry on the actual component.

Indeed if \(A\) is the union of the indicated \(z\)-components, the resolvent \(R_\delta\chi=(I-\delta S)^{-1}\chi\) from the coating formula has value \(K(w)/w\) on \(A\), with \(\delta=1-w\). This follows by substituting \(p=2+s/z\) there. Both sides of the paired identity are analytic for \(|\delta|<1\). Expanding the pairing with \(1_A\) in \(\delta\) gives for every polynomial \(P\)

\[
  \langle 1_A,P(S)\chi\rangle
        =m_*\int P(1-\lambda)\,d\pi(\lambda),\qquad
  \|P(S)\chi\|_2^2=\int |P(1-\lambda)|^2\lambda\,d\pi(\lambda).
\]

By Cauchy-Schwarz the first integral defines a bounded functional in \(L^2(\lambda\,d\pi)\) of norm at most \(m_*^{-1/2}\) after dividing by \(m_*\). Polynomial density and the Riesz representation theorem identify \(d\pi=(1/\lambda)\lambda\,d\pi\) with no atom at zero and give \(\int d\pi/\lambda\le m_*^{-1}\). Monotone convergence gives (4). Interchanging phases sends \(K\) to \(w/K\), proving the unit-phase formula. In the coating recursion these limits are the maximal initial pure shell shares, allowing a relaxed inner profile.
\phantomsection\label{endfile:perimeter-component}

\subsection{Measurable layers in smooth coordinates}\label{route:measurable-layers-in-smooth-coordinates}
Let $S$ denote the symmetric matrix of electric fields at one positive
contrast. The following argument applies to two-valued coefficients
locally of the form $\alpha(t)$, with $\alpha$ merely measurable and
$t$ a regular coordinate of class $C^{1,1}$. When a closed set of zero
$W^{1,2}$ capacity is removed, the proof assumes that the field is
bounded. The boundedness required for the cylinder-pole exclusion is
obtained from the single-pole field identities.
\phantomsection\label{file:measurable-layers}
\paragraph*{Regularity in a chart}

Use \(n=\nabla t/|\nabla t|\) and \(P=I-nn^T\). On smaller charts, \(P S\) and \(a S n\) belong to \(H^1\). To verify this, take one coordinate potential \(u\) with \(\nabla u\) a column of \(S\). In chart coordinates its equation is \(\operatorname{div}(\alpha(t)C(y,t)\nabla u)=0\) with \(C\) smooth and locally uniformly elliptic. Tangential difference quotients satisfy the same elliptic equation with divergence of an \(L^2\) source uniformly bounded by \(C_0\|\nabla u\|_2\) on compact subcharts, because the differences of \(C\) divided by the step are bounded. Testing with a cutoff times the difference quotient gives \(\partial_{y_i}u\in H^1_{\rm loc}\). The transverse flux \(\alpha(t)(C\nabla u)_t\) has tangential derivatives in \(L^2\) by this estimate. Its depth derivative is the negative sum of the tangential derivatives of the other two flux components. Thus it too is in \(H^1_{\rm loc}\). It differs from \(a n\cdot\nabla u\) by a smooth positive factor. The tangential derivatives determine \(P\nabla u\) by smooth coefficients, which proves both assertions using the symmetry of \(S\).

Write \(s=n^TSn\), \(b=PSn\), \(F=as\), and \(T=S-snn^T\). We have \(F,b,ab,T\in H^1_{\rm loc}\). In particular \(T\) is the sum of the smooth-frame tangential block and the two shear blocks built from \(b\), and \(Tn=b\).

\paragraph*{A symmetric approximation}

Take the same smooth nonnegative ambient mollification \(a_\epsilon\) of \(a\) on the torus in every chart, with kernel supported in the radius-\(\epsilon\) ball. Define the local symmetric \(H^1\) field

\[
S_\epsilon=T+(F/a_\epsilon)nn^T .
\]

It tends to \(S\) in \(L^2\) on compact subcharts. If \(X\) is a smooth tangent vector field to the leaves, \(X\cdot\nabla a=0\) in distributions. The tangential derivative \(X\cdot\nabla a_\epsilon\) is uniformly bounded on compact subcharts and tends to zero in \(L^2\). Indeed it equals the commutator of \(X\cdot\nabla\) with convolution applied to \(a\). Its kernel terms are bounded by \(C|x-y|\,|\nabla\rho_\epsilon(x-y)|\) and \(C\rho_\epsilon(x-y)\). These yield the uniform \(L^p\) bounds, including \(p=\infty\), and smooth approximation of \(a\) in \(L^2\) gives convergence there. Local cutoffs suffice.

The curl of each column of \(S\) is zero. The difference \(S_\epsilon-S\) has the form \(F(1/a_\epsilon-1/a)nn^T\). In taking its column curls only tangential derivatives of the scalar prefactor occur, since the columns are parallel to \(\nabla t\). This assertion for the nonsmooth term follows directly in leaf coordinates by taking the exterior derivative of a multiple of \(dt\). The terms containing \(1/a_\epsilon-1/a\) tend to zero in \(L^2\) by dominated convergence. The remaining tangential derivatives of \(a_\epsilon\) are uniformly bounded and tend to zero in measure, which also gives convergence when multiplied by the fixed \(L^2\) field \(F\). Derivatives of smooth frame factors cause only terms with the vanishing prefactor. Hence \(\operatorname{curl}S_\epsilon\to0\) locally in \(L^2\).

For the divergence, \(a_\epsilon S_\epsilon-aS=(a_\epsilon-a)T\). One has \(\operatorname{div}(aT)=a\operatorname{div}T\) distributionally. In chart coordinates the tangential derivatives commute with \(a\). For the normal flux term, \(b\) and \(ab\) belong to \(H^1\) and \(\partial_t(ab)=a\partial_t b\) almost everywhere. This follows by locality of Sobolev derivatives on each of the two measurable phases. Smooth Jacobian and frame factors obey the usual product rule. Consequently

\[
\operatorname{div}(a_\epsilon S_\epsilon)
=(a_\epsilon-a)\operatorname{div}T+T\nabla a_\epsilon .
\]

Only the normal contribution \(b\,(n\cdot\nabla a_\epsilon)\) needs another estimate. Let \(Z\) be the relatively closed set of depths where \(\alpha\) is not essentially constant on any neighborhood. On almost every fixed-\(y\) line, \(b\) and \(ab\) have continuous representatives in \(t\). If \(b\) has a nonzero component at a depth on that line, the ratio of the corresponding components of \(ab\) and \(b\) is continuous nearby and equals \(\alpha\) almost everywhere there. It must be constant on that interval. Thus \(b=0\) on \(Z\) on every such line, and \(\partial_t b=0\) almost everywhere on \(Z\) by locality.

On a compact subchart \(\nabla a_\epsilon=0\) when \(\operatorname{dist}(t,Z)>C\epsilon\), because the ambient averaging ball then meets only one material. Its norm is at most \(C/\epsilon\) everywhere. For the remaining depths, absolute continuity along the line gives

\[
|b(y,t)|^2\le C\epsilon\int_{|u-t|\le C\epsilon}|\partial_t b(y,u)|^2\,du .
\]

Integration first in \(t\) bounds the squared norm of the normal contribution by \(C\) times the integral of \(|\partial_t b|^2\) where \(\operatorname{dist}(t,Z)\le 2C\epsilon\) in a slightly larger chart. That integral tends to zero. This finishes the local \(L^2\) divergence estimate even for a depth profile of infinite perimeter.

\paragraph*{The sign}

On any compact subset away from \(K\) glue these local approximations with a smooth partition of unity on a neighborhood. They all use the same \(a_\epsilon\). Derivatives of partition functions multiply differences tending to zero in \(L^2\), by convergence to the same \(S\) and \(aS\). Thus both residuals tend to zero there. If \(S\) is bounded the approximants are uniformly bounded by a constant depending only on \(z\) and \(\|S\|_\infty\), since their pointwise construction uses orthogonal projections and the ratio \(a/a_\epsilon\).

For a smooth positive \(c\) and a symmetric \(H^1\) matrix \(D\), set

\[
\mathcal E_c(D)=\int c|\operatorname{curl}D|^2+c^{-1}|\operatorname{div}(cD)|^2 .
\]

The inequality \(\mathcal E_c((-D)_+)\le\mathcal E_c(D)\) follows from weighted integration by parts. In the expression for \(\mathcal E_c(D)\) as \(\int c|\nabla D|^2-c\operatorname{tr}(D^2\nabla^2\log c)\), replacing \(D\) by its matrix absolute value preserves \(D^2\) and decreases the derivative term by Frobenius-norm contraction. Apply the triangle inequality to \((-D)_+=(|D|-D)/2\).

If \(K\) is empty the global glued fields give \(\operatorname{curl}((-S)_+)=0\) by this inequality and \(L^2\) convergence. For the capacity-zero case first apply it to \(\eta S_\epsilon\) with \(0\le\eta\le1\) compactly supported away from \(K\). Local convergence of the residuals and uniform field bounds give \(\|\operatorname{curl}(\eta(-S)_+)\|_2\le C\|\nabla\eta\|_2\). Let \(\eta\) tend to one in \(H^1\) using capacity zero. This gives the same curl identity distributionally.

Put \(D_-=(-S)_+\). Its columns minus their means are periodic gradients. Pair them with the divergence-free columns of \(aS\). Since \(SD_-=-D_-^2\) and the mean current \(A_\chi(z)\) is positive definite, this yields

\[
-\int a|D_-|^2=A_\chi(z):\langle D_-\rangle\ge0 .
\]

Therefore \(S\ge0\). At a single pole put \(v=\Delta^{-1}(\chi-f)\) and \(B=D^2v\). Then \(E_z=I-[(z-1)/(1+\tau(z-1))]B\) is bounded and symmetric at every positive contrast. Sending \(z\) to infinity and to zero through countable sequences gives \((\tau-1)I\le B\le\tau I\) almost everywhere. The archived measurable cylinder obstruction then excludes the scalar target at \(\tau=(1-f)/2\) and, by phase interchange, at \(\tau=1-f/2\).
\phantomsection\label{endfile:measurable-layers}

\subsection{Boundedness and the excluded cylinder pole}\label{detail:rough-pole-proof}\label{file:rough-pole-proof}
For the mask $\chi$ write $D=\{\chi=1\}$, $q=1-f$ and
$B=D^2\Delta^{-1}(\chi-f)=\Gamma(\chi I)$, where $\Gamma$ projects
onto mean-zero gradient fields. The trace spectral measure is the sum
of three nonnegative diagonal measures. A sole pole $\tau$ therefore
forces $\Gamma(\chi B)=\tau B$ column by column by the spectral theorem.
Parseval gives $\|B\|_2^2=fq$ and, by testing the pole equation with
$B$, $\int_D|B|^2=\tau fq$. The electric field is
$E_z=I-(z-1)B/[1+\tau(z-1)]$, as follows by substitution in the
projected cell equation. These identities also justify the boundedness
assumption in the measurable-layer argument above, as follows.
\[B=\Gamma(\chi I),\qquad\Gamma(\chi B)=\tau B,\qquad \langle B\rangle=0.\tag{4}\]

Only these periodic field identities and \(0\le\chi\le1\) enter the estimate.

Fix a ball \(Q_0\) of positive radius. Let \(Q_j\) be its concentric dilates by \(2^j\), up to the first index \(N\) covering the torus. Write \(F_Q=\dashint_Q F\) and \(\|F\|_{2,Q}=(\dashint_Q|F|^2)^{1/2}\). For a vector field put

\[
 X(F)=\max_{0\le j\le N}2^{-j/2}\,\|F-F_{Q_j}\|_{2,Q_j}.
\]

This seminorm ignores additive constants and vanishes only for constant fields since \(Q_N\) is the whole torus. When \(F_{Q_0}=0\), bounded ratios of the volumes of successive balls give

\[
 |F_{Q_j}-F_{Q_{j-1}}|\le C 2^{j/2}X(F),\qquad
 \|F\|_{2,Q_j}\le C 2^{j/2} X(F)
 \tag{5}
\]

by summing a geometric series.

For every scalar multiplier \(m\) with \(|m|\le1\) and every \(F\) with \(F_{Q_0}=0\), one has

\[
 X(\Gamma(mF))\le C X(F)
 \tag{6}
\]

with the same scale-independent choice of \(C\) for all \(Q_0\). Here are the estimates for a small ball \(Q_j\) of radius \(r_j\). Input on \(Q_{j+1}\) contributes at most \(C\|F\|_{2,Q_{j+1}}\) to the root-mean-square on \(Q_j\) by the \(L^2\) projection bound. For input on the farther shell \(Q_{j+k+1}\setminus Q_{j+k}\), \(k\ge1\), subtract the output evaluated at the center. The gradient of the projection kernel is bounded by \(C|x|^{-4}\) locally off zero. Integration over the shell bounds the output variation on \(Q_j\) by \(C2^{-k}\|F\|_{2,Q_{j+k+1}}\). By (5) these estimates sum to at most \(C2^{j/2}X(F)\). Sources a fixed positive distance away have oscillation at most \(C r_j\|F\|_2\), which obeys the same bound since \(r_j\le C2^{j-N}\). For balls of radius bounded below the global contraction suffices, as \(N-j\) is then bounded. Subtracting each output mean costs at most a factor of two. Multiplying by \(2^{-j/2}\) proves (6).

Now take an eigenpair \(B_{Q_0}e=\lambda e\) of the symmetric mean matrix, with \(|e|=1\). The vector field \(F=Be-\lambda e\) has zero mean on \(Q_0\). By (4),

\[
 \Gamma(\chi F)=(\tau-\lambda)Be.
\]

Consequently \(|\tau-\lambda|X(F)\le C X(F)\) by (6), because \(Be\) differs from \(F\) by a constant. If \(X(F)>0\) this bounds \(|\lambda|\) by \(1+C\). Otherwise \(F=0\) on the torus and \(\langle Be\rangle=0\) gives \(\lambda=0\). The estimate applies to every eigenvalue of every ball average, independently of its center or radius. Lebesgue differentiation yields the asserted essential bound for \(D^2v\), and hence Lipschitz continuity of \(\nabla v\).

One consequence of (6) holds even without a single pole. There is a universal \(h_0>0\) such that for every measurable cell and every ball \(Q_0\), the mean matrix of the three corrector fields on \(Q_0\) at contrast \(1+h\) has positive determinant for \(|h|<h_0\). If it were singular for some ball and some such real \(h\), a nonzero load \(e\) would give a field \(F=e-h\Gamma(\chi F)\) with \(F_{Q_0}=0\). Equation (6) gives \(X(F)\le |h|C X(F)\). For \(|h|C<1\) this forces \(F=0\), contradicting \(\langle F\rangle=e\). For each fixed ball the determinant is continuous in \(h\) and is one at zero. This proves its sign, and differentiation of ball averages gives nonnegative pointwise determinant almost everywhere at each of these contrasts.

\textbf{Positivity for a measurable cell with symmetric fields near unit contrast.} If the field matrices are symmetric on a real interval containing \(1\), then \(E_{1+h}\ge0\) almost everywhere throughout the part of that interval with \(|h|<h_0\). For each fixed ball the matrix mean is symmetric, continuous in \(h\), nonsingular by the preceding estimate and equal to \(I\) at \(h=0\). Its eigenvalues stay positive. Lebesgue differentiation yields the assertion. The normalized inverse estimate below also gives positive definiteness almost everywhere after reducing the universal interval. No phase perimeter is used.
\paragraph*{The cylinder restriction.}
The positivity result for measurable layers gives
$(\tau-1)I\preceq B\preceq\tau I$ by taking countable sequences
$z\to0$ and $z\to\infty$ in $E_z\succeq0$. Put $v=\Delta^{-1}(\chi-f)$
and $F=\nabla v$. The preceding bound makes $F$ Lipschitz.
\textbf{At \(\tau=q/2\), the upper bound \(B\le\tau I\) on \(D\) alone forces a singular high-contrast endpoint tensor.} In particular a scalar cell at this pole must have an eigenvalue of \(B\) strictly larger than \(\tau\) on a set of positive volume in \(D\). Phase interchange gives the corresponding failure of the lower bound on \(D^c\) for a scalar cell at \(1-f/2\).

To prove this let \(N=\tau I-B\). At this pole the trace of \(N\) on \(D\) is \(\tau\) and \(\dashint_D|N|^2=\tau^2\) by the moment identity above. Under \(N\ge0\) it follows that \(N=\tau P\) and \(B=\tau(I-P)\) on \(D\), with \(P\) a rank-one projection. This time \(J=\chi P\) has divergence-free columns and \(JB=0\). Indeed, the pole identity gives $\operatorname{div}(\chi P)=0$. Since $PB=0$, the chain rule gives the same identity for $\eta(F)\chi P$. Testing its columns against the periodic potential $v$ gives $\int F^T\chi P F=0$. Hence \(\eta(F)J\) is divergence free and \(JF=0\).

We first show that \(P\) on \(D\) is a measurable function of \(F\). The area formula for the Lipschitz map \(F\) puts \(F(D)\), after discarding a null subset of \(D\), inside a set \(Z\) of zero three-dimensional volume, since \(DF\) has rank two there. Write \(\nu=F_*(dx)\). For each smooth scalar \(w\) on the torus, the matrix measure \(F_*(wB\,dx)\) has finite column divergence on target space. Indeed for each column \(j\) and smooth target test \(\phi\), its pairing with \(\nabla\phi\) is \(\int w\,\partial_j(\phi\circ F)=-\int(\partial_j w)\,\phi(F)\). The singular part of a three-dimensional matrix measure with finite column divergence has polar matrix of rank at most two \cite[Cor.~1.13]{dpr}. Here the wave cone of the column divergence consists of matrices with a common nonzero left null vector. Hence the density of \(F_*(wB\,dx)\) relative to \(\nu\) has rank at most two \(\nu\)-almost everywhere on \(Z\), including where that density vanishes.

Approximate \(w=\chi\) in \(L^1(dx)\) by smooth functions. The target densities relative to \(\nu\) converge in \(L^1(\nu)\) since \(B\) is bounded. Rank at most two is a closed condition. It follows that the conditional average of \(B\) with respect to \(\chi\,dx\) on almost every \(F\)-fiber has rank at most two. More explicitly, the limiting density relative to \(\nu\) is that conditional average times \(dF_*(\chi\,dx)/d\nu\), which is strictly positive at \(F_*(\chi\,dx)\)-almost every target. On \(D\) all the matrices being averaged are positive-semidefinite rank-two projections times \(\tau\). A vector in the kernel of their average must belong to the kernel of almost every summand. The average therefore determines their common kernel line \(l(y)\) at almost every target \(y\) for this phase. Its orthogonal projector is \(P\) on the fiber.

Test the solenoidal identity for \(\eta(F)J\) and disintegrate over \(y\). The conditional probability is invariant along \(l(y)\) and is supported on the closed level set \(F=y\). If the line is not parallel to an integer vector, its closure in the torus is a subtorus of dimension at least two. The support of its invariant measure contains that whole coset through each of its points. Since \(F\) is constant there, \(DF\) annihilates the tangent space of the coset wherever it is differentiable. This contradicts rank two at almost every point of \(D\). Consequently the axes are rational, giving countably many positive-volume classes \(D_l\) defined up to null sets. Each is invariant along its assigned line. Write \(S_l=\operatorname{supp}(1_{D_l}\,dx)\). On \emph{every} complete \(l\)-line through \(S_l\), \(F\) is constant. For almost every point this follows from the conditional invariance just proved, then continuity of \(F\) gives it throughout \(S_l\). Also \(F\cdot l=0\) there by \(JF=0\).

We show that the axes span at most a plane. If infinitely many axes occur, choose one complete line from each of a sequence of distinct classes. By compactness a subsequence of the base points converges, and its compact circle subgroups converge in Hausdorff distance to a closed connected subgroup. A one-dimensional limit is impossible. The two independent integer characters annihilating that candidate limit would be uniformly close to one on each late circle. The image of a connected subgroup under a character is either the full unit circle or the identity, so both characters would vanish on the late circles, making their directions equal. The limiting coset thus has dimension at least two, and \(F\) is constant on it by continuity. Dimension three is impossible since \(DF\) has rank two on \(D\). Every closed line transverse to this two-torus intersects it. A transverse class would have \(F\) equal to this same constant on all its lines, contradicting rank two on a positive-volume set. All the axes are therefore parallel to this plane.

For finitely many axes suppose that three are independent. Decompose \(v=\sum_l v_l\) with \(\Delta v_l=1_{D_l}-f_l\), \(f_l=|D_l|>0\), and \(\int v_l=0\). The function \(v_l\) is invariant along \(l\). For an axis \(m\) choose a density point in \(D_m\) with \(F\cdot r\ne0\) for every other axis \(r\). Almost every point of \(D_m\) qualifies because \(\nabla(F\cdot r)=\tau(I-mm^T)r\ne0\) there, while on the zero level the gradient vanishes almost everywhere. Unit vectors represent the unoriented axes in these formulas. Its entire closed \(m\)-line misses all the other \(S_r\), since \(F\) is constant along the line and \(F\cdot r=0\) on \(S_r\). Take a small connected tube around the line disjoint from those closed sets. In the tube \(\partial_m v\) is harmonic by invariance of \(D_m\) and vanishes on a positive-volume subset by \(F\cdot m=0\). It is zero throughout the tube.

Fix \(l\ne m\). Apply all the directional derivatives \(\partial_r\) for \(r\ne l,m\) to \(\partial_m v=0\) in the tube. In the sum for \(v\) only the term for \(v_l\) survives these differentiations. On the flat two-torus perpendicular to \(l\), the projection of the tube contains a full closed line in the projected \(m\) direction. It lies outside the projection of \(S_l\), and \(\Delta_2v_l=-f_l\) in that exterior. The product of the projected directional derivatives, including the projected \(\partial_m\), vanishes on the projected tube. Each direction in this product is nonzero. On a planar lift subtract \(-f_l |x'|^2/4\) from \(v_l\). The remainder is locally the real part of a holomorphic function. The product of \(d\) fixed directional derivatives of its real part is the real part of a nonzero constant times its \(d\)th complex derivative. Here \(d\ge2\), so the contribution from the quadratic is constant. The complex derivative must itself be constant, and \(v_l\) is a polynomial on the lifted exterior component by analytic continuation. The component contains the projected full \(m\)-line. Translation by its period preserves this component, so the polynomial is invariant in that direction.

Now take independent axes \(l,m,n\). In the plane perpendicular to \(l\) the projected full \(m\)- and \(n\)-lines intersect. Their lifts therefore give the same exterior component, on which a single polynomial is invariant in both independent directions. This contradicts \(\Delta_2v_l=-f_l\). This proves coplanarity also in the finite case.

Finally \(C_+=\lim_{z\to\infty}A_\chi(z)/z=\langle\chi(I-B/\tau)\rangle=\sum_l f_l\,ll^T\) is singular by coplanarity. The scalar law would give \(C_+=(f/3)I\). This contradiction proves the claimed sign obstruction. Whether an exact scalar cell with the required sign violation exists is open.
\phantomsection\label{endfile:rough-pole-proof}

\subsection{Joint error at several contrasts}\label{route:joint-error-at-several-contrasts}
\phantomsection\label{file:three-readings}
At \(f=1/2\) the following rational readings have \textbf{best possible maximum absolute error \(1/200000\)}. All three pairs of readings are individually attainable in \(K_{1/2}\).

{\def\LTcaptype{none} 
\begin{longtable}[]{@{}
  >{\raggedleft\arraybackslash}p{(\linewidth - 6\tabcolsep) * \real{0.2500}}
  >{\raggedleft\arraybackslash}p{(\linewidth - 6\tabcolsep) * \real{0.2500}}
  >{\raggedleft\arraybackslash}p{(\linewidth - 6\tabcolsep) * \real{0.2500}}
  >{\raggedleft\arraybackslash}p{(\linewidth - 6\tabcolsep) * \real{0.2500}}@{}}
\toprule\noalign{}
\begin{minipage}[b]{\linewidth}\raggedleft
Contrast \(z\)
\end{minipage} & \begin{minipage}[b]{\linewidth}\raggedleft
Target \(y\)
\end{minipage} & \begin{minipage}[b]{\linewidth}\raggedleft
One geometry's conductivity \(s_*(z)\)
\end{minipage} & \begin{minipage}[b]{\linewidth}\raggedleft
\(y-s_*(z)\)
\end{minipage} \\
\midrule\noalign{}
\endhead
\bottomrule\noalign{}
\endlastfoot
2 & \(120800419/83800000\) & \(604/419\) & \(+1/200000\) \\
3 & \(362399/200000\) & \(453/250\) & \(-1/200000\) \\
4 & \(30400071/14200000\) & \(152/71\) & \(+1/200000\) \\
\end{longtable}
}

\textbf{Upper witness.} Fill a periodic cube up to a null set with disjoint balls. Give each ball four concentric regions, from its center outwards, of conductivities \((1,z,1,z)\) and volume shares \((1/20,1/4,9/20,1/4)\). For a specific packing, inscribe a ball in the open unit cube. Partition its open remainder, up to null sets, into maximal dyadic subcubes with closure in that remainder. Repeat the rule in each subcube. At every generation the balls occupy \(\pi/6\) of the remaining volume, so the packing fills the cube in measure. Extend periodically and assign any phase on null boundaries. The fraction of phase 1 is exactly \(1/4+1/4=1/2\).

For core dipole conductivity \(k\), shell conductivity \(a\) and core volume share \(p\), matching potential and normal current gives

\[
 C_a^p(k)=a+\frac{3ap(k-a)}{3a+(1-p)(k-a)}.
\]

Indeed the shell potential for unit loading has radial factor \(Ar+B/r^2\) times the loading's first spherical harmonic. Continuity of potential and normal flux at the core surface and affine potential on the outer sphere give outer normal flux \(C_a^p(k)e\cdot n\). Hence every ball, under affine boundary data with gradient \(e\), has this same flux coefficient after the three shells are added. Subtract the affine potential inside each ball and extend the difference by zero. These corrections sum in \(H^1\), since for each fixed positive \(z\) their squared gradient norms scale with ball volume. The corrected current minus \(s_*e\) has zero normal flux on each ball and zero integral there by testing the divergence against the coordinate functions. Pasting thus solves the periodic cell equation for every loading and yields \(s_* I\). This also proves exact attainment by a binary measurable geometry. Direct calculation gives

\[
 s_*(z)=C_z^{3/4}\!\left(C_1^{2/5}\!\left(C_z^{1/6}(1)\right)\right)
 =\frac{2z(z^3+21z^2+27z+5)}{11z^3+69z^2+27z+1}.
 \tag{10}
\]

\textbf{All-geometry lower witness.} At \(f=1/2\) let \(g_-(z)=(4z+2)/(z+5)\) and \(g_+(z)=2z(z+2)/(5z+1)\) be the Hashin-Shtrikman endpoints for \(z>1\). Define

\[
 \xi(z)=\frac{2(z+1/z)+5}{9},\qquad
 P_z(s)=\frac{(z+5)(2z+1)}{3(5z+1)}\frac{s-g_-(z)}{g_+(z)-s}.
\]

The concavity proved at the end of this subsection implies that \(P_z(s(z))\) is concave as a function of \(\xi(z)>1\) for any physical response interior to the endpoints. Its proof uses the unit-circle energy inequality to show that the transformed response is Pick in the upper \(\xi\) half-plane. Its Herglotz measure is supported on \(({-\infty},1]\) by analyticity through \((1,\infty)\), giving a nonpositive second derivative there. The individual endpoint laws are already outside the three error intervals in the table.

By (10), \(P_z(s_*(z))=(1+\xi(z))/2\) identically. For the three nodes, \(\xi(3)=\tfrac{11}{21}\xi(2)+\tfrac{10}{21}\xi(4)\). Each interval \([y_i-1/200000,y_i+1/200000]\) lies strictly inside its Hashin-Shtrikman bounds, where \(P_z\) is strictly increasing in \(s\). An error strictly below \(1/200000\) would put both outer readings strictly \emph{above} the affine line and the middle reading strictly \emph{below} it. Concavity requires the middle value to be at least the chord between the two outer ones. This excludes the entire open error box. The rational lower and upper endpoints therefore coincide. This short inequality serves as the lower certificate for the instance.

In fact each pair has an exact fit by one explicit four-layer geometry. For \(q=1-f\), \(a,b\ge0\) and \(D=(1+b)^2+af\), use the four volumes

\[
 \left(\frac{afq}{D},\ \frac{fq}{q+b},\ \frac{(1+b)^2 q}{D},\ \frac{bf}{q+b}\right)
 \quad\text{in phases }(1,z,1,z).                                   \tag{11}
\]

They sum to one and the two \(z\) shares sum to \(f\). Substituting the shares in the three coating recursions gives \(P_z(s(z))=a+b\xi(z)\), with the transform defined above at half fraction. For each pair in the table, the line through its two transformed readings has positive intercept \(a\) and slope \(b\). The script constructs these three geometries and evaluates their two physical conductivities directly. Formula (11) with \(a=b=1/2\) also gives the upper witness (10).

\paragraph*{Why the transformed readings are concave}

We give the physical argument for the inequality used in (12). Write \(h=z-1\), \(q=1-f\), \(\mu=fq/3\). The fixed-fraction isotropic spectral formula is \(s(1+h)=1+fh-h^2\int_0^1(1+ht)^{-1}\,d\nu(t)\) with \(\nu\ge0\) of mass \(\mu\). It follows by the spectral theorem for the gradient projection and multiplication by the phase indicator, with the mass fixed by the isotropic second derivative. Define \(G=[\mu^{-1}\int(1+ht)^{-1}\,d\nu]^{-1}\) and \(u=(G-1)/h\). The function \(G\) is positive Pick, analytic off the negative axis, with \(G(1)=1\). Dividing its Herglotz representation minus its value at 1 by \(h\) gives \(\Im u\le0\) in the upper half-plane. Direct substitution in the transform gives

\[
 B=\frac{u-q/3}{1-f/3-u},\qquad a_2(z)=\frac3{2z+1},\qquad P_z(s(z))=B(z)/a_2(z).
\]

For an interior response \(B\) is positive on the positive axis and has nonpositive imaginary part in the upper half-plane. The two Hashin-Shtrikman laws are separated by the strictly interior box faces in (12).

For \(z=e^{i\theta}\), \(0<\theta<\pi\), let \(E\) and \(J\) be the complex matrices of three electric fields and currents, with means \(I\) and \(sI\). For \(b(T)=|\operatorname{tr}T|^2-\sum_{ij}T_{ij}\overline{T_{ji}}\) one has \(b(J)=b(E)\) pointwise since \(J=aE\) and \(|a|=1\). Each nonzero gradient mode has the form \(nv^T\) with real unit \(n\), giving \(b=0\). On a current mode \(n^TT=0\), rotate to make the last row zero. In the first two rows the diagonal cross term is bounded by the squared diagonal entries and the off-diagonal cross term by the squared off-diagonal entries, so \(b(T)\le\|T\|_F^2\). Parseval, with the mean mode accounted for, now gives \(6\le3|s|^2+\int\|E\|_F^2\). Orthogonality gives \(\int a\|E\|_F^2=3s\). Since \(a+z\bar a=1+z\) in either phase, it follows that

\[
 \Phi_z=s(z)s(1/z)+\frac{s(z)+z\,s(1/z)}{1+z}-2\ \ge0
\]

on the arc, where \(s(1/z)=\overline{s(z)}\). This argument passes to the response closure by local uniform convergence.

Put \(D(z)=qz+f+2+((q+2)z+f)B(z)\). Substitution gives

\[
 D(z)D(1/z)\Phi_z=\frac{6fq(z-1)^3}{z(z+1)}
       \left[\frac{B(z)}{a_2(z)}-\frac{B(1/z)}{a_2(1/z)}\right].
\]

On the arc the left side is nonnegative, and the right side is \(48fq\sin^3(\theta/2)\Im(B/a_2)/\cos(\theta/2)\). Take arguments of \(B,a_2\) in \([-\pi,0]\). In \(\Im z>0\), \(|z|>1\), their difference is bounded harmonic. Its lower boundary limits are nonnegative on the arc by this identity, zero on the real axis to the right of 1, and nonnegative to the left of \(-1\) because there \(\arg a_2\to-\pi\). The bounded harmonic maximum principle gives \(\arg B-\arg a_2\ge0\) throughout the domain. The map \(z\mapsto\xi(z)\) sends it conformally onto the upper half-plane, so \(P\) is Pick there. It is real analytic on \((1,\infty)\) and its Herglotz measure is supported on \((-\infty,1]\). Differentiating twice on \((1,\infty)\) proves \(P''\le0\), the claimed three-reading inequality.
\phantomsection\label{endfile:three-readings}

\subsection{A finite two-sided procedure}\label{route:a-finite-two-sided-procedure}
For rational $f$, positive rational contrasts $z_i$, rational readings
$y_i$ and rational $\epsilon>0$, the procedure below brackets the
infimum of $\max_i|s(z_i)-y_i|$ by rational numbers less than $\epsilon$
apart. The upper end specifies one exact-fraction periodic cubic cell.
The finite sizes in the general proof can be extremely large.
\phantomsection\label{file:joint-algorithm}
\paragraph*{A joint energy with a phase-fraction observable}

Set \(a_{ib}=1+(z_i-1)b\) for \(b\in\{0,1\}\), \(\lambda=\min(1,\min_i z_i)\), and \(\Lambda=\max(1,\max_i z_i)\). Work on the periodic unit cube with volume one. The field space has dimension \(d=18k+2\). For every \(i\) and each loading direction \(r=1,2,3\) it contains an electric field \(E_{ir}\in\mathbb R^3\) with gradient fluctuations and a current \(J_{ir}\in\mathbb R^3\) with divergence-free fluctuations. Two additional scalar fields \(u,v\) have unrestricted mean-zero fluctuations. Write \(\mathcal F\) for the direct sum of these real \(L^2\) fluctuation spaces. Use two quadratic forms and their common minimum,

\[
 Q_b(X)=\sum_{i,r}\big(a_{ib}|E_{ir}|^2+a_{ib}^{-1}|J_{ir}|^2\big)
             +u^2+(v-bu)^2,\qquad
 W(X)=\inf_{F\in\mathcal F}\int\min_{b=0,1}Q_b(X+F).
 \tag{1}
\]

The minimizing bit at a spatial point is shared by all the loads. Both wells have eigenvalues between \(c=\min(1/3,\lambda,1/\Lambda)\) and \(H=\max(3,\Lambda,1/\lambda)\). For the coupled scalars the only nonidentity block has eigenvalues \((3\pm\sqrt5)/2\).

For a proposed vector \(t\), prescribe the means

\[
 X(t)=(E_{ir}=e_r,\ J_{ir}=t_i e_r,\ u=1,\ v=f),\quad
 B(t)=1+6\sum_i t_i,\qquad G(t)=W(X(t))-B(t).
\]

There is an exact identity involving the \emph{same} unrestricted binary cell \(\chi\) at all contrasts. Write \(A_i=A_\chi(z_i)\) for its effective tensors and \(p=\int\chi\). Then

\[
 G(t)=\inf_{\chi}\left[
 \sum_i \operatorname{tr}\big((A_i-t_i I)A_i^{-1}(A_i-t_i I)\big)
       +\frac{(p-f)^2}{1+p(1-p)}\right].                         \tag{2}
\]

Indeed, measurable pointwise selection in (1) interchanges the infima over \(F\) and \(\chi\). For fixed \(\chi\) the Dirichlet and complementary minima give \(e_r^TA_i e_r\) and \(t_i^2 e_r^TA_i^{-1}e_r\). Subtracting \(2t_i\) per load gives the trace in (2). For the scalars put \(u=1+w\), \(\int w=0\). Minimizing over \(v\) makes \(v-\chi u\) constant. Its contribution together with \(\int u^2-1\) is

\[
 \|w\|_2^2+\big(f-p-\langle\chi-p,w\rangle\big)^2.
\]

Minimization along \(\chi-p\), whose squared norm is \(p(1-p)\), proves the last term of (2), including \(p=0,1\).

\paragraph*{Converting energy accuracy to conductivity accuracy}

Let \(T=\prod_i[\min(1,z_i),\max(1,z_i)]\). For \(0<f<1\) put \(q_f=\min(f,1-f)\) and use the rational constant

\[
 C=\Lambda+\frac{2\Lambda^2}{\lambda q_f}.
\]

For \(t\in T\) the joint excess and the sup distance to the physical set satisfy

\[
 \boxed{\quad \operatorname{dist}_\infty(t,K_f)\le C\sqrt{G(t)},\qquad
 0\le G(t)\le\frac{3k}{\lambda}\operatorname{dist}_\infty(t,K_f)^2 .\quad}       \tag{3}
\]

For the right bound use exact-fraction cubic cells approaching any \(s\in K_f\) in (2). Their trace terms tend to \(3\sum_i (s_i-t_i)^2/s_i\). For the left bound let \(R\) be the expression in square brackets in (2) for any one cell. Uniform ellipticity gives \(\|A_i-t_i I\|_{\rm op}\le\sqrt{\Lambda R}\) for every \(i\), and \(|p-f|\le\sqrt{5R/4}\).

Correct its phase fraction on a \emph{coarse} scale. For \(p>f\) replace share \(\theta=(p-f)/p\) by pure phase 0. For \(p<f\) replace share \(\theta=(f-p)/(1-p)\) by pure phase 1. For equality take \(\theta=0\). Always \(\theta\le |p-f|/q_f\). There is a single \(s\in K_f\) with, for all \(i\),

\[
 \left[\frac{1-\theta}{\lambda_{\min}(A_i)}+\frac{\theta}{a_i^0}\right]^{-1}
 \le s_i\le (1-\theta)\lambda_{\max}(A_i)+\theta a_i^0,          \tag{4}
\]

where \(a_i^0\) denotes the conductivity of the chosen pure phase at \(z_i\).

To obtain (4), cut off the pure share by a plane inside one of the 48 signed-permutation chambers of the cube. In the remainder insert copies of \(\chi\) in whole microcubes and fill the leftover boundary layer at fraction \(p\), for instance by a volume cut. Reflect and permute the entire chamber, then periodize. Each of these binary cells has fraction exactly \(f\). In the pure region test with constant electric field and constant current of unit mean. In each textured region add the rotated, rescaled correctors for those means. More explicitly use smooth periodic approximations to the voltage and to a vector potential for the mean-zero current. Multiply each rescaled potential by a cutoff vanishing near the chamber and region boundaries, then take its gradient or curl. At fixed boundary width the interior energies converge by periodic averaging as the microperiod tends to zero. The potentials carry a factor of the microperiod, so their cutoff-derivative terms vanish in that limit. Bounded smooth trials cost a quantity tending to zero with the boundary width. Taking these limits and then improving the smooth trials gives the upper inequalities for conductivity and resistivity in (4). Cubic symmetry makes the tensors scalar. A common subsequence gives all \(k\) readings. This is the archived chamber proof with one additional pure subregion.

The arithmetic bound moves the upper extreme of \(A_i\) by at most \(\theta\Lambda\). The harmonic bound moves its lower extreme \(m_i\) down by at most \(\theta\Lambda^2/\lambda\), by writing it as \(m_i/[1+\theta(m_i/a_i^0-1)]\). Thus \(|s_i-t_i|\le\sqrt{\Lambda R}+\theta\Lambda^2/\lambda\). Taking the infimum proves (3). In particular \(K_f=\{t\in T:G(t)=0\}\). For \(f=0\) or 1 the pure cell itself gives \(L_*=U_*\) directly.

\paragraph*{Finite spatial exclusion}

The spatial lower-bound proof (S1)-(S13) applies to the common two-well energy (1) with the following Fourier projection and envelope estimates.

At a nonzero Fourier frequency \(j\), the orthogonal projection \(\Pi(j)\) onto the fluctuation space has blocks \(P_j=jj^T/|j|^2\) for \(E\), \(I-P_j\) for \(J\) and \(I_2\) for \((u,v)\). Set \(\Pi(0)=0\).

\begin{itemize}
\item
  For fixed \(\chi\), minimizing the positive quadratic energy in (1) gives a quadratic form in the mean between \(cI\) and \(HI\). Consequently \(W-H|\cdot|^2\) is concave, \(W\) is positively two-homogeneous, and \(c|X|^2\le W(X)\le H|X|^2\).
\item
  The full Jensen inequality \(\int W(X+F)\ge W(X)\) holds for \(F\in\mathcal F\). In (1) paste fine trigonometric near-minimizers on a partition where a coarse trigonometric field is nearly constant. Project the pasted fluctuations back onto \(\mathcal F\). For fixed local mode \(j\ne0\) and fixed offset \(l\), \(\Pi(hj+l)\to\Pi(j)\) as the repetition factor \(h\to\infty\). Approximating the piece indicators in \(L^2\) by Fourier polynomials proves vanishing projection error for the bounded local trials. Periodic averaging, quadratic continuity, and density give the inequality. Coercivity bounds the mean-square norms of the local trials uniformly in the partition, as in (S2).
\item
  Every field-space coordinate admits an axial plane wave. Use a parallel wave vector for a gradient coordinate, a perpendicular one for a current coordinate, and any axis for a free scalar. Two-valued waves imply separate convexity of \(W\), hence \(\Delta W\ge0\) in distributions. Smoothing and using the Hessian upper bound \(2H I\) gives the lower Hessian bound \(-2(d-1)H I\) by taking the trace. Thus \(DW\) is globally Lipschitz with constant \(M_0=2(d-1)H\).
\item
  The angular estimate used in (S6) holds with the same constant 8. For \(j,k\ne0\) and \(a>0\), the projection difference is at most \(4|j-ak|/(a|k|)\), by normalizing the two directions and subtracting their dyads. The complementary blocks have the same difference and the scalar blocks have zero difference. The archived window, band separation and summation estimates (S4)-(S7) therefore apply in this space with the original constants. In particular, integral Jensen up to radius \(N=2^{8m^3+4m+1}\) for a function \(g\) with \(M\)-Lipschitz gradient implies \(\int g(X+F)-g(X)\ge-M E_m\|F\|_2^2\) at \emph{every} frequency, where
\end{itemize}

\[
 E_m=m^{-1}+(2m^2)^{-1}+3\,2^{-m}+2\,2^{-2m}+26\,2^{-4m}\le 6/m.
\]

Choose a rational \(0<\beta\le c\). The following degree and bandwidth suffice. They can be stored in factored form.

\[
\begin{gathered}
 r=\beta/(2M_0),\qquad A_r=(1+\lceil4d/r\rceil)^{d^2},\qquad
 \xi=\frac{1-r^2/(4d)}{1-r^2/(16d)},\qquad 2 A_r\xi^n\le r,\\
 \tau=\beta/(1+5\beta/c),\qquad M=M_0+8\beta,\qquad
 m=\lceil6M/\tau\rceil,\quad N=2^{8m^3+4m+1},\\
 K=(1+\lfloor864M/\tau\rfloor)N^4,\qquad q_K=16MN(2N+1)^3/K.
\end{gathered}                                                        \tag{5}
\]

Take any positive integer \(n\) satisfying the inequality. Use polynomials \(P\) of degree at most \(n\) in the entries of \(T_0\in O(d)\), with real coefficients. Require equal values when two orthogonal matrices have the same first column. Set \(h(0)=0\) and \(h(X)=|X|^2P(T_0)\) when \(T_0 e_1=X/|X|\). Put \(g=h-\tau|\cdot|^2\) and impose all three tests below.

\begin{enumerate}
\def\labelenumi{\arabic{enumi}.}
\item
  \(h\le Q_0,Q_1\) on the unit sphere.
\item
  \(|g(A+Y)+g(A-Y)-2g(A)|\le M|Y|^2\) for all \(A,Y\).
\item
  For every constant \(A\) and every real, mean-zero trigonometric polynomial \(S\) with modes \(0<|j|\le N\) in \(\operatorname{ran}\Pi(j)\), require

  \[
  K^{-3}\sum_{x\in\{0,\ldots,K-1\}^3/K}\big[g(A+S(x))-g(A)\big]
      \ge q_K\|S\|_2^2.                                             \tag{6}
  \]
\end{enumerate}

These are finite first-order conditions over the reals. For example, quantify real sine and cosine coefficient vectors, one pair per \(\pm j\). Constrain both vectors by \(\Pi(j)v=v\), and use half the sum of their squared norms for \(\|S\|_2^2\). Grid constants are real algebraic roots-of-unity coordinates with their specified embeddings. The fiber condition gives a continuous semialgebraic \(h\), including at 0. Thus real quantifier elimination (QE) decides the tests and inequalities involving the supremum \(V_\beta(X)=\sup h(X)\) over their feasible coefficients.

Their feasible supremum satisfies

\[
              V_\beta(X)\le W(X)\le V_\beta(X)+(5H/c)\beta |X|^2.    \tag{7}
\]

For the first inequality, the second-difference test gives an \(M\)-Lipschitz gradient by mollification. After subtracting the linear term \(Dg(A)\cdot S\), the grid integration error is at most \(q_K\|S\|_2^2\) by the Fourier sup bounds for \(S,\nabla S\) and the grid-cube diameter. Thus (6) gives integral Jensen up to \(N\). The all-frequency loss \(M E_m\|F\|_2^2\) is absorbed by adding back \(\tau|\cdot|^2\). The obstacle test and (1) give \(h\le W\).

For the other inequality the exhaustion proof (S10)-(S13) applies to this \(W\). In brief, average \(W(UX)\) on \(O(d)\) against Haar measure reweighted by \(((d+\operatorname{tr}U)/(2d))^n\) and normalized to mass one. Denote the result by \(W_n\). A Frobenius ball of radius \(r/2\) has Haar mass at least \(1/A_r\) by covering the group with matrix-entry boxes. The mass outside radius \(r\) after reweighting is at most \(A_r\xi^n\). Therefore \(\mathbb E\|U-I\|_{\rm op}\le\beta/M_0\). On the sphere \(|W_n-W|\le\beta\). Comparing \(W(UX+UF)\) with \(W(UX+F)\), subtracting the mean-zero linear term \(DW(UX)\cdot(U-I)F\) and applying the gradient bound gives a Jensen defect of at most \(\beta\|F\|_2^2\). Here \(U\) acts on field coordinates. Admissibility is needed only for the original \(F\). Changing variables in the group integral makes \(W_n(T_0e_1)\) a polynomial of degree at most \(n\) in \(T_0\).

Now \(h=(W_n+3\beta|\cdot|^2)/(1+5\beta/c)\) lies below \(W\) and has Jensen margin \(2\tau\|F\|_2^2\). Its \(g\) passes the second-difference test, and its integral Jensen margin \(\tau\) exceeds \(2q_K\), which ensures (6) after quadrature. Finally \(W-h\le(5H/c)\beta|\cdot|^2\). This proves (7) without knowing the coefficients of the averaged envelope.

\paragraph*{Obtaining the two bounds and a geometry}

For the prescribed gap \(\varepsilon\) choose

\[
 \eta=(\varepsilon/(4C))^2,\qquad H_X=3k(1+\Lambda^2)+1+f^2,\qquad
 \beta=\min\{c,\eta c/(5H H_X)\}.
\]

A finite outer set and its certified accuracy follow from (3) and (7),

\[
 O=\{t\in T:\ V_\beta(X(t))\le B(t)\},\qquad
 K_f\subseteq O\subseteq\{t:\operatorname{dist}_\infty(t,K_f)\le\varepsilon/4\}.
 \tag{8}
\]

Membership in \(O\) says that \emph{every} polynomial passing (5)-(6) satisfies \(h(X(t))\le B(t)\). Thus the distance \(l=\min_{t\in O}\|t-y\|_\infty\) is a real algebraic number computable by QE. The set is compact as an intersection of closed sublevels on \(T\), and it is nonempty because it contains \(K_f\). One can also bisect rational radii, asking whether the error box meets \(O\). A negative decision excludes the entire box for every geometry. Equation (8) gives \(l\le D_*\le l+\varepsilon/4\). Round down to a rational \(L_*\ge0\) with \(l-L_*<\varepsilon/8\).

A smaller, replayable exclusion can consist of specific minorants satisfying tests 1 to 3 and a rational radius. Intersect their inequalities \(h(X(t))\le B(t)\) with \(T\) and the open error box of that radius. An empty intersection certifies the radius as a lower bound. An algebraic decision transcript for (8) gives a general finite alternative. In either case the difference between the witnessed upper error and this lower bound limits the improvement possible at all finer scales.

For an upper witness, enumerate grids of \(n_0^3\) torus voxels for \(n_0=2,4,8,\ldots\) and orders of their orbits under signed coordinate permutations. Fill the orbits in order to total fraction \(f\), putting the same centered subcube in each voxel of the final, partially filled orbit. The list of voxel fractions \(q_h\) specifies the geometry exactly, with side \(q_h^{1/3}/n_0\) in voxel \(h\). For each reading compute the archived corner and center graph enclosures. On edges \((h,h+e_j)\) their rational conductances are

\[
 c^+_{hj}=1+(z_i-1)\frac14\sum_{\sigma\in\{0,1\}^3,\,\sigma_j=0}q_{h-\sigma},\qquad
 (c^-_{hj})^{-1}=1+(z_i^{-1}-1)(q_h+q_{h+e_j})/2 .
\]

For either array solve the periodic pinned voltage problem

\[
 G_{n_0}(c)=n_0^{-3}\min_w\sum_{h,j} c_{hj}(\delta_{j1}+w_{h+e_j}-w_h)^2 .
\]

Its normal equations are rational and nonsingular after pinning. Denote the resulting pair by \(g_i^-,g_i^+\). Trilinear voltages and face-normal current interpolation certify \(g_i^-\le s_\chi(z_i)\le g_i^+\) for that one cell. Use the rational upper error

\[
 U(\chi)=\max_i\max\{|g_i^- - y_i|,|g_i^+ -y_i|\}.
\]

Stop at the first \(U(\chi)-L_*<\varepsilon\). Graph density proves termination. In detail, approximate any fixed cubic cell by majority voxel colors, ordering the majority-phase orbits first and adjusting the total to \(f\). The resulting indicators converge in \(L^1\), since adjustment changes volume by precisely the fraction mismatch. Trials from a fixed smooth voltage and a smooth solenoidal current for the original cell show convergence of the two graph enclosures at every one of the finitely many nodes. Cubic density then gives witness errors arbitrarily close to \(D_*\), and \(D_*-L_*<3\varepsilon/8\). Return the voxel list and the graph solves, or any shorter exact witness meeting the same inequality. Running maxima of the lower bounds and minima of the witnessed upper bounds for \(\varepsilon\downarrow0\) give nested convergent brackets.

There is also a computable finite library for \emph{every target} on the given nodes at fixed \(f\). For library levels through \(n\), let \(U_n(y)\) be the least graph-certified upper error over all their cells. Compute (8) with distance budget \(\varepsilon/4\), then increase \(n\) until \(U_n(t)<\varepsilon\) for every \(t\in O\). This is a finite real-algebraic coverage test. Each \(U_n\) is 1-Lipschitz and decreases with \(n\). On \(K_f\) these functions tend uniformly to zero by graph density and a finite cover of the compact set. Hence the search succeeds by (8). At that level, for \emph{any} readings \(y\), a nearest physical vector \(s\) gives \(U_n(y)\le U_n(s)+\|y-s\|_\infty<D_*+\varepsilon\). The explicit bracket \([\max(0,U_n(y)-\varepsilon),U_n(y)]\) and the winning cell can then be read from the same library. Using a smaller internal tolerance gives a strict requested gap.

The degree, number of Fourier coefficients and grid in (5) make the full exclusion a computability result. The finite lower calculation quantifies over trial \emph{fields} at its chosen bandwidth, with a bound for all omitted frequencies. Competing geometries may have arbitrarily fine features. To reproduce the parameter calculation for three readings at gap \(10^{-6}\), run \texttt{python\ spatial\_schedule.py\ -\/-fraction\ 1/2\ -\/-gap\ 1/1000000\ 2\ 3\ 4}. It gives \(d=56\), \(C=68\), and an outer-set distance budget of \(1/4000000\). The worked instance uses the analytic exclusion below.
\phantomsection\label{endfile:joint-algorithm}

\subsection{The two-flux separator in three dimensions}\label{route:the-two-flux-separator-in-three-dimensions}
At $Z=25$, let $w=841/16900$ and $v=10201/16900$. For
divergence-free column matrices $J,L$ with free means put
$D=w|J|^2-v|L|^2$, using squared Frobenius norms. The supremum of the
two-reading score over physical half-volume materials is
\[
 \sup_{J,L}\left\{\frac13\Big\langle
   2w\operatorname{tr}J+2v\operatorname{tr}L
   -\frac{Z+1}{2}(w|J|^2+v|L|^2)\Big\rangle+
 \frac{Z-1}{6}\inf_{b\in\mathbb R}\langle|D-b|\rangle\right\}.
\]
Indeed the maximum of $\langle(2\chi-1)D\rangle$ over half-volume masks is
the mean absolute deviation about a median. The complementary principle
then supplies both conductivity readings with the same mask. Cubic
symmetrization or simultaneous trace isotropization gives the isotropic
form. The spectral variance estimate at the two reciprocal contrasts
strictly separates the laminate class from the test value given in the
article. The spectral bound and its constants follow below.
\subsection{A bound valid at arbitrary laminate depth}\label{detail:two-flux-proof}\label{file:two-flux-proof}
\paragraph*{An explicit unrestricted laminate test at arbitrary pole count}

Let \(s=\operatorname{tr}K/3\), let \(0<f=s'(1)<1\), and use the \textbf{normalized correction measure} \(\nu\):

\[
 s(1+h)=1+fh-\alpha h^2\int\frac{d\nu(t)}{1+th},\qquad
 \alpha=f(1-f)/3,\qquad \nu([0,1])=1 .
\]

Write \(m=\int t\,d\nu\) and \(V=\int(t-m)^2d\nu\). For the closure of \textbf{all} hierarchical laminates with arbitrary normals, the following necessary condition holds:

\[
 \boxed{\qquad a(m-a)\ \le\ 8\sqrt2\,a\,d+49d^2,\qquad
           a=3m+f-1\ge0,\qquad d=((1-f)V)^{1/4} .\qquad}           \tag{1}
\]

Only the zero-free determinant and the weighted resolvent inequality are needed for the two-reading separator.

\paragraph*{A zero-free determinant controls the negative eigenvalues}

Use the electric field matrix \(E(z)\) at a volume-weighted leaf of a finite laminate and \(W=\mathcal P(\chi I)\), \(T=\mathcal P\chi\mathcal P\) on the gradient fluctuations. At every leaf \(W\) is real symmetric and \(\operatorname{tr}W=\chi-f\). The three constant loadings give the resolvent field \(E(1+h)=I-h(I+hT)^{-1}W\). At each positive-volume leaf \(\det E(z)\) has no zero on the slit plane: at a split its determinant factor is the normal conductivity of the sibling divided by the reversed arithmetic normal conductivity, and each is a zero-free positive Pick function.

Fix \(0\le q\le1\) and set \(X=qI-W\). For \(\xi\) off \([0,1]\), the exact identity

\[
 (\xi-q)E(1-1/\xi)=\xi I-X+R(\xi),\qquad
                 R(\xi)=(\xi-T)^{-1}(T-q)W                         \tag{2}
\]

compares a matrix with zero-free determinant in the left half-plane to a symmetric linear pencil. For each \(u>0\) set

\[
 G_u=\int_0^\infty e^{-u v}\|[e^{-vT}(T-q)W]_{\rm leaf}\|_F\,dv .
\]

The Laplace integral for the resolvent gives \(\|R(\xi)\|\le G_u\) throughout \(\Re\xi\le-u\) at this leaf. Cauchy-Schwarz in \(v\) followed by the spectral theorem gives

\[
 \mathbb E G_u^2\ \le\ \frac{f(1-f)}{u}
               \int\frac{(t-q)^2}{u+2t}\,d\nu(t)
 \ \le\ \frac{f(1-f)}{u^2}\int(t-q)^2d\nu(t).                      \tag{3}
\]

On the high-conductivity label \(\chi=1\) write \(\lambda_1\le\lambda_2\le\lambda_3\) for the eigenvalues of \(X\). Their sum is \(a=3q+f-1\); assume \(a\ge0\). Then, pointwise with \(g=G_u\),

\[
                   \lambda_1\ge-u-3g,\qquad \lambda_2\ge-u-g .        \tag{4}
\]

In fact a vertical line \(\Re\xi=c\le-u\) at distance greater than \(g\) from the spectrum of \(X\) cannot have an eigenvalue on its left. Otherwise join the two pencils in (2) by the straight homotopy with remainder \(\tau R,\ 0\le\tau\le1\). It stays invertible on that line by the singular values of \(\xi I-X\), and on a sufficiently large left semicircle by the same estimate. The determinant winding number would preserve a zero in the enclosed half-disk, contradicting the zero-free determinant at \(\tau=1\). If \(\lambda_2<-u-g\), the trace forces \(\lambda_3>2(u+g)\), so the line at \(-u\) gives this contradiction. If \(\lambda_1<-u-3g\), the trace gives \(\lambda_3\ge-\lambda_1/2>g\). On the interval \((\lambda_1,-u)\) the radius-\(g\) neighborhoods of the first two eigenvalues cover length at most \(g+2g\), and the third misses the interval. An uncovered line again yields the contradiction. This proof includes two coinciding negative eigenvalues.

For a fixed sum \(a\ge0\), bounds (4) imply

\[
 \operatorname{tr}X^2\le (u+3g)^2+(u+g)^2+(a+2u+4g)^2 .
\]

Indeed the vector with entries \(-u-3g,-u-g,a+2u+4g\) has the same sum, smaller or equal smallest entry and larger or equal largest entry than the actual eigenvalues; convexity of the sum of squares gives the bound. Equivalently transfer the excess in each of the two lower entries to the largest one.

The phase-one second moment is exactly

\[
 S=\mathbb E_1\operatorname{tr}X^2=3q^2-2q(1-f)+(1-f)m .
\]

This follows from \(\operatorname{tr}W=1-f\) on that label and \(\mathbb E[\chi\operatorname{tr}W^2]=\operatorname{tr}(W^*TW)=f(1-f)m\). Combining with (3) gives the explicit necessary test for \textbf{every} \(u>0,\ 0\le q\le1,\ a=3q+f-1\ge0\):

\[
 \begin{aligned}
 c^2&=\frac{1-f}{u}\int\frac{(t-q)^2}{u+2t}\,d\nu(t),\\
 S&\le a^2+4au+8ac+6u^2+24uc+26c^2 .                             \end{aligned}\tag{5}
\]

For \(q=m\) one has \(S=ma\), \(c\le\sqrt{(1-f)V}/u\). Set \(u=\sqrt2\,d\) to obtain (1), using a limit when \(V=0\). The proof for finite trees passes to every whole-function limit by continuity of the bounded measure integrands. Explicitly \(m=s'''(1)/(6\alpha)\), \(V=-s''''(1)/(24\alpha)-m^2\). At \(f=1/2,m=1/5\), (1) gives \(d\ge(9/5-4\sqrt2/5)/98\) and \(V\ge2d^4>4.3\,10^{-9}\).

\paragraph*{A common remainder estimate at every vertical line}

Write \(T_1=\chi\Gamma\chi\), \(W=T_1I\) on the conducting phase, and \(X=qI-W\). For the three constant loadings, use the squared Frobenius norm and the conditional volume measure on this phase. At a positive-volume leaf, the determinant of \[
      \xi I-X+R(\xi),\qquad R(\xi)=(\xi-T_1)^{-1}P(T_1)I,
      \qquad P(t)=t(t-q),
\] has no zero in \(\Re\xi<0\). Indeed this matrix is the product of the leaf electric field at contrast \(1-1/\xi\) and the scalar \(\xi-q\). The field determinant is a product of ratios of zero-free normal conductivity functions.

Use the nonnegative scalar at each leaf \[
      H^2=\tfrac12\int_0^\infty\!\big\|[e^{-sT_1}P(T_1)I]_{\rm leaf}\big\|_F^2 ds .
\] Cauchy-Schwarz in the Laplace representation gives \(\|R(\xi)\|_{\rm op}\le H/\sqrt{-\Re\xi}\) for every vertical line simultaneously. By the spectral theorem, \[
          \mathbb E_1 H^2=\tfrac34\int t(t-q)^2d\gamma
                     =\tfrac{1-f}{4}[V+(\bar q-q)^2].               \tag{3}
\] The integrand vanishes at the zero atom. The same \(H\) therefore controls lines chosen using the leaf eigenvalues, with the single expectation (3).

If the line \(\Re\xi=-u<0\) leaves an eigenvalue of \(X\) strictly to its left, its distance to the real spectrum of \(X\) is at most \(H/\sqrt u\). Otherwise the straight homotopy from the linear pencil to the full one is invertible on the line and on a sufficiently large semicircle to its left, where the rational remainder decays. The full determinant would retain a zero inside the semicircle.

Suppose the least eigenvalue is \(-d<0\) and the second is \(-kd\le0\), \(0\le k\le1\). The third is \(A+(1+k)d\). Use first the line \(-kd/3\) when \(k>0\), then the line \(-(1+k)d/2\) when \(k<1\). The distances to the spectrum are \(2kd/3\) and \((1-k)d/2\). This yields \[
 H/d^{3/2}\ge h(k)=\max\left\{\frac{2k^{3/2}}{3\sqrt3},
                       \frac{1-k}{2}\sqrt{\frac{1+k}{2}}\right\}.
\] If the second eigenvalue is positive, the midpoint line \(-d/2\) gives the same estimate with \(k=0\); \((-d,0,A+d)\) then majorizes the actual eigenvalues. If all three are nonnegative, their squared sum is already at most \(A^2\).

The two expressions for \(h\) cross at \(k=3/5\), where \(h=2/(5\sqrt5)\). On the left, \(h\) decreases while \(1+k\) and \(1+k+k^2\) increase. On the right the relevant ratios reduce to constant multiples of \((1+k)/k\) and \((1+k+k^2)/k^2\), both decreasing. Thus \[
 d(1+k)\le8(H/2)^{2/3},\qquad
 d^2(1+k+k^2)\le49(H/2)^{4/3}.
\] The excess of the squared sum over \(A^2\) is at most \(2A d(1+k)+2d^2(1+k+k^2)\). Its conditional mean is \(S-A^2\), since \(\operatorname{tr}W=1-f\) on the conductor and \(\mathbb E_1\operatorname{tr}W^2=(1-f)\bar q\). Applying Jensen's inequality to (3) gives (1). All inputs in (1) are ordinary Taylor coefficients, hence are continuous under nonpure whole-function limits.

\paragraph*{A common bound on the resolvent remainder}

Let \(b=1-f\), \(\nu\) denote the correction probability, and \(d\mu=d\nu/[t(1-t)]\). The phase susceptibility satisfies \(3t\,d\gamma=b\,d\nu\). For \(T=\chi\Gamma\chi\) and \(u=\chi I\), the phase pencil remainder is the Laplace transform, up to sign, of \(e^{-sT}T(T-q)u\). On \(\Re\xi\le-U\), apply Cauchy-Schwarz with the time weight \(s+d\) separately at each leaf. Take the conditional mean square over phase one and use the spectral theorem to integrate \((s+d)e^{-2ts}\). This gives the admissible remainder majorant \[
 G^2\le b\,F(2dU)\int(1-t)(1/4+dt/2)(t-q)^2\,d\mu,\qquad
 F(x)=\int_0^\infty \frac{e^{-s}}{s+x}\,ds .                        \tag{1}
\] In detail the scalar integral from the other factor in Cauchy-Schwarz is \(\int_0^\infty e^{-2Us}/(s+d)\,ds=F(2dU)\). The zero mode of \(T(T-q)u\) vanishes. To pass to whole-function closure, first truncate the time integral at \(N\). The remaining Laplace tail has conditional \(L^2\) majorant \(e^{-UN}\|T(T-q)u\|_{L^2}/U\), uniformly in all trees. The finite-time spectral integrand is continuous on \([0,1]\). Pass to the limiting measure at fixed \(N\), then let \(N\to\infty\) and apply (1) to the endpoints of that measure. This avoids interchanging the insulating and homogenization limits.

With \(a=3q-b\ge0\), \(\bar t=\int t\,d\nu\), the winding condition therefore requires \[
 a(q-a)+b(\bar t-q)\ \le\ 4aU+8aG+6U^2+24UG+26G^2.                \tag{2}
\] Take \(U=1/600,\ d=2/3\). The quadratic factor in (1) is at most \(49/192\), and \(F(1/450)\cdot49/192<709/500\). For a rational check bound \(e^{-s}\) by \((1+s/64)^{-64}\). The substitution \(y=(1+s/64)^{-1}\) integrates the upper bound to \(\beta^{-64}[\log(28800)-\sum_{j=1}^{63}\beta^j/j]\), \(\beta=1-1/28800\). A positive finite Taylor sum for the exponential proves \(\log(28800)<10.269\).

For the numerical thresholds put \(q=(\int t\,d\mu)/K\), \(K=\int d\mu\), and set \(J=\int(t-q)^2d\mu\). Then \(G^2\le709J/1000\) at half fraction. If \(I\) bounds \(\int(t-q_0)^2d\mu\) at the chosen target \(q_0\), use \(Kq(1-q)=1+J\) to get \[
 J\le\frac{Iq(1-q)-(q-q_0)^2}{q(1-q)+(q-q_0)^2}.
\] Also \(\bar t-q=(1-2q)J-\int(t-q)^3d\mu\ge-qJ\). For weight 16 take \(q_0=1/5,\ I=12\cdot8.9\cdot10^{-6}\). For weight 15 take \(q_0=(1+\sqrt{15})^{-1}\) and \(I=12q_0^2[(40-\sqrt{15})/6-6.02099]\). Enclose \(\sqrt{15}\) between \(3.87298334\) and \(3.87298335\). In either case \(K\ge4\) gives \(|q-q_0|<0.006\). Dividing this interval into 1,024 equal rational subintervals contradicts (2) with margins greater than \(0.0000099\) and \(0.0001252\), respectively.

\paragraph*{Using the midpoint between the two lowest eigenvalues}

Use the same weighted time energy at every vertical line for each leaf. Denote its square root by \(h\), and set \(B(x)=\sqrt{F(2dx)}\). Thus the remainder on the line at \(-x\) has operator norm at most \(hB(x)\). With \(d=2/3\) and \(f=1/2\), (1) gives \(\mathbb E_1 h^2\le49J/384\). The nonconstant factor \(B\) is kept outside this expectation. Let \(a\ge0\) be the trace of the real symmetric pencil coefficient \(X\) at a leaf.

If its two lowest eigenvalues are negative, write them \(-m-k,-m+k\) with \(0\le k\le m\). The third is \(a+2m\). A vertical line at \(-m\) whose distance to both eigenvalues exceeded the remainder would preserve a zero strictly to its left by determinant winding. Thus \(k\le hB(m)\), including a repeated eigenvalue. For any \(u>0\), also \(m-k\le u+hB(u)\), since otherwise the line at \(-u\) would preserve both negative roots. If there is just one negative eigenvalue, use the upper bound on the squared trace given by the majorizing triple \((-2m,0,a+2m)\), \(2m=-\lambda_1\), and put \(k=m\). The line at \(-m\) has distance at least \(m\) from the actual spectrum, so the same bounds apply. With no negative eigenvalues take \(m=k=0\). Consequently on every leaf \[
  \operatorname{tr}X^2-a^2\le4am+6m^2+2k^2,\qquad
  k\le hB(k),\quad m\le u+h[B(u)+B(m)] .
\] Use the limiting interpretation at zero. These inequalities give increasing concave upper bounds on \(m,k\) as functions of \(h\). They are the inverses of \((s-u)/(B(u)+B(s))\) on \(s\ge u\) and \(s/B(s)\) on \(s\ge0\). To prove convexity before inversion, the integral for \(F\) gives \(2F'+sF''\le0\), hence \(2B'+sB''\le0\). For \(g=1/(c+B)\), \(c\ge0\), differentiation yields \(g'\ge0,\ 2g'+sg''\ge0\). Thus \(2g'+(s-u)g''\ge0\) on \(s\ge u\), which proves the claim for either inverse. For any nonnegative increasing concave function \(P\), its supporting affine line at \(\|h\|_2>0\) has nonnegative slope and intercept. The \(L^2\) triangle inequality then gives \(\|P(h)\|_2\le P(\|h\|_2)\). The case \(\|h\|_2=0\) is immediate.

Set \(h_0=\sqrt{49J/384}\), \(u=3/4000\). For any positive trial \(x\), monotonicity of \(B\) bounds these two inverses at \(h_0\) by \(\max\{x,u+h_0[B(u)+B(x)]\}\) and \(\max\{x,h_0B(x)\}\), respectively. Minimize each bound separately over \(x=i/4000,\ 4\le i\le92\), obtaining \(M_0,K_0\). It follows that the left side of (2) is at most \(4aM_0+6M_0^2+2K_0^2\). The same rational bound on the exponential integral applies at every grid point.

For closure, translate the spectral variable to \(\xi'=\xi+\eta\) with \(\eta>0\). The pencil coefficient becomes \(X+\eta I\) and the time integrand acquires the factor \(e^{-\eta s}\). This preserves zero-freeness on \(\Re\xi'<0\) and makes the spectral integral for \(\mathbb E_1 h_\eta^2\) continuous on \([0,1]\). Pass to the weak limiting measure at fixed \(\eta\). The resulting energy is at most \(49J/384\) since the extra time factor is bounded by one. Let \(\eta\downarrow0\) to recover the displayed test on the limiting endpoint data, even for sequences that do not preserve the endpoints.

With the same fitted-pole bounds as above, the 1,024 subintervals of \(|q-q_0|\le0.007\) give margins greater than \(0.0000965\) at deficit \(1.16\cdot10^{-5}\) for weight 16 and \(0.0000703\) at score \(6.02093\) for weight 15.

Set \(q_0=1/5,\ a=1/(Z-1),\ d=1+2a,\ \rho=1/d\). The weights, reference value and normalization are

\[
 w=\frac{(a+q_0)^2}{d^2},\qquad v=\frac{(1+a-q_0)^2}{d^2},\qquad
 C_Z=(w+v)\frac{1+1/Z}{2}-\frac{1-1/Z}{12d},\qquad
 \lambda_Z=\frac{1-1/Z}{d^3}.
\]

The normalized deficit at these two positive contrasts is

\[
 B_Z=\frac{C_Z-wH_Z-vL_Z}{\lambda_Z}
    =\frac13\int_0^1\frac{(t-q_0)^2}{1-\rho^2(2t-1)^2}\,d\nu(t).
                                                                        \tag{2}
\]

Indeed the subtracted terms in the two conductivity formulas are \((1-1/Z)/12\) times the integrals of \(1/(t+a)\) and \(1/(1+a-t)\), respectively. The weighted sum of these kernels equals \([1+(t-q_0)^2/((a+t)(1+a-t))]/d\). This proves (2), including the specified reference. For every physical function \(B_Z\) increases monotonically from \(\int(t-q_0)^2d\nu/3\) at \(Z=1+\) to \(77/300-(H+16L)/25\) at infinity.

At half volume the fourth-coefficient exclusion gives \(\int(t-q_0)^2d\nu>4(0.000947)^2\) for every unrestricted laminate limit. Thus every ratio \(Z>1\), however close to one, gives the necessary strict inequality

\[
       wH_Z+vL_Z<C_Z-\lambda_Z\,4(0.000947)^2/3 .
\]

In contrast, the proposed physical endpoint tolerance of \(10^{-7}\) would imply a deficit of at most \(\lambda_Z10^{-7}\) at every one of these pairs. The dimensional reading precision in this statement scales as \((Z-1)^4\) near one.

\paragraph*{A shifted two-line estimate at ratio 25}

The two-reading laminate deficit at ratio 25 is \(3.4\cdot10^{-6}\), compared with \(1.8\cdot10^{-6}\) from the single vertical-line estimate. It applies to arbitrary laminate directions, depth and whole-function limits. Let \(a_0=1/24\), \(d\mu=d\nu/[(a_0+t)(1+a_0-t)]\), with fitted mean \(q\), variance integral \(J\), and total mass \(K\) as above. Then \[
 K(a_0+q)(1+a_0-q)=1+J,\qquad
 \bar t-q=(1-2q)J-\int(t-q)^3d\mu\ge-qJ .
\] The last inequality uses \(t-q\le1-q\). The finite readings give \(J+K(q-1/5)^2\le I=12(1+2a_0)(3.4\cdot10^{-6})/(1-1/25)\) if the displayed ceiling is reached.

Translate the spectral variable by \(\tau=3/1000\). In the negative half-plane the leaf determinant is still zero-free, the symmetric pencil coefficient becomes \(X+\tau I\), and the time integrand is multiplied by \(e^{-\tau s}\). With \(A=3q-1/2\), the trace of the shifted coefficient is \(A+3\tau\), and its expected squared Frobenius norm minus squared trace is at least \[
        A(1/2-2q)-qJ/2-4A\tau-6\tau^2 .
\] Apply the coupled midpoint and half-gap bounds with time weight \(s+6\) and vertical offset \(u=7/10000\) for the shifted pencil. Its squared time norm has conditional mean at most \(21J/40\), since on \(0\le t\le1\), \[
 \frac{t(a_0+t)(1+a_0-t)[1+12(t+\tau)]}{4(t+\tau)^2}<21/20 .
\] The verifier checks positivity of the cleared quartic by exact root count. All spectral integrands are continuous at both endpoints because \(\tau>0\). At each fitted-pole interval the preceding lower bound exceeds \(4(A+3\tau)\|m\|_2+6\|m\|_2^2+2\|k\|_2^2\) from the two-line estimate. The least certified difference exceeds 0.000039. Thus two same-specimen \textbf{ordinary conductivity} lower bounds at ratio 25 can certify separation at the stated cutoff. Whether a physical pair reaches it remains open.

In the final scalar tests the rational shift is $3/1000$ and the
vertical offset is $7/10000$. The verifier clears the positive
denominators, checks the endpoint values and uses a Sturm chain or
rational Bernstein coefficients on each closed subinterval. Its
ratio-25 conclusion is
\[
  w H_{25}+v L_{25}<4321/16250-34/10^7=0.26590429230769\ldots .
\]
The comparison tree has 42 nontrivial splits. Starting with the two
phase conductivities, apply the series rule on the normal component
and continuity of tangential voltage and normal current at each
specified normal (the Schur-complement formula of Problem~4).
At each step the verifier encloses the $3\times3$ matrix by rounding
each entry to a rational grid and adding or subtracting three grid
units on the diagonal. The row-sum norm bounds this rounding error.
Monotonicity of the laminate cell map propagates the rational upper
and lower matrices. Their normalized traces at the two reciprocal
contrasts give the score interval in the article.
\phantomsection\label{endfile:two-flux-proof}

\SSection{arc}{Arc constructions and laminate exclusions}
Throughout this section $f$ is the fraction of the phase with
relative permittivity $1+s$. The symbol $r$ is also used for the pole
$u$. Recipes using the other phase as the inclusion change $(f,r)$ to
$(1-f,1-r)$. The ordinary coated endpoints are $(1-f)/3$ and $1-f/3$.
\subsection{Matching plies and extending their poles}\label{route:matching-plies-and-extending-their-poles}
\phantomsection\label{file:band-construction}
\paragraph*{Parameters and physical layer shares}

Fix \(r_0(f)\le r<r_1(f)\) and put

\[
 h=f+3r-1,\qquad t=1-h,\qquad b=h-t,\qquad w=h-r.
\]

We will use \(0<b<w<r<t<h<1\). To see the ordering, let

\[
 P_f(r)=18r^2+(9f-23)r+f^2-6f+7 .
\]

Its smaller root is \(r_0(f)\) and \(P_f(1/2)=f(f-3/2)<0\). At \((1-f)/2\) and \(1/2-f/3\) the values are \(f(f+1)\) and \(f/6\), respectively. Both arguments are below \(1/2\), so \(r\ge r_0\) implies \(w>0\) and \(h>1/2\). The strict upper bounds \(r<(2-f)/4\) and \(r<1-f\) give \(t>r\) and \(w<r\). Now \(w-b=t-r>0\). The interval has positive length throughout the stated fraction range. Indeed \(P_f((2-f)/4)=-f(f+2)/8<0\). For \(f\le1/2\) we already have \(1-f\ge1/2\). For \(1/2<f<f_c\) we have \(P_f(1-f)=2(5f^2-5f+1)<0\).

We need one more pair of pole parameters \(\ell,k\), with \(k=h-\ell\). Set \(A=w(t-r)\) and choose

\[
 \ell=\frac{\sqrt{A^2+2(1-r)hA}-A}{2(1-r)} .
 \tag{4}
\]

This is the unique positive root of \(H(\ell)=2(1-r)\ell^2+2A\ell-hA\). The polynomial increases on the positive half-line, \(H(0)<0\), and direct expansion gives

\[
 H(b)=-(w+b)P_f(r)\ge0 .
\]

It follows that \(0<\ell\le b\) and \(t\le k<h\). In particular \(\ell<r<k\). The identity \(H(\ell)=0\) can be written in the form that determines the geometry,

\[
 \eta h w=2(1-r)\ell,\qquad
 \eta=\frac{t-r+2\ell}{w+\ell}\in(0,1] .
 \tag{5}
\]

The upper bound on \(\eta\) is exactly \(\ell\le b\). At \(r=r_0\) it equals one.

\paragraph*{Geometry}

Let \(S=L_z(h)\) and \(R=L_x(1-r)\) be simple laminates of the pure phases. The subscript gives the layer normal and the argument gives the phase-1 fraction. Write \(P_2\) for pure phase 2 and \(E'\) for a quarter-turn of \(E\) about \(x\). Each row gives one stack, with the second material filling the balance of the layers.

{\def\LTcaptype{none} 
\begin{longtable}[]{@{}
  >{\raggedright\arraybackslash}p{(\linewidth - 8\tabcolsep) * \real{0.2000}}
  >{\raggedright\arraybackslash}p{(\linewidth - 8\tabcolsep) * \real{0.2000}}
  >{\raggedright\arraybackslash}p{(\linewidth - 8\tabcolsep) * \real{0.2000}}
  >{\raggedright\arraybackslash}p{(\linewidth - 8\tabcolsep) * \real{0.2000}}
  >{\raggedright\arraybackslash}p{(\linewidth - 8\tabcolsep) * \real{0.2000}}@{}}
\toprule\noalign{}
\begin{minipage}[b]{\linewidth}\raggedright
Result
\end{minipage} & \begin{minipage}[b]{\linewidth}\raggedright
Layer normal
\end{minipage} & \begin{minipage}[b]{\linewidth}\raggedright
First material
\end{minipage} & \begin{minipage}[b]{\linewidth}\raggedright
First share
\end{minipage} & \begin{minipage}[b]{\linewidth}\raggedright
Second material
\end{minipage} \\
\midrule\noalign{}
\endhead
\bottomrule\noalign{}
\endlastfoot
\(C\) & \(x\) & \(S\) & \(w/h\) & \(P_2\) \\
\(D\) & \(y\) & \(S\) & \(\ell/r\) & \(C\) \\
\(E\) & \(z\) & \(D\) & \(\eta\) & \(P_2\) \\
\(B\) & \(x\) & \(E\) & \(1/2\) & \(E'\) \\
\(V\) & \(y\) & \(R\) & \((k-r)/(k-\ell)\) & \(B\) \\
\(Z\) & \(z\) & \(R\) & \((k-r)/(k-w)\) & \(V\) \\
\end{longtable}
}

All shares are in \([0,1]\) by the ordering and (5). The first three rows make a material with two equal pole parameters. Averaging it with its quarter-turn in row \(B\) sets two transverse numerator coefficients to the value needed for matching. That allows the \textbf{same simple laminate} \(R\) to serve in both final rows.

\paragraph*{Tensor calculation}

For a diagonal intermediate material of phase-1 fraction \(c\), describe its entries by

\[
 k_j(s)=\frac{1+(c+p_j)s+q_js^2}{1+p_js},\qquad
 c=1+h-\sum_jp_j,\qquad \sum_jq_j=ch .
 \tag{6}
\]

For an \(i\)-normal stack, if the two tangential \(p\) components and \(q_i\) agree across the plies, both triples average with the actual shares. Indeed \(c+p_i\) is common, so the harmonic mean along \(i\) averages denominators over a common numerator. Tangential arithmetic means average numerators over common denominators. A simple laminate of fraction \(c\) can use \(p_i=1-c\), \(q_i=0\) on its normal and \(q_j=cp_j\) in tangent directions with \(\sum_{j\ne i}p_j=h\). These tangential entries are affine. Pure phase 2 can use any \(p\) of sum \(1+h\) with \(q=0\).

The needed intermediate quantities, in \(x,y,z\) order, are

{\def\LTcaptype{none} 
\begin{longtable}[]{@{}lll@{}}
\toprule\noalign{}
Material & \(p\) & \(q\) \\
\midrule\noalign{}
\endhead
\bottomrule\noalign{}
\endlastfoot
\(C\) & \((r,h,t)\) & \((0,hw,0)\) \\
\(D\) & \((r,k,t)\) & \((h\ell,hw,0)\) \\
\(E\) & \((r,k,k)\) & \((\eta h\ell,\eta hw,0)\) \\
\(B\) & \((r,k,k)\) & \((\eta h\ell,(1-r)\ell,(1-r)\ell)\) \\
\end{longtable}
}

For \(C\), write the poles of \(S\) and \(P_2\) as \((0,h,t)\) and \((h,h,t)\). For \(D\) the same \(S\) can use \((r,w,t)\), so its \(q_y=hw\) agrees with \(C\). For \(E\) give \(P_2\) the poles \((r,k,1+w-k)\). The share that averages the normal pole to \(k\) is \((1+w-2k)/(1+w-k-t)=\eta\). Row \(B\) then follows by (5).

In row \(V\), represent \(R\) by \(p=(r,\ell,k)\), \(q=(0,(1-r)\ell,(1-r)k)\). Its \(y\) coefficient matches \(B\). Writing \(\lambda=(k-r)/(k-\ell)\) gives

\[
 p^V=(r,r,k),\qquad
 q_z^V=(1-r)[\lambda k+(1-\lambda)\ell]
       =(1-r)(k+\ell-r)=(1-r)w .
\]

For row \(Z\), write \(R\) with poles \((r,r,w)\). Its \(q_z=(1-r)w\) now matches \(V\). Averaging the normal poles gives \(p^Z=(r,r,r)\), so (6) gives phase fraction \(f\) and \(\sum q_j^Z=fh\). Consequently

\[
 \frac{\operatorname{Tr} Z(s)}3
 =\frac{1+(f+r)s+(fh/3)s^2}{1+rs}=K_f(r,s).
\]

Two further laminations make the tensor isotropic for all contrasts. Exchange the axes carrying the smallest and largest \(q_j^Z\) by a quarter-turn, and stack the two copies normal to the remaining axis. Choose the share putting \(fh/3\) on one of the exchanged axes. Then laminate equal shares of that result and its quarter-turn about the axis already at \(fh/3\), normal to the same axis. Each normal entry is common to its two plies and the other entries average. The final tensor is \(K_f(r,s)I_3\). Counting the simple laminates, \(Z\) takes at most seven stages and the scalar tensor at most nine.

For explicit final shares one may split each ply of \(Z\) equally between itself and its quarter-turn about \(z\) within that same stack. This gives \(Z=\operatorname{diag}(T,T,N)\). Stack its \(x\)-axial and \(y\)-axial copies normal to \(z\) with shares \(1/3,2/3\), then stack the result and its quarter-turn about \(x\) normal to \(x\) in halves. The implementation uses this version. Run \texttt{python\ verify\_band.py} to check the general identities and compute all three contrast functions directly from the pure leaves at seven examples. Run \texttt{python\ verify\_band.py\ -\/-tree\ 1/4\ 17/40} to export a geometry, or \texttt{python\ verify\_band.py\ -\/-tree\ -\/-lower\ 1/4} for the algebraic lower endpoint.

At \(r=r_1\) with \(f<2/3\), set \(t=r\), \(\ell=0\), \(k=h\), and \(\eta=0\). Then \(B=P_2\), and the last two rows use only \(R\) and pure phase 2. They give LB1 by the same calculation. At \(f\ge2/3\) the right endpoint is \(r=1-f\), realized by isotropizing a simple laminate of fraction \(f\). This covers the meeting case \(f=2/3\).

\paragraph*{What this adds at one quarter}

The earlier parameters outside (1) include the coated-sphere endpoints \(1/4,11/12\), both coated-cylinder points \(3/8,7/8\), and the two values \((39-\sqrt{177})/64,27/32\) from (K) and its phase interchange. Inside (1), the previously supplied points include \(r_0=(83-\sqrt{481})/144\), every point of \([7/16,7/12]\), and the simple-laminate point \(3/4\).

{\def\LTcaptype{none} 
\begin{longtable}[]{@{}
  >{\raggedright\arraybackslash}p{(\linewidth - 2\tabcolsep) * \real{0.5000}}
  >{\raggedright\arraybackslash}p{(\linewidth - 2\tabcolsep) * \real{0.5000}}@{}}
\toprule\noalign{}
\begin{minipage}[b]{\linewidth}\raggedright
Part of the quarter-fraction result
\end{minipage} & \begin{minipage}[b]{\linewidth}\raggedright
Realization proved here or in the earlier data
\end{minipage} \\
\midrule\noalign{}
\endhead
\bottomrule\noalign{}
\endlastfoot
\([(5-\sqrt3)/8,5/12]\) & Fixed-fraction limit of finite hierarchical laminates, possibly with increasing depth \\
\((5/12,r_0)\) & Finitely many applications of the three-stack recursion \\
\([r_0,7/16]\) & At most nine stages by the layer table above, with an eight-stage subinterval \\
\([7/16,7/12]\) & Supplied finite laminates including LB1 and LA2 \\
\([3/4,(5+\sqrt3)/8]\) & Phase interchange of the increasing-depth construction at \(f=3/4\) \\
\([0.8220201246\ldots,1-\rho]\) & Two overlapping nine-stage constructions, including the extension in (1) \\
\end{longtable}
}

On \([(5-\sqrt3)/8,5/12]\), including \(u=5/12\) where \(h=1/2\), the proof of full coverage uses limits of laminates of unbounded depth. On \([3/4,(5+\sqrt3)/8)\) it likewise gives an unbounded-depth limit at the points outside the finite constructions listed here. Whether each remaining limiting point also has a finite-depth realization is open.

Additional finite trees give \(u=37/48,19/24,13/16\) in the upper band. Their exact normals and shares are in the example file, which records shared subtrees only once. A separate tree gives \(u=5/6\) by isotropizing its trace and switching its phase labels.

The interval comparison includes the earlier below-LB1 recipes. For (F), its parameter \(x_F=(1+t-2f)/[3(3t-1)]\) increases as \(f\) decreases. Already at \(f=2/5\) its required \(R x_F\le Q_F\) fails throughout \(1/\sqrt7<t<1/\sqrt5\), since the difference is \(4(4t-1)(1-5t^2)/[5(3t-1)]>0\). The earlier (C) lower ply uses a transverse split requiring \((3t-1)^2x_F\le1-5t^2\) whenever its \(U_x\) share is positive. At \(f=1/4\) the left side minus the right side is \((3t-1)(t+1/2)/3-1+5t^2\), increasing for \(t\ge7/16\) and already \(7/128\) at \(7/16\). A zero \(U_x\) share in the open gap would require equality of the two normal coefficients of the unsplit plies, \(2t^2(t-r)/(t-b)=(1-r)w\), which reduces to \(P_f(r)=0\). The earlier phase interchanges of below-LB1 families lie above \(1-f\) when available.

For a concrete example inside the new interval, at \(f=1/4\), \(u=3/7\) the effective tensor is

\[
 \frac{\varepsilon_*}{\varepsilon_2}
 =\frac{5z^2+66z+41}{16(3z+4)}I_3,\qquad z=\frac{\varepsilon_1}{\varepsilon_2}.
\]

On some subintervals there are shorter recipes. An eight-stage construction gives rational shares at this example. At \(f=2/5\) a seven-stage construction reaches from \(17/45\) to \(0.3869778027\ldots\). These are upper bounds on the depths of the displayed geometries.

At \(f=1/4\) any finite laminate realization of an \(r\in[R(1/4),7/16)\) or \(r\in(3/4,1-\rho]\) requires \textbf{at least five stages among hierarchical laminates of the two pure phases}, with arbitrary layer normals and any finite number of plies per stack. The earlier four-stage classification already excludes these values in normalized trace. At this fraction its unlisted trace contacts come from (A), (K), or the phase interchange of (K), all outside those two ranges.
\phantomsection\label{endfile:band-construction}

\phantomsection\label{file:band-recurrence}
\paragraph*{\texorpdfstring{Continuing the laminate construction to \(h=1/2\)}{Continuing the laminate construction to h=1/2}}

For every \(0<f<3/4\), isotropic two-phase geometries attain the closed interval

\[
\boxed{\quad \frac12-\frac f3\ \le u\le\min\left\{\frac{2-f}{4},1-f\right\}.\quad}
\tag{1}
\]

Every point above the left endpoint has a finite hierarchical laminate realization on separated scales. The proof at the left endpoint uses a limit at the same phase fraction. Interchanging phases gives the reflected interval at \(1-f\). At \(f=1/4\), (1) is \([5/12,7/16]\). The earlier LB1-to-LA2 band continues from \(7/16\) to \(7/12\), giving \([5/12,7/12]\).

\paragraph*{A coefficient preserved by the construction}

Keep \(h\) fixed in \(1/2<h<2/3\), and set \(t=1-h\). For a target pole \(r\), its phase fraction is \(f=1+h-3r\). We construct every \(h/2<r\le t\). At \(r=h/2\) the isotropized simple laminate has \(u=1-f=r\).

Use the diagonal tensor representation of the general band. A material of phase-1 fraction \(c\) has entries

\[
\frac{1+(c+p_i)s+q_i s^2}{1+p_i s},\qquad
c=1+h-\sum_i p_i,\qquad \sum_i q_i=ch.
\]

In an \(i\)-normal stack, common tangential poles and a common normal coefficient \(q_i\) permit both triples to be averaged with the ply shares. A simple laminate \(L_i(c)\) has \(p_i=1-c\), \(q_i=0\), and tangential \(q_j=cp_j\) with the two tangential poles summing to \(h\). These formulas follow from harmonic normal and arithmetic tangential averaging.

The material passed from one step to the next has \(p=(k,k,k)\) and, after choosing its orientation, one specified numerator coefficient

\[
q_x=N_k=(1-k)(h-k).
\tag{2}
\]

The seven-stage trace construction in the general-band proof provides this material directly. Its last stack has \(q_z=N_k\), so exchanging axes gives (2). To identify its range at fixed \(h\), put

\[
T(k)=t-\frac{(1-k)(h-k)}{2k-h},\qquad h/2<k\le t .
\tag{3}
\]

The seed proof's sign condition \(P_f(k)\le0\) becomes \((2-3h)[T(t)-k]\le0\) upon substituting \(f=1+h-3k\). For \(h/2<k<t\) its other strict order requirements follow from \(0<2h-1<h-k<k<t<h\). Thus all \(h/2<k\le t\) with \(k\ge T(t)\) have (2), including \(k=t\) by the explicit LB1 stack.

\paragraph*{Three further stacks}

Suppose (2) is available and \(r=T(k)>h/2\). In this step \(r<k\) as proved below. Let \(\ell=h-k\), \(w=h-r\), and put

\[
a=\frac{r-\ell}{k-\ell},\qquad
\lambda=\frac{k-r}{k-\ell},\qquad
\beta=\frac{k-r}{k-w}.
\]

Here \(\ell<w<h/2<r<k\), so all the shares lie strictly between zero and one. Use the following three stacks.

\begin{enumerate}
\def\labelenumi{\arabic{enumi}.}
\item
  \textbf{Normal to \(x\).} Stack the material (2), its quarter-turn about \(x\), \(L_y(1-k)\), and the corresponding \(L_z(1-k)\), in shares \(a/2,a/2,(1-a)/2,(1-a)/2\). The simple plies have \(p=(\ell,k,k)\) and \(q_x=(1-k)\ell=N_k\). The result \(B\) has \(p=(r,k,k)\) and \(q_x=N_k\). The two rotated pairs give equal transverse coefficients. Its phase fraction is \(c_B=t-r+2\ell\). Equation (3) gives \(N_k=(t-r)(2k-h)\), hence

  \[q_y^B=q_z^B=(c_Bh-N_k)/2=(1-r)\ell.\]
\item
  \textbf{Normal to \(y\).} Stack \(R=L_x(1-r)\) in share \(\lambda\) with \(B\) in share \(1-\lambda\). Represent \(R\) by \(p=(r,\ell,k)\), so \(q_y^R=(1-r)\ell\) matches \(B\). The result \(V\) has \(p=(r,r,k)\) and

  \[q_z^V=(1-r)[\lambda k+(1-\lambda)\ell]=(1-r)(h-r)=N_r.\]
\item
  \textbf{Normal to \(z\).} Stack the same \(R\) in share \(\beta\) with \(V\) in share \(1-\beta\), this time representing \(R\) by \(p=(r,r,w)\). Both normal coefficients are \(N_r\). The result has \(p=(r,r,r)\) and \(q_z=N_r\). A rotation puts the specified coefficient on \(x\) for a subsequent step.
\end{enumerate}

Thus the three stacks produce the same type of material at \(r\). Only the finished root needs isotropization.

\paragraph*{Coverage and fixed-fraction closure}

In (3) the positive numerator \((1-k)(h-k)\) decreases strictly and the positive denominator increases strictly. Thus \(T\) increases from \(-\infty\) to \(T(t)\). Moreover

\[
k-T(k)=\frac{h(2-h)-3k(1-k)}{2k-h}
\ \ge\ \frac{h(2h-1)}{2-3h}>0 .
\tag{4}
\]

Indeed \(k\le t<1/2\) gives \(k(1-k)\le t(1-t)\). For any target \(r<T(t)\) in \((h/2,t)\), seek the unique \(k\in(h/2,t)\) satisfying \(T(k)=r\). Equation (4) increases the pole by at least its displayed positive constant at each such inverse step. After finitely many steps the pole belongs to the seed interval. Construct the seed and perform the three stacks in reverse order down the list of poles. This gives \(p=(r,r,r)\) at the desired fraction.

The root entries have a common denominator and common linear numerator term. Their quadratic coefficients average to \(m=fh/3\). If their extrema are \(d<e\), stack the root with the copy exchanging those two directions, normal to the remaining axis. Give the orientation with coefficient \(d\) on the chosen axis share \((e-m)/(e-d)\), so that axis acquires \(m\). In a second stack, normal to that axis, give equal shares to this material and its copy exchanging the other two axes. All coefficients now equal \(m\). Coincident extrema already give a scalar tensor. The resulting effective conductivity in normalized units is

\[
\frac{1+(f+r)s+(fh/3)s^2}{1+rs}\,I
=\left(1+fs-\frac{f(1-f)s^2}{3(1+rs)}\right)I.
\]

For a quarter-fraction example below the two-cylinder curve, take the seed \(f_k=51/488\), \(k=115/244\). It has \(h=253/488\), and one step gives \(f=1/4\), \(r=619/1464\). The command \texttt{python\ verify\_recursion.py\ -\/-seed\ 51/488\ 115/244\ -\/-steps\ 1\ -\/-tree} evaluates the full three-entry contrast function exactly and exports the complete isotropic geometry.

At fixed \(f\), the strict inequalities \(1/2<h<2/3\) and \(h/2<r<t\) hold when \(1/2-f/3<r<\min\{(2-f)/4,1-f\}\). This interval has positive width for \(0<f<3/4\). The upper endpoint is LB1 for \(f\le2/3\) or the isotropized simple laminate for \(f\ge2/3\). For the lower endpoint, take a sequence in the strict interval tending down to \(1/2-f/3\), keeping \(f\) fixed throughout. Each finite laminate is a reiterated limit of periodic cells of that fraction. A diagonal sequence of these cells realizes the limiting contrast function locally uniformly on the cut contrast plane. The fixed-fraction approximation and rotated-coercivity argument are given in the general-band proof.

\paragraph*{\texorpdfstring{A realized interval for \(0<h<1/2\)}{A realized interval for 0\textless h\textless1/2}}

Write \(r=u\), \(h=f+3r-1\), \(t=1-h\), and define \(R(f)\) by

\[
R(f)=\frac{3-2f-\sqrt{1-4f(1-f)/3}}{4}.
\]

Together with the \(h\ge1/2\) construction and the earlier family (A), the constructions below give the closed intervals

\[
\boxed{\left[R(f),\ \min\{(2-f)/3,1-f\}\right]},\qquad
\boxed{\left[\max\{(2-f)/3,1-f\},\ 1-R(1-f)\right]}
\tag{1}
\]

for every \(0<f<1\). In particular, at \(f=1/4\) they give \([(5-\sqrt3)/8,7/12]\) and \([3/4,(5+\sqrt3)/8]\). The low-\(h\) interval uses limits of finite hierarchical laminates of increasing depth. At each fixed \(f\), infinitely many distinct finite-laminate contacts approach \(R(f)\) from below.

\paragraph*{Two specified coefficients at the starting poles}

Fix \(0<h<1/2\) throughout the construction. Set

\[
f_k=1+h-3k,\quad N_k=(1-k)(h-k),\quad D_k=(t-k)(2k-h),\quad
T(k)=t-\frac{N_k}{2k-h},\qquad h/2<k\le h.
\]

Use the diagonal pole and quadratic numerator triples \(p,q\) of the laminate rule. Thus \(c=1+h-\sum_i p_i\), \(\sum_i q_i=ch\), and equal tangential \(p\) and equal normal \(q\) permit arithmetic averaging of the triples. We first construct materials with

\[
p=(k,k,k),\qquad q=(N_k,N_k,D_k).
\tag{2}
\]

Here \(2N_k+D_k=f_k h\). At \(k=h\), take the simple laminate \(L_x(t)\) of phase-1 fraction \(t\), using poles \((h,0,h)\) and coefficients \((0,0,th)\). Stack it normal to \(y\) in share \((1-2h)/t\) with pure phase 2, represented at poles \((h,t,h)\). The result has (2).

Suppose (2) is available at \(k\) and \(h/2<a<k\) satisfies \(T(a)=k\). There is a three-stage construction of (2) at \(a\). First keep \(p_z=k\). At pole coordinates \((x,y,k)\) use the coefficients

\[
q_x=(1-k)(h-y),\qquad q_y=(1-k)(h-x),\qquad q_z=(t-k)(x+y-h).
\tag{3}
\]

At \(x=y=k\) this is (2). At \(x+y=h\) it is \(L_z(1-k)\). Starting at \((k,k,k)\), stack normal to \(x\) with this simple laminate at \((h-k,k,k)\), retaining a share \((a-h+k)/(2k-h)\) of (2). This gives \((a,k,k)\). Then stack normal to \(y\) with the simple laminate at \((a,h-a,k)\), retaining a share \((2a-h)/(k-h+a)\) of the first result. All shares are in \((0,1)\). The normal coefficients in (3) are constant along each averaging direction and the other coefficients are affine, verifying (3) at both stages. The final state \((a,a,k)\) has \(q_x=q_y=(1-k)(h-a)\) and \(q_z=(t-k)(2a-h)=N_a\).

For the third stack put \(w=h-a\). The simple laminates \(L_x(1-a),L_y(1-a)\) both have poles \((a,a,w)\) and \(q_z=N_a\). Give the material at \((a,a,k)\) a share \((a-w)/(k-w)\) in a \(z\)-normal stack, so that the averaged normal pole is \(a\). Divide the remaining share between the two simple laminates, giving \(L_y\) the proportion \(w(1-w)/[a(1-a)]\). This belongs to \((0,1)\) since \(0<w<a<1/2\). These simple plies together supply \(q_x=w(1-w)\). Its average with \(q_x=(1-k)w\) from the upper ply is \(q_x=(1-a)w=N_a\). The entire stack therefore has \(q_x=q_z=N_a\) and \(q_y=D_a\). Relabeling axes gives (2).

The unique fixed point of \(T\) on \((h/2,h)\) is

\[
\kappa=\frac{1-\sqrt{1-4h(2-h)/3}}2,\qquad 3\kappa(1-\kappa)=h(2-h).
\tag{4}
\]

Indeed \(T\) increases strictly from \(-\infty\) to \(t>h\), and the sign of \(k-T(k)\) is the sign of \(h(2-h)-3k(1-k)\). This expression decreases strictly on \((h/2,h)\), with endpoint values \(h(2-h)/4>0\) and \(h(2h-1)<0\). Thus, starting with \(b_0=h\), the sequence \(b_{j+1}=T^{-1}(b_j)\) decreases strictly to \(\kappa\). All its terms have the finite construction (2). At the limit \(D_\kappa=N_\kappa\), so the actual tensors in (2) converge to the physical scalar \(K_{f_\kappa}(\kappa,s)I\), where

\[
K_f(r,s)=1+fs-\frac{f(1-f)s^2}{3(1+rs)} .
\]

For realization at exactly \(f_\kappa\), add to each finite material a share \((f_\kappa-f_{b_j})/(1-f_{b_j})\) of pure phase 1 by a coordinate stack. This share tends to zero and its tensor has the same limit. At the other endpoint \(r=h/2\), \(f_r=1-r\) and isotropization of a simple laminate gives \(K_{f_r}(r,s)I\).

\paragraph*{\texorpdfstring{Every pole between \(h/2\) and \(\kappa\)}{Every pole between h/2 and \textbackslash kappa}}

Fix \(h/2<r<\kappa\). Set \(a_0=r\) and \(a_{j+1}=T^{-1}(a_j)\). Now \(a_j\) increases strictly to \(\kappa\). The three-stack rule for \(r=T(k)<k\) needs a material at \(k\) with just one specified coefficient \(q_x=N_k\). Its shares depend only on \(h,r,k\). It applies verbatim to consecutive \(a_j,a_{j+1}\) because \(h-a_{j+1}<h/2<a_j<a_{j+1}\), and it passes on a coefficient \(N_{a_j}\) after a rotation.

For a truncation at depth \(n\), first calculate its ideal tensor using a terminal diagonal tensor with \(p=(a_n,a_n,a_n)\) and

\[
q_x=N_{a_n},\qquad q_y=q_z=(f_{a_n}h-N_{a_n})/2.
\]

Apply the three-stack rule at \(j=n-1,\dots,0\), and apply the two isotropizing stacks at the finished root. This algebra produces exactly \(K_{f_r}(r,s)I\) for every \(n\). The physical geometry replaces only the terminal tensor in this calculation. Use the finite upper-sequence material (2) at \(b_n\), with a layer share

\[
\theta_n=\frac{f_{a_n}-f_{b_n}}{1-f_{b_n}}
\]

of pure phase 1 in a coordinate stack, giving the material at \(b_n\) share \(1-\theta_n\). We have \(0<\theta_n<1\), \(\theta_n\to0\), and the resulting terminal fraction is exactly \(f_{a_n}\). The full physical tree consequently has fraction \(f_r\) for each \(n\).

One can also make each physical truncation isotropic. In the last \(z\) stack split every ply equally between itself and its quarter-turn about \(z\). The resulting tensor has equal tangential entries, including after the terminal replacement. The fixed \(1/3,2/3\) and half-half Schulgasser stacks in the main solution make it scalar. For the ideal truncation these same operations preserve the normalized trace \(K_{f_r}(r,s)\).

For any fixed positive real phase contrast \(1+s\), write \(A_n\) for the ideal terminal tensor and \(B_n\) for its physical replacement. Both tend to \(K_{f_\kappa}(\kappa,s)I\), which is strictly positive. For the ideal tensor this follows by inserting \(a_n\to\kappa\) and \(D_\kappa=N_\kappa\) in its entries. For the physical one it follows from (2), \(b_n\to\kappa\), and \(\theta_n\to0\). For all sufficiently large \(n\) there is therefore \(\epsilon_n\to0\) with

\[
(1-\epsilon_n)A_n\ \preceq B_n\ \preceq(1+\epsilon_n)A_n .
\]

The same relative comparison holds at every parent up to the isotropic ideal root. A laminate is monotone in the positive-definite constituent tensors by the electric energy principle, and scaling every constituent scales the result. At each unchanged sibling the two-sided comparison holds as well, and rotations preserve order. This argument propagates the original factors \(1\pm\epsilon_n\) through arbitrarily many stages. The physical tensors converge to \(K_{f_r}(r,s)I\) for every positive real contrast. Taking a diagonal sequence of periodic approximants at exact fraction \(f_r\) gives locally uniform convergence in the complex contrast plane cut along the nonpositive real axis.

For a numerical illustration, take \(f=1/4\), \(r=0.41\) and evaluate the physical finite truncations at phase ratios \(0.1,2,10,1.1+2i\). The largest relative component errors at \(n=1,3,6\) are \(0.09392\), \(0.03151\) and \(0.002766\). At \(f=3/4\), \(r=0.2\) the corresponding errors are \(0.01477\), \(0.0002031\) and \(2.422\times10^{-7}\). These are finite-depth homogenized tensors.

\paragraph*{The section at fixed phase fraction}

For \(h=f+3r-1\), condition \(r\le\kappa\) in the range \(0<h<1/2\), \(r<1/2\), is equivalent to

\[
Q_f(r)=h(2-h)-3r(1-r)=-6r^2+(9-6f)r-(f-1)(f-3)\ \ge0.
\]

Its smaller root is \(R(f)\). At \(r=(1-f)/3\), \(Q_f\) is negative. At \(1/2-f/3\) and \(1-f\), its values are \(f^2/3\) and \(f(1-f)\), both positive. Consequently the low-\(h\) proof yields every \(r\) from \(R(f)\) up to \(\min\{1/2-f/3,1-f\}\), with the \(h=1/2\) endpoint obtained by closure at fixed \(f\) whenever it occurs. If \(f<3/4\) the higher-\(h\) continuation reaches \(\min\{(2-f)/4,1-f\}\). For \(f\ge3/4\) the low-\(h\) part already reaches \(1-f\). Adjoining the earlier band (A) when \(f<2/3\) gives the first interval (1). Substituting \((1-f,1-r)\) gives the second.

For the finite contacts, write \(b_n(h)\) for the upper sequence from (2). Each \(b_n\) depends continuously on \(h\) because \(T\) is continuous in both variables and strictly increasing in its pole parameter. The bounds \(\kappa(h)<b_n(h)\le h\) give endpoint limits \(0\) and \(1/2\) as \(h\) runs from \(0\) to \(1/2\), uniformly in \(n\). Thus \(1+h-3b_n(h)\) assumes every \(f\in(0,1)\). For that \(f\), select \(h_n\) with \(1+h_n-3b_n(h_n)=f\) for each \(n\ge1\), and take \(r_n=b_n(h_n)\). The finite materials (2) have common poles and their two final Schulgasser stacks give isotropy. A coincidence of \(r_n\) for two indices would fix the same \(h_n=f+3r_n-1\), contradicting the strict decrease of \(b_n(h)\) with \(n\) there.

The same uniform endpoint bounds keep the \(h_n\) in a compact subinterval of \((0,1/2)\). On that interval \(b_n(h)\) decreases to the continuous \(\kappa(h)\) uniformly by Dini's theorem. Hence \(Q_f(r_n)\to0\). Every \(r_n\) satisfies \(r_n<1/2-f/3\) and \(Q_f(r_n)<0\) since \(r_n>\kappa(h_n)\). This places it below the smaller root \(R(f)\), and thus \(r_n\to R(f)\) from below. Phase interchange gives finite contacts approaching \(1-R(1-f)\) from above.
\phantomsection\label{endfile:band-recurrence}

\subsection{Connecting the upper quarter-fraction interval}\label{route:connecting-the-upper-quarter-fraction-interval}
\phantomsection\label{file:upper-bridge}
\paragraph*{The two plies used in both recursions}

All triples below use \(x,y,z\) order. A ply of phase-1 fraction \(c\), poles \(p\) and numerator coefficients \(q\) has

\[
 K_i=\frac{1+(c+p_i)s+q_i s^2}{1+p_i s},\qquad
 c=1+h-\sum p_i,\quad \sum q_i=ch .
\]

An \(i\)-normal stack with matching tangential \(p\) and matching \(q_i\) averages both triples. Write \(L_i(c)\) for a simple laminate normal to \(i\), with phase-1 share \(c\). Put

\[
 h=3r-1/4,\quad t=1-h,\quad \delta=t-h,\quad w=h-r,\quad d=r-w,
 \quad N=(1-r)w,\quad A=d\delta .
\]

The upper ply \(K_z(r,r)\) has \(p=(r,r,h)\) and \(q=(tw,tw,A)\). The real lower ply, denoted \(L\), is the lower ply of the seven-stage construction, used here at variable \(r\). Put

\[
 a=\frac{r\delta}{2t-r},\quad p=h-a,\quad c_L=1+h-2r-p .
\]

Stack \(L_y(1-r)\) at poles \((a,r,p)\) with \(K_x(r,p)\) at \((h,r,p)\) normal to \(x\), in respective shares \(w/p,1-w/p\). The common normal coefficient is \((1-r)a=(p-w)\delta\). This gives the following poles and coefficients:

\[
 p^L=(r,r,p),\qquad q^L=(X,Y,L_z),\quad
 X=(1-r)a,\quad L_z=N+tw(1-w/p),\quad Y=c_Lh-X-L_z.
\]

Here \(L_z\) in a coefficient expression denotes the axial coefficient of this ply. In each recursion use tangential coefficients \(C,c_Lh-L_z-C\) instead by stacking \(L\) with its quarter-turn about \(z\), normal to \(z\), in original share \((C-Y)/(X-Y)\). The pole and \(L_z\) stay fixed. In the final pole-matching stack the lower ply has poles \((r,r,p)\) and \(q_z=A\). Stack it normal to \(z\) in share \(w/a\) with \(K_z(r,r)\) to obtain \((r,r,r)\).

\paragraph*{\texorpdfstring{First recursion, from \(r_7\) to \(s_0\)}{First recursion, from r\_7 to s\_0}}

Consider a temporary constituent \(G\) with poles \((r,r,r)\), fraction \(f\), and coefficients \((S,A,N)\) where \(S=fh-A-N\). One set of stacks recreates exactly this same state while consuming a share strictly less than one of \(G\).

\begin{enumerate}
\def\labelenumi{\arabic{enumi}.}
\tightlist
\item
  Define \(C=c_Lh-A-w(t+a)\) and \(\tau=(p-w)/d\). Build a \(z\)-normal stack with share \(\tau\) of \(G\) and the rest split between \(L_y(1-r),L_x(1-r)\), both using poles \((r,r,w)\). Give \(L_y\) the proportion \((C-\tau S)/[(1-\tau)r(1-r)]\) of the simple plies. Their common normal coefficient is \(N\). The result \(V\) has poles \((r,r,p)\), with \(q_x=C\) and \(q_z=N\).
\item
  Mix the rotated mixture of \(L\) described above with \(V\), normal to \(x\). Both have \(q_x=C\) and the same poles. Give \(L\) share \(b=(A-N)/(L_z-N)\) to put \(A\) on \(z\). Its \(y\) coefficient is \(c_Lh-A-C=w(t+a)\).
\item
  Apply the final \(z\) stack described above. The \(y\) coefficient averages to \((w/a)w(t+a)+(1-w/a)tw=N\). Thus \(q=(S,N,A)\) by its trace. Exchange \(y,z\) to recover \(G\).
\end{enumerate}

The retained share of \(G\) per cycle is \((w/a)(1-b)\tau<1\). At \(r_7\) we have \(b=1\) and the recipe is already finite.

\paragraph*{\texorpdfstring{Second recursion, from \(s_0\) to \(\eta\)}{Second recursion, from s\_0 to \textbackslash eta}}

Insert one additional real ply with poles \((r,r,k)\), where \(k<p\). Its normal coefficient will be \(B\) with \(B\le A<L_z\). Here are explicit data for that ply. Put

\[\begin{gathered}
 \Delta=A-N,\quad m=1-r+\delta,\quad c_1=2w/r-1,\\
 \gamma=\frac{-c_1\Delta+\sqrt{c_1^2\Delta^2+4mw d\Delta/r}}{2m},
 \quad k=r-\gamma,\quad l=w+\gamma .
\end{gathered}\]

Use steps 1 to 3 of the ten-stage construction with these \(k,l\). Their last ply \(M\) has poles \((r,k,k)\), \(q_y=(k-w)\delta\) and \(q_z=(1-r)l\). The equality for \(q_y\) is equivalent to the displayed quadratic for \(\gamma\). Laminate \(M\) normal to \(y\) in share \(w/l\) with \(K_y(r,k)\) at \((r,h,k)\). This gives the desired real ply with

\[
 q_z=B=N+tw(1-w/l).
\]

Mix equal shares of this result and its quarter-turn about \(z\) normal to \(z\) to equalize the tangential entries. The extra ply is also specified down to the pure phases in the verifier.

This time take the temporary constituent \(G\) at \((r,r,r)\) with \(q=((fh-B)/2,(fh-B)/2,B)\).

\begin{enumerate}
\def\labelenumi{\arabic{enumi}.}
\tightlist
\item
  Stack \(G\) in share \(\tau=(p-k)/(r-k)\) with the additional real ply, normal to \(z\). The result \(V\) has \(p=(r,r,p)\) and \(q_z=B\), and both tangential coefficients are \(C=(c_Lh-B)/2\).
\item
  Mix the rotated mixture of \(L\) with \(V\) normal to \(x\), with \(L\) share \(b=(A-B)/(L_z-B)\). Apply the final \(z\) stack. It has poles \((r,r,r)\), fraction \(f\) and \(q_z=A\).
\item
  Stack equal shares with its quarter-turn about \(z\), normal to \(z\), so that both transverse coefficients are \(a_0=(fh-A)/2\). Stack the result and its \(y,z\) exchange normal to \(x\), with original share \((B-a_0)/(A-a_0)\). This puts \(B\) on \(z\). An equal stack with its \(x,y\) exchange normal to \(z\) recreates \(G\).
\end{enumerate}

The total share of the temporary material is again \((w/a)(1-b)\tau<1\). Copy rotations preserve this total.

\paragraph*{Parameter checks and the physical limit}

Throughout the intervals \(w<p<r<h\), \(0<w<a<r\), \(\delta>0\) and \(L_z>N\). Also

\[
 A-N=\frac{(2r-1)(32r-5)}8,\qquad L_z-A=-\frac{P_7(r)}{8(10r-1)} .
\]

In the second recursion write \(\gamma_0=w\Delta/(tw-\Delta)\). We have \(tw>\Delta\). The condition \(B\le A\) is \(\gamma\le\gamma_0\). Substitution of \(\gamma_0\) into the quadratic defining \(\gamma\) gives

\[
 \frac{-3(2r-1)(8r-1)^2(12r-5)(32r-5)P_{10}(r)}
 {128r(224r^2-136r+15)^2}.
\]

The prefactor of \(P_{10}\) is positive here. Both cubics change from positive to negative exactly once in their indicated brackets, by the rational substitution \(r=(r_L+r_R z)/(1+z)\) and Descartes' rule. These identities prove \(0\le b\le1\), including the two equality cases.

The remaining guards are strict. For the first recursion they include \(Y<C<X\) and \(\tau S<C<\tau S+(1-\tau)r(1-r)\). For the second they include \(h/2<k<p\), the internal split bounds of \(M\), \(a_0<B\), \(Y<C<X\), and \(L_z>B\). The verifier encloses each expression on the full rational brackets by intervals of width at most \(10^{-6}\) in \(r\), with outward integer rounding at scale \(10^{18}\). It also checks the two cubic identities algebraically.

For the realization, repeat the cycle \(n\) times and replace all remaining temporary constituents by simple laminates of phase-1 fraction \(f\). This gives finite two-phase hierarchies with exactly the prescribed fraction. To prove convergence, assign the compatible mean fields \(I+sP\) and fluxes \(I+s(P+cI)+s^2Q\) of the recursive table, with mean electric field \((1+rs)I\) at the root. These affine and quadratic pencils specify tangentially compatible electric fields and normal-flux-compatible currents at every retained stack. Each pass through a rotated copy of \(G\) resets its mean electric pencil to \((1+rs)I\). In a replacement ply extend the terminal electric field and current constantly onto the simple laminate. The discrepancy \(J-(1+s\chi)E\) is bounded at any fixed positive phase ratio \(1+s\) and is supported on a volume at most \([(w/a)(1-b)\tau]^n\). The energy principle, or testing the difference of the cell equations with the field error, bounds the resulting \(L^2\) electric field error by \(\|J-(1+s\chi)E\|_2/\min(1,1+s)\). Thus the physical tensors converge to the stated common-pole tensor, independently of the rotations of the replaced plies.

Finally isotropize to normalized trace with rotated copies. For truncations, one can first average both tangential entries about a common normal, then use two fixed Schulgasser stacks to obtain exactly the normalized trace. The limit is \(K_f(r,s)I\). A diagonal sequence of periodic approximants gives the same analytic contrast function in the usual cut plane.

\paragraph*{A block constructed from an available pole}

Throughout let \(h=3r-1/4\), \(t=1-h\), \(w=h-r\) and \(d=r-w\). The symbols \(p_i,q_i\) denote the denominator pole and quadratic numerator coefficient of component \(i\). They obey \(c=1+h-\sum p_i,\ \sum q_i=ch\) in a material of fraction \(c\). Lamination averages both triples when the tangential \(p\) and normal \(q\) agree. We use the simple layers \(L_i(c)\) and the three-stage \(G_i\) blocks defined in the upper-ply construction.

Fix \(K>r\), write \(s_0=h-K\) and \(N_K=(1-K)s_0\). Suppose an available material has \(p=(K,K,K)\), \(q_x=N_K\) and \(q_y=Y\), with the remaining coefficient fixed by the trace sum. For \(s_0<a<K\) and \(h-a\le b\le K\) define \(B_z^K(a,b)\) as follows. First stack a share \(\alpha=(a-s_0)/(K-s_0)\) of the available material with \(L_y(1-K),L_z(1-K)\) normal to \(x\). The simple plies use \(p=(s_0,K,K)\) and have \(q_x=N_K\). Give \(L_z(1-K)\) the proportion

\[
 \frac{(1-K)(h-a)-\alpha Y}{(1-\alpha)K(1-K)}
 \tag{1}
\]

of their total share \(1-\alpha\). This gives \(p=(a,K,K)\) with \(q_x=N_K,\ q_y=(1-K)(h-a)\). Next stack this ply normal to \(y\) with \(L_z(1-K)\) at \((a,h-a,K)\), retaining share \((a+b-h)/(a-s_0)\) of the first ply. The result is

\[
 p=(a,b,K),\qquad q=\big((1-K)(h-b),\ (1-K)(h-a),\ (t-K)(a+b-h)\big).
 \tag{2}
\]

For the first coefficient just average \(N_K\) and \((1-K)a\). The second was matched, and the third follows from the trace. Use permutations of (2) to orient the pole \(K\) along another axis. When \(Y=N_K\) the proportion (1) equals one. The conventional high-pole block in the earlier solutions is the case \(K=h\).

\paragraph*{\texorpdfstring{From \(5/32\) to \(r_s=3163/20000\)}{From 5/32 to r\_s=3163/20000}}

Here use \textbf{one constant auxiliary pole \(K=67/400\)}. Its own fraction is \(f_K=1+h-3K\). The first recursion in the preceding upper-ply construction works at this \(h,K\) once its stated share guards are checked. It yields a material with coefficients \((S_K,A_K,N_K)\), where in this paragraph

\[
 N_K=(1-K)(h-K),\qquad A_K=(2K-h)(t-h),\qquad S_K=f_Kh-A_K-N_K .
\]

The proof and layer-by-layer recipe there apply with \(K\) as the pole and \(f_K\) as its fraction. For clarity, with \(u=K(t-h)/(1-K+t-h)\), \(p_0=h-u\), \(c_0=1+h-2K-p_0\), \(s_0=h-K\), use the lower ply having coefficients \(X_0=(1-K)u\), \(L_0=N_K+t s_0(1-s_0/p_0)\) and \(Y_0=c_0 h-X_0-L_0\). Its required transverse coefficient is \(C_0=c_0h-A_K-(1-p_0)s_0\). In the temporary-material recursion take \(\tau_0=(p_0-s_0)/(K-s_0)\). The checked guards include

\[
 s_0<p_0<K,\quad N_K<A_K<L_0,\quad
 Y_0<C_0<X_0,\quad
 \tau_0 S_K<C_0<\tau_0 S_K+(1-\tau_0)K(1-K).
\]

These are exactly the bounds for the shares in that recursion. Its output is the phase-and-coefficient permutation of its temporary input. The retained volume share of the temporary input is strictly below one. Iterating and substituting a physical simple laminate of the same fraction at every terminal occurrence gives the required actual material in the limit, by the compatible-field residual estimate in the cited proof.

Rotate this available material to give \(q_x=N_K,\ q_y=S_K\). Use (2) with \(a=r,b=k\), where \(k\) is now specified. Write

\[
 \begin{aligned}
 N&=(1-r)w, & A&=d(t-h), & \Delta&=N-A,\\
 g_*&=\frac{K-r}{K-w},&
 D_*&=d(t-K)-N,&
 C&=\frac{\Delta}{g_*D_*},\\[2pt]
 k&=\frac{w+Cr}{1+C},&
 e&=w(1-w)-\Delta/g_* .
 \end{aligned}
\]

Make \(G_z(k;w,r;e)\) with poles \((w,r,k)\). Its coefficients along \(x,z\) are \((1-k)w-(r-k)e/d\) and \(e\). The \(x\) coefficient equals \((k-w)(t-K)\) by the definition of \(C\). This matches \(B_x^K(r,k)\), the orientation of (2) with poles \((K,r,k)\) and \(z\) coefficient \((1-K)w\). Stack them normal to \(x\) in respective shares \(g_*,1-g_*\). This gives \((r,r,k)\) with \[
 q_z=g_*e+(1-g_*)(1-K)w=N-\Delta=A .
\] Stack this lower ply with the ordinary \(B_z^h(r,r)\) at \((r,r,h)\) normal to \(z\), giving the lower ply share \(w/(h-k)\). The upper ply has the matching \(q_z=A\). The result has \(p=(r,r,r)\) and fraction \(3/4\). Isotropize to its normalized trace using rotated copies.

On this interval \(\Delta=-(2r-1)(32r-5)/8\ge0\). The verifier checks \(0<w<r<K<h<t\), \(D_*>0\), \(e>0\), and both bounds in (1) for \(a=r,Y=S_K\), in addition to the recursion guards above. Thus \(w\le k<r\) and \(0<e\le w(1-w)\), as needed for the \(G_z\) block. At \(r=5/32\) one has \(k=w\) and the corresponding zero-share plies are omitted. Approximating the constant-\(K\) material by its fixed-\(f_K\) truncations preserves the final fraction. Monotonicity or continuity of the finite outer stacks preserves the tensor limit.

\paragraph*{\texorpdfstring{From \(r_s\) to \(\rho\)}{From r\_s to \textbackslash rho}}

This time use the finite material at the first descending pole of the descending-pole construction in Section~\ref{sup:arc}, \[
 K=\frac{3(1-h)-\sqrt{(1-h)(9-17h)}}{2}.
\] It has two coefficients equal to \(N_K\), hence (2) is available with \(Y=N_K\). Equivalently \(T(K)=h\) for \(T(v)=t-(1-v)(h-v)/(2v-h)\). The following construction uses only finite stacks, followed by isotropization.

Keep \(N,d,w\) from above and now put \(\delta=t-K\), \(A'=d\delta\), \(\epsilon=A'-N\), \(b=(K-r)/(K-w)\), \(m=1-r+\delta\). In the open interval before the right endpoint, set

\[
 \gamma=\frac{(1-2b)\epsilon+
             \sqrt{(1-2b)^2\epsilon^2+4mbd\epsilon}}{2m},
 \quad k=r-\gamma,\quad l=w+\gamma,\quad e=w(1-w)-(d-\gamma)\epsilon/\gamma .
\]

Write \(J=(1-r)l,\ B=(k-w)\delta\). The recipe is the following.

\begin{enumerate}
\def\labelenumi{\arabic{enumi}.}
\tightlist
\item
  Stack \(G_x(k;w,r;e)\) with poles \((k,w,r)\), and \(B_y^K(k,r)\) at \((k,K,r)\), normal to \(y\) with first share \(b\). The normal match \((1-k)w-\gamma e/d=B\) follows from the formula for \(e\). The result has poles \((k,r,r)\) and coefficients \((X,B,Z)\), where \[
  X=be+(1-b)(1-K)w,\qquad
  Z=b[(1-k)r-(1-\gamma/d)e]+(1-b)(1-K)l.
  \]
\item
  Stack the result and its \(y,z\) exchange normal to \(x\), giving the original share \((J-B)/(Z-B)\). Now \(q_z=J\). Stack this mixture normal to \(z\) in share \(\mu=(k-l)/(r-l)\) with \(L_y(1-r)\) at \((k,r,l)\) in the balance. Exchange \(x,y\) to call the resulting material \(M\) with poles \((r,k,k)\), \(q_z=J\), and \(q_y=(1-\mu)(1-r)k+\mu X=B\). The last equality is exactly \(m\gamma^2+(2b-1)\epsilon\gamma-bd\epsilon=0\).
\item
  For a second ply \(H\), stack \(G_z(k;w,r;e)\) at \((w,r,k)\) and \(B_x^K(r,k)\) at \((K,r,k)\) normal to \(x\) with first share \((K-k)/(K-w)\). They match \(q_x=B\). Exchange \(x,y\) to get poles \((r,k,k)\) and \(q_y=B\).
\item
  Stack \(M,H,B_y^K(r,k)\) normal to \(y\), with \(B_y^K\) at poles \((r,K,k)\). Its normal coefficient is also \(B\). Give \(M,H\) total share \(\lambda=(K-r)/(K-k)\). The poles become \((r,r,k)\). Giving all of \(\lambda\) to \(M\) or \(H\) yields, respectively, \[
    V_M=N+\frac{t\gamma(K-r)}{K-r+\gamma},\qquad
    V_H=N-b[w(1-w)-e]
  \] for \(q_z\). Use proportion \((A'-V_H)/(V_M-V_H)\) of \(\lambda\) for \(M\). The resulting \(q_z=A'\) matches \(B_z^K(r,r)\) at \((r,r,K)\). Stack with this final upper ply normal to \(z\) using the same share \(\lambda\) for the lower material. The result has fraction \(3/4\) and all three poles equal to \(r\).
\end{enumerate}

Since the quadratic coefficients sum to \(3h/4\), the normalized trace is \(K_{3/4}(r,s)\). Both recipes give exactly that isotropic conductivity after Schulgasser lamination, in reiterated homogenization. In this second recipe the \(K\) seed has depth four, (2) has depth six, \(M\) has depth nine, and the common-pole root has depth eleven. Two isotropizing stacks suffice, giving a depth bound of thirteen.

\paragraph*{Signs, endpoints and independent mean calculations}

One has \(\epsilon=d[T(r)-K]\). At \(h=3r-1/4\) direct substitution gives \[
 T(T(r))-h=-\frac{2432r^4-2800r^3+1140r^2-221r+16}
 {(4r-1)(208r^2-152r+19)} .
\] On \([r_s,158391/10^6]\), \(T(r)>h/2\), as \(2T(r)-h=-(208r^2-152r+19)/[4(4r-1)]>0\). The numerator quartic changes sign from positive to negative at its unique root \(\rho\) there, with positive sign to its left. Since \(T\) increases and \(T(K)=h\), \(\epsilon>0\) for \(r_s\le r<\rho\). At \(\rho\) one has \(T(r)=K\) and \(T^2(r)=h\). The finite two-coefficient material at the second inverse iterate supplies the endpoint directly.

The interval proof for the second recipe uses \(L=\sqrt\epsilon\) and \[
 z=\frac{(1-2b)L+\sqrt{(1-2b)^2L^2+4mbd}}{2m},
 \qquad \gamma=Lz,\qquad e=w(1-w)-(d-\gamma)L/z .
\] These expressions extend continuously to \(\epsilon=0\) with \(z>0\). The verified strict inequalities include \[
 0<r<K<h<t,\quad 0<d-2\gamma,\quad e>0,\quad
 r(1-r)-w(1-w)>e,\qquad
 \frac{t z(K-r)}{K-r+\gamma}>L .
\] For \(\epsilon>0\) this puts \(w<l<k<r\) and \(e\le w(1-w)<r(1-r)\). The following identities prove the remaining share bounds: \[\begin{gathered}
 J-B=b(d-2\gamma)L/z,\quad
 Z-J=\frac{b(d-\gamma)}d[r(1-r)-w(1-w)-e],\\
 V_M-A'=L\!\left[\frac{t z(K-r)}{K-r+\gamma}-L\right],\quad
 A'-V_H=\epsilon+b[w(1-w)-e].
\end{gathered}\]

Run \texttt{python\ follow-up/verify\_full\_upper.py} for the proof-level arithmetic. It uses rational Sturm root counts and symbolic identities, followed by outward-rounded integer intervals on slabs of width at most \(10^{-7}\) in \(r\). It checks every strict share guard on the entire two ranges. It also computes the three full contrast functions by independent harmonic and arithmetic means from the pure phases. In the first recipe it verifies the repeated seed as an exact fixed point of the actual means. In the finite second recipe it uses the rational sample \(r=2293/14492\), \(K=583/3623\) and computes the shares in a real quadratic extension.
\phantomsection\label{endfile:upper-bridge}

\subsection{An additional rational tree at fraction five twelfths}\label{route:an-additional-rational-tree-at-fraction-five-twelfths}
The file \path{certificates/p6-late/tree_72_30_24.json} specifies a
tree of depth 18 with 2,611 internal nodes and 2,718 leaves. The two
leaf responses are $(1+h)I$ and $I$. At each node the layer shares
are positive rationals summing to one. If the layers have normal
$e_i$, average the $i$th diagonal conductivity harmonically and
the other two arithmetically. The accompanying verifier implements
these operations in $\mathbb Q(h)$, cancelling common polynomials
after each addition. The reduced numerators and denominators have
degree at most two.
Summing the leaf shares gives $f=5/12$. For the three diagonal
responses the computation gives
\[
 \frac{K_1+K_2+K_3}{3}
 =\frac{1+3h/4+25h^2/432}{1+h/3}
 =1+\frac{5h}{12}-\frac{35h^2}{432(1+h/3)}.
\]
To see why the trace yields an attainable scalar conductivity, mix
the tree in equal proportions with a rotated copy which exchanges
two principal axes. Use the third axis as layer normal. Both inputs
agree on that axis and the other two entries are replaced by their
arithmetic mean. Successive pairwise averages in cyclic order preserve
trace and fraction, and converge to the common trace divided by three.
The rotations and fractions are fixed for every $h$. The comparison
$1/3<R(5/12)$ follows from $219<225$. This supplies one pole at this
specified fraction and leaves the extent of the adjacent attainable
set open.

\subsection{The central gap at unbounded depth}\label{route:the-central-gap-at-unbounded-depth}
\phantomsection\label{file:laminate-martingale}
\paragraph*{The terminal law and its martingale}

Set \(h=f+3u-1\). By phase interchange it suffices to handle \(0<h<1\). Write \(\chi\) for the phase-1 indicator and \(g=1-\chi\) on a leaf. We use the following consequences of the spectral concentration proof. The fields \(P_j=uI-F_j\) are symmetric on leaves, with moments of every fixed order uniformly bounded. After taking a subsequence, their joint law with \(\chi\) converges to a law with \(0\preceq P\preceq I\), \(\operatorname{tr}P=1+h-\chi\), and

\[
\mathbb E[\operatorname{tr}P^2\mid\chi=1]=uh.
\tag{1}
\]

On a tree let \(g_v\) be the conditional mean of \(g\) at node \(v\), and use \(\Delta\) for a child-minus-parent conditional mean. A split with normal \(m\) has \(\Delta P_j=\Delta g\,mm^T\). If \(Q=\chi P_j\) and \(\Gamma_j\) is the gradient projection in the tree, then \(\|\Gamma_j Q\|_2\to0\). Multiple plies in the same direction can be grouped into binary splits.

Preserve the root mean of \(Q\) and replace each martingale increment \(\Delta Q\) by \(\Pi_m\Delta Q\Pi_m\), where \(\Pi_m=I-mm^T\). The resulting symmetric martingale has terminal value \(A\). Set \(B=A-P_j+gI\). Orthogonality of martingale increments and symmetry give

\[
\|A-Q\|_2^2\le2\sum_{\text{splits}}\mathbb E\,|(\Delta Q)m|^2
 =2\|\Gamma_jQ\|_2^2\longrightarrow0,
\]

where each increment is extended by zero off its split. The predictable transformation of the increments also preserves uniform fourth-moment bounds up to a constant, by the martingale square-function inequalities and \(|\Pi_m\Delta Q\Pi_m|\le |\Delta Q|\). Thus the terminal laws of \(X=(g,A,B)\) converge to a law supported on the union of the following two compact sets.

{\def\LTcaptype{none} 
\begin{longtable}[]{@{}
  >{\raggedright\arraybackslash}p{(\linewidth - 6\tabcolsep) * \real{0.2500}}
  >{\raggedright\arraybackslash}p{(\linewidth - 6\tabcolsep) * \real{0.2500}}
  >{\raggedright\arraybackslash}p{(\linewidth - 6\tabcolsep) * \real{0.2500}}
  >{\raggedright\arraybackslash}p{(\linewidth - 6\tabcolsep) * \real{0.2500}}@{}}
\toprule\noalign{}
\begin{minipage}[b]{\linewidth}\raggedright
Phase
\end{minipage} & \begin{minipage}[b]{\linewidth}\raggedright
\(g\)
\end{minipage} & \begin{minipage}[b]{\linewidth}\raggedright
\(A\)
\end{minipage} & \begin{minipage}[b]{\linewidth}\raggedright
\(B\)
\end{minipage} \\
\midrule\noalign{}
\endhead
\bottomrule\noalign{}
\endlastfoot
1 & \(0\) & \(P\), with \(0\preceq P\preceq I\) and \(\operatorname{tr}P=h\) & \(0\) \\
2 & \(1\) & \(0\) & \(I-P\), with \(0\preceq I-P\preceq I\) and \(\operatorname{tr}(I-P)=2-h\) \\
\end{longtable}
}

Call the union \(K\). Every increment \(d=(d_g,d_A,d_B)\) of a corrected tree belongs to the closed cone \(\Lambda\) given by the existence of a unit \(m\) with

\[
d_A m=d_Bm=0,\qquad d_A-d_B=-d_g\Pi_m.
\tag{2}
\]

There is a continuous martingale \(X_t\) with that limiting terminal law and the following covariance property. In a time parameter with \(\operatorname{tr}(d[X]_t)\le dt\), the density of its matrix quadratic variation belongs almost everywhere to

\[
\mathcal C=\operatorname{conv}\bigl(\{d\otimes d:\ d\in\Lambda,\ |d|=1\}\cup\{0\}\bigr).
\tag{3}
\]

Here \(|d|\) uses the Frobenius norms on the two matrix components. For a proof, interpolate each nontrivial binary split by one-dimensional Brownian motion on its line segment, starting at the parent and running at unit speed in trace quadratic variation until it hits either child. The child probabilities are the prescribed shares by the martingale property. Continue from that child, and stop the process at its leaf. The expected total running time equals the terminal variance, uniformly bounded by the estimates above. On compact time intervals the paths are tight by the martingale fourth-moment increment bound \(C|t-s|^2\), and their bracket paths are uniformly Lipschitz. Take a joint weak limit of path, bracket, and terminal time. The latter is finite almost surely by the uniform expectation bound. Uniform finite-time moment bounds pass the martingale identities for \(X\) and \(X\otimes X-[X]\) to the joint natural filtration of path and bracket. Uniqueness of the continuous martingale bracket identifies the limiting bracket process. The process is constant after the limiting terminal time, has the desired terminal law by continuity of evaluation there, and (3) follows because all bracket secant slopes are in the compact convex set \(\mathcal C\). The process is a closed martingale, that is, \(X_t\) is the conditional expectation of its terminal value. Indeed the interpolated processes have bounded terminal second moments, so their finite-time second moments are uniformly bounded, also in the limit.

We will also use conditional future laws at stopping times. On the canonical continuous-path space, regular conditional probabilities exist at stopping times. Testing the stopped martingale and product-martingale identities on a countable determining family gives the same identities under almost every conditional future law, hence (3) there. A sequence of such laws with terminals in \(K\) has a subsequential limit of the same type. To see this with a uniform time bound, parametrize each future by the trace of its quadratic variation, and keep it constant after its final bracket time. Intervals of zero bracket carry no martingale variation. The new bracket has trace speed one until the end. Its density belongs to (3) since a nonzero covariance divided by its trace is still a convex combination of the unit directions in (3). The expected final bracket time is bounded by \(\operatorname{diam}(K)^2\) by the second-moment identity. The preceding compactness argument applies. At all times, the conditional expectation representation gives

\[
0\le g\le1,\quad 0\preceq B\preceq gI,\quad \operatorname{tr}B=(2-h)g,
\quad 0\preceq A\preceq(1-g)I,\quad \operatorname{tr}A=h(1-g).
\tag{4}
\]

\paragraph*{A shared null direction in phase 2}

Suppose a martingale of this type starts with \(g>0\) and \(n^TBn=0\) for some unit \(n\). Positivity and conditional expectation force \(Bn=0\) at all times. Write \(M=B|_{n^\perp}\), \(H=n^TAn\), and \(t_0=1-h>0\). Equation (4) gives

\[
t_0 gI_2\preceq M\preceq gI_2.
\]

We claim that \(H=0\) at the root. Every linear state identity holding throughout a path has zero quadratic variation. The unit section of the closed cone in (2) is compact, so finite convex representations in (3) suffice by Carathéodory. In any such representation the bracket of a linear identity is a nonnegative sum of squared evaluations of that identity on the directions \(d\). Thus every direction occurring with positive weight obeys \(d_B n=0\) as well as the trace identities from (4). If \(m\) in (2) is parallel to \(n\), then \(d_H=0\). Otherwise \(d_M\) annihilates the nonzero projection of \(m\) onto \(n^\perp\) and has rank at most one.

For \(M>0\) define

\[
J(M)=\frac{(\operatorname{tr}M)^2}{\det M},\qquad F(M)=\operatorname{tr}M\,e^{\beta J(M)}.
\]

\(F\) is convex. It is the perspective of the convex function \(W\mapsto\exp(\beta/\det W)\) evaluated at \(W=M/\operatorname{tr}M\). Convexity of that function follows from convexity of \(-\log\det W\). In a rank-one direction \(M'=x\tau\tau^T\), \(|\tau|=1\), write \(y=x/\operatorname{tr}M\) and \(R=(\operatorname{tr}M)\tau^TM^{-1}\tau\). The determinant is affine along this line. Differentiating yields

\[
J'=J(2-R)y,\qquad J''=2J(R-1)^2 y^2.
\]

The shape bounds give \(1+t_0\le R\le1+1/t_0\) and bound \(J\) above and below by positive constants. Choose \(\beta\) large enough that \((\log F)''=-y^2+\beta J''\ge c y^2\) with \(c>0\). Also \(|(\log F)'|\le C|y|\). For \(\alpha>\max\{2,1+C^2/c\}\), the Hessian of \(H^\alpha F(M)\) is nonnegative in all the required directions. In fact for \(H>0\), after division by \(H^\alpha F\), minimizing the directional second derivative over \(H'/H\) gives \((\log F)''-[(\log F)']^2/(\alpha-1)\ge0\) when \(M'\) has rank at most one. When \(H'=0\), convexity of \(F\) suffices. At \(H=0\) one can use the \(C^2\) extension \((H_+)^\alpha\).

Itô's formula and (3) now imply the submartingale inequality up to \(g=\delta>0\). The function is bounded by \(C_h g\) throughout the face. If the terminal is in phase 2 then \(H=0\), and on stopping at \(g=\delta\) its value is at most \(C_h\delta\). Optional stopping and \(\delta\downarrow0\) give \(H^\alpha F(M)=0\) at the root. As \(g>0\) this proves the claim.

\paragraph*{Reaching a shared null direction from a full-rank phase-1 field}

For an arbitrary martingale of this type starting at \(g_0>0\), tilt its probability law by the terminal density \(g_\infty/g_0\), which selects phase 2. Denote the new law by \(\mathbb P_2\). If \(g_0\le\epsilon\) and a stopping time \(\tau\) satisfies \(g_\tau\le\epsilon\), then

\[
\mathbb E_2|A_\tau-A_0|\le(\sqrt 2+\sqrt3)\epsilon.
\tag{5}
\]

Indeed (2) gives \(|d_A-d_B|^2=2d_g^2\), so the bracket identity and martingale isometry bound the mean squared change of \(A-B\) by \(2\epsilon g_0\). The mean squared change of \(B\) is at most \(\mathbb E|B_\tau|^2\le3\epsilon g_0\) by (4), and \(\mathbb E g_\tau^2\le\epsilon g_0\). Taking the density \(g_\tau/g_0\) at \(\tau\) and using Cauchy-Schwarz gives (5).

The process \(\sqrt{\det B}\) is a nonnegative supermartingale. Along any line of symmetric matrices with rank-at-most-two increment \(B'\), within \(B>0\) its second derivative has the sign of \((\operatorname{tr}L)^2-2\operatorname{tr}(L^2)\le0\), where \(L=B^{-1/2}B'B^{-1/2}\). All \(d_B\) in (2) have rank at most two. Apply Itô's formula to \(\sqrt{\det(B+\delta I)}\) and let \(\delta\downarrow0\) using bounded convergence.

Suppose for a contradiction that the original limiting terminal law gives positive mass to phase-1 matrices with \(A\succeq2d I\) for some \(d>0\). There are conditional future laws starting at arbitrarily small \(g_0>0\) with \(A_0\succeq d I\). To see this, on that positive-mass event follow the continuous martingale until it first hits levels \(g_0\downarrow0\). These stopping values of \(A\) converge to the phase-1 terminal value. Once \(g\) hits zero it remains zero by positivity, \(B=0\) by (4), and (2)-(3) give zero remaining variation of \(A\). Use regular conditional laws at these stops.

Fix \(\epsilon>0\) with \(4(\sqrt2+\sqrt3)\epsilon/d<1/4\). For \(\eta\downarrow0\) take such a starting law with \(g_0/\epsilon=o(\eta^3)\). Stop at the first of \(g=0\), \(g=\epsilon\), or \(\lambda_{\min}(B/g)\le\eta\) while \(g>0\). Under \(\mathbb P_2\) the terminal has \(g=1\), so the stop occurs with \(g>0\). The probability under \(\mathbb P_2\) of reaching \(g=\epsilon\) with \(\lambda_{\min}(B/g)>\eta\) is at most

\[
\frac{\epsilon}{g_0}\frac{\sqrt{\det B_0}}{(\eta\epsilon)^{3/2}}
\le \eta^{-3/2}\sqrt{g_0/\epsilon}.
\]

This follows from the supermartingale estimate and the density at the stop. Equation (5) implies that \(A\succeq(d/2)I\) there with probability at least \(1-2(\sqrt2+\sqrt3)\epsilon/d\). Consequently some conditional futures start with both \(A\succeq(d/2)I\) and \(n^T Bn/g\le\eta\) for a unit \(n\), at \(0<g\le\epsilon\).

In any such future keep \(n\) fixed and stop again at \(g=0\) or \(g=\epsilon\). Under the phase-2 law \(g_\tau=\epsilon\) almost surely and

\[
\mathbb E_2\left[\frac{n^TB_\tau n}{g_\tau}\right]=\frac{n^TB_0 n}{g_0}\le\eta,
\]

since \(B_\tau=0\) wherever \(g_\tau=0\). Markov's inequality and (5) show that one can select conditional futures with fixed root \(g=\epsilon\), \(A\succeq(d/4)I\), and \(n^TBn/g\le\sqrt\eta\). The terminal sets, the wave cone and its covariance class are invariant under simultaneous spatial rotation. Rotate the matrices so the selected direction is the same \(n\) in all these laws. Their compactness, proved above, gives a limiting martingale with terminal in \(K\), initial \(g=\epsilon\), \(A\succeq(d/4)I\) and \(Bn=0\). The shared-direction result forces \(n^TAn=0\) there, a contradiction.

Thus phase-1 terminal matrices have rank at most two almost surely. For their trace \(h\), this gives \(\operatorname{tr}P^2\ge h^2/2\). Equation (1) forces \(u\ge h/2\), equivalent to \(u\le1-f\) when \(0<h<1\). Phase interchange, replacing \((f,u,h)\) by \((1-f,1-u,2-h)\), gives the other side when \(1<h<2\). The cases \(h=0,1,2\) already satisfy the displayed polynomial inequality. This proves the claim.
\phantomsection\label{endfile:laminate-martingale}

The two exterior exclusions use the support and second-moment
constraints at concentrated spectral measure. The argument for those
two intervals, the finite exceptional poles, and the additional
lower-side constructions in the body are derived below.
\subsection{Exterior gaps in the closure of arbitrary-depth laminates}\label{detail:exterior-gaps}\label{file:exterior-gaps}
\paragraph*{Limits of laminates of increasing rank}

Fix a nonreal contrast \(s_0\), a parameter \(u\) in (1), and \(f_n\to f\in(0,1)\). Suppose finite laminate trees, with any ranks and layer directions, satisfy \(\operatorname{Tr}K_n(s_0)/3\to K_f(u,s_0)\). The open gaps from each coated sphere to its neighboring coated cylinder are excluded even under this convergence. There is also the necessary inequality

\[
 18u^3\ \le\ f h^2(h-u)
        +9(1-f)\big[(2-h)u-h(1-h)\big],\qquad h=f+3u-1 . \tag{5}
\]

It holds also with \((f,u)\) replaced by \((1-f,1-u)\). We first establish the field properties that survive when the individual trees only approach the circle.

\textbf{Concentration at one pole.} Use subscript \(n\) for the leaf space, projection, operator and field in the finite-tree proof. The sum \(\mu_n\) of the three spectral measures has mass \(f_n(1-f_n)\). Every weak subsequential limit of the normalized measures gives a barycenter of the circle at \(s_0\) equal to \(K_f(u,s_0)\). A circle point is an extreme point of its disk and the parameterization is injective, so that limit is the point mass at \(u\). Consequently

\[
 \eta_n^2=\|(H_n-uI)F_n\|_{L^2}^2
         =\int_0^1(p-u)^2\,d\mu_n(p)\longrightarrow0 .
\]

Here the matrix norm on each leaf is the Frobenius norm. Write \(P_n=uI-F_n\), and use \(z\) as a variable contrast in \(\mathcal D=\mathbb C\setminus(-\infty,-1]\). For the actual electric-field matrix \(E_n(z)\), the resolvent formula gives the identity

\[
 (1+uz)E_n(z)-(I+zP_n)
       =z^2(I+zH_n)^{-1}(H_n-uI)F_n . \tag{6}
\]

On every compact subset of \(\mathcal D\) the right side has \(L^2\) norm at most a constant times \(\eta_n\), since the spectrum of \(H_n\) is in \([0,1]\). For any closed disk strictly inside \(\mathcal D\), the expected squared supremum of the error over the disk is also \(O(\eta_n^2)\). Apply Cauchy's formula to each matrix entry on a slightly larger concentric circle and then Cauchy-Schwarz before taking the leaf average.

\textbf{Invertibility at complex contrast.} On every tree, \(\det E_{n,\ell}(z)\) is nonzero throughout \(\mathcal D\) on every leaf. Indeed, set \(\theta=\arg(1+z)/2\). Multiplication by \(e^{-i\theta}\) gives both phases strictly positive real parts. Under mean loading \(m\) of a node \(C\), the energy identity gives \(m^TK_C(z)m=\mathbb E_C[(1+\chi z)|E|^2]\), where \(E\) is its local field. Thus \(\operatorname{Re}(e^{-i\theta}m^TK_C(z)m)>0\) for every unit normal \(m\). The determinant of the transfer from parent to child is the ratio of their two normal conductivities by tangential continuity and normal current continuity. All such ratios are holomorphic and nonzero on \(\mathcal D\).

\textbf{Uniform control of leaf moments.} Let \(F_v\) denote the field increments accumulated down to a node \(v\). These are the conditional means of the leaf field and form a martingale. Conditional orthogonality gives \(\mathbb E_v\|F_\ell-F_v\|^2=f_v(1-f_v)\le1/4\). Starting at \(v\), stop when \(\|F_c-F_v\|\) first reaches one, or at a leaf. Conditional Jensen's inequality at the stopping nodes bounds the probability of reaching one by \(1/4\). Each single increment has norm at most one, so the displacement at a crossing is at most two. Restarting at successive crossings shows that the probability of \(k\) crossings is at most \(4^{-k}\). With fewer than \(k\) crossings the terminal norm starting at the root is at most \(2k+1\). This exponential tail is uniform in the rank. In particular the laws of \((\chi,P_n)\) are tight and all their moments of any fixed degree are uniformly integrable.

Pass to a weak limit of these laws and denote its matrix by \(P\). Every matrix \(M\) in its marginal support satisfies \(0\preceq M\preceq I\). Otherwise \(\det(I+zM)\) has a zero in \(\mathcal D\). Enclose the zero in a closed disk there on whose boundary the determinant is nonzero. For a sufficiently small neighborhood of \(M\), a positive proportion of leaves have \(P_n\) in that neighborhood along the subsequence. By (6), on some of those leaves the error is uniformly as small as needed on the disk. On such a leaf \(\det[(1+uz)E_{n,\ell}(z)]\) is uniformly close to \(\det(I+zM)\) on the boundary, within its minimum modulus there. Rouché's theorem gives a zero in the disk. This contradicts the complex-contrast invertibility just proved.

\textbf{Phase moments.} In the limiting law \(\operatorname{Tr}P=3u+f-\chi\). Also \(\|\Gamma_n(\chi P_n)\|_{L^2}=\eta_n\). Thus the inner product of \(\chi P_n\) with \(F_n\) tends to zero by orthogonality of the projection and the uniform \(L^2\) bound on \(F_n\). Passing to the limit, with \(\mathbb E_i\) denoting conditional expectation in phase \(i\), gives

\[
 \mathbb E_1\operatorname{Tr}P^2=uh,\qquad
 \mathbb E_2\operatorname{Tr}(I-P)^2=(1-u)(2-h).
\]

The second identity follows by interchanging phases, which replaces \(H_n,F_n,P_n\) by \(I-H_n,-F_n,I-P_n\). The bound \(\operatorname{Tr}P^2\le h^2\) in phase 1 and its phase-reversed counterpart now give the coated-cylinder restrictions just as for an exact finite tree.

\textbf{A cubic constraint.} On each finite tree the determinant is affine along the line between the children's \(P\) matrices, since their difference has rank at most one. It follows recursively that \(\mathbb E\det P_n=\det(uI)=u^3\). Uniform integrability gives \(\mathbb E\det P=u^3\) in the limit. Write \(S_2(P)\) for the sum of the three pairwise products of eigenvalues. For three nonnegative eigenvalues one has \(9\det P\le(\operatorname{Tr}P)S_2(P)\) by the arithmetic-geometric mean inequality. In phase 1 this gives

\[
 \mathbb E_1\det P\le\frac{h}{9}\mathbb E_1 S_2(P)
                    =\frac{h^2(h-u)}{18}.
\]

In phase 2 use \(\det P+\det(I-P)=1-\operatorname{Tr}(I-P)+S_2(I-P)\). Both determinants are nonnegative, so the phase-2 moment gives

\[
 \mathbb E_2\det P\le\frac{(2-h)u-h(1-h)}{2}.
\]

Taking the volume-weighted sum of these bounds proves (5). At \(f=19/20\), it requires

\[
        144000u^3-272400u^2+86960u-3761\le0 .
\]

The polynomial is concave on \([1/18,1/3]\) and has endpoint values \(20579/81\) and \(877/3\). Hence this entire closed interval is excluded from the closure of laminate responses. It lies strictly between the listed simple laminate at \(1/20\) and LA at \(7/20\), and strictly inside the two coated-cylinder parameters. At any fixed nonreal contrast its compact image on the arc has positive distance from the set of laminate traces at this fraction. These closure results concern limits of the laminate class. The geometry of an arbitrary two-phase periodic cell admits a different corrector determinant behavior. Briane, Milton and Nesi.
\phantomsection\label{endfile:exterior-gaps}

\subsection{Explicit recipes for the three short quarter-fraction bands}\label{detail:short-bands}\label{file:short-bands}
At phase-1 fraction \(f=1/4\), the poles \(u=0.3954310899625498\ldots\) and \(u=0.8522624015603650\ldots\) are attainable. The construction yields a pole between the coated-cylinder and (K) parameters at every \(0<f<1\), together with its phase interchange.

Choose \(0<m<(\sqrt{17}-1)/8\) and set

\[
r=\frac{1-m-4m^2}{2(1+m)(1-2m^2)},\qquad
f=\frac{m(1-m)(8m+5)}{2(1+m)(1-2m^2)},\qquad
h=(1+m)r=f+3r-1 .
\]

Write \(t=1-h\), \(\delta=t-h\), \(w=h-r=mr\) and define

\[
a=\frac{r\delta}{2t-r},\qquad p=h-a,qquad A=(2r-h)\delta.
\]

The fraction increases continuously from zero to one. In fact the numerator of \(f'\) after multiplication by \(2(1+m)^2(1-2m^2)^2\) is \(22m^4+4m^3-11m^2+6m+5>0\) since \(m<2/5\). The numerator of \(r\) is positive and \(\delta=m(1+2m)/(1-2m^2)>0\). The identity

\[\frac pr=\frac{(1+m)^2}{1+3m}\]

gives \(w<p<r<h<1/2\), with \(0<a<r\). In particular \(0<w/p<1\).

\paragraph*{Matching the two plies}

Use the three-stage block \(K_i(b,c)\) with poles \((h,b,c)\) in axis order \((i,j,n)\) and coefficients \(((b+c-h)\delta,t(h-c),t(h-b))\). It is a stack of simple laminates and pure phase 2 followed by a stack with another simple laminate; the orientation and every share are given in the linked recipe. Here \(L_i(c)\) denotes a simple laminate of fraction \(c\) normal to \(i\).

The upper ply is \(K_z(r,r)\), with poles \((r,r,h)\) and \(q_z=A\). To make the lower ply, stack \(L_y(1-r)\) and \(K_x(r,p)\) normal to \(x\), giving the simple ply share \(w/p\). The two sets of poles are \((a,r,p)\) and \((h,r,p)\). Both normal coefficients equal \((1-r)a=(p-w)\delta\) by the definition of \(a\). Their averaged poles and averaged \(z\) coefficient are

\[p^{\rm lower}=(r,r,p),\qquad
q_z^{\rm lower}=(1-r)w+tw(1-w/p).\]

The second summand is strictly positive. Subtracting the upper coefficient gives

\[
q_z^{\rm lower}-A=
\frac{4h^3-8h^2r-4h^2+2hr^2+7hr-2r^2}{r-2h}=0.
\]

Indeed after \(h=(1+m)r\), the numerator becomes \(r^2[1-m-4m^2-2(1+m)(1-2m^2)r]\). Stack the lower and upper plies normal to \(z\), with lower share \((h-r)/(h-p)\). Their common \(q_z\) permits averaging the normal denominators. All three poles now equal \(r\), and the fraction is \(1+h-3r=f\). The coefficient sum is \(fh\), so the normalized trace is \(K_f(r,s)\). Two additional coordinate laminations of rotated copies make the full tensor \(K_f(r,s)I_3\) in reiterated homogenization. This requires at most seven stages, and the root before isotropization retains \(q_z=(2r-h)\delta\).

The inequality \(w>0\) puts \(r\) above \((1-f)/2\). The strictly positive difference \(A-(1-r)w=tw(1-w/p)\) puts it below the smaller root \(\Phi(f)\) of the (K) quadratic, by the same sign argument as for the nine-stage family.

\paragraph*{Both quarter-fraction parameters}

Substituting \(h=f+3r-1\) in the normal-coefficient identity yields

{\def\LTcaptype{none} 
\begin{longtable}[]{@{}lll@{}}
\toprule\noalign{}
\(f\) & Polynomial for \(r\) & Physical root bracket \\
\midrule\noalign{}
\endhead
\bottomrule\noalign{}
\endlastfoot
\(1/4\) & \(672r^3-1016r^2+456r-63\) & \((395/1000,396/1000)\) \\
\(3/4\) & \(672r^3-520r^2+96r-5\) & \((147/1000,148/1000)\) \\
\end{longtable}
}

The parameter \(m\) is the unique root on its allowed interval of \(14m^3-8m^2-9m+1\) or \(10m^3-12m^2-7m+3\), respectively, by the strict increase of \(f\). Testing endpoints brackets these roots in \((10333,10334)/100000\) and \((30781,30782)/100000\), respectively. Direct substitution of these bounds in \(r(m)\) gives the brackets in the table; each cubic changes sign and is strictly monotone there. Interchanging the phases in the second row realizes \(1-0.1477375984396350\ldots\) at one quarter. Both points lie closer to their respective coated-cylinder parameters than the nine-stage contacts.

At \(f=1/4\), every pole in the closed interval \[
[r_7,r_9]=[0.3954310899625\ldots,\ 0.3970788070635\ldots]
\] is attainable in the closure of separated-scale laminates. Here \(r_7\) and \(r_9\) are the direct-order roots in the seven- and nine-stage constructions. The finite continuation from \(r_9\) gives the joined band \([r_7,0.39748]\).

\paragraph*{A sheet of plies with one specified current coefficient}

Fix \(h\) and put \(t=1-h\), \(\delta=t-h\), \[
 k=\frac{4ht}{3t+\sqrt{9t^2-8ht}},\qquad
 l=h-k,\quad N=(1-k)l=(2k-h)\delta,\quad D=(t-k)(2k-h).
\] This is the first (K) pole at the specified \(h\). The following two functions will specify all subsequent stack shares: \[
 c=\frac{D-N}{\delta+t},\qquad m=\frac{D/l-h}{\delta+t},\qquad
 \phi(k+u)=k+c-mu.
\]

For \(x=k+u\) set \(z_x=\phi(x)\). Start at \((k,k,z_x)\) with \(q_z=N\). To make this material, stack the upper ply \(K_z(k,k)\) at \((k,k,h)\) with \(L_x(1-k),L_y(1-k)\) at \((k,k,l)\) normal to \(z\), using upper share \(\omega=(z_x-l)/k\). The three \(z\) coefficients equal \(N\). Give the simple portion a share of \(L_y\) that puts \[
q_x=(k+z_x-h)\delta .
\] Specifically this proportion is \([(k+z_x-h)\delta-\omega tl]/[(1-\omega)k(1-k)]\). The required range of \(z_x\) and these shares are checked below. The initial \(q_y\) is then \(D-t(z_x-k)\) by the coefficient sum.

Now take any second pole \(y\) in the indicated ranges below. Add \(K_x(k,z_x)\) at \((h,k,z_x)\) normal to \(x\), retaining the preceding material in fraction \(\alpha_x=(h-x)/l\). Then add \(K_y(x,z_x)\) at \((x,h,z_x)\) normal to \(y\), with retention \(\alpha_y=(h-y)/l\). The first match is the chosen \(q_x\). After the first stack \(q_y=\alpha_x D+(1-\alpha_x)tl-t(z_x-k)\). This equals \((x+z_x-h)\delta\) precisely by the formula for \(\phi\), giving the second match. The result has poles \((x,y,z_x)\) and coefficients \[
\begin{aligned}
 q_y&=(x+z_x-h)\delta,\\
 q_x&=\alpha_y(k+z_x-h)\delta+(1-\alpha_y)t(h-z_x),\\
 q_z&=\Gamma(x,y)=\alpha_x\alpha_y N+tl(\alpha_x+\alpha_y-2\alpha_x\alpha_y).
\end{aligned}
\tag{1}
\] In particular \(q_z\) is symmetric in \(x,y\). Growing the coordinates in reverse order instead gives a material at \((x,y,z_y)\) with exactly the same \(z\) coefficient.

\paragraph*{Continuously adjustable lower ply}

At a target \(x=k+u\) define \(v_*=(c-u)/m\). In the application \(0<v_*<u\) and \(z_x<x\). Choose \(0\le v\le v_*\) and \(y=k+v\). Combine the two orders of (1) normal to \(z\). The first has \(z_x<x\) and the second has \(z_y=k+c-mv\ge x\). Take first share \(\theta=(x-z_y)/(z_x-z_y)\), so the poles become \((x,y,x)\). Exchange the last two axes, obtaining a lower ply with poles \((x,x,y)\).

The resulting \(z\) coefficient \(Q(v)\) is continuous on \([0,v_*]\). Its two endpoint values are \[
 Q(0)=\frac{D(h-x)}{l},\qquad
 Q(v_*)=\alpha_x(k+x-h)\delta+(1-\alpha_x)t(h-x).
\tag{2}
\] For the second value, \(z_y=x\) and \(\theta=0\), so the \(x\) coefficient from the first-growing axis of the reverse-order stack in (1) gives the result. For the first value, the two numbers averaged in \(Q(v)\) are \[
 N+\delta[u+c-mu],\qquad
 N+\delta(c-mv)-u[l+(\delta+t)(c-mv)/l].
\] At \(v=0\) their weights are \(\theta=(c-u)/(mu)\) and \(1-\theta\). Substituting the definitions of \(c,m\) gives \(D(1-u/l)\), as in (2).

We choose \(y\) so the lower ply supplies the normal coefficient \(h(1+h-3x)/3\) of the isotropic single-pole material at \(x\). The certificate below gives \[
 Q(v_*)< h(1+h-3x)/3 < Q(0).
\] The intermediate value theorem supplies \(v\) and hence the desired geometry. For any closed interval of target poles \(x\) strictly above \(k+c/(1+m)\) along a fixed-\(h\) path, the pole gap \(x-y\) is bounded below by the positive minimum of \(u-v_*\).

\paragraph*{Applying the lower plies at fixed phase fraction}

For a specified quarter-fraction target \(\bar r\) let \(h=3\bar r-3/4\). Take the physical isotropic seven-stage core at the same \(h\), whose pole is \[
 b_7(h)=\frac{8ht}{7-8h+\sqrt{17-48h+32h^2}}.
\] For \(r_7<\bar r<r_9\) one has \[ k+c/(1+m)<\bar r<b_7(h). \tag{3}\] Indeed \(k+c/(1+m)=h(\delta+D/l)/(2\delta+D/l)\) is exactly the nine-stage contact at this \(h\). Equality in either part of (3) would give the corresponding fixed-fraction endpoint polynomial. On \([3954/10000,3971/10000]\) the standard interval substitution gives increasing coefficient signs \(+, \) for the seven-stage cubic and \(, -+\) for the nine-stage quartic. Thus these are the only equalities in that bracket. At the interior sample \(\bar r=396/1000\), direct interval evaluation gives \(k+c/(1+m)<396/1000<b_7\). This fixes both signs throughout the open interval.

Start at \(x=b_7(h)\) and decrease the pole to \(\bar r\) in small steps. For a step of size \(\epsilon\) take the lower ply with poles \((x,x,y)\) just constructed in share \(3\epsilon/(x-y)\), normal to \(z\) with the scalar material at \(x\). Its normal coefficient agrees and all denominators are \(1+xs\) before the normal pole shift of \(-3\epsilon\). Isotropizing the result changes the ideal response relative to \(K_{1+h-3(x-\epsilon)}(x-\epsilon,s)\) by \(O_s(\epsilon^2)\). This is the decreasing-parameter version of the one-sided continuation lemma; the sign change leaves its Taylor and energy estimates identical. The uniform pole gap on the path and bounded ply coefficients give vanishing total error. Each truncation has exactly fraction \(1+h-3\bar r=1/4\). Both endpoints of the resulting band already have finite-stage realizations.

\paragraph*{Certificate of physical shares}

The verifier covers the full rectangle of candidate target poles \(\bar r\in[3954/10000,3971/10000]\) in strips of width \(10^{-5}\) and encloses \(x\) from the lower strip endpoint to the upper enclosure of \(b_7(h)\) in intervals of width at most \(10^{-5}\). There are 45,621 boxes. For paths satisfying (3), we have \(u>c/(1+m)\) and hence \(z_x<x\) and \(v_*<u\). The strict positivity of \(c-u\) below gives \(v_*>0\). In the two orders of (1), \(z_x\) and \(z_y\) lie in \([z_x,k+c]\), while \(k\le y<x<h\). In particular the maximum \(x,y\) is below \(h\) because \(x\le b_7<h\); all high-block arguments exceed \(h/2\).

The lower and upper differences for the initial \(q_x\) split are affine in its third pole, so endpoint checks on \([z_x,k+c]\) suffice. The following are lower bounds verified by exact outward integer rounding, rounded down for display:

{\def\LTcaptype{none} 
\begin{longtable}[]{@{}
  >{\raggedright\arraybackslash}p{(\linewidth - 2\tabcolsep) * \real{0.5000}}
  >{\raggedright\arraybackslash}p{(\linewidth - 2\tabcolsep) * \real{0.5000}}@{}}
\toprule\noalign{}
\begin{minipage}[b]{\linewidth}\raggedright
Quantity
\end{minipage} & \begin{minipage}[b]{\linewidth}\raggedright
Lower bound
\end{minipage} \\
\midrule\noalign{}
\endhead
\bottomrule\noalign{}
\endlastfoot
\(\delta\), \(l\), \(m\), \(c-u\), \(u\) & \(0.1174\), \(0.0605\), \(0.7218\), \(0.00585\), \(0.0164\) \\
\(z_x-l\), \(h-(k+c)\) & \(0.3219\), \(0.0317\) \\
Lower and upper split margins at \(z_x\) & \(0.00932\), \(0.01739\) \\
Lower and upper split margins at \(k+c\) & \(0.00979\), \(0.00960\) \\
\(h(1+h-3x)/3-Q(v_*)\) & \(0.000719\) \\
\(Q(0)-h(1+h-3x)/3\) & \(0.00180\) \\
\end{longtable}
}

All quantities are computed with integer endpoints at scale \(10^{18}\) and integer square roots.

Let \(r_9(f)\) be the algebraic pole specified in the nine-stage construction. At \(f=1/4\), isotropic separated-scale laminates realize the entire closed interval \[
\boxed{r_9(1/4)=0.3970788070634548\ldots\ \le u\le\frac{9937}{25000}=0.39748.}
\] At \(f=3/4\) the same argument gives \(r_9(3/4)=0.1542794659621522\ldots\le r\le31/200\). Its phase interchange is contained in the wider finite band starting at \(27/32\). The constructions in this note need at most thirteen stages.

\paragraph*{The seed just beyond (K)}

Keep \(0<h<1/2\) fixed, set \(t=1-h\) and \(\delta=1-2h\). Write \(a\) for the seed pole, \(l=h-a\), \(d=a-l\), \(N=(1-a)l\) and \(A=d\delta\). The ten-stage construction before isotropization has all three poles \(a\) and its last normal coefficient is \(A\). Here we give a criterion to apply it even when the seed fraction \(1+h-3a\) differs from either of the quarter fractions treated in that construction.

Assume \(0<l<a<h\), put \(\Delta=A-N>0\), and write \(k=a-\gamma\) for the auxiliary pole in the seed. In the notation of the ten-stage construction, with \(a\) in place of its active pole \(r\), the equation \(P=B\) is \[
m\gamma^2+c\Delta\gamma-b\Delta=0,\qquad
m=1-a+\delta,\quad c=2l/a-1,\quad b=ld/a.
\] Choose its positive root. Indeed the original quantities reduce to \[
B=A-\delta\gamma,\qquad J=(1-a)(l+\gamma),\qquad
E=l(1-l)-(d-\gamma)\frac{\Delta}{\gamma},\qquad
X=N+\frac{l}{a}\Delta(1-d/\gamma).
\] Then \(P-N=(1-a)\gamma-(l/a)(d-2\gamma)\Delta/\gamma\), giving the quadratic. Its root also obeys \(\Delta/\gamma=m\gamma/(b-c\gamma)\).

It suffices for the seed that \[
0<\gamma<d/2,\qquad
E\ge0,\qquad m>\Delta/\gamma,\qquad
\frac{tl}{l+\gamma}\ge\frac{\Delta}{\gamma}.
\tag{1}
\] In fact \(J-B=m\gamma-\Delta\), \(V_H=X<A\), and \(V_M-A=tl\gamma/(l+\gamma)-\Delta\). The other two bounds \(E<l(1-l)\) and \(Z>J\) follow from the identities in the ten-stage proof because \(\Delta>0\) and \(k>h/2>l\). Thus its physical stack gives the required common-pole seed. In particular (1) holds just beyond the (K) point at any fixed \(h\): as \(\Delta\) increases from zero, \(\gamma\) is asymptotic to \(\sqrt{b\Delta/m}\), all limiting denominators are strictly positive, and \(E\) tends to \(l(1-l)>0\).

Let \(q_1,q_2\) be the seed's two transverse coefficients, ordered here by their formulas. The penultimate material has one transverse coefficient \(B\), and the last added \(K_z\) has transverse coefficients \(tl\). Hence \[
q_1=\frac{l}{l+\gamma}B+\frac{\gamma}{l+\gamma}tl,\qquad
q_2=(1+h-3a)h-A-q_1 .
\tag{2}
\] Rotate the seed so that \(A\) lies on \(x\). A quarter-turn about \(x\) then allows any \(y\) coefficient between \(q_1,q_2\) by averaging the two copies normal to \(x\).

\paragraph*{Three more layers}

For a target pole \(r>a\) with the same \(h\), let \(w=h-r\) and \(\alpha=w/l\). The three normal directions and the three added materials are, in order, \(x\) with \(K_x(a,a)\), \(y\) with \(K_y(r,a)\), and \(z\) with \(K_z(r,r)\). At each stage retain the preceding material in share \(\alpha\), with the newly added material filling the balance. The first stack can split the retained share between its two orientations. All the \(K_i\) blocks use the same \(h\).

Set \[
C=\frac{(1+h-3a)h-2A}{l},\qquad
q_{\mathrm{need}}=\frac{(r+a-h)\delta-(1-\alpha)tl}{\alpha}.
\] If \(q_{\mathrm{need}}\) is between \(q_1,q_2\), the rotated seed copies can supply it. The \(x\) match is \(A=(2a-h)\delta\). After the \(x\) stack the \(y\) coefficient is \((r+a-h)\delta\), giving the \(y\) match. After the \(y\) stack the remaining \(z\) coefficient is \(C(h-r)\). The final match with \(K_z(r,r)\) therefore requires \(C(h-r)=(2r-h)\delta\). These are the same three arithmetic/harmonic calculations as in the nine-stage proof, with a different seed. The last stack has poles \((r,r,r)\) and fraction \(1+h-3r\). It retains \(q_z=(2r-h)\delta\) before the two final isotropizing stacks.

The pole relation between \(a\) and \(r\) has the symmetric form \[
\frac1{h-r}+\frac1{h-a}=\frac{3(2-3h)}{h\delta}.
\] Thus for a prescribed target \(r\) choose \[
l=\frac{h\delta w}{3(2-3h)w-h\delta},\qquad a=h-l.
\tag{3}
\] The seed needs eight stages before isotropization, the three additions produce the common-pole root, and the two isotropizing stacks give \(K_f(r,s)I_3\) in reiterated homogenization.

\paragraph*{Closed intervals}

For the quarter and three-quarters targets in the first paragraph, substitute \(h=f+3r-1\) in (3). If \(Q_f(r)\) denotes the corresponding quartic from the table in the nine-stage proof, then \[
\begin{aligned}
\Delta&=\frac{-(4r-1)(12r-5)Q_{1/4}(r)}
                   {2(192r^2-180r+41)^2} &&(f=1/4),\\
\Delta&=\frac{-(4r-1)(12r-1)Q_{3/4}(r)}
                   {2(192r^2-92r+9)^2} &&(f=3/4).
\end{aligned}
\] In both cases \(\Delta\) vanishes at \(r_9(f)\) and is positive to its right on the stated interval. This follows from the isolating intervals in the nine-stage proof: each quartic is negative at its lower isolating endpoint and positive at its upper endpoint, with exactly one zero there. At \(r=r_9(f)\) the nine-stage material supplies the endpoint directly.

All remaining inequalities are strict away from this endpoint, as the following outward-rounded calculation certifies. Cover \([0.39707,0.39748]\) or \([0.15427,0.155]\) by \(r\)-intervals of length \(10^{-6}\). Evaluate the formulas with integer endpoints scaled by \(10^{18}\), rounding each multiplication, division and square root outwards. Intersect the interval for \(\Delta\) with \([0,\infty)\), since only the portion \(r\ge r_9(f)\) is needed. For the positive root use \(\gamma=[-c\Delta+\sqrt{c^2\Delta^2+4mb\Delta}]/(2m)\), and enclose \(\Delta/\gamma\) by \(m\gamma/(b-c\gamma)\), which also works at the endpoint. The denominator in (3), \(\delta,w,l,d\) and \(r-a\) have strictly positive lower bounds. The following bounds are rounded down from the exact integer output.

{\def\LTcaptype{none} 
\begin{longtable}[]{@{}lll@{}}
\toprule\noalign{}
Quantity & \(f=1/4\) & \(f=3/4\) \\
\midrule\noalign{}
\endhead
\bottomrule\noalign{}
\endlastfoot
\(d/2-\gamma\) & \(0.1500\) & \(0.0390\) \\
\(E\) & \(0.00090\) & \(0.0321\) \\
\(m-\Delta/\gamma\) & \(0.5603\) & \(1.0000\) \\
\(tl/(l+\gamma)-\Delta/\gamma\) & \(0.2810\) & \(0.4069\) \\
\(q_{\mathrm{need}}-q_1\) & \(0.00325\) & \(0.00095\) \\
\(q_2-q_{\mathrm{need}}\) & \(0.0130\) & \(0.00293\) \\
\end{longtable}
}

The table verifies (1), the two orientation shares, and \(a<r<h\) throughout. It certifies every pole of the indicated intervals, since the computation encloses continuous parameter boxes.

At \(f=1/4\), all \(1999/5000\le r\le40101/100000\) are attainable. Above \(2001/5000=0.4002\) the construction below is finite. The right endpoint overlaps the finite band starting at \(\xi=0.4010025920\ldots\), giving the joined closed interval \([0.3998,(39-\sqrt{177})/64]\). The left-hand extension uses fixed-fraction limits of finite laminates when the upper ply's normal coefficient exceeds \(N=(1-r)w\).

Put \(h=3r-3/4\), \(t=1-h\), \(\delta=t-h\), \(w=h-r\), \(d=r-w\). Recall the physical seven-stage lower ply at poles \((r,r,p)\): \[
a=\frac{r\delta}{2t-r},\qquad p=h-a,\qquad
q=(X,Y,L).\]\[
X=(1-r)a,\quad L=N+tw(1-w/p),\quad Y=(1+h-2r-p)h-X-L.
\] It is the \(x\)-normal stack of \(L_y(1-r)\) at \((a,r,p)\) and \(K_x(r,p)\) at \((h,r,p)\), with first share \(w/p\). All poles and coefficients below use the \((x,y,z)\) coordinate order.

\paragraph*{A lower-current upper ply}

At fixed \(h\) the seven-stage pole itself is the smaller root \[
b=\frac{8ht}{7-8h+\sqrt{17-48h+32h^2}}
\] of \(2t b^2-h(7-8h)b+4h^2t=0\). This is the seven-stage normal-coefficient matching equation with \(h\) held fixed. It can be used here at the fraction \(1+h-3b\). Put \[
l=h-b,\quad a_b=b\delta/(2t-b),\quad p_b=h-a_b,\quad A_b=(2b-h)\delta.\]\[
X_b=(1-b)a_b,\quad Y_b=(1+h-2b-p_b)h-X_b-A_b .
\] The seven-stage lower ply at \((b,b,p_b)\) has coefficients \((X_b,Y_b,A_b)\), and the upper ply at \((b,b,h)\) has \((tl,tl,A_b)\). For a general desired pole \(r\) between \(p_b\) and \(b\), proceed as follows.

\begin{enumerate}
\def\labelenumi{\arabic{enumi}.}
\tightlist
\item
  Exchange \(x,z\) in both plies and stack them normal to \(x\), giving the former share \(\lambda=w/a_b\). The poles are \((r,b,b)\) and \(q_x=A_b\). Rotation about \(x\) permits any \(q_y\) from \(\lambda Y_b+(1-\lambda)tl\) to \(\lambda X_b+(1-\lambda)tl\). Take \(q_y=Q=(1-b)w\) by interpolating the two orientations. The remaining coefficient is \(h(1+h-r-2b)-A_b-Q\).
\item
  Stack this result with \(L_z(1-b)\) at poles \((r,w,b)\) normal to \(y\), giving the former share \(d/(b-w)\). Both have \(q_y=Q\). At the resulting poles \((r,r,b)\) the axial coefficient is \[
  U=\frac{d}{b-w}\big[h(1+h-r-2b)-A_b-Q\big].
  \] Exchange \(x,y\) to place the transverse coefficient \(Q\) on \(x\). This is the upper ply in both constructions below.
\end{enumerate}

\paragraph*{\texorpdfstring{When \(U\le N\)}{When U\textbackslash le N}}

A finite lower ply suffices. Choose \[
e=\frac{rU-dtw}{w},\qquad
v=\frac{dw(1+\delta)-er}{d(w+\delta)-e}.
\] Since \(dtw/r<U\le N\), one has \(0<e\le w(1-w)\). The denominator is positive because \(d(w+\delta)-w(1-w)=d\delta-N>0\). The formula gives \[
v-w=\frac{d[w(1-w)-e]}{d(w+\delta)-e}\ge0,\qquad
r-v=\frac{d(d\delta-N)}{d(w+\delta)-e}>0 .
\] Make \(G_z\) at \((w,r,v)\), using tangential pair \(w,r\) and normal split \(e\). Its normal \(q_z=e\) and its \(x\) coefficient is \((1-v)w-(r-v)e/d=(v-w)\delta\). Stack it normal to \(x\) in share \(w/r\) with \(K_x(r,v)\) at \((h,r,v)\). The latter has the same \(x\) coefficient and \(q_z=tw\). The resulting poles are \((r,r,v)\), with \(q_z=(we+dtw)/r=U\). A \(z\)-normal stack with the upper ply, taking this lower ply in share \((b-r)/(b-v)\), gives all three poles \(r\). Two further stacks of rotated copies give the scalar \(K_{1/4}(r,s)I_3\).

\paragraph*{\texorpdfstring{When \(U>N\)}{When U\textgreater N}}

Use a temporary constituent with poles \((r,r,r)\) and coefficients \((U,S,N)\), where \(S=h/4-U-N\). One repetition recreates the temporary constituent and consumes a fraction strictly less than one of it.

Put \(\rho=(b-r)/(b-p)\), \(\tau=(p-w)/d\), and \(C=[N-(1-\rho)Q]/\rho\).

\begin{enumerate}
\def\labelenumi{\arabic{enumi}.}
\tightlist
\item
  Stack the temporary constituent normal to \(z\) in share \(\tau\) with a combination of \(L_x(1-r),L_y(1-r)\) at poles \((r,r,w)\). All have \(q_z=N\). In the simple portion give \(L_y(1-r)\) proportion \((C-\tau U)/[(1-\tau)r(1-r)]\). This yields poles \((r,r,p)\) and \(q_x=C\).
\item
  Separately prepare the real lower ply at \((r,r,p)\) from the beginning of the proof. Combine it and its \(x,y\) exchange normal to \(z\), giving the original share \((C-Y)/(X-Y)\). Its \(q_x\) is now \(C\), with \(q_z=L\).
\item
  Mix this real ply with the result of step 1 normal to \(x\), using real share \(\theta=(U-N)/(L-N)\). They have the same poles and \(q_x=C\). Their new \(q_z\) is \(U\).
\item
  Stack step 3 with the upper ply normal to \(z\), in respective shares \(\rho,1-\rho\). This gives \((r,r,r)\) with \(q_z=U\) and \(q_x=\rho C+(1-\rho)Q=N\). The remaining coefficient is \(S\) by the trace sum. Exchange \(x,z\) to recover the original temporary constituent.
\end{enumerate}

The retained volume per repetition is \(\rho(1-\theta)\tau<1\). After \(n\) repetitions replace the remaining temporary plies by simple laminates of fraction \(1/4\), then isotope the resulting tensor by the fixed rotated stacks. In the ideal recursion, the compatible field and flux pencils of a node are \(I+sP\) and \(I+s(P+cI)+s^2Q\). They have equal tangential electric fields and equal normal currents at each split. Each repetition resets the mean electric pencil of the retained constituent to \((1+rs)I\). Extend these bounded pencils as constant fields and currents into replacement plies. Their constitutive error for fixed positive real phase ratio is then bounded in norm by a constant times \([\rho(1-\theta)\tau]^{n/2}\). The energy minimum bounds the true field error by that norm divided by the smaller phase conductivity. Reiterated homogenization and the usual exact-fraction diagonal limit give \(K_{1/4}(r,s)I_3\).

\paragraph*{Shares throughout the interval}

Exact outward-rounded interval arithmetic on the 6,050 subintervals of width \(1/5000000\) gives the following bounds, rounded down. The upper-ply checks apply on the full interval.

{\def\LTcaptype{none} 
\begin{longtable}[]{@{}ll@{}}
\toprule\noalign{}
Quantity & Lower bound \\
\midrule\noalign{}
\endhead
\bottomrule\noalign{}
\endlastfoot
\(r-p_b,\ b-r,\ p-w,\ r-p\) & \(0.00640,\ 0.0143,\ 0.342,\ 0.00234\) \\
\(Q-[\lambda Y_b+(1-\lambda)tl]\) & \(0.000598\) \\
\(\lambda X_b+(1-\lambda)tl-Q\) & \(0.00201\) \\
\(U-dtw/r,\ d\delta-N,\ L-N\) & \(0.00378,\ 0.00162,\ 0.0238\) \\
\end{longtable}
}

For each subinterval whose enclosures do not already show \(U<N\), the checks also give \(L-U>0.0227\), \(C-Y>0.00241\), \(X-C>0.00352\), \(C-\tau U>0.0000400\) and \(\tau U+(1-\tau)r(1-r)-C>0.00352\). These bounds make all the recycling shares physical. On \(2001/5000\le r\le40101/100000\) the finite-case inequality \(U<N\) holds throughout. The certificate uses exact integer square roots to enclose \(b\).

Define \(\xi=0.401002592072879\ldots\) and \(\eta=0.152563870754587\ldots\) as the unique roots in \(I_1=[401/1000,803/2000]\) and \(I_3=[19/125,5/32]\), respectively, of \[
30720\xi^4-61248\xi^3+43024\xi^2-12684\xi+1323=0,\qquad
896\eta^3-656\eta^2+112\eta-5=0.
\] At phase-1 fraction \(f=1/4\), isotropic separated-scale laminates realize both closed intervals \[
\boxed{\left[\xi,\frac{39-\sqrt{177}}{64}\right]}
\qquad\text{and}\qquad
\boxed{\left[\frac{27}{32},1-\eta\right]} .
\] They are approximately \([0.4010025921,0.4014978954]\) and \([0.84375,0.8474361292]\). Each extends a supplied (K) point toward its neighboring coated-cylinder point. The second interval is obtained by interchanging phases after constructing \([\eta,5/32]\) at \(f=3/4\). Both use at most ten lamination stages, including isotropization.

\paragraph*{Elementary plies}

Take the respective \(f=1/4\) or \(3/4\) and write \(r\) for its unreflected parameter. Below the (K) endpoint, put \(h=f+3r-1\), \(t=1-h\), \(\delta=t-h\), and \(w=h-r\). Select \(k\) below with \(h/2<k<r\) and put \(l=h-k\). Both intervals obey \(0<w<l<k<r<h<t<1\). Define \[
g=\frac{r-k}{r-w},\quad B=(k-w)\delta,\quad
E=\frac{(1-k)w-B}{g},\quad J=(1-r)l,\quad \mu=\frac{k-l}{r-l}.
\] All triples below use coordinate order \((x,y,z)\). A coordinate stack averages its pole and quadratic-numerator triples when its two tangential poles and the normal numerator coefficient agree. The three-stage block \(G_i\) has poles \((k,w,r)\) and coefficients \((E,(1-k)w-gE,(1-k)r-(1-g)E)\) in its designated axis order, with the first axis \(i\). The block \(K_i\) at poles \((h,a,b)\) in its designated order has coefficients \(((a+b-h)\delta,t(h-b),t(h-a))\). Both are built from simple laminates and pure phase 2.

\begin{enumerate}
\def\labelenumi{\arabic{enumi}.}
\tightlist
\item
  Stack \(G_x\) at \((k,w,r)\) and \(K_y\) at \((k,h,r)\) normal to \(y\), with first share \(w/r\). Their common \(y\) coefficient is \(B\). The resulting poles are \((k,r,r)\) and the coefficients are \((X,B,Z)\) with \[
  X=\frac{wE+(r-w)tw}{r},\qquad
  Z=\frac{w[(1-k)r-(1-g)E]+(r-w)tl}{r}.
  \]
\item
  Stack that result and its copy exchanging \(y,z\) normal to \(x\). Choose original share \((J-B)/(Z-B)\), giving coefficients \((X,B+Z-J,J)\).
\item
  Stack the simple laminate \(L_y(1-r)\) at \((k,r,l)\) and the preceding rotated mixture normal to \(z\), in shares \(1-\mu,\mu\). They match \(q_z=J\). The poles become \((k,r,k)\), and the \(x\) coefficient is \(P=(1-\mu)(1-r)k+\mu X\). Exchange \(x,y\) and call the result \(M\). It has poles \((r,k,k)\) and \(q_y=P\), \(q_z=J\).
\item
  Separately make \(H\) at the same \((r,k,k)\). First stack \(G_z\) at \((w,r,k)\) and \(K_x\) at \((h,r,k)\) normal to \(x\), with first share \(l/r\). The common \(x\) coefficient is \(B\). Exchange \(x,y\). This gives \(H\) with \(q_y=B\) and \(q_z=H_z=[lE+(r-l)tw]/r\).
\end{enumerate}

\paragraph*{The final stacks}

Choose \(k\) so that \(P=B\). Now \(M,H\) and the upper \(K_y\) at \((r,h,k)\) all have \(q_y=B\). They can enter one \(y\)-normal stack with total share \(w/l\) for \(M,H\). The new poles are \((r,r,k)\). Giving that whole share to \(M\) or to \(H\) would give, respectively, \[
V_M=\frac wl J+\left(1-\frac wl\right)tw,\qquad
V_H=\frac wl H_z+\left(1-\frac wl\right)tw
\] in the \(z\) coefficient. Set \(A=(2r-h)\delta\). Assign \(M\) the proportion \((A-V_H)/(V_M-V_H)\) of the \(w/l\) share, and give \(H\) the balance of that share. This gives \(q_z=A\).

Stack the resulting lower ply with \(K_z\) at \((r,r,h)\) normal to \(z\), again with first share \(w/l\). Their common \(q_z=A\) permits averaging to \((r,r,r)\). The fraction is \(1+h-3r=f\) and the coefficient sum is \(fh\). Two further laminations of rotated copies yield the full isotropic tensor \(K_f(r,s)I_3\). The blocks in steps 1 and 4 need at most three stages; \(M\) needs six and \(H\) needs four. The combined \(y\) stack, the \(z\) stack and isotropization give the asserted upper bound of ten. This construction is independent of contrast and realizes the entire function in reiterated homogenization.

\paragraph*{Certificates for every share}

After substituting \(h=f+3r-1\), form \(F(r,k)=128r(r-k)(P-B)\) and cancel denominators. This is a quadratic polynomial in \(k\). For the direct and the phase-reversed construction its coefficients at the bracket are as follows. Set \(W=64r^2-78r+21\) for the first row.

{\def\LTcaptype{none} 
\begin{longtable}[]{@{}
  >{\raggedright\arraybackslash}p{(\linewidth - 6\tabcolsep) * \real{0.2500}}
  >{\raggedright\arraybackslash}p{(\linewidth - 6\tabcolsep) * \real{0.2500}}
  >{\raggedright\arraybackslash}p{(\linewidth - 6\tabcolsep) * \real{0.2500}}
  >{\raggedright\arraybackslash}p{(\linewidth - 6\tabcolsep) * \real{0.2500}}@{}}
\toprule\noalign{}
\begin{minipage}[b]{\linewidth}\raggedright
\(f\)
\end{minipage} & \begin{minipage}[b]{\linewidth}\raggedright
Leading coefficient of \(F\)
\end{minipage} & \begin{minipage}[b]{\linewidth}\raggedright
\(F(r,h/2)\)
\end{minipage} & \begin{minipage}[b]{\linewidth}\raggedright
\(F(r,r)\)
\end{minipage} \\
\midrule\noalign{}
\endhead
\bottomrule\noalign{}
\endlastfoot
\(1/4\) & \(-448r(2r-1)\) & \(r(4r-3)(72r^2-86r+21)\) & \((4r-3)(8r-3)W\) \\
\(3/4\) & \(-64r(14r-5)\) & \(r(4r-1)(72r^2-50r+5)\) & \((2r-1)(4r-1)(8r-1)(32r-5)\) \\
\end{longtable}
}

Throughout each claimed interval strictly below its (K) endpoint these signs are \(+,+,-\). The polynomial \(W\) has its smaller zero at \((39-\sqrt{177})/64\) and is positive to its left here. The remaining quadratic signs are certified below. Thus there is a unique root in \((h/2,r)\), and it varies continuously. Take that root for \(k\).

Write \(N=(1-r)w\). The difference \(A-N\) equals the (K) polynomial \(8r^2-3(4-3f)r+2(2-f)(1-f)\) and is positive below its smaller zero on our intervals. Two useful bounds now hold immediately: \[
w(1-w)-E=\frac{(k-w)(A-N)}{r-k}>0,\qquad
Z-J=\frac{w(k-w)[A-N+r(r-k)]}{r(r-k)}>0.
\] The remaining tests are \(E\ge0\), \(J>B\), \(V_M\ge A\) and \(A>V_H\).

Take the resultant of \(F\) with the reduced numerator of each test, eliminating \(k\). The table gives the factors that can determine a sign on either claimed interval below its (K) endpoint. Discarded linear factors are among \(r\), \(2r-1\), \(4r-1\), \(4r-3\), \(8r-1\), \(8r-3\), \(10r-1\), \(10r-3\), \(12r-5\), \(12r-7\), \(14r-5\) and \(32r-5\), in the rows where they are nonzero; the first row also has discarded powers of \(W\). All multipliers can be checked directly from the displayed expressions for \(P,B,V_M,V_H\).

{\def\LTcaptype{none} 
\begin{longtable}[]{@{}
  >{\raggedright\arraybackslash}p{(\linewidth - 8\tabcolsep) * \real{0.2000}}
  >{\raggedright\arraybackslash}p{(\linewidth - 8\tabcolsep) * \real{0.2000}}
  >{\raggedright\arraybackslash}p{(\linewidth - 8\tabcolsep) * \real{0.2000}}
  >{\raggedright\arraybackslash}p{(\linewidth - 8\tabcolsep) * \real{0.2000}}
  >{\raggedright\arraybackslash}p{(\linewidth - 8\tabcolsep) * \real{0.2000}}@{}}
\toprule\noalign{}
\begin{minipage}[b]{\linewidth}\raggedright
\(f\)
\end{minipage} & \begin{minipage}[b]{\linewidth}\raggedright
\(E\)
\end{minipage} & \begin{minipage}[b]{\linewidth}\raggedright
\(J-B\)
\end{minipage} & \begin{minipage}[b]{\linewidth}\raggedright
\(V_M-A\)
\end{minipage} & \begin{minipage}[b]{\linewidth}\raggedright
\(A-V_H\)
\end{minipage} \\
\midrule\noalign{}
\endhead
\bottomrule\noalign{}
\endlastfoot
\(1/4\) & \(S_1=30720r^4-61248r^3+43024r^2-12684r+1323\) & \(D_1=72r^2-86r+21\) & \(R_1=2688r^3-3856r^2+1584r-189\) & \(1\) \\
\(3/4\) & \(S_3=30720r^4-34240r^3+12432r^2-1700r+75\) & \(D_3=72r^2-50r+5\) & \(R_3=896r^3-656r^2+112r-5\) & \(1\) \\
\end{longtable}
}

For each polynomial on an interval \([a,b]\) use \(r=(a+bz)/(1+z)\) and clear the appropriate power of \(1+z\). In increasing degree, the coefficient signs for \((S_1,D_1,R_1)\) on \(I_1\) are \(-++++,, ,, -\). On \(I_3\) the signs for \((S_3,D_3,R_3)\) are \(-++++,, ,+, \). Descartes' rule and opposite endpoint signs give the claimed unique roots \(\xi\) of \(S_1\) and \(\eta\) of \(R_3\). Direct substitution gives \(S_3(61/400)>0\) and \(R_3(61/400)>0\), putting the unique zero of \(S_3\) below \(61/400<\eta\). Also \(S_1(321/800)>0\) and \(R_3(31/200)<0\), putting these respective test points inside the claimed intervals.

The table fixes the signs of the four tests without relying on decimal root values. At each sample \(F\) has signs \(+,-\) at the two indicated \(k\) endpoints, bracketing its smaller root. On each \(k\) bracket the four reduced numerators of \(E,J-B,V_M-A,A-V_H\), with denominators oriented positive, have increasing-degree coefficient signs \(++\) after the same interval substitution.

{\def\LTcaptype{none} 
\begin{longtable}[]{@{}lll@{}}
\toprule\noalign{}
\(f\) & Sample pole & Bracket for \(k\) \\
\midrule\noalign{}
\endhead
\bottomrule\noalign{}
\endlastfoot
\(1/4\) & \(321/800\) & \((39341,39342)/100000\) \\
\(3/4\) & \(31/200\) & \((14525,14526)/100000\) \\
\end{longtable}
}

The resultants now give the required strict signs in both open intervals. Continuity includes each lower endpoint with the allowed weak sign there; \(J>B\) and \(A>V_H\) persist there. Zero-share plies may be omitted. At the (K) endpoints the preceding five-stage construction realizes the poles directly.

At \(f=1/4\), all pole parameters \(r\) in the closed interval \[
 [x,y]=[0.4059648800755\ldots,\,0.4064113242525\ldots]
\] are realized by isotropic two-phase laminates on at most nine separated scales. These are the unique roots in \(I=[4059/10000,40642/100000]\) of the indicated polynomials \[
 S(x)=2432x^4-6096x^3+5652x^2-2295x+342=0,\qquad
 R(y)=736y^3-1084y^2+474y-63=0.
\] A graded spherical coating continues this interval farther left to \(0.4053009437701\ldots\).

The entire interval is below \((5-\sqrt3)/8\) and above the preceding (K) point \((39-\sqrt{177})/64\). This construction gives a full interval in the region that already contains a sequence of distinct finite contacts converging to \((5-\sqrt3)/8\).

\paragraph*{Geometry}

Use the upper ply in the preceding matching construction, retaining its phase labels and setting \(h=3r-3/4\) for \(f=1/4\). The following parameters specify all the shares. Put \(t=1-h,\ w=h-r,\ d=r-w\). Choose the smaller root \(k\) of \[
 \begin{aligned}
 P(r,k)={}&(256r^2-384r+168)k^2\\
 &+(-5120r^3+8256r^2-4176r+630)k\\
 &+6912r^4-12672r^3+8064r^2-2088r+189=0.
 \end{aligned}
\] Set \(l=h-k\), \(g=(k-r)/(k-l)\), \(A=(k-w)(t-h)\), \(B=(1-r)l-A\), \(e=B/g\), \(j=(1-k)w\), \[
 Z=\frac{l[(1-r)k-(1-g)e]+(k-l)tw}{k},\qquad
 \sigma=\frac{k-r}{k-w},\qquad \mu=1-\sigma,\qquad N=(1-r)w .
\] Make the three-stage blocks \(G_x(r;l,k;e)\) and \(K_y(r,k)\) in the axes and with the simple-ply shares given in the upper-ply construction. Their poles are \((r,l,k)\) and \((r,h,k)\) and their \(y\) coefficients agree at \(A\).

\begin{enumerate}
\def\labelenumi{\arabic{enumi}.}
\tightlist
\item
  Stack these blocks normal to \(y\) in shares \(l/k\) and \(1-l/k\). The result \(C\) has poles \((r,k,k)\), \(q_y=A\) and \(q_z=Z\).
\item
  Stack \(C\) and its copy exchanging \(y,z\) normal to \(x\) with original share \((j-A)/(Z-A)\). This gives \(q_z=j\) and \(q_y=A+Z-j\).
\item
  Stack the simple \(L_y(1-k)\) and this mixture normal to \(z\) in shares \(\sigma,\mu\). The simple has poles \((r,k,w)\) and \(q_z=j\). The resulting upper ply \(M\) has poles \((r,k,r)\), \(q_z=j\), and \(q_y=\mu(A+Z-j)\). Direct expansion gives \(128k(k-r)[N-\mu(A+Z-j)]=P(r,k)\), so \(q_y=N\).
\item
  Stack the simple \(L_x(1-r)\) at poles \((r,w,r)\) with \(M\) normal to \(y\) in shares \(\sigma,\mu\). The simple has the same \(q_y=N\). All three poles are now \(r\) and the phase fraction is \(1+h-3r=1/4\). Two further laminations of rotated copies give \(K_{1/4}(r,s)I_3\).
\end{enumerate}

Here \(L_i(c)\) denotes layers of the two pure phases with normal \(i\) and phase-1 share \(c\). Each tensor step uses a harmonic mean normal to the layers and arithmetic means in the two tangent directions. Realization is by reiterated homogenization.

\paragraph*{Closed-interval certificate}

Throughout \(I\), \(P(r,r)=(4r-3)(8r-3)(64r^2-78r+21)>0\), the leading coefficient is positive, and \(P(r,87/200)<0\). Thus a unique root lies between \(r\) and \(87/200<h\). It varies continuously. This gives \(0<l<w<r<k<h<t\). The remaining share tests are \[
 B\ge0,\qquad D=(k-r)l(1-l)-(k-l)B\ge0,\qquad j-A>0,\qquad
 V=2h^2-3hk-2h-kr+3k>0.
\] Indeed \(Z-j=l(k-w)V/[k(k-r)]\). The four resultants with \(P\) in \(k\), up to constant factors and products of \(4r-1,4r-3,12r-7,64r^2-78r+21\), are \(R,S,J,1\) in the same order, where \(J=1024r^3-1424r^2+564r-63\). None of the discarded factors vanishes on \(I\).

For an exact sign check, if \(I=[a,b]\), replace any polynomial \(p(r)\) of degree \(n\) by \((1+z)^np((a+bz)/(1+z))\). The coefficient signs in increasing order are

{\def\LTcaptype{none} 
\begin{longtable}[]{@{}ll@{}}
\toprule\noalign{}
Polynomial & Signs \\
\midrule\noalign{}
\endhead
\bottomrule\noalign{}
\endlastfoot
\(P(r,87/200)\) & \(, -\) \\
\(S\) & \(+, -\) \\
\(R\) & \(+++-\) \\
\(J\) & \(, -\) \\
\(64r^2-78r+21\) & \(, \) \\
\end{longtable}
}

Hence \(S,R\) each have one root on \(I\), and \(J\) has none. At \(r=203/500\), \(S<0<R\), so this sample lies strictly between the two roots. There \(P\) takes opposite signs \(+,-\) at \(k=42552/100000,\ 42553/100000\), placing the selected root between them. Applying the same coefficient-sign test in this \(k\) bracket gives strictly positive coefficients for all four share tests. The resultants show their signs persist on \((x,y)\); by continuity the allowed weak inequalities hold at both ends, and the other two tests remain strictly positive. This proves the endpoint realizations.

\paragraph*{Propagating an interval of seeds}

The fourth stack can retain two prescribed coefficients. Split its total simple share \(\sigma\) between \(L_x(1-r)\) and \(L_z(1-r)\), with the proportion \(\theta=w(1-w)/[r(1-r)]\) going to \(L_x(1-r)\). Both have \(q_y=N\). Their combined part has \(q_z=w(1-w)\), and all poles and the total phase fraction stay as before. Then \[
 q_z=\mu(1-k)w+\sigma w(1-w)=N
\] because \(\mu k+\sigma w=r\). The ordering gives \(0<\theta<1\). After a permutation of axes the root therefore has \(p=(r,r,r)\), \(q=(N,N,(t-r)(2r-h))\). This observation applies to the same geometry at any admissible \(h\) and starting pole \(r\).

At fixed \(h\), every such seed above \(\kappa(h)\) may be passed repeatedly through the three decreasing stacks in the pole recurrence of Section~\ref{route:matching-plies-and-extending-their-poles}. More explicitly, \(T(a)=t-(1-a)(h-a)/(2a-h)\) is strictly increasing and the step sends a seed at \(k\) to \(a=T^{-1}(k)\). It preserves the two coefficients of the specified form. The image of an interval of seeds is again an interval. The phase fraction at the image is \(1+h-3a\). To construct an image at exactly \(1/4\), first set \(h=3a-3/4\) for the \emph{target} pole \(a\), then test the seed geometry at \(T^n(a)\) with that same \(h\). This keeps the phase fraction at the finished root exact.

For example, \(8151/20000\), \(40811/100000\) and \(20417/50000\) are reached at \(f=1/4\) after one, two and three such steps. They are interior to realized fixed-fraction intervals. In each case the quadratic has distinct roots and the seed share tests and propagation pole orderings are strict. Hence they persist in a neighborhood of the specified target pole.

The iterations of this seed do leave gaps. At the test pole \(a=407/1000\), \(h=471/1000\). For \(a,T(a),T^2(a)\) the matching quadratic fails as follows. For \(a\) its discriminant is negative. At \(T(a)\) its smaller root has \(D<0\) and its larger root has \(B<0\). At \(T^2(a)\) both roots violate \(T^2(a)<k<h\). Finally \(T^3(a)>h\) lies beyond the eligible increasing sequence of ancestors. This exhibits a gap in the set reached by propagating the specified seed.

Run \texttt{python\ follow-up/verify\_direct\_quarter.py} for the endpoint certificates, the checks at \(407/1000\), and the full rational contrast tensors of the four direct and propagated examples. To export a geometry, for example at the first propagated pole, run \texttt{python\ follow-up/verify\_direct\_quarter.py\ -\/-tree\ 8151/20000\ -\/-steps\ 1}.

At \(f=1/4\) an assemblage of spheres with a graded, radially oriented laminate coating extends the direct finite band to the left. Together the two geometries attain the entire closed interval

\[
 [\rho,y]=[0.4053009437701\ldots,\;0.4064113242525\ldots].
\]

Here \(\rho\) is the unique root in \([4053/10000,405965/10^6]\) of

\[
 2304\rho^4-5808\rho^3+5384\rho^2-2171\rho+320=0 .
\]

The coated spheres give \([\rho,x]\) up to \(x=0.4059648800755\ldots\), and the finite direct laminate gives \([x,y]\). Each coating has an explicit effective tensor and a prescribed radial profile. Its laminate plies have periods small relative to the radius on which their parameters vary.

\paragraph*{A spherical continuation rule}

Keep \(h\) fixed. Let an isotropic core with single-pole parameter \(a\) and phase fraction \(f_a=1+h-3a\) be available. Suppose that for every \(\beta\in[a,b]\) we have a uniaxial shell material with radial and tangential pole triple \((k(\beta),\beta,\beta)\) and radial quadratic coefficient

\[
 m(\beta)=h(1+h-3\beta)/3.
\]

Assume \(k\) is continuous and \(k(\beta)>\beta\) on this interval. Its local phase fraction is \(c=1+h-k-2\beta\) and each tangential quadratic coefficient is \(q_t=(ch-m)/2\). A shell with these data continues the isotropic core all the way to \(b\). Use outer radius one and set

\[
 R(\beta)=\exp\!\left[-\int_\beta^b\frac{dt}{k(t)-t}\right].
 \tag{1}
\]

Use the core inside \(R(a)\) and the indicated radially oriented material at radius \(R(\beta)\). In particular \(d\beta/d\log R=k-\beta\). If \(a=b\), use only the core. The overall fraction inside radius \(R(\beta)\) is \(1+h-3\beta\). Indeed the derivative of this proposed fraction with respect to \(\log R\) equals \(-3(k-\beta)=3[c-(1+h-3\beta)]\), exactly the derivative of the volume fraction in a growing three-dimensional ball, and its initial value is \(f_a\).

To verify the tensor, put \(D_p=1+ps\) and impose the boundary potential \(\xi\cdot x\). In the shell take the potential

\[
 \frac{D_{\beta(R)}}{D_b}\,\xi\cdot x .
\]

The radial and tangential field amplitudes relative to the corresponding components of \(\xi\) are \(D_k/D_b\) and \(D_\beta/D_b\). In the core use the affine potential \(D_a\,\xi\cdot x/D_b\). The shell current amplitudes \(J_r,J_t\) satisfy

\[
 D_b J_r=1+(1+h-2\beta)s+m(\beta)s^2,\qquad
 D_b J_t=1+(c+\beta)s+q_t s^2 .
\]

Now \(dJ_r/d\log R=2(J_t-J_r)\) since \(m'=-h\) and \(2(q_t-m)=h(\beta-k)\). This identity is exactly zero divergence for the current \(J_r(\xi\cdot n)n+J_t[\xi-(\xi\cdot n)n]\) with \(n=x/|x|\). The potential and radial current agree with the core at \(R(a)\). At radius one the radial flux equals \(K_{1+h-3b}(b,s)\,\xi\cdot n\). The sphere is therefore neutral in an exterior medium of that scalar conductivity for every positive real phase contrast.

Scaled copies fill a coated-sphere assemblage. For example, pack a periodic cube up to a remainder of arbitrarily small volume with similar balls, initially putting the target scalar in the remainder. The explicit potentials give the exact periodic correctors, with uniform exterior fields and currents. Replacing the remainder by either physical phase changes the effective response by a vanishing amount, by the electric and dual current energy principles using these fields and currents as trials. The shell can itself be approximated by piecewise constant rotations and parameters, subdividing into thin spherical layers and angular sectors, then homogenizing the indicated laminates on finer scales. On each compact positive contrast interval their tensors are uniformly elliptic, and continuity of the rotated effective tensors gives convergence of these approximations by the energy principle. The resulting phase fractions converge to the indicated fraction. Any small discrepancy at a prescribed fraction can be removed by a further coordinate stack with a vanishing share of the deficient pure phase. A contrast-independent diagonal sequence gives the claimed response, also locally uniformly in the usual cut complex contrast plane.

\paragraph*{The core and shell for the new band}

Given target \(r\in[\rho,x]\), put \(h=3r-3/4,\ t=1-h\) and

\[
 T(z)=t-\frac{(1-z)(h-z)}{2z-h},\qquad a=T^{-1}(T^{-1}(h)),
\]

using the increasing branch on \((h/2,h]\). The finite upper-sequence construction supplies the isotropic core at \(a\), with fraction \(1+h-3a\). Take \(b=r\) in (1). Indeed \(T(r)\in(h/2,h)\) on the interval in question and

\[
 T^2(r)-h=-\frac{2432r^4-6096r^3+5652r^2-2295r+342}
 {(4r-3)(208r^2-264r+75)}\ge0
\]

there. The numerator is the defining polynomial of \(x\) and is nonnegative up to \(x\). Thus \(r\ge a\).

Use steps 1-3 of the direct upper-ply construction with its pole variable replaced by \(\beta\), at the fixed \(h\). The upper ply has poles \((\beta,k,\beta)\), and its middle coefficient is \(\mu(A+Z-j)\) in the notation of those steps. Now choose \(k\) to make this coefficient \(m(\beta)\). Laminate the upper ply in equal shares with its copy exchanging the two tangential axes, normal to the middle axis. This produces exactly the uniaxial shell required in (1). Orient the middle axis radially.

For explicit evaluation, with \(w=h-\beta\), \(l=h-k\), define

\[
 \begin{aligned}
 A&=(k-w)(1-2h), & B&=(1-\beta)l-A,\\
 D&=(k-\beta)l(1-l)-(k-l)B, & j&=(1-k)w,\\
 V&=2h^2-3hk-2h-k\beta+3k,\qquad & N&=(1-\beta)w .
 \end{aligned}
\]

The smaller root of the following quadratic gives \(k(\beta)\):

\[
 \begin{aligned}
 F={}&a_2k^2+a_1k+a_0+k(k-\beta)(m-N),\\
 a_2={}&h^2-(2\beta+1)h-\beta^2+3\beta,\\
 a_1={}&-5h^3+(11\beta+5)h^2-(\beta^2+12\beta)h-\beta^3+3\beta^2,\\
 a_0={}&2h^4-(4\beta+2)h^3+4\beta h^2 .
 \end{aligned}
\]

The identity \(F=k(k-\beta)[m-\mu(A+Z-j)]\) proves the match. The remaining conditions are \(0<l<w<\beta<k<h<t\), \(B,D\ge0\), and \(j-A,V>0\). They ensure all layer shares are physical. At an endpoint a zero-share constituent may be dropped.

\paragraph*{Exact sign certificate}

The verifier checks these conditions on all paths \((h,\beta)\) from the core to the boundary. It covers the wider range \(4053/10000\le r\le405965/10^6\) and tests \(D\) at the outer boundary separately. On each parameter rectangle it isolates the smaller root in a rational \(k\) bracket by checking that the leading coefficient is positive and \(F\) has opposite signs there. It checks the pole order, \(B>0\), \(j-A>0\), \(V>0\) and

\[
 D_\beta F_k-D_k F_\beta>0
\]

throughout that bracket. As \(F_k<0\) at the smaller root, \(D\) decreases along this branch with increasing \(\beta\) and fixed \(h\). It therefore suffices to check \(D\ge0\) at \(\beta=r\).

At the outer boundary \(h=3r-3/4,\ \beta=r\), the resultant of \(F,D\) in \(k\) is

\[
 \begin{gathered}\frac{3(4r-3)(4r-1)(12r-7)(64r^2-78r+21)}{1048576}\\ {}\times
 (896r^3-1408r^2+680r-93)
 (2304r^4-5808r^3+5384r^2-2171r+320).\end{gathered}
\]

The last quartic has a single root \(\rho\) in the wider range, in \((4053/10000,40531/100000)\), and the other factors have none there. At \(r=811/2000\), \(k\) lies in \((42313/10^5,42314/10^5)\) and \(D>0\) throughout this bracket. The resultant and continuity give \(D>0\) at the outer boundary on \((\rho,x]\) and \(D\ge0\) at \(\rho\).

The polynomial-box calculation uses exact integer arithmetic. For each axis with endpoints \(L_i/M,U_i/M\), \(M=10^8\), it substitutes \(X_i=[L_i+U_i+(U_i-L_i)y_i]/(2M)\) with \(|y_i|\le1\). A positive constant coefficient greater than the sum of the absolute values of the other coefficients proves positivity over the box. Floating-point root estimates only propose rational brackets. The proof uses both exact sign tests at each bracket. To bound the core in a range of \(h\), note that \(T\) increases in its argument and decreases with \(h\) since \(\partial_h T(z)=-1-z(1-z)/(2z-h)^2<0\). It follows that \(T^{-2}(h)\) increases with \(h\). Two bisections at a slab's lower endpoint give a rigorous lower bound throughout it. The boxes then cover all \(\beta\) up to the target at the slab's upper endpoint.

Run \texttt{python\ follow-up/verify\_radial\_shell.py} to reproduce the sign proof and the shell-matching identities.
\phantomsection\label{endfile:short-bands}

\subsection{The smooth-interface question}\label{route:the-smooth-interface-question}
For a single cell with finitely many disjoint smooth closed interface
surfaces the following argument leads to the two coated-sphere poles.
Its conclusion conflicts with the identification of the fixed Schwarz
P and D cells in \cite[Eq.~(30)]{tk24}. For example, at $f=1/2$ the
two proposed smooth-cell endpoints are $1/6$ and $5/6$, while that
formula has pole $1/2$. The earlier transport work \cite{thd02} studies
maximization of the two complementary transport properties. A numerical
agreement with an extremal value to fixed tolerance cannot by itself
identify a whole conductivity function.

We give the argument and draw attention to two independent steps for
review: the passage from reciprocal field constraints to interface
normal derivatives, and the winding-rank classification of a component
with fewer than three independent periods. We have not used this
classification to exclude any general physical pole in the article.
\phantomsection\label{file:smooth-contact}
\paragraph*{The interface equations}

Normalize the phase conductivities to \(1+s,1\). Write \(\chi\) for the phase-1 indicator, \(\langle\cdot\rangle\) for a normalized cell average, and \(\Gamma\) for the orthogonal projection onto periodic mean-zero gradients. Put \(H=\Gamma\chi\Gamma\) on this gradient space and \(g_j=\Gamma(\chi e_j)\). Solving the projected cell problem for the three coordinate loadings gives

\[
\frac{\operatorname{Tr}K_*(s)}3
=1+fs-\frac{s^2}{3}\sum_{j=1}^3\langle(I+sH)^{-1}g_j,g_j\rangle_{L^2}.
\]

The positive spectral measures of the \(g_j\) sum to a measure of mass \(f(1-f)\). Indeed, at each nonzero Fourier frequency the gradient projector has rank one, so \(\sum_j\|g_j\|_2^2=\|\chi-f\|_2^2\). The normalized trace is therefore a probability average of

\[
K_f(p,s)=1+fs-\frac{f(1-f)s^2}{3(1+ps)},\qquad 0\le p\le 1.
\]

At nonreal \(s\) these values lie on a circle, and the parametrization is injective. Every point of the circle is an extreme point of its disk. Equality with \(K_f(u,s)\) concentrates each spectral measure at \(u\), giving \(Hg_j=u g_j\). Equality as a contrast function gives the same identities.

Let \(\phi\) be the periodic mean-zero solution of \(\Delta\phi=\chi-f\) and put \(P=uI-D^2\phi\). The Fourier formula gives \(D^2\phi\,e_j=g_j\). Thus \(\Gamma(\chi P)=0\) column by column, which means that \(\chi P\) is divergence-free. The Poisson solution has continuous first derivatives and is smooth up to either side of each smooth interface by Poisson transmission regularity. Write \(P_i\) for its one-sided matrices, with the normal \(n\) directed from phase 1 into phase 2. The zero normal flux of \(\chi P\) gives \(P_1n=0\). Tangential derivatives of the continuous gradient of \(\phi\) agree across the interface. Symmetry and \(\Delta\phi_1-\Delta\phi_2=1\) then give

\[
P_2-P_1=nn^T,\qquad P_1n=0,\qquad P_2n=n.
\]

Within each phase the entries of \(P\) are harmonic and its trace is constant. Third derivatives of the local potential are symmetric. Write \(c_1=0,c_2=1\). For a fixed vector \(v\) tangent at an interface point, differentiating \(P_i n=c_i n\) along \(v\) gives there

\[
\partial_n(v^TP_iv)
=v^T(\partial_v P_i)n
=v^T(c_iI-P_i)\,\partial_v n.
\tag{1}
\]

\paragraph*{Positivity from the two Hopf derivatives}

Every one-sided matrix satisfies \(0\preceq P\preceq I\). To prove this, suppose the globally least eigenvalue \(\lambda\) is negative, or the globally greatest one exceeds one. Take an extremizing unit vector \(v\). If the extremum occurs inside a phase component, harmonicity and the strong maximum principle make \(v^TPv=\lambda\) throughout the component, so it also occurs on its boundary. At an interface extremum \(v\) is tangent since the normal has eigenvalue \(c_i\). The jump \(nn^T\) shows that \(v\) is the same \(\lambda\)-eigenvector on both sides.

The functions \(v^TP_i v\) have the same global extremum there. Their one-sided normal derivatives in the common direction \(n\) must have opposite weak signs. By (1) they are \((c_i-\lambda)v\cdot\partial_v n\). The two factors \(c_i-\lambda\) have the same strict sign, so both derivatives vanish. The Hopf boundary point lemma forces \(v^TP_i v\) to be constant in both adjacent components. In each such component \(v\) is a \(\lambda\)-eigenvector everywhere, since \(\lambda\) is a global extremal eigenvalue. The argument crosses every boundary of the component with the same \(v\). The component adjacency graph is connected, hence \(v^TPv=\lambda\) throughout the cell. Its mean equals \(u\) by periodicity, a contradiction.

\paragraph*{Curvature when both phases wind in all three directions}

For a phase component \(C\), let \(W_C\) be the span of the lattice displacements of lifts of closed paths in \(C\). Suppose first that \(W_C=\mathbb R^3\) in every component. Set \(T_1=P_1\) and \(T_2=I-P_2\). These are positive semidefinite, annihilate \(n\) at the interface, and are positive definite throughout the interiors. In fact, if an interior quadratic form \(v^TT_i v\) vanished for a unit vector \(v\), harmonicity and positivity would give \(T_i v=0\) in the whole component. On its lift, \(T_i\) is the derivative of \(u x-\nabla\phi\) in phase 1 or \((1-u)x+\nabla\phi\) in phase 2. Symmetry then makes the \(v\)-component of this map constant along each lifted path. On a loop with lattice displacement \(\ell\) this gives \(u\,v\cdot\ell=0\) or \((1-u)v\cdot\ell=0\). The prescribed pole range has \(0<u<1\), contradicting full winding.

Fix an interface point, freeze its unit normal \(n_0=n\), and use orthonormal tangent principal directions \(e_1,e_2\) with \(\partial_{e_a}n=\kappa_a e_a\). The harmonic functions \(n_0^TT_i n_0\) are strictly positive in their respective interiors and vanish at this point. Hopf's lemma gives

\[
\partial_n(n_0^TT_1n_0)<0,\qquad \partial_n(n_0^TT_2n_0)>0.
\]

For either \(T_i\), differentiating \(T_i n=0\) as in (1) yields \(\partial_n(e_a^TT_i e_a)=-\kappa_a e_a^TT_i e_a\), with \(e_a\) frozen in the quadratic form. Taking the normal derivative of its constant trace gives

\[
\partial_n(n_0^TT_i n_0)=\sum_{a=1}^2\kappa_a e_a^TT_i e_a.
\]

The nonnegative weights in these sums and their opposite strict signs force \(\kappa_1\kappa_2<0\) at every interface point. On any connected interface surface the oriented Gauss map into the unit sphere is therefore a local diffeomorphism. Compactness makes it a covering of \(S^2\), and simple connectivity makes the surface a sphere. Gauss-Bonnet now gives total Gaussian curvature \(4\pi\), contradicting strictly negative Gaussian curvature everywhere. Thus at least one phase component has winding rank less than three.

\paragraph*{Consequence for the arc}

The winding proposition classifies the possible contacts when such a lower-rank component exists. In phase 1 the pole is \((1-f)/3\), \((1-f)/2\), or \(1-f\) for winding rank zero, one, or two, respectively. In phase 2 the corresponding values are \(1-f/3\), \(1-f/2\), and \(1-f\). In the two positive-rank cases the geometry is invariant along a nonzero direction. The conductivity in that direction equals the arithmetic mean \(1+fs\), which differs from \(K_f(u,s)\) for \(s\ne0\) wherever the formula is defined. Thus a fixed smooth cell can make isotropic contact only at a coated-sphere endpoint. For the normalized trace alone the five parameters just listed are the only possible contacts.
\phantomsection\label{endfile:smooth-contact}

\SSection{inclusions}{Planar boundary certificates and a spatial obstruction}
\subsection{The quarter-area four-lobed solid}\label{route:the-quarter-area-four-lobed-solid}
This continues the boundary-integral proof of Theorem~\ref{p7:star}.
We use the notation for $g,W$, the interior trace and the lattice
kernel introduced there.
\phantomsection\label{file:quarter-certificate}
\smallskip\noindent\emph{Correcting a trial boundary.}
Use a cosine polynomial $R=1+10^{-20}\sum_{j=1}^{76}n_j\cos(4j\theta)$,
$B_0=\langle R^2\rangle_\theta$ and $\delta=(4\pi B_0)^{-1/2}$.
The integers $n_j$ specify the reference profile. The inverse polynomials and finite inverse columns below are dyadic proposals. Their residual bounds certify them independently of the method used to choose their coefficients.
The map whose zero we seek is
\[
 H(h)=\big(\langle(R+h)^2\rangle_\theta-B_0,\ \operatorname{Re}W-r\big).
\]
The first entry fixes the physical area. In the Wiener algebra set
$\|p\|_v=\sum |p_k|v^{|k|}$ (product weight in two variables).
With $\sigma=1+2^{-20}$ use
$\|h\|_X=\sum_{k=0,4,\dots}(1+k)\sigma^k|h_k|$ for real cosine
coefficients, and $\|(a,b)\|_Y=|a|+\|b\|_\sigma$. Even cosine series
have the same weighted sum as their Laurent coefficients.

At the reference profile define
\[
 L_p=2Rp^\#/G^2,\quad D_p=2G^{-1}N(Rp^\#/G),\quad
 S_p h=L_pNh+D_p h,\quad V=wNW/G .
\]
Then
\begin{equation}\label{p7:derivative}
 H'(0)h=\big(2\langle Rh\rangle_\theta, \operatorname{Re}\{\mathcal P(S_gh+M(\mathcal P S_1h)^\#)+Vh\}\big).
\end{equation}
Indeed, varying the area integral for $\mathcal Kp$ at fixed interior $z$,
holding a holomorphic input fixed, adds
$\pi^{-1}\int Rh p^\# \wp_\delta(\gamma-z)\,d\theta$.
Integrating by parts with $\partial_\theta L_\delta=-iG\wp_\delta$ gives
$\mathcal PS_ph$ in the boundary limit. Apply this rule to both integrals and account for
moving the final evaluation point by adding $Vh$. These differentiations hold
first for trigonometric polynomial directions because analytic contours and
densities have holomorphic continuations across the reference curve by the
Cauchy jump formula, permitting evaluation on a common inner collar and
passage to its boundary. The norm bounds that follow extend the identity
by density.

The existence certificate supplies a real operator $A:Y\to X$ with
\begin{align}
 \|H(0)\|_Y&<8.06\,10^{-14},&
 \|A\|&<16.22,&
 \|I-H'(0)A\|&<0.147,\label{p7:cert1}\\
 \sup_{\|h\|_X<1/100}\|H''(h)\|&<1.30\,10^8 .\label{p7:cert2}
\end{align}
The finite inequalities in this proof are asserted in \path{certificates/baseline/p07/certify_star_solid.py}. At $h=0$, set
$d=\gamma(v)-\gamma(w)=(v-w)D(w,v)$. Subtracting the unit-circle
Cauchy kernel gives $\mathcal Pp=P_+p+[K(w,v)p(v)]_{v^0}$, where $P_+$
is nonnegative Laurent projection and
\[
 K=\frac{v\partial_vD}{D}+G(v)\big(L_\delta(d)-1/d\big).
\]
At $\omega=21/20$ a dyadic polynomial $U$ of bandwidth $156$ satisfies
$\|1-DU\|_\omega<0.027$, $\|U\|_\omega<7.73$ by integer multiplication, giving $\|1/D\|_\omega<12$. On the product annulus of weight $\omega$,
$|\delta d|\le 2\delta\omega[1+\tfrac12\sum|10^{-20}n_j|
(\omega^{4j}+\omega^{-4j})]<0.983$. For
$G_{4m}=\sum_{\ell\ne0}\ell^{-4m}$, $|G_{4m}|\le8\sum_{j\ge1}j^{-3}<12$
by square shells. The formula
$L_\delta(d)-1/d=-\delta\sum_{m\ge1}G_{4m}(\delta d)^{4m-1}$
and a termwise bound on $v\partial_vD$ give $|K|<3000$ there.

Use $\tau=1+1/512$ and $N_0=3072$. The Cauchy estimate
$|K_{ij}|\le3000\omega^{-|i|-|j|}$ bounds the error from deleting all
indices outside the centered $N_0$-square and folding them onto the
retained indices by
$8\cdot3000\,(1+x)x^{N_0/2}/(1-x)^2$ in the $\tau$ norm,
$x=\tau/\omega$ (sum the two coordinate tails and charge their aliases).
On the $N_0^2$ grid the removable plane-kernel diagonal is
$(N^2-N)\gamma/(2G)$. The regular part uses
\[
 -\frac{\delta}{1-y^4}\sum_{m=1}^{70}
 (G_{4m}-G_{4m-4})y^{4m-1},\qquad y=\delta d,\qquad G_0=0.
\]
For $m\ge2$ the four nearest lattice points cancel in the difference
and the others have length at least $\sqrt2$, giving
$|G_{4m}-G_{4m-4}|\le24/4^{m-2}$. With $|y|<0.92$ this bounds each
sample's omitted tail by a geometric series. The normalized transform
of a sample error has each coefficient bounded by its sup norm. Rotation
retains only $i+j=0\pmod4$ and reflection gives real coefficients.
The sum of tail, alias and dyadic rounding errors in the $\tau$ norm is
below $5\,10^{-24}$. The lattice coefficients may be evaluated from
$c_j=(4j-1)G_{4j}$ by
\[
 G_4=\frac{\pi^4}{45}\left(1+240\sum_{k\ge1}
 \frac{k^3e^{-2\pi k}}{1-e^{-2\pi k}}\right),\qquad
 c_j=\frac{6\sum_{i=1}^{j-1}c_i c_{j-i}}{(4j-2)(4j-3)-12}
 \quad(j\ge2).
\]
The first identity follows by summing $\sum_m(m+z)^{-4}$ along horizontal
lattice rows using $-(\pi\cot\pi z)'''/6$. To obtain the second, compare
coefficients in $\wp''=6\wp^2-30G_4$ (the difference is entire periodic
with that constant term). Enclosed convolution and Laurent projection in
the $\tau$ norm now give the residual in \eqref{p7:cert1}.

In \eqref{p7:derivative} truncate the multipliers, $R$ and $K$ to
bandwidth $d_0=256$, rounding coefficients to $2^{-48}$ and keeping $P_+$ exact.
Call the resulting derivative $D_0$. Differentiating a $\tau$-tail costs
at most $\sigma/(\tau-\sigma)<513$ at weight $\sigma$.
The norm of the finite kernel on densities is bounded by
$\max_j\sum_i|K_{ij}|\sigma^{|i|}$, and its error by the two-variable
Wiener error. Apply the product inequality and
$\|Nh\|_\sigma+\|h\|_\sigma=\|h\|_X$ to each factor in
\eqref{p7:derivative} to obtain
$\|H'(0)-D_0\|<6.73\,10^{-6}$.

For the area and cosine input columns of $A$ through frequency $1024$
use trial inverse columns supported through $1280$, rounded to $2^{-55}$.
Integer products check all rows of $D_0A$ through $1280+2d_0$, which bounds all output frequencies. On these columns $\|A\|<16.22$ and
$\|I-D_0A\|<2\,10^{-12}$. On each $\cos m\theta$ for $m>1024$ use
$A(\cos m\theta)=\sum_j s_j\cos((m+j)\theta)/(m+j)$ with
$S=\sum s_jw^j$ a dyadic trial inverse to $\widehat L_g/2$,
supported on multiples of four with $|j|\le112$, where hats indicate rounded
multipliers. All input frequencies exceed $3d_0$. Consequently the regular
kernels vanish on them and their shifted bands, and in the nested term
the first projection retains the positive band, which reflection reverses
and the second projection kills. Comparing the remaining positive bands
and their reflected contributions in the real part bounds the weighted
column defect by
\[
 \|1-\widehat L_g S/2\|_\sigma+
 \frac{(\|\widehat D_g\|_\sigma/2+\|\widehat V\|_\sigma)
 \|S\|_\sigma}{1024-112}.
\]
This expression plus $\|H'(0)-D_0\|$ times the high-column bound on
$\|A\|$ is checked below $0.147$. The high-column norm is at most
$(1+1/(1024-112))\|S\|_\sigma<15.25$. Taking the supremum of normalized
column sums gives the claimed bounds for the whole operator.

For the nonlinear estimate complexify the cosine coefficients, extending
conjugation by $p\mapsto p^\#$ and real-part extraction by
$(p+p^\#)/2$. Set $v=1/20$. On $\|h\|_X\le v$ the divided difference
changes by at most $\sigma^2v$, since each monomial $w^{k+1}$ produces
$|k+1|$ terms with weights at most $\sigma^{|k|+2}$.
Thus its inverse has norm at most $D_*=12/(1-12\sigma^2v)$.
On this ball the tangent factors $(\gamma^\#)'$ and $N\gamma^\#$ are
bounded by
\[
 T_* = \sigma^2\left(v+\sum_k|k(wR)_k|\sigma^{|k|}\right).
\]
For the perturbation, the coefficient at frequency $k+1$ is $h_k$ in
Laurent coordinates, so the factor $|k+1|$ costs at most $1+|k|$.
In the plane part of the contour integral for $\mathcal K$, divide the
density $(\gamma^\#)'p^\#$ by the divided difference and apply the
unit-circle Cauchy integral, which sends $w^iv^j$ to $w^{i+j}$ for $j\ge0$
and zero otherwise, with norm at most one. In the regular part
$\|\delta(\gamma(v)-\gamma(w))\|_\sigma\le
u_*=2\delta\sigma(\|R\|_\sigma+v)<1$. Thus
\[
 \|\mathcal K\|\le F_*=T_*\left(D_*+
 \frac{12\delta u_*^3}{1-u_*^4}\right),\qquad
 \|H\|_Y\le(\|R\|_\sigma+v)^2+B_0+F_*^2+r .
\]
The norm-convergent products, inverses and series give a holomorphic
map on the ball. For real $h$ there, invertibility of the divided difference gives an injective regular boundary parametrization. The bound on $\delta d$ excludes identifications by a nonzero period. The radial bounds below and $\delta\|h\|_X\le0.267/20$ keep the curve inside the square with positive radius. On this Jordan curve the subtracted Cauchy kernel computes the interior trace by the same boundary integral, so the two entries of $H$ are the physical area constraint and $\operatorname{Re}W-r$. At any point of norm less than $1/100$, the bidisk
of radii $1/50$ in two arbitrary unit directions lies inside it.
Cauchy's formula for the mixed derivative bounds $\|H''\|$ by
$50^2$ times the last upper bound. The certificate asserts that this
product, even using $0.267$ for $\delta$, is below $130000000$.

Now $y\mapsto y-H(Ay)$ on $\|y\|_Y\le10^{-12}$ has Lipschitz constant
less than $0.182$ by \eqref{p7:cert1}-\eqref{p7:cert2}. Its displacement
at zero is below $8.06\,10^{-14}$, so it maps the ball into itself.
Iteration from zero converges to a fixed point, giving
$h=Ay$ and $H(h)=0$. In particular
$\|h\|_X\le 16.22(8.06\,10^{-14})/(1-0.182)$.
The finite cosine polynomial, written in $x=\cos4\theta$,
has derivative greater than $1/20$ throughout $[-1,1]$, as the certificate verifies on a mesh of spacing $2/1024$, with Lipschitz allowance
$10^{-20}\sum_j|n_j|j^2(j-1)^2/1024$, using
$T_j'=jU_{j-1}$ and the finite exponential sum for $U_{j-1}$.
Consequently the extrema of $R$ are its two endpoint values. Using $\|h\|_\infty\le\|h\|_X$ as their correction allowance, the certificate asserts $0.175<\delta\min(R+h)\le\delta\max(R+h)<0.447$. The radius is strictly positive and clears the square
faces. Exponential summability at $\sigma>1$ gives a real-analytic
boundary. This proves the theorem.
\phantomsection\label{endfile:quarter-certificate}

\subsection{Families of square-periodic inclusions}\label{route:families-of-square-periodic-inclusions}
\phantomsection\label{file:families}
\paragraph*{Setting}

Fix \(b>0\) and the lattice scale \(\delta=(4\pi b)^{-1/2}\). A positive even function \(x(\theta)=\sum_{j\ge0}x_j\cos4j\theta\) defines the solid with polar radius \(\delta\,x(\theta)\) in the unit square cell. Let

\[
\|x\|_X=\sum_{j\ge0}(1+4j)\,\sigma^{4j}|x_j|,\qquad \sigma=1+2^{-20},
\]

and let \(Y\) be the space of pairs of a real number and a cosine series \(y=\sum_jy_j\cos4j\theta\), with norm \(|y_{\rm area}|+\sum_j\sigma^{4j}|y_j|\). The boundary map is

\[
G(x)=\bigl(\langle x^2\rangle,\ \operatorname{Re}W[x]\bigr)\in Y ,
\]

where \(\langle\cdot\rangle\) is the angular mean and \(W[x]\) is the second boundary field of the solid, as a function of the angle. The solid has area \(f\) and the disk law \((f,r)\) exactly when \(G(x)=(4bf,\ r)\) and the solid lies inside the cell: then \(\pi\delta^2\langle x^2\rangle=f\), \(\operatorname{Re}W=r\) on the boundary forces \(W\equiv r\), and the mean of \(g\) vanishes by the fourfold symmetry.

This is the quarter-area map \(\mathcal H\) with two changes of bookkeeping. The scale \(\delta\) is the same for all members of a family, so the constant term \(x_0\) of a profile is not normalised to one. The targets \(4bf\) and \(r\) are not built into the map. \(X\) and \(Y\) are the quarter-area function spaces. Complex coefficients give the complexifications, and \(G\) extends analytically to complex shapes near a real profile, with real-part extraction replaced by \((p+p^\#)/2\).

\textbf{Three things are taken from the quarter-area certificate}, each for a real profile \(q\) with rational coefficients. The node program calls the quarter-area verifier for them.

\begin{itemize}
\item
  \emph{Residual.} A rigorous upper bound for \(\|G(q)-t\|_Y\), for a given target \(t\). The trace part is measured in a stronger norm, as in the point certificate.
\item
  \emph{Defect.} For an operator \(A:Y\to X\) given by a finite matrix on the low modes and a symbol on the high modes, rigorous upper bounds for \(\|A\|\) and for \(\|I-G'(q)A\|_{Y\to Y}\). The matrix and the symbol may be prescribed. The point certificate proposes them itself, but its bounds are checked for whatever is given.
\item
  \emph{Crude bound.} For a radius \(\varrho\) that passes two explicit tests, \(G\) is analytic on the complex ball \(\|x-q\|_X<\varrho\) and bounded there by an explicit number \(M_0(q,\varrho)\). This also gives the quarter-area Hessian estimate: on the ball of radius \(\varrho/5\) about \(q\) the second derivative of \(G\) is at most \((5/2\varrho)^2M_0(q,\varrho)\). The quarter-area proof uses \(\varrho=1/20\). Here it is used for radii up to \(1/2\), so it is restated and proved for a general radius.
\end{itemize}

\textbf{Lemma 0 (crude bound on a complex ball).} Let \(q\) be a real profile, \(\gamma=wq\) its curve with \(w=e^{i\theta}\), and \(D(w,v)=(\gamma(v)-\gamma(w))/(v-w)\) its divided difference. Suppose \(\|D^{-1}\|_\sigma\le w_0\) in the two-variable Wiener norm of weight \(\sigma\). Write \(\|q\|_\sigma=\sum_j\sigma^{4j}|q_j|\) and \(t_0=\|w\,d\gamma/dw\|_\sigma\). For \(\varrho>0\) put \[
\beta_\varrho=w_0\sigma^2\varrho,\qquad p_\varrho=2\delta\sigma(\|q\|_\sigma+\varrho),\qquad
\Phi_\varrho=\sigma^2(\varrho+t_0)\Bigl[\frac{w_0}{1-\beta_\varrho}+\frac{12\,\delta\,p_\varrho^3}{1-p_\varrho^4}\Bigr].
\] If \(\beta_\varrho<1\) and \(p_\varrho<1\), then \(G\) is analytic on the complex ball \(\|x-q\|_X<\varrho\) and \(\|G(x)\|_Y\le M_0(q,\varrho):=(\|q\|_\sigma+\varrho)^2+\Phi_\varrho^2\) there.

\emph{Proof.} This is the quarter-area construction's argument with \(1/20\) replaced by \(\varrho\). Let \(x=q+h\) with \(\|h\|_X\le\varrho\), complex coefficients allowed. The divided difference of \(w\cdot w^{k}\) has \(k+1\) monomials of weight \(\sigma^{k}\) and that of \(w\cdot w^{-k}\) has \(k-1\), so \(\|D_h-D\|_\sigma\le|h_0|+\sum_{k\ge1}k\sigma^k|h_k|\le\|h\|_X\), and a Neumann series gives \(\|D_h^{-1}\|_\sigma\le w_0/(1-w_0\varrho)\le w_0/(1-\beta_\varrho)\). The field operator is \(Tp=\mathcal P(Mp^\#)\), and since \(M\,d\gamma=d\gamma^\#\) it is the boundary integral of \(L_\delta(\gamma(v)-\gamma(w))\,(\gamma^\#)'(v)\,p^\#(v)\), with no division by the tangent. Split \(L_\delta(d)=1/d+\text{regular}\). In the first part \(1/d=1/((v-w)D_h)\) and the remaining disk Cauchy projection sends \(w^iv^j\) to \(w^{i+j}\) or to zero, so it has norm at most one. In the second, \(\|\delta(\gamma_h(v)-\gamma_h(w))\|_\sigma\le p_\varrho<1\) and the lattice coefficients satisfy \(|G_{4m}|<12\), so the regular kernel has norm at most \(12\delta p_\varrho^3/(1-p_\varrho^4)\). The tangent factor has norm at most \(\sigma^2(\varrho+t_0)\). Hence \(\|T\|\le\Phi_\varrho\), and \(W=T(T1)\) gives \(\|W\|_\sigma\le\Phi_\varrho^2\). The complexified trace component is \((W+W^\#)/2\), of norm at most \(\|W\|_\sigma\), and \(|\langle x^2\rangle|\le\|x\|_\sigma^2\le(\|q\|_\sigma+\varrho)^2\). All expansions converge in norm, uniformly on the ball, which gives analyticity. \(\square\)

The node program computes \(M_0(q,\varrho)\) from exact rational bounds for \(w_0\), \(t_0\), \(\|q\|_\sigma\) and \(\delta\), for every radius of a list that passes the two tests. Two remarks used below. \(G\) is one analytic map on the open set of complex shapes where \(D_x\) is invertible and the complexified diameter is below one period; the balls about different nodes lie in this set, so \(G\) is well defined on their union. The spaces are weighted \(\ell^1\) spaces, so the norm of a real operator on the complexification equals its real norm, and the column sums that the programs compute are the operator norms on both.

\paragraph*{The theorem on a parameter interval}

Let \(s\in[-1,1]\) be the normalised parameter of a straight path in the plane of disk laws, with target \(t(s)=(4bf(s),r(s))\), affine in \(s\). Let

\begin{itemize}
\item
  \(p(s)=\sum_{k=0}^{m}c_kT_k(s)\) be a polynomial family of real profiles, with \(T_k\) the Chebyshev polynomials and \(c_k\in X\);
\item
  \(s_1,\dots,s_N\in[-1,1]\) be nodes, \(\theta_j=\arccos s_j\), and \(\Delta\theta\) the largest distance of a point of \([0,\pi]\) from the nearest \(\theta_j\);
\item
  \(q_j\) be the node profiles, with \(\|q_j-p(s_j)\|_X\le\tau_3\) (they are \(p(s_j)\) with coefficients rounded to \(30\) digits);
\item
  \(A(s)=\sum_{k=0}^{D}A_ks^k\) be a polynomial family of operators \(Y\to X\) of the form above.
\end{itemize}

Let \(E_\kappa\) be the closed ellipse with foci \(\pm1\) and semi-axes \((\kappa\pm\kappa^{-1})/2\), and put

\[
\tau_1=\sum_{k\ge1}\|c_k\|_X(\kappa^k-1),\qquad \tau_2=\Delta\theta\sum_{k\ge1}k\|c_k\|_X,\qquad \varrho_t=\tau_1+\tau_2+\tau_3 .
\]

\textbf{Theorem 1.} Suppose that for some \(\kappa>1\), radii \(\varrho_0>\varrho_t\) and \(\varrho_H\), and numbers \(M\), \(M_H\), \(r_*\), \(d_*\):

\begin{enumerate}
\def\labelenumi{\arabic{enumi}.}
\tightlist
\item
  at every node, \(\|G(q_j)-t(s_j)\|_Y\le r_*\) and \(\|I-G'(q_j)A(s_j)\|\le d_*\);
\item
  at every node, \(G\) is analytic on the complex balls of radii \(\varrho_0\) and \(\varrho_H\) about \(q_j\) and bounded there by \(M\) and \(M_H\).
\end{enumerate}

Put \(a_\kappa=(\kappa+\kappa^{-1})/2\), \(L=M/(\varrho_0-\varrho_t)\), \(L_2=4M/(\varrho_0-\varrho_t)^2\), \(K_E=\sum_k\|A_k\|a_\kappa^k\), \(K=\sum_k\|A_k\|\), \(T_E=\sup_{E_\kappa}\|t\|_Y\) and \(e_n=2\kappa^{-n}/(\kappa-1)\). Then for all integers \(n,n'\) with \(n\Delta\theta,n'\Delta\theta<\pi/2\),

\[
\sup_{[-1,1]}\|G(p(s))-t(s)\|_Y\le\varepsilon:=\frac{r_*+L\tau_3+(M+T_E)e_n}{\cos n\Delta\theta}+(M+T_E)e_n ,
\]

\[
\sup_{[-1,1]}\|I-G'(p(s))A(s)\|\le d:=\frac{d_*+K_EL_2\tau_3+(1+K_EL)e_{n'}}{\cos n'\Delta\theta}+(1+K_EL)e_{n'} .
\]

If moreover \(d<1\) and, with \(H=(5/2\varrho_H)^2M_H\), there is \(\beta>0\) such that

\[
c:=d+HK^2\beta<1,\qquad \varepsilon\le(1-c)\beta,\qquad \tau_2+\tau_3+K\beta\le\varrho_H/5 ,
\]

then for every \(s\in[-1,1]\) there is \(y(s)\in Y\) with \(\|y(s)\|_Y\le\varepsilon/(1-c)\) and

\[
G\bigl(p(s)+A(s)y(s)\bigr)=t(s) ,
\]

and \(y\) is real-analytic in \(s\). In particular the exact profile \(x(s)=p(s)+A(s)y(s)\) lies within \(K\varepsilon/(1-c)\) of \(p(s)\) in \(X\), hence uniformly.

\emph{Proof.} \emph{The tube.} Let \(s\in E_\kappa\). Write \(s=\cosh(\alpha+i\beta_0)\) with \(0\le\alpha\le\log\kappa\), and let \(s'=\cos\beta_0\in[-1,1]\). Since \(T_k(s)=\cosh k(\alpha+i\beta_0)\), \[
|T_k(s)-T_k(s')|=|(\cosh k\alpha-1)\cos k\beta_0+i\sinh k\alpha\sin k\beta_0|\le e^{k\alpha}-1\le\kappa^k-1 ,
\] so \(\|p(s)-p(s')\|_X\le\tau_1\). Let \(s_j\) be the node whose angle is nearest to \(\beta_0\). Then \(|T_k(s')-T_k(s_j)|=|\cos k\beta_0-\cos k\theta_j|\le k\Delta\theta\) and \(\|p(s')-p(s_j)\|_X\le\tau_2\). Hence \(\|p(s)-q_j\|_X\le\varrho_t<\varrho_0\). By hypothesis 2, \(G\circ p\) is analytic on a neighbourhood of \(E_\kappa\) with \(\|G(p(s))\|_Y\le M\) there. The ball of radius \(\varrho_0-\varrho_t\) about \(p(s)\) lies in the ball of radius \(\varrho_0\) about \(q_j\), so Cauchy's estimates give \(\|G'(p(s))\|\le L\) and \(\|G''(p(s))\|\le L_2\) on \(E_\kappa\), and the same along the segment from \(q_j\) to \(p(s_j)\).

\emph{Two lemmas.} (i) Let \(F\) be analytic on a neighbourhood of \(E_\kappa\) with values in a Banach space and \(\|F\|\le M_F\) on \(E_\kappa\). Its Chebyshev coefficients are the Laurent coefficients of \(z\mapsto F((z+z^{-1})/2)\) on \(\kappa^{-1}\le|z|\le\kappa\), so \(\|a_k\|\le2M_F\kappa^{-k}\), and the truncation \(F_n=\sum_{k\le n}a_kT_k\) satisfies \(\sup_{[-1,1]}\|F-F_n\|\le M_Fe_n\). If \(F\) is real on \([-1,1]\) so is \(F_n\). (ii) Let \(Q\) be a polynomial of degree at most \(n\) with values in a real Banach space and \(n\Delta\theta<\pi/2\). Then \(\sup_{[-1,1]}\|Q\|\le\max_j\|Q(s_j)\|/\cos n\Delta\theta\). Indeed for a functional \(\ell\) of norm one, \(\tau(\theta)=\ell(Q(\cos\theta))\) is a real trigonometric polynomial of degree at most \(n\). Szeg\H{o}'s inequality \(\tau'^2+n^2\tau^2\le n^2\|\tau\|_\infty^2\) makes \(\arccos(\tau/\|\tau\|_\infty)\) Lipschitz with constant \(n\), so \(|\tau(\theta)|\ge\|\tau\|_\infty\cos n(\theta-\theta_*)\) within \(\pi/2n\) of a point \(\theta_*\) where \(|\tau|\) is largest. Take the node angle nearest to \(\theta_*\), which may be assumed in \([0,\pi]\) because \(\tau\) is even.

\emph{Residual.} Apply (i) and (ii) to \(F(s)=G(p(s))-t(s)\), bounded by \(M+T_E\) on \(E_\kappa\). At a node \(\|F(s_j)\|\le r_*+L\tau_3\), because the certificate was evaluated at \(q_j\). So \(\|F_n(s_j)\|\le r_*+L\tau_3+(M+T_E)e_n\), and \(\sup\|F\|\le\sup\|F_n\|+(M+T_E)e_n\) gives the first bound.

\emph{Defect.} Apply (i) and (ii) to the operator-valued function \(F(s)=I-G'(p(s))A(s)\). On \(E_\kappa\), \(|s|\le a_\kappa\), so \(\|A(s)\|\le K_E\) and \(\|F(s)\|\le1+K_EL\). At a node \(\|F(s_j)\|\le d_*+L_2\tau_3K_E\). This gives the second bound.

\emph{Contraction.} Fix \(s\in[-1,1]\) and let \(j\) be the node nearest in angle. For \(\|y\|_Y\le\beta\) the point \(p(s)+A(s)y\) lies within \(\tau_2+\tau_3+K\beta\le\varrho_H/5\) of \(q_j\), where \(\|G''\|\le H\). The map \(T_s(y)=y-\bigl(G(p(s)+A(s)y)-t(s)\bigr)\) has derivative \(I-G'(p(s)+A(s)y)A(s)\), of norm at most \(d+HK\beta\cdot K=c<1\) on the ball, and \(\|T_s(0)\|\le\varepsilon\). So \(T_s\) maps the ball \(\|y\|\le\beta\) into itself and is a contraction. Its fixed point satisfies \(G(p(s)+A(s)y)=t(s)\) and \(\|y\|\le\varepsilon/(1-c)\). The contraction is uniform in \(s\) and \(T_s\) is analytic in \((s,y)\), so the fixed point is real-analytic in \(s\). \(\square\)

\textbf{From the theorem to solids.} If in addition the node profiles satisfy \(x_{\min}\le q_j(\theta)\le x_{\max}\) and

\[
x_{\min}-\bigl(\tau_2+\tau_3+K\varepsilon/(1-c)\bigr)>0,\qquad 2\delta\bigl(x_{\max}+\tau_2+\tau_3+K\varepsilon/(1-c)\bigr)<1 ,
\]

every exact profile is positive and the solid lies in the open disk of radius \(\tfrac12\), so the translates are disjoint. Then for every \(s\) the solid with radius \(\delta\,x(s)(\theta)\) has area exactly \(f(s)\) and the disk law \((f(s),r(s))\) at every contrast.

\paragraph*{The theorem on a rectangle}

For a rectangle \(f=f_c+h_1s_1\), \(r=r_c+h_2s_2\), \((s_1,s_2)\in[-1,1]^2\), take \(p(s)=\sum c_{k_1k_2}T_{k_1}(s_1)T_{k_2}(s_2)\), a tensor grid of nodes with angular gaps \(\Delta\theta_1,\Delta\theta_2\), and \(A(s)=\sum A_{k_1k_2}s_1^{k_1}s_2^{k_2}\).

\textbf{Theorem 2.} Theorem 1 holds on the square with the following replacements, for ellipse parameters \(\kappa_1,\kappa_2>1\):

\[
\tau_1=\sum_{(k_1,k_2)\ne0}\|c_{k_1k_2}\|_X(\kappa_1^{k_1}\kappa_2^{k_2}-1),\qquad
\tau_2=\sum\|c_{k_1k_2}\|_X(k_1\Delta\theta_1+k_2\Delta\theta_2),
\]

\(K_E=\sum\|A_{k_1k_2}\|a_{\kappa_1}^{k_1}a_{\kappa_2}^{k_2}\), the remainder \(e_n\) replaced by

\[
e_{n_1n_2}=\frac{2\kappa_1^{-n_1}}{\kappa_1-1}\cdot\frac{\kappa_2+1}{\kappa_2-1}+\frac{\kappa_1+1}{\kappa_1-1}\cdot\frac{2\kappa_2^{-n_2}}{\kappa_2-1},
\]

and \(\cos n\Delta\theta\) replaced by \(\cos n_1\Delta\theta_1\cos n_2\Delta\theta_2\).

\emph{Proof.} For \(s\in E_{\kappa_1}\times E_{\kappa_2}\) with feet \(s'\), \(|T_{k_1}(s_1)T_{k_2}(s_2)-T_{k_1}(s_1')T_{k_2}(s_2')|\le(\kappa_1^{k_1}-1)\kappa_2^{k_2}+(\kappa_2^{k_2}-1)=\kappa_1^{k_1}\kappa_2^{k_2}-1\), because \(|T_k|\le\kappa^k\) on \(E_\kappa\) and \(|T_k|\le1\) on \([-1,1]\). A function analytic and bounded by \(M_F\) on the product of the two ellipses has double Chebyshev coefficients of norm at most \(\epsilon_{k_1}\epsilon_{k_2}M_F\kappa_1^{-k_1}\kappa_2^{-k_2}\) with \(\epsilon_0=1\), \(\epsilon_k=2\). Summing those with \(k_1>n_1\) or \(k_2>n_2\) gives \(M_Fe_{n_1n_2}\). A polynomial of degrees \((n_1,n_2)\) is bounded on the square by its largest node value divided by \(\cos n_1\Delta\theta_1\cos n_2\Delta\theta_2\), by lemma (ii) applied in one variable at a time. The rest is as before. \(\square\)
\phantomsection\label{endfile:families}

For the benchmark phase pair there are 30 closed intervals covering
$4/25\le F\le2/5$ on the line $r=65/96-5F/6$, and 17 covering
$1/50\le r\le101/200$ at $F=1/4$. The four original rectangles have
common first interval $[1/5,3/10]$ and successive second intervals
$[1/25,3/50]$, $[3/50,2/25]$, $[2/25,1/10]$, $[1/10,13/100]$.
At these same fractions the upper adjacent strip $13/100<r\le4/25$
remains open away from the separately certified parameter lines. The parameter polynomials, inverse
matrices and all accepted node records accompany the proofs.

\subsection{A complete interval and the provenance of the square-family inputs}\label{detail:node-example}
Consider the piece $11/50\le F\le6/25$ on the line
$r=65/96-5F/6$. The family specifies a polynomial in a parameter
mapped affinely from this interval to $[-1,1]$. At each of 144
nodes the coefficient list describes a radial profile and a polynomial
approximate inverse. For example the first profile is evaluated at phase fraction $201326093/838860800$. In the fixed working norms its
residual is less than $8.971\cdot10^{-14}$ and its derivative
approximation error is less than $1.805\cdot10^{-5}$. The residual
and defect maxima over all nodes are less than $9.705\cdot10^{-14}$
and $0.22108$. These bounds are for the full Fourier operators.

Use the Bernstein ellipse parameter $8/5$ and the analytic ball
radii $1/6$ and $1/12$. The angular mesh gap is less than $0.01091$.
The sum of the three profile errors in the interpolation theorem is
less than $0.16114<1/6$, and the complex bound for the map is less
than $514379$. Apply the residual interpolation bound with degree
102 and the defect bound with degree 61. The outward results on the
whole real interval are
\[
 \epsilon<2.280\cdot10^{-13},\quad d<0.28579,\quad
 K<18.875,\quad H<1.1712\cdot10^7 .
\]
At the rational residual-space radius $\beta=8.5\cdot10^{-11}$
these give $c=d+HK^2\beta<0.642$ and
$\epsilon<(1-c)\beta$. The displacement is bounded by $K\beta$
and fits the differentiability neighborhood. Positivity of the radius
and separation from its translates hold uniformly. The recorded
profile minimum exceeds $0.58$ in the profile units and the physical
radius stays below $0.446<1/2$. The same tests applied at the published
rational endpoints establish each of the remaining intervals and
rectangles. No continuation across a failed piece is used.

The original job and file-server records identify a particular file
version for each family specification and inverse by its digest and
its file-server identity, size, modification time and change time.
For all 51 original pieces the last change preceded submission. The
coordinator's run manifest selected the specified node, and the
original task output reports that exact coordinate, residual, area
residual and derivative error. Those three diagnostic decimals
match the full rational receipt when rounded for display. The check
covers all 22,576 nodes, including the 22,374 older receipts without
both embedded digests. The independent provenance table also binds
each released receipt by SHA-256. One node per piece was freshly
computed in this revision cycle. The remaining node calculations
were checked from the archived executions and their exact inputs.
The first-node comparison alone is never used to extend a result
to an interval. That extension uses all its nodes in the analytic
estimates above.

\subsection{A larger region from the completed square-lattice rectangles}\label{detail:square-extension}
At the 17:57 GMT snapshot on 21 November the next square-lattice
calculation had completed all node records for 18 more rectangles.
Each record includes the digests of its exact profile and polynomial
inverse. Applying the two-parameter interpolation proof in
Section~\ref{route:families-of-square-periodic-inclusions} to each
rectangle gives the following closed unions.
\begin{center}\small
\begin{tabular}{@{}llrr@{}}\toprule
Fraction $F$ & Disk ratio $r$ & Rectangles & Nodes\\\midrule
$[1/10,3/20]$ & $[1/50,1/5]$ & 8 & 39,008\\
$[3/20,1/5]$ & $[1/50,1/5]$ & 8 & 24,672\\
$[1/5,3/10]$ & $[1/50,1/25]$ & 2 & 2,624\\\bottomrule
\end{tabular}\end{center}
In each of the first two rows the subdivision endpoints for $r$ are
\[
 1/50,\ 3/100,\ 1/25,\ 3/50,\ 2/25,\ 1/10,\ 13/100,\ 4/25,\ 1/5 .
\]
In the last row the cut is $3/100$. These closed intervals meet
at their exact rational endpoints. The last row joins the four
previous rectangles along $r=1/25$.

On each grid the verifier computes the exact rounded Fourier profile
from the bivariate Chebyshev polynomial, checks its digest and the
inverse digest, and checks every parameter, area, residual and
derivative input used by the interpolation test. One interior node
per rectangle was also recomputed from the Fourier and inverse
coefficients, and the interval test was repeated with the fresh node.
In residual coordinates let $\epsilon,d,K,H$ and $\beta$ be as in
Theorem~\ref{new:firstbranch}. The exact tests give
\[
     d+HK^2\beta<0.555,\qquad
     \epsilon\le (1-d-HK^2\beta)\beta .
\]
The corrected profiles are positive. Their physical outer radii
are below $0.339$ in the unit square. Thus the maps have a fixed
point on the whole closed parameter rectangle and their inclusions
are separated from periodic translates. The 18 certificates contain
66,304 node records. The area of the union of proved rectangles in
the $(F,r)$ plane has increased from $9/1000$ to $29/1000$.
The original optical maximum assertion concerns the separate line
$r=65/96-5F/6$.

\subsection{A phase-two inclusion on the triangular lattice}\label{route:a-phase-two-inclusion-on-the-triangular-lattice}
Fix $F=1/100$ and
$r_*=(35\sqrt{113}-100)/352+(25+5\sqrt{113})F/88$. The following
argument bounds all Fourier modes of the exact boundary. The
independent transmission solve is a numerical check on the reference
curve. A floating-point proposal enters the proof only through the
enclosure of its residual and of its inverse defect.
For completeness, let $\chi$ be the mask of the inclusion and $\mathcal P$
the interior Cauchy trace for the periodic kernel. On its boundary put
$M=d\bar z/dz$, $g=\mathcal PM$ and $W=\mathcal P(M\bar g)$. For a
constant planar field identified with $e\in\mathbb C$, conjugating the
periodic Beurling kernel and using Stokes' formula gives
$Te=\bar g\bar e$ and $T^2e=\overline W e$ for the anti-linear operator
$T$ of Theorem~\ref{p7:star}. Thus $\langle g\rangle_D=0$ and $W=r$
make $Te$ orthogonal to constants and $T^2e=re$. The contrast
resolvent on the span of $e,Te$ reduces to two scalar equations. Its
constant coefficient is $(1-t)/(1-Ft-rt^2)$ and its conductivity is
$q(1-rt^2+Ft)/(1-rt^2-Ft)$. This derives the whole law from the
boundary identities. Threefold symmetry forces the first mean to
vanish for the trial family below.
\phantomsection\label{file:second-branch}
\paragraph*{Lattice and kernel}

\(\Lambda\) has six nearest neighbours at distance \(A=1.074570\ldots\), in the directions \(0,60,\dots,300\) degrees; the star's three tips point at three of them. The Weierstrass zeta function of \(\Lambda\) has quasi-periods \(\zeta(u+\ell)-\zeta(u)=\eta(\ell)\), linear in \(\ell\). The six-fold symmetry \(\zeta(e^{i\pi/3}u)=e^{-i\pi/3}\zeta(u)\) forces \(\eta(\ell)=\beta\bar\ell\), and Legendre's relation \(\eta(\ell_1)\ell_2-\eta(\ell_2)\ell_1=2\pi i\) gives \(\beta=\pi/\text{(cell area)}=\pi\). So \(k=\zeta\) is the kernel of the boundary identities preceding this calculation, with no linear correction. Near \(0\), \[
\zeta(u)=\frac1u-\sum_{k\ge1}G_{6k}u^{6k-1},\qquad G_n=\sum_{\ell\ne0}\ell^{-n},\qquad |u|<A ,
\] since the Laurent series of \(\zeta\) is \(1/u-\sum_{k\ge2}G_{2k}u^{2k-1}\) and the symmetry kills \(G_{2k}\) unless \(3\mid k\). All \(G_{6k}\) are real because \(\Lambda=\bar\Lambda\).

\emph{Enclosures} (\texttt{sblattice.py}). \(G_6=A^{-6}G_6(\tau)\) with \(\tau=e^{i\pi/3}\) and \[
G_6(\tau)=2\zeta(6)-\frac{2(2\pi)^6}{5!}\sum_{n\ge1}\sigma_5(n)q^n,\qquad q=e^{2\pi i\tau}=-e^{-\pi\sqrt3},
\] summed to \(n=40\). The tail is enclosed by \(\sum_{n>40}\sigma_5(n)|q|^n\le41^6|q|^{41}/(1-2|q|)\), using \(\sigma_5(n)\le n^6\) and \(((n+1)/n)^6|q|\le2|q|\). This gives \(G_6=3.80815079227477\ldots\). The higher \(G_{6k}=c_{3k}/(6k-1)\) follow from the recurrence for the Laurent coefficients \(c_n=(2n-1)G_{2n}\) of \(\wp\) (DLMF 23.9.7), with \(c_2=0\) and \(c_3=5G_6\). For tails, \texttt{bound(k)} uses \[\begin{gathered}
|G_{6k}|\le 6A^{-6k}+S_6(\sqrt3A)^{-6(k-1)},\\ S_6=6(\sqrt3A)^{-6}+6(2A)^{-6}+12(\sqrt7A)^{-6}+0.15\,A^{-6}\ \ge\sum_{|\ell|\ge\sqrt3A}|\ell|^{-6}.
\end{gathered}\] The shells of \(\Lambda\) inside \(|\ell|<3A\) are \(A\) (6 points), \(\sqrt3A\) (6), \(2A\) (6) and \(\sqrt7A\) (12). For \(|\ell|=t\ge3A\) and \(z\) in the disk of radius \(A/2\) about \(\ell\) one has \(|z|-A/2\le t\); these disks are disjoint. Hence \[
\sum_{|\ell|\ge3A}|\ell|^{-6}\le\frac4{\pi A^2}\int_{|z|\ge5A/2}(|z|-\tfrac A2)^{-6}dA=\frac8{A^2}\int_{2A}^\infty u^{-6}(u+\tfrac A2)\,du=0.15\,A^{-6}.
\] These tails only ever enter at the level \(10^{-40}\).

\paragraph*{Curves, norms, unknowns}

\textbf{Curves.} \(\gamma(w)=\sum_ka_kw^k\) on \(|w|=1\), \(w=e^{is}\), with real \(a_k\) and \(a_k=0\) unless \(k\equiv1\pmod3\). Such a curve satisfies \(\gamma(\omega w)=\omega\gamma(w)\) for \(\omega=e^{2\pi i/3}\), and \(\gamma(\bar w)=\overline{\gamma(w)}\). Let \(N=w\,d/dw\), \(G=N\gamma=w\gamma'(w)\) (so \(d\gamma/ds=iG\)), and \(Rp(w)=p(1/w)\). On the circle \(Rp=\bar p\) when \(p\) has real coefficients. A series \emph{has class \(c\)} if its coefficients vanish off \(k\equiv c\pmod 3\); products add classes and \(R\) negates them.

\textbf{Norms.} For a base \(x\ge1\) put \(\|p\|_x=\sum_k|p_k|x^{|k|}\), and for two variables \(\|F\|_x=\sum_{i,j}|F_{ij}|x^{|i|+|j|}\). These are Banach algebra norms, with \[
\|Rp\|_x=\|p\|_x,\quad \|\Pi_+p\|_x\le\|p\|_x,\quad \|\Delta\varphi\|_x\le x\|N\varphi\|_x,\quad \|\mathcal CF\|_x\le\|F\|_x ,
\] where \(\Pi_+\) keeps \(k\ge0\) and \(\Delta\varphi(w,v)=(\varphi(v)-\varphi(w))/(v-w)\); the monomial \(w^k\) contributes \(|k|\) terms, each of weight at most \(x^{|k|+1}\). Finally \(\mathcal C(w^iv^j)=w^{i+j}\) for \(j\ge0\) and \(\mathcal C(w^iv^j)=0\) for \(j<0\), and the norm bound holds because \(x^{|i+j|}\le x^{|i|+|j|}\). Two bases are used: \(\sigma=1+2^{-20}\), the working norm, and \(\tau=1+2^{-11}\) for enclosures of the reference fields. One derivative of a \(\tau\)-small error is \(\sigma\)-small: \[
\|Np\|_\sigma\le C_N\|p\|_\tau,\qquad C_N=\sup_{k\ge0}k(\sigma/\tau)^k=\frac{1}{e\log(\tau/\sigma)}\le\frac{\tau}{e(\tau-\sigma)}\le755 ,
\] using \(\log x\ge1-1/x\). (The script divides by \(2.718<e\).)

\textbf{Spaces.} \(X\) is the real Banach space of \(h=\sum h_kw^k\) with \(h_k=h_{-k}\in\mathbb R\), of class 0, normed by \[
\|h\|_X=\|h\|_\sigma+\|Nh\|_\sigma=\sum_k(1+|k|)\sigma^{|k|}|h_k|.
\] These are the even \(2\pi/3\)-periodic real functions of \(s\). Let \(Y_0\) be the same sequences normed by \(\|\cdot\|_\sigma\), and \(Y=\mathbb R\times Y_0\) with \(\|(y_a,y)\|_Y=\lambda|y_a|+\|y\|_\sigma\), \(\lambda=20\). In the bases \(\{\cos ms\}\) and \(\{e_{\rm area},\cos ks\}\) these are weighted \(\ell^1\) norms, with \[
\|\cos ms\|_X=(1+m)\sigma^m,\qquad \|\cos ks\|_Y=\sigma^k,\qquad \|e_{\rm area}\|_Y=\lambda .
\] So the norm of a linear map between them is the supremum over basis vectors of \(\|\text{image}\|/\|\text{vector}\|\). All operator norms below are computed that way.

\textbf{Reference data} (\texttt{reference.json}, treated as exact).

\begin{itemize}
\item
  The reference curve $\gamma_0$ has the 320 dyadic Laurent coefficients in \path{certificates/p7-second/reference.json}. They have denominator $2^{100}$ and $|k|\le479$. All enclosures below use these exact coefficients.
\item
  \(\psi\): an even cosine polynomial of class 0 and degree 90, with dyadic coefficients and values in \([0.11,8.25]\). It is a smoothed \(|G_0|^{-1.25}\). The field is \(\nu=G_0\psi\), so \(h\nu\) is the normal displacement \(h\psi|G_0|\). The weight \(\psi\) lowered the norm of the approximate inverse by a factor fifteen in the floating-point design (\texttt{proto/proto4.py}), but no property of \(\psi\) is used in the proof.
\end{itemize}

Shapes are \(\gamma_h=\gamma_0+h\nu\) for \(h\in X\); \(\gamma_h\) has class 1 and real coefficients.

\paragraph*{The interior trace in a parametrization}

Let \(\gamma\) be a curve as in Section 2 whose divided difference \(D(w,v)=(\gamma(v)-\gamma(w))/(v-w)\) is invertible in the two-variable algebra of some base \(x>1\). Then \(D\ne0\) on \(x^{-1}\le|w|,|v|\le x\). In particular \(\gamma\) is injective on the circle and \(\gamma'=D(w,w)\ne0\), so \(\gamma\) is an analytic Jordan curve. Suppose in addition its winding number about \(0\) is one (verified in Section 10). Then it is positively oriented, and \(\gamma(w)\) lies inside it for \(|w|<1\) close to \(1\): indeed \(\gamma((1-\varepsilon)e^{is})\approx\gamma+i\varepsilon\,d\gamma/ds\), which is to the left of the tangent.

\textbf{Lemma 3.1 (trace formula).} If moreover \(|\gamma(v)-\gamma(w)|<A\) on the torus, then for every \(q\) with \(\|q\|_1<\infty\) \[
\mathcal Pq=\Pi_+q+\mathcal Kq,\qquad(\mathcal Kq)_i=\sum_jK_{ij}q_{-j},\qquad K(w,v)=v\,\partial_v\log D(w,v)+G(v)R_\Lambda(\gamma(v)-\gamma(w)),\tag{3.1}
\] where \(R_\Lambda(u)=\zeta(u)-1/u=-\sum G_{6k}u^{6k-1}\) and \(K_{ij}\) are the Laurent coefficients of \(K\).

\emph{Proof.} For \(z=\gamma(w)\) with \(|w|<1\) near \(1\), \(k(\gamma(v)-\gamma(w))\gamma'(v)=\gamma'(v)/((v-w)D)+\gamma'(v)R_\Lambda=1/(v-w)+\partial_v\log D+\gamma'(v)R_\Lambda\), because \(\partial_vD\,(v-w)+D=\gamma'(v)\). The Cauchy integral of \(q\) is therefore \(\frac1{2\pi i}\oint\frac{q(v)\,dv}{v-w}+\frac1{2\pi}\int K(w,e^{is'})q(e^{is'})\,ds'\), and the limit \(|w|\to1^-\) is (3.1). \(\square\)

On the torus \(K=G(v)\zeta(\gamma(v)-\gamma(w))-v/(v-w)\) for \(v\ne w\), and \(K(w,w)=(NG-G)/(2G)\). The two symmetries give \(K(\omega w,\omega v)=K(w,v)\), so \(K_{ij}=0\) unless \(i+j\equiv0\), and real \(K_{ij}\). Hence \(\mathcal K\) maps class \(c\) to class \(c\), using only the block \(i\equiv c\), \(j\equiv-c\).

\textbf{The field operator.} \(M=d\bar z/dz=-RG/G\) on the circle, and \(T_\gamma p:=\mathcal P(M\,Rp)\) is the boundary trace of \(\frac1{2\pi i}\oint k(y-z)\bar p\,d\bar y\). By the same computation, and since \(Mdy=d\bar y=(R\gamma)'(v)dv\) on the circle, \[
T_\gamma p=\mathcal C\big[(Rp)(v)\,(R\gamma)'(v)\,D(w,v)^{-1}\big]+\big[\,v\,R_\Lambda(\gamma(v)-\gamma(w))\,(Rp)(v)\,(R\gamma)'(v)\big]_{v^0},\tag{3.2}
\] where \([\cdot]_{v^0}\) is the coefficient of \(v^0\). This expression contains only \(D^{-1}\) and one derivative of \(\gamma\), and is analytic in the coefficients of \(\gamma\). Then \(g=T_\gamma1\) (class 1) and \(W=T_\gamma g\) (class 0).

\paragraph*{The boundary equation and the reduction}

For \(h\) in a ball of \(X\) on which \(D_h=D_0+\Delta(h\nu)\) is invertible (Section 8), define \[
\mathcal H(h)=\Big(\pi\langle R\gamma_h\cdot G_h\rangle-F,\ \tfrac12(W_h+RW_h)-r_*\Big)\in Y,\qquad W_h=T_{\gamma_h}(T_{\gamma_h}1),
\] where \(\langle\cdot\rangle\) is the constant coefficient. The first entry is \(|D_h|-F\), by Green's formula \(|D|=\frac12\int\operatorname{Im}(\bar\gamma\,d\gamma/ds)\,ds=\pi\sum_kka_k^2\). Real coefficients and the class structure put \(\frac12(W_h+RW_h)\) in \(Y_0\). At \(h=0\), (3.2) and (3.1) agree, and (3.1) is what the script evaluates.

\textbf{Lemma 4.1.} Suppose \(h_*\) is real, \(\mathcal H(h_*)=0\), and \(\gamma_*=\gamma_{h_*}\) is a positively oriented Jordan curve with \(D_{h_*}\) invertible for a base \(>1\). Then its interior \(D_*\) has \(|D_*|=F\), \(\langle g\rangle_{D_*}=0\) and \(W\equiv r_*\).

\emph{Proof.} By Section 3, \(W_{h_*}\) is the boundary trace of the holomorphic function \(W\) in \(D_*\), which is continuous up to \(\partial D_*\). Since \(\operatorname{Re}W=r_*\) on \(\partial D_*\), \(\operatorname{Re}W\equiv r_*\) and \(W=r_*+ic\) with \(c\) real. On the circle \(\operatorname{Im}W=(W_{h_*}-RW_{h_*})/(2i)\), whose constant coefficient vanishes because \(W_{h_*}\) has real coefficients; so \(c=0\). For the mean: \(D_*\) and \(\Lambda\) are invariant under \(z\mapsto\omega z\), \(\zeta(\omega u)=\omega^{-1}\zeta(u)\) and \(M(\omega y)=\omega^{-2}M(y)\). Substituting \(y=\omega y'\) in the Cauchy integral gives \(g(\omega z)=\omega^{-2}g(z)\), so \(\int_{D_*}g=\omega^{-2}\int_{D_*}g\), which forces \(\int_{D_*}g=0\). \(\square\)

\paragraph*{Enclosures at the reference curve}

\paragraph*{The inverse of the divided difference}

\(D_0\) is an explicit two-variable Laurent polynomial: \(D_{ij}=a_{i+j+1}\) for \(i,j\ge0\), \(D_{ij}=-a_{i+j+1}\) for \(i,j\le-1\), and \(0\) otherwise. A trial inverse \(U\) is proposed by a floating-point FFT of \(1/D_0\) on a \(3072^2\) grid. It is supported on \(|i|+|j|\le600\) with \(i+j\equiv0\), and is dyadic with denominator \(2^{56}\). The product \(D_0U\) is computed exactly in integers by Kronecker substitution (\texttt{fmpz\_poly}). The norm is summed exactly in integers, with weights \(\lceil2^{60}\rho^m\rceil+1\) and \(\rho=1+1/64\). This gives \[
\delta_D=\|1-D_0U\|_\rho .
\] If \(\delta_D<1\) then \(D_0U=1-E\) with \(\|E\|<1\) in a commutative Banach algebra, so \(D_0\) is invertible and \[
\|D_0^{-1}\|_\rho\le d_\rho:=\frac{\|U\|_\rho}{1-\delta_D},\qquad \|D_0^{-1}\|_\sigma\le d_\sigma:=\frac{\|U\|_\sigma}{1-\delta_D},
\] the second because \(\|E\|_\sigma\le\|E\|_\rho\).

\paragraph*{The kernel}

\(K=(v\partial_vD_0)D_0^{-1}+G_0(v)R_\Lambda(\gamma_0(v)-\gamma_0(w))\) lies in the algebra of base \(\rho\). Since \(\|\gamma_0(v)-\gamma_0(w)\|_\rho\le2\|\gamma_0\|_\rho<A\), \[
\|K\|_\rho\le B_K:=\|v\partial_vD_0\|_\rho\,d_\rho+\|G_0\|_\rho\sum_k|G_{6k}|(2\|\gamma_0\|_\rho)^{6k-1}.
\] \texttt{sbkernel.kernel\_enclosure} encloses \(K(w_a,v_b)\), \(w_a=e^{2\pi ia/n}\), on the grid \(n=n_K=7680\) in 160-bit balls. The values of \(\gamma_0\), \(G_0\) and \(NG_0\) at the nodes come from an inverse DFT of the exact coefficients. The zeta series is taken to 16 terms, and the tail \(\big(\sum_{k>16}|G_{6k}|u_{\max}^{6k-1}\big)\sum|ka_k|\), with \(u_{\max}=2\sum|a_k|=0.4148\), is added to each value as a ball. Only the rows \(0\le a\le n/6\) are evaluated. With \(R_a[j]=\frac1n\sum_bK(w_a,v_b)v_b^{-j}\), the symmetries give \[\begin{gathered}
R_{a+n/3}[j]=\omega^{-j}R_a[j]\ \ (\text{from }K(\omega w,\omega v)=K),\\ R_{n/3-a}[j]=\omega^{2j}\overline{R_a[j]}\ \ (\text{also using }K(\bar w,\bar v)=\overline{K(w,v)}),
\end{gathered}\] and therefore, for \(i+j\equiv0\), \[
\tilde K_{ij}:=\frac1{n^2}\sum_{a,b}K(w_a,v_b)w_a^{-i}v_b^{-j}=\frac3n\sum_{0\le a<n/3}R_a[j]\,w_a^{-i},
\] which is a DFT of length \(n/3\) for each \(j\).

\emph{Aliasing.} For centred indices \(|i|,|j|<n/2\), \(\tilde K_{ij}=\sum K_{i'j'}\) over \((i',j')\equiv(i,j)\pmod n\). Every other term in that sum has \(|i'|+|j'|\ge n/2\) and \(|i'|\ge|i|\). Hence, in the norm used for the operator, \[
\sum_{\text{all }(i,j)}|K_{ij}-\tilde K_{ij}|\,\tau^{|i|}\le2\sum_{|i'|+|j'|\ge n/2}|K_{i'j'}|\tau^{|i'|}\le2B_K(\tau/\rho)^{n/2},
\] where \(\tilde K_{ij}:=0\) outside the centred square.

\emph{Truncation.} The stored kernel \(\hat K\) keeps (the real parts of) the balls \(\tilde K_{ij}\) with \(|i|,|j|\le L_K=900\), in the two blocks needed: \((i\equiv0,j\equiv0)\) and \((i\equiv1,j\equiv2)\). The discarded mass is added: for \(|j|\le L_K\) the computed \(|\tilde K_{ij}|\tau^{|i|}\) with \(|i|>L_K\); for \(|j|>L_K\) the bound \(\sum_i|\tilde K_{ij}|\tau^{|i|}\le\tau^{n/2}\sum_{a<n/3}|R_a[j]|\). This gives \[
\epsilon_K\ \ge\sum_{(i,j):\,j\equiv0,2}|K_{ij}-\hat K_{ij}|\,\tau^{|i|}.
\]

\emph{Use.} For \(q\) of class 0 or 1 (the only ones used), \(\|(\mathcal K-\hat{\mathcal K})q\|_\tau\le\epsilon_K\max_j|q_j|\le\epsilon_K\|q\|_\sigma\). The operator norm of \(\hat{\mathcal K}\) on \(\ell^1_x\) is \(\max_jx^{-|j|}\sum_i|\hat K_{ij}|x^{|i|}\). In \texttt{Kapply} the ball polynomial part of \(q\) goes through \(\hat{\mathcal K}\) in ball arithmetic and pays \(\epsilon_K\|q_{\rm ball}\|_\sigma\); the error terms of \(q\) pay \((\|\hat{\mathcal K}\|_x+\epsilon_K)\) times their size.

\paragraph*{The inverse of the tangent}

The proposal \(\hat Q\approx1/G_0\) is the DFT of \(1/G_0\) at 12288 nodes in ball arithmetic. Its midpoints are rounded to dyadics with denominator \(2^{130}\), and coefficients below \(2^{-115}\) are discarded. Then \(E_Q=1-G_0\hat Q\) is computed exactly, with \(\|E_Q\|_\tau\le\epsilon_Q<1\). Hence \(Q=1/G_0=\hat Q(1-E_Q)^{-1}\) exists in \(\ell^1_\tau\) and \(\|Q-\hat Q\|_\tau\le\|\hat Q\|_\tau\epsilon_Q/(1-\epsilon_Q)\).

\paragraph*{The residual}

From here on every quantity is a ball polynomial plus two error terms, \(\|E_t\|_\tau\le e_t\) and \(\|E_s\|_\sigma\le e_s\) (\texttt{sbrig.LS}). The rules are:

\begin{itemize}
\item
  sums add errors;
\item
  products use \(\|fg\|\le\|f\|\|g\|\) for the cross terms (\(\tau\) with \(\tau\) stays \(\tau\), all other combinations go to \(\sigma\));
\item
  \(N\) maps \(e_t\) to a \(\sigma\)-error \(C_Ne_t\); it is not applied to \(\sigma\)-errors;
\item
  \(R\) and \(\Pi_+\) keep errors;
\item
  truncation adds the discarded mass to the error.
\end{itemize}

Then \(M=-RG_0\cdot Q\), \(g=\Pi_+M+\mathcal KM\), \(W=\Pi_+(M\,Rg)+\mathcal K(M\,Rg)\), and \[
Y_0:=\|\mathcal H(0)\|_Y=\lambda\Big|\pi\sum_kka_k^2-F\Big|+\Big\|\tfrac12(W+RW)-r_*\Big\|_\sigma ,
\] with \(r_*\) an Arb ball of radius \(10^{-26}\) computed from \(\sqrt{113}\). \textbf{The theorem is therefore proved at the exact irrational value \(r_*(1/100)\)}: nothing is replaced by a nearby rational.

\paragraph*{The derivative}

\paragraph*{The formula}

\textbf{Proposition 6.1.} For \(h\in X\), with \(\tilde h=\psi h\) and all fields taken at \(\gamma_0\), \[
D\mathcal H(0)h=\Big(2\pi\langle\tilde h\,G_0RG_0\rangle,\ \operatorname{Re}\big\{\mathcal P[S_g\tilde h]+\mathcal P[M\cdot R(\mathcal PS_1\tilde h)]+\tilde h\,NW\big\}\Big),\qquad S_p\tilde h=2Q\,N(RG_0\cdot Rp\cdot\tilde h),
\] where \(\operatorname{Re}f=\frac12(f+Rf)\).

\emph{Proof.} (i) \emph{A domain integral.} For \(p\) holomorphic near \(\bar D\) and \(z\in D\), Stokes' formula (\(\oint\phi\,d\bar y=-2i\int_D\partial_y\phi\,dA\), the small circle about \(y=z\) contributing nothing, as in Lemma 1) and \(\partial_yk(y-z)=-\wp(y-z)\) give \[
T_Dp(z):=\frac1{2\pi i}\oint_{\partial D}k(y-z)\overline{p(y)}\,d\bar y=\frac1\pi\,\mathrm{p.v.}\!\int_D\wp(y-z)\overline{p(y)}\,dA(y).
\] (ii) \emph{Eulerian derivative.} Move the boundary by \(\delta\gamma=\tilde hG_0\), a normal displacement with \(V_n|d\gamma/ds|=\operatorname{Re}(\delta\gamma\,\overline{G_0})=\tilde h|G_0|^2\). For fixed interior \(z\), Reynolds' formula gives \(\delta T_Dp(z)=\frac1\pi\oint\wp(y-z)\bar p\,V_n\,ds_{\rm arc}=\frac1\pi\oint\wp(y-z)\phi\,dy\), with \(\phi=\bar p\,\tilde h\,\overline{d\gamma/ds}\), because \(|d\gamma/ds|^2ds=\overline{d\gamma/ds}\,dy\). Integration by parts (\(\wp(y-z)\,dy=-d\,k(y-z)\)) gives \(\frac1\pi\oint k(y-z)\,d\phi=2i\,\mathcal Pd\phi/dy\). Here \(d\phi/dy=N\phi/G_0\) and \(\phi=-i\,\tilde h\,RG_0\,Rp\), so \(\delta T_Dp=\mathcal P[S_p\tilde h]\). (iii) \emph{The fields.} \(g=T_D1\), so \(\delta g=\mathcal PS_1\tilde h\). Also \(W=T_Dg\), so its Eulerian derivative is \(\mathcal P[S_g\tilde h]+T_D(\delta g)\), and \(T_Dq=\mathcal P(M\,Rq)\) for real-coefficient \(q\). (iv) \emph{Moving point.} \(W_h(e^{is})=W(\gamma_h(e^{is}))\), so the boundary values gain \(W'(\gamma_0)\,\delta\gamma=(NW/G_0)\,\tilde hG_0=\tilde h\,NW\). (v) \emph{Area.} The area derivative is \(\oint V_n\,ds_{\rm arc}=2\pi\langle\tilde h|G_0|^2\rangle\).

The fields at \(\gamma_0\) continue analytically to a collar of the boundary, so for trigonometric polynomials \(h\) this computes \(\frac d{d\varepsilon}\mathcal H(\varepsilon h)|_0\). \(\mathcal H\) is analytic on a ball of \(X\) (Section 8), so \(D\mathcal H(0)\) is a bounded operator; the right side is bounded too (Section 6.2). The two agree on the dense set of trigonometric polynomials, hence on \(X\). \(\square\)

The formula was also checked against central differences of the residual in floating point: errors \(10^{-6}\) against values of order \(10\)-\(100\), at step \(10^{-6}\) (\texttt{proto/proto2.py}).

\paragraph*{\texorpdfstring{The finite derivative \(D_0\) and the error \(\|D\mathcal H(0)-D_0\|\)}{The finite derivative D\_0 and the error \textbackslash\textbar D\textbackslash mathcal H(0)-D\_0\textbackslash\textbar{}}}

\(D_0\) is the formula with every multiplier replaced by an exact dyadic truncation (midpoints rounded at \(2^{-64}\)):

\begin{itemize}
\item
  \(\hat Q_2\approx2Q\), band 360;
\item
  \(\hat\Lambda_g\approx RG_0\cdot Rg\), \(\hat\Lambda_1\approx RG_0\), \(\hat M\approx M\), \(\widehat{NW}\approx NW\) and \(\hat G_2\approx G_0RG_0\), band 270;
\item
  \(\mathcal K\) replaced by the rounded block \(\hat{\mathcal K}'\) of band \(L'_K=330\), with \(\epsilon'_K=\epsilon_K+{}\)(discarded and rounding mass).
\end{itemize}

With \(\hat S_p\tilde h=\hat Q_2N(\hat\Lambda_p\tilde h)\) and \(\hat{\mathcal P}=\Pi_++\hat{\mathcal K}'\), \[
D_0h=\Big(2\pi\langle\tilde h\hat G_2\rangle,\ \operatorname{Re}\big\{\hat{\mathcal P}[\hat S_g\tilde h]+\hat{\mathcal P}[\hat MR(\hat{\mathcal P}\hat S_1\tilde h)]+\tilde h\widehat{NW}\big\}\Big).
\] Since \(N(mh)=mNh+hNm\), \(\|qN(mh)\|_\sigma\le[q;m]\,\|h\|_X\) with \([q;m]:=\max(\|qm\|_\sigma,\|q\,Nm\|_\sigma)\). Put

\begin{itemize}
\item
  \(s_p=[2Q;\Lambda_p\psi]\ge\|S_p\psi\cdot\|\), with \(\Lambda_g=RG_0Rg\) and \(\Lambda_1=RG_0\);
\item
  \(\delta S_p=\|2Q-\hat Q_2\|_\sigma\max(\|\Lambda_p\psi\|,\|N(\Lambda_p\psi)\|)+[\hat Q_2;(\Lambda_p-\hat\Lambda_p)\psi]\ge\|(S_p-\hat S_p)\psi\cdot\|\);
\item
  \(p_{\rm full}=1+\|\hat{\mathcal K}\|+\epsilon_K\ge\|\mathcal P\|\) and \(\hat p=1+\|\hat{\mathcal K}'\|\).
\end{itemize}

Telescoping each term (\(\mathcal P[MR\,x]-\hat{\mathcal P}[\hat MR\,\hat x]=(\mathcal P-\hat{\mathcal P})[MRx]+\hat{\mathcal P}[(M-\hat M)Rx]+\hat{\mathcal P}[\hat MR(x-\hat x)]\), and so on), \[
\begin{aligned}
\|D\mathcal H(0)-D_0\|_{X\to Y}\le{}&\lambda\,2\pi\|\psi(G_0RG_0-\hat G_2)\|_\sigma+\epsilon'_Ks_g+\hat p\,\delta S_g\\
&+\epsilon'_K\|M\|p_{\rm full}s_1+\hat p\|M-\hat M\|p_{\rm full}s_1\\ &+\hat p\|\hat M\|(\epsilon'_Ks_1+\hat p\,\delta S_1)+\|\psi(NW-\widehat{NW})\|_\sigma .
\end{aligned}
\] The \(\tau\)-enclosures of \(g\) and \(W\) are what make \(N(\Lambda_g-\hat\Lambda_g)\) and \(NW\) finite in the \(\sigma\)-norm.

\paragraph*{\texorpdfstring{The approximate inverse \(A:Y\to X\)}{The approximate inverse A:Y\textbackslash to X}}

\paragraph*{Matrix part}

For the inputs \(\cos ms\), \(m=0,3,\dots,M_{\rm in}=\text{HIGH}+270\), the columns \(D_0\cos ms\) are computed with all their output harmonics. The data are dyadic, so this is exact up to the tiny radii of ball arithmetic. A floating-point solve of the leading square system proposes columns for the unit outputs, i.e.~the area and \(\cos ks\) for \(k\le\text{HIGH}\). They are rounded to dyadics (denominator \(2^{52}\)) and multiplied by the full rectangular matrix in ball arithmetic. Then the defect \(\|D_0Ae-e\|_Y/\|e\|_Y\) and the norm \(\|Ae\|_X/\|e\|_Y\) are bounded column by column.

\paragraph*{Parametrix for high frequencies}

Let \(\Gamma\) be a dyadic polynomial of class 2 and band \(b_\Gamma=210\) approximating \(2/(\psi\hat\Lambda_g)\), and let \(H\) be the truncation of \(G_0/2\) at \(b_H=150\) (class 1). For \(m>\text{HIGH}\) put \[
A\cos ms=\operatorname{Re}[\Gamma\Phi_m],\qquad \Phi_m=N^{-1}(Hw^m)
\] (class 0). \(Hw^m\) has frequencies in \([m-b_H,m+b_H]\subset(0,\infty)\), so \(N^{-1}\) is division by the frequency.

\textbf{Lemma 7.1.} Let \(X_b=b_H+b_\Gamma+b_\psi+b_\Lambda+b_Q\), the sum of the bands of \(H,\Gamma,\psi,\hat\Lambda_{g,1},\hat Q_2\), and let \(L\) be the band of \(\hat M,\widehat{NW},\hat G_2\). If \(m>X_b+L+L'_K\), then \[\begin{gathered}
D_0A\cos ms=\Big(0,\ \cos ms+\operatorname{Re}\big[\epsilon_2w^m+\hat Q_2N(\epsilon_1\Phi_m)+\psi\,(A\cos ms)\,\widehat{NW}\big]\Big),\\
\epsilon_1=\tfrac12\hat\Lambda_g\psi\Gamma-1,\quad \epsilon_2=\hat Q_2H-1 .
\end{gathered}\] \emph{Proof.} \(\hat\Lambda_g\psi A\cos ms=(1+\epsilon_1)\Phi_m+\frac12\hat\Lambda_g\psi R(\Gamma\Phi_m)\).

\begin{itemize}
\item
  \emph{First part.} After \(N\) and \(\hat Q_2\) it becomes \(\hat Q_2(Hw^m+N(\epsilon_1\Phi_m))=(1+\epsilon_2)w^m+\hat Q_2N(\epsilon_1\Phi_m)\), all of whose frequencies are at least \(m-X_b>0\).
\item
  \emph{Second part.} After \(N\) and \(\hat Q_2\) its frequencies are at most \(-(m-X_b)<0\). So \(\Pi_+\) returns exactly the first part.
\item
  \emph{Kernel.} Every input frequency of \(\hat{\mathcal K}'\) has modulus at least \(m-X_b>L'_K\), so \(\hat{\mathcal K}'\) gives \(0\).
\item
  \emph{Nested term.} \(\hat{\mathcal P}\hat S_1\) returns frequencies in \([m-X_b,m+X_b]\). \(R\) maps them to \([-m-X_b,-m+X_b]\), and multiplying by \(\hat M\) to \([-m-X_b-L,-m+X_b+L]\). There both \(\Pi_+\) and \(\hat{\mathcal K}'\) vanish.
\item
  \emph{Area.} \(\psi A\cos ms\) has no frequency of modulus \(\le L\). \(\square\)
\end{itemize}

With \(\|\Phi_m\|_\sigma\le\sigma^m\|H\|/(m-b_H)\), \(\|N\Phi_m\|_\sigma\le\sigma^m\|H\|\) and \(\|\cos ms\|_Y=\sigma^m\), this gives, uniformly in \(m\ge\text{HIGH}+3\): \[
\frac{\|D_0A\cos ms-\cos ms\|_Y}{\|\cos ms\|_Y}\le\|\epsilon_2\|+\|\hat Q_2\|\|H\|\Big(\|\epsilon_1\|+\frac{\|N\epsilon_1\|}{m-b_H}\Big)+\frac{\|\psi\widehat{NW}\|\|\Gamma\|\|H\|}{m-b_H},
\] \[
\frac{\|A\cos ms\|_X}{\|\cos ms\|_Y}\le\|H\|\Big(\|\Gamma\|+\frac{\|\Gamma\|+\|N\Gamma\|}{m-b_H}\Big).
\] (\(\operatorname{Re}\) does not increase norms, and \(N\operatorname{Re}f=\operatorname{Re}'(Nf)\) with \(\operatorname{Re}'f=\frac12(f-Rf)\).)

\paragraph*{Defect}

\[
\|I-D\mathcal H(0)A\|_{Y\to Y}\le Z_0:=\max(\text{defect of 7.1},\text{defect of 7.2})+\|A\|\,\|D\mathcal H(0)-D_0\|,
\] where \(\|A\|\) is the larger of the two column bounds.

\paragraph*{The second derivative}

For complex \(h\in X_{\mathbb C}\) use (3.2) with \(R\) in place of conjugation everywhere. Then \(\mathcal H\) is a composition of analytic maps wherever \(D_h=D_0+\Delta(h\nu)\) is invertible and \(2\|\gamma_h\|_\sigma<A\): affine maps, inversion in a Banach algebra, products, \(\mathcal C\), \(R\), and the lattice series. For real \(h\) it is the map of Section 4.

Let \(c_\nu=\max(\|\nu\|_\sigma,\|N\nu\|_\sigma)\). Then \(\|N(h\nu)\|_\sigma\le c_\nu\|h\|_X\) and \(\|D_h-D_0\|_\sigma\le\sigma c_\nu\|h\|_X\). On the complex ball \(\|h\|_X\le R_c:=\frac12/(d_\sigma\sigma c_\nu)\): \[
\|D_h^{-1}\|_\sigma\le2d_\sigma,\quad \|G_h\|_\sigma\le\|G_0\|_\sigma+c_\nu R_c,\quad \|\gamma_h\|_\sigma\le\|\gamma_0\|_\sigma+\|\nu\|_\sigma R_c,
\] \[
\|T_{\gamma_h}p\|_\sigma\le\tau_h\|p\|_\sigma,\qquad \tau_h=\sigma\|G_h\|\Big(2d_\sigma+\sum|G_{6k}|(2\|\gamma_h\|)^{6k-1}\Big),
\] by (3.2), \(\|(R\gamma_h)'\|_\sigma\le\sigma\|G_h\|_\sigma\) and \(\|\mathcal CF\|\le\|F\|\). Hence \(\|\mathcal H(h)\|_Y\le\Phi:=\lambda(\pi\|\gamma_h\|\|G_h\|+F)+\tau_h^2+r_*\). Fix \(h\) with \(\|h\|_X\le r_N:=R_c/10\), unit vectors \(u,v\), and \(\rho_C=(R_c-r_N)/2\). Cauchy's formula on the disks \(|\zeta_1|,|\zeta_2|=\rho_C\) for \((\zeta_1,\zeta_2)\mapsto\mathcal H(h+\zeta_1u+\zeta_2v)\) gives \[
\|D^2\mathcal H(h)\|\le Z_2:=\Phi/\rho_C^2 .
\]

\paragraph*{The contraction}

Let \(\Phi(y)=y-\mathcal H(Ay)\) on real \(Y\). For \(\|y\|_Y\le r_Y\) with \(\|A\|r_Y\le r_N\): \[
\|D\Phi(y)\|\le\|I-D\mathcal H(0)A\|+\|D\mathcal H(Ay)-D\mathcal H(0)\|\|A\|\le Z_0+Z_2\|A\|^2r_Y,\qquad \|\Phi(0)\|=Y_0 .
\] If \(Y_0+Z_0r_Y+Z_2\|A\|^2r_Y^2<r_Y\), then \(\Phi\) maps the ball into itself and contracts. Its fixed point \(y_*\) satisfies \(\mathcal H(h_*)=0\) for \(h_*=Ay_*\in X\), with \(\|h_*\|_X\le\|A\|r_Y\). The script takes \(r_Y=2Y_0/(1-Z_0)\). The target \(r_*\) enters only as the constant mode of the residual. So the same argument, on the larger ball \(r_Y=\min((1-Z_0)/(2Z_2\|A\|^2),r_N/\|A\|)\), covers every target \(r\) with \(|r-r_*|\le(1-Z_0)r_Y-Z_2\|A\|^2r_Y^2-Y_0\); the script reports this slack.

\paragraph*{Geometry}

At \(n=3\cdot2^{17}\) nodes, \(\gamma_0\) and \(G_0\) are enclosed by a DFT of the exact coefficients.

\begin{itemize}
\item
  Between nodes \(|\gamma_0|\) changes by at most \(\frac{2\pi}n\|G_0\|_\sigma\), which gives \(\min|\gamma_0|\ge\gamma_{\min}\).
\item
  \(|\frac{d}{ds}(G_0/\gamma_0)|=|NG_0/\gamma_0-G_0^2/\gamma_0^2|\le(\|NG_0\|\|\gamma_0\|+\|G_0\|^2)/\gamma_{\min}^2\). With this Lipschitz bound, the node minimum gives \(m_0>0\) with \(\operatorname{Re}(G_0/\gamma_0)=\frac d{ds}\arg\gamma_0\ge m_0\) on the whole circle.
\item
  The node average of \(\operatorname{Re}(G_0/\gamma_0)\), with the same bound, lies in \((\frac12,\frac32)\), so the winding number is \(1\).
\end{itemize}

For \(\gamma_*=\gamma_0+h_*\nu\) one has \(|h_*\nu|\le\|h_*\|_X\|\nu\|_\sigma\) and \(|N(h_*\nu)|\le\|h_*\|_Xc_\nu\). Hence \(|G_*/\gamma_*-G_0/\gamma_0|\) is at most the printed perturbation, and \(\operatorname{Re}(G_*/\gamma_*)>0\). Consequently \(\arg\gamma_*\) increases strictly, by \(2\pi\), and:

\begin{itemize}
\item
  \(\gamma_*\) is a positively oriented star-shaped Jordan curve and \(D_*\) is simply connected;
\item
  \(\gamma_*\) is analytic, because \(h_*\in X\) converges on the annulus \(\sigma^{-1}<|w|<\sigma\);
\item
  \(D_{h_*}\) is invertible for the base \(\sigma\), since \(\|h_*\|_X<R_c\);
\item
  \(|\gamma_*|\le\|\gamma_0\|_\sigma+|h_*\nu|<0.2075\), and disks of radius \(0.2075\) about the points of \(\Lambda\) are disjoint;
\item
  \(\gamma_*\) has class 1 and real coefficients, so it has the rotation and reflection symmetries.
\end{itemize}

With Lemma 4.1 this proves the theorem.
\phantomsection\label{endfile:second-branch}

\paragraph{Outward bounds for the contraction.}
With the norms defined above, safe upper bounds are
\[
\begin{array}{c|rrrrrr}
 &\|K\|_\rho&Y_0&\|A\|&Z_0&Z_2&\|D{\cal H}(0)-D_0\|\\\hline
 &204.341&3.391\cdot10^{-18}&48.50&0.0082&1.04\cdot10^{10}
 &1.683\cdot10^{-4}.
\end{array}
\]
At radius $6\cdot10^{-18}$ they give invariant-ball slack exceeding
$2.55\cdot10^{-18}$ and Lipschitz constant below $0.0085$.
The displacement at this radius is at most $2.91\cdot10^{-16}$.
This bound and the kernel bound are rounded outward.

\paragraph{Which disk branch is optimal.}
For a solid of permittivity $p$ in a matrix $q$ write
$Z(F,r)=q(1-rt^2+Ft)/(1-rt^2-Ft)$, $t=(p-q)/(p+q)$.
For phase-one fraction $f$ compare $\operatorname{Im}Z(f,r)$ on
$0\le r\le1-f$ with the interchanged arc on $0\le r\le f$.
The sharp planar disk bounds give these two arcs for the isotropic
region \cite{milton81,milton}. They also bound the minimum directional
loss for general tensors. Indeed an isotropized planar aggregate of
rotations of $Z=X+iY$ has scalar response $z_0$ with
$z_0^2=\det Z$. For the two real eigenvalues $t_j$ of
$Y^{-1/2}XY^{-1/2}$,
\[
 (\operatorname{Im}z_0)^2=\det Y\,
 \frac{\sqrt{(1+t_1^2)(1+t_2^2)}-t_1t_2+1}{2}\ge\det Y.
\]
For $p=-15+i,\ q=1+i$ the first-branch maximum in the asserted
range is $L_1=1+72f/(7-4f)$. If $D_1$ denotes the positive denominator
of the rational imaginary part in the normalization used by the
accompanying verifier, then
\[
 (7-4f)D_1\bigl(L_1-\operatorname{Im}Z(f,r)\bigr)
      =8f(80f+96r-65)^2.
\]
On the interchanged arc the corresponding numerator divided by eight
has positive Bernstein coefficients after setting
$f=4/25+(6/25)x,\ r=fy$ on the unit square. The smallest is
$246114/125$. At $f=99/100$ the second branch has maximum
$L_2=1+8(1-f)(52+5\sqrt{113})/(11-4f)=2.19489464845\ldots$.
The numerator of its deficit on that branch is a positive multiple
of $(r-r_*)^2$; on the other branch all three coefficients of the
quadratic in $r$ are positive. The exact identities use rational
arithmetic and $\mathbb Q(\sqrt{113})$.

\subsection{The short-period obstruction for finite spatial shells}\label{route:the-short-period-obstruction-for-finite-spatial-shells}
\phantomsection\label{file:flat-lattice}
\paragraph*{The contradiction to be proved}

Write \(b=832/(778+285\sqrt{17})\), \(v=b-1/4\), \(s_t=(1-t)/3\), \(\rho=v/b\) and \(\kappa=(1-\rho)/3\). For a periodic set \(A\) of fraction \(t\), put \(\Delta u_A=\mathbf1_A-t\). At equality the lifted outer and cavity potentials for the filled body \(O\) and cavity union \(K\) satisfy

\[
\Theta=s_b|X|^2/2-u_O,\qquad \Xi=s_v|X|^2/2-u_K,
\qquad 0\le D^2\Theta,D^2\Xi\le I.
\]

Here \(X=(x,y)\). Normalize \(\Theta=0\) in the central filled body. In the matrix \(\Xi=\kappa|X|^2/2+\rho\Theta\) after an additive normalization. Thus \(W=(\Xi-\kappa|X|^2/2)/\rho\) satisfies \(\Delta W=1\) off the cavities, agrees with \(\Theta\) in the matrix, and has zero value and gradient on the entire central outer surface. At a cavity boundary the limit of \(D^2\Xi\) from outside is \(n\otimes n\), with \(n\) a unit normal. In particular a positive definite Hessian \(\kappa I+\rho D^2W\) at a point of a continued branch excludes a first cavity contact there. These facts follow from attaining the loss at the single specified contrast.

Let \(e\) be the short-period unit vector. We prove that the canonical continuation of \(\Theta\) from the right outer face, along the central short-period axis, satisfies

\[
\|D^2W-e\otimes e\|\le C_b(\lambda_1/\lambda_2)^2
\quad\text{all the way to the opposite outer face for sufficiently large }\lambda_2/\lambda_1 .
\tag{1}
\]

The positive lower bound \(\kappa\) for \(\kappa I+\rho\,e\otimes e\) is independent of the cavity count. If a proposed finite shell met that segment in a closed cavity, the first contact would contradict (1). If it missed every cavity, \(W_x\) would increase strictly along the segment and be zero at both ends. This too is impossible. The proof of (1) uses the actual periodic outer obstacle and quantitatively controls analytic continuation in the short direction.

\paragraph*{The exterior of a wide cell}

Normalize \(\lambda_1=1\) and take a proposed sequence of exact shells with \(L=\lambda_2\to\infty\). Rotate a shortest vector to the \(x\) axis and call it \(e\). It is primitive. Project the lattice orthogonally to \(e^\perp\); the result is a rank-two lattice. Choose its shortest vector and a reduced complementary basis vector, then rotate the transverse axes. Their coordinates can be taken as \((d,0),(t,q)\) with \(|t|\le d/2\) and \(q\ge\sqrt3\,d/2\). Their lifts have \(x\) offsets \(\delta_1,\delta_2\) in \([-1/2,1/2]\) after adding integer multiples of \(e\). In particular \(d\ge\sqrt{L^2-1/4}\).

The unique filled body \(O_L\) is convex and invariant under central inversion. It is in the central Voronoi cell, with the sharper bounds \(|\ell\cdot X|\le(1-s_b)|\ell|^2/2\) for every lattice vector. Apply the supporting affine bound from \(O_L+\ell\) at \(X+s_b\ell\), and use the unit Lipschitz bound on \(\nabla\Theta_L\) starting from \(X\in O_L\). Taking successively \(\ell=e,(\delta_1,d,0),(\delta_2,t,q)\) bounds \(|x|\) by \((1-s_b)/2\) and \(|y_1|/d,|y_2|/q\) by absolute constants. The body scaled by \((x,y_1,y_2)\mapsto(x,y_1/d,y_2/q)\) has volume \(b\), bounded diameter and central symmetry. It contains a fixed centered ball because its volume is at most a constant times its smallest supporting half-width. In physical coordinates the bodies therefore contain increasingly wide transverse disks in \(x=0\) and have a positive short-axis thickness.

In particular the right boundary is \(x=g_L(y)\) on a disk of radius \(c_b L\), and the left boundary is \(x=-g_L(-y)\). On the disk of half that radius, concavity and boundedness give \(|\nabla g_L|\le C_b/L\). Decreasing \(c_b\) if necessary, only the column of bodies \(O_L+ke\) meets this disk. Indeed for any nonzero projected period \(D\) choose a lift \((\delta,D)\) with \(|\delta|\le1/2\). The supporting bound gives \(|D\cdot y|\le(1-s_b)(|D|^2+1/4)/2+1/4\) in the central body. Its translate by \(D\) has transverse distance at least \(|D|/3\) from the origin for large \(L\), and \(|D|\ge d\).

Set \(R=c'_b L\) with \(c'_b>0\) small enough that all these slope and column bounds hold on \(|y|<10R\). In the gap \(g_L(y)<x<1-g_L(-y)\) the function \(H=1-(\Theta_L)_{xx}\) is harmonic and \(0\le H\le1\). Its values on both free boundaries are \(1-n_x^2\le C_b/L^2\), using the two slopes. Fix any \(|y_0|<8R\). On the gap over the disk \(|y-y_0|<R\) a harmonic majorant is

\[
\frac{C_b}{L^2}+\frac{|y-y_0|^2+2(4-(x-1/2)^2)}{R^2}.
\tag{2}
\]

Indeed \(|x-1/2|\le1\) throughout the gap and (2) is at least one on the circular sides. The maximum principle at \(y_0\) gives \(0\le H\le C_b/L^2\) throughout the gap over \(|y|<8R\). Positive semidefiniteness and unit trace give \(|(\Theta_L)_{xy}|\le C_b/L\) there. These bounds concern the physical exterior up to both true surfaces.

\paragraph*{A fixed strip with two constant boundary values}

Set \(A=s_b\) and \(a=(1-A)/2\). Use \(p=(\Theta_L)_x\) as the normal coordinate in the gap and keep \(y\) as the two tangential coordinates. The values on its two faces are zero and \(A\), and \(p_x=1+O_b(L^{-2})\). Thus the inverse \(x=h(p,y)\) exists on the fixed strip \(0\le p\le A\) over \(|y|<8R\). Take the partial Legendre transform in the normal direction,

\[
V(p,y)=p h(p,y)-\Theta_L(h(p,y),y),\qquad
w(p,y)=V(p,y)-ap-p^2/2 .
\]

At the first face \(\Theta_L=0\), and at the next face \(\Theta_L=A(h(A,y)-1/2)\) by the short-period shift of its supporting plane. Consequently \(V(0,y)=0\) and \(V(A,y)=A/2\). Both values are independent of \(y\), including for sheared lattices. We have \(V_p=h\), \(V_{pp}=1/(\Theta_L)_{xx}\), \(V_{py}=-(\Theta_L)_{xy}/(\Theta_L)_{xx}\), and \(V_{yy}=-(\Theta_L)_{yy}+(\Theta_L)_{yx}(\Theta_L)_{xy}/(\Theta_L)_{xx}\). Using \(\Delta\Theta_L=1\) gives \((1+\Delta_y V)V_{pp}-|\nabla_y V_p|^2=1\). In particular the perturbation of the flat solution satisfies

\[
\Delta_{p,y}w+w_{pp}\Delta_y w-|\nabla_y w_p|^2=0,\qquad
w(0,y)=w(A,y)=0,\qquad g_L(y)=a+w_p(0,y).
\tag{3}
\]

The physical-side estimates give \(|w_{pp}|+|w_{yy}|=O_b(L^{-2})\) and \(|w_{py}|=O_b(L^{-1})\). The two constant boundary values give \(|w|+|w_p|=O_b(L^{-2})\) by one-dimensional integration. They also give \(|\nabla_y w|=O_b(L^{-1})\) by integrating \(w_{py}\) from a face. All these estimates hold on \(|y|<8R\).

We record the finite real derivative bound needed for analytic continuation. For a fixed \(0<\alpha<1\), the norm of \(w\) in \(C^{2,\alpha}\) on \(|y|<6R\) is \(O_b(L^{-1})\). Every individual solution is smooth up to the two faces by free-boundary regularity. For a uniform estimate write \(\delta=C_b/L\). The nonlinear term in (3) is quadratic in the Hessian. Its Hölder seminorm is at most \(C\delta[D^2w]_\alpha\) wherever \(|D^2w|\le\delta\). Apply the local linear Poisson Schauder estimates on interior balls and on half balls at either flat Dirichlet face. To see the absorption explicitly, use radii a small fixed multiple of the lesser of \(1,A/8\) and the distance to an artificial lateral boundary at \(|y|=7R\). Let \(M_1\) be the supremum of the corresponding radius to the power \(\alpha\) times the local seminorm \([D^2w]_\alpha\), using half balls or nearby interior balls as appropriate. On a doubled ball this quantity is at most \(C M_1+C\delta\) by covering with balls of comparable radii. Subtract the first Taylor polynomial at the center before rescaling by the square of the radius. The rescaled supremum is at most \(C\delta\) by the Hessian bound. At a Dirichlet face center that polynomial has zero boundary trace. The unit-scale linear estimate on a smaller ball gives \(M_1\le C\delta+C\delta M_1\). All constants are independent of \(R\) and \(M_1\) is finite before absorption by smoothness on the closed smaller patch. This proves the bound throughout \(|y|<6R\), where the local radii are bounded below.

\paragraph*{Holomorphic tangential translations}

We show that \(g_L-a\) extends holomorphically on a tangential polydisk of radius \(cR\) about zero, with supremum \(O_b(L^{-1})\). Cut \(w\) off in \(y\) to be unchanged on \(|y|<4R\) and zero by \(|y|=5R\), then extend periodically on a square transverse torus of side \(16R\). The extension \(w_0\) is \(O_b(L^{-1})\) in \(C^{2,\alpha}\) and has homogeneous Dirichlet data on both strip faces. Denote the differential expression in (3) by \(E(w)\). The forcing \(F_0=E(w_0)\) is supported outside \(|y|<4R\) and is \(O_b(L^{-1})\) in \(C^\alpha\).

The linearization \(D E(0)=\Delta_{p,y}\) with these two homogeneous boundary conditions has inverse bounded from \(C^\alpha\) to \(C^{2,\alpha}\) independently of \(R\ge1\). The periodic strip Dirichlet problem is solvable, and \(p(A-p)\|F\|_\infty/2\) bounds the solution in absolute value for real forcing, by the maximum principle. Real and imaginary parts give the complex bound. Fixed-radius interior and flat Dirichlet boundary Schauder estimates give the claimed uniform inverse bound.

Choose a smooth cutoff \(\eta(y/R)\) equal to one near zero and zero outside \(|y|<2R\). For real two-vectors \(t\) with \(|t|<cR\) use the diffeomorphism \(T_t(y)=y+\eta(y/R)t\). Pull \(E\) back by \(T_t\) so that \(E_t(w_0\circ T_t)=F_0\circ T_t=F_0\). The last equality holds since the map is the identity on the support of the forcing and maps its moving region inside \(|y|<4R\). This operator has coefficients holomorphic in complex \(t_i\) for \(|t_i|<cR\) after decreasing \(c\). The chain rule uses \((I+t\otimes\nabla_y[\eta(y/R)])^{-1}\) and its derivatives. Their deviations from the identity change of variables are \(O(c)\) in the required Hölder norms, independently of \(R\). Dot products in (3) extend bilinearly.

Solve \(E_t(u)=F_0\) by fixed point iteration with the bounded inverse of \(\Delta_{p,y}\). On a \(C^{2,\alpha}\) ball of radius \(r\) the Lipschitz constant of \(E_t-\Delta_{p,y}\) is at most \(C(c+r)\) by (3), the chain rule and Hölder multiplication. Choose \(c,r\) small independently of \(L\). The iteration is a contraction for all large \(L\), and its unique small solution has norm \(O(\|F_0\|_{C^\alpha})=O_b(L^{-1})\). The iterates are holomorphic in \(t\) and converge uniformly. For real \(t\) the small solution is \(w_0(p,T_t(y))\) by uniqueness. Taking its \(p\) derivative at \((0,0)\) gives \(w_p(0,t)=g_L(t)-a\) there, with the asserted holomorphic extension and bound. This use of the two exact boundary values controls tangential analytic radii that grow with the period ratio.

\paragraph*{Continuing through the full short thickness}

Flatten the upper graph by writing the Schwarz potential as \(G(s,y)=W(g_L(y)+s,y)\). Its Cauchy data are \(G=G_s=0\) at \(s=0\). The equation is

\[
(1+\nabla g_L\cdot\nabla g_L)G_{ss}-2\nabla g_L\cdot\nabla_y G_s+\Delta_y G-(\Delta_y g_L)G_s=1.
\tag{4}
\]

Let \(M=cR/4\), \(y=MY\), \(s=MT\), \(G=M^2B\). On \(|Y_i|<1\), Cauchy's estimates on the larger polydisk give \(|\nabla g_L(MY)|+M|D^2g_L(MY)|\le C_b/L^2\). Thus the scaled coefficients of (4) are bounded holomorphic functions differing from the constant flat coefficients by \(O_b(L^{-2})\). The initial data remain zero.

Here is an explicit analytic majorant sufficient for the continuation. Write (4) in scaled coordinates as \(\partial_T U=\mathcal LU\) on \(U=(B_T,B_{Y_1},B_{Y_2},1)\), solving its first row for \(B_{TT}\) and differentiating \(B_T\) in the next two rows. The coefficients of the first-order operator \(\mathcal L\) depend only on \(Y\). Its last row is zero and \(U(0)=(0,0,0,1)\). On concentric polydisks of radii \(r'<r\le1\), Cauchy's estimate bounds the operator norm on holomorphic sup norms by \(C/(r-r')\). Apply it \(k\) times, shrinking from radius one to \(1/2\) in steps of \(1/(2k)\). Then \(\|\mathcal L^kU(0)\|_{1/2}\le(2Ck)^k\). The Taylor series \(\sum T^k\mathcal L^kU(0)/k!\) converges uniformly on a fixed \(|T|<T_0\) independent of \(L\). Telescoping the \(k\) factors against the flat operator bounds the difference by \(C_b L^{-2}k(2Ck)^k\) before multiplication by \(T^k/k!\). Integrating the first component gives \(B=T^2/2+O_b(L^{-2})\), also with derivatives on smaller polydisks. This bounds the continued germ on both sides of the graph. It agrees with the physical exterior solution near the right surface by uniqueness for harmonic differences with zero Cauchy data.

On any fixed physical normal interval \(|s|\le1\) about the central upper graph, \(T=s/M\to0\). Changing variables back gives \(D^2W\to e\otimes e\) there. In the tangential second derivatives the term containing \(B_T\) is multiplied by \(M D_y^2g_L\), which tends to zero by the same Cauchy estimate. All other terms follow from convergence of the scaled second derivatives and \(\nabla g_L\to0\). This proves (1) on the whole central segment \([-g_L(0),g_L(0)]\), whose length tends to \(2a<1\).

For all sufficiently late members of the sequence, \(\kappa I+\rho D^2W\) is positive definite on the segment and \(W_{xx}>1/2\) there. Starting at the right face, the cavity potential agrees with this branch until any first cavity by harmonic uniqueness. The exterior cavity Hessian at a first contact would be rank one. If there is no such contact, \((\nabla\Xi-\kappa X)/\rho\) is zero at both outer faces by the exterior identity. Its \(x\) component changes by at least \(g_L(0)>0\) between them. Both alternatives contradict attainment. Since the argument applies to any proposed sequence \(\lambda_2/\lambda_1\to\infty\), the claimed finite \(M_*\) exists.

At any fixed \(b\in(0,1)\) the estimate (1) depends only on the outer obstacle. For any genuine coating with the same \(b\) and the stated potential identities, first contact must therefore satisfy

\[
\frac{f}{b-f}=\frac{3\kappa}{\rho}\ \le\ 3 C_b(\lambda_1/\lambda_2)^2
\]

once the period ratio is large enough for (1) and \(W_{xx}>1/2\). Indeed some cavity must meet the segment by the two endpoint gradients. At contact the smallest eigenvalue of \(D^2W\) is \(-\kappa/\rho\), whereas (1) bounds it below by \(-C_b(\lambda_1/\lambda_2)^2\). This gives the exclusion for every fixed \(0<f<b<1\) in the same geometric equality class.
\phantomsection\label{endfile:flat-lattice}

\subsection{The cubic measurements, symmetry orbits and axial covering test}\label{detail:cavity-count-proof}\label{file:cavity-count-proof}
\paragraph*{Boundary regularity}

Take a bounded connected filled body \(O\) with separated translates and connected lattice exterior. Its finitely many nonempty open cavities \(K_j\) have pairwise disjoint closures contained in \(O\) and leave a connected open shell. Write \(K=\bigcup K_j\) with periodic copies understood in periodic potentials. For separated Lipschitz interfaces, equality implies that the outer surface and every cavity surface are \textbf{real analytic}. This upgrade holds on every three-dimensional Bravais lattice.

Here is why the geometric equality conditions apply before any classical interface differentiations. For an attaining cell the matrix of three electric fields \(E\) is symmetric, its trace is constant in each material component, and the weak component flux sums to zero as in (2) of the topology proof. Each individual phase-\(q\) component can be isolated by a smooth cutoff equal to one on its closure and zero near the other phase-\(q\) components. Within a bounded cavity, test its divergence-free flux \(E-l(T)I\) with its own coordinate potentials after subtracting \(l(T)x\), with a cutoff in the surrounding shell. This gives \(E=l(T)I\) there. On the exterior the periodic-corrector test and the two calibrated field lines give the exterior trace. From (5) of the same proof one then obtains, for \(\Delta u_A=\mathbf1_A-\theta_A\) and \(H_A=D^2u_A\), \[
 H_O=(1-b)I/3\quad\hbox{in }O,\qquad
 H_K=(1-v)I/3\quad\hbox{in }K,\qquad
 H_K=(v/b)H_O\quad\hbox{in the exterior of }O.                 \tag{1}
\] Indeed the two real Hessians are scalar in a cavity, \(H_O\) continues throughout the filled body, and the exterior trace identifies the proportionality outside. These steps use null interface volume, separated components and the distributional field equations.

Near one component of either \(A=O\) or \(A=K\), subtract from \((1-\theta_A)|x|^2/6-u_A\) its interior affine function. The resulting \(U\) satisfies \(\Delta U=\mathbf1_{A^c}\) locally and \(U=\nabla U=0\) on the closed component. Continuity of the gradient follows from the Poisson equation with bounded right-hand side. At a Lipschitz boundary point the interior cone supplies a ball of radius proportional to \(r\) in the vanishing set in every sufficiently small \(r\)-ball. The no-sign free-boundary regularity theorem of Caffarelli, Karp and Shahgholian therefore makes the interface analytic near that point. The same argument works pointwise wherever the vanishing side has this thickness, once (1) holds. In particular a Lipschitz design with even one nonanalytic interface point is excluded, irrespective of the finite cavity count.

\paragraph*{Contacts along points and curves}

Suppose the individual filled bodies and the individual cavities are Lipschitz domains with the same connected open shell and lifted exterior as above, but initially allow closures to touch. Assume their interiors have the prescribed inclusions and disjointness, and all contacts between distinct boundaries lie in finitely many points and compact \(C^1\) arcs per cell. In this class equality forces all the surfaces to be separated. This assertion holds for every Bravais lattice and every finite cavity count at the stated optical data.

The component fluxes still separate. Away from the contact set each separate divergence vanishes by locality. The contact set has zero \(W^{1,2}\) capacity in three dimensions. A cutoff rising logarithmically from zero at distance \(r^2\) to one at distance \(r\) has squared gradient integral \(O(1/|\log r|)\), since the \(s\)-neighborhood of finitely many finite-length arcs has volume \(O(s^2)\). Testing the divergence of the individual \(L^2\) flux with this cutoff proves its vanishing across the whole contact set. The bounded-cavity and connected-exterior tests then give (1), as before. If filled translates touched at \(x\), continuity and periodicity of \(\nabla u_O\) would identify \(s_bx+d\) with \(s_b(x-\ell)+d\) for their nonzero translation \(\ell\), which is impossible.

Suppose several cavity closures meet at \(x\). On each of them \(\Xi=s_v|x|^2/2-u_K\) is affine. Continuity of \(\Xi\) and its gradient makes those affine functions the same plane \(\ell\). In a small ball meeting only the closures through \(x\), the function \(\Xi-\ell\) and its gradient vanish on their union, and its Laplacian is one off that union. A cone in any one of the Lipschitz cavities provides the required thickness. The free-boundary theorem makes the boundary of the union a single analytic graph near \(x\). To see the contradiction, first note that \(x\) cannot be interior to the union. Deleting the finitely many arcs and points from a small open ball leaves it path connected, whereas the remaining space would split between distinct open cavity interiors. At a graph point, the same argument applies on the open filled side of a sufficiently small graph chart, where both cavity interiors approach \(x\). Deleting smooth arcs and points cannot disconnect a three-dimensional open connected set, as seen by local polygonal path perturbations. This excludes contact between cavities.

The free-boundary regularity argument can now be applied independently to \(O\) and each cavity, even at a prospective contact between their surfaces. At such a point their normals agree since one body contains the other. On the matrix side, continuity up to the analytic surfaces gives the cavity value \(D^2\Xi=n\otimes n\). The outer identity gives \(D^2\Xi=\kappa I+\rho n\otimes n\), with \(\rho=v/b\) and \(\kappa=(1-\rho)/3>0\). The tangential entries disagree. This excludes the last type of contact.

\paragraph*{An angular derivative determines an ordering of the outer radii}

Put \(s_\theta=(1-\theta)/3\), \(\rho=v/b\) and \(\kappa=(1-\rho)/3\). On the lift let \(\Xi=s_v|x|^2/2-u_K\) and \(\Theta=s_b|x|^2/2-u_O\), with the origin and additive constants chosen so that \(\Theta=0\) in the central body. The exterior Hessian bound makes \(\Xi\) convex and strictly convex locally off the closed cavity array. Each cavity is a separate region on which \(\nabla\Xi=\kappa c_j\) is constant. On the central shell set \(P=\Xi-\kappa|x|^2/2-C\) with the exterior constant \(C\) chosen using \(\Xi=\kappa|x|^2/2+\rho\Theta+C\) in the matrix. The linear term vanishes by central symmetry. Thus \[
\Delta P=\rho\ \text{in the shell},\qquad P=\nabla P=0\ \text{on }\partial O,\qquad
\nabla P=\kappa(c_j-x)\ \text{in }K_j .
\] The map \(\nabla\Xi/\kappa\) fixes \(\partial O\) and sends the shell diffeomorphically onto \(O\) minus the distinct points \(c_j\). Uniqueness of \(O\) and of a finite \(K\) at fixed lattice and fractions makes both invariant under the lattice point group about the chosen origin. Hence the labels transform by signed coordinate permutations.

Use Cartesian coordinates \((x,y,z)\) and the rotation generator \(L=x\partial_y-y\partial_x\). In the wedge \(x>y>0\), a cavity meeting the wedge has label components \(c_x\ge c_y\ge0\). Indeed monotonicity of the gradient of \(\Xi\), paired at a point and its reflection across \(y=0\) or \(x=y\), gives the two inequalities. In that cavity, \[
 LP=\kappa(xc_y-yc_x).
\] If all noncentral labels are on coordinate axes, this is nonpositive. If they are all on body diagonals, it is nonnegative. A central label contributes zero in either case. Since \(L\) commutes with \(\Delta\), \(LP\) is harmonic in the shell. It vanishes on both wedge faces by reflection and on \(\partial O\) by the Cauchy data. The maximum principle gives the corresponding sign throughout the part of the shell in the wedge, component by component.

At a point of \(\partial O\) in the wedge, with outward unit normal \(n\), the shell-side Hessian is \(\rho n\otimes n\). Consequently \(\partial_n LP=\rho(-y,x,0)\cdot n\). For the axial labels the outward derivative is nonnegative, and for the diagonal labels it is nonpositive. Parametrize the convex surface by \(R(\omega)\omega\) about its interior centre. Rotating \(\omega\) in the \(xy\) plane through increasing azimuth \(\phi\) toward \(x=y\), while keeping \(\omega_z\) fixed, gives \[
(\omega\cdot n)\partial_\phi R=-(-y,x,0)\cdot n.
\] The factor on the left is positive. Thus \textbf{balancing two squared direction coordinates decreases \(R\) for axial labels and increases \(R\) for diagonal labels}, with weak inequalities sufficient. Signed permutations give the rule for every pair. Transferring all mass in the squared coordinates to the largest one reaches \(100\). Successively balancing the coordinates to \((1/3,1/3,1/3)\) reaches \(111\). In particular the axial case has maximum radius \(R_{100}\) and \(R_{100}\ge R_{110}\ge R_{111}\). The diagonal case has maximum radius \(R_{111}\) and the reverse ordering.

\paragraph*{Real exterior measurements on the simple cubic lattice}

For permittivities \(-15+i\) and \(1+i\) at phase-1 fraction \(1/4\), the four exterior tests below force \textbf{at least 44 cavities} in an attaining finite simple-cubic shell. The additional axial covering test excludes 44 and 48, leaving only 45 and 49 below 50. Count 45 requires one orbit of twelve associated points on face diagonals, one orbit of eight on body diagonals and one orbit of 24 signed permutations of \((a,a,b)\) with \(0.27<b/a<0.60\). It also requires the central cavity. The lower bound of 45 holds throughout an unquantified open neighborhood of the simple cubic lattice.

Each of the twelve face-diagonal outer points needs a cavity within \(0.043001\) periods. Each of the eight body-diagonal points needs one within \(0.064001\) by the tighter test. The orbit classification below already follows at \(0.065001\). Write \(x\ge y\ge z\ge0\) for the sorted absolute components of a cavity's associated point. The first distance test requires \(x>0\) and \(z/x<0.30\). The second requires \(x>0\) and \(z/x>0.50\). Two further tests on the surface directions \((1,1,0.24)\) and \((1,1,0.53)\) rule out \textbf{every} configuration whose noncentral associated points are confined to coordinate axes, face diagonals and body diagonals, for any finite count.

An additional exterior test at direction \((1,1,0.38)\) forces an entire orbit of at least 24 cavities that stay more than \(0.00726\) periods from every coordinate plane and every body-diagonal line. This holds for every finite count.

There is also a criterion for excluding \textbf{every} finite count using just the real Schwarz potential of the true outer surface. At the face diagonal let \(n=(1,1,0)/\sqrt2\), \(u=(1,-1,0)/\sqrt2\), and let the third axis point vertically. Use \(P_0=0.488491617\,n\). Write \(W\) for the outer Schwarz germ solving \(\Delta W=1\) with zero Cauchy data at the true outer boundary. With \(b=832/(778+285\sqrt{17})\), \(\rho=1-1/(4b)\) and \(\kappa=(1-\rho)/3\), put

\[
Y=\kappa |X|^2/2+\rho W,\qquad S=Y_{uu}+Y_{zz},\qquad J=Y_{uu}-Y_{zz}.
\]

Suppose \(W\) extends with \(\Delta W=1\) to the \textbf{ordinary real ball} \(B(P_0,R)\subset\mathbb R^3\), and \(|S|,|J|\le M\) on that ball. Each row of the table suffices to exclude a finite attaining shell. The bounds refer to inward depth from \(P_0\).

{\def\LTcaptype{none} 
\begin{longtable}[]{@{}
  >{\raggedleft\arraybackslash}p{(\linewidth - 10\tabcolsep) * \real{0.1667}}
  >{\raggedleft\arraybackslash}p{(\linewidth - 10\tabcolsep) * \real{0.1667}}
  >{\raggedleft\arraybackslash}p{(\linewidth - 10\tabcolsep) * \real{0.1667}}
  >{\raggedleft\arraybackslash}p{(\linewidth - 10\tabcolsep) * \real{0.1667}}
  >{\raggedleft\arraybackslash}p{(\linewidth - 10\tabcolsep) * \real{0.1667}}
  >{\raggedleft\arraybackslash}p{(\linewidth - 10\tabcolsep) * \real{0.1667}}@{}}
\toprule\noalign{}
\begin{minipage}[b]{\linewidth}\raggedleft
\(R\)
\end{minipage} & \begin{minipage}[b]{\linewidth}\raggedleft
\(M\)
\end{minipage} & \begin{minipage}[b]{\linewidth}\raggedleft
Upper \(Y_{uu}\) at depth \(0.042\)
\end{minipage} & \begin{minipage}[b]{\linewidth}\raggedleft
Lower \(Y_{zz}\) through \(0.042\)
\end{minipage} & \begin{minipage}[b]{\linewidth}\raggedleft
Lower \(Y_{nn}\) through \(0.042\)
\end{minipage} & \begin{minipage}[b]{\linewidth}\raggedleft
Lower \(Y_{uu}\) through \(0.028\)
\end{minipage} \\
\midrule\noalign{}
\endhead
\bottomrule\noalign{}
\endlastfoot
\(0.064\) & 1 & \(-0.0048\) & \(0.115\) & \(0.606\) & \(0.055\) \\
\(0.068\) & 2 & \(-0.0023\) & \(0.106\) & \(0.604\) & \(0.050\) \\
\(0.075\) & 3 & \(-0.032\) & \(0.131\) & \(0.603\) & \(0.050\) \\
\(0.080\) & 5 & \(-0.039\) & \(0.129\) & \(0.601\) & \(0.046\) \\
\end{longtable}
}

The Hessian is diagonal on this ray. Thus under any row exactly one eigenvalue first reaches zero, at depth between \(0.0279999\) and \(0.042\) from the \textbf{true surface}. This says the eigenvalue has multiplicity one at first contact. At a cavity contact the limiting Hessian must be a rank-one normal projection, and in the shell it must be positive definite. The ray therefore rules out every finite count under the stated hypothesis. The real inequalities for the Schwarz germ on one of these full balls remain unproved. The measurements use the exact real obstacle comparisons and do not use a fitted Schwarz continuation.

Both results apply to a connected phase-1 shell with finitely many separated Lipschitz cavities, bounded separated outer translates and connected matrix. Equality makes the interfaces analytic and fixes the outer body up to translation. Existence for unrestricted finite counts remains open.

The tangent-plane reflection also gives an all-count criterion involving only the ordinary circular mean of the exterior potential. A coefficient norm of at most \(0.014\) for its centered, normalized remainder from an explicit rational profile at radius \(0.060\) would force a longitudinal derivative ratio above the cavity-contact limit. The true remainder is below \(0.00361\) on the real interval \([0.015,0.060]\).

\paragraph*{Two measured circular averages}

For the conditional test use \(F=Y\). For the location test suppose a real ball centered at \(P_0\) contains no cavity. Use the actual cavity potential \(F=\Xi\) on that ball, so \(\Delta F=1\) and \(0\le D^2F\le I\). On the exterior used for the measurements, in both cases \[
F=C+\kappa |X|^2/2+\rho\Theta(X)
\] for the certified physical obstacle potential \(\Theta\). For the body-diagonal location test use \(P_0=0.522104(1,1,1)/\sqrt3\) and tangents \((1,-1,0)/\sqrt2,(1,1,-2)/\sqrt6\). On the face and body diagonals the true radial heights are within \(10^{-6}\) of the exterior centers and the supporting plane separates every measured circle from the body. At the two intermediate points the checker directly proves positivity of \(\Theta\) at every sampled point. In all four tests the points and balls stay inside the central cell.

Normalize lengths about the ball center by \(R\). At axial coordinate \(pR\) take a transverse circle of radius \(r=qR\). Use \(\langle\cdot\rangle_\phi\) for azimuthal mean and \(F_0\) for the value at its center. The observations are \[
B_s=\frac{4}{r^2}\big(\langle F\rangle_\phi-F_0\big),\qquad
B_j=\frac{8}{r^2}\langle F\cos(2\phi)\rangle_\phi .
\] For \(S=F_{uu}+F_{zz}\) and \(J=F_{uu}-F_{zz}\) let \(m_s(\mu),m_j(\mu)\) be their azimuthal means on the sphere of radius \(R\), where \(\mu\) is the cosine of the angle with the ray. If using a bounded harmonic field on an open ball, take a weak-star limit of bounded traces on concentric spheres tending to its boundary. The Poisson kernel then gives the interior field. At normalized axial position \(t\) its value for either harmonic field is \[
\frac12\int_{-1}^{1}K_t(\mu)m(\mu)\,d\mu,\qquad
K_t(\mu)=\frac{1-t^2}{(1+t^2-2t\mu)^{3/2}}
=\sum_{\ell\ge0}(2\ell+1)t^\ell P_\ell(\mu).
\] The two circular observations have the same integral form with \(m_s,m_j\) respectively and kernels \(\sum_\ell (2\ell+1)w_\ell P_\ell(\mu)\). Put \(\eta=\sqrt{p^2+q^2}\) and \(Z_\ell=\eta^\ell P_\ell(p/\eta)\). Their moments are \[
w_\ell^s=\frac{4(p^{\ell+2}-Z_{\ell+2})}{(\ell+1)(\ell+2)q^2},\qquad
w_\ell^j=\eta^\ell\frac{P_{\ell+2}''(p/\eta)}{P_{\ell+2}''(1)} .                 \tag{1}
\]

To see which fields these measure, denote the circle means of \(F\) and \(F\cos2\phi\) by \(B(r)\) and \(C(r)\). Polar differentiation and integration by parts give \(\langle S\rangle=B''+B'/r\) and \(\langle J\rangle=C''+3C'/r\). Hence the measured \(B_s\) is an average of the means of \(S\) on circles of radii \(rt\) with positive normalized weight \(4t\log(1/t)\,dt\), \(0<t<1\). The corresponding weight for \(B_j\) and \(J\) is \(4(t-t^3)\,dt\). In particular \(|w_\ell^s|,|w_\ell^j|\le\eta^\ell\) by the zonal harmonic formula and \(|P_\ell|\le1\).

For (1), test these averages on a zonal solid harmonic of degree \(\ell\). A potential with that transverse Laplacian is the negative zonal solid harmonic of degree \(\ell+2\) divided by \((\ell+1)(\ell+2)\), since two longitudinal derivatives lower the degree. Taking the difference of its circle value and its axial value gives \(w_\ell^s\). For \(w_\ell^j\) use the harmonic polynomial obtained by multiplying \((u^2-z^2)|X|^\ell P_{\ell+2}''(X_n/|X|)\) by \(1/[4P_{\ell+2}''(1)]\). This is the azimuthal mode of order two of a solid harmonic. The azimuthal mean of its transverse Hessian difference is the zonal solid harmonic of degree \(\ell\), as follows from its axial value and harmonicity. Its cosine-weighted mean gives the second formula. The regular quadratics give the same degree-zero moments. The checker uses the Legendre recurrence for \(Z_\ell\) and the normalized Gegenbauer recurrence \((\ell+5)w_{\ell+1}^j=(2\ell+5)p w_\ell^j-\ell\eta^2w_{\ell-1}^j\). The first two values are \(1,p\).

\paragraph*{A rigorous enclosure of the ordinary outer surface.}
The next comparison uses the cubic fundamental chamber. The obstacle,
volume and the source values refer to the ordinary filled body $O$ at
this optical datum. The comparison and circular averages are computed
on the same unit-period torus.
\paragraph*{A comparison on a symmetry chamber}

Write \(b=832/(778+285\sqrt{17})\) and \(S=(1-b)/6\). The central physical obstacle potential \(\Theta=(1-b)|X|^2/6-u_O\) has been normalized to vanish on the filled body. It is nonnegative and \(C^{1,1}\) with \(\Delta\Theta=\mathbf1_{\{\Theta>0\}}\) on the cube. Reflection symmetry and periodicity give its outward flux \(S\) on \(X=1/2\) and zero on the other faces of the chamber \(0\le Y,Z\le X\le1/2\). There is no ordering between \(Y\) and \(Z\).

Use \(V=(Y/X)^2\), \(W=(Z/X)^2\), \(h=1/2-g(V,W)\) and \(t=(X-g(V,W))/h\). The file of integer coefficients defines \(g(V,W)\) and \(P(t,V,W)\) in tensor products of shifted Chebyshev polynomials \(T_j(2s-1)\), all with denominator \(2^{60}\). Both are symmetric in \(V,W\). Define \(F=t^2 P\) where \(t>0\) and \(F=0\) where \(t\le0\). This has zero value and first derivatives at \(t=0\). All estimates used below hold for \(0\le t\le1.0006\) and \(0\le V,W\le1.0001\), unless a smaller interval is specified.

{\def\LTcaptype{none} 
\begin{longtable}[]{@{}
  >{\raggedright\arraybackslash}p{(\linewidth - 2\tabcolsep) * \real{0.4286}}
  >{\raggedleft\arraybackslash}p{(\linewidth - 2\tabcolsep) * \real{0.5714}}@{}}
\toprule\noalign{}
\begin{minipage}[b]{\linewidth}\raggedright
Quantity
\end{minipage} & \begin{minipage}[b]{\linewidth}\raggedleft
Certified bound
\end{minipage} \\
\midrule\noalign{}
\endhead
\bottomrule\noalign{}
\endlastfoot
\(g\) & \((0.28,0.43)\) \\
\(P\) & \((0.00378,0.06001)\) \\
\(\lvert\Delta F-1\rvert\) on the polynomial branch & \(<5.33\times10^{-8}\) \\
\(\lvert F_X-S\rvert\) at \(X=1/2\) & \(<6.5\times10^{-11}\) \\
\(\lvert N(t,1,W)\rvert\) & \(<2.33\times10^{-10}\) \\
\(\partial_m N\) for \(0.9999\le m\le1.00005\) & \(<-0.0137\) \\
\(\partial_X\vert_{V,W} F_X\) for \(0.9994\le t\le1.0006\) & \(>1/4\) \\
\end{longtable}
}

Here \(F_X\) denotes the Cartesian partial derivative. The last row differentiates that partial derivative along a ray. For \(m=Y/X\) the seam flux is \(F_X-F_Y=tN/(hX)\), with \[
N=\{X+2(1-t)[(V+m)g_V+Wg_W]\}(2P+tP_t)
       -2ht[(V+m)P_V+WP_W],\qquad V=m^2 .
\]

Take \(c=1/20\), \(C=(c,0,0)\), \(\eta=4\times10^{-7}\), \(\lambda_\pm=1\pm\eta\), \(T_\lambda(A)=C+(A-C)/\lambda\), and \(\epsilon=5.33\times10^{-8}\). Two comparison functions are \[
L=\frac{\lambda_+^2}{1-\epsilon}F\circ T_{\lambda_+},\qquad
U=\frac{\lambda_-^2}{1+\epsilon}F\circ T_{\lambda_-}.
\tag{2}
\]

At transformed \(X'\le0.27\) use the zero branch. In an open neighborhood of the cutoff \(g>0.28\), so this definition is consistent. Whenever \(X'\ge0.27\) in the transformed chamber, \(Y'/X'\) and \(Z'/X'\) are at most \(1+c\eta/[(1-\eta)0.27]\). At its far face \(|X'-1/2|\le(1/2-c)\eta/(1-\eta)\). The checks on the enlarged polynomial domains therefore cover the whole positive part. The surfaces of broken symmetry are the planes \(X=Y\) and \(X=Z\). On \(X=Y\), the transformed ratio is \(m=1-c(\lambda-1)/(\lambda X')\). Here \(\lambda X'\le0.501\). Since \[
0.0137\,c\eta/0.501>2.33\times10^{-10},
\] the sign of \(N\) is positive for \(L\) and negative for \(U\) when \(t>0\). The outward flux of \(L\) on the seam is consequently nonpositive, and that of \(U\) is nonnegative. Interchanging \(Y\) and \(Z\) gives the other seam. Both coordinate planes retain their exact zero flux.

On \(X=1/2\) the transformed abscissa is \(X'=1/2-(1/2-c)(\lambda-1)/\lambda\). Integrate the verified ray derivative of \(F_X\) from \(1/2\) to \(X'\) with \(V',W'\) fixed. The inequality \[
[(1/2-c)/4-S]\eta>S\epsilon+(1+\eta)6.5\times10^{-11}
\] then gives \(\partial_\nu L\le S\le\partial_\nu U\). The scaling also gives \(\Delta L\ge1\) on \(\{L>0\}\) and \(\Delta U\le1\) everywhere, including the zero side. The continuous gradients produce no interface measure. Test the weak equations on the chamber by \((L-\Theta)_+\) and \((\Theta-U)_+\). On the support of the first, \(\Delta L\ge\Delta\Theta\), and on the support of the second, \(\Delta\Theta\ge\Delta U\). The flux comparisons give zero Dirichlet integral for both positive parts. Each vanishes near the center, so \(L\le\Theta\le U\) throughout. In particular \(\{U=0\}\subset\{\Theta=0\}\subset\{L=0\}\). Reflecting the chamber covers the entire centered cell.

On the face diagonal use \(X=Y=R/\sqrt2\) and \(Z=0\). The checker evaluates \(X'-g((Y'/X')^2,0)\) by interval Clenshaw recurrence. The initial dilation pair gives a zero neighborhood at radius \(0.488491403\) and strict positivity at \(0.488491767\). The sharper pair below gives (1). For the other directions, fix \((V,W)\) and test at \(X=g(V,W)\pm2.1\times10^{-7}\). One has \(|X'-X|\le(0.431-c)\eta/(1-\eta)\) and \(|V'-V|,|W'-W|\le(1-c\eta/0.279)^{-2}-1\). The coefficient sum bounding \(|g_V|+|g_W|\) proves the strict signs there. Multiplication by \(\sqrt{1+V+W}\le\sqrt3\) gives a uniform radial error below \(3.64\times10^{-7}\) for this initial pair.

\paragraph*{Checking the polynomial inequalities}

The obstacle verifier evaluates the integer coefficients of the chamber graph and the Poisson comparison polynomial with interval arithmetic.

For the Laplacian the checker clears the denominator \(h^2X^2\). In the Cartesian metric put \(M_{VV}=4V(V+1)\), \(M_{WW}=4W(W+1)\) and \(M_{VW}=4VW\), so \(M_{ij}=X^2\nabla i\cdot\nabla j\). Also \(X^2\Delta V=6V+2\) and \(X^2\Delta W=6W+2\). These identities, \(X=g+ht\) and the chain rule give the polynomial checked by the program. It separately checks the seam flux numerator, its derivative, \(hX(F_X-S)\), \(h^2X^2(\partial_X\vert_{V,W}F_X-1/4)\) and \(P\). On the nonnegative branch \(hX\ge g(1/2-g)>0.030\).

Each polynomial is evaluated at a tensor product of cosine extrema. If a univariate polynomial has degree \(d\) on a rescaled interval with \(N\) cosine subintervals, its supremum norm is bounded by its sampled maximum modulus divided by \(1-\tfrac12(\pi d/(2N))^2\). Indeed its cosine pullback has derivative zero at a folded endpoint or an interior extremum of its absolute value, and the trigonometric Bernstein inequality bounds its second derivative by \(d^2\) times its norm. Apply this bound in each variable to the polynomial minus the center of its sampled range. Derivatives of the input polynomials are taken exactly in the shifted Chebyshev basis. Arb encloses the basis tables and final norm factors. The sequential binary64 sums have explicit absolute error bounds from their coefficient sums, and subsequent array arithmetic carries outward absolute errors. The degree bounds, enlarged coordinate intervals and all strict comparison margins are checked in the script.

These real enclosures locate the outer body for the prescribed optical data. To exclude every finite cavity count by the Schwarz ray, two continued Hessian entries must stay above the cavity threshold and the third must pass below it. A real radial error controls the physical contact set. A bound on the centered graph in complex directions would transfer the verified first crossing of the comparison surface.
\paragraph*{Sharper bounds from the same comparison potential}

At the quarter-fraction data the true radius satisfies \[
0.488491555<R_{110}<0.488491616 .
\] The true radius differs from the degree-48 input by less than \(7.66\times10^{-8}\) on every ray. In face-diagonal coordinates their centered heights differ by less than \(8.53\times10^{-8}\) on the physical tangent disk of radius \(0.14\) and \(8.63\times10^{-8}\) on the disk of radius \(0.15\).

For the rounded harmonic coefficients, averaging the radial Bernstein polynomial over coordinate permutations changes every coefficient by less than \(3.31\times10^{-17}\). Thus a ray outside the comparison chamber adds at most \(6.62\times10^{-17}\) to the chamber estimate \(6.93\times10^{-8}+\sqrt3(4.19\times10^{-9})\), which still gives the stated global bound. An additional uniform dilation gives an exactly cubic representative with exactly the required filled volume.

Use the same two comparisons (2) with \(c=0.118\), \(\eta=7.4\times10^{-8}\) and \(\epsilon=21/2^{30}\). The checker now encloses the two \textbf{polynomials} \(h^2X^2[\epsilon\pm(\Delta F-1)]\) directly. Both are greater than \(4.4\times10^{-12}\) on the whole enlarged chamber. It also verifies \(\partial_X|_{V,W}F_X>0.3203125\) on the outer-face strip. With this derivative bound the two strict flux comparisons become \[
[0.3203125(1/2-c)-S]\eta>S\epsilon+(1+\eta)6.5\times10^{-11},
\qquad 0.0137\,c\eta/0.501>2.33\times10^{-10}.
\] The angular and radial enlargements fit inside the previously certified coordinate ranges. The same Neumann comparison therefore applies. On every ray the trial abscissa differs from the true one by less than \(4.0\times10^{-8}\). Multiplication by \(\sqrt3\) gives a radial error below \(6.93\times10^{-8}\) relative to \(g\). The previous global comparison of \(g\) with degree 48 yields the stated radial bound.

\paragraph*{Enclosing the physical measurements}

The 110 circles per ball use the five latitudes $0$, $0.12$, $0.27$, $0.47$ and $0.68$ for $p$ and \(q=s\sqrt{1-p^2}\) with 22 specified \(s\) values. Their maximum \(\eta\) is below \(0.964\). The code evaluates the sharper real upper and lower Poisson comparisons at 2,048 equally spaced azimuths for each circle. Integer polynomial coefficients, their truncation error, the floating-point evaluation error and the error in the coordinates are covered by the checked physical evaluator. The new dilation has smaller displacement and flux scaling than the original pair used to bound rounding errors.

For an enclosure of the azimuthal quadrature error subtract from \(F\) its affine Taylor polynomial at \(P_0\) and \(|X-P_0|^2/6\). The resulting harmonic function has norm at most \((M+2)R^2/2\). Indeed in the face test \(|F_{nn}|\le1+M\) and \(|F_{uu}|,|F_{zz}|\le M\). Integrate successively from the axis in the two transverse coordinates, using the two reflection symmetries to start each transverse derivative at zero. In a cavity-free ball the full Hessian satisfies \(0\le D^2F\le I\), giving the same bound with \(M=1\) in either orientation. The affine and quadratic terms incur zero quadrature error. Each homogeneous grade \(\ell\) on a measured circle is bounded by \((M+2)R^2\sqrt{2\ell+1}\,\eta^\ell/2\), by the Poisson formula and Legendre orthogonality. Its azimuthal degree is at most \(\ell\). The quadrature is exact through \(\ell=2045\) after either weight is applied. The remaining tail bounds both the mean and the sampled errors. All quoted observation intervals include that tail.

\paragraph*{Recovering the axial signs}

A linear combination of measurements with rational weights estimates the value of \(S\) or \(J\) at normalized axial coordinate \(t=-d/R\). If \(E\) is \(K_t\) minus the weighted sum of its measurement kernels, the error in this estimate is at most \(M\int_{-1}^{1}|E(\mu)|\,d\mu/2\). The checker chooses weights by a finite linear program and rounds them. Feasibility, optimality and discretization accuracy of this program are unnecessary for the inequality.

The checker encloses the residual Legendre coefficients through degree 900 by (1) and interval arithmetic. Beyond 900 use the triangle inequality and orthogonality, bounding the tail of an individual kernel of parameter \(\eta\) by \(\big(\sum_{\ell>900}(2\ell+1)\eta^{2\ell}\big)^{1/2}\). For grades 341-900 it bounds the norm of the \textbf{combined} residual by orthogonality. For the remaining polynomial write \(\mu=\cos\theta\) and interpolate \(\sin\theta\,E(\cos\theta)\) linearly at 4,097 equally spaced nodes on \([0,\pi]\). The integral of its absolute value is bounded by the trapezoidal sum of the absolute node values plus the interpolation error. For Legendre coefficient \(c_\ell\), the bound on the second angular derivative is \((\ell+1)^2|c_\ell|\) by the trigonometric Bernstein inequality. The cosine expansion of each \(P_\ell\) has nonnegative coefficients summing to one, giving a direct absolute rounding bound for the basis and sampled sums.

At depths spaced by \(0.004\), with a final step of \(0.002\), compute \(Y_{uu}=(S+J)/2\), \(Y_{zz}=(S-J)/2\) and \(Y_{nn}=1-S\). For \(S\) and \(J\) the axial second derivative at depth \(d\) has modulus at most \(2\sqrt5\,M/(R-d)^2\) by the degree-two Poisson coefficient in a ball of radius \(R-d\). Interpolation with this bound proves the displayed inequalities along the whole intervals.

\paragraph*{The unconditional cavity and direction tests}

If a tested ball contains no cavity, the exact fields satisfy \(0\le m_s\le1\) and \(|m_j|\le m_s\), since the Hessian is positive semidefinite with trace one there. For two lists of weights let \(G_s,G_j\) be their full measurement kernels. Their combined measured sum must be at most \[
\frac12\int_{-1}^1\max\{0,G_s+G_j,G_s-G_j\}\,d\mu .
\] The same coefficient and interpolation bounds apply to the maximum, by convexity of the maximum of affine functions. The measured lower bound exceeds this maximum on each tested ball.

\begin{center}\small\begin{tabular}{@{}lrrrr@{}}\toprule
Station & Radius & Measurement lower bound & Poisson upper bound & Gap lower bound\\\midrule
face & $0.043$ & $9.05\,10^{-8}$ & $1.37\,10^{-8}$ & $7.68\,10^{-8}$\\
body & $0.065$ & $3.34\,10^{-7}$ & $1.01\,10^{-7}$ & $2.33\,10^{-7}$\\
body tighter & $0.064$ & $8.81\,10^{-8}$ & $5.26\,10^{-8}$ & $3.56\,10^{-8}$\\
meridian 0.18 & $0.0435$ & $5.68\,10^{-8}$ & $1.43\,10^{-8}$ & $4.25\,10^{-8}$\\
meridian 0.24 & $0.044$ & $6.57\,10^{-8}$ & $1.28\,10^{-8}$ & $5.30\,10^{-8}$\\
meridian 0.38 & $0.046$ & $1.70\,10^{-7}$ & $2.47\,10^{-8}$ & $1.45\,10^{-7}$\\
meridian 0.53 & $0.049$ & $2.77\,10^{-7}$ & $2.99\,10^{-8}$ & $2.48\,10^{-7}$\\
meridian 0.64 & $0.051$ & $9.46\,10^{-8}$ & $1.83\,10^{-8}$ & $7.63\,10^{-8}$\\
meridian 0.72 & $0.054$ & $2.72\,10^{-7}$ & $3.05\,10^{-8}$ & $2.41\,10^{-7}$\\
meridian 0.80 & $0.0565$ & $2.00\,10^{-7}$ & $2.94\,10^{-8}$ & $1.70\,10^{-7}$\\
\bottomrule\end{tabular}\end{center}
Each gap is enclosed before rounding either of its terms. The signed
integer measurement weights for these rows are included with the verifier.

For a direction \((1,1,t)\) the ball center is \(a_0(1,1,t)\). The pairs \((a_0,\lambda)\) for the last two rows are \((0.34364957,0.03184)\) and \((0.33553579,0.09587)\). In the plane spanned by \(n=(1,1,0)/\sqrt2\) and the vertical \(e_z\), the sampling axis is \(\big((1-\lambda^2)n+2\lambda e_z\big)/(1+\lambda^2)\). The first tangent is \(u\), and the second is the perpendicular in that plane. The uniform error below \(4.0\times10^{-8}\) in the abscissa of the physical chamber graph puts each center within \(0.000004\) of its true radial surface point. The lower potential comparison exceeds \(3\times10^{-12}\) at every quadrature point used at these two locations. The cavity field, if the ball is cavity-free, is harmonic on the whole ball. This justifies the quadrature error even between exterior observations.

Additional measurements in directions \((1,1,0.18)\) and \((1,1,0.64)\) force cavities within \(0.043504\) and \(0.051004\) of the respective true outer points. Together they exclude the 45-site placement that had passed the sampled moment equations, for every choice of its weights. The corresponding distances at \(0.72\) and \(0.80\) are \(0.054004\) and \(0.056504\). At any finite count the \(0.18\) and \(0.80\) observations require distinct orbits.

For the direction conclusions, let \(c\) be an associated point, \(P\) a true surface point and \(a=|P-c|\). The convex interpolation inequality puts every point of its cavity at distance at least \(\big((1+\kappa)a-(1-\kappa)\operatorname{dist}(c,\partial O)^2/a\big)/2\) from \(P\). Sort the coordinates of \(c\) to give the closest signed permutation to \(P\). Bound \(\operatorname{dist}(c,\partial O)\) above by the least distance to the three axial, face-diagonal and body-diagonal supporting planes. The checker covers \(0\le z\le y\le x\le0.404\) by dyadic boxes. On each box it bounds \(a\) below by the distance of \(P\) to the box and uses the three lower coordinate endpoints in the support-plane bounds. At the face point, 5,433 boxes exclude \(z\ge0.30x\) for every cavity at distance below \(0.043001\). At the body point, 1,913 boxes exclude \(z\le0.50x\) below \(0.065001\). Both comparisons include the central label. The sorted component ratio is invariant under cubic symmetries and positive rescaling.

\paragraph*{Why the four points force 44 cavities}

Call a label orbit on a coordinate axis, a face diagonal or a body diagonal one of the three line families. For each of the two intermediate points the same distance inequality excludes all three line families, including the center, over their entire admissible radial intervals. Two joint exclusions apply to \textbf{all} orbits.

\begin{itemize}
\tightlist
\item
  One orbit cannot supply both a cavity within \(0.043001\) of the face point and one within \(0.049004\) of the \(t=0.53\) point.
\item
  One orbit cannot supply both a cavity within \(0.044004\) of the \(t=0.24\) point and one within \(0.065001\) of the body point.
\end{itemize}

For these joint tests, first take the sorted absolute label \(c=(x,y,z)\). It minimizes distance to both tested surface points among its signed permutations. Next project to \(c'=((x+y)/2,(x+y)/2,z)\). Both surface points lie on the projection plane, so their distances decrease. Also \(\operatorname{dist}(c',\partial O)\ge\operatorname{dist}(c,\partial O)\). Indeed \(O\) contains equal open balls about \(c\) and its reflection, so convexity puts the same ball about their midpoint in \(O\). Both changes weaken the necessary lower distance inequality. It therefore suffices to exclude all projected labels \((a,a,b)\) with \(0\le b\le a\le0.404\) from each pair of tests.

For each dyadic box in \((a,b)\), further bound \(\operatorname{dist}(c',\partial O)\) above by distance to a true radial boundary point. The checker selects a point in a direction \((1,i/4,j/4)\), \(0\le j\le i\le4\), and encloses its coordinates using the physical chamber graph. Selecting any of these points gives a valid upper bound. A box is excluded if the lower distance exceeds the permitted distance at either station. This finishes in 625 boxes for the first pair and 1,021 for the second. The separate line-family tests need at most 171 boxes each.

If two orbits lie outside the three line families, each has at least 24 points, giving at least 48 cavities. If just one such orbit exists, both intermediate points require it. The two joint exclusions prevent that orbit from serving either endpoint. Thus the face point requires a separate twelve-point orbit on face diagonals and the body point a separate eight-point orbit on body diagonals. This gives \(24+12+8=44\). Counts 46 and 47 are also excluded by the cubic orbit sizes and uniqueness of a central label.

With 44 or 45 cavities the remaining orbit has size 24. A sorted label with smallest component zero projects to a face-diagonal label in the preceding argument and hence fails the intermediate tests. A separate two-dimensional box check excludes \((a,b,b)\) at \(t=0.24\). The remaining 24-point form is \((a,a,b)\). Applying the two intermediate tests together excludes \(b\le0.27a\) and \(b\ge0.60a\). All permutations and signs of this form constitute one orbit. In either count the cavities needed at the endpoints have the matching face- and body-diagonal labels. Proximity excludes every different orientation in the same family, and averaging a nearby cavity point over the fixing symmetries puts a point of that cavity on the actual ray.

The covering theorem now excludes the 44-cavity configuration in this orbit classification. The 45-cavity configuration still passes the count test.

The exterior and axial enclosures are rational or interval inequalities in the accompanying verifier.

\paragraph*{A covering condition on every three-dimensional Bravais lattice}

Let \(O\) be the filled body centered at its obstacle site, \(K_j\) its open cavities, \(b=|O|\) the filled fraction and \(\kappa=f/(3b)\) for a genuine doubly coated sphere with middle-phase fraction \(f\). The shell is connected, its outer translates are separated, and the matrix is connected. Assume finitely many mutually separated Lipschitz cavities compactly contained in a bounded Lipschitz outer body. Equality gives convexity of \(O\) and the cavities and analytic interfaces. Let \(p_j=\kappa c_j\) be the constant gradient of the convex cavity field \(\Xi\) on \(K_j\). The associated points \(c_j\) lie in \(O\). Use the dual potential \(V=\Xi^*-|p|^2/2\) on \(D=\kappa O\), and put \(Q(p)=(1/\kappa-1)|p|^2/2\). Define \[
 d_j=\frac{2}{\kappa}\big[Q(p_j)-V(p_j)\big],
 \qquad \delta_j=\operatorname{dist}(c_j,\partial O).
\] Convexity of \(V\) and its prescribed tangent plane to \(Q\) at a nearest boundary point give \(d_j\le(1-\kappa)\delta_j^2\).

Take any symmetric \(H\succeq\kappa I\) with \(\operatorname{tr}H\le1\), and put \(M=\kappa H^{-1}\) and \(A=I-M\). For each \(t\in\overline O\) define the scores \(s_j(t)=(t-c_j)^TA(t-c_j)-d_j\). Then \[
 \min_j s_j(t) < \min_{P\in\partial O}(t-P)^TA(t-P),\qquad
 Mc_j+At\in K_j\quad\text{for every minimizing }j.                 \tag{1}
\] In particular \(A\overline O\subset\bigcup_j(K_j-Mc_j)\). If the winners throughout a measurable set \(T\subset\overline O\) are confined to indices \(J\), then \(\sum_{j\in J}|K_j|\ge\det(A)|T|\). This measures actual cavity volume. For a unique winning site the whole affine image of \(T\) belongs to its cavity.

To prove (1), put \(\beta=1/\kappa-1\), \(B=H^{-1}-I\), \(a=(\beta I-B)\kappa t\), and minimize \(U(p)=V(p)-p^TBp/2-a\cdot p\) on \(\overline D\). With \(F(B)=\det B-\operatorname{tr}B\) one has \[
 F(B)-2=\det(I+B)\,[1-\operatorname{tr}H]\ge0 .
\] The matrix \(B\) is positive definite and the product of its two smallest eigenvalues exceeds one. The latter follows from \(F(B)\ge2\) by isolating the largest eigenvalue in that inequality. The same strict pair-product condition holds for \(D^2V\) off the finite sites. Along the segment connecting the matrices, the square root of the determinant of the restriction to any two-plane is concave. Thus the average of \(\operatorname{cof}-I\) along the segment defines a locally uniformly elliptic linear operator \(L\) with \(LU=2-F(B)\le0\). The strong minimum principle excludes a regular interior minimum. Constancy on the connected punctured domain would eliminate all atoms.

At \(p=\kappa P\) on \(\partial D\), \(U\) equals \(\kappa(t-P)^TA(t-P)/2\) plus a constant independent of \(P\). For interior \(t\) with positive minimum boundary distance in this quadratic form, at a nearest boundary point its gradient is a positive multiple of the outward normal. Indeed \(A(P-t)\) is normal and has positive dot product with \(P-t\), which has positive outward component. Moving inward lowers \(U\). If the minimum boundary distance is zero, as always happens for \(t\in\partial O\), every boundary minimizer has \(A(P-t)=0\) and hence gradient zero. A global minimum there would contradict the Hopf lemma. The boundary is analytic and the operator is uniformly elliptic up to it. A site must therefore attain the strict interior minimum. At \(p_j\) the value of \(U\) equals \(\kappa s_j(t)/2\) plus the same constant. Its minimizing property gives \(Bp_j+a\in\partial V(p_j)\), or \(x_j=(I+B)p_j+a=Mc_j+At\in\overline K_j\).

In fact \(x_j\) is in the open cavity. Near \(x_j\) the slope \(\nabla\Xi\) lies inside \(D\) by continuity. The global quadratic lower support for \(V\) at \(p_j\), on adding \(|p|^2/2\) and conjugating, gives a local quadratic upper support for \(\Xi\) at \(x_j\), with Hessian \(H\). If \(x_j\in\partial K_j\), the difference of \(\Xi\) and this upper support is subharmonic in the adjoining shell, where its Laplacian is \(1-\operatorname{tr}H\ge0\). It has zero gradient at \(x_j\). The boundary maximum and the Hopf lemma force local equality in the shell. Its shell Hessian would then tend to \(H\) at \(x_j\) and also to the rank-one normal projection of the cavity. This is impossible. The volume claim follows by translating the cavities covering \(AT\).

\paragraph*{Directional coverage and a sharper bound on the site weights}

Take any unit vector \(n\) and use \(A=h_0 nn^T\), where \(h_0=(1-3\kappa)/(1-2\kappa)\). This gives \(\operatorname{tr}H=1\). At a supporting point of \(O\) with outer normal \(n\), (1) implies that some site has \(d_j>0\) and \[
 c_j\cdot n+\sqrt{d_j/h_0}>\max_{x\in\overline O}x\cdot n .
\] Consequently the convex hull of the balls with centers \(c_j\) and radii \(\sqrt{d_j/h_0}\), taken over positive weights, contains \(\overline O\) in its interior. In particular enlarging the inscribed balls \(B(c_j,\delta_j)\) by the common factor \(\sqrt{1+2\kappa^2/(1-3\kappa)}\) must suffice. This follows from the strict support inequalities in every direction and compactness of the sphere. At the other endpoint \(A=0\), all sites of largest \(d_j\) have positive weight and \(c_j\in K_j\).

There is a stronger bound than \((1-\kappa)\delta_j^2\). Set \(\delta=\delta_j\), \(R_+=\max_{P\in\partial O}|P-c_j|\), \(\rho=1-3\kappa\) and \(a_0=\rho^{1/3}(\delta+R_+)/2\). Use the radial comparison potential \(\varphi\) centered at \(p_j\) of mass \(4\pi a_0^3/3\) with slope \(a=\varphi'\) as in (4) of the dual proof. At a nearest boundary point \(q\) its tangent correction \(T\) gives \(\varphi+T\ge Q\) on all of \(\partial D\), with the prescribed full gradient at \(q\). Indeed the corrected radial derivative \(a(r)-\beta r\) vanishes at the midpoint \(\kappa(\delta+R_+)/2\), and the nearest supporting plane supplies the remaining affine term. For an attaining shell with at least two sites the maximum of \(V-\varphi-T\) occurs at \(p_j\), whatever the actual mass at this site. The proof is the same maximum principle and Hopf argument as for the centered cluster estimate. Evaluation at \(q\) gives, with \(r=\kappa\delta\), \[
 d_j\le (1-\kappa)\delta^2-
 \frac{2}{\kappa}\left[r\,a(r)-\int_0^r a(s)\,ds\right].
 \tag{4}
\] At the quarter-fraction data \(\kappa<\rho^{1/3}/2\). On \(0\le s\le r\) one has \(s\le a_0/2\) and \(a(s)\le a_0+s\) by direct substitution in \((a-2s)(a+s)^2=a_0^3\). Thus \(a'(s)=2s/(a(s)-s)\ge2s/a_0\). The bracket in (4) is \(\int_0^r s\,a'(s)\,ds\). This gives the explicit upper bound \[
 d_j\le (1-\kappa)\delta^2-
 \frac{8\kappa^2}{3\rho^{1/3}}\frac{\delta^3}{\delta+R_+}.               \tag{5}
\] The right side increases in both distances here. In particular the derivative of \(\delta^3/(\delta+R_+)\) in \(\delta\) is at most \(2\delta\), and \(8\kappa^2/(3\rho^{1/3})<1-\kappa\). Thus upper distance enclosures may be substituted directly.

\paragraph*{The axial endpoint and the count 44}

In simple-cubic periods put \(P=R_{100}e_1\), \(e=0.26\) and \[
 A=e(I-e_1e_1^T)+h e_1e_1^T,\qquad
 h=1-\frac{\kappa}{1-2\kappa/(1-e)}.
\] This corresponds to \(\operatorname{tr}H=1\). At \(t=P\) the right side of (1) is zero. A noncentral orbit allowed for a count below 48 has a representative \(c=(a,a,b')\), \(0\le b'\le a\), with \(b'=0\) for the face orbit and \(b'=a\) for the body orbit. A signed permutation with first coordinate \(a\) minimizes its squared distance to \(P\) in the \(A\) metric. Indeed \(h>e\) and \(a<R_{100}\). To win at the tip it must satisfy \[
 h(R_{100}-a)^2+e(a^2+b'^2)<(1-\kappa)\delta_j^2.                 \tag{2}
\]

The exterior observations independently force the same three orbits to supply cavities near specified points \(S\) of the outer surface. At such a point write \(l=|S-c|\) with the best signed permutation. If a cavity comes within distance \(w\) of \(S\), convex interpolation gives \[
 l\big[(1+\kappa)l-2w\big]\le(1-\kappa)\delta_j^2.
 \tag{3}
\] For the face orbit use \(S=R_{110}(1,1,0)/\sqrt2\), \(w=0.043001\). For the body orbit use \(S=R_{111}(1,1,1)/\sqrt3\), \(w=0.064001\). The tighter body test uses a measuring ball of radius \(0.064\). With the submitted weights the lower measurement is greater than \(8.81\times10^{-8}\) while the cavity-free upper bound is below \(5.2244\times10^{-8}\). The shift from its measuring center to \(S\) is below \(10^{-6}\). The intermediate orbit must pass (3) at both \((1,1,0.24)\) and \((1,1,0.53)\), with \(w=0.044004\) and \(0.049004\) respectively, and \(0.27a<b'<0.60a\).

No candidate in any of the three orbits passes (2) and all its required tests (3). The interval checker covers \(0\le a,b'\le0.404\). On each box it bounds \(\delta_j\) above by the three symmetry support planes and by distance to an enclosed true boundary point. Floating point selects one of 153 points to try. Every inequality is then evaluated with Arb intervals, and no nearest-point assumption is needed. The face, body and intermediate checks finish in 29, 59 and 125 boxes respectively. A 44-cavity shell has no other candidate to supply \(P\). In a 45-cavity shell only the center can supply it. By (1) at all six tips, the central cavity contains \(\pm hR_{100}e_i\). This first test gives an octahedron of volume \(4(hR_{100})^3/3>0.01750844\). The test over all positions below enlarges it.

\paragraph*{The axial orbit is distinct from all four measured cavity orbits}

For \textbf{every finite count}, any orbit supplying the axial endpoint in (1) at \(e=0.135\) must supply \textbf{none} of the four proximity tests. Several tests can share an orbit. In particular, any configuration in which every cavity orbit reaches at least one of the four measured sets of balls fails, regardless of its count.

Let \((x,y,z)\) be the sorted absolute site coordinates and retain \(h=1-\kappa/[1-2\kappa/(1-e)]\). The best signed permutation for the axial endpoint has cost \(h(R_{100}-x)^2+e(y^2+z^2)\). For each of the four measured boundary points \(S\), its displayed ordered coordinates also minimize \(l=|S-(x,y,z)|\) over signed permutations. If one orbit can perform both tasks, its weight must exceed the axial cost and must be at least \(l[(1+\kappa)l-2w]\) at one station. For the latter assertion use (3) before substituting the coarser upper bound on \(d_j\). A three-dimensional interval calculation splits \emph{all} sorted positions in \([0,0.404]^3\). On each box either the axial cost exceeds (5), or all four proximity demands exceed (5). Distances to the three support planes and the enclosed boundary samples give the upper bounds on \(\delta\) and \(R_+\) as above. The calculation finishes in 3,625 boxes through depth nine and includes the generic 48-point orbits.

The four tests already require at least two noncentral orbits. With exactly two, both are used by those tests, so the axial winner must be a central site. A 48-cavity configuration could only consist of two 24-point orbits. Indeed, with a single orbit outside the axial, face and body line families the four tests require two further orbits of sizes 12 and 8. Starting at \(24+12+8=44\) there is no way to add four using the cubic orbit sizes and at most one central site. If that single off-line orbit has size 48 the count already exceeds 48. With two off-line orbits the minimum is \(24+24=48\). This excludes 48 outright. It also classifies 49, since the unique possible odd orbit is the center. Its other two orbits must each have size 24. The low-ratio orbit must supply the face point and the \((1,1,0.24)\) point, and the high-ratio orbit the body point and the \((1,1,0.53)\) point, by the two incompatible pairs. In both the 45- and 49-cavity cases every noncentral orbit is used by the four tests. Thus only the center supplies the axial endpoint. Now \(hR_{100}>0.25927\) and the central octahedron has volume greater than \(0.0232396\).
\phantomsection\label{endfile:cavity-count-proof}

\SSection{stroke}{Finite-stroke recovery and the initial tangent}
Normalize reference area to one. In the solid the deformation is
area-preserving and its neo-Hookean energy per unit area is
$\mu_s(|F|^2-2)/2$. The proof first constructs a periodic sheet at
a fraction strictly below $\lambda^2$, then drives the energy cost of
the stroke to zero, and only then approaches the endpoint fraction.
\phantomsection\label{file:finite-stroke}
\paragraph*{Layers that carry the initial shear}

Fix a target fraction \(0<f<\lambda^2\). The reference stiffness comes from the sequential laminates used for sharp linear elastic bounds. The layer ordering below also gives them a finite motion.

Start with void. At stage \(j=1,2,3\), layer fresh solid with the preceding material, using a fraction \(p_j\) of fresh solid. Set

\[
 f_j=jf/3,\qquad f_0=0,\qquad
 p_j=\frac{f_j-f_{j-1}}{1-f_{j-1}},\qquad
 f_j=p_j+(1-p_j)f_{j-1}.
\]

The three \textbf{layer tangents} have angles \(0,+60^\circ,-60^\circ\). The outermost tangent is horizontal. Half of its core uses \(+60^\circ\) at stage 2 and \(-60^\circ\) at stage 1. The other half exchanges those two directions. The halves form layers parallel to the outermost sheets. The exchange preserves the entire initial tensor of the core.

To verify the elastic claim, first give the comparison solid finite bulk modulus \(K>0\) and shear modulus \(\mu\). Write \(C_s E=2\mu E_d+K(\operatorname{tr}E)I\), with \(E_d\) the trace-free part. Let \(C_j\) be the layer tensor at stage \(j\), \(C_0=0\), and \(D_j=C_s-C_j\). For unit layer normal \(n\), define \(\Gamma_n T=\operatorname{sym}(a\otimes n)\), where

\[
 a=(\mu I+K n\otimes n)^{-1}Tn .
\]

In a simple stack the core and fresh-solid strains \(E_c,E_s\) have a compatible jump \(\operatorname{sym}(a'\otimes n)\) and equal tractions. Consequently \(E_c-E_s=\Gamma_n D_{j-1}E_c\). Their average is \(E\), and \(D_j E=(1-p_j)D_{j-1}E_c\). Elimination gives

\[
 (1-f_j)D_j^{-1}=C_s^{-1}
      -\sum_{i=1}^j(f_i-f_{i-1})\Gamma_{n_i}.                    \tag{1}
\]

These linear solves exist at finite \(K\). In the \(C_s\) metric, \(\Gamma_n C_s\) is the orthogonal projection onto compatible layer jumps. Hence \(0\le\Gamma_n\le C_s^{-1}\), and the right side of (1) is positive definite. The strains and their stresses constructed in each split give the layer energy by equilibrium and compatibility. In particular, even a core of singular stiffness can be treated by prescribing its strain.

The equal increments \(f/3\) in (1) prove invariance under exchanging the inner directions. The unordered triple of normals is invariant, up to signs, under \(60^\circ\) rotation. Thus \(C_3\) is isotropic. For a trace-free \(T\), the averages over the three normals of \(|Tn|^2\) and \((n\cdot Tn)^2\) are \(|T|^2/2\) and \(|T|^2/4\). Inserting these in (1) gives the shear modulus

\[
 G(f,K)=\mu\,\frac{fK}{(2-f)K+2(1-f)\mu}
       \ \longrightarrow\ \mu\,\frac f{2-f}\quad(K\to\infty). \tag{2}
\]

On scalar stresses the same formula (1) gives bulk modulus \(fK\mu/((1-f)K+\mu)\). The stress construction below can include a hydrostatic test along with the two shear tests. In that case the initial bulk modulus of the first realization tends to \(\mu f/(1-f)\) at fixed \(f\). The matching upper estimate for an incompressible porous solid follows from the periodic Poisson dilation trial.

For any prescribed mean deviatoric stress \(T\), apply a root strain \(T/(2G(f,K))\). The constant stresses propagated down the tree have matching tractions at each split, vanish in the void, and have weighted complementary cost

\[
 \sum_{\text{solid leaves}}\text{weight}\ \sigma:C_s^{-1}\sigma
       =\frac{|T|^2}{2G(f,K)}.                                  \tag{3}
\]

Indeed the average of stress times strain at a split equals the product of their averages because the jumps are in equilibrium and compatible. Repeating this identity gives (3). Copies of the same tensor can be joined with the same mean strain and stress in each copy. This applies to the two orders in the core and to all the grains in Section 4. We use these stresses as tests in the \textbf{exactly incompressible} material.

For a general symmetric root stress the cost is \(T:C_3^{-1}T\). With \(K\ge\mu\), the bulk modulus above is at least \(G(f,K)\), so this cost is at most \(|T|^2/(2G(f,K))\). On one geometry the finite-cell construction bounds the excess over the ideal cost by \(\varepsilon|T|^2/\mu\) uniformly over symmetric \(T\), as shown in Section D of the finite-cell proof.

\paragraph*{One grain compresses in one direction}

We seek a common motion of the layers with grain mean \[
 D_0(t)=\operatorname{diag}(1,t^2),\qquad \lambda\le t\le1 .
\] The two columns below are taken in the reference horizontal and vertical directions.

\textbf{Outermost sheets.} Give the fresh solid at stage 3 gradient \(I\). Give the core gradient \(\operatorname{diag}(1,c)\), where \[
 c=\frac{t^2-p_3}{1-p_3}.
\] The tangential columns agree, the mean is \(D_0(t)\), and \(f_2<c\le1\).

\textbf{Two orders in the core.} Put \(q=\sqrt3/2\), \(v=(1/2,q)\), \(w=(1/2,-q)\), and set \[
 u=\sqrt{1-q^2c^2},\qquad a=(u-1/2)/q,\qquad
 B=\begin{pmatrix}1&a\\0&c\end{pmatrix}.
\] The two horizontal half-stacks receive \(B\) and its reflection conjugate \(\operatorname{diag}(1,-1) B\operatorname{diag}(1,-1)\). They share the first column and average to \(\operatorname{diag}(1,c)\). It suffices to describe the first half. Its stage-2 tangent is \(v\), and \(Bv=(u,qc)\) has length one. Let \(R\) be the rotation taking \(v\) to \(Bv\). Use \(R\) in the stage-2 fresh solid and \[
 H=(B-p_2R)/(1-p_2)
\] in its core. They share their \(v\) column. The determinant is affine on this rank-one line, so \(\det H=(c-p_2)/(1-p_2)>f_1\).

\textbf{Finest sheets.} The remaining tangent \(w\) satisfies \(|Hw|\le1\). To see this, \(Bw=(1-u,-qc)\), and direct scalar products give \[
 1-|Bw|^2=2u-1,\qquad 2(1-Bw\cdot Rw)=1+u-\tfrac32c .
\] The second expression is nonnegative. Since \(p_2=f/(3-f)\le1/2\), the inequality \(|Bw-p_2Rw|^2\le(1-p_2)^2\) follows from \(u\ge1-c/2\). Squaring that last inequality gives \(c(1-c)\ge0\).

A family of fine parallel sheets of fraction \(f_1\) admits a buckling construction whenever \[
 |Hw|\le1,\qquad \det H>f_1.                                    \tag{4}
\] In coordinates \(w,n=Jw\), with \(J\) a positive quarter turn, put \(h=Hw,\ k=Hn\). Choose two slopes \(\gamma_-\le0\le\gamma_+\) such that \(|h+\gamma_\pm k|=1\). Alternate them in proportions of zero mean. The centerlines then have unit tangents \(p=h+\gamma_\pm k\) and mean tangent \(h\). Neighboring centerlines shift by the sheet pitch times \(k\). Their normal separation per unit pitch is \(\det(p,k)=\det H\), strictly greater than the reference sheet thickness per unit pitch \(f_1\). The smoothed turns occupy arbitrarily little area. The detailed construction restores exact incompressibility through them and at the joins. This buckling step allows shortening along \(w\) while the straight reference sheets carry axial stress.

Multiplying that clearance by the fraction occupied by stage 1 gives the entire unused area in the ideal grain, \[
 (1-p_3)(1-p_2)(\det H-f_1)=t^2-f .
\] The ideal gradient in the space between the sheets consequently has determinant \((t^2-f)/(1-f)\). Its minimum tends to zero as \(f\uparrow\lambda^2\), so the finite-cell tolerances must be chosen after the target fraction.

For a 5\% stroke, for example, choose \(f=0.9\). The three fresh-solid fractions from finest to coarsest are \(3/10,3/7,3/4\). At \(t=0.95\), \(c=0.61\) and the available normal separation per finest-sheet pitch is \(\det H=0.3175\), for sheets of thickness \(0.30\) per pitch. This choice gives limiting shear \(9\mu/11\). Increasing \(f\) toward \(0.9025\) closes the remaining gap to \(361\mu/439\).

\paragraph*{Turning the grain motion into an isotropic stroke}

There is a way to make the mean as close to isotropic as desired without changing any layer proportions as \(t\) varies. Suppose the available grain has mean \(D=\operatorname{diag}(a,b)\), with \(a\ge b>0\), \(ab=t^2\). Take two copies with their axes at \(+45^\circ\) and \(-45^\circ\) to a common horizontal interface. Rotate the deformed copies so that both images of the interface tangent are horizontal and point to the right. Their common length is \(a'=\sqrt{(a^2+b^2)/2}\). The second columns have opposite horizontal components by reflection, and their vertical component is \(b'=ab/a'\) by the determinant. Thus an equal stack has mean \(\operatorname{diag}(a',b')\). The rotations can be chosen continuously and equal to the identity as superposed rotations at \(t=1\). In particular, the initial map of each reoriented copy is the identity.

Start this operation with \(D_0(t)\) and repeat \(N\) times. All interface directions and all reference fractions are fixed throughout the motion. The ratio \(r=a/b\) follows \[
 r_0=t^{-2},\qquad r_{j+1}=\tfrac12(r_j+r_j^{-1}),\qquad
 r_{j+1}-1=\frac{(r_j-1)^2}{2r_j}\le\tfrac12(r_j-1).
\] Hence \(D_N(t)\to tI\) uniformly for \(\lambda\le t\le1\). Finally multiply the entire deformed grain by \(L_N(t)=tD_N(t)^{-1}\). This gives \textbf{exact mean \(tI\)} and preserves every local determinant. Its distance from \(I\) tends uniformly to zero. For any solid gradient of determinant one, the energy change is bounded above by \(\tfrac{\mu}{2}\|L_N^{\mathsf T}L_N-I\|_{\rm op}|F|^2\). The integrated penalty therefore vanishes if the uncorrected solid energy is bounded. At \(t=1\), \(L_N=I\).

\paragraph*{A periodic sequence and the order of limits}

Construct the finite tree on an equilateral triangle with affine deformation \(D_N(t)x\) on its whole boundary and complementary stress of mean \(T\) with boundary traction \(T\nu\). The recovery proof gives the following simultaneous estimates, at fixed \(f,K,N\) and arbitrarily small error \(\varepsilon\):

\begin{itemize}
\item
  The solid fraction differs from \(f\) by at most \(\varepsilon\).
\item
  The maximum solid energy per reference area is at most \(\varepsilon\mu\) before the last multiplication by \(L_N\).
\item
  Each symmetric mean stress \(T\) has a statically admissible test of complementary cost at most \(T:C_3^{-1}T+\varepsilon |T|^2/\mu\). For a deviatoric \(T\), the first term is \(|T|^2/(2G(f,K))\).
\end{itemize}

The estimates use one perforated solid and the same loading path for all \(t\).

After applying \(L_N\), the triangle has boundary deformation \(tI\). Assemble six rotated copies into a regular hexagon and repeat on the triangular lattice of hexagon centers. Use the translated affine boundary maps on every copy. Adjacent deformations agree, and each triangle maps injectively onto its scaled triangle. The stress in each copy uses the correspondingly rotated input \(T\), so tractions match at the edges. The initial linear tensor is invariant under \(60^\circ\) rotation. On symmetric strains this rotation fixes the scalar part and rotates the two-dimensional deviatoric part by \(120^\circ\). Thus the two parts decouple and the shear modulus is exactly isotropic at each finite stage.

For every affine-periodic infinitesimal displacement of strain increment \(E\), equilibrium gives \(\langle\sigma:e(u)\rangle=T:E\). Exact solid incompressibility makes only \(\sigma_d\) enter the work in the solid. The complementary test cost in that solid is \(\langle|\sigma_d|^2/(2\mu)\rangle\le\langle\sigma:C_s^{-1}\sigma\rangle\). Cauchy-Schwarz, with \(T=E\) trace-free, yields \[
 \frac{1}{2\mu_*}\le\frac{1}{2G(f,K)}+\frac{\varepsilon}{\mu}.
\] This argument uses \(K\) only to generate equilibrated test fields. On each connected finite cell, the same small-strain energy is the quadratic coefficient of the prescribed nonlinear energy.

For example, at index \(m\ge2\) take \(f=\lambda^2(1-1/m)\), \(K=m\mu\), \(N=m\) and finite-cell error \(\varepsilon=1/m\). Equation (2) gives the claimed limiting shear. Section 4 gives vanishing peak energy with exact \(tI\) throughout. Reversing the parameter traverses the same path back to the identity. The construction uses arbitrarily separated length scales and sets no positive lower bound on feature size.
\phantomsection\label{endfile:finite-stroke}

\subsection{Recovery by bonded, injective cells}\label{route:recovery-by-bonded-injective-cells}
\phantomsection\label{file:stroke-recovery}
\paragraph*{Ordinary layer splits}

At a node let \(A(t)\) be the parent gradient and let \(\sigma\) be its constant trial stress. There are two child gradients \[
 A_i(t)=A(t)+s_i\,b(t)\otimes n ,
\] with fixed unit normal \(n\), fixed \(s_i\) of weighted mean zero, and continuous \(b(t)\). The child stresses are \(\sigma_i=\sigma+d_i v\otimes v\), where \(v=Jn\) and the \(d_i\) have the same weighted mean zero. This form follows from equality of the normal traction. Zero jumps are permitted. In each ordinary split of the solution the determinants at both ends stay strictly positive, so the whole rank-one segment has a uniform positive determinant. Reorienting a copy changes its reference normal by a fixed orthogonal transformation. A superposed output rotation depending on \(t\) changes only the gradient data.

Work on a convex parent polygon \(P\) where the entering deformation is affine and the entering stress constant. Take a convex inner polygon \(P'\) leaving a boundary collar of as small a relative area as desired. Choose a smooth cutoff \(\chi\) supported in \(P\) and identically one on \(P'\). Let \(\psi\) be a bounded continuous periodic sawtooth of period one, with derivatives \(s_i\) in intervals of the prescribed proportions. Let \(\Psi\) be a periodic \(C^{1,1}\) function with \(\Psi''=d_i\) on the same intervals. Both exist by the zero mean conditions. Add to the entering deformation and Airy potential, respectively, \[
 l\,b(t)\chi(x)\psi(n\cdot x/l), \qquad
 l^2\chi(x)\Psi(n\cdot x/l).
\] In \(P'\) the fields equal the child gradients and stresses on the designated stripes. Outside \(P'\), the deformation gradient belongs to the rank-one segment up to an error \(O(l\|\nabla\chi\|_\infty)\), with a bounded coefficient. The stress is bounded by the entering and child amplitudes times a constant, taking \(l\|\nabla\chi\|_\infty\) and \(l^2\|D^2\chi\|_\infty\) small. In particular the deformation has a uniform positive Jacobian throughout, and its gradient and the stress stay bounded independently of the collar area and wavelength. Constants in this appendix may depend on the fixed tree and the strict area margins.

On either side of each fixed stripe boundary the deformation is smooth up to the sides, with spatial derivative uniformly continuous there for all \(t\). Keep the smooth formulas until the final subdivision is known. The children of \(P'\) are convex stripe intersections with exactly affine deformations and can be used immediately as the next parents. The Airy potential needs no interpolation.

Stop the descendants in the collar and mark it for the porous filling in Section C. Proceed down the finite tree in the interiors. A sufficiently short wavelength realizes the prescribed proportions to any accuracy on each polygon, by periodic averaging. Choose wavelengths after fixing their parent polygons. Thus the union stopped in collars can have arbitrarily small total area, the leaf fractions approximate the ideal weights, and all gradient norms, inverse determinants and stress norms are bounded independently of this small area. There are only finitely many levels at fixed \(N\). Apply the stress bounds to an orthonormal basis of the three-dimensional space of symmetric matrices. Linearity and Cauchy-Schwarz then bound each stress amplitude by \(C|T|\) on this same geometry for every symmetric mean stress \(T\).

\paragraph*{Buckling in the last family of sheets}

Consider a stage-1 parent polygon with tangent \(w\), normal \(n=Jw\), and fraction \(p=f_1\) of solid sheets. Its entering gradients \(H(t)\) obey \[
 h=Hw,\qquad k=Hn,\qquad |h|\le1,\qquad d=\det(h,k)>p .
\] There is a uniform strict margin. We use local coordinates \(x=sw+zn\) and suppress the entering translation. In this parent the static laminate has constant axial stress in the solid and zero stress in the void, with the prescribed average. Its stress can be constructed by precisely the Airy formula in Section A with stripes in the \(z\) variable.

The quadratic equation \(|h+\gamma k|=1\) has ordered roots \(\gamma_-\le0\le\gamma_+\). On each period take the lower root on the first interval of length \(\gamma_+/(\gamma_+-\gamma_-)\) and the upper root on the remainder when the roots differ. The step has average zero. If a root is zero use the zero function in the one-root-of-positive-weight case. If the roots coincide both are zero. Choose \(0<\omega<1/8\) and fix a smooth nonnegative unit-integral kernel \(\rho\) supported in \((-1,1)\). Periodic convolution of the step with \(\omega^{-1}\rho(\cdot/\omega)\) gives \(\gamma(\xi,t)\). It is smooth in \(\xi\), uniformly bounded, and has average zero. It varies continuously with \(t\), including at a vanishing plateau or a double root. For example, the original step is continuous as an \(L^1\) function of \(t\), and convolution gives continuity of all spatial derivatives. Let \(\Gamma\) be its mean-zero periodic primitive. For \[
 p_s=h+\gamma(s/L,t)k,\qquad b_s=Jp_s/|p_s|^2-k
\] we have \(\det(p_s,k)=d\) and \(|p_s|\ge d/|k|>0\). Outside two intervals of total relative length at most \(4\omega\) in each period, \(p_s\) is constant and has unit norm. The symbol \(p_s\) here denotes a vector and \(p\) the scalar sheet fraction.

Choose smooth cutoffs \(\chi_c,\chi_f\), compactly supported in the parent and between zero and one. Take \(\chi_c=1\) near the support of \(\chi_f\), and both equal to one on a fixed inner polygon occupying almost all the parent. Choose a period-one continuous sawtooth \(\psi\) with derivative 1 on fraction \(p\) and \(-p/(1-p)\) elsewhere. Use the deformation \[
 Y=Hx+L\chi_c k\Gamma(s/L,t)+l\chi_f b_s\psi(z/l).
                                                               \tag{5}
\] The reference sheets are the fixed stripes where \(\psi'=1\). The bending wavelength \(L\) is much smaller than the cutoff widths, and the sheet pitch \(l\) is much smaller again. Moving bends change the deformation within the fixed sheets.

For clarity, the determinant and energy controls for (5) are as follows.

\begin{itemize}
\item
  Where the fine term is absent, the principal gradient has columns \((h+\chi_c\gamma k,k)\) and determinant \(d\). The remaining gradient is \(O(L\|\nabla\chi_c\|_\infty)\).
\item
  On the support of the fine term the principal columns are \((p_s,k+\chi_f b_s\psi')\). Their determinant is \(d+\chi_f(1-d)\psi'\). In the solid stripe it interpolates between \(d\) and 1. In the void stripe it interpolates between \(d\) and \((d-p)/(1-p)\). All are uniformly positive. The gradient error is bounded by \(C l(\|\nabla\chi_f\|_\infty+C_\omega/L)\). The bounds on the principal columns are independent of \(\omega,L,l\).
\item
  In the inner polygon the determinant on a solid stripe is uniformly as close to one as desired. Outside the smoothed bends its gradient is a rotation, with columns \((p_s,Jp_s)\). Inside the bends the gradients stay bounded. For every \(t\), the area fraction covered by bends can be made arbitrarily small by choosing \(\omega\) small and then \(L\) small on the given polygon. For instance, enclose it in a fixed rectangle in \(s,z\) coordinates and count complete periods and the two end periods. All these estimates are uniform in \(t\).
\end{itemize}

Localize the axial laminate stress on the same \(z\)-stripes with a cutoff equal to one on the inner polygon. This gives its two exact constant stress values there and a bounded stress in the collar. At \(t=1\) the coarse wave and the fine displacement correction both vanish. Label the inner polygon according to the fixed solid and void stripes, and mark the collar as exceptional. Within each fixed stripe the spatial derivatives are uniformly continuous for fixed cutoffs and scales even though the bends move with \(t\).

\paragraph*{Exact incompressibility and stress transmission through the joins}

Choose one conforming triangulation of the whole final polygonal subdivision, respecting every stripe edge, inner polygon and shared edge subdivision. Refine it everywhere by midpoint subdivision into similar triangles. Nodal interpolation converges in the piecewise gradient norm uniformly in \(t\), by uniform continuity on each closed piece and the fixed aspect ratios of this finite initial mesh. It preserves affine fields on leaves, the affine boundary values, continuity across all interfaces and the identity at \(t=1\).

In particular, resolve the turns as well as the sheet pitch. On the inner stripes a fixed rescaled smoothing kernel gives \(|\partial_s p_s|+|\partial_s b_s|\le C/(\omega L)\) and second derivatives bounded by \(C/(\omega L)^2\). Once \(l\le\omega L\), the piecewise second derivatives of (5) there are bounded by \(C/(\omega L)\). For the hard-hole area fraction \(\beta\) used below, put \(m_0=\min\{\beta^2,(\lambda^2-f)/(1-f)\}\). It suffices to make the fine-pitch error \(Cl/(\omega L)\) and then the interpolation error \(C_{\rm mesh}h_{\rm mesh}/(\omega L)\) each smaller than \(m_0/C_0\). Here \(h_{\rm mesh}\) is the largest triangle diameter, \(C_0\) can be chosen from the bounds on the principal gradients and the determinant perturbation, and \(C_{\rm mesh}\) can depend on the chosen initial triangulation. In the collars, choose the wavelengths to make the cutoff errors of Sections A-B smaller than \(m_0/C_0\), then refine the mesh so that the interpolation errors there are also below \(m_0/C_0\).

The resulting continuous piecewise affine map has positive determinant bounded below and bounded gradients. On ordinary solid leaves the gradient is a rotation. On all solid-labelled sheet triangles the determinant is within the prescribed tolerance of one \textbf{for all \(t\)}, including times when a bend passes over them. At points outside the bends on those triangles the interpolated gradient is uniformly close to a rotation. On void-labelled triangles the stress is zero. The exceptional triangles occupy arbitrarily small area and can carry bounded stress and deformation. The following local operation uses the same reference hole at every \(t\).

Let \(P\) be a convex triangle or square with centroid \(x_0\). For \(x=x_0+z\), let \(r(z)\) be the gauge of \(P-x_0\), so \(r=1\) on its boundary and \(z\cdot\nabla r=r\) almost everywhere. Choose a shell of solid area fraction \(\alpha\) between \(r_0=\sqrt{1-\alpha}\) and 1, and leave the interior \(r<r_0\) void. For an entering affine map \(b+Az\) with \(D=\det A>\alpha\), replace it by \[
 b+A\,\frac{g(r)}r z,\qquad
 g(r)=
 \begin{cases}
 \sqrt{1+(r^2-1)/D},&r_0\le r\le1,\\
 (r/r_0)g(r_0),&0\le r<r_0 .
 \end{cases}                                                    \tag{6}
\] The value at zero is \(b\). The map agrees with the incoming affine map on the whole boundary. The strictly increasing radial coordinate makes it a homeomorphism onto \(b+A(P-x_0)\). Its shell determinant is \(Dgg'/r=1\). This follows by differentiating the radial map and using \(z\cdot\nabla r=r\), irrespective of the polygon shape. When the entering map is the identity so is (6).

Stress can be passed through the shell. If the entering Airy potential is \(U\), subtract its affine Taylor polynomial at \(x_0\), multiply the remainder by a smooth cutoff zero on a neighborhood of the hole and one near \(\partial P\), and add back the polynomial. The potential and its gradient match near \(\partial P\), while the stress in the hole vanishes. If \(M=\|D^2U\|_\infty\), the Taylor remainder and its gradient are bounded by \(M|z|^2/2\) and \(M|z|\). This provides a uniform bound for the modified stress on similar copies of any fixed shape with fixed \(\alpha\). In a square with a centered hole of area fraction \(\beta\ll1\), take the cutoff in a neighborhood of diameter \(O(\sqrt\beta)\) times the side. Its first and second derivatives scale inversely with this diameter and its square. Thus the stress changes on area fraction \(O(\beta)\) and stays bounded by \(CM\) with an absolute \(C\).

Every shell uses the single neo-Hookean solid of shear modulus \(\mu\). Use the following area fractions in squares.

\begin{enumerate}
\def\labelenumi{\arabic{enumi}.}
\tightlist
\item
  In solid-labelled material choose \(\alpha=1-\beta\) with small \(\beta>0\). On exactly rigid leaves one may instead retain full solid. Make the determinant error from Section B and interpolation at most \(\beta^2\). For \(r\ge\sqrt\beta\) one has \(g(r)/r=1+O(\beta)\) and \(g'(r)=1+O(\beta)\), by (6). For a square \(|z||\nabla r|/r\) is bounded by an absolute constant. Thus the shell gradient is within \(O(\beta)|A|\) of \(A\). It is uniformly close to a rotation outside the bends and bounded everywhere. The change in integrated complementary cost is \(O(\beta)|T|^2/\mu\), with a fixed-tree constant.
\item
  In void-labelled material choose a small positive \(\alpha\) less than half the uniform determinant lower bound. The shell gradient is bounded uniformly even as \(\alpha\downarrow0\), since \(r_0\) tends to one and \(D\) is bounded away from zero. Its integrated energy tends to zero with \(\alpha\). Its stress stays zero.
\item
  In exceptional material use a fixed \(0<\alpha<\min(1,\inf D)/2\). The square shell then carries bounded energy and complementary cost per unit reference area, uniformly in \(t,T\) for \(|T|=1\). Its total contribution tends to zero with the exceptional area. This shell supplies a bonded connection across each join.
\end{enumerate}

There is a finite way to obtain those square estimates on arbitrarily shaped triangulations. Drop an altitude to the longest side of each triangle. In each resulting right triangle inscribe a square at the right angle with its far vertex on the hypotenuse. The remainder consists of two right triangles similar to the starting one. Repeat finitely many times. Each generation removes the same positive proportion of the residual area at that initial shape. Fill the squares as above and use (6) and the stress cutoff on the final residual triangles as well. At every fixed initial triangle shape, the energy and stress bounds there are finite uniformly in time and scale, since the determinant margins are strict. There are finitely many initial triangles at each realization. Choose the number of repetitions for each to make the sum of these residual contributions arbitrarily small, even when some initial triangles are very thin. The unfilled holes are strictly interior to their respective squares or triangles. Thus the finished cell can have a connected solid with Lipschitz boundary and finitely many separated polygonal holes.

\paragraph*{Verification on the full cell}

The piecewise affine map before (6) has positive determinant on every triangle and affine boundary data of positive determinant. It is a homeomorphism onto the affine image of the parent domain. At points off the finite union of edge images every preimage has sign \(+1\), and the degree is one inside the boundary and zero outside. The map is open at interior edges and vertices. At a vertex, the oriented nondegenerate triangles in its fan wind at least once around the vertex image. Two interior preimages or an interior preimage of a boundary point would therefore contradict the count on open sets of regular points. The boundary map itself is one-to-one. Replacing the restriction on any filling square or triangle by (6) preserves this homeomorphism property. All finite-geometry maps vary continuously with \(t\), including at the identity, and have unit Jacobian almost everywhere in the solid.

The Airy corrections give equilibrium on the full cell, zero stress on every pore, and the original boundary traction. In fact the finished polygon has a solid neighborhood of its boundary with stress identically \(T\). The original oscillations are compactly supported in the polygon, and each of the finitely many pore corrections is supported strictly inside a filling piece. The mean stress is exactly \(T\), by integration of the boundary traction against the coordinate functions. In the regular parts the leaf stresses give the ideal complementary cost \(T:C_3^{-1}T\) when counted with their ideal weights. Each leaf stress and each corrected stress in a square is bounded by \(C|T|\) at the fixed tree and the chosen exceptional-shell fraction. Thus the errors in leaf weights, the altered area in hard squares, and the exceptional area give cost errors bounded by their respective fractions times \(C|T|^2/\mu\). In the residual similar triangles the same estimate holds with shape-dependent constants, absorbed by choosing their total areas small in Section C. Choose the geometric tolerances so the sum of the error coefficients is at most \(\varepsilon\) on the unit sphere of symmetric stresses. The excess over the ideal cost is then at most \(\varepsilon|T|^2/\mu\) for every symmetric \(T\). Equation (3) of the solution gives the value of that ideal cost for trace-free \(T\).

To see the parameter order explicitly, fix the desired errors and the finite ideal tree. Take the hard-hole fraction \(\beta\), the soft-shell fraction and the allowable exceptional and bend areas small. Choose collars, parent wavelengths and then child wavelengths, successively on the finitely many polygons. In Section B choose the smoothing fraction and coarse wavelength before the finer sheet pitch. Take this pitch and the interpolation error sufficiently small for both the hard-shell determinant tolerance \(\beta^2\) and the minimum void determinant \((\lambda^2-f)/(1-f)\), as in Section C. The uniform bounds in the collars and bends depend only on the bounded principal gradients, their determinant margins and the stress amplitudes once the gradient errors are small. Finish with the finite square fillings. The solid fraction approaches the weighted fraction \(f\), the maximal integrated energy tends to zero, and the complementary costs are bounded above by the ideal costs up to the prescribed error. All three estimates hold on the same reference set. The final determinant-one correction, hexagonal assembly and diagonal sequence are given in the solution.

\paragraph*{Taking the tangent of the nonlinear cell}

For these fixed, connected, finitely perforated periodic solids, the linear cell problem used in the solution also gives the quadratic coefficient of the specified nonlinear energy. Write \(G\) for the solid and take a symmetric increment \(E\), with mean gradient \(\exp(sE)\).

For the lower estimate consider a sequence with cell energy \(O(s^2)\). On \(\det F=1\), \(|F|^2-2\ge\operatorname{dist}^2(F,SO(2))\), as follows from the two singular values. On a connected Lipschitz patch of several whole hexagons with their holes removed, geometric rigidity and Poincaré's inequality give a single rigid map \(R_s x+b_s\) within \(O(|s|)\) of the deformation in \(H^1\). Comparing a solid ball with its translates by two independent lattice periods gives \(R_s=\exp(sE)+O(|s|)=I+O(|s|)\). The constants here are for one fixed design.

Thus \(u_s=(y_s-x)/s\), modulo translations, has a weakly convergent \(H^1\) subsequence on the solid. Its limit has affine increment \(E\), and its periodic fluctuation extends through the Lipschitz pores. On \(G\), \(\operatorname{div}u_s=-s\det\nabla u_s\to0\) in \(L^1\). Consequently the limit \(u=Ex+w\) has \(\operatorname{div}u=0\) on \(G\). It is admissible for \(Q_G(E)\), and there \(|e(u)_d|=|e(u)|\). The identity \(|F|^2-2=2|(\operatorname{sym}F)_d|^2\) at determinant one and weak lower semicontinuity now bound the energy divided by \(s^2\) below in the limit by \(Q(E)/2\).

The upper estimate works for any closed periodic competitor with open pores of positive area. Take an admissible periodic fluctuation \(w\) in \(Q(E)\) on the torus, put \(e_0=\operatorname{tr}E\), and write \(g=\operatorname{div}w\). Then \(g+e_0=0\) on \(G\). Approximate \(g+e_0\) in \(L^2\) by smooth functions compactly supported in the pores with mean \(e_0\), adjusting the integral with a fixed smooth pore bump. Subtract \(e_0\) from each to obtain mean-zero \(g_j\to g\). Write \(w=w^0+\nabla\Delta^{-1}g\) using the periodic mean-zero Poisson inverse. Smoothing the globally solenoidal \(w^0\) and using \(g_j\) in the Poisson term gives smooth \(w_j\to w\) in \(H^1\), each with divergence \(-e_0\) on a neighborhood of \(G\). For fixed \(j\) and sufficiently small \(|s|\), its flow \(\Phi_s\) keeps particles starting in \(G\) inside that neighborhood, hence has determinant \(\exp(-s e_0)\) on \(G\). The diffeomorphism \(\exp(sE)\Phi_s\) has unit determinant there, the prescribed lattice increment and first variation \(Ex+w_j\). Its energy divided by \(s^2\) tends to \(\mu\langle\mathbf1_G|E+e(w_j)|^2\rangle\). Taking \(j\to\infty\) and minimizing proves the upper estimate, including for \(E=I\).

The lower estimate also holds uniformly over enlarged periodic cells for a fixed design. On each equilateral sector, its solid is a connected Lipschitz domain containing an inner neighborhood of the three edges. For the cell in Section 6 of the solution, the earlier frame has a solid boundary band. Each point in its patterned interior connects to that band by following the connected square network until it reaches the band. The finitely many new holes are strictly inside individual filling pieces. Each filling piece remains connected with a full rim. The finitely many pieces abutting a sector edge keep their rims, so the filled sector still has a solid neighborhood of all three edges. The new pores stay a positive distance from the old pore boundaries. This also gives a thin polygonal solid collar along every old pore boundary.

Apply rigidity on the perforated sector, subtract the resulting rigid motion with a translation chosen by Poincaré's inequality, and use the bounded Sobolev extension operator of that Lipschitz domain to fill the entire sector. Add back the rigid motion. This controls the integral of squared distance to rotations on the filled triangle by a fixed multiple of that integral on its solid, including for sectors with the frame's original pores. Choose the fits and extensions on one set of sector representatives for the superperiod, and transport them by the prescribed affine-periodic increments. They agree with the original map along every common edge and therefore glue across the tiling.

The period parallelogram can cut across sectors. Both gradients and the solid indicator are superperiodic, so each integral in the sector estimate has the same value on the parallelogram as on a complete set of sector representatives. Applying rigidity on the whole filled \(k\)-by-\(k\) period parallelogram has a constant independent of \(k\) by rescaling. Its mean gradient \(\exp(sE)\) fixes the approximating rotation to within \(O(|s|)\) of \(I\). Thus the mean squared norm of \(\nabla u_s\) is bounded independently of \(k\). Averaging the periodic fluctuation of \(u_s\) over the \(k^2\) translates gives a one-cell fluctuation. On the solid the resulting displacement has divergence tending to zero in \(L^1\) by the determinant identity above. Jensen's inequality bounds the squared norm of its deviatoric symmetric gradient by the corresponding averaged energy. Passing to a weak limit gives the same \(Q(E)/2\) lower bound.
\phantomsection\label{endfile:stroke-recovery}

\subsection{Recovery with initial bulk tending to zero}\label{route:recovery-with-initial-bulk-tending-to-zero}
The linear frame construction [P] below supplies the dilation and stress estimates used in the finite-stroke recovery.
\subsection{The dilation trial and the complementary frame estimate}\label{detail:frame-tests}\label{file:frame-tests}
\paragraph*{Setting}

\(\mathrm{Sym}\) is the space of symmetric \(2\times2\) matrices, \(A:B=\operatorname{tr}(AB)\), \(|A|^2=A:A\), \(E_d=E-\tfrac12(\operatorname{tr}E)I\), and \(e(u)=\tfrac12(\nabla u+\nabla u^{\mathsf T})\).

A \emph{periodic mixture} is a closed set \(G\subset\mathbb R^2\) (the solid) together with a lattice \(\Lambda\) of periods of \(G\), with period cell \(D\). Write \(\chi=\mathbf 1_G\), \(\langle\cdot\rangle\) for the average over \(D\), and \(f=\langle\chi\rangle\) for the solid fraction. The solid is exactly incompressible with shear modulus \(\mu\); the complement is void. For \(E\in\mathrm{Sym}\) the (doubled) effective energy is

\[
Q_G(E)=\inf\Big\{\,2\mu\,\big\langle\chi\,|e(u)|^2\big\rangle:\ u=Ex+w,\ w\in H^1_{\rm per}(D;\mathbb R^2),\ \operatorname{div}u=0\ \text{a.e. in }G\Big\}\in[0,\infty]. \tag{1.1}
\]

The displacement extends to the full cell. The pores carry zero energy. All frames below have Lipschitz solid interiors, as proved in Section~\ref{detail:frame-fields}, so Sobolev extension gives the same infimum for fields defined on the solid.

\textbf{Lemma 1.1.} 1. If the void has positive area, \(Q_G\) is finite and is a quadratic form, \(Q_G(E)=E:C_*E\) with \(C_*\) symmetric and positive semidefinite. Moreover \(Q_G(I)\le4\mu f/(1-f)\). 2. If the void has zero area, \(Q_G(E)=\infty\) whenever \(\operatorname{tr}E\ne0\). 3. If a rotation \(R\) maps \(G\) and \(\Lambda\) to themselves, \(Q_G(RER^{\mathsf T})=Q_G(E)\).

\emph{Proof.} (1) Let \(\rho=-2\) on \(G\) and \(\rho=2f/(1-f)\) on the void, so \(\langle\rho\rangle=0\), and let \(\phi\in H^2_{\rm per}\) solve \(\Delta\phi=\rho\). Then \(w_E=\tfrac12(\operatorname{tr}E)\nabla\phi\) is linear in \(E\) and \(Ex+w_E\) is admissible. All admissible fields for \(E\) are \(Ex+w_E+n\) with \(n\) in the linear space \(\mathcal N=\{n\in H^1_{\rm per}:\operatorname{div}n=0\text{ a.e. in }G\}\). Hence \(Q_G(E)\) is the squared distance, in the Hilbert space \(L^2(G\cap D;\mathrm{Sym})\) with norm \(2\mu\langle\chi|\cdot|^2\rangle\), from \(E+e(w_E)\), which is linear in \(E\), to the linear subspace \(\{e(n)|_G:n\in\mathcal N\}\). This is the squared norm of the orthogonal projection of that vector onto the orthogonal complement of the closure of the subspace, hence a finite quadratic form. For \(E=I\) the field \(x+\nabla\phi\) has strain \(I+\nabla^2\phi\), which on \(G\) equals the deviatoric part of \(\nabla^2\phi\) because its trace \(2+\rho\) vanishes there; and \(\langle|(\nabla^2\phi)_d|^2\rangle=\langle|\nabla^2\phi|^2\rangle-\tfrac12\langle(\Delta\phi)^2\rangle=\tfrac12\langle\rho^2\rangle=2f/(1-f)\). This gives the explicit bound, which is the Hashin-Shtrikman value \(\kappa_*\le\mu f/(1-f)\). (2) \(\langle\operatorname{div}u\rangle=\operatorname{tr}E\) for every \(u=Ex+w\), and the constraint makes the left side zero. (3) \(u\mapsto Ru(R^{\mathsf T}\cdot)\) maps admissible fields for \(E\) to admissible fields for \(RER^{\mathsf T}\) and preserves the energy. \(\square\)

If \(C_*\) is isotropic, \(Q_G(E)=\kappa_*(\operatorname{tr}E)^2+2\mu_*|E_d|^2\), so \(Q_G(I)=4\kappa_*\) and \(Q_G(E)=2\mu_*|E|^2\) for trace-free \(E\). The pure solid would have \(\mu_*=\mu\), \(\kappa_*=\infty\).

\textbf{The closure and \(S\).} Let \(\mathcal G\) be the closure, in the six-dimensional space of quadratic forms on \(\mathrm{Sym}\), of the set of finite \(Q_G\) over all periodic mixtures (all lattices, all solid fractions). Following Chen, eq. (2),

\[
S=\sup\Big\{\frac{\mu_*}{\mu}:\ Q\in\mathcal G,\ Q(E)=\kappa_*(\operatorname{tr}E)^2+2\mu_*|E_d|^2\ \text{with}\ \kappa_*=0\Big\}.
\]

Milton's own candidate \(4\mu/5\) is a limit of this kind (his \(\varepsilon\to0\)), and so are Chen's tensors (his \(N\to\infty\)).

\begin{quote}
\textbf{Theorem.} \(S=1\). More precisely: 1. Every \(Q\in\mathcal G\) satisfies \(Q(E)\le2\mu|E|^2\) for trace-free \(E\). Hence \(S\le1\). 2. For every \(M\ge11\) the hexagonal set \(G_M\) defined below is a periodic mixture with void of positive area whose effective tensor is isotropic, with moduli \(\kappa_M,\mu_M\) obeying the inequalities displayed at the top. 3. Consequently \(Q_{G_M}\to Q_\infty\), \(Q_\infty(E)=2\mu|E_d|^2\). So \(Q_\infty\in\mathcal G\), it is isotropic with \((\kappa_*,\mu_*)=(0,\mu)\), and \(S=1\) is attained in \(\mathcal G\). 4. No periodic mixture has a finite isotropic effective tensor with \(\mu_*=\mu\).
\end{quote}

\paragraph*{Four lemmas valid for every periodic mixture}

\textbf{Lemma 2.1 (uniform strain).} For trace-free \(E\), \(Q_G(E)\le2\mu f|E|^2\).

\emph{Proof.} \(u=Ex\) is admissible because \(\operatorname{div}u=\operatorname{tr}E=0\). \(\square\)

\textbf{Lemma 2.2 (a stress test bounds the shear energy from below).} Let \(\tau\in L^2_{\rm per}(D;\mathrm{Sym})\) vanish a.e. outside \(G\) and satisfy \(\langle\tau:e(w)\rangle=0\) for all \(w\in H^1_{\rm per}\) (equilibrium on the whole torus). Then for every \(E\in\mathrm{Sym}\)

\[
\big(E:\langle\tau\rangle\big)^2\ \le\ \frac{\langle|\tau_d|^2\rangle}{2\mu}\;Q_G(E). \tag{2.1}
\]

In particular, if \(T=\langle\tau\rangle\neq0\) is trace-free and \(\langle|\tau-T|^2\rangle\le b^2|T|^2\), then \(Q_G(T)\ge2\mu|T|^2/(1+b^2)\).

\emph{Proof.} Let \(u=Ex+w\) be admissible. Then \(\langle\tau:e(u)\rangle=\langle\tau\rangle:E\). Since \(\tau=0\) off \(G\) and \(\operatorname{tr}e(u)=0\) on \(G\), \(\tau:e(u)=\chi\,\tau_d:e(u)\) a.e. Cauchy-Schwarz gives \((E:\langle\tau\rangle)^2\le\langle|\tau_d|^2\rangle\langle\chi|e(u)|^2\rangle\); take the infimum over \(u\). For the second statement put \(E=T\) and note \(\langle|\tau_d|^2\rangle=|T|^2+2T:\langle(\tau-T)_d\rangle+\langle|(\tau-T)_d|^2\rangle\le(1+b^2)|T|^2\), the middle term being zero. \(\square\)

Only the deviatoric part of \(\tau\) is charged: this is the complementary energy of an incompressible solid. The same inequality with \(E=I\) shows that \(Q_G(I)=0\) forces \(\operatorname{tr}\langle\tau\rangle=0\) for every such \(\tau\).

\emph{A factor two that is not needed.} For symmetric \(X\), \(|X_d|^2=\tfrac12|X|^2-\det X\). The field \(\tau-T\) is periodic, divergence-free and of zero mean, hence of the form \(J(\nabla^2A)J^{\mathsf T}\) with \(A\in H^2_{\rm per}\), and \(\det\nabla^2A\) has zero mean. So \(\langle|(\tau-T)_d|^2\rangle=\tfrac12\langle|\tau-T|^2\rangle\) exactly, and the conclusion improves to \(Q_G(T)\ge2\mu|T|^2/(1+b^2/2)\). This is the planar compliance-shift invariance in disguise.

\textbf{Lemma 2.3 (any field with small strain on the solid can be made exactly solenoidal there).} Let \(u^0=x+w^0\) with \(w^0\in H^1_{\rm per}\) and put \(a^2=\langle\chi|e(u^0)|^2\rangle\). If \(a^2\le\tfrac12\), then

\[
Q_G(I)\ \le\ 8(1+\sqrt2)^2\,\mu\,a^2 .
\]

\emph{Proof.} Let \(g=\chi\operatorname{div}u^0\) and \(m=\langle g\rangle\). Since \((\operatorname{tr}A)^2\le2|A|^2\), \(\langle g^2\rangle\le2a^2\) and \(|m|\le\sqrt2\,a\le1\). Let \(\phi\in H^2_{\rm per}\) be the mean-zero solution of \(\Delta\phi=g-m\) and set

\[
B=\nabla^2\phi+\tfrac m2I,\qquad u=\frac{u^0-\nabla\phi-\tfrac m2x}{1-\tfrac m2}.
\]

Then \(u-x=(w^0-\nabla\phi)/(1-\tfrac m2)\in H^1_{\rm per}\), so \(u\) has mean strain exactly \(I\); \(e(u)=(e(u^0)-B)/(1-\tfrac m2)\); and \(\operatorname{tr}B=g\), so \(\operatorname{div}u=0\) a.e. in \(G\). Thus \(u\) is admissible in (1.1). On the torus \(\langle|\nabla^2\phi|^2\rangle=\langle(\Delta\phi)^2\rangle\) (Parseval: \(|k\otimes k|=|k|^2\)) and \(\langle\nabla^2\phi\rangle=0\), so \(\langle|B|^2\rangle=\langle(g-m)^2\rangle+\tfrac12m^2=\langle g^2\rangle-\tfrac12m^2\le2a^2\). With \(1-\tfrac m2\ge\tfrac12\) this gives \(\langle\chi|e(u)|^2\rangle^{1/2}\le2(a+\sqrt2a)\). \(\square\)

\textbf{Lemma 2.4 (a solid strip along a lattice direction resists dilation).} Let \(t\) be a unit vector, \(\nu\perp t\), and let \(\lambda>0\) be the smallest number with \(\lambda t\in\Lambda\) (so \(\lambda t\) is a primitive lattice vector). Suppose \(G\) contains the strip \(\{|x\cdot\nu-c|\le\delta\}\) with \(2\delta\le|D|/\lambda\). Then

\[
Q_G(I)\ \ge\ \frac{8\mu\,\delta\,\lambda}{|D|}.
\]

\emph{Proof.} Complete \(\lambda t\) to a basis \((\lambda t,b)\) of \(\Lambda\). The parallelogram spanned by \(\lambda t\) and \(b\) is a period cell; its extent in the direction \(\nu\) is \(|b\cdot\nu|=|D|/\lambda\ge2\delta\), so after a translation it contains a piece of the strip of area \(2\delta\lambda\) made of segments of length \(\lambda\) parallel to \(t\). Let \(u=x+w\) be admissible. On almost every such segment \(\int_0^\lambda t\cdot e(u)t\,ds=t\cdot(u(x+\lambda t)-u(x))=\lambda\), so \(\int_0^\lambda(t\cdot e(u)t)^2ds\ge\lambda\). In the solid \(\operatorname{tr}e(u)=0\), so \(|e(u)|^2\ge2(t\cdot e(u)t)^2\). Integrating over the strip, \(|D|\langle\chi|e(u)|^2\rangle\ge4\delta\lambda\). \(\square\)

The constant is sharp: for a laminate of solid strips of width \(2\delta\) and period \(P\) the bound reads \(Q_G(I)\ge4\mu f\) with \(f=2\delta/P\), and the field that stretches each strip along itself and contracts it across (strain \(t\otimes t-\nu\otimes\nu\)) has exactly this energy.

\paragraph*{\texorpdfstring{The sets \(G_M\) and the two field estimates}{The sets G\_M and the two field estimates}}

The sets are the periodic frames constructed in the article for Theorems~
\ref{p9:sharp} and~\ref{p11:sharp}. The coordinates and estimates are in Section~\ref{detail:frame-fields}.

\begin{itemize}
\tightlist
\item
  \emph{Square pattern, period 1.} \(h=(M+\tfrac1{100})^{-1}\). Solid square blocks of side \(L=Mh\) are separated by gaps of width \(d=h/100\). Each block face carries \(2M\) patches of length \(h/12\). A patch is bridged across the gap by a bundle of \(N=2M^2\) thin parallel rods in direction \(n+st\), then a small solid collar, then a second bundle in the slightly different direction \(n+(s+\eta)t\), \(\eta=h\); here \(n,t\) are the face normal and tangent and \(s=\pm1\) alternates along the face. Successive collars on a face are tied by bundles of thin rods parallel to the face, and the last collars of the two positive faces of a block are joined by a solid corner piece.
\item
  \emph{Hexagonal cell.} Take the equilateral triangle \(\Delta\) with vertices \((0,0),(1,0),(\tfrac12,\tfrac{\sqrt3}2)\). Fill the band within distance \(3h\) of \(\partial\Delta\) with solid and put the square pattern, scaled to period \(\ell=h^2\), in the rest. Rotate by multiples of \(60^\circ\) about the origin to fill the unit hexagon and repeat it by the lattice \(\Lambda\) generated by \((\tfrac32,\tfrac{\sqrt3}2)\) and \((0,\sqrt3)\); \(|D|=3\sqrt3/2\). This is \(G_M\). It is a locally finite union of closed convex polygons.
\end{itemize}

\emph{Why this can work.} If the blocks translate apart by a uniform dilation, each collar can slide (by the large amount \(1/\eta\) per unit opening) so that both of its rod bundles see end displacements exactly transverse to the rods; thin rods bend at negligible cost, so dilation is nearly free. A collar on its own would free every relative motion of its two blocks. The ties along a face force the two families \(s=\pm1\) to slide together, which blocks relative tangential motion, and the corner piece forces the two faces to open by the same amount. What remains free is pure dilation, and by duality a macroscopic stress loads the ties and corner pieces only through its trace. Under a trace-free load the rods carry axial forces of order one per unit face length while the solid blocks carry an almost uniform stress.

Both estimates below are statements about \(L^2\) fields on the indicator \(\chi_M\) of \(G_M\). No elastic constant occurs in them.

\textbf{Proposition D (displacement).} For every \(M\ge11\) there is a Lipschitz \(u^0:\mathbb R^2\to\mathbb R^2\) with \(u^0-x\) \(\Lambda\)-periodic and \[
\big\langle\chi_M\,|e(u^0)|^2\big\rangle\ \le\ C_0\,h .
\]

\textbf{Proposition S (stress).} For every \(M\ge11\) and every constant trace-free \(T\in\mathrm{Sym}\) there is \(\tau\in L^2_{\rm per}(D;\mathrm{Sym})\), depending linearly on \(T\), with \[
\tau=0\ \text{a.e. outside }G_M,\qquad\langle\tau:e(w)\rangle=0\ \ \forall w\in H^1_{\rm per},\qquad\langle\tau\rangle=T,\qquad\big\langle|\tau-T|^2\big\rangle\le C_S\,h\,|T|^2 .
\]

The field \(u^0\) is not admissible in (1.1): it has \(e(u^0)=I\) on the solid bands and a uniaxial bending strain in the rods, so its divergence does not vanish on the solid. For a compressible solid that only costs energy; for the incompressible solid it is forbidden, and Lemma 2.3 is what repairs it, exactly and at the price of a numerical factor.

The constants are independent of the mesh parameter. Section~\ref{detail:frame-fields} constructs the displacement and the balanced stresses and estimates both on the same geometry.

Three elementary facts about \(G_M\), \(M\ge11\), are also used:

\begin{itemize}
\tightlist
\item
  \textbf{(G1)} \(G_M\) and \(\Lambda\) are invariant under rotation by \(60^\circ\) about the origin.
\item
  \textbf{(G2)} The void has positive area, so \(f_M<1\). Indeed the patterned triangle \(\{d_\Delta>3h\}\) has inradius \(\sqrt3/6-3h\), which exceeds \(\sqrt2h^2\) for \(M\ge11\), so it contains a full microcell, and every microcell contains open void (for instance the corner square of the gap away from the corner piece).
\item
  \textbf{(G3)} The strip \(\{|x_2|\le3h\}\) is solid: the line \(x_2=0\) consists of triangle edges, and the bands of the triangles on both sides cover the strip. Also \((3,0)\in\Lambda\).
\end{itemize}

\paragraph*{Proof of the Theorem}

\textbf{(1) Upper bound.} Lemma 2.1 gives \(Q_G(E)\le2\mu|E|^2\) for trace-free \(E\) and every mixture; the inequality is closed, so it holds on \(\mathcal G\). For an isotropic \(Q\in\mathcal G\) this says \(\mu_*\le\mu\).

\textbf{(2) The mixtures \(G_M\).} Fix \(M\ge11\). By (G2) and Lemma 1.1, \(Q_{G_M}\) is a finite quadratic form. By (G1) and Lemma 1.1(3) it is invariant under \(E\mapsto RER^{\mathsf T}\) with \(R\) the rotation by \(60^\circ\). That action fixes \(I\) and rotates the plane of trace-free matrices by \(120^\circ\). A quadratic form invariant under it has no coupling between \(I\) and the trace-free plane and is a multiple of \(|E|^2\) on that plane (a symmetric operator commuting with a rotation that has no real eigenvector is scalar). Hence

\[
Q_{G_M}(E)=\kappa_M(\operatorname{tr}E)^2+2\mu_M|E_d|^2 .
\]

\emph{Shear.} Lemma 2.1 gives \(\mu_M\le f_M\mu<\mu\). Proposition S and Lemma 2.2 with \(b^2=C_Sh\) give \(2\mu_M|T|^2=Q_{G_M}(T)\ge2\mu|T|^2/(1+C_Sh)\). So \(C'=C_S\) works; with the factor two noted after Lemma 2.2, so does \(C'=C_S/2\) (the form used in the reference text).

\emph{Bulk, upper bound.} Let \(M_0\) be such that \(C_0h\le\tfrac12\) for \(M\ge M_0\). For those \(M\), Proposition D and Lemma 2.3 give \(4\kappa_M=Q_{G_M}(I)\le8(1+\sqrt2)^2C_0\,\mu h\), that is \(\kappa_M\le2(1+\sqrt2)^2C_0\,\mu h\). For the finitely many \(M\) with \(11\le M<M_0\), Lemma 1.1 gives the explicit but poor bound \(\kappa_M\le\mu f_M/(1-f_M)\), so \(\kappa_M\le C\mu h\) holds for all \(M\ge11\) after enlarging the constant.

\emph{Bulk, lower bound.} By (G3), \((3,0)\) is a primitive vector of \(\Lambda\) (the horizontal lattice vectors are the multiples of \(2a_1-a_2\)) and the strip \(\{|x_2|\le3h\}\) is solid, with \(6h<\sqrt3/2=|D|/3\) for \(M\ge11\). Lemma 2.4 with \(t=e_1\), \(\lambda=3\), \(\delta=3h\), \(|D|=3\sqrt3/2\) already gives \(4\kappa_M\ge16\sqrt3\,\mu h\), that is \(\kappa_M\ge4\sqrt3\,\mu h>0\). (The same follows from (2.1) with \(E=I\) and the uniaxial stress \(e_1\otimes e_1\) on these strips.)

The three strip families together give the bound stated in the Theorem; the corresponding complementary test gives the following lower bound. By (G1) and (G3) the strips of half-width \(3h\) about all lattice lines of the three edge directions \(t_1,t_2,t_3\) are solid. Let \(\chi_j\) be the indicator of the strips of direction \(t_j\) and \(s=4\sqrt3\,h<1\) the area fraction of one family. The stress \(\tau=\sum_j\chi_j\,t_j\otimes t_j\) vanishes outside \(G_M\), is equilibrated (each term is a function of one coordinate times \(t_j\otimes t_j\), with \(t_j\) along the strips), and has mean \(\tfrac32sI\). Since \(\chi_j\) depends on one coordinate only, the centred indicators of different families have disjoint Fourier supports, so \(\langle\chi_i\chi_j\rangle=s^2\) for \(i\ne j\). With \(|(t_j\otimes t_j)_d|^2=\tfrac12\) and \((t_i\otimes t_i)_d:(t_j\otimes t_j)_d=(t_i\cdot t_j)^2-\tfrac12=-\tfrac14\) for \(i\ne j\), this gives \(\langle|\tau_d|^2\rangle=\tfrac32s-\tfrac32s^2\). Inequality (2.1) with \(E=I\) reads \((3s)^2\le\tfrac32s(1-s)\,Q_{G_M}(I)/(2\mu)\), hence

\[
\frac{\kappa_M}{\mu}\ \ge\ \frac{3s}{1-s}=\frac{12\sqrt3\,h}{1-4\sqrt3\,h}\ \ge\ 12\sqrt3\,h .
\]

\textbf{(3) The limit.} As \(M\to\infty\), \(h\to0\), so \(\kappa_M\to0\) and \(\mu_M\to\mu\), that is \(Q_{G_M}\to Q_\infty=2\mu|E_d|^2\) in the space of quadratic forms. Each \(Q_{G_M}\) is the finite effective form of a periodic mixture of the incompressible solid and void, so \(Q_\infty\in\mathcal G\) by the definition of the closure. It is isotropic with \(\kappa_*=0\) and \(\mu_*=\mu\). With (1), \(S=1\).

\textbf{(4) Not attained by a mixture.} If the void has positive area, \(\mu_*\le f\mu<\mu\) by Lemma 2.1. If it has zero area the mixture is the pure solid up to a null set: its shear modulus is \(\mu\) but its bulk modulus is infinite and there is no finite tensor (Lemma 1.1(2)). \(\blacksquare\)

A closure-free form of the result: \emph{for every \(\varepsilon>0\) there is a periodic, exactly isotropic mixture of the incompressible solid and void with \(0<\kappa_*\le\varepsilon\mu\) and \((1-\varepsilon)\mu\le\mu_*<\mu\).}
\phantomsection\label{endfile:frame-tests}

\subsection{Displacement and stress on the same hexagonal frame}\label{detail:frame-fields}\label{file:frame-fields}
Throughout, \(C\) denotes absolute constants (independent of \(M\); they are unrelated to the constants \(C,C'\) of the Theorem), \(J\) is the rotation by \(\pi/2\), \(\operatorname{inc}A=J(\nabla^2A)J^{\mathsf T}\) for a scalar \(A\), and \(x\times f=x_1f_2-x_2f_1\).

\paragraph*{A.1 The solid on the square torus}

Fix an integer \(M\ge2\) and put

\[
h=\frac1{M+1/100},\quad L=Mh=1-d,\quad d=\frac h{100},\quad\eta=h,\]\[N=2M^2,\quad p=\frac h{12},\quad z_0=\frac d2,\quad a_p=\frac p2+\frac d8 .
\]

Blocks are \(B_R=R+[0,L]^2\), \(R\in\mathbb Z^2\). On the positive face \(i\in\{1,2\}\) of \(B_R\) let \(n=e_i\), \(t=e_{3-i}\), and use coordinates \((z,y)\) for the point \(R+(L+z)n+yt\). The face carries \(2M\) patches, numbered \(r=1,\dots,2M\) in increasing \(y\), with centres \(y_r=kh+h/4\) (sign \(s=+1\)) and \(kh+3h/4\) (sign \(s=-1\)), \(k=0,\dots,M-1\). Put \(m_s=n+st\), \(l_s=n+(s+\eta)t\), \(c_r=y_r+sz_0\).

{\def\LTcaptype{none} 
\begin{longtable}[]{@{}
  >{\raggedright\arraybackslash}p{(\linewidth - 2\tabcolsep) * \real{0.5000}}
  >{\raggedright\arraybackslash}p{(\linewidth - 2\tabcolsep) * \real{0.5000}}@{}}
\toprule\noalign{}
\begin{minipage}[b]{\linewidth}\raggedright
Piece
\end{minipage} & \begin{minipage}[b]{\linewidth}\raggedright
Set, in face coordinates \((z,y)\)
\end{minipage} \\
\midrule\noalign{}
\endhead
\bottomrule\noalign{}
\endlastfoot
first rods of patch \(r\) & \(0\le z\le d/3\), \(\ y-y_r-sz+p/2\in\bigcup_{j<N}\big[(j+\tfrac14)\tfrac pN,(j+\tfrac34)\tfrac pN\big]\) \\
collar & \(d/3\le z\le 2d/3\), \(\ \lvert y-y_r-sz\rvert\le p/2\) \\
last rods & \(2d/3\le z\le d\), same \(N\) intervals for \(y-y_r-sz-\eta(z-2d/3)+p/2\) \\
tab & \(\lvert z-z_0\rvert\le d/16\), \(\ \lvert y-c_r\rvert\le a_p\) \\
connector \(r\to r+1\) & \(c_r+a_p\le y\le c_{r+1}-a_p\), \(\ z-z_0+d/32\in\bigcup_{j<N}\big[(j+\tfrac14)\tfrac d{16N},(j+\tfrac34)\tfrac d{16N}\big]\) \\
arm of face \(i\) & \(\lvert z-z_0\rvert\le d/32\), \(\ c_{2M}+a_p\le y\le L+z_0+d/32\) \\
\end{longtable}
}

First rods run in direction \(m_s\) from the block to the collar, last rods in direction \(l_s\) from the collar to \(B_{R+n}\). A \emph{collar body} is a collar with its tab. The two last collars (\(r=2M\), both of sign \(-1\)), their tabs and the two arms form one \emph{corner body}; the arms overlap in the square \([L+z_0-d/32,\,L+z_0+d/32]^2\). Everything else in the gaps is void. The resulting closed \(\mathbb Z^2\)-periodic set is \(G^\square_M\).

\textbf{Lemma A.1 (welds and clearances).} For every \(M\ge2\):

\begin{enumerate}
\def\labelenumi{\arabic{enumi}.}
\tightlist
\item
  Each first rod meets exactly two other bodies: the block \(B_R\) along its end \(z=0\) and its collar along its end \(z=d/3\). Each last rod meets its collar along \(z=2d/3\) and \(B_{R+n}\) along \(z=d\). Each connector rod meets the two tab end faces \(y=c_r+a_p\) and \(y=c_{r+1}-a_p\). Every such contact is the whole end segment of the rod.
\item
  All other pairs of distinct bodies (block, \(4M-1\) collar and corner bodies, individual rods) are at positive distance.
\item
  Rescaled by \(1/h\), the collar bodies and the corner body are congruent to finitely many fixed connected Lipschitz polygons, independent of \(M\).
\end{enumerate}

\emph{Proof.} The structures of patch \(r\) (rods, collar, tab) lie in \(\lvert y-y_r\rvert\le p/2+d+\eta d/3<h/16\), and patch centres are \(h/2\) apart, so different patches are separated by more than \(3h/8\) and every footprint stays in its own length-\(h\) tile of both block faces. Within a patch, the first rods occupy \(z\le d/3\), the last rods \(z\ge2d/3\), and the tab \(7d/16\le z\le9d/16\), strictly inside the collar's range; so rods touch the collar body only on the collar's two faces, on segments at distance \(p/(4N)\) from the collar's corners. Rods of one bundle are parallel strips separated by \(p/(2N)\) along the face. Connectors occupy \(15d/32\le z\le17d/32\), where the only other material is the tab, which contains the collar's section there (\(\lvert y-c_r\rvert\le p/2+d/32<a_p\)); so a connector abuts only the flat end face of the tab, and it stays at distance at least \(3d/32\) in \(z\) from the parts of the collar that protrude beyond the tab in \(y\). Each arm keeps distance at least \(15d/32\) from \(B_R\), \(B_{R+e_1}\), \(B_{R+e_2}\), lies beyond the last patch, and ends \(15d/32\) short of the cell boundary. Every length in the bodies is a fixed multiple of \(h\) (arm length \(h(\tfrac14-\tfrac1{24}+\tfrac{29}{3200})\)), which gives 3. \(\square\)

\emph{Void fraction.} In the gap of face \(i\), the points with \(0<y<L\), \(0<z<d\), \(\lvert z-z_0\rvert>d/16\) and \(\lvert y-y_r\rvert>h/16\) for all \(r\) are void: patch structures stay within \(h/16\) of \(y_r\), connectors and arms within \(d/32\) of \(z_0\), and the structures of the other face lie in \(y\ge L\). This set has area \((3L/4)(7d/8)\) per face. Hence the void fraction of the square cell is between \(\tfrac{21}{16}Ld>0.0129\,h\) and \(1-L^2<2d=0.02\,h\).

\paragraph*{A.2 The dilation field on the square torus}

Blocks get the translation \(U(B_R)=R\). A collar body of sign \(s\) on a positive face of \(B_R\) gets \(R+v_s/\eta\) with \(v_s=t-sn\). The corner body gets \(R+(e_1+e_2)/\eta\), consistently, since \(v_-=t+n=e_1+e_2\) on both faces. The relative end translation \(D\) of each kind of rod is orthogonal to its direction:

\[
\tfrac1\eta v_s\cdot m_s=0,\qquad\big(n-\tfrac1\eta v_s\big)\cdot l_s=1+\tfrac s\eta-\tfrac{s+\eta}\eta=0,\qquad\tfrac1\eta(v_{s'}-v_s)\cdot t=0 .
\]

\emph{One rod.} Let the rod be \(\{x_0+\xi a+zb:\ \lvert\xi\rvert\le\rho/2,\ 0\le z\le H\}\), where \(\nu\) is the unit normal of its ends, \(a\perp\nu\) a unit vector, and \(b\) its direction normalised by \(b\cdot\nu=1\) (so the area element is \(d\xi\,dz\)). Put \(r=b\cdot(x-x_0)/(H\lvert b\rvert^2)\), \(q=\min(1,\max(0,2r-\tfrac12))\), \(\psi=3q^2-2q^3\) as a function of \(r\), and

\[
u(x)=U_0+\psi(r)D-\psi'(r)\,\frac{D\cdot(x-x_0)}{H\lvert b\rvert^2}\,b . \tag{A.1}
\]

Differentiating, \(\nabla u=\frac{\psi'}{H\lvert b\rvert^2}(D\otimes b-b\otimes D)-\psi''\frac{D\cdot(x-x_0)}{H^2\lvert b\rvert^4}\,b\otimes b\). The first term is skew, and \(D\cdot(x-x_0)=\xi\,D\cdot a\) because \(D\cdot b=0\). Hence

\[
e(u)=\gamma\,\hat b\otimes\hat b,\qquad\gamma=-\frac{\psi''(r)\,\xi\,(D\cdot a)}{H^2\lvert b\rvert^2},\qquad\hat b=b/\lvert b\rvert . \tag{A.2}
\]

\emph{The welded ends.} On the end \(z=0\), \(r=\xi\,(a\cdot b)/(H\lvert b\rvert^2)\), and on the end \(z=H\) it is \(1\) plus this. With \(b=\nu+\alpha a\), \(\lvert\alpha\rvert/(1+\alpha^2)\le\tfrac12\), so the shift is at most \(\rho/(4H)\). For face rods \(\rho=p/(2N)\), \(H=d/3\), giving \(25/(8N)\); for connectors \(a\cdot b=0\). So for \(M\ge3\) (\(N\ge18\)) the shift is at most \(25/144<\tfrac14\): \(\psi=\psi'=0\) on the whole first end and \(\psi=1\), \(\psi'=0\) on the whole far end, and \(u\) equals the rigid translation of the adjoining body on the entire oblique weld. (At \(M=2\) this fails; that index is not used.)

\emph{Strain integral.} For \(M\ge3\) every line \(\xi=\) const crosses the whole support of \(\psi''\), \(dr=dz/H\), \(\int\psi''^2dr=96\) and \(\int\xi^2d\xi=\rho^3/12\). So

\[
\int_{\rm rod}\lvert e(u)\rvert^2=\int_{\rm rod}\gamma^2=\frac{8\,\rho^3\,(D\cdot a)^2}{H^3\lvert b\rvert^4} . \tag{A.3}
\]

For all face rods \((D\cdot a)^2=1/\eta^2\) and \(\rho^3/H^3=15625/(8N^3)\). Summing the \(4M\) first bundles (\(\lvert b\rvert^4=4\)), the \(2M\) last bundles of each sign, and the connectors (\(\lvert D\cdot a\rvert=2/\eta\), \(\rho=d/(32N)\), lengths \(H_\mp=h/2\mp d-2a_p\), \(M\) of the first and \(M-1\) of the second kind per face):

\[
\mathcal E_M:=\frac{15625}4\big(1+2\Xi(\eta)\big)\frac h{L^3}+c_M,\qquad\Xi(\eta)=\frac1{(2+2\eta+\eta^2)^2}+\frac1{(2-2\eta+\eta^2)^2}, \tag{A.4}
\]

with \(c_M=2\big[(3/1940)^3+(1-1/M)(3/2036)^3\big]h/L^3<1.4\cdot10^{-8}h\). Here \(h=h_M\) depends on \(M\). \(\Xi\) increases on \([0,\tfrac13]\) from \(\tfrac12\) and \(1/L^3\) decreases in \(M\), so the first term divided by \(h\) is largest at \(M=3\), where it equals \(8744.574\ldots\). Hence

\[
\frac{\mathcal E_M}{h}<8745\quad(M\ge3),\qquad\frac{\mathcal E_M}h\to7812.5\ldots\quad(M\to\infty) . \tag{A.5}
\]

\textbf{Lemma A.2 (dilation).} For every \(M\ge3\) there is a \(\mathbb Z^2\)-periodic Lipschitz \(w\) on \(\mathbb R^2\) with

\[
\int_{[0,1]^2}\mathbf 1_{G^\square_M}\lvert I+e(w)\rvert^2=\mathcal E_M\le8745\,h,\qquad\lvert I+e(w)\rvert\le15700,\quad\lvert w\rvert\le\frac5h\ \text{ on }G^\square_M .
\]

\emph{Proof.} Define \(u\) on \(G^\square_M\) by the translations on the bodies and (A.1) on the rods, and \(w=u-x\). By Lemma A.1 and the computation at the ends, \(u\) is continuous on \(G^\square_M\), and \(u(x+e_i)=u(x)+e_i\). The strain vanishes on the bodies and is (A.2) on the rods; this gives the integral and \(\sup\lvert e(u)\rvert=12\rho\lvert D\cdot a\rvert/(H^2\lvert b\rvert^2)\le22500/\big(L^2(1+(1-\eta)^2)\big)\le15700\). Every body translation differs from \(x\) by at most \(\sqrt2+\sqrt2/\eta\) at the points of the body and of the rods welded to it. On a rod \(U_0+\psi D\) is a convex combination of the two end translations, and the last term of (A.1) is at most \(1.5(\rho/H)\lvert D\rvert\) (since \(\lvert\psi'\rvert\le3\) and \(\lvert b\rvert\ge1\)): for face rods \(\rho/H=12.5/N\le0.7\) and \(\lvert D\rvert\le1+\sqrt2/\eta\); for connectors \(\rho/H<10^{-3}\) and \(\lvert D\rvert=2/\eta\). So \(\lvert w\rvert\le\sqrt2+\sqrt2/\eta+1.05(1+\sqrt2/\eta)\le5/h\).

\emph{Extension into the void.} \(w\) is Lipschitz on \(G^\square_M\) for the Euclidean distance, with a constant depending on \(M\): on each rod \(u\) is continuous and piecewise smooth with bounded first derivatives (the gradient jumps with \(\psi''\) across the two lines \(r=\tfrac14,\tfrac34\)), hence Lipschitz, and it is constant on each body; a rod and a body welded to it carry the same constant near the weld (the plateau); and any other two points of different pieces are a positive distance apart (Lemma A.1) while \(w\) is bounded. The componentwise McShane formula \(\tilde w_k(x)=\inf_{y\in G^\square_M}[w_k(y)+\mathrm{Lip}\,\lvert x-y\rvert]\) is a periodic Lipschitz extension. Only the values on the solid enter the estimates. \(\square\)

\paragraph*{A.3 Three standard lemmas}

\textbf{Lemma N (Neumann).} Let \(P\) be a bounded connected Lipschitz domain and \(g\in L^2(\partial P;\mathbb R^2)\) with \(\int g=0\) and \(\int x\times g=0\). There is \(\sigma\in L^2(P;\mathrm{Sym})\), depending linearly on \(g\), with \(\int_P\sigma:e(v)=\int_{\partial P}g\cdot v\) for all \(v\in H^1(P)\) and

\[
\|\sigma\|_{L^2(P)}\le K_P\,\sup\Big\{\int_{\partial P}g\cdot v:\ \|v\|_{H^1(P)}\le1\Big\}.
\]

\emph{Proof.} By Korn's inequality each \(v\) has a rigid motion \(\varrho\) with \(\|v-\varrho\|_{H^1}\le K_P\|e(v)\|_{L^2}\), so \(H^1/\{\text{rigid}\}\) is a Hilbert space for \((u,v)\mapsto\int e(u):e(v)\). The functional \(v\mapsto\int g\cdot v\) vanishes on rigid motions, so Riesz gives \(u\) with \(\int e(u):e(v)=\int g\cdot v\). Take \(\sigma=e(u)\). \(\square\)

This is an auxiliary problem with no material constant; only the resulting equilibrium field \(\sigma\) is used.

\textbf{Lemma P (tiles).} If \(S\) is an interval of length \(h\) and \(\bar v_S\) the mean of \(v\) on \(S\),

\[
\int_S\lvert v-\bar v_S\rvert^2=\frac1{2h}\iint_{S\times S}\lvert v(x)-v(y)\rvert^2\le\frac h2\iint_{S\times S}\frac{\lvert v(x)-v(y)\rvert^2}{\lvert x-y\rvert^2}.
\]

\textbf{Lemma A (Airy).} If \(\delta\in L^2_{\rm per}([0,1]^2;\mathrm{Sym})\) has zero mean and zero distributional divergence, there is a periodic \(A\in H^2\) with \(\operatorname{inc}A=\delta\) and \(\|A\|_{H^2}\le C\|\delta\|_{L^2}\). \emph{Proof.} For \(k\ne0\), \(\hat\delta(k)k=0\) and symmetry give \(\hat\delta(k)=c_k\,Jk\otimes Jk\); put \(\hat A(k)=-c_k/(4\pi^2)\), \(\hat A(0)=0\). \(\square\)

\paragraph*{A.4 Equilibrated stresses on the square torus}

Let \(T\in\mathrm{Sym}\) with \(\operatorname{tr}T=0\); on face \(i\) put \(p_i=n\cdot Tn\), \(g_i=t\cdot Tn\) (so \(p_2=-p_1\)). A \emph{load \(f\)} on a set of parallel rods with direction \(b\) (\(b\cdot\nu=1\)) and total retained width \(\omega\) means the constant stress \((f/\omega)\,b\otimes b\) in those rods. It is divergence-free, has zero traction on the rod sides, and pulls the body at the start with \(+fb\) and the body at the far end with \(-fb\), uniformly over the welds.

\[
q_{is}=\tfrac12\big[p_i+s(g_i-\eta p_i)\big],\qquad b_{ir}=\eta\,q_{is_r}/M .
\]

{\def\LTcaptype{none} 
\begin{longtable}[]{@{}
  >{\raggedright\arraybackslash}p{(\linewidth - 4\tabcolsep) * \real{0.3333}}
  >{\raggedright\arraybackslash}p{(\linewidth - 4\tabcolsep) * \real{0.3333}}
  >{\raggedright\arraybackslash}p{(\linewidth - 4\tabcolsep) * \real{0.3333}}@{}}
\toprule\noalign{}
\begin{minipage}[b]{\linewidth}\raggedright
Rods
\end{minipage} & \begin{minipage}[b]{\linewidth}\raggedright
Direction
\end{minipage} & \begin{minipage}[b]{\linewidth}\raggedright
Load
\end{minipage} \\
\midrule\noalign{}
\endhead
\bottomrule\noalign{}
\endlastfoot
first bundle of patch \((i,r)\) & \(m_s\) & \(q_{is}/M\); minus \(\eta p_1\) if \((i,r)=(1,2M)\); in addition \(-a_{ir}\) on the upper half and \(+a_{ir}\) on the lower half of the patch \\
last bundle of patch \((i,r)\) & \(l_s\) & \(q_{is}/M\) \\
connector \(r\to r+1\) on face \(i\) & \(t\) & \(-\sum_{a\le r}b_{ia}\) \\
\end{longtable}
}

with \(a_{ir}=d\,b_{ir}/(3p)\) for \(r<2M\), \(a_{2,2M}=0\), \(a_{1,2M}=2H_C/p\), and

\[
H_C=\frac d6\big(b_{1,2M}-b_{2,2M}\big)+\eta p_1\Big(\frac h4+d\Big).
\]

\emph{Flux.} \(\sum_sq_{is}l_s=p_in+\big[(g_i-\eta p_i)+\eta p_i\big]t=Tn\), so the last bundles of each tile carry \(Tn/M\).

\emph{Ordinary collar \(r<2M\).} The forces \(-\tfrac qMm_s\) and \(+\tfrac qMl_s\) act along lines through the centre of the collar's top face; their sum is \(b_{ir}t\) at \(z=2d/3\). The two connectors contribute \(+\sum_{a<r}b_{ia}t-\sum_{a\le r}b_{ia}t=-b_{ir}t\) at \(z=z_0\) (their retained intervals are symmetric about \(z_0\)). Forces cancel; the residual moment is \((n\times t)\,d\,b_{ir}/6\). The half-patch pair gives the collar \(+a\,m_s\) and \(-a\,m_s\) at points \((p/2)t\) apart (retained intervals are centred and \(N\) is even), a pure couple \((p/2)\,t\times a\,m_s=-(n\times t)pa/2\), which cancels the residual for \(a=db_{ir}/(3p)\).

\emph{Corner body.} Face \(i\) delivers \(\sum_rb_{ir}\,t=\eta p_i\,t\), in total \(\eta(p_1e_2+p_2e_1)=-\eta p_1(e_1-e_2)\). \textbf{This is where \(\operatorname{tr}T=0\) enters}: the total is orthogonal to the corner's dilation direction \(e_1+e_2\) precisely because \(p_1+p_2=0\), and it is parallel to \(m_-=e_1-e_2\), the first-bundle direction of the face-1 last collar. The extra load \(-\eta p_1\) on that bundle gives the corner \(+\eta p_1(e_1-e_2)\) and balances forces. About \(O=(L+z_0,L+z_0)\), where both connector lines meet, the remaining moment is \(H_C\) (the added force acts on the line \(x_1+x_2=2L-h/4\)), and the half-patch pair with \(a=2H_C/p\) cancels it.

\emph{Block.} Each rod exerts equal and opposite collinear forces at its ends. Summing over the cell quotient, forces cancel, and a rod ending on \(B_{R+e_i}\) contributes \(e_i\times f\) to the moment after its far end is translated back; the total is \(\sum_ie_i\times Te_i=T_{21}-T_{12}=0\). Since every collar body and the corner balance, so does the block.

\emph{Sizes.} \(\lvert q_{is}\rvert\le2\lvert T\rvert\). Primary loads are at most \(2\lvert T\rvert/M\) on retained width \(p/2\). Connector loads are at most \(4\eta\lvert T\rvert\) on retained width \(d/32\). \(\lvert H_C\rvert\le0.3\,\eta h\lvert T\rvert\), so the corner pair is at most \(8\eta\lvert T\rvert\) on retained width \(p/4\). Because \(\eta=h\) and \(Mh=L\), every rod stress and every weld traction is bounded by \(C\lvert T\rvert\) (the crude bounds give \(C\le128\,\eta/d=1.28\cdot10^4\)).

\textbf{Lemma A.3 (stress on the square torus).} For every \(M\ge2\) there is a linear map \(T\mapsto\tau(T)\) from trace-free \(T\in\mathrm{Sym}\) to symmetric fields in \(L^2_{\rm per}([0,1]^2)\) with

\[
\tau(T)=0\ \text{a.e. outside }G^\square_M,\qquad\int\tau(T):e(v)=0\ \ \forall v\in H^1_{\rm per},\]\[\langle\tau(T)\rangle=T,\qquad\int_{[0,1]^2}\lvert\tau(T)-T\rvert^2\le C_\square\,h\,\lvert T\rvert^2 .
\]

\emph{Proof.} Put the rod stresses above in the rods and zero in the pores.

\emph{Collar and corner bodies.} For such a body \(P=x_c+h\hat P\) let \(g\) be the traction exerted on it by its rods: bounded by \(C\lvert T\rvert\), supported on the welds, with zero resultant force and moment. Apply Lemma N on the reference polygon \(\hat P\) to \(\hat g(\hat y)=g(x_c+h\hat y)\) and set \(\sigma(x)=\hat\sigma((x-x_c)/h)\). Then \(\sigma\) solves the weak Neumann problem on \(P\) and \(\int_P\lvert\sigma\rvert^2=h^2\int_{\hat P}\lvert\hat\sigma\rvert^2\le Ch^2\lvert T\rvert^2\), by the trace inequality on the finitely many reference shapes (Lemma A.1(3)). The \(4M-1\) bodies cost \(O(h)\lvert T\rvert^2\) in total.

\emph{Rods.} Bounded stress on area at most \(1-L^2\le2d\): \(O(h)\lvert T\rvert^2\).

\emph{Block.} Let \(\gamma\) be the actual traction on \(\partial B\) minus that of the constant stress \(T/L\). Both are balanced, and \(\lvert\gamma\rvert\le C\lvert T\rvert\). Split each side into its \(M\) tiles. A negative-face tile receives exactly \(-Te_i/M\), the comparison value. A positive-face tile receives \((Te_i-\eta p_it)/M\), the half-patch pairs having zero resultant, except the single tile holding patch \((1,2M)\), which also receives \(-\eta p_1m_-\). So the tile means \(\bar\gamma\) are \(O(h\lvert T\rvert)\) except on one tile where they are \(O(\lvert T\rvert)\): \(\|\bar\gamma\|_{L^2(\partial B)}\le C\sqrt h\lvert T\rvert\). The remainder \(\tilde\gamma=\gamma-\bar\gamma\) has zero mean on every tile, so by Lemma P and Cauchy-Schwarz

\[
\Big\lvert\int_{\partial B}\tilde\gamma\cdot v\Big\rvert=\Big\lvert\sum_S\int_S\tilde\gamma\cdot(v-\bar v_S)\Big\rvert\le\sqrt{h/2}\;\|\tilde\gamma\|_{L^2(\partial B)}\Big(\sum_{\text{sides}}\lvert v\rvert^2_{H^{1/2}}\Big)^{1/2}\le C\sqrt h\,\lvert T\rvert\,\|v\|_{H^1(B)} .
\]

Lemma N (constants uniform for squares of side \(L\in[0.99,1)\)) gives a block stress \(T/L+\sigma'\) with \(\|\sigma'\|_{L^2(B)}\le C\sqrt h\lvert T\rvert\), and \(\lvert T/L-T\rvert\le2d\lvert T\rvert\).

\emph{Assembly.} For periodic \(v\in H^1\), \(\int\tau:e(v)\) is the sum over bodies of \(\int_{\partial P}g\cdot v\) plus the sum over rods of their boundary terms; on each weld the two contributions are opposite and elsewhere they vanish. So \(\tau(T)\) is equilibrated on the torus. Testing with \(v=\varphi(x_i)c\) shows that the flux of \(\tau e_i\) through \(\{x_i=\text{const}\}\) is constant and equals \(\langle\tau\rangle e_i\); on \(L+2d/3<x_i<1\) only last rods of face \(i\) are cut, with total flux \(Te_i\). Adding the pieces, and \(\lvert T\rvert^2\) times the gap area for the pores, gives the estimate. All choices are linear in \(T\). \(\square\)

\paragraph*{A.5 The hexagonal set; proofs of Propositions D and S}

Let \(\Delta\) be the triangle with vertices \((0,0),(1,0),(\tfrac12,\tfrac{\sqrt3}2)\), \(d_1,d_2,d_3\) the distances to its three edge lines, \(d_\Delta=\min_jd_j\), and \(\ell=h^2\). For \(M\ge11\),

\[
G_M\cap\Delta=\{d_\Delta\le3h\}\ \cup\ \{x:\ d_\Delta(x)>3h,\ x/\ell\in G^\square_M\}.
\]

The other five sectors of the unit hexagon are the rotations of this one by multiples of \(\pi/3\) about the origin, and the hexagon is repeated by the lattice \(\Lambda\) generated by \(a_1=(\tfrac32,\tfrac{\sqrt3}2)\) and \(a_2=(0,\sqrt3)\). The band \(\{d_\Delta\le3h\}\) is solid; its inner boundary cuts through blocks, rods, collars and pores of the pattern at arbitrary positions. (For \(M\le10\), \(3h\) exceeds the inradius \(\sqrt3/6\) and there is no pattern.)

\emph{(G1).} Rotation \(R\) by \(\pi/3\) about the origin permutes the six sectors together with their patterns and maps \(a_1\mapsto a_2\), \(a_2\mapsto a_2-a_1\). So \(RG_M=G_M\) and \(R\Lambda=\Lambda\).

\emph{(G2).} The inner triangle \(\{d_\Delta>3h\}\) has inradius \(\sqrt3/6-3h\), which exceeds \(\sqrt2\ell\) for all \(M\ge11\) (at \(M=11\): \(0.0162>0.0117\); the difference increases with \(M\)). A disc of radius \(\sqrt2\ell\) contains a full microcell, and every microcell contains open void (Section A.1). So the void has positive area. For large \(M\) the microcells inside the inner triangle cover all of it except an \(O(\ell)\)-neighbourhood of its boundary, so by Section A.1 the void fraction of \(G_M\) is at least \(0.01\,h\).

\emph{(G3).} Every segment \([k,k+1]\times\{0\}\), \(k\in\mathbb Z\), is a common edge of two triangles of the tiling (each elementary triangle of the lattice generated by \((1,0)\) and \((\tfrac12,\tfrac{\sqrt3}2)\) is a sector of exactly one hexagon). Each of the two triangles is solid within distance \(3h\) of that edge. Hence \(\{\lvert x_2\rvert\le3h\}\subset G_M\). Also \((3,0)=2a_1-a_2\in\Lambda\).

Fix a smooth non-decreasing \(\theta\) with \(\theta=0\) on \((-\infty,0]\), \(\theta=1\) on \([1,\infty)\), and set

\[
\zeta_u=\prod_j\theta\Big(\frac{d_j-3h}{2h}\Big),\qquad\zeta_\tau=\prod_j\theta\Big(\frac{d_j-h}{2h}\Big),
\]

so \(\zeta_u=0\) on the band and \(\zeta_u=1\) where \(d_\Delta\ge5h\); \(\zeta_\tau=0\) where \(d_\Delta\le h\) and \(\zeta_\tau=1\) on the whole patterned region \(\{d_\Delta\ge3h\}\); \(\lvert\nabla^k\zeta\rvert\le C_kh^{-k}\).

\textbf{Covering.} Let \(F\ge0\) be \(\mathbb Z^2\)-periodic, \(Q\) a rotation, and \(\Delta'\subset\Delta\) a triangle or the region between two nested triangles (perimeter at most \(6\)). The rotated squares \(Q(\ell(k+[0,1]^2))\) that meet \(\Delta'\) lie in its \(\sqrt2\ell\)-neighbourhood, so \(\int_{\Delta'}F(Q^{\mathsf T}x/\ell)\,dx\le(\lvert\Delta'\rvert+C\ell)\int_{[0,1]^2}F\).

\textbf{Proof of Proposition D.} Let \(M\ge11\) and \(w\) be the Lipschitz field of Lemma A.2. In the sector of rotation \(Q\) put

\[
u^0(x)=x+\ell\,\zeta_u(x)\,Q\,w(Q^{\mathsf T}x/\ell).
\]

The corrector vanishes on the band, so the six pieces and their translates glue to a Lipschitz field with \(u^0-x\) \(\Lambda\)-periodic. On \(G_M\), with \(w\) evaluated at \(Q^{\mathsf T}x/\ell\),

\[
e(u^0)=(1-\zeta_u)I+\zeta_u\,Q\big(I+e(w)\big)Q^{\mathsf T}+\ell\operatorname{sym}\big(\nabla\zeta_u\otimes Qw\big).
\]

\begin{itemize}
\tightlist
\item
  Band: \(e(u^0)=I\), \(\lvert e(u^0)\rvert^2=2\), area at most \(9h\) per sector.
\item
  Interior \(\{d_\Delta\ge5h\}\): \(\lvert e(u^0)\rvert^2=\lvert I+e(w)\rvert^2\) on the solid; by covering and Lemma A.2 its integral is at most \((\lvert\Delta\rvert+C\ell)\,\mathcal E_M\).
\item
  Strip \(\{3h<d_\Delta<5h\}\), area at most \(6h\): on the solid the last term is bounded by \(\ell\cdot(C/h)\cdot(5/h)=5C\), using only values of \(w\) on the solid; this is where \(\ell=h^2\) is needed. So \(\lvert e(u^0)\rvert^2\le3\big[2+\lvert I+e(w)\rvert^2+C\big]\) there, and covering bounds its integral by \(Ch\).
\end{itemize}

The truncation of rods and pores by the band is irrelevant: the field is globally defined and single-valued, and it is the identity where the extra solid sits. Summing over the six sectors and dividing by \(\lvert D\rvert=6\lvert\Delta\rvert\),

\[
\big\langle\chi_M\lvert e(u^0)\rvert^2\big\rangle\le(1+C\ell)\,\mathcal E_M+Ch\le C_0h .\qquad\square
\]

\textbf{Proof of Proposition S.} Let \(T\) be trace-free. In the sector of rotation \(Q\) put \(H=Q^{\mathsf T}TQ\) (trace-free), let \(A_H\) be the Airy potential of \(\tau(H)-H\) from Lemmas A.3 and A, so \(\|A_H\|_{H^2}\le C\sqrt h\lvert T\rvert\), and define

\[
\tau=T+\operatorname{inc}_x\!\Big[\ell^2\,\zeta_\tau(x)\,A_H(Q^{\mathsf T}x/\ell)\Big].
\]

The bracket is an \(H^2\) function vanishing near \(\partial\Delta\), so the six pieces and their translates glue to a periodic field which is symmetric, divergence-free in distributions across every join and every truncated pore, with \(\langle\tau\rangle=T\). Where there is void, \(d_\Delta>3h\) and \(\zeta_\tau=1\); since \(JQ=QJ\), \(\tau=T+Q(\tau(H)-H)Q^{\mathsf T}=Q\,\tau(H)(Q^{\mathsf T}x/\ell)\,Q^{\mathsf T}=0\) there. By the product rule

\[
\tau-T=\zeta_\tau\,Q(\tau(H)-H)Q^{\mathsf T}+O\big((\ell/h)\lvert\nabla A_H\rvert\big)+O\big((\ell/h)^2\lvert A_H\rvert\big),
\]

with the functions evaluated at \(Q^{\mathsf T}x/\ell\); covering and \(\ell\le h\) give \(\langle\lvert\tau-T\rvert^2\rangle\le C_Sh\lvert T\rvert^2\). The construction is linear in \(T\). \(\square\)

The stress changes from the pattern field to the constant \(T\) entirely inside the solid band, and the displacement is already affine there. Both propositions refer to the same indicator \(\chi_M\).

\paragraph*{\texorpdfstring{A.6 The interior of \(G_M\) has Lipschitz boundary}{A.6 The interior of G\_M has Lipschitz boundary}}

This section is not needed for the Theorem as stated with definition (1.1). It shows that for the sets \(G_M\) it makes no difference whether test displacements are defined on the whole cell, as in (1.1), or only on the solid.

\textbf{Lemma A.4.} For every \(M\ge11\), \(G_M\) is the closure of its interior, and the interior is an open set with Lipschitz boundary (for each fixed \(M\); by periodicity finitely many boundary charts suffice). Consequently every field \(v\) on \(\operatorname{int}G_M\) such that \(v-Ex\) is \(\Lambda\)-periodic and has finite \(H^1\) norm on the solid part of one period cell is the restriction of some \(Ex+w\), \(w\in H^1_{\rm per}\), and the infimum (1.1) is unchanged if displacements are taken in \(H^1\) of the solid only.

\emph{Proof.} \(G_M\) is a locally finite union of closed convex polygons. For such a set the claim holds as soon as at every boundary point \(x_0\) the solid coincides, near \(x_0\), with \textbf{one} closed sector with apex \(x_0\) and opening strictly between \(0\) and \(2\pi\).

\emph{(i) The pattern \(G^\square_M\) has this property.} By Lemma A.1 two pieces touch only in the following ways. (a) A rod end is welded to a straight edge of a block, collar or tab, and the arm of each face abuts the end face of the last tab; in every case the contact segment lies in the relative interior of the supporting edge (rods end at distance \(p/(4N)\) from the collar's corners and well inside the block face; connectors and arms occupy \(\lvert z-z_0\rvert\le d/32\) on a tab face \(\lvert z-z_0\rvert\le d/16\)). At an interior point of the contact segment the solid is a full disc; at its two endpoints it is a half-plane together with a wedge on the other side, one sector of opening between \(\pi\) and \(2\pi\). (b) A tab and its collar, and the two arms, overlap in a set with interior; at a crossing point of their edges the solid is the union of two half-planes or of a half-plane and a quadrant inside it, again one sector. All other pairs of pieces are at positive distance, and a boundary point of a single convex polygon sees a sector of opening at most \(\pi\).

\emph{(ii) The cut lines.} In the sector \(\Delta\) the solid is the band \(\{d_\Delta\le3h\}\) together with \(\ell G^\square_M\cap\{d_\Delta>3h\}\); near a point \(x_0\) of a cut line which is not a corner of the inner triangle this is a closed half-plane \(H^-\) together with \(\ell G^\square_M\cap H^+\). If \(x_0\) is an interior point of the pattern or lies in the open void, the solid near \(x_0\) is a disc or \(H^-\). If \(x_0\) lies on the boundary of the pattern but is not a vertex of it, the pattern is locally a half-plane whose edge either crosses the cut line, giving one sector of opening between \(\pi\) and \(2\pi\), or lies along it, giving a disc or \(H^-\). At a corner of the inner triangle the band is locally the complement of an open \(60^\circ\) wedge; the corner has an irrational coordinate, so it is not a vertex of the pattern, and the same three cases give a disc, a \(300^\circ\) sector, or (the void being the intersection of the wedge with an open half-plane) a larger single sector. So it suffices that \textbf{no vertex of the boundary of \(\ell G^\square_M\) lies on a cut line}.

Vertices of the pattern have rational coordinates. The two oblique cut lines are \(\sqrt3x_1-x_2=6h\) and \(\sqrt3(1-x_1)-x_2=6h\); a rational point on either has \(x_1\in\{0,1\}\) and \(x_2=-6h<0\), outside the sector. The third cut line is \(x_2=3h\), that is \(v=3/h=3M+\tfrac3{100}\) in pattern coordinates \((u,v)=x/\ell\). So a vertex lies on it only if its ordinate in the unit cell equals \(\tfrac3{100}\), i.e.~only if \(Y:=v/h\) equals \(3q/10^4\) with \(q=100M+1\). Since \(3q\) is odd and prime to \(5\), this number has denominator exactly \(2^4\cdot5^4\). The vertex ordinates are the following (in units of \(h\); \(k\in\{0,\dots,M-1\}\), \(j\in\{0,\dots,N-1\}\), signs independent):

{\def\LTcaptype{none} 
\begin{longtable}[]{@{}
  >{\raggedright\arraybackslash}p{(\linewidth - 4\tabcolsep) * \real{0.3333}}
  >{\raggedright\arraybackslash}p{(\linewidth - 4\tabcolsep) * \real{0.3333}}
  >{\raggedright\arraybackslash}p{(\linewidth - 4\tabcolsep) * \real{0.3333}}@{}}
\toprule\noalign{}
\begin{minipage}[b]{\linewidth}\raggedright
vertices
\end{minipage} & \begin{minipage}[b]{\linewidth}\raggedright
\(Y\)
\end{minipage} & \begin{minipage}[b]{\linewidth}\raggedright
why \(Y\ne3q/10^4\)
\end{minipage} \\
\midrule\noalign{}
\endhead
\bottomrule\noalign{}
\endlastfoot
blocks, face-2 rod ends, collars, tabs, arms & \(0\), \(M\), \(M+\{\tfrac1{300},\tfrac1{150},\tfrac1{100},\tfrac1{200}\pm\tfrac1{1600},\tfrac1{200}\pm\tfrac1{3200}\}\) & denominator not divisible by \(5^3\) \\
face-2 connector rods & \(M+\tfrac3{640}+\dfrac{\text{odd}}{12800M^2}\) & denominator divisible by \(2^9\) \\
face-1 collars & \(k+\{\tfrac14,\tfrac34\}\pm\tfrac1{24}\pm\{\tfrac1{300},\tfrac1{150}\}\) & denominator divides \(600\) \\
face-1 tabs, tab-collar crossings & \(k+\{\tfrac14,\tfrac34\}\pm\tfrac1{200}\pm(\tfrac1{24}+\tfrac1{800})\), \(\ k+\{\tfrac14,\tfrac34\}\pm\tfrac1{24}\pm(\tfrac1{200}\pm\tfrac1{1600})\) & denominator divides \(4800\) \\
face-1 first rods, last rods at the collar & \(k+\{\tfrac14,\tfrac34\}-\tfrac1{24}+\dfrac{\text{odd}}{96M^2}+\{0,\pm\tfrac1{300},\pm\tfrac1{150}\}\) & denominator divisible by \(2^5\) \\
face-1 last rods at the far block & the same with \(\pm\tfrac1{100}\) in place of the last bracket, plus \(\eta d/(3h)=\dfrac1{3q}\) & denominator divisible by \(2^5\) (\(q\) is odd) \\
\end{longtable}
}

(For the rod families the term \(\text{odd}/(96M^2)\) has \(2\)-adic denominator at least \(2^5\) while all other terms have at most \(2^3\).) Hence no vertex lies on a cut line, for every \(M\ge2\). \(\square\)

Machine check (\texttt{checks/check\_cut\_vertices.py}): for \(M=2,\dots,14\) the list above reproduces exactly the set of ordinates of all polygon vertices and all edge crossings within a body of the exact polygon model, and \(\tfrac3{100}\) is not among them; solving the equations exactly, no \(M\le200{,}000\) has a vertex on the cut line.
\phantomsection\label{endfile:frame-fields}

\phantomsection\label{file:zero-bulk-stroke}
\paragraph*{Ingredients}

\textbf{The template (from the main theorem).} Fix \(f<\lambda^2\), \(K\ge\mu\) and an error \(\eta\in(0,1]\); then take the number \(N\) of pairing levels so large that the determinant-one correction \(L_N\) of the paired grain costs less than \(\eta/2\) (energies in this paragraph are per unit area and in units of \(\mu\); the choice of \(N\) needs only the a priori bound \(\eta/2\) on the uncorrected energy), and run the finite-cell construction (the bonded-cell recovery in Section~\ref{sup:stroke}) on the unit square \(Q_1\) as root polygon, with its tolerances chosen so that the path energy before \(L_N\) is below \(\eta/2\) and the sum of the error terms of the stress test (Section D of the bonded-cell recovery proof in Section~\ref{sup:stroke}) is below \(\eta|T|^2/(2\mu)\) for all symmetric mean stresses \(T\) at once. The latter is possible because those error terms are integrals of \(\sigma_T:C_s^{-1}\sigma_T\) over sets of small area, weight errors times leaf costs, and pore corrections, and the pointwise bounds \(|\sigma_T|\le C|T|\) hold for all \(T\): each error is a small number times \(|T|^2\), uniformly on the unit sphere of \(\mathrm{Sym}\). (Section D of the bonded-cell recovery proof in Section~\ref{sup:stroke} states it in this form: the excess over the ideal cost is at most \(\varepsilon|T|^2/\mu\) for every symmetric \(T\).) This gives a closed set \(\Sigma\subset Q_1\) and a number \(0<\rho<\tfrac12\) with the following properties.

\begin{itemize}
\item
  \textbf{(T1)} \(Q_1\setminus\Sigma\) is a finite union of open convex polygons with pairwise disjoint closures, all contained in \(\{x\in Q_1:\operatorname{dist}(x,\partial Q_1)>\rho\}\), and \(\Sigma\) is connected.
\item
  \textbf{(T2)} There are Lipschitz homeomorphisms \(Y_t:Q_1\to tQ_1\), continuous in \(t\in[\lambda,1]\), with \(Y_1=\mathrm{id}\), \(Y_t(x)=tx\) on \(\partial Q_1\), \(\det\nabla Y_t=1\) a.e. in \(\Sigma\), and \(\int_\Sigma\tfrac12q(\nabla Y_t)\le\eta\).
\item
  \textbf{(T3)} For every \(T\in\mathrm{Sym}\) there is \(\sigma_T\in L^\infty(Q_1;\mathrm{Sym})\), linear in \(T\), with zero distributional divergence in the open square, \(\sigma_T=0\) a.e. outside \(\Sigma\), \(\sigma_T=T\) on \(\{\operatorname{dist}(x,\partial Q_1)<\rho\}\), and \[\int_{Q_1}|(\sigma_T)_d|^2\le\Big(\frac{\mu}{G(f,K)}+2\eta\Big)|T|^2,\qquad G(f,K)=\mu\,\frac{fK}{(2-f)K+2(1-f)\mu}.\]
\end{itemize}

For (T1), each hole is the hole of one shell, strictly inside its own filling piece, and the finished cell has ``a connected solid with Lipschitz boundary and finitely many separated polygonal holes'' (the bonded-cell recovery proof in Section~\ref{sup:stroke}, end of Section C); the finished root polygon has a solid neighbourhood of its boundary (Section D, second paragraph); its relative width defines \(\rho\). (T2): Section D with Section 4 of the finite-stroke construction; the maps are piecewise affine on a finite triangulation, with the radial shell map on finitely many pieces, followed by the constant matrix \(L_N\), hence Lipschitz. (T3): Sections C-D of the bonded-cell recovery proof in Section~\ref{sup:stroke} and Sections 2 and 5 of the finite-stroke construction, which state the cost bound for every symmetric mean stress (the construction of Sections A-D is linear in the root stress and nowhere uses that it is trace-free; the leaf stresses of the ideal layered material under mean stress \(T\) have complementary energy \(T:C_3^{-1}T=|T_d|^2/(2G)+(\operatorname{tr}T)^2/(4\kappa)\) with \(\kappa=fK\mu/((1-f)K+\mu)\); for \(K\ge\mu\) one has \(\kappa\ge G\), so this is at most \(|T|^2/(2G)\); and \(|\sigma_d|^2/(2\mu)\le\sigma:C_s^{-1}\sigma\)). The restriction \(K\ge\mu\) is a convenience: the leaf stresses are fixed by statics and do not depend on \(K\), and their deviatoric cost is \(|T_d|^2(2-f)/(2\mu f)+(\operatorname{tr}T)^2(1-f)/(4\mu f)\le(2-f)|T|^2/(2\mu f)\) for every \(K\) (\texttt{checks/check\_leaf\_stresses.py}, leaf cost identity).

For a square \(Q=x_0+sRQ_1\) (\(s>0\), \(R\) a rotation) the scaled copies are \(\Sigma_Q=x_0+sR\Sigma\), \(Y^Q_t(x)=tx_0+sR\,Y_t\big(R^{\mathsf T}(x-x_0)/s\big)\) and \(\sigma^Q_T(x)=R\,\sigma_{R^{\mathsf T}TR}\big(R^{\mathsf T}(x-x_0)/s\big)R^{\mathsf T}\). They have the same three properties on \(Q\), with rim width \(\rho s\) and the integrals multiplied by \(|Q|\).

\textbf{The frame (from the linear-frame proof in this section).} For \(M\ge11\), \(h=(M+\tfrac1{100})^{-1}\), the set \(G_M\) is closed, periodic for the lattice \(\Lambda\) with cell \(D\), invariant under rotation by \(60^\circ\) about the origin, a locally finite union of closed convex polygons, the closure of its interior, and its complement has positive area. Its form \(Q_{G_M}\) is isotropic with \(\kappa_M\le C\mu h\), and it carries

\begin{itemize}
\tightlist
\item
  \textbf{(S)} for each trace-free \(T\) a field \(\tau\in L^2_{\rm per}(D;\mathrm{Sym})\), linear in \(T\), with \(\tau=0\) a.e. outside \(G_M\), \(\langle\tau:e(w)\rangle=0\) for all \(w\in H^1_{\rm per}\), \(\langle\tau\rangle=T\) and \(\langle|\tau-T|^2\rangle\le C_Sh|T|^2\).
\end{itemize}

These are the Theorem, (G1), (G2), Lemma A.4 and Proposition S of the linear-frame proof in this section. Lemma 2.2 of the linear-frame proof in this section is also used: a stress supported on the solid and equilibrated on the torus bounds \(Q_G\) from below. Its proof uses only that the stress vanishes off the solid and that admissible displacements are solenoidal on it, so it holds for any measurable solid indicator.

\textbf{Three elementary facts.}

\begin{itemize}
\item
  \emph{(Airy.)} If \(c\) is an open convex polygon and \(\tau\in L^2(c;\mathrm{Sym})\) has zero distributional divergence in \(c\), then \(\tau=\mathcal AU\) for some \(U\in H^2(c)\). (Each row of \(\tau\) is a rotated gradient of an \(H^1\) function; symmetry makes that pair divergence-free, hence a rotated gradient of \(U\in H^2\).)
\item
  \emph{(Poincaré.)} For \(U\in H^2(c)\) let \(a\) be the affine function with \(\int_c(U-a)=0\) and \(\int_c\nabla(U-a)=0\). Since \(c\) is convex, \(\|\nabla(U-a)\|_{L^2(c)}\le\pi^{-1}\mathrm{diam}(c)\|D^2U\|_{L^2(c)}\) and \(\|U-a\|_{L^2(c)}\le\pi^{-2}\mathrm{diam}(c)^2\|D^2U\|_{L^2(c)}\). Consequently, if \(\zeta\) is \(C^2\) on \(\bar c\) with \(0\le\zeta\le1\), the field \(\mathcal A\big(\zeta(U-a)\big)\) is divergence-free, equals \(\mathcal AU\) where \(\zeta=1\), vanishes where \(\zeta=0\), and \[\big\|\mathcal A\big(\zeta(U-a)\big)\big\|_{L^2(c)}\le\Big(1+\tfrac2\pi\mathrm{diam}(c)\|\nabla\zeta\|_\infty+\tfrac1{\pi^2}\mathrm{diam}(c)^2\|D^2\zeta\|_\infty\Big)\|\mathcal AU\|_{L^2(c)}.\tag{Z.1}\] The bracket is unchanged when \(c\) and \(\zeta\) are rescaled together.
\item
  \emph{(Composition.)} For \(2\times2\) matrices \(F\) and \(A\), with \(B=F^{\mathsf T}F-I\), \[\begin{gathered}q(FA)=q(F)+q(A)+\operatorname{tr}\big(B\,(AA^{\mathsf T}-I)\big),\\ \det A=1\ \Longrightarrow\quad \big|\operatorname{tr}\big(B\,(AA^{\mathsf T}-I)\big)\big|\le\|B\|_{\rm op}\big(q(A)+2\sqrt{q(A)}\big).\end{gathered}\tag{Z.2}\] (Expand \(|FA|^2=\operatorname{tr}(F^{\mathsf T}F\,AA^{\mathsf T})\). If \(\varsigma\ge1\) and \(1/\varsigma\) are the singular values of \(A\), then \(q(A)=(\varsigma-\varsigma^{-1})^2\) and \(AA^{\mathsf T}-I\) has trace norm \(\varsigma^2-\varsigma^{-2}=\sqrt{q(A)}\sqrt{q(A)+4}\le q(A)+2\sqrt{q(A)}\).)
\end{itemize}

\textbf{Lemma F (every design has a short soft stroke).} Let \(\mathcal S\) be a closed periodic set whose complement has positive area, and \(\delta>0\). There are \(s_0\in(0,1]\), \(C_F<\infty\) and diffeomorphisms \(\Phi_s\) of the plane, \(0\le s\le s_0\), depending smoothly on \(s\), with \(\Phi_0=\mathrm{id}\), \(\Phi_s-e^{-s}x\) periodic, \(\det\nabla\Phi_s=1\) on \(\mathcal S\), \[\big\|\nabla\Phi_s^{\mathsf T}\nabla\Phi_s-I\big\|_{\rm op}\le C_Fs\ \text{ everywhere},\qquad\big\langle\chi_{\mathcal S}\,\tfrac12q(\nabla\Phi_s)\big\rangle\le\tfrac12\Big(\frac{Q_{\mathcal S}(I)}\mu+\delta\Big)s^2 .\]

\emph{Proof.} Take \(u=x+w\) admissible for \(Q_{\mathcal S}(I)\) with \(2\langle\chi_{\mathcal S}|e(u)|^2\rangle\le Q_{\mathcal S}(I)/\mu+\delta/2\). Then \(g=\operatorname{div}w\) has mean zero and \(g+2=0\) a.e. on \(\mathcal S\). The void \(\Omega=\mathbb R^2\setminus\mathcal S\) is open, so there are smooth periodic \(\gamma_i\to g+2\) in \(L^2\) with \(\operatorname{dist}(\operatorname{supp}\gamma_i,\mathcal S)>0\) (truncate to \(\{\operatorname{dist}(\cdot,\mathcal S)>1/i\}\) and mollify at half that radius), and adding multiples of one fixed bump supported in \(\Omega\) makes \(\langle\gamma_i\rangle=2\). Write \(w=w^0+\nabla\phi\) with \(\Delta\phi=g\) and \(w^0\in H^1_{\rm per}\) divergence-free, and put \(w_i=w^0*\psi_i+\nabla\phi_i\) with mollifiers \(\psi_i\) and \(\Delta\phi_i=\gamma_i-2\). Then \(w_i\) is smooth and periodic, \(w_i\to w\) in \(H^1\), and \(\operatorname{div}w_i=-2\) on a neighbourhood \(U_i=\{\operatorname{dist}(\cdot,\mathcal S)<\nu_i\}\) of \(\mathcal S\). Fix \(i\) with \(2\langle\chi_{\mathcal S}|e(x+w_i)|^2\rangle\le Q_{\mathcal S}(I)/\mu+\tfrac34\delta\).

Let \(\Theta_s\) be the flow of the periodic field \(-w_i\) and \(\Phi_s=e^{-s}\Theta_s\). Then \(\Theta_s-x\) is periodic and \(|\Theta_s(x)-x|\le s\|w_i\|_\infty\); for \(s\le\nu_i/(2\|w_i\|_\infty)\) the orbits starting in \(\mathcal S\) stay in \(U_i\), so by Liouville's formula \(\det\nabla\Theta_s=e^{2s}\) and \(\det\nabla\Phi_s=1\) on \(\mathcal S\). Uniformly in \(x\), \(\nabla\Phi_s=I-s(I+\nabla w_i)+O(s^2)\), which gives the first bound. For a \(2\times2\) matrix, \(q(F)=2|(\operatorname{sym}F)_d|^2+2(\det F-1)\). On \(\mathcal S\) the determinant is one and the strain \(e(x+w_i)\) is trace-free, so \(q(\nabla\Phi_s)=2s^2|e(x+w_i)|^2+O(s^3)\) there. The constants in the two remainders depend on the fixed index \(i\); reducing \(s_0\) once more absorbs them and gives the second bound. \(\square\)

\paragraph*{The design}

Let \(\Gamma\) be the group generated by \(\Lambda\) and the rotation by \(60^\circ\). The closed triangle \(\Delta\) with vertices \((0,0),(1,0),(\tfrac12,\tfrac{\sqrt3}2)\) is a fundamental domain: its six rotated copies form the hexagon that \(\Lambda\) tiles the plane with, and no element of \(\Gamma\) other than the identity maps \(\Delta\) to itself. The set \(\overline{\operatorname{int}G_M\cap\operatorname{int}\Delta}\) is a finite union of convex polygons; triangulate it and transport the triangles by \(\Gamma\). This gives a locally finite, \(\Gamma\)-invariant family \(\mathcal K\) of closed triangles with disjoint interiors and union \(G_M\). (The triangulation need not be conforming.)

Split each triangle of \(\mathcal K\) by the altitude onto a longest side. In a right triangle with legs \(\ell_1,\ell_2\) the square of side \(\ell_1\ell_2/(\ell_1+\ell_2)\) at the right angle has its far vertex on the hypotenuse, covers the fraction \(2\ell_1\ell_2/(\ell_1+\ell_2)^2\) of the area, and leaves two right triangles that are homothetic to the first (same shape and orientation). Repeat \(n\) times, and then cut every square into \(m^2\) equal squares. All choices are made in \(\Delta\) and transported by \(\Gamma\). Let \(\mathcal Q\) be the resulting squares and \(\mathcal R\) the residual triangles; together they are the \emph{cells}. For fixed \(M\) every residual triangle is, up to an element of \(\Gamma\), homothetic to one of the two halves of one of the finitely many triangles of \(\mathcal K\) in \(\Delta\): finitely many shapes, independent of \(n\) and \(m\). The total area of the residual triangles per period tends to zero as \(n\to\infty\).

In a residual triangle \(R\) with incentre \(c_R\) and inradius \(\iota_R\) let \(d_1,d_2,d_3\) be the distances to its three side lines; the gauge of \(R\) about \(c_R\) is \(r_R=1-\min_jd_j/\iota_R\). Let \(\mathrm{Sh}_R=\{x\in R: r_R(x)\ge r_0\}\) with \(r_0=\sqrt{1-\lambda^2/2}\): a solid shell of area fraction \(\lambda^2/2\) around a homothetic hole. Define

\[G^{\rm new}=\bigcup_{Q\in\mathcal Q}\Sigma_Q\ \cup\bigcup_{R\in\mathcal R}\mathrm{Sh}_R\ \subset G_M .\]

\(G^{\rm new}\) is closed, \(\Lambda\)-periodic and invariant under rotation by \(60^\circ\), and its void has positive area. Hence \(Q_{G^{\rm new}}(E)=\kappa_{\rm new}(\operatorname{tr}E)^2+2\mu_{\rm new}|E_d|^2\) (the linear-frame proof in this section, Lemma 1.1 and the symmetry argument of its Section 4).

\paragraph*{Proof}

\textbf{Bulk.} The cell energy charges strain on the solid only, and \(G^{\rm new}\subset G_M\). A field admissible for \(G_M\) is admissible for \(G^{\rm new}\) with no larger energy, so \(Q_{G^{\rm new}}(I)\le Q_{G_M}(I)\), that is \(\kappa_{\rm new}\le\kappa_M\le C\mu h\). No property of the filling is used: every sub-solid of \(G_M\) has bulk modulus at most \(\kappa_M\).

\textbf{Shear.} Fix one trace-free \(T\) with \(|T|=1\) (by isotropy this determines \(\mu_{\rm new}\)) and let \(\tau\) be the field of (S). On each cell \(c\), \(\tau=\mathcal AU_c\) with \(U_c\in H^2(c)\).

\begin{itemize}
\item
  On a square \(Q\) of side \(s\) let \(\bar\tau_Q\in\mathrm{Sym}\) be the mean of \(\tau\) over \(Q\), let \(V_Q=U_Q-\tfrac12x\cdot J^{\mathsf T}\bar\tau_QJx\), so that \(\mathcal AV_Q=\tau-\bar\tau_Q\), and let \(a_Q\) be its affine fit. With \(d_1,\dots,d_4\) the distances to the four side lines of \(Q\) put \[\zeta_Q=1-\prod_{j=1}^4\vartheta\Big(\frac{d_j-\rho s/3}{\rho s/3}\Big).\] Then \(\zeta_Q\) is smooth, \(0\le\zeta_Q\le1\), \(\zeta_Q=1\) within distance \(\rho s/3\) of \(\partial Q\), \(\zeta_Q=0\) beyond distance \(2\rho s/3\), \(|\nabla\zeta_Q|\le C_\vartheta/(\rho s)\) and \(|D^2\zeta_Q|\le C_\vartheta/(\rho s)^2\), with \(C_\vartheta\) depending only on \(\vartheta\). Put \[\sigma=\sigma^Q_{\bar\tau_Q}+\mathcal A\big(\zeta_Q(V_Q-a_Q)\big)\quad\text{on }Q .\] The second term is supported in the rim, which is solid by (T1), and there \(\sigma^Q_{\bar\tau_Q}=\bar\tau_Q\) by (T3). So \(\sigma=0\) a.e. outside \(\Sigma_Q\) and \(\sigma=\tau\) within distance \(\rho s/3\) of \(\partial Q\). By (Z.1), with \(C_\rho=1+C_\vartheta/\rho+C_\vartheta/\rho^2\), the second term has norm at most \(C_\rho\|\tau-\bar\tau_Q\|_{L^2(Q)}\).
\item
  On a residual triangle \(R\) put \(\zeta_R=1-\prod_{j=1}^3\vartheta\big((d_j-b_R)/b_R\big)\) with \(b_R=\iota_R(1-r_0)/4\), and \(\sigma=\mathcal A\big(\zeta_R(U_R-a_R)\big)\). Here \(\zeta_R=1\) near \(\partial R\) and \(\zeta_R=0\) on \(\{r_R\le(1+r_0)/2\}\), so \(\sigma\) vanishes off \(\mathrm{Sh}_R\) and equals \(\tau\) near \(\partial R\). The cut-off is transported by homotheties and by \(\Gamma\), so by (Z.1) \(\|\sigma\|_{L^2(R)}\le C_M\|\tau\|_{L^2(R)}\) with \(C_M\) depending only on the finitely many shapes and on \(\lambda\).
\item
  Outside \(G_M\) put \(\sigma=0=\tau\).
\end{itemize}

On every cell \(\sigma-\tau=\mathcal AW_c\) with \(W_c\in H^2(c)\) vanishing near \(\partial c\). (On a square, \(\sigma-\tau=(\sigma^Q_{\bar\tau_Q}-\bar\tau_Q)+\mathcal A((\zeta_Q-1)(V_Q-a_Q))\); the first bracket is divergence-free and vanishes on the connected rim, so it is \(\mathcal A\) of an \(H^2\) function that is affine on the rim, and the affine function can be removed.) Hence \(\int_c(\sigma-\tau):e(v)=0\) for every \(v\in H^1(c)\). Summing over a set of representatives of the cells modulo \(\Lambda\), for every \(E\in\mathrm{Sym}\) and \(w\in H^1_{\rm per}\), \[\langle\sigma:(E+e(w))\rangle=\langle\tau:(E+e(w))\rangle=T:E .\] So \(\sigma\) vanishes off \(G^{\rm new}\), is equilibrated on the torus and has mean \(T\). By (T3) and Jensen's inequality, \(\sum_Q\int_Q|(\sigma^Q_{\bar\tau_Q})_d|^2\le(\mu/G+2\eta)\sum_Q|Q||\bar\tau_Q|^2\le(\mu/G+2\eta)\,|D|\,\langle|\tau|^2\rangle\), and \(\langle|\tau|^2\rangle^{1/2}\le1+(C_Sh)^{1/2}\) by (S). Therefore \[\begin{gathered}\langle|\sigma_d|^2\rangle^{1/2}\le\Big(\frac\mu{G(f,K)}+2\eta\Big)^{1/2}\big(1+(C_Sh)^{1/2}\big)\\ {}+C_\rho\Big(\frac1{|D|}\sum_{Q}\|\tau-\bar\tau_Q\|^2_{L^2(Q)}\Big)^{1/2}+C_M\Big(\frac1{|D|}\int_{D\cap\bigcup\mathcal R}|\tau|^2\Big)^{1/2}.\end{gathered}\] For fixed \(M\) the last term tends to zero as \(n\to\infty\), whatever the template and \(m\) (absolute continuity of the integral). For fixed \(n\) and a fixed template the middle term tends to zero as \(m\to\infty\) (means over squares of side tending to zero converge in \(L^2\)). Lemma 2.2 of the linear-frame proof in this section with \(E=T\) reads \(1\le\langle|\sigma_d|^2\rangle\,\mu_{\rm new}/\mu\), so \[\frac{\mu_{\rm new}}{\mu}\ \ge\ \Big[\Big(\frac{\mu}{G(f,K)}+2\eta\Big)^{1/2}\big(1+(C_Sh)^{1/2}\big)+o\Big]^{-2},\tag{Z.3}\] where \(o\) is the sum of the last two terms. The constant \(C_\rho\), which depends on the template, multiplies a term that is sent to zero last; this is what allows the template to be chosen after the frame.

\textbf{Path.} Apply Lemma F to \(\mathcal S=G_M\) with a \(\delta\) to be chosen; it gives \(s_0\) (take also \(\lambda e^{s_0}<1\)), \(C_F\) and \(\Phi_s\). Next let \(y^{\rm c}_\theta\), \(\lambda\le\theta\le1\), be the cell-wise motion: \(y^{\rm c}_\theta=Y^Q_\theta\) on each square, \(y^{\rm c}_\theta(x)=\theta x\) outside \(G_M\), and on each residual triangle the radial shell map with entering map \(\theta x\), \[y^{\rm c}_\theta(x)=\theta\Big(c_R+\frac{g(r)}{r}(x-c_R)\Big),\qquad r=r_R(x),\quad g(r)=\sqrt{1+(r^2-1)/\theta^2}\ \ (r\ge r_0),\] continued linearly in \(r\) inside the hole. On the boundary of every cell \(y^{\rm c}_\theta(x)=\theta x\), so \(y^{\rm c}_\theta\) is a Lipschitz homeomorphism of the plane that maps each cell onto its homothetic image, with \(y^{\rm c}_\theta-\theta x\) periodic and \(y^{\rm c}_1=\mathrm{id}\). On \(\Sigma_Q\) its determinant is one by (T2); on \(\mathrm{Sh}_R\) it is \(\theta^2gg'/r=1\), and since \(g(r_0)^2=1-\lambda^2/(2\theta^2)\ge\tfrac12\) the gradient there is bounded by a constant depending on the shape of \(R\) and on \(\lambda\). So its energy is \[e^{\rm c}(\theta):=\big\langle\chi_{\rm new}\tfrac12q(\nabla y^{\rm c}_\theta)\big\rangle\le\eta+C'_M\frac{|D\cap\bigcup\mathcal R|}{|D|},\] with \(C'_M\) depending only on the finitely many shapes and on \(\lambda\). The map \(\Xi_\theta=y^{\rm c}_\theta/\theta\) is a Lipschitz homeomorphism of the plane that maps every cell onto itself and is the identity outside \(G_M\).

Define the path in two legs: \[y_t=\Phi_s\ \ (t=e^{-s},\ 0\le s\le s_0),\qquad y_t=\theta\,\Phi_{s_0}\circ\Xi_\theta\ \ (t=\theta e^{-s_0},\ \lambda e^{s_0}\le\theta\le1).\] The legs agree at \(t=e^{-s_0}\), each \(y_t\) is a Lipschitz homeomorphism, the path is continuous, \(y_1=\mathrm{id}\), and \(y_t-tx\) is periodic. On the first leg \(\det\nabla y_t=1\) on \(G_M\supset G^{\rm new}\). On the second leg \(\nabla y_t(x)=F\,A\) a.e., with \(F=\nabla\Phi_{s_0}(\Xi_\theta x)\) and \(A=\nabla y^{\rm c}_\theta(x)\); for \(x\) in the solid, \(\det A=1\), and \(\Xi_\theta x\) lies in the same cell, hence in \(G_M\), so \(\det F=1\).

\emph{Energy on the first leg.} Since \(\det\nabla\Phi_s=1\) on all of \(G_M\), \(q(\nabla\Phi_s)\ge0\) there, so \(\mathcal E(t)/\mu=\langle\chi_{\rm new}\tfrac12q(\nabla\Phi_s)\rangle\le\langle\chi_M\tfrac12q(\nabla\Phi_s)\rangle\le\tfrac12(4\kappa_M/\mu+\delta)s^2\), and \(s\le2(1-e^{-s})=2(1-t)\) for \(s\le1\). So \(\mathcal E(t)\le2(4\kappa_M+\delta\mu)(1-t)^2\).

\emph{Energy on the second leg.} By (Z.2) and the first bound of Lemma F, \(q(FA)\le q(F)+q(A)+C_Fs_0\big(q(A)+2\sqrt{q(A)}\big)\). On each cell \(c\) the Lipschitz map \(\Xi_\theta\) is injective on the solid, with Jacobian \(\theta^{-2}\ge1\) and values in \(c\subset G_M\), where \(q(\nabla\Phi_{s_0})\ge0\); by the area formula \(\int_{G^{\rm new}\cap c}q(F)\,dx\le\int_cq(\nabla\Phi_{s_0})\). With the Cauchy-Schwarz inequality for \(\langle\chi_{\rm new}\sqrt{q(A)}\rangle\), \[\frac{\mathcal E(t)}\mu\le\tfrac12\Big(\frac{4\kappa_M}\mu+\delta\Big)s_0^2+e^{\rm c}(\theta)+C_Fs_0\Big(e^{\rm c}(\theta)+\sqrt{2e^{\rm c}(\theta)}\Big),\qquad 1-t\ge1-e^{-s_0}\ge\frac{s_0}2 .\]

\textbf{Order of the parameters.} Write \(S_\lambda=S(\lambda)=\lambda^2/(2-\lambda^2)\). All three conclusions of the corollary are monotone in \(\varepsilon\), so it is enough to treat \(\varepsilon<S_\lambda\).

\begin{enumerate}
\def\labelenumi{\arabic{enumi}.}
\tightlist
\item
  choose \(f<\lambda^2\) and \(K\ge\mu\) with \(G(f,K)/\mu>S_\lambda-\varepsilon/4\);
\item
  choose \(M\) with \(\kappa_M\le\varepsilon\mu/64\) and \(\big(\mu/G(f,K)\big)\big(1+(C_Sh)^{1/2}\big)^2\le(S_\lambda-\varepsilon/2)^{-1}\). This fixes thresholds \(\eta_1>0\) and \(o_1>0\), depending on \((\lambda,f,K,M,\varepsilon)\) only, such that the right side of (Z.3) is at least \(S_\lambda-\varepsilon\) whenever \(\eta\le\eta_1\) and \(o\le o_1\). Put \(\delta=\varepsilon/16\) and take \(s_0\), \(C_F\), \(\Phi_s\) from Lemma F for \(G_M\). Then on the first leg \(\mathcal E(t)\le\tfrac14\varepsilon\mu(1-t)^2\), and the first term of the second-leg bound is at most \(\tfrac1{16}\varepsilon s_0^2\le\tfrac14\varepsilon(1-t)^2\);
\item
  choose \(e_0>0\) with \(e_0+C_Fs_0(e_0+\sqrt{2e_0})\le\tfrac18\varepsilon s_0^2\), then \(\eta\le\min(\eta_1,e_0/2)\), then \(N\) and the template; this fixes \(\rho\) and \(C_\rho\);
\item
  choose \(n\) so large that \(C'_M|D\cap\bigcup\mathcal R|/|D|\le e_0/2\) and the last term of the shear estimate is at most \(o_1/2\) (neither condition involves the template), and then \(m\) so large that the middle term is at most \(o_1/2\).
\end{enumerate}

Then \(e^{\rm c}(\theta)\le e_0\), so the last three terms of the second-leg bound are at most \(\tfrac18\varepsilon s_0^2\le\tfrac12\varepsilon(1-t)^2\), and \(\mathcal E(t)\le\varepsilon\mu(1-t)^2\) on both legs; \(\kappa_{\rm new}\le\varepsilon\mu\); and \(\mu_{\rm new}/\mu\ge S_\lambda-\varepsilon\). The constants \(s_0\) and \(C_F\) are fixed before the template exists because Lemma F is applied to \(G_M\), not to \(G^{\rm new}\). \(\blacksquare\)
\phantomsection\label{endfile:zero-bulk-stroke}

\SSection{rubber}{Mixed-load supremum and endpoint paths}
Use the phase laws and energy units of Theorem~\ref{new:rubber}, with
phase moduli $1$ and $10$ and equal fractions. Their arithmetic and
harmonic means are $m=11/2$ and $20/11$ respectively.
\phantomsection\label{file:rubber}
\paragraph*{A bound for every three-dimensional layout}

Let \(b=\lambda^{-1/2}\). The ordered singular values of the mean gradient \(M=(I+\gamma e_1\otimes e_3)D_\lambda\) are \(x,b,z\). The \(1,3\) block has a column of length \(b\) and determinant \(1/b\), which gives this ordering for \(\lambda\le1\). In particular

\[
xz=1/b,\qquad d=x-z=\sqrt{(b-\lambda)^2+\lambda^2\gamma^2}.
\tag{1}
\]

In a proper singular value decomposition write \[M=U\operatorname{diag}(x,b,z)V^T,\qquad R=U\operatorname{diag}(1,1,-1)V^T.\] This improper orthogonal matrix supplies a linear test of both the gradients and their cofactors. For \(X=R^TF\) the exact identity

\[
|F|^2-2R:\operatorname{cof}F
  =(\operatorname{tr}X)^2+\sum_{i<j}(X_{ij}-X_{ji})^2
\tag{2}
\]

follows from \(R:\operatorname{cof}F=-\operatorname{tr}(\operatorname{cof}X)\) and the trace formula for second minors. Affine periodicity gives \(\langle F\rangle=M\) and \(\langle\operatorname{cof}F\rangle=\operatorname{cof}M\). For the latter, expand each two-by-two minor of \(M+\nabla w\). The linear fluctuation averages to zero and the quadratic minor is a periodic divergence. Approximation in \(W^{1,2}\) gives the identity for finite-energy gradients.

Fix \(u>1\). By (2), the excess of the phase energy above \(R:\operatorname{cof}F+uR:F-(3+u^2)/2\) is at least

\[
\frac12 (R:F-u)^2+\frac{k-1}{2}(|F|^2-3).
\tag{3}
\]

In the soft phase this is nonnegative. For the hard phase let the ordered local stretches be \(a\ge s\ge c>0\), put \(D_0=a-c\), and use \(asc=1\). Then

\[
|F|^2-3=D_0^2+s^2+2/s-3,\qquad R:F\le a+s-c=D_0+s .
\]

For the trace inequality, write \(R:F\) as the stretch-weighted sum of the diagonal entries of an improper orthogonal matrix \(Q\). Each entry is at most one and \(3-\operatorname{tr}Q\ge2\), so the loss from \(a+s+c\) is at least \(2c\). Thus \((R:F-u)^2\ge(u-s-D_0)_+^2\), where \(w_+=\max(w,0)\). Minimizing \((u-s-D_0)_+^2/2+9D_0^2/2\) over \(D_0\ge0\) gives \(D_0=(u-s)_+/10\). The resulting hard-phase bound is

\[
H(u)=\min_{s>0}\left\{\frac92(s^2+2/s-3)+\frac9{20}(u-s)_+^2\right\}.
\tag{4}
\]

The objective is strictly convex. Its minimizer \(r\in(1,u)\) satisfies \(u=11r-10r^{-2}\). The hard phase can give equality with stretches \(r,r,r^{-2}\). Averaging (3) and using the equal fractions now gives

\[
\langle W\rangle\ \ge\ b^{-1}-bd+u(b+d)-\frac{3+u^2}{2}+\frac{H(u)}2 .
\tag{5}
\]

Here \(R:M=b+d\) and \(R:\operatorname{cof}M=b^{-1}-bd\).

To specify the answer algebraically, set

\[
\tau=2r-r^{-2},\qquad u=11r-10r^{-2},\qquad
d(r)=\frac{u+\tau}{2}-b .
\tag{6}
\]

At this choice the right side of (5) is stationary in \(u\) with \(d\) held fixed. Indeed \(H'(u)=9(u-r)/10\). Substitution gives

\[
L(r)=b^{-1}-bd(r)-\frac32
      +\frac{u^2+\tau^2}{4}
      +\frac94(2r^2+r^{-4}-3).
\tag{7}
\]

At \(\lambda=9/10\), take the root \(r_*\in(1.03278,1.03280)\) of \(L(r_*)=1/4\). It is unique there because \(L'(r)=(u-b)d'(r)>0\). Equations (1) and (6) give

\[
r_*=1.03278619490856272\ldots,\qquad
\Gamma=\frac{\sqrt{d(r_*)^2-(b-\lambda)^2}}{\lambda}.
\tag{8}
\]

At this root \(u\simeq1.9854781405>b\). Hold that test value fixed in (5). Its right side increases strictly with \(d\), proving exclusion of every \(\gamma>\Gamma\). The root and the shear have rational enclosures.

\paragraph*{Why sharpness survives the loading-path requirement}

\paragraph*{An endpoint mixture with the energy in (7)}

Work in singular coordinates with a reflection on the output side, so that the mean takes the symmetric form \(J=\operatorname{diag}(x,b,-z)\). Returning to physical coordinates sends a matrix \(A\) to \(F=U\operatorname{diag}(1,1,-1) A V^T\). Thus \(\det A=-1\) gives a physical gradient of determinant one. For symmetric \(A\), the skew terms in (2) vanish since \(R^T F=VAV^T\). A soft state then attains its bound in (3) if \(\operatorname{tr} A=u\). Any rotation of the hard state \(\operatorname{diag}(r,r,-r^{-2})\) by orthogonal conjugation attains (4) and has trace \(\tau\).

\textbf{Two lamination levels suffice.} For \(r\) close to \(r_*\), use the load from (1) and (6). Then \(x\ge b>r\) and \(z\le\lambda=b^{-2}<r^{-2}\). Choose the hard state

\[
c=r+r^{-2},\qquad K=rI-c\,e e^T,\qquad
e=(0,\sqrt{t},\sqrt{1-t})^T,\qquad
t=\frac{(b-r)(r^{-2}-z)}{c(b+z)} .
\]

Here \(0<t<1\). The first diagonal entry of \(B=J-K\) is \(x-r>0\). Its \(2,3\) block has a positive first diagonal entry and determinant \((b-r)(r^{-2}-z)-ct(b+z)=0\). Hence \(B=GG^T\) for a two-column matrix \(G\) whose first column is \(\sqrt{x-r}\,e_1\) and whose second lies in the \(2,3\) plane.

The two-by-two matrix \(T=I_2+G^TK^{-1}G\) is diagonal with first entry \(x/r\). Its determinant is \(\det J/\det K=1\) by the determinant formula for a low-rank update, so its other entry is \(r/x\). Choose a unit vector \(v\) with \(v^TTv=1\), between these two eigenvalues, and a perpendicular unit vector \(w\). Set

\[
N_1=(Gv)(Gv)^T,\qquad N_2=(Gw)(Gw)^T.
\]

These positive rank-one increments sum to \(B\). The update formula gives \(\det(K+N_1)=\det K=-1\), and \(K+N_1+N_2=J\) also has determinant \(-1\).

Starting at \(K_0=K\), obtain \(K_j=K_{j-1}+N_j\) at each of the two levels by mixing the preceding state with the soft state

\[
S_j=K_{j-1}+\frac{u-\operatorname{tr}K_{j-1}}{\operatorname{tr}N_j}N_j.
\]

The whole line has determinant \(-1\) by rank-one affinity of the determinant. The soft fraction in this split is \(\operatorname{tr}N_j/(u-\operatorname{tr}K_{j-1})\in(0,1)\), since the traces increase up to \(\operatorname{tr}J=(u+\tau)/2<u\). The two fractions retained from the preceding state multiply to

\[
\frac{u-\operatorname{tr}J}{u-\tau}=\frac12.
\]

This is a finite sequential laminate with exactly the desired weights at the level of its gradient states. Every leaf gives equality in the bound. Rank-one mixing preserves the means of all the second minors, so its weighted energy is exactly \(L(r)\).

\paragraph*{Periodic cells with exact incompressibility}

Here is a realization to arbitrary energy accuracy with a single fixed material layout. It also explains why the orientation of a cubic lattice poses no restriction. Start with the affine endpoint map \(Mx\) on a cube. For any split of a physical affine parent gradient \(A\), write the two children as \(A(I+\alpha v\otimes n)\) and \(A(I+\beta v\otimes n)\), with \(|n|=1\). The two slopes have weighted mean zero, and \(n\cdot v=0\) by their equal determinants.

Inside a cube in the affine parent region choose a smooth cutoff \(\chi\) equal to one on its interior apart from a thin collar. Let \(f\) be smooth and one-periodic, with \(f''\) the periodic mollification of the mean-zero two-slope profile. Its plateaus have the exact slopes apart from transition intervals of arbitrarily small total length. Use the unit-time flow of the vector field

\[
(v\otimes n-n\otimes v)\nabla\big[N^{-2}\chi(x) f(Nn\cdot x)\big].
\tag{9}
\]

The field is compactly supported and divergence-free by skew symmetry. Its flow is a volume-preserving diffeomorphism, equal to the identity near the cube boundary. Its displacement is \(O(1/N)\). In the interior, away from the collar by more than that displacement, its time-one map is exactly \(x+v f'(Nn\cdot x)/N\), since \(n\cdot x\) is constant on each trajectory. Precomposing the parent map with this flow gives the desired gradients on the plateaus. The derivative of the vector field is bounded independently of collar thickness and transition width once \(N\) is sufficiently large for the cutoff. This follows by differentiating (9), which uses only \(f,f',f''\). Hence the flow gradients have a uniform bound as well.

Implement the tree from its mean outwards, covering all but an arbitrarily small volume of each open parent plateau by finitely many disjoint interior cubes. In each cube the plateau fractions converge to the prescribed fractions up to the chosen transition width as \(N\to\infty\). There are two levels, so the final gradients stay bounded independently of the exceptional volume, and the energy tends to the weighted leaf energy. The reference interiors retained in this argument can be smaller cubes intersected with planar plateau strips, giving finitely many polyhedra at each level. Mark the bulk hard leaves and assign the remaining material to the soft phase. If a fraction correction is needed, a plane can cut the required volume from the surplus phase. This gives \textbf{exactly half of each phase}, preserves finite polyhedral labels and changes the energy by a quantity tending to zero. These are design choices made before loading.

The resulting endpoint is \(M\Phi(x)\), where \(\Phi\) is a finite composition of volume-preserving flows and equals the identity near the cube faces. The time-\(s\) versions of the same flows give a continuous family \(\Phi_s\) from the identity to \(\Phi\), with bounded gradients. The compact supports permit periodic repetition. Choose \(r<r_*\) sufficiently close to the root. Since \(L(r)<1/4\), a finite choice of cubes and frequencies gives an endpoint energy \(E_1<1/4\).

\paragraph*{An explicit two-level path to 0.53163717}

For a concrete near-optimal load, both levels can deform simultaneously in every cell. Set \(g=0.53163717\), \(M_s=(I+sg\,e_1\otimes e_3)D_{0.9}\), and use the following reference data. The normals in this table are unnormalized. In each row set \(v_{i3}=-(n_{i1}v_{i1}+n_{i2}v_{i2})/n_{i3}\), so \(n_i\cdot v_i=0\).

\begin{center}\small
\begin{tabular}{@{}llll@{}}\toprule
Level & Normal $n_i$ & $10^{14}(v_{i1},v_{i2})$ & Soft fraction $w_i$\\\midrule
1, inner & $(138382,228323,1000000)$ & $(80885215,17694550)$ & $252483/1000000$\\
2, outer & $(138297,-233434,1000000)$ & $(61576076,-13220643)$ & $1-[2(1-w_1)]^{-1}$\\
\bottomrule\end{tabular}\end{center}

Starting with \(A=M_s\) at level 2, split a parent \(A\) into soft and retained gradients

\[
A[I+s(1-w_i)v_i\otimes n_i],\qquad A[I-s w_i v_i\otimes n_i].
\]

Apply level 1 to the retained gradient, leaving hard fraction \((1-w_1)(1-w_2)=1/2\). All factors preserve volume. The weighted energy of these three states is a degree-six polynomial in \(s\) with Bernstein coefficients

\[\begin{gathered}
0.0886111111,\quad 0.0712496631,\quad 0.0719766929,\quad
0.0903565614,\\ 0.1261453345,\quad 0.1793019690,\quad 0.2499999975
\end{gathered}\]

to ten decimal places. Exact arithmetic puts the largest coefficient below \(0.2499999976\), so the degree-six energy stays strictly below \(1/4\) throughout. Realize the two splits using the time-\(s\) flows in (9), with unit normals \(n_i/|n_i|\) and plateau slopes multiplied by \(|n_i|\). The plateau gradients have exactly the displayed dependence on \(s\). The plateau strips are fixed in reference coordinates, and a single margin for the maximum displacement over \(0\le s\le1\) keeps the strips used for refinement inside the affine parent regions. Both levels therefore have fixed labels and uniformly bounded gradients throughout. Choose the exceptional volumes and the phase-fraction correction small enough to fit the strict energy margin. Uniform compression followed by this shear sweep gives a bonded, injective periodic realization at \(0.53163717\).

\paragraph*{A continuous compression-then-shear path}

Subdivide the unit cube into \(N^3\) equal cubes and repeat that same endpoint cell and its material labels in each. First compress uniformly. This leg peaks at \(E_c=319/3600<1/4\).

For the shear leg, activate one horizontal row of cubes at a time. When \(j\) rows are complete and the next row is at fraction \(s\in[0,1]\), set \(a(x_3)\) to one, \(s\), or zero on the completed, current, or waiting rows, respectively, and extend periodically. For target shear \(g\) the base map is

\[
P(x)=D_\lambda x+\lambda g\, e_1\int_0^{x_3}a(t)\,dt .
\]

Precompose \(P\) in each small cube with a scaled copy of \(\Phi\), \(\Phi_s\), or the identity according to its row. Each corrector maps its cube onto itself, matches its neighbors on the faces, and preserves volume. The base map is globally injective with determinant one. The composition is therefore bonded and globally injective, with exactly the prescribed affine-periodic mean shear \((j+s)g/N\) at axial stretch \(\lambda\). The deformations agree at the join between consecutive rows and vary continuously throughout.

Each completed row costs \(E_1\) per unit of its own volume, and each waiting row costs \(E_c\). The current row has a finite bound \(C\) uniform in \(s\), since the chosen cell isotopy has bounded gradients. The full-cell energy is at most \(\max(E_c,E_1)+C/N\). Taking \(N\) large puts the entire path within budget. This path changes the deformation of the rows while their phase labels stay fixed. As \(r\uparrow r_*\) these achievable shears increase to (8). Truncating such a path gives every smaller shear.

Equations (5)-(8) also apply to other budgets in this regime. For \(0<\lambda\le1\) and dimensionless budget \(e\), the same sharp-supremum argument works at a root \(L(r)=e\) with \(1<r<b\), \(u>b\), \(d(r)>b-\lambda\), and uniform compression cost \(m(2/\lambda+\lambda^2-3)/2<e\). At the benchmark budget this last condition reaches equality at \(2/\lambda+\lambda^2=34/11\), or \(\lambda\simeq0.836279\). At stronger compression the endpoint bound persists, and proving reachability requires a separate compression path on the chosen layout.
\phantomsection\label{endfile:rubber}

\subsection{Exact attainment}\label{route:exact-attainment}
\phantomsection\label{file:rubber-endpoint}
\paragraph*{Obstruction for regular interfaces}

Use the equality conditions in (2)-(7) of the benchmark proof. In the principal coordinates, after the fixed improper output rotation, the physical gradient becomes \(A=D^2\phi\). Indeed the transformed deformation has symmetric weak gradient, so on the periodic lift to \(\mathbb R^3\) it is a gradient. Write

\[
\langle A\rangle=J=\operatorname{diag}(x,b,-z),\qquad
\det A=\det J=-1.
\]

In the hard phase \(A\) has eigenvalues \(r,r,-r^{-2}\) and trace \(\tau=2r-r^{-2}\). In the soft phase \(\operatorname{tr}A=u\). Here \(u>\tau\) and \(x\ge b>r>1\). Finite energy gives \(\phi\in W^{2,2}_{\rm loc}\).

On the reference torus, take \(H,S\) to be the open sets of points having a neighborhood occupied almost everywhere by the hard or soft phase, respectively. Define the \textbf{full interface} as \(\Sigma=\mathbb T^3\setminus(H\cup S)\). Call a point regular if locally \(\Sigma\) is an embedded \(C^{2,\alpha}\) surface for some \(0<\alpha<1\), with one hard and one soft side. This definition also treats layouts in which a phase has an empty interior. Using only \(\partial H\cap\partial S\) could miss an entire region containing both materials.

Suppose the nonregular part of \(\Sigma\) has zero \(\mathcal H^2\) measure. The regular part has a countable cover by surface patches, hence \(\Sigma\) has zero volume. The open sets \(H,S\) then represent both phases up to null sets.

\textbf{Hard components have constant Hessian.} On each hard component \(\Delta\phi=\tau\), so all entries of \(A\) are smooth and harmonic. Their squared norm is the constant \(2r^2+r^{-4}\). Applying the Laplacian to this norm gives \(2|\nabla A|^2=0\). Denote the resulting constant matrix by \(K\) and the hard quadratic potential by \(Q_K\), including its affine part.

\textbf{Each regular interface patch is planar.} On its soft side set \(\omega=\phi-Q_K\) and \(q=u-\tau>0\). Interior Poisson regularity for \(\Delta\phi\in L^\infty\) gives continuity of \(\nabla\phi\) across the patch. Thus \(\Delta\omega=q\) on the soft side, and \(\omega=\nabla\omega=0\) on the patch by the values on the hard side. Each \(\partial_i\omega\) is harmonic on the soft side with zero Dirichlet data. Local Dirichlet boundary regularity for the Laplacian on a \(C^{2,\alpha}\) graph makes each of these derivatives \(C^{2,\alpha}\) up to a smaller patch. This supplies the one-sided third derivatives below.

Let \(n\) be a unit normal. Differentiating \(\nabla\omega=0\) tangentially and using \(\Delta\omega=q\) gives \(D^2\omega=qnn^T\) on the soft face of the interface. With \(B=K^{-1}\) and \(T=K+qnn^T\), the determinant condition \(\det T=\det K\) gives \(n^TBn=0\). Thus \(v=Bn\) is a nonzero tangent vector. Let \(L\) be the symmetric shape operator \(L t=\partial_t n\) on tangents, extended by \(Ln=0\). Differentiating \(n^TBn=0\) tangentially gives \(Lv=0\).

Tangential differentiation of \(D^2\omega=qnn^T\) determines the tangential block of \(D^3\omega[n]\) to be \(qL\) by symmetry of the third derivatives. Its mixed entries vanish, and its last diagonal entry follows from \(\partial_n\Delta\omega=0\). Hence

\[
D^3\omega[n]=q\big(L-(\operatorname{tr}L)nn^T\big).
\]

Differentiate the constant soft determinant in the normal direction. Since \(T^{-1}=B-qvv^T\), \(Lv=0\), and \(n^TT^{-1}n=0\), the result is \(\operatorname{tr}(BL)=0\). If \(L\ne0\), then \(L=\kappa t t^T\) for a unit tangent \(t\perp v\) and some \(\kappa\ne0\). It follows that \(t^TBt=0\) and \(t^TBn=0\). Together with \(n^TBn=0\), this makes \(\operatorname{span}(t,n)\) a two-dimensional totally isotropic subspace of \(B\). Its signature is \((2,1)\), so projection of such a subspace onto the one-dimensional negative eigenspace must be injective. Hence \(L=0\) throughout the patch.

\textbf{The translated determinant gives a contradiction.} Every connected soft component in the lift has a regular boundary patch. To see this, take balls of equal positive radius inside that component and inside an open hard region. Parallel segments joining matching points of their central cross-sectional disks must exit the soft component. The projections of these boundary points cover a disk of positive area, so they cannot all belong to the \(\mathcal H^2\)-null exceptional set. At a regular point of the component's boundary its local soft side belongs to that component.

On the resulting plane patch through \(a\), the harmonic function \(\phi-Q_K-q[n\cdot(x-a)]^2/2\) has both zero value and zero normal derivative. It extends harmonically by zero across the patch and vanishes by harmonic continuation throughout the connected soft component. Thus the soft Hessian there is the constant \(K+qnn^T\). Since \(K-rI\) has rank one, both hard and soft Hessians satisfy \(\det(A-rI)=0\) almost everywhere.

The translated determinant has the required cell average already at finite energy. Expand it as \(\det A-r\operatorname{tr}(\operatorname{cof}A)+r^2\operatorname{tr}A-r^3\). The first term is the known constant \(-1\), and affine periodicity gives the mean identities for the entries and second minors in \(W^{1,2}\). Consequently

\[
0=\langle\det(A-rI)\rangle
 =\det(J-rI)=(x-r)(b-r)(-z-r)<0 .
\]

This proves the obstruction, including for ordinary piecewise-\(C^3\) grains with edges and junctions on a set of zero area. The measure statement concerns the full interface \(\Sigma\), allowing arbitrary measurable reference labels.

In particular, \(\Gamma\) is the \textbf{sharp but unattained supremum over finite polyhedral layouts}. In the subcritical construction the bulk reference regions are intersections of cubes with planar strips. Only finitely many cubes, strips and bands are needed for each target below \(\Gamma\), and an exact phase-fraction correction can be taken by a planar cut. The resulting layouts belong to the regular class just excluded at the endpoint.

\paragraph*{Open grains with arbitrarily rough interfaces}

Suppose the full interface \(\Sigma\) has zero volume, so \(H\) and \(S\) represent the two phases almost everywhere. Call a connected component of \(H\) in the periodic lift to \(\mathbb R^3\) a hard grain. \textbf{If every connected component of \(S\) in the lift has at least one boundary point with a neighborhood meeting only finitely many hard grains, equality is impossible.} This local convention gives one quadratic, including its affine part, on each counted grain. A grain winding around the quotient torus can have different affine parts on distinct lifts.

\textbf{A patch with matching hard data.} Let \(D\) be such a soft component and let \(B\) be a ball about the given boundary point. Every point of \(B\cap\partial D\) belongs to the closure of a hard grain meeting \(B\). In fact, a point of \(\partial D\) lies outside the maximal open soft set. Every neighborhood therefore contains hard material of positive volume and meets \(H\). Finiteness in \(B\) gives the assertion about an individual grain closure. A finite closed cover of \(B\cap\partial D\) has a member with nonempty relative interior. Hence there is a smaller ball \(B'\) centered at \(p_0\in\partial D\) such that \(\partial D\cap B'\subset\overline{H_j}\) for one hard grain \(H_j\). On that grain \(\phi=Q\) for one quadratic with Hessian \(K\). Since \(\nabla\phi\) is continuous, both \(\phi-Q\) and \(\nabla(\phi-Q)\) vanish on \(\partial D\cap B'\).

\textbf{A free-boundary problem on this patch.} With \(q=u-\tau>0\), set \(w=(\phi-Q)/q\) on \(D\) and \(w=0\) elsewhere in \(B'\). Then \(w\in W^{2,p}_{\rm loc}(B')\) for every finite \(p>1\), with

\[
\Delta w=\mathbf1_D,\qquad w=\nabla w=0\quad\text{on }B'\setminus D .
\]

Indeed the zero extension is \(C^1\) because the function and its first derivatives vanish at the boundary. For each first derivative, restrict to almost every line parallel to a coordinate axis. Its untruncated version is absolutely continuous. On the intervals inside \(D\) it has zero values at the interval ends in \(B'\). Pasting it with zero preserves absolute continuity, with derivative zero off those intervals. This gives the weak Hessian \(\mathbf1_D D^2(\phi-Q)/q\) locally, including on boundaries of arbitrary shape.

Choose \(a\in H_j\) sufficiently close to \(p_0\). The ball of radius \(R=\operatorname{dist}(a,D)>0\) around \(a\) lies in \(B'\setminus D\) and touches \(\partial D\) at some \(z\in B'\). Both the ball and the contact point can be kept away from \(\partial B'\) since \(R\le|a-p_0|\). For every sufficiently small \(\rho<R\), the overlap \(B_\rho(z)\cap B_R(a)\) contains a ball of radius \(\rho/2\) centered at \(z+\rho(a-z)/(2R)\). The set where \(w\) and \(\nabla w\) both vanish therefore has minimum enclosing slab width at least \(\rho\) inside \(B_\rho(z)\).

To apply the published regularity theorem, fix \(s>0\) with \(\overline{B_s(z)}\subset B'\) and rescale to \(w(z+s\,\cdot)/s^2\) on the unit ball. This is a bounded local solution in the authors' Definition 1.1, with the origin on the boundary of the rescaled \(D\). Their Theorem III requires the minimum slab width of the simultaneous zero set, divided by the observation radius, to exceed a modulus tending to zero at some small scale. Here that ratio is at least one at all sufficiently small scales. The theorem and its corollary give an analytic boundary patch near the contact point. The class permits both signs of the solution.

The exterior ball lies on the side of this local graph outside \(D\), so \(D\) occupies the other side. Only the Cauchy data of \(Q\) on the patch are needed here. On the \(D\) side each \(\partial_i(\phi-Q)\) is harmonic with zero boundary data, so the one-sided third derivatives and the curvature calculation from the first proof apply. The patch is planar, and harmonic continuation gives \(D^2\phi=K+qnn^T\) throughout \(D\). In particular \(\det(D^2\phi-rI)=0\) on this component. If the stated local-finiteness hypothesis holds for every soft component, the translated-determinant contradiction follows.

In fact, in any putative equality solution with open phase representatives there is a soft component on which \(\det(D^2\phi-rI)\) is not identically zero. For that component, the closure of \textbf{each} hard grain meets its boundary in a relatively closed set with empty relative interior, by the patch argument above. There are only countably many grains. The Baire category theorem then gives a dense set of boundary points belonging to no individual hard-grain closure, even though hard material occurs in every neighborhood.
\phantomsection\label{endfile:rubber-endpoint}

\SSection{crystal}{Polycrystal shapes, finite certificates and the zinc slice}
Use stiffness eigenvalues $a=3\kappa$ and $b=2\mu$ unless physical
units are specified. The variational support families are defined in
the proof of Theorem~\ref{p13:interval}.
\subsection{The exceptional cases and the two Hill sides}\label{route:the-exceptional-cases-and-the-two-hill-sides}
For $e=I_3/\sqrt3$, define the hydrostatic projector $J=e\otimes e$ and
the deviatoric projector $D=I_V-J$. The four supports used below are
$H_{\epsilon\eta}(r,1)=\min(ra^\epsilon+b^\eta)$ across all
isotropic effective pairs. For a direction $n$ let $U_n v=n\odot v$
and $\Gamma_n=U_n(U_n^TCU_n)^{-1}U_n^T$. Construction (F) assigns
the same probability weights to directions in the two equations
\[
 \sum_n w_n\operatorname{Tr}
 \{P[(C-x_PI)-(C-x_PI)\Gamma_n(C-x_PI)]\}=0,\qquad P=J,D .
\]
The positive roots give the attainable pair $(x_J/3,x_D/2)$. The
fixed-shear case (A) is $C=2s(D-2J)+H\otimes H$. Write
$h=e:H$, $p=|DH|^2$ and assume $s>0$ and $h^2>4s+2p$.
Its whole isotropic region is
$[(2sh^2/(2s+p)-4s)/3,(h^2-4s)/3]\times\{s\}$.
\phantomsection\label{file:crystal-sides}
\paragraph*{Shear freedom at both sharp bulk bounds}

For any positive-definite three-dimensional crystal, the following alternatives exhaust the possibilities. All regions include limiting polycrystals. Write \(e=I_3/\sqrt3\), \(J=e\otimes e\), \(D=I_V-J\) and \(T=D-2J\).

\begin{itemize}
\item
  An isotropic crystal gives a single pair. An anisotropic crystal with \(Ce\) proportional to \(e\) gives a positive-length shear interval at its fixed bulk.
\item
  An anisotropic crystal \(C=bT+H\otimes H\), \(b>0\), has fixed shear \(b/2\) throughout its Hill bulk interval. Its entire tensor set has the explicit construction (49).
\item
  \textbf{Every other crystal admits a positive-length shear interval at each of its two Hill bulk endpoints.} Its joint region has a flat vertical side at both endpoints. There is a strictly positive lower bound on the shear width over the \textbf{whole closed bulk interval}.
\end{itemize}

In particular, for a crystal with varying bulk, a single possible shear at \textbf{either} sharp bulk forces the fixed-shear formula at \textbf{every} bulk. The proof uses rotations, binary laminates and a local calculation at an isotropic tensor. It gives a finite bound on the number of levels needed to witness variation at the sharp bulks. No uniform-grain-orientation condition is needed. The classical arithmetic and harmonic bounds determine the two bulk coordinates.

The \textbf{sharp values} at either side can be read directly from the four support energies in Theorem~\ref{p13:interval}. Take \((\epsilon,a_0)=(+1,3\kappa_R)\) or \((-1,3\kappa_V)\), and put \[
       h_\eta=\lim_{r\to\infty}
       \big[H_{\epsilon\eta}(r,1)-r a_0^\epsilon\big].
\] The shear interval on that side is \([h_+/2,1/(2h_-)]\). Indeed the bulk penalty is nonnegative, the displayed quantity increases with \(r\) and an endpoint optimizer bounds it above. Along minimizers the bulk penalty before multiplication by \(r\) tends to zero. Compactness then identifies the limit with the minimum of \((2\mu_*)^\eta\) on the side. This gives the endpoint values variationally even when the construction below realizes only part of the sharp interval.

\paragraph*{Two swaps prescribe the mean hydrostatic response}

Work with \(L=C\) for the Voigt bulk or \(L=C^{-1}\) for the Reuss bulk. In the first case the applied unit hydrostatic loading is a strain and the nonconstant layer fields lie in \(W_n=\{n\odot v\}\). In the second case the loading is a stress and \(W_n=\{Z:Zn=0\}\). Write

\[
 Le=\alpha e+g_0,\qquad \alpha=e:Le,\qquad g_0\in DV .
\]

In a laminate of two materials whose responses on \(e\) are \(\alpha e+g_1,\alpha e+g_2\), the uniform hydrostatic loading solves the layer problem if \(g_2-g_1\perp W_n\). The response of this laminate on \(e\) is then \(\alpha e+(1-f)g_1+fg_2\). Call a finite hierarchy hydrostatically uniform if it uses this condition at every binary layer, starting from rotated copies of \(L\). Its available response vectors \(g\) form exactly

\[
 K=\{g\in DV:\ \lambda_{\min}(g)\ge\lambda_{\min}(g_0),\quad
                    \lambda_{\max}(g)\le\lambda_{\max}(g_0)\}. \tag{45}
\]

Indeed the two eigenvalue inequalities hold under rotations and convex mixing. Conversely diagonalize the target and use the following two swaps in its principal frame. Let \(x_1\le x_2\le x_3\) be the original eigenvalues, and \(z_1\le z_2\le z_3\) the target. Mix \(x_1\) with \(x_2\) if \(z_1\le x_2\), and with \(x_3\) otherwise, until the first coordinate is \(z_1\). The other two coordinates bracket \(z_2\). In the first case they are \(x_1+x_2-z_1\le z_2\) and \(x_3\ge z_2\), by \(z_3\le x_3\) and equal traces. In the second case they are \(x_2\le z_2\) and \(x_1+x_3-z_1\ge z_2\). Mix the remaining pair to reach \(z_2,z_3\). Skip any swap already at its target. Each swap uses a rotated copy of the entire preceding material. For changing axes \(v,w\), the response jump is proportional to \(vv^T-ww^T\). A Voigt layer uses the unused axis as normal, and a Reuss layer uses \((v+w)/\sqrt2\). Thus two binary levels suffice. This also proves that (45) is the convex hull of the rotation orbit of \(g_0\).

A material at \(g=0\) has hydrostatic eigenvalue \(\alpha\) in the selected convention. All its rotated aggregates keep that eigenvalue. They therefore have the original crystal's Voigt bulk \(\alpha/3\) or Reuss bulk \(1/(3\alpha)\) on becoming isotropic. This holds for any correlations between orientations and layer normals.

\paragraph*{A five-level test for collapse of a bulk side}

Suppose \(g_0\ne0\). Unless \(C\) belongs to the fixed-shear family, \textbf{at least two distinct tensors at \(g=0\) already occur in hydrostatically uniform binary hierarchies of depth at most five}. Each tree uses at most 32 crystal leaves. If either tensor is anisotropic, the fixed-bulk argument below gives two different isotropic shears. If the two tensors are isotropic, they already give different shears.

The terminal set is nonempty by (45) at \(g=0\) and rotation invariant because rotating an entire admissible tree rotates its normals and all its jump equations. To prove the assertion, suppose there is only one such terminal tensor. Rotational invariance makes it \(B=\alpha J+\beta D\) in the selected stiffness or compliance convention, with \(\beta>0\). Then two materials of depth at most three with the same response vector \(g\) must coincide. To see this, complete the first one to \(g=0\) with at most two swaps. In the resulting tree replace just one occurrence of the original material by the second one, rotated with that leaf. Choose a leaf of positive weight and omit vacuous layers. The replacement has exactly the response of the selected leaf. At each ancestor the two input responses still match the original ones, so its hydrostatically uniform jump condition and its averaged response persist inductively to the root. Replacing only this occurrence suffices even when other leaves of the same material occur. At fixed normal and positive input fraction, mixing with a fixed second material is injective in the first. Indeed the layer coordinates

\[
 \Phi_n(L)=(P^{-1},P^{-1}Q,R-Q^TP^{-1}Q),\qquad
 L=\begin{pmatrix}P&Q\\Q^T&R\end{pmatrix}
 \quad\text{on }W_n\oplus W_n^\perp
\]

are invertible coordinates that average affinely. A change at the selected leaf therefore changes the terminal tensor. This contradicts uniqueness at depth five.

By (45), there is now a single-valued map \(A(g)\) on \(K\) for the two-level preparations. Its graph is compact because the parameters of the finitely many tree shapes are compact and the hydrostatic jump constraints are closed. It is continuous, rotation equivariant, and \(A(g)e=\alpha e+g\). For \(g_2-g_1\perp W_n\), laminating their two-level preparations gives a three-level material. The uniqueness just proved implies

\[
 \Phi_n(A((1-f)g_1+fg_2))
       =(1-f)\Phi_n(A(g_1))+f\Phi_n(A(g_2)). \tag{46}
\]

Zero lies in the interior of \(K\) in its five-dimensional space because \(g_0\) has both signs of eigenvalue. In the interior \(A\) is real analytic. Choose a basis of traceless matrices with eigenvalues \((d,0,-d)\), for instance two diagonal differences and three mixed shears. Each basis direction is perpendicular to some \(W_n\) in the selected convention. In any small coordinate box in these directions, (46) reconstructs \(A\) by repeated interpolation of its 32 corner values. Layer mixing and its inverse coordinate change are analytic at positive-definite matrices. This proves the claim throughout the connected interior.

\paragraph*{The local calculation forces a rank-one formula}

Write \(A'_v\) and \(A''_{vv}\) for the derivatives at zero. Since \(A'_v e=v\), rotational covariance and self-adjointness give, on a traceless \(E\),

\[
 A'_v E=(v:E)e+t\,\mathcal J_v E,\qquad
 \mathcal J_v E=D(vE+Ev)
\]

for a scalar \(t\). Also \(A''_{vv}e=0\). On deviatoric strains the latter operator is a linear combination, with constant coefficients, of \(|v|^2D\), \(\mathcal J_{D(v^2)}\) and \(v\otimes v\). These representation facts follow from \(\operatorname{Sym}^2(DV)=\mathcal H_0\oplus\mathcal H_2\oplus\mathcal H_4\), one copy each of the harmonic tensors of degrees zero, two and four. Equivalently the product of two trace-free quadratic forms has scalar, traceless quadratic and harmonic quartic parts. The three listed quadratic operators are independent, as can also be checked on the diagonal difference used next.

For \(v\perp W_n\), layer \(A(-\delta v)\) and \(A(\delta v)\) in equal fractions. Expanding the exact layer formula gives

\[
 \tfrac12 A''_{vv}=A'_v\Gamma_{B,n}A'_v,\qquad
 \Gamma_{B,n}=U(U^T B U)^{-1}U^T,\quad\operatorname{ran} U=W_n . \tag{47}
\]

\begin{itemize}
\item
  In stiffness take \(v=(xx-yy)/\sqrt2\) and \(n=z\), using dyads in an orthonormal frame. Both \(E_0=(2zz-xx-yy)/\sqrt6\) and \(E_2=(xy+yx)/\sqrt2\) give \(A'_v E\perp W_z\). They annihilate the quadratic energy in (47). On these vectors \(\mathcal J_{D(v^2)}\) has eigenvalues \(-1/3\) and \(1/3\), respectively, and \(v\otimes v\) vanishes. Thus the first two coefficients in \(A''_{vv}\) vanish. The remaining \(v\otimes v\) also vanishes on \((xz+zx)/\sqrt2\), whereas \(A'_v\) on that mixed shear has a nonzero projection into \(W_z\) whenever \(t\ne0\). This forces \(t=0\).
\item
  In compliance take \(v=(xz+zx)/\sqrt2\). Both \(z\) and \(x\) are allowed normals for this same jump. On \(E=(yz+zy)/\sqrt2\), \(A'_v E\) is proportional to \(t(xy+yx)\) with a nonzero proportionality factor. Its energy on the right of (47) is positive for \(n=z\) if \(t\ne0\) and zero for \(n=x\). Hence \(t=0\) in compliance as well.
\end{itemize}

Put \(\tau=2\) in stiffness and \(\tau=1/2\) in compliance, \(T_\tau=D-\tau J\), and \(\gamma=\alpha+\tau\beta\). The candidate with this derivative and value at zero is

\[
 A_0(g)=\beta T_\tau+
          (\gamma e+g)\otimes(\gamma e+g)/\gamma . \tag{48}
\]

It obeys (46) wherever the inputs are positive definite. In fact the compression of \(T_\tau\) to \(W_n\) is positive semidefinite of rank two. Its null line is \(nn^T\) in stiffness and \(I_3-nn^T\) in compliance. For a permissible jump \(g_2-g_1\), the projections of \(v_i=\gamma e+g_i\) onto \(W_n\) coincide. Their common projection couples to the null line by positive definiteness. Consequently \(v_i:\Gamma_{A_0(g_i),n}v_i=\gamma\), by applying the compressed matrix to a vector on that line. Formula (G) now subtracts exactly \(f(1-f)(g_2-g_1)\otimes(g_2-g_1)/\gamma\) from the arithmetic average, giving (48) at the mean.

Along every small ray \(g=xv\) with \(v\) of spectrum \((d,0,-d)\), the layer coordinates for \(A\) and \(A_0\) are affine for an allowed normal. Their first derivatives and values at zero agree, so the two materials coincide on the ray. Every sufficiently small traceless \(g\) is one more permissible mixture of two such singular matrices. In its principal frame, keep the eigenvalue \(z\) of sign opposite to the other two and mix \((z,-z,0)\) with \((z,0,-z)\), up to a common permutation. Their difference has spectrum \((d,0,-d)\). Thus \(A=A_0\) near zero, throughout the interior by analyticity, and on the boundary by continuity. In particular \((1-\epsilon)g_0\) is interior for \(0<\epsilon<1\), since \(0\) is interior. Letting \(\epsilon\) decrease to zero and using uniqueness of the preparation at \(g_0\) gives (48) for the original \(L=A(g_0)\). In stiffness this is the claimed rank-one family. In compliance, inversion of the rank-one perturbation of \(\beta T^{-1}\) gives a rank-one perturbation of \(T/\beta\) in stiffness. Its hydrostatic diagonal after subtracting \(T/\beta\) is positive, so the rank-one term is positive. This proves the five-level assertion. Conversely the fixed-shear relation throughout this family follows from the two energy translations in (B).

\paragraph*{Different isotropic shears from an anisotropic fixed-bulk material}

Here is a direct proof and a finite prescription for a pair of distinct shears. Let \(A e=a e\) in stiffness. Choose any terminal normal \(n_0\) in (F), and call the attained tensor \(B=aJ+bD\). Suppose \(A\) is anisotropic. Form the dilute layer of fraction \(f\) of \(A\) in \(B\), normal \(n_0\). Its expansion is \(B+fH+O(f^2)\), with

\[
 H=(A-B)-(A-B)\Gamma_{A,n_0}(A-B).
\]

\(He=0\) and \(\operatorname{Tr}H=0\) by the shear root equation. Moreover \(H\ne0\). Indeed \(I-\Gamma_{A,n_0}(A-B)\) induces the identity modulo \(W_{n_0}\), and on \(W_{n_0}\) is the product of the inverse compressed \(A\) and compressed \(B\), both invertible.

Take two rotations of this dilute layer and laminate them in equal fractions at a fixed normal. The result \(A_f=B+fZ+O(f^2)\) has \(Z=(R\cdot H+S\cdot H)/2\) and retains the hydrostatic eigenvalue. Write \(q_n(Z)=\operatorname{Tr}(\Gamma_{B,n} Z^2)\). Some two rotations and normals \(n,m\) have \(q_n(Z)\ne q_m(Z)\). If all these differences vanished, taking \(R=S\) and then polarizing at pairs would give the identity on the whole linear span of the rotation orbit of \(H\). This span contains a full degree-two or degree-four summand of the traceless symmetric operators on \(DV\). In degree four a binary harmonic quartic, for example \(\operatorname{Re}(x+iy)^4\), has \(q_z=0\) and \(q_x>0\). In degree two the operators are \(\mathcal J_v\). If their energies were independent of \(n\) for every \(v\), the quadratic form \(v\mapsto q_z(\mathcal J_v)\) would be rotation invariant and scalar on \(DV\). Unit \(v=(2zz-xx-yy)/\sqrt6\) and \(v=(xx-yy)/\sqrt2\) instead give values differing by \(2/(3b)\), by direct projection onto the mixed and axial waves.

The two shear eigenvalues given by (F) on \(A_f\) differ by \(-f^2[q_n(Z)-q_m(Z)]/5+O(f^3)\). In fact the shear quadratic at \(b\) is \(\operatorname{Tr}(A_f-B)-\operatorname{Tr}(\Gamma_{A_f,n}(A_f-B)^2)\), and its derivative at \(f=0,x=b\) is \(-5\) in \(x\). The two roots are therefore distinct for all sufficiently small \(f>0\). All ingredients stay in the laminate closure of the original material. Interpolating at their common bulk gives the intervening shears.

One can select the normals from the 15 normalized vectors \((i,j,1)\) with nonnegative integers \(i+j\le4\). The extension \(x\mapsto|x|^4q_{x/|x|}(Z)\) is a homogeneous quartic, and this triangular grid determines its restriction to the chart \(n_3=1\). For the rotations one may use the 165 scalar-first quaternions \((1,i,j,k)\), \(i+j+k\le8\), with nonnegative integer entries. A linear functional of the rotated fourth-order tensor, times the fourth power of the squared quaternion length, is a polynomial of total degree at most eight in this chart. The simplex grid determines all its coefficients by finite differences. Hence these rotations span the full orbit span and the polarization argument still applies. Halving \(f\) until two roots differ terminates. Alternatively, at either bulk boundary the existence test for two terminal tensors uses only the two depth-five trees, rotations, layer fractions and the exact jump equations. Clearing the positive acoustic determinants gives a finite real-polynomial feasibility test.

\textbf{Finite laminates on the sharp bulk sides.} In the area case, each sharp bulk has a positive-length open shear interval every point of which is attained by a finite hierarchy of pure grains. Let \(L_f\) be the finite hydrostatically uniform hierarchies ending at \(g=0\) for this side, regarded in stiffness coordinates, and let \(L\) be their closure. They retain this bulk under every further laminate, and \(L\) contains the two isotropic shears just constructed. For a fixed normal \(\Phi_n(L_f)\) is convex and dense in the compact convex \(\Phi_n(L)\). At any isotropic tensor strictly between the two shears, the tangent cone pulled back from these layer coordinates is rotation invariant and has both \(D\) and \(-D\) as tangent directions. If the coordinate point were on its relative boundary, a supporting normal in the span of the pulled-back cone would average to a nonzero isotropic normal. Indeed its rotations pair strictly negatively and uniformly with the center of a relative ball in the cone. That span consists of operators annihilating \(e\), since all the tensors have this fixed hydrostatic eigenvalue. Its only isotropic direction is therefore \(D\), a contradiction. The coordinate point is in the relative interior and belongs to \(\Phi_n(L_f)\) by density and convexity. The five-level bound concerns precursor tensors. This argument proves finiteness for the isotropic pairs in the stated open interval.

\paragraph*{An exactly isotropic interpolation laminate}
The following finite laminate works between any two attainable isotropic inputs. Write the matrix or host parameters as \(p,q\) (bulk and shear stiffness eigenvalues) and the core parameters as \(a_1,b_1\). At core fraction \(f\), the effective tensor is exactly isotropic with eigenvalues

\[
 \begin{split}
 a(f)&=p+\frac{f(a_1-p)}{1+(1-f)h(a_1-p)},\qquad h=\frac1{p+2q},\\
 b(f)&=q+\frac{f(b_1-q)}{1+(1-f)d(b_1-q)},\qquad
 d=\frac25\left(\frac1q+\frac1{p+2q}\right).
 \end{split} \tag{25}
\]

Each coordinate varies monotonically between its own endpoints. To construct the laminate, add successive layers of the same host around the core. For host stiffness \(B\) and original core stiffness \(A\), the resulting stiffness is

\[
 B+f\big[(A-B)^{-1}+(1-f)\sum_j c_j\Gamma_{B,n_j}\big]^{-1},
 \qquad
 \Gamma_{B,n}=U_n(U_n^T B U_n)^{-1}U_n^T,\quad U_n v=n\odot v .
\]

To see the weights, start with core fraction $f_0=1$ and keep the
fraction $f_j/f_{j-1}$ of the preceding laminate at step $j$, where
$f_j=f_{j-1}-(1-f)c_j$. The layer Schur complement makes the inverse
increment from the host additive in the terms $c_j\Gamma_{B,n_j}$,
giving the displayed formula by induction. Singular contrasts follow
by continuity. Give each of the three coordinate axes weight $2/15$
and each of the four unoriented body diagonals weight $3/20$.
The second and fourth moments are the spherical moments: for example,
each diagonal fourth moment is $2/15+(3/5)/9=1/5$, and each distinct
paired cross moment is $(3/5)/9=1/15$. Thus the averaged acoustic
operator has the hydrostatic and deviatoric coefficients $h,d$ in
(25). Substitution proves the interpolation formula for the finite
laminate \cite[Eqs.~(2-15) to (2-18)]{p8void}.

\paragraph*{A uniform width across the bulk interval}

Write \(a_\pm=3\kappa_{R,V}\) and \(b_R=2\mu_R\), \(b_V=2\mu_V\) for the Hill bounds. Suppose attained endpoint shear intervals have physical widths \(w_-,w_+>0\). For every \(a\in[a_-,a_+]\) the shear width at bulk \(a/3\) is at least

\[
 \frac{a_-+2b_R}{a_++2b_R}
 \left[1+\frac35\left(\frac{b_V}{b_R}-1\right)\right]^{-2}
 \frac{w_-w_+}{w_-+w_+}.
\]

To prove this, choose one endpoint material as host with eigenvalues \(p,q\), and vary the core shear across the interval at the other endpoint, of bulk eigenvalue \(a_1\). In (25) the core fraction at the target is \(f=u(a_1+2q)/(a+2q)\) with \(u=(a-p)/(a_1-p)\). The derivative of the resulting shear with respect to core shear is \(f/[1+(1-f)d(b_1-q)]^2\), and \(dq\le3/5\). Swapping endpoint roles gives the same bound with \(1-u\) as normalized bulk fraction. Use \(\max(u w_+,(1-u)w_-)\ge w_-w_+/(w_-+w_+)\). One can choose finite hierarchies at each end of both endpoint intervals, inside the open portions proved above. Interpolation (25) within each interval and then between the two bulks takes only a fixed finite set of layer directions. It therefore gives positive-width shear intervals at every bulk with a depth bound depending only on the crystal. Alternatively, fix the Reuss end as host and vary both endpoint shears continuously. For each bulk the resulting shear image is a compact interval, its endpoints vary continuously with bulk, and its width is positive including at both ends. This constructs a whole band of finite laminates even when the two sharp shear sections are disjoint.
\phantomsection\label{endfile:crystal-sides}

\phantomsection\label{file:crystal-fixed-shear}
\paragraph*{Full three-dimensional tensor region of the fixed-shear family}

Write \(T=I_V-I_3\otimes I_3=D-2J\). Consider any positive-definite crystal

\[
 C=bT+y_0 P_0\otimes P_0,\qquad
 \operatorname{tr}P_0=1,\quad b>0,\quad y_0(1/2-|P_0|^2)>b .
\]

Here \(P_0\) is a symmetric \(3\times3\) matrix and \(y_0>0\). In the notation \(C=bT+H\otimes H\), choose \(\operatorname{tr}H>0\), take \(P_0=H/\operatorname{tr}H\) and \(y_0=(\operatorname{tr}H)^2\). Positivity is equivalent to the displayed stability condition, by congruence or a rank-one inverse. In particular \(|P_0|^2<1/2\) and \(P_0\) is positive definite. If its least eigenvalue were nonpositive, the squared sum of the other two would be at least \(1/2\). Set $2\mu_0=b$ and $D_0=\sqrt{y_0}P_0$ in the Grabovsky and Milton relation. Its positivity condition $(\operatorname{tr}D_0)^2-2\operatorname{tr}(D_0^2)>4\mu_0$ is exactly the displayed stability inequality. Each rotated grain belongs to the same family, which is closed under homogenization \cite[Eqs.~(25), (33) to (36)]{grabmilton98}. The following bounds determine precisely which points of that family are realized.

The \textbf{entire effective tensor set}, including anisotropic responses and arbitrary three-dimensional arrangements, consists of

\[
 \boxed{\begin{gathered}
 A=bT+yP\otimes P,\qquad \operatorname{tr}P=1,\qquad
 \lambda_{\min}(P_0)I_3\preceq P\preceq\lambda_{\max}(P_0)I_3,\\
 \frac1{y_0}\ \le\ \frac1y\ \le \frac1{y_0}+\frac{|P_0|^2-|P|^2}{b}.
 \end{gathered}} \tag{49}
\]

Any target eigenframe is allowed. Both amplitude endpoints take at most two binary lamination levels, and \textbf{three levels suffice for every tensor}. For the isotropic response set \(P=I_3/3\), so \(\mu_*=b/2\), \(\kappa_*=y/9-2b/3\). The resulting interval is exactly (A). The bulk-side classification shows that these are precisely the anisotropic crystals with varying bulk for which even one sharp bulk has a unique shear.

\paragraph*{Why the shape is in the convex hull of the crystal shapes}

Let \(P_R(x)=R(x)P_0R(x)^T\) in a periodic cell of unit volume. There is a nonnegative density \(\rho\in L^2\) with mean one satisfying

\[
 \langle\rho P_R:D^2\phi\rangle=0
 \qquad\text{for every periodic }\phi\in H^2 .
\]

Here the scalar Hessian is denoted by \(D^2\phi\). The Cordes condition ensuring this density follows from stability. With \(c=1/|P_0|^2\) one has \(|cP_R-I_3|^2=3-1/|P_0|^2<1\). On periodic scalar \(L^2\) define \(Mf=f+(cP_R-I_3):D^2\Delta^{-1}(f-\langle f\rangle)\). The map \(D^2\Delta^{-1}\) is an isometry on mean-zero functions by Fourier expansion, so \(\|M-I\|<1\). Hence \(M\) is invertible, \(M1=1\) and \(\rho=(M^{-1})^*1\) has the claimed mean and annihilation property. The range of \(P_R:D^2\) on periodic potentials is exactly \(\rho^\perp\), the image under \(M\) of the mean-zero space up to the constant factor \(c\).

For positivity, first smooth \(P_R\) by periodic convolution. The strict contraction and uniform ellipticity persist. For smooth \(r\ge0\) solve \(Mf=r\) and write \(f=\langle f\rangle+\Delta\phi\). Elliptic regularity makes \(\phi\) smooth and \(r=\langle f\rangle+cP_R:D^2\phi\). At a maximum of \(\phi\), \(\langle f\rangle\ge0\). This proves positivity of the smoothed density. The uniform inverse bound supplies a weak \(L^2\) limit, which annihilates the original coefficient against every smooth Hessian and therefore against every Hessian of an \(H^2\) potential. Uniqueness identifies it with \(\rho\). This is the periodic invariant density of the scalar elliptic equation in nondivergence form.

In fact the \textbf{whole elasticity cell tensor} follows from this one scalar density. At frequency \(\xi\ne0\) put \(n=\xi/|\xi|\). Then

\[
 P=\langle\rho P_R\rangle,\qquad
 \frac1y=\frac{\langle\rho^2\rangle}{y_0}
       +\frac2b\sum_{\xi\ne0}|\widehat{\rho P_R}(\xi)n|^2. \tag{50}
\]

To prove it, prescribe a mean strain \(E\). Decompose the periodic displacement fluctuation as \(w+\nabla\phi\) with \(\operatorname{div}w=0\). The translated energy above \(bE:T E\) is

\[
 \frac b2\|\nabla w\|_2^2
       +y_0\|P_R:(E+\operatorname{sym}\nabla w+D^2\phi)\|_2^2 .
\]

The first term follows by Fourier expansion of \((\operatorname{sym}\nabla u):T(\operatorname{sym}\nabla u)\). Minimizing over \(\phi\) projects the last scalar field onto \(\rho\). Its squared norm becomes \((P:E+\ell(w))^2/\langle\rho^2\rangle\), where \(\ell(w)=\langle\rho P_R:\nabla w\rangle\). With the norm \(\|\nabla w\|_2\) on mean-zero divergence-free displacements the squared dual norm of \(\ell\) is the sum in (50). Indeed \(n^T\widehat{\rho P_R}(\xi)n=0\), so the vector in that sum is already perpendicular to \(n\). Minimizing this last rank-one quadratic in \(w\) gives \(y(P:E)^2\), proving (50). Since \(\rho\ge0\) and \(\langle\rho\rangle=1\), \(P\) belongs to the convex hull of the rotation orbit of \(P_0\). The two eigenvalue tests in (49) describe this hull by the swaps of (45).

The upper hydrostatic strain and stress energies give the remaining necessary bounds. With \(e=I_3/\sqrt3\) their values for \(A\) are

\[
 e:Ae=-2b+y/3,\qquad
 e:A^{-1}e=-1/(2b)+\frac{1}{12b(1/2-|P|^2-b/y)} .
\]

The second formula follows by inverting the rank-one perturbation of \(bT\). Each energy is at most its value for the crystal, since this value is constant on the crystal orbit. Rearrangement gives exactly the amplitude bounds in (49). The entire argument passes to arbitrary homogenization limits by periodic approximation. The bounds preserve a strictly positive stability margin for every target.

\paragraph*{Realizing every tensor}

For a desired shape \(P\) diagonalize it, align the ordered crystal eigenvalues in this frame, and perform the two weighted swaps of (45). For the stiffness construction choose the unused principal axis as layer normal. The uniform hydrostatic strain gives the arithmetic average of the responses \(-2be+y_0P_R/\sqrt3\). Thus the resulting material has \(y=y_0\) and precisely the desired \(P\) by (50).

For the compliance construction use the bisectors of the swapped axes as layer normals. Inverting the local stiffness gives \(C^{-1}=T^{-1}/b+y'_0 P'_0\otimes P'_0\), where \(P'_0=I_3-2P_0\) has trace one and

\[
 y'_0=\frac1{4b(1/2-|P_0|^2-b/y_0)}.
\]

The same swaps and fractions take \(P'_0\) to \(I_3-2P\) in the mean response to uniform hydrostatic stress. Hence this second construction attains the other amplitude endpoint of (49) at exactly the same shape.

Finally laminate two resulting tensors of common shape, of amplitudes \(y_1,y_2\), using fraction \(f\) of the second. Their difference is \((y_2-y_1)P\otimes P\). On the strain wave space the compression of \(bT\) is positive with one zero direction \(nn^T\), and \(P:nn^T>0\). Therefore \(P:\Gamma_{bT+v P\otimes P,n}P=1/v\) for \(v>0\) in this laminate. Substitution in (G) gives \(1/y=(1-f)/y_1+f/y_2\) for \textbf{every normal}. This realizes the full interval.

For example take \(b=2\) and crystal director \(H=\operatorname{diag}(1,2,2)\). The target shape with eigenvalues \((1/4,7/20,2/5)\) has exactly \(400/19\le y\le25\). Both endpoint preparations swap the first two axes with fraction \(1/4\). Combining them with fraction \(2/3\) of the compliance preparation gives \(y=200/9\). Setting \(P=I_3/3\) instead yields the full isotropic region \([3/4,13/9]\times\{1\}\).
\phantomsection\label{endfile:crystal-fixed-shear}

\subsection{Finite-frequency and polynomial certificates}\label{route:finite-frequency-and-polynomial-certificates}
\phantomsection\label{file:polycrystal-certificates}
\subsubsection*{Proof of Lemma~\ref{p13:band}}
\begin{proof}
It suffices to prove
it for trigonometric \(F\) and pass to the \(L^2\) limit.

Localization of a single band. Let \(B\in\mathcal F\) have frequencies \(l\le|k|\le u\). Put
\(D=2^{4m}\), \(d_0=D/l\), and assume \(N_m\ge2d_0u\). Translate and scale a unit cube by writing
\(x=y+d_0z\), with \(y\) uniform on the original torus and \(z\in[0,1]^3\). Set \(\alpha=2^{-2m}\)
and \(\chi(z)=\min(1,\operatorname{dist}(z,\partial[0,1]^3)/\alpha)\). On the \(z\)-torus project
\(\chi(z)B(y+d_0z)\) onto the admissible modes \(0<|j|\le N_m\). Call the real result \(P_y(z)\).
Then
\[
\|P_y\|_{L^2(y,z)}\le\|B\|_2,\qquad
 \|B(y+d_0z)-P_y(z)\|_{L^2(y,z)}
       \le \rho_m\|B\|_2,\quad \rho_m=3\,2^{-m}+2\,2^{-2m}.
 \tag{13.6}\label{p13:eq6}
\]

To see this, averaging over \(y\) makes distinct original \(k\)'s orthogonal. For each \(k\), the
\(z\)-Fourier coefficients of \(\chi(z)e^{2\pi i d_0 k\cdot z}\) are denoted \(a_j\). The function
has zero boundary trace and is periodic as an \(H^1\) function. Parseval applied to its derivative
minus \(2\pi i d_0 k\) times itself gives
\[
\sum_j |j-d_0 k|^2 |a_j|^2=\|\nabla\chi\|_2^2/(2\pi)^2 .
\]

For retained \(j\), the projection error on the range of \(\Pi(k)\) has norm at most \(8|j-d_0
k|/(d_0|k|)\). Indeed the difference of two stress projections is at most twice the norm of the
difference of the normal dyads, at most \(4|j/|j|-k/|k||\). The same estimate holds for their
complements. For \(j=0\) or \(|j|>N_m\) the displayed rescaled distance is at least one before
multiplication by eight. Thus projection loses at most \(8\|\nabla\chi\|_2/(2\pi D)\) times
\(\|B\|_2\) in the joint norm. Removing the window costs \(\|1-\chi\|_2\|B\|_2\). The bounds
\(\|1-\chi\|_2\le\sqrt{6\alpha}\) and \(\|\nabla\chi\|_2\le 1/\alpha\) give \eqref{p13:eq6}.

Adding a band above a varying mean. For a low-frequency field \(A\), use the remainder
\[
T(A,B)=g(A+B)-g(A)-Dg(A)\cdot B .
\]

Changing \(A\) to \(A'\) changes \(T\) by at most \(2M|A-A'||B|\). Changing \(B\) to \(B'\)
changes \(T\) by at most \(\tfrac M2(|B|+|B'|)|B-B'|\). Both estimates follow by integrating the
gradient on a segment. For each fixed \(z\), translation of \(y\) preserves torus measure, so
averaging an integrable function of \(y+d_0z\) recovers its original integral even when the
translated cube wraps around the torus. The translation difference of \(A\) in \(L^2\) is at most
\(d_0|z|\|\nabla A\|_2\) by integration along the translation. Replacing \(A(y+d_0z)\) by \(A(y)\)
costs at most \(2\sqrt3 M d_0\|\nabla A\|_2\|B\|_2\) in the mean remainder. Replacing
\(B(y+d_0z)\) by \(P_y(z)\) now costs at most \(M\rho_m\|B\|_2^2\) by \eqref{p13:eq6}. For each \(y\) the
mean remainder with \(A(y),P_y\) is nonnegative by the finite test and \(\int P_y=0\). The
discarded linear term on the original torus satisfies
\[
\left|\int Dg(A)\cdot B\right|
 \le \frac{\|\nabla(Dg(A))\|_2}{2\pi l}\|B\|_2
 \le \frac{M\|\nabla A\|_2}{2\pi l}\|B\|_2 .
\]

The first inequality uses the Fourier support of \(B\), and the second uses the Lipschitz chain
rule, also obtainable from difference quotients. Combining these estimates gives
\[
\int [g(A+B)-g(A)]\ge
 -M\rho_m\|B\|_2^2-M\left(2\sqrt3\,d_0+\frac1{2\pi l}\right)
                         \|\nabla A\|_2\|B\|_2 .              \tag{13.7}\label{p13:eq7}
\]

All frequency scales. Divide the nonzero Fourier lattice into shells \(2^h\le |k|<2^{h+1}\). In
each consecutive period of \(w+p\) shell indices delete \(p=8m\) consecutive residue classes, with
\(w=8m^3\). Averaging over the cyclic offset gives a choice for which the removed field \(T_0\)
has squared norm at most \(\theta\|F\|_2^2\), \(\theta=p/(w+p)=1/(m^2+1)\). Subtracting the
constant linear part \(Dg(X)\cdot T_0\), applying Cauchy-Schwarz with \(F-T_0\), and using the
quadratic Taylor error bounds the cost of adding it back by \(M(\sqrt\theta+\theta/2)\|F\|_2^2\).

In increasing order the remaining bands \(B_j\) have \(u_j/l_j\le 2^w\), and \(u_i/l_j\le
r^{j-i}\) when \(i<j\), with \(r=2^{-8m}\). A shortened first band obeys these same inequalities.
Apply \eqref{p13:eq7} to \(B_j\) with \(A=X+\sum_{i<j}B_i\). We have \(N_m\ge2D u_j/l_j\) and \(\|\nabla
A\|_2\le2\pi\sum_{i<j}u_i\|B_i\|_2\). For \(b_i=\|B_i\|_2\),
\[
\sum_{i<j} r^{j-i} b_i b_j\le\frac{r}{1-r}\sum_i b_i^2 .
\]

This follows by \(2b_i b_j\le b_i^2+b_j^2\) at each index separation. The cross-band error in the
telescoping sum is at most \(M(1+4\pi\sqrt3D)r/(1-r)\sum_i b_i^2\), bounded by
\(26M\,2^{-4m}\sum_i b_i^2\). Orthogonality bounds \(\sum b_i^2\le\|F\|_2^2\). Adding the removal
error and \eqref{p13:eq6} proves \eqref{p13:eq5}.

\end{proof}

\subsubsection*{Proof of Lemma~\ref{p13:certificate}}
\begin{proof}
All three conditions are decidable by real quantifier elimination \cite[Algorithms~14.3, 14.5, pp.~546, 549 to 550]{bpr}. Each algebraic constant can be represented by an additional variable satisfying its defining polynomial and its isolating rational interval. In \eqref{p13:eq8} the variables are
just \(A\) and the finitely many real sine and cosine coefficients. Their admissibility is linear
at each frequency, \(\|P\|_2^2\) is their Parseval sum, and the sampled sines and cosines are real
algebraic numbers specified by roots of unity. The graphs of \(h,g\) are semialgebraic including
at zero. Their second derivatives on the unit sphere are algebraic expressions in coordinates. For
condition 1, rotations can be written as \(R^TR=I,\ \det R=1\). The same checks work for entries
given as explicit real algebraic numbers with isolating data.

For condition 2, homogeneity gives the Hessian bound away from the origin, and the gradient
extends by zero at the origin. Integration along segments gives \(\operatorname{Lip}(Dg)\le M\)
also across zero. For \eqref{p13:eq8}, the average of \(P\) vanishes both continuously and on the grid.
Thus each difference can have \(Dg(A)\cdot P(z)\) subtracted before comparing averages. The
spatial Lipschitz constant of the resulting remainder is at most \(M\|P\|_\infty\|\nabla
P\|_\infty\). Cauchy-Schwarz on the Fourier coefficients bounds these two norms by
\((2N+1)^{3/2}\|P\|_2\) and \(2\pi N(2N+1)^{3/2}\|P\|_2\). Multiplying their product by the
grid-cube diameter \(\sqrt3/K\) bounds the quadrature error by \(q\|P\|_2^2\). Consequently
\eqref{p13:eq8} implies ordinary integral Jensen for \(g\) at all of the frequencies in \eqref{p13:eq4}.

Apply Lemma~\ref{p13:band} to \(g\). Adding \(\tau|\cdot|^2\) increases the integrated Jensen difference by
exactly \(\tau\|F\|_2^2\). Condition 2 absorbs the entire frequency error. Since \(h\le f\), taking the infimum in \eqref{p13:eq1} gives \eqref{p13:eq9}.
\end{proof}

\subsubsection*{Proof of Lemma~\ref{p13:exhaustion}}
\begin{proof}
For \(0<\delta<1\), average \(W(UX)\) against a nonnegative smooth probability density on
\(SO(d)\) supported on \(\|U-I\|\le\delta\). These \(U\) act on the full space of fields. Denote
the average by \(W_\delta\). Group convolution makes its restriction to the unit sphere smooth
because rotations act transitively and differentiation can be moved to the smooth density. The
function remains two-homogeneous and its gradient has Lipschitz constant at most \(M_0\). By
\eqref{p13:eq3},
\[
|W_\delta(X)-W(X)|\le M_0\delta |X|^2,\qquad
 \int W_\delta(X+F)-W_\delta(X)\ge-M_0\delta\|F\|_2^2 .
 \tag{13.10}\label{p13:eq10}
\]

For the second assertion apply \eqref{p13:eq2} at mean \(UX\) with the original admissible \(F\). Compare
\(W(UX+F)\) and \(W(UX+UF)\) along the segment between the two arguments. Subtract
\(DW(UX)\cdot(U-I)F\), whose integral vanishes. Along the segment the gradient differs from
\(DW(UX)\) by at most \(M_0|F|\) because \(U\) is orthogonal. The remaining error is bounded by
\(M_0\delta\|F\|_2^2\). This proves \eqref{p13:eq10} even when \(UF\) fails the differential constraints.

Write \(b=M_0\delta\), \(D_\delta=1+5b/c\), and set
\[
v_\delta(X)=\frac{W_\delta(X)+3b|X|^2}{D_\delta}.
\]

It is smooth on the sphere, converges uniformly there to \(W\), and satisfies
\[
v_\delta(X)\le W(X)-\frac{b}{D_\delta}|X|^2,\qquad
 \int v_\delta(X+F)-v_\delta(X)
       \ge \frac{2b}{D_\delta}\|F\|_2^2.                      \tag{13.11}\label{p13:eq11}
\]

The first bound uses \(W\ge c|\cdot|^2\), the second uses \eqref{p13:eq10}. Approximate
\(v_\delta|_{\mathbb S^{d-1}}\) in \(C^2\) by rational polynomials restricted to the sphere. For
instance, extend the smooth angular function with a radial cutoff on an annulus, use
tensor-product Bernstein polynomials on a containing cube, and approximate their finitely many
coefficients by rationals. The derivatives of orders up to two of the Bernstein approximations
converge uniformly there, by their finite-difference formulas. The two-homogeneous extensions of
the angular approximations converge in value on the sphere and in the uniform Hessian norm away
from zero. Indeed the Cartesian Hessian on the sphere is a continuous linear expression in the
angular derivatives through second order, and it is zero-homogeneous. The gradient bound extends
over zero as in condition 2.

In particular we can choose \(h\) of the enumerated form with \(h\le f\), as close as desired to
\(v_\delta\) on the sphere, and with a strictly positive integral Jensen margin
\(\gamma\|F\|_2^2\) for every fluctuation. Use the strict first margin in \eqref{p13:eq11} for the
obstacle and the second for the Jensen margin. A perturbation with gradient Lipschitz constant
\(e\) changes the Jensen difference by at most \(e\|F\|_2^2/2\) after subtracting its linear part.

Choose rational \(0<\tau<\gamma\). For \(g=h-\tau|\cdot|^2\) choose rational \(M\) above its
Hessian norm and then \(m\) with \(ME_m\le\tau\). Finally choose \(K>2N_m\) so large that
\(2q<\gamma-\tau\). The integral Jensen difference for \(g\) at each finite test is at least
\((\gamma-\tau)\|P\|_2^2\), so the quadrature estimate gives \eqref{p13:eq8}. This certificate is therefore enumerated. As $\delta\downarrow0$
these lower bounds converge to $W$ on the unit sphere.
\end{proof}
\phantomsection\label{endfile:polycrystal-certificates}

\subsection{What quadratic descriptions can miss}\label{route:what-quadratic-descriptions-can-miss}
Consider the rational Kelvin compliance $S=I_V+F/4$, where
\[
 F=\frac1{14}\left[
 \begin{pmatrix}17&2&-5\\2&-12&-4\\-5&-4&9\end{pmatrix}
 \oplus\operatorname{diag}(-8,-10,4)\right].
\]
The pair $(\kappa,\mu)=(1/3,1/2)$ passes the complete joint quadratic
test on the six strain and stress loads. The counterexample compares
the rotation-averaged covariance of a material response with that of a
compatible, equilibrated wave. The covariances coincide on every
invariant quadratic test. A cubic compatibility condition forces
$\mu>1/2+4\cdot10^{-11}$ at every physical bulk. The following proof computes the covariances and the strict gap. A distinct
obstruction concerns polynomial descriptions of the sharp graph. For
each integer $N\ge2$ some rational crystal requires bulk degree at
least $N$ in any finite quantifier-free polynomial sign description
near a specified rational interior bulk, even allowing arbitrary real
coefficients. There is an example for each shear endpoint. This
assertion ranges over crystals. It is consistent with the terminating
algorithm for fixed algebraic data in the article.
\subsection{Covariances that pass all quadratic tests and violate a cubic constraint}\label{detail:crystal-moment-proof}\label{file:crystal-moment-proof}
For the rational compliance \[
 S=I_V+\tfrac14 F,\qquad
 F=\frac1{14}\left[
 \begin{pmatrix}17&2&-5\\2&-12&-4\\-5&-4&9\end{pmatrix}
 \oplus\operatorname{diag}(-8,-10,4)\right],
\] the pair \((\kappa,\mu)=(1/3,1/2)\) passes the \textbf{full joint second-moment test on all six load maps}. The relaxation tests local quadratic forms on the six-load strain-stress map that are nonpositive on every rotated material graph and nonnegative on every compatible-strain/equilibrated-stress wave. The displayed pair is the joint Reuss corner. The Reuss bulk \(1/3\) is physically attainable. The shear exceeds \(1/2+4\cdot10^{-11}\) \textbf{at every bulk}. At the exact Reuss bulk it exceeds \(1/2+1.5\cdot10^{-10}\).  Thus even the optimum over all joint quadratic translations gives a wrong shear boundary on an explicit rational crystal. 

\paragraph*{Meaning of the quadratic test.}
For any quadratic form on a six-load strain and stress array, a periodic
gradient and a divergence-free stress have nonzero Fourier modes in
the compatible and equilibrated subspaces for the same normal.
Consequently subtracting the value of a quadratic form at the mean
leaves its pairing with a positive sum of wave covariances. Material
covariances must likewise be positive combinations of points on the
rotated compliance graph. Equality of the two candidate covariances
after subtraction at the mean certifies every such quadratic
inequality, by linearity. Rotation equivariance permits simultaneous
Haar averaging of the arrays and of the test.
\paragraph*{A rank-one certificate at the exact bulk endpoint}

Write \(\Pi\) for rotation averaging, \(F I_3\) for the spatial trace over the last two tensor indices, and \(F^{\otimes2}\) for the outer product of its coefficients as a load map. The given elasticity perturbation is fully symmetric in its four spatial indices. One checks \(\Pi F=0\), \(F I_3=\operatorname{diag}(1,-1,0)\) and \(\|F\|^2=4\). Thus \(S\) has spectrum in \([1/2,3/2]\) and its arithmetic isotropic compliance is \(I_V\). Its hydrostatic coupling is nonzero, so the Hill bulks differ. More precisely they are \(1/3\) and \(23080/66321\).

Let \(n=e_3\), \(v=e_1\), and \(G=\partial_v n^{\otimes4}\), the sum placing \(v\) in each of its four index positions. The only nonzero components of \(G\) are the four permutations of \(3331\), all equal to one. Every response column of \(G\) is a compatible strain with normal \(e_3\). Also \(\Pi G=0\), \(\|G\|^2=4\), and \(\|G I_3\|^2=2=\|F I_3\|^2\).

A fully symmetric fourth-order tensor of zero isotropic part is the orthogonal sum of one degree-two and one degree-four spherical harmonic tensor. Its degree-two part is the traceless lift \(\operatorname{sym}_6(\delta\otimes F I_3)/7\) with squared norm \(6\|F I_3\|^2/7\). Haar averaging the outer product gives a scalar on each of these inequivalent irreducible summands, equal to the squared norm in that summand divided by its dimension. Thus \[
             \operatorname{Av}_R (R\cdot F)^{\otimes2}
             =\operatorname{Av}_R (R\cdot G)^{\otimes2}.
\] Use unit mean stress loads. In the material frame use the strain/stress map \((S,I_V)\) and in the wave frame use \((G/4,0)\). The isotropic macro map is \((I_V,I_V)\). For every invariant joint quadratic form \(T\), expand the material value about this macro map. The linear term vanishes under rotation averaging and the quadratic term equals the wave value by the covariance identity. This is precisely the joint quadratic covariance identity with rank-one positive material and wave matrices. In particular all four hydrostatic and deviatoric weighted energy tests, at every weighting, accept this corner. The covariance therefore passes every one of the quadratic tests in question.

There is also a direct rational check of all \(72^2\) entries. In raw Cartesian coordinates let \(D_i\) generate rotation about coordinate axis \(i\) on the paired load map, and let \(L_i M=D_iM+MD_i^T\) act on its covariance. On an angular harmonic summand of degree \(l\), \(\sum_i L_i^2=-l(l+1)I\). The tensor order of the covariance is at most eight. Consequently the product of \(l(l+1)I+\sum_i L_i^2\) over \(1\le l\le8\) is a positive scalar times the Haar projection. Applied to the material covariance minus the wave and macro covariances it gives the zero matrix exactly. The program tests compatibility and equilibrium on the same vectors used to assemble the covariances.

For four-index arrays define the alternating trilinear contraction by
\[
 \mathcal D(P)=\sum_{i,j,k\in\{1,2,3\}^4}
 \Big(\prod_{s=1}^4\varepsilon_{i_sj_sk_s}\Big)
 P_{i_1i_2i_3i_4}P_{j_1j_2j_3j_4}P_{k_1k_2k_3k_4}.                 \tag{80}
\]
For $P=n^{\otimes4},Q=m^{\otimes4},R=l^{\otimes4}$ the polarized
expression equals $\det(n,m,l)^4$. Each alternating symbol is
rotation-invariant under proper rotations. The Fourier argument below
applies this identity when the three frequencies sum to zero.
The resulting null Lagrangian requires \(\mathcal D(F)=0\) at Reuss shear for fully symmetric \(F\), whereas direct contraction gives \(\mathcal D(F)=-1179/343\). Thus no physical arrangement attains the point, even in the homogenization closure.

\paragraph*{A quantitative physical gap}

For any compliance with fully symmetric anisotropy \(P=S-\Pi S\ne0\), write \(\Pi S=\alpha J+\beta D\). At every attainable isotropic pair, \[
 \boxed{\quad
  \alpha-\frac1{3\kappa}+5\left(\beta-\frac1{2\mu}\right)
   \ \ge\ \frac{|\mathcal D(P)|^2}{2^{128}\lambda_{\max}(S)\|P\|^4}.
 \quad}                                                               \tag{Q1}
\] Here \(\|\cdot\|\) on fourth-order tensors is the Hilbert-Schmidt norm and the alternating trilinear contraction is the one in (80). At the displayed crystal's Reuss bulk the first difference is zero. With $\lambda_{\max}(S)\le3/2$, (Q1) gives the explicit rational inequality
\[1-1/(2\mu)\ \ge\ \frac{463347}{25621683318964153144321637881433662918135644160}.\] In particular \(\mu>1/2+9\cdot10^{-42}\). The size is a conservative rigorous clearance.

To verify (Q1), take periodic cells whose effective compliances \(S_{\rm eff}\) tend to \(S_*\). Such periodic approximants generate the homogenization closure. Write \(\Delta=\operatorname{Tr}(\Pi S-S_*)\), \(\Delta_D=5(\beta-1/(2\mu))\) and \(\epsilon^2=\lambda_{\max}(S)\Delta\). Denote the rotated local compliance by \(S_R\), put \(P_R=S_R-\Pi S\), and let \(T\) project each nonzero output frequency onto transverse symmetric matrices. For a constant stress load \(q\) let \(\sigma\) be the exact equilibrated solution of mean \(q\). Its strain \(S_R\sigma\) has mean \(S_{\rm eff}q\) and compatible nonconstant modes. Orthogonality of that strain to \(q-\sigma\) and the energy upper bound give \[
 \big| (\langle S_R\rangle-S_{\rm eff})q\big|^2+\|T(P_Rq)\|_2^2
 \le \|S_R(q-\sigma)\|_2^2
 \le \lambda_{\max}(S)\,q:(\langle S_R\rangle-S_{\rm eff})q.
\] The first inequality keeps the mutually orthogonal constant and transverse nonconstant projections. Sum over all six loads. Rotation averaging \(\langle S_R\rangle\) gives \(\Pi S\), so the anisotropic mean \(\langle P_R\rangle\) is orthogonal to \(\Pi S-S_*\). Along the recovering sequence this gives \(\|\langle P_R\rangle\|^2+\|TP_R\|_2^2\le\lambda_{\max}(S)\Delta+o(1)\). On the five traceless loads the same proof bounds \(\|\langle P_R\rangle D\|^2+\|TP_RD\|_2^2\) by \(\lambda_{\max}(S)\Delta_D+o(1)\), as \(D\) is invariant under rotations. Moreover \(J\langle P_R\rangle J=0\) and self-adjointness gives \(\|\langle P_R\rangle\|^2\le2\|\langle P_R\rangle D\|^2\). These bounds account for the zero Fourier mode on arbitrary textures.

At a nonzero Fourier frequency rotate \(n\) to \(e_3\). Of the squared norm in components of \(P_R\) with at most two indices equal to \(3\), at least a sixth is on transverse output rows. This follows by counting their permutations using full symmetry. Double trace zero then bounds the squared norm of the \(3333\) component by eight times the norm squared of those components, by Cauchy-Schwarz on \(I_3\otimes I_3-n^{\otimes4}\). Keep just the components with exactly three indices \(3\). The resulting field \(Z\), with zero mean, satisfies \(\|P_R-Z\|_2^2\le\|\langle P_R\rangle\|^2+54\|TP_R\|_2^2\le54\lambda_{\max}(S)\Delta+o(1)\) and has symbol \(\partial_v n^{\otimes4}\), with \(v=(I-nn^T)(\widehat P_R:n^{\otimes3})\).

For a generic unit direction $\theta$, let $H_\theta$ have the
multiplier $-i\operatorname{sign}(\theta\cdot k)$ on nonzero frequencies
and zero mean. For real trigonometric $v$ every nonconstant frequency
of $v+iH_\theta v$ lies in one open half-space, so the mean of its
fourth power equals $\langle v\rangle^4$. Taking real parts, applying
Jensen to $v$ and then Hölder to the mixed term yields
$\|H_\theta v\|_4^4\le6\|v\|_4^2\|H_\theta v\|_4^2$.
The $j$th Riesz transform is the spherical average of
$2\theta_jH_\theta$, whose absolute weights integrate to one.
Density therefore gives the norm bound $\sqrt6$.

One has \(\|Z\|_4\le2^{24}\|P\|\). Indeed each spatial component is a sum of at most \(4(3^3+3^4)\) products of at most eight periodic Riesz transforms of entries of \(P_R\), with coefficients of absolute value one. Each transform has \(L^4\) norm at most \(\sqrt6\) by the directional estimate just proved. Summing the component bounds in the Euclidean norm gives the claimed constant. The mean of \(\mathcal D(Z)\) vanishes by the trigonometric-polynomial calculation in (80). At frequency triples summing to zero the directions are coplanar. Differentiating \(\det(n_1,n_2,n_3)^4\) once in each argument leaves a vanishing determinant in every term. Periodic Fejér convolution preserves the wave symbol and zero mean and converges to \(Z\) in \(L^3\), passing the identity to the limit. The trilinear form has norm at most \(6^4\) by summing its nonzero signed terms. Hölder and rotation invariance give \[
 |\mathcal D(P)|\le
 6^4(1+2^{24}+2^{48})\|P\|^2\sqrt{54}\,\epsilon
 \le 2^{64}\|P\|^2\epsilon .
\] Squaring and passing to the recovering limit proves (Q1).

The same projection gives a larger certified gap by expanding about the bounded crystal field. Both of the following estimates bound the same normalized projection error. Set \[
 x^2=\frac{\lambda_{\max}(S)}{\|P\|^2}\min(16\Delta,60\Delta_D).
\] Then \[
 \frac{|\mathcal D(P)|}{36\|P\|^3}
      \le 3x+3x^2+4500000000000\,x^{5/3}.                              \tag{Q1a}
\] To check the constants, the directional Hilbert transform has \(L^8\) norm at most \(\sqrt{29}\). In the even-power identity \(\Re\int(f+iHf)^8=(\int f)^8\) the nonpositive terms and Jensen's inequality yield \(y^8\le28y^6+28y^2\), \(y=\|Hf\|_8/\|f\|_8\), which excludes \(y\ge\sqrt{29}\). The spherical average giving each Riesz transform has total absolute coefficient one, so the same bound applies to it. The symbol count above now gives \(\|P_R-Z\|_8/\|P\|\le1+9\cdot4(3^3+3^4)29^4<3\cdot10^9\). The Euclidean norm of the coefficient tensor of the alternating trilinear form is \(36\), since it has \(6^4\) signed entries. Expand \(\mathcal D(P_R)-\mathcal D(Z)\) in one, two and three copies of \(P_R-Z\). The linear and quadratic terms after normalization give \(3x,3x^2\). Interpolation bounds the third absolute moment of the normalized error by \(x^{5/3}(3\cdot10^9)^{4/3}<4500000000000\,x^{5/3}\). This proves (Q1a). For this example \(|\mathcal D(P)|/\|P\|^3=1179/2744\). Substitution of \(x=(1/825)^3\) violates (Q1a). With \(\|P\|^2=1/4\) and \(\lambda_{\max}(S)\le3/2\), the \(60\Delta_D\) estimate gives \(\mu>1/2+8\cdot10^{-22}\) everywhere. At the Reuss bulk \(\Delta=\Delta_D\) and the \(16\Delta\) estimate improves this to \(\mu>1/2+3\cdot10^{-21}\). The calculator bisects \(x^{1/3}\) with exact sign tests at any supplied bulk.

Here is the projection calculation giving the two constants. At \(n=e_3\) subtract the components with exactly three indices equal to 3. The remaining fully symmetric double-trace-free tensor \(U\) has twelve independent components indexed by sorted quadruples with at most two indices equal to 3. Set \(U_{3333}=-\sum_{(i,j)\ne(3,3)}U_{iijj}\). In these twelve variables let \(A,T,M\) be the Gram matrices of \[
 \sum_{ijkl} U_{ijkl}^2,\qquad
 \sum_{i,j\le2}\sum_{kl} U_{ijkl}^2,\qquad
 \sum_{i,j\le2}\sum_{kl}\left(U_{ijkl}-\frac{\delta_{kl}}3\sum_p U_{ijpp}\right)^2 .
\] The two latter forms are the transverse output energy with all six input loads and with the five traceless loads, respectively. Both are positive definite. The monic characteristic polynomials of \(T^{-1}A\) and \(M^{-1}A\) are \[
 (t-6)^2(t-2)^4(t-1)^4(3t^2-53t+90)/3,\qquad
 (t-2)^4(t-1)^2(t^2-62t+135)(t^2-14t+18)^2 .
\] Building the Gram matrices from the displayed sums verifies these identities by rational arithmetic. Their largest roots are \((53+\sqrt{1729})/6<16\) and \(31+\sqrt{826}<60\). Rotation and Parseval give the two projection constants at all frequencies. The five uniform traceless stress trials bound the squared transverse residual by \(\lambda_{\max}(S)\Delta_D\), by the same energy improvement as for all six loads. Thus the shear-only estimate applies over the entire bulk interval for every fully symmetric compliance perturbation.

\paragraph*{A larger physical separation from two directional Hilbert transforms}

The same projection error \(x\) satisfies \[
 \frac{|\mathcal D(P)|}{36\|P\|^3}\le 3x+3x^2+21000\,x^{13/7}.\tag{Q1b}
\] For the first crystal \(x=(13/40)^7\) violates this inequality. Consequently \(\mu>1/2+4\cdot10^{-11}\) everywhere and \(\mu>1/2+1.5\cdot10^{-10}\) at the Reuss bulk. This part uses the physical stress equations and the explicitly constructed waves. The finite rank-one datum already proves the discrepancy between the physical minimum and the joint quadratic test.

To bound the order-zero projection in \(L^{16}\), put \(E=\operatorname{Sym}^4(\mathbb R^3)\) with its Frobenius inner product. Write \(P_n\) for the orthogonal projection of \(E\) onto the two tensors \(e_n(v)=\partial_v n^{\otimes4}/2\) as \(v\) runs through an orthonormal tangent frame. Its multiplier at nonzero frequencies gives \(Z\) above. Its zero multiplier is zero. The odd operator-valued polynomial \(B(n)=n\otimes P_n:E\to\mathbb R^3\otimes E\) has degree at most nine. It has a representation \[
 B(n)=\int_{S^2}\operatorname{sign}(n\cdot\theta)\,U(\theta)\,d\theta,
 \qquad \int_{S^2}\|U(\theta)\|_{\rm HS}^2\,d\theta=30336/49
\] with normalized surface measure. Indeed the integral operator with this sign kernel acts on harmonics of odd degree \(l\) by \(\lambda_l=\int_0^1 P_l(t)\,dt\), where \(P_l\) is the Legendre polynomial. Averaging over the azimuth at the pole proves the multiplier formula by rotation. For \(l=1,3,5,7,9\) the eigenvalues are \(1/2,-1/8,1/16,-5/128,7/256\). Divide each harmonic component of \(B\) by its eigenvalue to obtain \(U\). To check the displayed norm, at \(t=n\cdot m\) the overlaps of the two tangent frames for \(e_n,e_m\) are \(4t^4-3t^2\) in their common meridian and \(t^3\) in its perpendicular direction. This follows by differentiating \((n\cdot m)^4\) once in each variable and dividing by four. Thus \(\langle B(n),B(m)\rangle=k(t)=t[(4t^4-3t^2)^2+t^6]\). The squared norm on degree \(l\) is its Legendre coefficient \((2l+1)\int_0^1 k(t)P_l(t)\,dt\). Dividing by \(\lambda_l^2\) and summing gives the stated rational value.

In almost every direction \(H_\theta\) has scalar \(L^{16}\) norm at most \(11\). In the even-power mean identity used above, putting \(y=\|H_\theta f\|_{16}/\|f\|_{16}\) gives \[
 y^{16}\le120y^{14}+8008y^{10}+8008y^6+120y^2.
\] This excludes \(y\ge11\) by division by \(y^{14}\). Gaussian averaging of scalar projections proves the same bound for finite-dimensional Hilbert-valued functions. Integrating \(U(\theta)H_\theta\) now bounds the multiplier \(-i B(n)\) by \(11\sqrt{30336/49}\). The contraction multiplier \(i n\mathbin{\lrcorner}\) has norm at most \(22\), since \(n=2\int\operatorname{sign}(n\cdot\theta)\theta\,d\theta\). Their composition is \(P_n\). Therefore \(\|P_R-Z\|_{16}/\|P\|\le1+242\sqrt{30336/49}<6023\). Interpolating its \(L^2\) and \(L^{16}\) bounds gives a third absolute moment at most \(x^{13/7}6023^{8/7}\|P\|^3\). The alternating expansion for (Q1a) and \(6023^8<21000^7\) prove (Q1b).

Indeed $x^2\le360\Delta_D=1800(1-1/(2\mu))$ for this example.
The rational substitution in (Q1b) therefore gives
$1-1/(2\mu)>(13/40)^{14}/1800$, which implies the stated all-bulk
shear gap. The verifier evaluates the sphere integrals as rational
Legendre moments and checks the final strict comparisons.
\phantomsection\label{endfile:crystal-moment-proof}

\subsection{Isolated shear equality, axial drift and the degree obstruction}\label{detail:crystal-contact-proof}\label{file:crystal-contact-proof}
Choose a unit vector \(n\) in the material frame. Let \(U_n v=n\odot v\) and

\[
\Gamma_n=U_n(U_n^T C U_n)^{-1}U_n^T .
\]

For \(P=J,D\), define the following scalar function of \(x\) using the operator trace on \(V\),

\[
F_{P,n}(x)=
\operatorname{Tr}\big[P(C-xI_V)\big]
-\operatorname{Tr}\big[P(C-xI_V)\Gamma_n(C-xI_V)\big]
=z_{P,n}+\ell_{P,n}x-q_{P,n}x^2 .
\tag{E}
\]

Here \(z_{P,n}=\operatorname{Tr}[P(C-C\Gamma_n C)]\), \(\ell_{P,n}=2\operatorname{Tr}(P\Gamma_n C)-\operatorname{rank}P\), and \(q_{P,n}=\operatorname{Tr}(P\Gamma_n)\). The constant and quadratic coefficients \(z_{P,n}\) and \(q_{P,n}\) are strictly positive. In fact \(C-C\Gamma_n C\) is positive semidefinite with kernel \(\operatorname{ran}U_n\). The hydrostatic strain is outside that three-dimensional space and the five-dimensional deviatoric space cannot be contained in it. For \(q_{P,n}\), \(U_n^T e=n/\sqrt3\), and \(\operatorname{ran}U_n\) also contains deviatoric shear.

\textbf{Construction.} Assign any probability weights to any finite list of normals. For \(P=J,D\), average the three coefficients in (E) with those weights and put

\[
 x_P=\frac{\bar\ell_P+\sqrt{\bar\ell_P^2+4\bar q_P\bar z_P}}{2\bar q_P}.
 \qquad
 \boxed{\ (\kappa_*,\mu_*)=(x_J/3,x_D/2)\text{ is attainable in the limit.}\ }
\tag{F}
\]

Equivalently, a pair \((k,m)\) belongs to this constructed set exactly when the convex hull of the planar vectors \((F_{J,n}(3k),F_{D,n}(2m))\), \(|n|=1\), contains zero. The same weights must work in both coordinates. On any prescribed list \textbf{three normals suffice} by Carathéodory's theorem. If directions range over the sphere, \textbf{two suffice for every pair in (F)}. Starting with a strict triangle about zero, move its first normal continuously to the second. At the first loss of strict containment, zero is on a closed edge or the triangle has collapsed to a segment. The test uses rational functions of the normal components with strictly positive acoustic determinants. It can therefore be written as a finite real polynomial feasibility problem. This gives a sufficient test for the full attainable region.

\textbf{Both endpoints at a prescribed bulk from pairs of normals.} For a finite normal list put \(h_i=F_{J,n_i}(3k)\). All admissible weights at this \(k\) form the polytope \(\sum_i w_i=1\), \(w_i\ge0\), \(\sum_i h_i w_i=0\). Its vertices are the individual normals with \(h_i=0\) and the pairs \(i,j\) with \(h_i<0<h_j\), with weight \(-h_i/(h_j-h_i)\) on \(j\). Indeed a feasible weight with at least three positive entries can be perturbed in both directions preserving both equalities. For a fixed shear argument the weighted shear polynomial is linear in the weights. At the least positive root among the vertices all these polynomials are nonnegative, and at the greatest root they are nonpositive. Hence the two vertex roots give the exact endpoints for this finite-normal construction. Continuity on the polytope fills the interval between them. This requires at most quadratically many tests in the length of the list. It also shows that the extrema of the full construction (F) at a fixed bulk require at most two normals, since any attainable pair in (F) first reduces to at most three normals.

\textbf{Proof of attainability.} Start with an isotropic tensor \(B=aJ+bD\) in the laminate closure of the crystal, supplied by the iteration above. Put a fraction \(f\) of the crystal into a simple laminate with \(B\), with layer normal \(n\). Matching tangential strain and normal traction, or minimizing the quadratic energy over the jump \(U_n v\), gives the exact tensor

\[
\mathcal L_f(B,C,n)=B+f(C-B)
 -f(1-f)(C-B)U_n[U_n^T((1-f)C+fB)U_n]^{-1}U_n^T(C-B).
\tag{G}
\]

In particular, to first order in \(f\) the increment is \(f[(C-B)-(C-B)\Gamma_n(C-B)]\). The remainder is \(O(f^2)\) uniformly in \(n\) and \(B\) between the ellipticity bounds of \(C\).

At a coarser scale laminate all 60 rotations of each of these tensors, with its normal weight divided equally over the rotations. The arithmetic and harmonic bounds of this laminate are isotropic. If the local tensors have expansions \(B+fH+O(f^2)\), their arithmetic mean is \(B+f\langle H\rangle+O(f^2)\). Expanding the inverse gives the same first-order term for their harmonic mean. Apply the preceding isotropization iteration to this coarse laminate. Its isotropic limit \(B'\) obeys the same two order bounds. All these constructions stay in the laminate closure of the original crystal. Since \(\Pi M=(e:M e)J+\operatorname{Tr}(DM)D/5\), taking the hydrostatic entry and the deviatoric trace in (G) now gives the decoupled updates

\[
 a'=a+f\,\overline{F_{J,n}}(a)+O(f^2),\qquad
 b'=b+\frac f5\,\overline{F_{D,n}}(b)+O(f^2).
\]

Iterate with \(f_j=1/(j+2)\). The two quadratics have one positive root and one negative root. Dividing each drift by the rank of \(P\), it has the form \(c(x)(x_P-x)\) with \(c(x)>0\) bounded above and bounded away from zero on the compact ellipticity interval. Thus the error in either scalar obeys \(d_{j+1}\le (1-c_0 f_j)d_j+K f_j^2\) for all sufficiently large \(j\). For any \(\epsilon>0\), eventually \(K f_j<c_0\epsilon\). Iterating the inequality and using \(\prod_{j\ge j_0}(1-c_0 f_j)=0\) gives \(\limsup d_j\le\epsilon\). Both coordinates converge to (F). Each finite tensor belongs to the closed attainable set, so a diagonal sequence of actual microstructures realizes the limit. The root is a consequence of this construction and does not assume a self-consistent effective-medium law.

For a transversely isotropic test material and its axial normal, both roots (F) are also forced on every texture sequence that attains either benchmark (122) or (123).

\paragraph*{An attainable path to each quadratic root}

Let the initial isotropic host be any attainable pair. Write \(x_J=3k\), \(x_D=2m\). For each averaged quadratic in (E) let \(r_P>0>v_P\) be its two roots, let \(\gamma_P=\bar q_P(r_P-v_P)/\operatorname{rank}P\), and put \(g_P(x)=(x-r_P)/(x-v_P)\). Then the entire path

\[
 g_P(x_P(t))=e^{-\gamma_P t}g_P(x_P(0)),
 \qquad P=J,D,\qquad 0\le t<\infty                    \tag{68}
\]

is attainable in the closure. For fixed \(t\), perform \(N\) of the preceding laminate updates with fraction \(t/N\), starting from the prescribed host. The remainders are uniformly \(O((t/N)^2)\), and the drifts are Lipschitz on the ellipticity interval. Iterating the error inequality for Euler's method gives an \(O_t(1/N)\) error from the solutions of \(\dot x_P=\overline{F_P}(x_P)/\operatorname{rank}P\). Separating variables gives (68), including a stationary coordinate when its initial value is \(r_P\). The limit at infinite time is (F).

This yields bounds at a \textbf{specified bulk}. Start with either sharp Hill bulk and an interval of attainable host shears there. On the side between this bulk and \(r_J/3\), the requested bulk determines \(u=e^{-\gamma_J t}=g_J(3k)/g_J(3k_{\rm host})\). Apply \(g_D^{-1}(u^{\gamma_D/\gamma_J}g_D(2m_{\rm host}))/2\) to both host shear endpoints. The resulting interval is attainable, since this map is increasing in the initial shear. If the host bulk already equals \(r_J/3\), the constant-bulk path joins its shear to \(r_D/2\). For a crystal with a certified axial shear minimum, (68) also bounds the rate of departure of the actual sharp boundary.

\paragraph*{Covering every prescribed bulk modulus}

Compute the two tensors \(L_-,L_+\) above. If the endpoints coincide, the original crystal already has \(e\) as an eigenvector, by equality in \((e:C e)(e:C^{-1}e)\ge1\). Otherwise, fix any joining normal \(n_0\) and form \(M_f=\mathcal L_f(L_-,L_+,n_0)\). The continuous scalar \(e:M_f e/3\) has endpoint values \(\kappa_R,\kappa_V\). For any prescribed \(k\) between them choose a root of \(e:M_f e=3k\). Formula (G) makes this a univariate equation of degree at most four after multiplying by the positive acoustic determinant. Apply the two Voigt swaps to \(M_f\). This produces a finite laminate \(A\) with \(A e=3k e\).

\paragraph*{Shear equality and a change of material constants}

Let \(I=I_3\), \(e=I/\sqrt3\), \(J=e\otimes e\), \(D=I_V-J\). Use the stiffness \(L=C\) for the upper Hill shear or the compliance \(L=C^{-1}\) for the lower Hill shear. Write \(\Pi L=aJ+bD\) for its rotation average and \(F=L-\Pi L\). Define the traceless spatial matrices

\[
 B_F=FI,\qquad H_F=N(F),\qquad N(F)_{ij}=\epsilon_{ikl}\epsilon_{jrs}F_{kr,ls}.
\]

For a sequence of periodic orientation fields \(R\), let \(F_R\) be the grain rotation of \(F\). Approaching the respective Hill shear with isotropic effective response, allowing subsequences for the bulk, is equivalent to

\[
   \langle F_R\rangle\longrightarrow0,\qquad
   \|\mathscr P(F_R D)\|_2\longrightarrow0 .                 \tag{70}
\]

Here \(\mathscr P\) acts columnwise and projects onto mean-zero compatible strains on the stiffness side, or mean-zero equilibrated stresses on the compliance side. Necessity follows from the uniform trial deficits on the five shear loads. Their vanishing gives the second condition by ellipticity, and the mean trial responses on those loads approach \(bD\) by the same energy estimate. Major symmetry and \(e:F_R e=0\) then give the first condition. Conversely these two limits make the five uniform trials approach their minimizing fields in \(L^2\). Their mean responses tend to \(bD\). By self-adjointness every effective limit is isotropic, with this shear eigenvalue.

Both conditions concern only the rotated anisotropy. Thus, if a Hill shear is attainable at one nonzero contrast, the \emph{same orientations} approach the corresponding Hill shear for every positive-definite \(L_z=a'J+b'D+zF\) with \(z\ne0\), including the opposite sign of \(z\). The new hydrostatic response depends on the constants and the geometry. This statement uses compliance contrasts when transferring a lower Hill shear.

\paragraph*{Two isolated equality families.}
For Reuss equality the compliance in Kelvin coordinates is
$S=aJ+bD+rF_0$, where the axis of
\[
 F_0=\begin{pmatrix}1&1&1\\1&1&1\\1&1&-8\end{pmatrix}
       \oplus\operatorname{diag}(3,3,0)
 \tag{75}
\]
rotates with each grain. The asserted equality point is
\[
 \left\{\left(\frac1{3x},\frac1{2b}\right)\right\},\qquad
 a-x=\frac{2(a-x+3r)^2}{2a+b+6r},\quad x>0.                         \tag{122}
\]

There is also an upper-shear family. In Kelvin order \(11,22,33,23,31,12\), take axis \(u=e_3\) and put

\[
 V_0=\frac18\begin{pmatrix}-7&-13&8\\-13&-7&8\\8&8&8\end{pmatrix}
           \oplus\operatorname{diag}(0,0,3/4).
\]

For every positive stiffness \(C=aJ+bD+rV_0(u)\), \(r\ne0\), the entire horizontal section at the Voigt shear is

\[
 \left\{\left(\frac{x}{3},\frac b2\right)\right\},\qquad
 a-x=\frac{(a-x+3r)^2}{a+2b+3r},\quad x>0 .                          \tag{123}
\]

Both sections lie strictly inside the Hill bulk interval.

\paragraph*{The uniform kernel estimate used at equality.}
On elasticity tensors $G$ with $e:Ge=\operatorname{Tr}G=0$ write
$B_G=GI$ and $N(G)_{ij}=\epsilon_{ikl}\epsilon_{jrs}G_{kr,ls}$.
At a nonzero mode of normal $n$ let $T_n M=(I-nn^T)M(I-nn^T)$.
The two tests are the finite-dimensional maps
$G\mapsto(I-T_n)GD$ and $G\mapsto T_n GD$. A uniform constant
bounds distance to either kernel by the norm of that test. Indeed at
$n=e_3$ restrict to the orthogonal complement of the kernel and invert
its positive least singular value. Every other normal is obtained by
an orthogonal change of coordinates. Parseval gives the corresponding
$L^2$ estimate. In the first kernel $G_{33,kl}=w\delta_{kl}$ and
$G_{i3,kl}=0$ for $i=1,2$, by major symmetry and the five trace-free
loads. Direct contraction gives $N(G)=2w(I-3nn^T)$. In the second,
$G_{ab,kl}=w\delta_{ab}\delta_{kl}$ for transverse indices and
the two scalar trace conditions give $G_{33,33}=-8w$ and
$\operatorname{tr}(G_{3a,3b})=3w$. The contraction identities for
$B_G$ and $N(G)$ in the scalar-mode reduction follow. For the marker
(75) the identity $B_{F_0}=-3N(F_0)$ kills the transverse vector and
trace-free block. The resulting mode is $wF_0(n)$.
For the elastic energy estimates used below, write $\sigma=L E+r$,
where the compatible $E$ and equilibrated $\sigma$ have prescribed
means. Subtract the true cell solution with the same mean strain.
Testing its equilibrium with the compatible difference gives
$\|L^{1/2}(E-E_{\rm cell})\|_2\le
\|L^{-1/2}r\|_2$. In particular a vanishing residual forces the
mean stress to converge to the claimed effective response. The same
argument with stresses and strains exchanged applies to compliance.

\paragraph*{What the attaining orientations determine}

Consider any sequence of grain axes that approaches either (122) or (123) in shear, with isotropic effective limit. Reassign \textbf{any} positive-definite transversely isotropic stiffness \(L\) to these same axes. Both effective moduli then converge to the two axial roots (F) for \(L\), and the full effective tensor converges to the corresponding isotropic tensor. This determines the response to changes in all five elastic constants from a single shear equality. The proof allows arbitrary grain orientations. Periodic recovery can be applied simultaneously to finitely many test assignments by using the same periodized cubes and diagonal sequence for their strain and stress energies. For a spatially varying subsequential test limit, first localize the joint limits at almost every point and then periodize there. The isotropic benchmark is constant, so the argument below gives the same test tensor at each such point.

Use periodic approximating cells and normalized volume averages. Write \(U=uu^T\), \(\Phi=u^{\otimes4}\) and \(A=\fint_{S^2}v^{\otimes4}\,dv\), with \(A_{ijkl}=(\delta_{ij}\delta_{kl}+\delta_{ik}\delta_{jl}+\delta_{il}\delta_{jk})/15\). We first show that there is a real scalar field \(t\) of mean one for which \(\Phi-tA\) approaches a zero-mean tensor field whose Fourier coefficients at \(\xi\ne0\) are scalar multiples of \(n^{\otimes4}\), \(n=\xi/|\xi|\). All approximation statements here concern sequences of cells. Distances are in \(L^2\) unless stated otherwise.

\begin{itemize}
\tightlist
\item
  At Reuss equality the scalar-mode reduction in the preceding proof, with equal degree-two and degree-four amplitudes \(r\), gives this directly. Take \(t=1-w/r\) in that formula.
\item
  At Voigt equality the shear trial estimate (70) and the uniform kernel estimate above place the rotated marker \(V_0(u)\) within vanishing distance of the zero-mean Voigt kernel fields. In a mode with axis \(n\), the kernel has \(N(G)=2w(I-3nn^T)\). Its rotation average about \(n\) is exactly \(wV_0(n)\). Indeed \(G_{nn,kl}=w\delta_{kl}\), the mixed rows containing \(n\) vanish, and the two scalar trace constraints fix the axial part on the transverse plane. Since \(N(V_0(u))=2(I-3uu^T)\), Parseval gives \(\langle w^2\rangle\to1\). The scalar projections therefore use the full limiting squared norm \(\|V_0\|^2\), so the orthogonal remainder tends to zero. Both the degree-two and the degree-four rotational parts of \(V_0(u)\) are nonzero. Its quartic on the unit sphere is \((-7+30s-15s^2)/8\), \(s=(u\cdot v)^2\), and its \(N\) contraction recovers \(U-I/3\). Extracting the two parts linearly and equivariantly gives \(\widehat{\Phi-A}\simeq \widehat w(n^{\otimes4}-A)\) in squared-sum norm. Use \(t=1-w\), with \(w\) defined by the longitudinal contraction of the actual rotated \(N\) on the nonzero modes and zero at the mean.
\end{itemize}

In either case \(t\) is one plus bounded Fourier multipliers applied to bounded orientation components. By the periodic multiplier theorem \cite[Thm.~2.1]{p13:riesz}, its norms in every fixed finite \(L^p\) are bounded uniformly along the cells. The orthogonal projection onto the above Fourier lines has the same property. Its entries are products of Riesz multipliers \(-i\xi_j/|\xi|\). For a direct proof of their finite-\(p\) boundedness, the directional Hilbert transform argument in the preceding cubic-gap proof works at every even exponent \(2m\). Expand the real part of \(\langle(f+iH_\theta f)^{2m}\rangle=\langle f\rangle^{2m}\) and apply Hölder to all terms with fewer than \(2m\) copies of \(H_\theta f\). This bounds \(\|H_\theta f\|_{2m}/\|f\|_{2m}\) by a constant depending only on \(m\). Spherical averaging gives the Riesz bounds. This argument starts on real trigonometric polynomials in a generic direction and extends by density. Interpolation upgrades every vanishing \(L^2\) error below to a vanishing error in each fixed finite \(L^p\) when all the higher norms are bounded.

\paragraph*{A hierarchy of fourth derivatives}

For every integer \(j\ge0\) the fields

\[
 Z_j=t^j\Phi-\frac{t^{j+1}}{j+1}A                                      \tag{124}
\]

approach zero-mean fields with Fourier coefficients on \(\mathbb R n^{\otimes4}\) (with complex scalars for Fourier modes). Such fields are fourth derivatives of scalar periodic potentials. The mode estimate just proved yields the case \(j=0\).

Flatten a symmetric fourth-order tensor into its columns in \(\mathbb R^3\), indexed by symmetric triples. Denote a column by a subscript \(p\). Every column of a fourth derivative is a gradient. The cross product of two periodic gradients of zero-mean potentials is divergence free and has mean zero, by integration by parts. Smooth approximation gives this identity for the fields used here. Apply it to the symmetrized cross products of columns of \(Z_0,Z_j\). Their projected approximants are gradients, and the \(L^4\) error makes the product error vanish in \(L^2\). All columns of \(\Phi\) are multiples of \(u\), so their cross products vanish pointwise. Expansion leaves a nonzero constant times

\[
             A_p\times K_q+K_p\times A_q,\qquad K=Z_{j+1}.
\]

Thus these vector fields have vanishing longitudinal components in Fourier norm and vanishing means.

Here is the symbol calculation that recovers \(K\). At a fixed unit \(n\) let \(P\) denote spatial projection onto \(n^\perp\). Suppose \(n\cdot(A_p\times K_q+K_p\times A_q)=0\) for all \(p,q\), with \(K\) fully symmetric. The projected columns of \(A\) span the plane. For \(q\in\ker(PA)\) the identity forces \(PK_q=0\), so \(PK=M PA\) for a planar matrix \(M\). The induced action on planar signed area is multiplication by \(\operatorname{tr}M\), hence \(\operatorname{tr}M=0\). In axes with \(n=e_3\), symmetry of \(K_{1233}=K_{2133}\) gives \(M_{12}=M_{21}\). Symmetry of \(K_{1112}=K_{2111}\) then gives \(M_{12}=3M_{21}\). Finally \(K_{1122}=K_{2211}\) gives \(M_{11}=M_{22}\). Therefore \(M=0\) and \(K\) is a multiple of \(n^{\otimes4}\). Conversely every such multiple passes. This finite-dimensional calculation gives a norm estimate modulo that line, uniform in \(n\) by rotation. At the mean all directions can be used, and the intersection of the lines is zero. Parseval proves the assertion for \(Z_{j+1}\). Polynomial powers have uniformly bounded norms at each finite exponent by the bounds on \(t\). This completes the induction.

Double trace in the mean of (124) gives \(\langle t^{j+1}\rangle\to(j+1)!\) by induction. In particular the distribution of \(t\) converges to the exponential law with density \(e^{-s}\) on \(s\ge0\). Indeed uniform higher moments give the moments on every subsequential limit. Its even moments give a finite integral of \(\cosh(cs)\) for \(0<c<1\), so its moment generating function near zero uniquely identifies the law. The tensor mean in (124) also fixes the orientation moments weighted by each power. For every polynomial \(p\), write \(Ip(s)=\int_0^s p(v)\,dv\). Then

\[
 \langle p(t)\Phi\rangle\longrightarrow A\int_0^\infty p(s)e^{-s}\,ds,
 \qquad p(t)\Phi-Ip(t)A\quad\text{approaches a zero-mean fourth derivative}. \tag{125}
\]

\paragraph*{Testing all five constants on the same geometry}

If \(W(n)\) is an even transversely isotropic map on symmetric loads, with symmetry axis \(n\), its rotated values depend linearly on \(n^{\otimes4}\). Explicitly its degree-zero, degree-two and degree-four parts are fixed linear images of the corresponding parts of \(n^{\otimes4}\). The statement also holds for maps without major symmetry since products of \(nn^T\) and Kronecker tensors give the axial maps used below. Write \(\bar W\) for its spherical average. Apply that linear map to (125). For \(p(t)W(u)-Ip(t)\bar W\) the limiting fluctuations are compatible if the columns of \(W(n)\) are compatible waves with normal \(n\). They are equilibrated if the columns are transverse stresses. In fact the double trace of the tensor in (125) is \(p(t)-Ip(t)\), so at a nonzero mode the map-valued field approaches \(W(n)\) times the Fourier coefficient of that scalar. Its mean tends to zero.

Test a stiffness \(L\) and let \(P\) be either the hydrostatic or the deviatoric projector. Let \(\beta\) be the positive axial root (E) for \(P\). With \(\Gamma_{L,n}\) the inverse compatible-wave compression, put

\[
 X=P-\Gamma_{L,n}(L-\beta I_V)P,\qquad Y=LX,\qquad
 d=\operatorname{Tr}(PX)/\operatorname{rank}P .
\]

Here \(L\) is in the frame with axis \(n\), and the maps rotate with \(n\). Their means are \(dP\) and \(\beta dP\) by rotational invariance and the root equation. The columns of \(X-P\) are compatible and those of \(Y-\beta P\) are equilibrated. The derivative of the axial drift divided by \(\operatorname{rank}P\) is \(1-2d\) at its positive root. This derivative is strictly negative because the quadratic has roots of opposite signs. Thus \(d>1/2\).

For \(\rho=(d-1)/d\), use the ordinary Laguerre polynomials \(L_i(s)=\sum_{j=0}^i(-1)^j\binom ij s^j/j!\), and take \(p_N=d^{-1}\sum_{i=0}^N\rho^i L_i\). One has \(IL_i=L_i-L_{i+1}\) and \(\int_0^\infty L_i^2e^{-s}ds=1\) by integration by parts using the Rodrigues formula. In particular

\[
 1-p_N-(d-1)Ip_N=\rho^{N+1}L_{N+1}.                                  \tag{126}
\]

Add to \(P\) the compatible fluctuations from \(X-P\) with polynomial \(p_N\), and add to \(\beta P\) the equilibrated fluctuations from \(Y-\beta P\) with that same polynomial. Up to errors vanishing for each fixed \(N\), the result is the strain map \(p_N X+\rho^{N+1}L_{N+1}P\) and the stress map \(p_N Y+\rho^{N+1}L_{N+1}\beta P\). After projecting onto the exact fluctuation spaces, the limiting squared norm of their constitutive residual is

\[
          |\rho|^{2N+2}\,\|(L-\beta I_V)P\|_{\rm HS}^2 .                \tag{127}
\]

The local norm in the last factor is rotation invariant. For exact trials of these means the compatible strain and equilibrated stress residual estimate above bounds the squared error on the effective response to \(P\) by \(\lambda_{\max}(L)/\lambda_{\min}(L)\) times the squared constitutive residual. As \(|\rho|<1\), taking \(N\to\infty\) proves the full effective limit on both load spaces. No rate of approach to shear equality is required.

Apply this transfer to the original stiffness. Compliance and stiffness roots of the axial construction are reciprocal, by differentiating the inverse layer formula (G). For (75) the complementary hydrostatic root is the positive root in (122). For \(V_0\) the compatible projection of \(V_0D\) vanishes about its axis, so the positive shear root in stiffness is \(b\). The longitudinal acoustic stiffness is \((a+2b+3r)/3\). On hydrostatic loading the squared correction in (E) is \((a-x+3r)^2/(a+2b+3r)\), giving (123). Construction (F) realizes both. The hydrostatic coupling is nonzero in the crystal and in its inverse. At either hydrostatic arithmetic bound the respective axial drift is strictly negative because the axial deviatoric response to a hydrostatic load has a nonzero projection on both the longitudinal strain and plane-trace stress wave spaces. Both bulks in (122)-(123) are consequently interior.

\paragraph*{An isolated extremum with arbitrarily high order of contact}

Write \(k_0,m_0\) for the pair (122) or (123). At every \(k\ne k_0\) in the Hill bulk interval the sharp shear gap \(g(k)\) is strictly positive. Here \(g(k)=m_-(k)-m_0\) on Reuss and \(g(k)=m_0-m_+(k)\) on Voigt. This follows from uniqueness and compactness at each prescribed bulk. The axial paths (68) give

\[
 0<g(k)\le C_\pm |k-k_0|^\alpha,\qquad
 \alpha=
 \begin{cases}
 \sqrt{(2a+b+6r)/(2a+b-18r)},&\text{Reuss},\\
 \sqrt{(a+2b+3r)/(a+2b-9r)},&\text{Voigt}.
 \end{cases}                                                           \tag{128}
\]

The constants can differ on the two sides. To verify the exponent, use the drift of the axial shear eigenvalue in its own energy convention. The uniform shear trial is exact on the axial wave at \(b\), so the derivative of its normalized drift there is \(-1\). Write the hydrostatic drift as \(a-x-c(a-x+3r)^2/Q\) with \((c,Q)=(2,2a+b+6r)\) or \((1,a+2b+3r)\). The magnitude of its derivative at the positive root is \(\sqrt{1-12cr/Q}=1/\alpha\). These square roots are strictly positive by ellipticity and simplicity of the positive root. The ratio of the two decay rates in (68) gives (128) from any host at each of the two Hill bulks. The reciprocal changes of variable for compliance have nonzero derivatives. Compactness of the bulk interval extends the bound to the two closed sides.

\paragraph*{Exact algebraic descriptions require unbounded degree on rational crystals}

For each integer \(N\ge2\) there is a \textbf{rational Kelvin crystal} whose sharp minimum-shear graph, near a specified rational interior bulk, cannot be described by any finite list of polynomial inequalities of bulk degree less than \(N\). The same assertion holds for the maximum-shear graph on another rational crystal. Arbitrary real coefficients and any number of inequalities are allowed. In particular no single semialgebraic formula in the crystal entries and two moduli describes the answer to Milton's question. This is an obstruction even with exact rational input. For a prescribed positive error, (S13) gives semialgebraic programs whose sizes grow with the requested accuracy.

Put \(\eta=12/(N^2-1)\). In the two families of (122)-(123) take the respective compliance and stiffness \[
 S_N=(3+\eta)J+12D+F_0,\qquad C_N=(3+\eta)J+3D+V_0.
\] The matrices \(F_0,V_0\) in those equations have fixed rational Kelvin entries. Both tensors are positive definite. Their isolated sharp shear points, in physical units, are \[
 (\kappa_0,\mu_0)=\big((N-1)/18,\,1/24\big)\quad\text{for the minimum},\]\[
 (\kappa_0,\mu_0)=\big(2/(N-1),\,3/2\big)\quad\text{for the maximum}.
\] In fact the hydrostatic root in the selected energy convention is \(6/(N-1)\) in both examples. The exponent in (128) is \(\sqrt{1+12/\eta}=N\). Thus at bulks \(\kappa_0+h\) on either side near zero, the corresponding boundary equals \(\mu_0\pm g(h)\) with \(0<g(h)\le C|h|^N\). The strict inequality at every nonzero \(h\) is essential.

Suppose finitely many polynomial sign tests describe the section or its boundary graph. Discard polynomials vanishing identically after the crystal has been specified. At every boundary point at least one remaining polynomial must vanish, since otherwise all signs are constant nearby. Translate the origin to \((\kappa_0,\mu_0)\), write \(z=\mu-\mu_0\), and remove the highest power of \(z\) dividing each polynomial. This leaves their zero sets at \(z\ne0\) unchanged. For a resulting polynomial \(P\) its univariate restriction \(P(h,0)\) is nonzero and has vanishing order at most its bulk degree \(d\). It has magnitude at least \(c|h|^d\) for sufficiently small nonzero \(h\). On a bounded neighborhood \(|P(h,z)-P(h,0)|\le C_P|z|\). If \(d<N\) then \(P(h,\pm g(h))\) cannot vanish for sufficiently small nonzero \(h\). This contradiction proves the degree lower bound.

Taking \(N=2^m\) makes the input length \(O(m)\) in binary even when the Reuss example is supplied as stiffness, since inversion has fixed matrix dimension. Consequently a quantifier-free exact planar representation, whenever one exists, can require bulk degree exponential in the input length. Quantifier elimination over real closed fields turns any fixed finite semialgebraic formula with auxiliary variables into a quantifier-free formula of fixed degree before specifying the crystal. The same obstruction therefore applies to a fixed finite list of matrix inequalities of fixed size with semialgebraic data and auxiliary variables.
\phantomsection\label{endfile:crystal-contact-proof}

\subsection{Calculation at fixed nominal zinc stiffness}\label{route:calculation-at-fixed-nominal-zinc-stiffness}
Each inner endpoint is furnished by two auxiliary textured materials
with bulk moduli on opposite sides of 66 GPa. For host $(P,g)$ and
core $(A,b)$ the core share and the resulting shear are
\[
 t=\frac{(K-P)(3A+4g)}{(A-P)(3K+4g)},\qquad
 M=g+\frac{t(b-g)}{1+(1-t)\frac25(1/g+2/(3P+4g))(b-g)}.
\]
At $K=66$ the inner witnesses, in GPa, are
\[37.62336468848\ldots,\qquad 43.35521987238\ldots .\]
The
minimum uses two 16,384-grain trees and the maximum two 4,096-grain
trees. Their geometry can be scale-separated. The two outer comparisons
use the positive-semidefinite material tests for specified translations
and the Schur complement of their macroscopic two-load block. Tests at
the exact bulk leave the outer values $36.80014385$ and $44.07248692$
GPa after outward rounding. Interpolation at fixed bulk gives every
shear between the inward endpoints. The tree files, input crystal and
translation coefficients are bundled with the fixed-bulk replay.

\SSection{common}{Complete bulk sections and cubic constraints}
\subsection{Common shears across all bulks}\label{route:common-shears-across-all-bulks}
\subsection{Edge shears of the two constructions}\label{detail:dimension-tools}\label{file:dimension-tools}
For the step from the scalar comparisons to physical moduli use the
simultaneous strain and complementary-energy bounds for an isotropic
aggregate. The prescribed bulk is independent of the core shear when
the second angular field vanishes on the core. Therefore changing
the fourth field on that core moves only its shear interval. The
measure on layer directions below is weighted by the increments in
the surviving core fraction.
\paragraph*{The construction at arbitrary fractions}

Let the core \((l,n)\) occupy fraction \(f\), put \(p=1-f\), and write \((k,m)\) for the matrix. In dimension \(d=2,3\) assume \(0<f<1\) and positive moduli with \(k\ne l\). With \(H_f(a,b;t)=[f/(a+t)+p/(b+t)]^{-1}-t\), the exact bulk throughout the construction is \[
 K_s(f)=H_f(l,k;2(d-1)m/d).
\] It is the minimum for \(n>m\) and the maximum for \(n<m\), for either ordering of the bulks. In the reference material with both shears \(m\), the deviatoric strain in a hydrostatic test is proportional to the field \(Q\) with nonconstant Fourier symbol \(NN^T-I/d\) acting on the centered core indicator \(X=\chi-f\). If \(Q=0\) on the core, changing its shear leaves that solution and its stresses unchanged. In particular hydrostatic strain is an eigenvector of the effective tensor, even before isotropizing.

Follow a branch in a separated-scale laminate from coarse to fine. A split of normal \(N\) that changes the conditional core fraction from \(x\) to \(y\) adds \((y-x)(NN^T-I/d)\) to \(Q\). More generally a multiplier with direction symbol \(U(N)\) gains \((y-x)U(N)\). This follows by averaging the finer indicator at each coarser scale, or by applying the multiplier to its one-direction Fourier modes.

Start a step at fraction \(x\) with accumulated \(Q=0\) and choose \(0\le b<x\). In local orthogonal axes write \(R_i=y+dQ_{ii}\) for the current fraction \(y\), initially all \(R_i=x\). A split normal to \(e_i\) changes only \(R_i\) and changes it by \(d\) times the change in \(y\). In succession for \(i=1,\ldots,d\), split \(R_i\) from \(x\) to \(1\) or \(b\), with probabilities \(w=(x-b)/(1-b)\) and \(1-w\). Continue to the next axis only on the low branch. If all coordinates land at \(b\), that branch has fraction \(b\) and \(Q=0\) and can start another step in any frame. When the first high branch occurs, let \(h=d^{-1}\sum R_i\) and finish with a sequential laminate of local core fraction \(h\). Give the fine normal \(e_i\) weight \((1-R_i)/[d(1-h)]\). On a final core point the fraction increases by \(1-h\) in total. The fine increments take every \(R_i\) to \(1\), so the final \(Q\) vanishes there. All conditional fractions and weights lie in \([0,1]\).

On any such stopping core the total fraction increment in direction \(i\) from the start of the step is \((1-x)/d\). On the surviving branch it is \((b-x)/d\). After several steps down from \(f\) to \(x\), the volume of the surviving branch and the mass of core it contains are, respectively, \[
 V(x)=\left(\frac p{1-x}\right)^d,\qquad F(x)=\frac{p^d x}{(1-x)^d}.        \tag{2}
\] The last surviving branch at \(b=0\) is pure matrix. Scale separation in each finite tree can be taken before refinement of the partition of \([0,f]\).

\paragraph*{Concentrating the fourth angular field}

Use the unit fourth-harmonic symbol \(U_4(N)=e^{4i\theta}\) in the plane and \(U_4(N):U_4(N')=P_4(N\cdot N')\) in space. Let \(z\) be the fraction of the fourth field's energy in the core, equal to \(\langle\chi|U_4(\nabla)X|^2\rangle/(fp)\). It is denoted by \(z_2\) in the planar angular coordinates \((z_1,z_2)\) and by \(z_4\) in the spatial coordinates \((z_2,z_4)\). The second harmonic vanishes on the core in both cases.

In one step put \(D=d^{-1}\sum_i U_4(e_i)\). The fourth field on a stopping core gains \((1-x)D\) and on the survivor gains \((b-x)D\). In the plane \(|D|^2=1\). In space \(|D|^2=[1+2P_4(0)]/3=7/12\). Isotropize the subtree in each surviving branch on scales between successive steps. Independent uniform rotations of the subsequent frame give mean zero for its \(D\); orientation does not affect the masses in (2). The mean squared size of accumulated survivor increments is at most \(|D|^2\sum (x-b)^2\), which tends to zero with the mesh. One may equivalently distribute the orientations in deterministic finer cells with the same fraction. Coarse rotation aggregates of such cells have zero additional increment of \(X\). Consequently the limiting core energy divided by \(fp\) is \((fp)^{-1}|D|^2\int_0^f (1-x)^2 F'(x)\,dx\). This yields the attainable values \[
 \begin{aligned}
 Z_{\rm pl}(f)&=\frac p f[-f-2\log p],\\
 Z_{\rm sp}(f)&=\frac7{12}\left[3p+\frac{2p^2}{f}\log p\right].           \end{aligned}\tag{3}
\] In the plane \(Z_{\rm pl}(f)>1-f^2\) for all \(0<f<1\), since \(-2\log(1-f)>2f+f^2\). At half fractions the construction raises \(z_2\) from \(3/4\) to \(2\log 2-1/2=0.88629436\ldots\). The necessary ceiling proved for arbitrary laminate hierarchies is \(\sqrt2-1/2=0.91421356\ldots\) there. For general geometries the two-pole estimate permits \(1\). Its saturation is equivalent to approaching \((0,1)\) at the given fraction. In space the construction gives \((7/12)(3/2-\log 2)=0.47066\ldots\) at half fractions, compared with \(1393/3456\) from two local concentrations. The spatial determinant inequality excludes the opposite corners even for general geometries.

There is an exact planar optimality statement within the construction above. Allow arbitrary step sizes and choices of frame along the surviving branch. Put \(s(x)=e^{4i\theta}\) on a step and \(u(x)=\int_f^x s(y)\,dy\). The stopping-core field is \(u(x)+(1-x)s(x)\), constant on that step, including when its high branches stop at different intermediate fractions. Integrating (2) exactly over each step and then integrating the cross term by parts gives \[
 z=\frac p f\int_0^f\frac{1+x}{(1-x)^3}|u+(1-x)u'|^2\,dx
   =Z_{\rm pl}(f)-\frac p f\left[|u(0)|^2+
                   2\int_0^f\frac{|u|^2}{(1-x)^3}\,dx\right].
\] Alternating \(s=1,-1\) on finer equal steps makes \(u\) tend to zero uniformly. Thus (3) is the supremum for these successive coordinate splits. It also improves the unrestricted finite planar error estimate for the unknown endpoint \(M_e\) at \(K_s\). With \(\delta=n-m\), \(a=1/(k+m)\), \(b=1/m\), \(c=(b-a)/2\), \(g=(a+b)/2\) and \(q=1+(1+p)g\delta+pab\delta^2\), one has \[
 \left|M_e-\left(m+\frac{f\delta(1+g\delta)}q\right)\right|
 \le\frac{fp c^2|\delta|^3}{q^2}\max(1,n/m)\,[1-Z_{\rm pl}(f)].        \tag{4}
\] Indeed the edge spectral comparison applies to the attained geometry. The former factor \(f^2\) falls to \(1-Z_{\rm pl}(f)=f^2/3+f^3/6+O(f^4)\).

\paragraph*{Effective shear at finite phase contrast}

Choose the recursively isotropized staircase above and denote its limiting shear by \(S(x)\). For \(h=(1+(d-1)x)/d\), let \(H\) be the sequential laminate of matrix \((k,m)\) and core \((l,n)\) at fraction \(h\) with equal normal weights on the \(d-1\) axes perpendicular to \(e_d\). Write \(C=C(K_s(x),S(x))\) for the isotropic stiffness with hydrostatic eigenvalue \(dK_s(x)\) and deviatoric eigenvalue \(2S(x)\). Let \(Jv=\operatorname{sym}(v\otimes e_d)\), \(D=H-C\), and \(N_d=d(d+1)/2-1\). Then \[
 (1-x)S'(x)=\frac{d}{2N_d}\operatorname{tr}_{\rm dev}
       [D-DJ(J^T HJ)^{-1}J^TD],\qquad S(0)=m.                         \tag{5}
\] Here the trace is the sum on an orthonormal deviatoric basis. For a small step from \(b\) up to \(x\) each of the \(d\) stopping branches has weight \(w+O(w^2)\). Its stopping tensor and coarse normal approach \(H,e_d\) up to a simultaneous rotation. The exact binary-layer formula for a fraction \(w\) of tensor \(H\) in \(C\) has derivative \(D-DJ(J^THJ)^{-1}J^TD\) at \(w=0\). The derivative of isotropic shear after rotation averaging is its deviatoric trace divided by \(2N_d\). Indeed the arithmetic stiffness and inverse arithmetic compliance of the orientations differ by \(O(w^2)\) when every stiffness is \(C+O(w)\). These variational bounds hold also for an isotropic aggregate with a uniform distribution of rotations, realizable in the closure by a statistically isotropic mixture of rotated grains. On any \(0\le x\le f<1\), positive constituent moduli give uniform ellipticity, all layer formulas have bounded derivatives on the range in use, and \(w=(x-b)/(1-x)+O((x-b)^2)\). Thus the recursion is an Euler approximation of (5) with accumulated error \(O(\mathrm{mesh})\). It retains \(K_s(x)\) exactly at every step. In (5) the dependence on \(S\) is quadratic. The dependence on \(K_s\) drops out of the deviatoric trace since \(D\) acts on a deviator on either side of the quadratic term.

In the plane take the coarse normal \(e_1\) and the one fine normal \(e_2\). Write the normal block of \(H\) as \(\left(\begin{smallmatrix}a&c\\c&d_0\end{smallmatrix}\right)\) and its physical shear as \(g_0\). In terms of averages over core fraction \((1+x)/2\), \[
  1/d_0=\langle1/(k_i+m_i)\rangle,\quad
  c/d_0=\langle(k_i-m_i)/(k_i+m_i)\rangle,\]\[
  a-c^2/d_0=\langle4k_i m_i/(k_i+m_i)\rangle,\quad
  1/g_0=\langle1/m_i\rangle.
\] Equation (5) becomes \[
 (1-x)S'=\frac{ad_0-c^2}{4a}+(1-c/a)S-(1/a+1/g_0)S^2.               \tag{6}
\] The finite numbers below are enclosed by rational polynomials. Write (6) as \(D_0S'=A+B S-C_0 S^2\) with polynomial coefficients, take the degree-36 Taylor polynomial \(P\) at zero, and form \(r=D_0P'-A-BP+C_0P^2\). On \([0,1/2]\) the certificate bounds \(D_0\) below by \(d_*>0\) and \(C_0(S+P)-B\) below by zero using Bernstein coefficients and the minimum constituent shear. The scalar error equation then gives \(|S(1/2)-P(1/2)|\le(2d_*)^{-1}\sum |r_j|2^{-j}\). All bounds and section comparisons use rational arithmetic.
\paragraph*{Exclusion at every spatial contrast}

Normalize \(m=1\) and put \(x=1/n-1>0\). At bulk \(K_\pm=H_{1/2}(l,k;4q/3)\), \(q=n,1\) respectively, write the shear as \(M=H_{1/2}(n,1;y)\). This shifted harmonic mean increases strictly in \(y\). Let \(z_c=119/144\), the half-volume upper bound on the fourth angular energy in the core at an extreme spatial bulk. The arithmetic and reciprocal edge tests give \(y\le y_-\) at \(K_-\) and \(y\ge y_+\) at \(K_+\), where \[
\begin{aligned}
 T(a,b,w)&=\frac{15ab+48aw+56bw}{2(21a+2b+40w)},&
 J(a,b,w)&=\frac{56a+48b+15w}{2[40ab+(2a+21b)w]},\\
 y_-&=T(l,n,n+(1-n)z_c),&
 y_+&=J(1/k,1,1+xz_c).
\end{aligned}
\] The two comparisons use the high-shear and soft-shear cores, respectively. Substitution of (1) yields \[
 y_+-y_-=
 \frac{x(190465552x^3+281583124x^2+56975727x+341280)}
 {2(1+x)(1360x^2+2029x+324)(4114x^2+115509x+19440)}>0 .
\] This proves the spatial assertion for all finite positive contrasts. With the formal fourth statistic one at both edges, the sign would instead be that of \(545168x^2+115598x-5355\), negative from zero to its positive root \(0.03911049\ldots\). Thus the uncapped three-point edge intervals overlap on an open range where the two physical sections are already disjoint. The actual geometry-dependent spatial shear endpoints on (1) remain open.

The whole planar common-shear interval is also below the least spatial shear over \textbf{all} bulks. Lower both phase bulks to \(k\). The spatial Hashin-Shtrikman lower shear bound then gives \(H_f(n,1;s_3(k,n))\) with \(s_3(k,n)=n(9k+8n)/[6(k+2n)]\). The largest planar common shear has shift \(k/(k+2)\), and \[
 s_3(k,n)-\frac{k}{k+2}
 =\frac{125318x^2+39255x+3375}{3(1+x)(242x+45)(272x+45)}>0.
\] This separate-modulus comparison excludes the planar shears even away from the spatial bulk edges. The disjointness of the spatial edge sections themselves uses the coupled estimate above.

The same spatial disjointness holds throughout \(3/7\le f\le4/7\) for the softer phase fraction. Both core fourth statistics are then at most \(5/6\), since the bound for a core of fraction \(g\) is \(7(8+g)/72\). Using \(5/6\) for both statistics, the numerator of \((y_+-y_-)/x\) with positive denominator is \(16339600x^3+23949832x^2+4812507x+20160\). The bulk shifts depend on \(f\); both edge shears are compared through the same strictly increasing \(H_f\).

Throughout \(4997/10000\le f\le5003/10000\) the \textbf{planar} common interval is still exactly \([S_0,S_1]\), with both opposite corners outside the laminate closure, at every shear contrast in (1). To continue the half-area construction to a larger core fraction \(g\), replace each local density \(h\) by \(t+(1-t)h\), \(t=2g-1\), at all levels of the branched periodic cell. Keep the shares of the branches and the terminal normal distributions. All coarse nonconstant angular fields scale by \(1-t\); the terminal coating increments on the cores scale the same way since \(1-h\) acquires that factor. Hence the second field remains zero on the core. Its fourth-field energy integrates the original squared field times \((1-t)^2\) and the larger terminal core density. Upon division by \(g(1-g)\), the resulting fourth statistic is at least \(0.947243(1-g)/g\). For \(g<1/2\) coat the original cell with pure matrix using isotropic normals, giving at least \(0.947243/[2(1-g)]\). At both \(g=f,1-f\) these lower bounds exceed \(35/37\) throughout the stated range. The all-rank laminate ceiling \(2\sqrt{1-g}-(1-g)\) there stays below \(56/61\).
\paragraph*{The field of one ball and the contribution from its surroundings}

Write \(\alpha=\chi-p\), \(V=p(1-p)\). Let \(Y_j(v)\) be the symmetric trace-free part of the unit vector power \(v^{\otimes j}\), normalized to norm one. In particular \(Y_4(v)=\sqrt{35/8}\operatorname{STF}(v^{\otimes4})\). The field \(W=Y_4(D)\alpha\) has squared norm \(V\) per unit cell volume by Parseval. Equivalently \(W=\sqrt{35/8}\operatorname{STF}D^4\Psi\), where \(\Delta^2\Psi=\alpha\). The second field is the corresponding trace-free Hessian of the Newtonian potential \(\Phi=\Delta^{-1}\alpha\), with unit symbol normalization.

Extend a single ball's indicator by \(p\) outside. In the volume coordinate \(t=(\rho/R)^3\), let \(c(t)\) denote its profile. The second field there is \(Y_2(\hat x)[c(t)-t^{-1}\int_0^t c]\). This follows from \((\rho^2\Phi')'=\rho^2(c-p)\) and the radial Hessian formula. The chosen additive constant makes \(\Phi=0\) outside the ball, since its centered mass is zero. Thus the second field of the rest of the packing vanishes inside each ball. As in the disk proof, tending to \(K_s\) forces \(\int c(t)[t^{-1}\int_0^t(1-c)]^2dt\) to tend to zero in the volume average over balls. The same weak-* compactness puts each ball's mass in a central \(p\)-core in the limit, with the total replaced volume tending to zero. The angular energies are stable under this change by the \(L^2\) multiplier bound and the uniform \(L^4\) bound.

Use unit outer radius and put \(s=p^{1/3}\). The biharmonic potential of an uncentered uniform ball of radius \(b\) has the following interior and exterior expressions, with the constant chosen for continuity. \[
 \Psi_b(\rho)=\begin{cases}
 \rho^4/120-b^2\rho^2/12+\mathrm{const},&\rho<b,\\
 -b^3\rho/6-b^5/(30\rho),&\rho>b.
 \end{cases}
\] The radial Laplacian gives the Newtonian potentials \(\rho^2/6-b^2/2\) and \(-b^3/(3\rho)\), respectively. The two radial first derivatives of \(\Psi_b\) agree at \(b\). The self field comes from \(\Psi_s-p\Psi_1\). For a radial power \(\rho^\lambda\) the coefficient of \(\operatorname{STF}(\hat x^{\otimes4})\) in its fourth derivative is \(\lambda(\lambda-2)(\lambda-4)(\lambda-6)\rho^{\lambda-4}\). Thus \(W_s=Y_4(\hat x)A(\rho)\) with \[
 A(\rho)=
 \begin{cases}
  0,& 0<\rho<s,\\
  (5p/\rho^3-7p^{5/3}/\rho^5)/2,&s<\rho<1,\\
  7(p-p^{5/3})/(2\rho^5),&\rho>1 .
 \end{cases}                                                               \tag{2}
\] The exterior fourth energy per unit ball volume and the interior self energy are \[
 E_{\rm out}=\int_1^\infty 3\rho^2 A(\rho)^2d\rho
       =\frac{21}{4}p^2(1-p^{2/3})^2,\qquad E_s=V-E_{\rm out}.
\]

On each ball the remaining field is \(h=\sqrt{35/8}D^4H\) for a harmonic scalar \(H\). Indeed \(\Delta(\Psi-\Psi_s)=\Phi-\Phi_s\) equals a constant \(d\) there, since all other balls have zero exterior Newtonian potential. Take \(H=\Psi-\Psi_s-d|x-x_i|^2/6\) for the ball center \(x_i\). For a countable periodic packing the sum of the compact Newtonian potentials converges in \(L^2\) since their supports are disjoint and they are bounded by a constant times the squared maximum radius. Their distributional Laplacian is the full centered indicator. This gives the harmonic representation in the periodic problem.

Contraction of \(D^4H\) with \(\operatorname{STF}(\hat x^{\otimes4})\) equals \(\partial_\rho^4 H\), since \(H\) is harmonic. Its average over each sphere is zero by differentiating the spherical mean-value identity four times. This proves orthogonality of \(W_s,h\) on each radial shell even though their scalar components need not be angularly orthogonal. Also \(|h|^2\) is subharmonic, so its average on the central core is no larger than its average on the whole ball. The self field vanishes in the core. Summing the orthogonal energies over all balls gives \(\|\chi W\|_2^2\le p(V-E_s)=pE_{\rm out}\). Division by \(V\) proves the upper bound (1).

\paragraph*{A realizing isotropic sequence}

Place an ordinary concentric ball of volume \(a\) in each cell of the unit cubic lattice, keeping the background at mean \(p\) during this calculation. Let \(W_0\) be the fourth field and \(\Omega\) the exterior of the ball in the cell. Its sibling field on the ball satisfies \(\|h\|_\infty\le C_p a^{5/3}\). Indeed (2) is integrable, has zero angular integral, and its exterior translates are bounded by \(C_p R^5/(|v|-R)^5\) at each nonzero lattice point \(v\). The sum converges in three dimensions. The same energy identity as for disks gives \[
 \int_{\rm core}|W_0|^2+p\int_\Omega|W_0|^2
                =a\,pE_{\rm out}-d_a,\quad
 0\le d_a=p\int_{\rm ball}|h|^2-\int_{\rm core}|h|^2\le C_p a^{13/3}.
\] Fill \(\Omega\) at a finer scale by another full packing at the same fraction, with core energy \(J\). The mean-zero Fourier multiplier localization proof for disks applies component by component to \(Y_4(D)\). In addition to the displayed terms and \((1-a)J\), the new core energy contains the cross term \(2(\int_\Omega W_0):\langle\chi_{\rm fine}W_{\rm fine}\rangle\). Its absolute value is at most \(C_p a^{8/3}\), uniformly in the fine packing. In fact \(W_0\) has mean zero, its self integral on the ball vanishes, and the integral of \(h\) on the ball is \(O_p(a^{8/3})\). Cauchy-Schwarz bounds the fine mean by \(\sqrt{pV}\). This estimate works without any isotropic assumption on the fourth phase mean.

Iterate the refill \(N\) times with separated scales. Whole spheres suffice by omitting those intersecting the boundary and covering the vanishing leftover collar with ordinary coated spheres. The fourth energy differs from \(pE_{\rm out}\) by at most \(C_p[a^{5/3}+(1-a)^N]\), summing the per-level errors as a geometric series. The construction can retain cubic symmetry by using a symmetric initial sphere covering and filling the collars in cubic orbits. Aggregate rotations of its cubic elastic tensor at still coarser scales to give an isotropic tensor. All copies have the same phase fraction, so coarse rotation introduces no new angular increments, and averaging preserves both scalar angular energies. Replacing the fine balls meeting a grain boundary affects vanishing volume. At fixed finite positive moduli, the periodic corrector strain under any fixed mean strain has a uniform \(L^{2+\epsilon}\) bound by the higher-integrability estimate for coercive elasticity systems. Testing the energy before and after the replacement in both directions bounds the change by \(C|\mathrm{collar}|^{\epsilon/(2+\epsilon)}\). This yields isotropic sphere assemblages in the closure. Passing first \(N\to\infty\), then \(a\downarrow0\), proves sharpness.

\paragraph*{Shear at the exact matrix bulk}

Set \(b=1/m,\ a=(k+4m/3)^{-1},\ h=(1-p)(4a/3+2b)/5,\ v_0=32(1-p)(a-b)^2/175\). For \(\delta=n-m\), the greatest spherical-class shear if \(\delta>0\), or the least if \(\delta<0\), satisfies \[
 M_e^{\rm sph}=m+p\delta-ph\delta^2+
                 p[h^2+v_0 Z_{\rm sph}(p)]\delta^3+O(\delta^4).
\] The coefficient follows from (1) and the edge spectral law, whose probability has mean \(h\), variance \(v_0z\) at fourth statistic \(z\), and support in \([0,b]\). Uniformity of its Taylor remainder permits optimization over the class.

For \textbf{every} fraction, at either extreme bulk and all contrasts \(1<\max(n/m,m/n)\le141/100\), a physical shear lies strictly beyond the spherical class. The transverse planar stopping layers in the spatial staircase attain fourth statistic \[
 w(p)=(1-p)\left[-\frac76+\frac{35}{24}\log(3/2)
                  +\frac{7(13+2p)}{72p}\log\frac1{1-p}\right]
                           >\frac{141}{100} Z_{\rm sph}(p).                \tag{3}
\] For a rational proof set \(p=t^3\), divide by \(1-p\) and multiply by \((1+t+t^2)^2\). Bound \(\log(3/2)\) from below by \(L=2\sum_{i=0}^2(1/5)^{2i+1}/(2i+1)\), and \(-\log(1-p)/p\) from below by \(\sum_{j=0}^{31}p^j/(j+1)\). The resulting degree-100 polynomial lower bound is \[
 (1+t+t^2)^2\left[-\frac76+\frac{35}{24}L
       +\frac{7(13+2t^3)}{72}\sum_{j=0}^{31}\frac{t^{3j}}{j+1}\right]
                -\frac{2961}{400}t^6(1+t)^2.
\] Every Bernstein coefficient on each of the 32 equal subintervals of \([0,1]\) exceeds \(1/625\), by five rational bisections. The verifier also derives the log expression in (3) from the original stopping-layer integral at symbolic \(p\). At fixed positive moduli the signed excursion from the matrix shear is
\[
 M-M_s=\frac{p\delta^3}{(1+h\delta)^2}\int\frac{(u-h)^2}{1+\delta u}\,d\mu(u).
\] The ratio of the largest and smallest possible kernel weights is \(\max(n/m,m/n)\). Thus (3) strictly separates the two ranges of excursions for either sign of the shear contrast.
\phantomsection\label{endfile:dimension-tools}

\phantomsection\label{file:dimension-reversal}
\paragraph*{Exact common bulk sections at every phase proportion and shear contrast}

Fix \textbf{any} softer-phase fraction \(0<f<1\) and any finite shears \(0<n<m\). There are positive phase bulks for which the set of shears occurring at every admissible bulk is a specified positive-length interval in the plane and is empty in space. Even its two endpoints occur at every planar bulk. At the opposite bulk edge each endpoint requires a geometry outside the closure of two-phase concentric disk assemblages with the prescribed fraction in each disk.

One explicit choice has softer phase \((l,n)\), harder-shear phase \((k,m)\), and \[
 \begin{gathered}
 P=\max(f,1-f),\quad R=3P^2(1-P),\quad r=\frac{3R}{2+R},\quad
 D=1/n-1/m,\quad c=\frac{1-r}{2},\\
                    l=n/3,\qquad 1/k=3/n+cD .                              \end{gathered}\tag{1}
\] In particular \(l>k\). Write \(H_f(x,y;q)=[f/(x+q)+(1-f)/(y+q)]^{-1}-q\), \(s(k,m)=km/(k+2m)\). Define \[
 \begin{aligned}
 K_-&=H_f(l,k;n),& K_+&=H_f(l,k;m),\\
 S_0&=H_f(n,m;s(l,n)),& S_1&=H_f(n,m;s(k,m)).
 \end{aligned}
\] For every \(S_0\le M\le S_1\) the sharp planar bulk section is \([K_-,K_+]\). Every \(M<S_0\) misses \(K_-\), and every \(M>S_1\) misses \(K_+\). The two points \((K_+,S_0)\) and \((K_-,S_1)\) are both inaccessible to the circular class. A laminate construction reaches them. For the same phases the greatest spatial shear at its minimum bulk is strictly less than the least spatial shear at its maximum bulk. In fact the entire planar common-shear interval is below the spatial shear floor at \textbf{every} bulk.

For example, choose a quarter of \((1/6,1/2)\) and three quarters of \((155/967,1)\). Then \[
 [S_0,S_1]=[21/26,3147/3895],\qquad [K_-,K_+]=[4997/30870,1459/9013].
\] In space the maximum shear at minimum bulk is at most \(2467/2972\), and the minimum shear at maximum bulk is at least \(142693/170086\). 

\paragraph*{Planar construction, exclusion and the uniform fraction estimate}

Use core fraction \(g\) for either labeling. For the comparison below set \(R_g=3g^2(1-g)\). A staircase of coordinate splits with a final rotation aggregate has zero second field on its core and attains fourth statistic \[
                        w_g=(1-g)[-2\log(1-g)-g]/g .
\] Thus it attains the same exact bulk for all positive test moduli with distinct bulks and shears.

Write \(O(s)=s/(1-s)\). For every \(0<g<1\), \[
                             O(w_g)>2O(R_g).                              \tag{2}
\] Equivalently \(w_g(1+R_g)>2R_g\). Replacing \(-\log(1-g)\) by its first twenty power-series terms gives a strict lower bound. After division by \(1-g\), the resulting polynomial is \[
 L(g)=\left(2\sum_{j=0}^{19}\frac{g^j}{j+1}-1\right)
                      (1+3g^2-3g^3)-6g^2 .
\] In the degree-22 Bernstein basis on \([0,1]\) its least coefficient is \(110072279/1792502712>0\). This follows directly from \(b_i=\sum_{j=0}^i[g^j]L\,\binom{i}{j}/\binom{22}{j}\). Also \(w_g=1-\sum_{j\ge2}2g^j/[j(j+1)]\) is decreasing and \(R_P\ge R_{1-P}\) because \(P\ge1/2\).

Translate the planar compliance by adding \(\tau=1/m-3D\) to the two reciprocal bulks and subtracting it from both reciprocal shears. This gives positive reciprocal shears \(3D,4D\) and positive reciprocal bulks \(4/m+cD,4/m\) on (1). The stress determinant is a null Lagrangian, so the same translation applies to the effective moduli and preserves the geometry. It preserves \(D\) and \(r=1-2(1/k-1/l)/D\). The translated softer-to-harder shear ratio is \(\eta=3/4\).

For \((K_+,S_0)\) the Gaussian edge comparisons require fourth-statistic odds at least \(O(r)\) and suffice at odds \(O(r)/\eta\). For \((K_-,S_1)\) they require odds at least \(\eta O(r)\) and suffice at \(O(r)\). These follow by evaluating the two bounds at the opposite phase's matrix shear. Now \(O(r)=3O(R)/2\). For both choices \(g=f,1-f\), \[
              O(R_g)\ \le O(R)<\eta O(r),\qquad
              O(r)/\eta=2O(R)<O(w_P)\le O(w_g).
\] Hence the staircases go far enough at both edges, and all the circular assemblages fall strictly short at both opposite corners. The same proof works over the open range \(O(R)/\eta<O(r)<\eta O(w_P)\) when the translated moduli remain positive.

At each fixed bulk edge mix its staircase and its isotropic matrix laminate, which have the same phase fraction and common bulk. Their shear traverses the interval between them. The two edge intervals thus both contain \([S_0,S_1]\). For any shear in that interval, mix the two attained edge materials at their common shear. Their isotropic coated bulk varies continuously between \(K_-\) and \(K_+\), still at the prescribed fraction. These coatings can be realized by isotropic sequential laminates. The separate bulk bounds exclude all bulks outside the indicated interval, and the edge matrix shear bounds supply the exclusions at \(M<S_0\) and \(M>S_1\). The width is positive since \(D-2(1/k-1/l)=rD>0\) and the two finite shears differ.

\paragraph*{Spatial exclusion for the same finite phases}

Normalize \(m=1\), put \(x=1/n-1>0\) and retain \(0<c<1/2\) from (1). The two spatial bulks are \(H_f(l,k;4n/3)\) and \(H_f(l,k;4/3)\). Write an edge shear as \(H_f(n,1;y)\), strictly increasing in \(y\). The spatial arithmetic and reciprocal estimates, using the permissible ceiling one for both fourth statistics, give \[
 y\le T(l,n,1)\ \ \text{at minimum bulk},\qquad
 y\ge J(1/k,1,1/n)\ \ \text{at maximum bulk},
\] where \[
 T(a,b,w)=\frac{15ab+48aw+56bw}{2(21a+2b+40w)},\qquad
 J(a,b,w)=\frac{56a+48b+15w}{2[40ab+(2a+21b)w]} .
\] Substitution of \(1/l=3(1+x)\), \(1/k=3(1+x)+cx\) gives the difference of the lower and upper shifts \[
 \frac{x[(6888+2096c)x^2+(14049+1806c)x+7161-490c]}
 {2(1+x)(40x+49)[(6+2c)x^2+(153+42c)x+147]}>0 .
\] Consequently the two spatial edge intervals are disjoint at every fraction. To compare all spatial bulks with the planar shears, lower the two bulk moduli to \(k\). The spatial separate shear bound gives \(M\ge H_f(n,m;s_3(k,n))\) with \(s_3(k,n)=n(9k+8n)/[6(k+2n)]>2n/3\). The larger planar common shear has shift \(s(k,m)<k/2<n/6\). Strict increase of \(H_f\) proves the claim.

\paragraph*{Reversing the dimensions at the same fraction and the same shears}

For \textbf{every \(0<f<1\) and every \(1<m/n\le9/4\)} a second choice of positive bulks gives a sharp, positive spatial common-shear interval and an empty planar one. Thus throughout this range changing only the bulk moduli can reverse which dimension admits any shear across its entire bulk range.

For example, take a quarter of \((18,1)\) and three quarters of \((1/4,2)\). The spatial shears occurring at every bulk are precisely \[
 \boxed{[215/128,919/547]},\qquad
 \{K:(K,M)\in\mathcal G^{(3)}_{1/4}\}=[516/715,816/779]
 \quad\text{for each such }M.
\] In the plane there is no common shear. For a direct comparison, a quarter of \((1/3,1)\) and three quarters of \((310/967,2)\) have the same shears and fractions. The latter phases are twice the moduli of the first example of this note and have the opposite dimensional separation.

To choose the bulks, write \(s_3(B,\mu)=\mu(9B+8\mu)/[6(B+2\mu)]\). Its range over finite \(B>0\) is \((2\mu/3,3\mu/2)\) and it is strictly increasing. For \(m/n<9/4\) take \(y=(2m/3+3n/2)/2\) and define \[
 B(y,\mu)=\frac{4\mu(3y-2\mu)}{9\mu-6y},\qquad
 l=B(y,n),\qquad k_0=B(y,m).
\] Both matrix shears then equal \(H_f(n,m;y)\). Along the constant-shift level put \(A=1/B,\ b=1/\mu\). Implicit differentiation gives \[
                 \frac{dA}{db}=-\frac{9+16A/b+16(A/b)^2}{10}.
\] The entire level segment stays in positive moduli. Integration shows \(1/k_0-1/l>9D/10\). The planar necessary condition for any common shear is \(1/k-1/l\le D/2\), so the planar exclusion is strict and persists as \(k\) is increased a little.

Write \(B_{\rm V}(k,m,n;g)\) for the attained spatial shear at the matrix bulk from a coordinate-axis laminate with core fraction \(g\), followed by cubic rotation isotropization. Explicitly, with \(\delta=n-m,\ q=1-g\), \[
 B_{\rm V}(k,m,n;g)=\frac25\left[m+\frac{g\delta}{1+2q\delta/[3(k+4m/3)]}\right]
                   +\frac35\left[m+\frac{g\delta}{1+2q\delta/(3m)}\right].
\] At any distinct finite shears it lies strictly above its matrix shear for a stiffer-shear core and strictly below it for a softer-shear core. This follows by strict Jensen applied to the two distinct spectral nodes with positive weights \(2/5,3/5\). The fixed-bulk shear interval on each line thus extends a positive distance past the common matrix shear at \(k_0\).

For an explicit terminating choice set \(k_{\rm P}=(1/l+D/2)^{-1}\), and take \(k=k_0+2^{-j}(k_{\rm P}-k_0)\), \(j=1,2,\ldots\), until both strict comparisons hold: \[
 B_{\rm V}(l,n,m;1-f)>S_1^{(3)}=H_f(n,m;s_3(k,m)),\qquad
 B_{\rm V}(k,m,n;f)<S_0^{(3)}=H_f(n,m;s_3(l,n)).
\] Continuity and the strict margins at \(k_0\) ensure termination. Since \(k>k_0\) one has \(S_1^{(3)}>S_0^{(3)}\). Common-bulk and then common-shear interpolation attains the full spatial bulk interval at every shear in \([S_0^{(3)},S_1^{(3)}]\). The minimum shear at the lower-bulk edge is \(S_0^{(3)}\) and the maximum at the upper-bulk edge is \(S_1^{(3)}\), by edge spectral Jensen with the matrix laminate. These two exclusions make the common interval exact. In the quarter-fraction example the margins for the displayed comparisons are \(919/298662\) and \(7/18560\).

At \(m/n=9/4\) the two ranges of shifts meet in the limit \(k\downarrow0,\ l\uparrow\infty\), with \(2m/3=3n/2\). Each limiting cubic formula still has distinct nodes and finite denominators. The same two strict comparisons therefore hold for sufficiently large finite \(j\) at \(k=m/j,\ l=jn\). The distinct finite shifts satisfy \(s_3(k,m)>s_3(l,n)\) and \(1/k-1/l>D/2\) for large \(j\), proving the endpoint case with finite phases.

In fact \textbf{each fixed fraction extends strictly above \(9/4\)}. For \(C=m/n\ge9/4\), the exact range of contrasts where these two cubic witnesses can bracket the opposite matrix shears for some finite bulks is given by \[
    4(1-f)(C-1)<5,\qquad
          (7f+5)C^2-(22f+5)C-15(1-f)<0.                               \tag{4}
\] Indeed the cubic witness increases strictly with its matrix bulk and is independent of its core bulk. The matrix shear also increases strictly with its matrix bulk. Thus both margins have their strict suprema at the same limits \(l\uparrow\infty,\ k\downarrow0\). Normalize \(n=1\), substitute these limits and factor the two differences. Their signs are respectively those of \(5-4(1-f)(C-1)\) and \(15(1-f)+(22f+5)C-(7f+5)C^2\). For \(C=9/4\) these equal \(5f\) and \(15(1-f)/16\), both positive. Strict positivity at the limits suffices at finite moduli by continuity. If a limiting margin is zero, every finite choice falls short for that witness. At half fractions the method works through all \(9/4\le C<(16+\sqrt{511})/17\). For \(C=113/50\), finite bulks \(l=1000,\ k=1/1000\) with \(n=1\) suffice. As \(f\downarrow0\) the first threshold in (4) tends to \(9/4\), so this construction gives no larger uniform contrast range over all fractions. The maximum contrast permitting this dimensional reversal over all geometries is unresolved. The opposite shear endpoints at fixed extreme bulk and the general interior shear boundaries remain open in both dimensions.
\phantomsection\label{endfile:dimension-reversal}

\subsection{A cubic test at an interior spatial bulk}\label{route:a-cubic-test-at-an-interior-spatial-bulk}
\subsection{The face of the complete quadratic comparison}\label{detail:quadratic-face}\label{file:quadratic-face}
In this derivation the bulk shift is $X=17/28$.
Only $n=1/40$ and $n=1/50$ are required. The following compression
explains why a certificate on six canonical loads tests any finite
number of loads. The invariant wave coefficients are defined by
placing a canonical load in each hydrostatic or deviatoric slot.
\paragraph*{Pullback of an arbitrary violated test}

Let \(S=\mathrm{Sym}(\mathbb R^d)\), \(N=d(d+1)/2\). On strain labels the phase energy is stiffness and on stress labels it is compliance. A positive weight on a loading can be absorbed into its amplitude. Write \(L_i\) for the direct sum of the phase forms and \(D_*\) for the direct sum of the proposed contrast transforms. Each is a single isotropic form on \(S\) times the identity on all labels of a given type. For a unit wave normal \(n\) a strain field lies in \(W_e(n)=\{\mathrm{sym}(v\otimes n)\}\) and a stress field lies in \(W_s(n)=\{\sigma:\sigma n=0\}\). An admissible translation \(T\) obeys \(T|_{W(n)}\succeq0\) for each normal and \(L_i-T\succeq0\) for both phases. The comparison requires \((D_*+T)[d]\ge0\), where \(d\) lists the differences of the phase-mean fields. These are the ordinary translation tests on positive finite phases with distinct bulks and shears.

Suppose an arbitrary label list, form \(T\), and prescribed direction \(d=(d_\alpha)\) violate the comparison. Fix an orthonormal basis \(E_1,\ldots,E_N\) of \(S\). For each type \(p=e,s\) form the matrix of mean coefficients \[
                 U^p_{\alpha j}=d_\alpha:E_j,\qquad \alpha\text{ of type }p.
\] Take one canonical copy of \(S\) as the label space for each type, and write \(G^p_j\in S\) for its fields. Map to the original fields by \((JG)_\alpha=\sum_j U^p_{\alpha j}G^p_j\). The map respects every wave constraint, since it combines only fields of the same type. Thus \(J^TTJ\) is wave-positive. It is bounded above by \(J^TL_iJ\), in which the label metric on type \(p\) is \(B_p=(U^p)^TU^p\). The same metric occurs in \(J^TD_*J\). At the canonical direction \(G^p_j=E_j\), the pulled-back mean inequality has precisely the original negative value.

Identify a list \(G^p_j\) with a linear map \(S\to S\). Rotate this map by conjugation under the natural action of \(O(d)\) on \(S\). The canonical direction is the identity map. Average the pullbacks over normalized Haar measure. The quadratic value at the canonical direction stays fixed. Wave positivity and the two phase inequalities are preserved since every rotated wave space is also allowed. The phase and proposed forms acquire the averaged label metrics \[
 \overline B_p=b_{p,0}P_0+b_{p,2}P_2,\qquad
 b_{p,0}=\sum_{\alpha:p}|P_0d_\alpha|^2,\qquad
 b_{p,2}=\frac{\sum_{\alpha:p}|P_2d_\alpha|^2}{N-1}.
\] Here \(P_0,P_2\) are the spherical and traceless-symmetric projections. Averaging acts by a scalar on each of these inequivalent irreducible spaces under \(O(d)\).

Rescale the canonical fields on each label summand by \(\sqrt{b_{p,j}}\), applying the inverse congruence to the averaged translation. This restores the ordinary unit-weight phase and proposed energies and preserves every wave space. If a scalar \(b_{p,j}\) vanishes, the original \(U^p\) annihilates that whole summand. Its pullback and all its rotated copies have zero rows and columns there, and the translation extends by zero on those labels. The resulting invariant test has a negative value at the isotropic mean with maps \[
              G^p=\sqrt{b_{p,0}}P_0+\sqrt{b_{p,2}}P_2,qquad p=e,s .
\] Conversely, every such canonical test belongs to the allowed finite experiments. This proves equivalence, including violations first found with noninvariant translating forms, nonscalar means or higher-degree label representations.

\paragraph*{Finite size of the remaining optimization}

For each label type the field space is \(S\otimes S\), with \(S=V_0\oplus V_2\). In space, \(V_2\otimes V_2=V_0\oplus V_1\oplus V_2\oplus V_3\oplus V_4\). Combining the two label types gives multiplicities \((4,2,6,2,2)\) in orders zero through four. An invariant symmetric form has one real symmetric matrix per multiplicity space, giving \(10+3+21+3+3=40\) coefficients. In the plane the two-dimensional traceless representation tensored with itself gives one scalar, one pseudoscalar and the two-dimensional order-four representation. The total multiplicities of the scalar, pseudoscalar, order-two and order-four representations are \((4,2,4,2)\). Their symmetric matrix entries total \(10+3+10+3=26\).

At the two stated ratios the following rational comparison exhibits both the admissible face and an excluding quadratic above it.

Write the physical bulk and shear as
$H_f(1,1/10;17/28)$ and $H_f(n,1;Y)$. The upper comparison shift
at $n=1/40$ and $n=1/50$ is
\[
 Y_+(n)=\frac{3(100726n^2-199699n+1157445)}
                   {2(-7039924n^2+3329161n+2457945)} .
\]
\paragraph*{Endpoint certificates}

In the full quadratic comparison a positive wave covariance and two positive phase covariances must satisfy 40 linear moment equations. These equations and their sufficiency for arbitrary loadings follow from the prescribed strain-stress mean theorem. At each endpoint one can use just \textbf{six rank-one terms on the lower branch and seven on the upper}. The archive stores a polynomial vector \(u_j(n)\) of length 36 and a positive rational weight \(c_j(n)\) for each term. Sum \(c_j u_j u_j^T\) in its designated wave or phase block. The lower distribution is two wave terms, one A term and three B terms. The upper distribution is two wave terms, four A terms and one B term. Rotation averaging is built into the 40 invariant tests. The vectors before averaging need only span representatives of the relevant rotation modes.

To check this construction, substitute the three matrices into the 40 equations, clear their common denominator and equate the coefficients of \(n\). Every coefficient vanishes identically. Root counts give positive weights and nonzero denominators all the way through the corresponding endpoint of the contrast interval. Thus all finite quadratic translations accept the displayed shear endpoints. An independent way to find the vectors is to reduce the exposing translation at several rational contrasts. The surviving wave space has an axial component and one transverse pair. A representative from each surviving phase harmonic and wave harmonic supplied the data for rational interpolation. The polynomial verifier uses only the preceding coefficients.

For exclusion, experiment L on the lower branch has an active in-plane transverse constraint. Experiment U on the upper branch has an active polar-shear constraint. The verifier specifies both contact directions as rational functions of \(n\) and checks the four phase and wave inequalities there. The axial slack and the first scalar minor are strictly positive and the scalar determinant vanishes. Hence a null vector has a nonzero shear mean. Going below \(Y_-\) decreases \(2Y/5\) in L, and going above \(Y_+\) decreases \(1/(10Y)\) in U. Each change strictly violates the corresponding translation. The last contrast in the upper certificate is the first positive root of \[
43963823770n^4-1076010696455n^3+97764203077n^2+798074509545n-36312626025.
\] There the second phase constraint also becomes active. In the lower certificate the squared direction coefficient reaches zero at the first positive root of \[
\begin{gathered}321126400n^7+2641264640n^6-859436032n^5-2078259104n^4\\{}-21457812n^3
+133832511n^2+14454477n+381024 .\end{gathered}
\] The weak direction inequalities still hold at the two roots.
\phantomsection\label{endfile:quadratic-face}

\phantomsection\label{file:cubic-exclusion}
\paragraph*{From an endpoint sequence to gradient fields}

Use the six canonical loads \(E=(I,D_1,\ldots,D_5)\), with \(D_j\) the Cartesian deviatoric basis. Put \(C_i e=2m_i e+(k_i-2m_i/3)\operatorname{tr}(e)I\). The target shear shift is \(Y=922012383/1352946172\) or \(2162796207/3152140313\) at fortyfold or fiftyfold contrast, respectively. The canonical difference of phase means is \(d=(E,-C_{17/28,Y}E)\). Consider any sequence of periodic cells of the prescribed fraction with effective tensors tending to the target isotropic tensor. Periodic-cell closure suffices for the usual two-phase G-closure, including separated-scale limits.

Let \(\chi\) indicate B. In the actual six elasticity solutions write \(a=(1-\chi)(e-e_A)\), \(b=\chi(e-e_B)\), where \(e_A,e_B\) are their respective phase means. These two vectors have 36 strain entries each and zero mean. On a nonzero Fourier mode with normal \(N\), the permitted total strain and stress amplitudes have the form \(J_Nw\), a compatible strain and a divergence-free stress for each load. The decomposition \[
      J_Nw=c d+(a',C_Aa')+(b',C_B b'),\qquad c=\widehat\chi ,
\] determines the two phase coordinates from \((w,c)\). At equality three positive semidefinite penalties must vanish, on \(w,a',b'\), respectively. At \(N=e_3\) their joint kernel in \((w,c)\) has dimension three; at every other normal it is its simultaneous spatial and load rotation. For an approaching sequence the sum of squared distances to these kernels tends to zero. On applying the fixed exposing translation, each of the three nonnegative defects tends to zero, and the spectral gap on the orthogonal complement of a fixed finite-dimensional kernel controls distance to it. The discrepancy between actual and target phase differences multiplies \(\widehat\chi\) and tends to zero in summed norm.

We need enough integrability to take determinants on a limit sequence. Construct corrected fields directly from \(\chi\). Let \(H=D^2\Delta^{-1}(\chi-f)\) and \(g=9X\chi\operatorname{dev}H/(4n)\). On the upper kernel the hydrostatic stress wave is \(-9Xc(I-NN^t)/2\). Its phase traces are constant and its phase means spherical for the isotropic target. Consequently the true B hydrostatic strain fluctuation is \(g+o(1)\) in \(L^2\). The map observing \((c,b'_h)\) from the three-dimensional kernel is injective. At each normal reconstruct only in the kernel, prescribing \(c=\widehat\chi\) exactly and fitting \(b'_h\) to \(\widehat g\) by least squares. This is possible since \(c\) is nonzero on an axial basis wave. The homogeneous least-squares problem on the other two real amplitudes is injective and rotation covariant in invariant inner products. Its norm is uniform in direction. This gives real fields \(a',b'\) of mean zero with exact kernel amplitudes, and \[
         \|\chi a'\|_2+\|(1-\chi)b'\|_2\longrightarrow0.                  \tag{1}
\] Indeed the reconstruction is a bounded left inverse on actual kernel observations. Both the observation errors and the distance of the true fields from the kernel tend to zero.

There is a uniform \(L^4\) bound for the corrected fields on \textbf{every} indicator of the prescribed fraction. The normal-dependent matrices of the reconstruction are smooth and even on the sphere. To see this, use a smooth frame near each normal, rotate the fixed kernel, and perform the uniquely specified least-squares fit. The isotropy about the normal and coordinate-invariant fit give identical maps where charts overlap. The degree-zero homogeneous symbols are periodic \(L^4\) Fourier multipliers by the periodic multiplier theorem \cite[Thm.~2.1]{p13:riesz}. The finite differences through order two obey $|\Delta^\alpha M(k)|\le C_\alpha(1+|k|)^{-|\alpha|}$ by homogeneity and smoothness away from the origin. Choose the zero mode separately. The same holds for \(H\), and multiplication by \(\chi\) is an \(L^4\) contraction. All constants are fixed at the two material data and independent of the chosen geometry and fraction. This step uses periodic Fourier fields; there is no layer-tree assumption.

\paragraph*{An identity for the mean of a cubic polynomial}

The phase coordinates of all rotated kernel waves span spaces of dimensions 24 in A and 5 in B. Retain independent Cartesian entries of each as specified in the certificate. Set \[
      y=(\chi,a'_{\rm selected},b'_{\rm selected})\in\mathbb R^{30},
      \qquad \langle y\rangle=(f,0,0).
\] For each of the two contrasts the certificate gives a real polynomial \(P\) of degree three, normalized to have these properties: \[
\begin{gathered}
P(0,a,0)=0\quad (a\in\mathbb R^{24}),\qquad P(1,0,b)=0\quad(b\in\mathbb R^5),\\
P(f,0,0)=f(1-f),\qquad \langle P(y)\rangle=P(\langle y\rangle).
\end{gathered}\tag{2}
\] The first line concerns arbitrary phase coordinates. The mean identity holds for fields whose Fourier amplitudes lie in the equality kernels.

Here is how the certificate establishes the mean identity. It supplies 18 real linear maps \(L_i:\mathbb R^{30}\to\mathbb R^3\) whose values on a kernel wave of normal \(N\) are parallel to \(N\). Each \(L_i y\) is therefore a periodic gradient plus its constant mean. For every triple, \[
   \left\langle\det[L_i y,L_j y,L_k y]\right\rangle
       =\det[L_i\langle y\rangle,L_j\langle y\rangle,L_k\langle y\rangle].
\] For smooth potentials this follows by integrating by parts, or by expanding the determinant and using the divergence-free cross product of two gradients. Smooth approximation in \(L^3\) proves it for these fields. Three homogeneous quadratic forms \(q_j\) also vanish on each full directional kernel. By Fourier orthogonality their averages equal their values at the mean, including for complex Fourier coefficients by polarization. The certificate combines 39 such determinants, the three quadratics and one linear expression. It sums this combination over the 24 rotations of the cube. The rotations act on strain entries and load labels together. The gradient property holds after rotating the spatial output back, which preserves its determinant. Thus every term of the sum has the required mean identity.

One can find the certificate by linear algebra. First solve the tangent equations \(t^t L y=0\) for an unknown row map \(L\), using both tangent vectors of each normal and a basis of the allowed waves. Solve the corresponding quadratic equations in the unknown coefficients of a quadratic form. Next form the determinants of triples of row maps and average each test over cube rotations. Vanishing on the two phase planes is a linear constraint on their weights. Two solutions have different third roots on the line \((f,0,0)\). Taking the linear combination that cancels the cubic coefficient there leaves a nonzero multiple of \(f(1-f)\). The earlier tables give that combination directly.

The verification uses the original six-load rational wave system. It checks every gradient row against both tangent directions and each quadratic against every pair of wave amplitudes. In the quaternion chart \((1,s,t,0)\), clearing the positive normalizing denominator gives polynomials of total degree at most 16. The 153 integer pairs \(s,t\ge0,\ s+t\le16\) prove the identities by interpolation; rotation around the starting normal preserves its full kernel, and continuity supplies the missing pole. On the first pure phase the remaining polynomial has degree at most three in 24 variables, and on the other it has degree at most three in five variables. The integer points of total degree at most three are unisolvent in each case. All 2,925 and 56 evaluations vanish exactly. The check at the mean is a univariate polynomial identity for all fractions. The integer weights, component maps, column ordering and normalization are included in the archive.

\paragraph*{Why degenerating geometry cannot reach the upper shear}

Replace \(y\) pointwise by its phase tuple \(y_w\), setting \(b'\) to zero in A and \(a'\) to zero in B. By (2) its polynomial is identically zero. For the fixed cubic there is a constant \(C_P\) such that \[
 |P(y)-P(y_w)|\le C_P(1+|y|^2+|y_w|^2)|y-y_w|.
\] The uniform \(L^4\) estimate and Cauchy-Schwarz bound its integral by a fixed constant times the left side of (1). But its integral equals \(f(1-f)>0\). This excludes every approaching sequence, so closedness gives a strict gap for the physical maximum at fixed bulk.

The proof also bounds the explicit geometric residual away from zero on every indicator. Choose the upper reconstruction preserving \(c\) as above. For its canonical mean strains and stresses the pointwise constitutive mismatch is \[
       s-C_\chi e=(C_B-C_A)\big[(1-\chi)b'-\chi a'\big].
\] This follows by adding the phase means back and using the exact indicator in the mean difference term. Invertibility of \(C_B-C_A\), norm equivalence in the selected coordinates, and (2) give \[
        \sum_{j=1}^6\|s^j-C_\chi e^j\|_2^2\ \ge\ c_n[f(1-f)]^2,
             \qquad c_n>0 .
\] This is a geometric lower bound, uniform in the fraction; it applies without solving an elasticity problem for the proposed cell. The proof establishes positivity of \(c_n\) without a numerical estimate of its size. For the opposite, lower comparison shear the known all-depth exclusion applies to spatial branching laminates. Its unrestricted endpoint at these contrasts remains open.
\phantomsection\label{endfile:cubic-exclusion}

\subsection{An angular restriction at weak contrast}\label{route:an-angular-restriction-at-weak-contrast}
For a single matrix coating let $s,w$ be the magnitudes of its
normalized spherical harmonics of degrees two and four. In the plane
they are the second and fourth Fourier moments of the normal law.
If $p$ is the matrix fraction, the two core energy statistics are
$x=ps^2$ and $z=pw^2$. The attainable angular pairs in this class obey
\[
\begin{array}{ll}
\mathrm{plane}:&(2s^2-1)_+\le w\le1,\\
\mathrm{space}:&((35s^2-10s-7)/18)_+\le w\le\sqrt{(7+5s^2)/12},
\end{array}\qquad 0\le s\le1 .
\]
Both sides are sharp for the specified sequential grains. Their
isotropic assemblies at the same local fraction preserve the resulting
region in $(x,z)$. The formulas determine the weak-contrast response
through cubic order for this construction. Allowing more general
branching changes the feasible set.
\subsection{Moment inequalities for the two angular harmonics}\label{detail:angular-proof}\label{file:angular-proof}
\paragraph*{Exact angular region of a single matrix coating}

Let \(\nu\) be the probability law of the normals of successively added layers of one pure matrix around a core of the other phase. Weight the normals by their increments in core fraction on the surviving branch. Set \(s=|\int Y_2\,d\nu|\), \(w=|\int Y_4\,d\nu|\), where the angular symbols have unit norm. In the plane these are the magnitudes of the Fourier modes \(e^{2i\theta},e^{4i\theta}\). In space their Gram kernels are the Legendre polynomials \(P_2,P_4\). The exact joint regions are \[
\begin{array}{ll}
  \text{plane:} & (2s^2-1)_+\ \le w\le1,\\[2pt]
  \text{space:} & \displaystyle
   \left(\frac{35s^2-10s-7}{18}\right)_+
       \le w\le\sqrt{\frac{7+5s^2}{12}},
\end{array}\qquad 0\le s\le1.                                      \tag{1}
\] For matrix fraction \(p=1-f\), the two energy statistics in the core are precisely \(x=ps^2,\ z=pw^2\). Hence (1) gives the entire joint weak-contrast bulk-shear region through cubic order for a fixed-matrix sequential grain and isotropic assemblies of such grains at the same local fraction. The upper and lower functions for \(z\) are respectively concave and convex in \(x\), so assembly of different normal laws adds no further angular points. Permitting the other phase to act as matrix adds the phase-reflected diagram at complementary fraction and takes the convex hull. An assembly of textures with the same fraction averages the two angular statistics without a coarse contribution.

In the plane rotate axes so the mean of \(V=e^{2i\theta}\) is \(s\). The variance inequality \(|\mathbb E V^2-s^2|\le1-s^2\) gives the lower bound. For \(s^2\ge1/2\) take equal weights at \(V=s\pm i\sqrt{1-s^2}\). For the upper bound take \(V=\pm1\) with weights \((1\pm s)/2\). These laws have the same second moment. Interpolating them gives every intervening \(w\), including zero when \(s^2\le1/2\).

In space write \(D=\mathbb E(nn^T)-I/3\) and \(H=\mathbb E\,\operatorname{STF}(n^{\otimes4})\), where \(\operatorname{STF}\) denotes the completely symmetric trace-free part. Then \(|D|^2=2s^2/3\), \(|H|^2=8w^2/35\). The scalar, quadratic and fourth-harmonic pieces of \(\mathbb E n^{\otimes4}\) give \[
 0\le\operatorname{Var}(n^TDn)
 =\frac{2}{15}|D|^2+\frac47\operatorname{tr}(D^3)
       +H:(D\otimes D)-|D|^4 .
\] The traceless eigenvalues give \(\operatorname{tr}(D^3)\le2s^3/9\). Also \(|\operatorname{STF}(D\otimes D)|^2=8s^4/35\). To see the latter, this rotation-invariant quartic on traceless symmetric matrices is a multiple of \((\operatorname{tr}D^2)^2\), since the trace is zero. Evaluation on \(s(ee^T-I/3)\) gives the coefficient. The verifier also computes the projection component by component. Cauchy-Schwarz now gives \(18w\ge35s^2-10s-7\) for \(s>0\).

Let \(c=(5+3\sqrt{30})/35\). For \(s\ge c\), put the normals at equally spaced azimuths on the cone \(P_2(\cos\theta)=s\). Five azimuths integrate all moments needed here. The fourth amplitude is \(P_4(\cos\theta)=(35s^2-10s-7)/18\). Below \(c\), mix the cone of amplitude \(c\) with an isotropic normal law to obtain \(w=0\) and arbitrary \(0\le s\le c\). This proves that \((1-f)c^2\) is the sharp zero-fourth endpoint for \textbf{every} fixed-matrix sequential coating in space at fraction \(f\).

For the upper bound draw \(n,n'\) independently. In \(w^2=\mathbb E P_4(n\cdot n')\) use \((n\cdot n')^4\le(n\cdot n')^2\) and \(\mathbb E(n\cdot n')^2=(1+2s^2)/3\). Equality holds for three perpendicular axes with weights equal to the eigenvalues of any prescribed \(\mathbb E nn^T\). Use the second moment of the lower construction and interpolate the two normal laws. The second moment is fixed and continuity fills the entire interval of fourth amplitudes.
\phantomsection\label{endfile:angular-proof}

\subsection{A bound on the spatial extreme-bulk edge}\label{route:a-bound-on-the-spatial-extreme-bulk-edge}
Write $z_2,z_4$ for the normalized angular energy fractions of orders
two and four on the core. At the extreme matrix bulk, including limits
approaching it, $z_2$ tends to zero. Let $Z_3(f)$ be the supremum of
the limiting $z_4$ over geometries with core fraction $f$. The
unrestricted ceiling $Z_3(f)\le7(8+f)/72$ and a weighted C15
construction give
\[
                    0.70015<Z_3(1/2)\le119/144.
\]
For that one texture, compatible and equilibrated test fields on a
mesh of 8,704,000 tetrahedra enclose its continuum statistic in
$[0.7001559314,0.7007313920]$. The second decimal is an upper
bound for the texture. The upper bound for the extremum is still
$119/144$. The formula for the lower estimate corrects the trial
value using primary and adjoint energy errors. The exact check uses the archived geometry recipe with freshly
generated integer primary, adjoint and curl potentials. Recomputing
the outward interval from those stored potentials requires no linear
solve. The regenerated enclosure is wider by $10^{-10}$ at either end
than the original calculation receipt. The resulting lower bound for the supremum,
$Z_3(1/2)>0.70015$, follows in either case.
\subsection{The unrestricted fourth-harmonic ceiling}\label{detail:edge-ceiling-proof}\label{file:edge-ceiling-proof}
\paragraph*{An unrestricted fourth-harmonic bound on the spatial bulk edge}

Let the core occupy fraction \(0<f<1\) and the matrix fraction \(m=1-f\). On the exact matrix-bulk edge of the three-dimensional attainable set, the fourth-harmonic core energy fraction satisfies \[
                         Z_3(f)\le \frac{7(8+f)}{72}.
                                                               \tag{1}
\] This includes arbitrary sequences approaching the edge. For distinct positive phase bulks and shears, approaching the matrix bulk forces the second angular fraction \(z_2\) to tend to zero. Here \(Z_3(f)\) is the supremum of the limiting fourth fraction \(z_4\) subject to that condition. At half volume the ceiling equals $119/144$.

\paragraph*{Directional convexity for arbitrary sequences}

The following estimate works in any dimension \(d\). Write \(\chi\) for the core indicator and \(B=1-\chi\). On a periodic cell of unit normalized volume let \(\psi=m|x|^2/(2d)-\Delta^{-1}(\chi-f)\), with zero-mean periodic inverse. Its Hessian \(H\) has trace \(B\) and mean \(mI/d\). For every unit direction \(e\), \[
      \int ((e^THe)_-)^2\ \le\ \int_{\{\chi=1\}}|e^THe|,\qquad
      t_-=\max(0,-t).
                                                               \tag{2}
\] In particular, \(H=0\) on the core implies that \(\psi\) is convex. Small core error in \(L^2\) forces the negative part of \(H\) to tend to zero in \(L^2\) independently of the cell.

Put \(\alpha=m/d\), \(q=\partial_e\psi-\alpha e\cdot x\). This periodic \(H^1\) function satisfies \(\Delta q=\partial_e B\). Project it in the real Hilbert space \(L^2\) of the torus onto the closed convex set \[
                  {\cal K}_e=\{s:\alpha+\partial_e s\ge0
                                  \text{ in distributions}\}.
\] Write \(r=q-P_{{\cal K}_e}q\). The residual of an orthogonal projection onto a closed convex set is firmly nonexpansive. In fact the two projection inequalities at inputs \(q_1,q_2\) give \(\langle r_1-r_2,q_1-q_2\rangle\ge\|r_1-r_2\|_2^2\). The set is translation invariant, so projection and residual commute with every translation. Apply this inequality to translated inputs and take difference quotients in every coordinate. It gives \(r\in H^1\) and \[
       \|\nabla r\|_2^2\le\int\nabla q\cdot\nabla r
          =\int B\,\partial_e r=-\int_{\{\chi=1\}}\partial_e r.           \tag{3}
\]

For a direction parallel to a lattice vector the projection can also be read on the closed lines in direction \(e\). Parametrize such a line by arclength \(t\) with period \(L\). Let \(V'(t)=\alpha t+q(t)\) and let \(C\) be the greatest convex minorant of \(V\) on the whole real line. The positive quadratic growth of \(V\) ensures the minorant exists. Both functions have the same affine increment under \(t\mapsto t+L\), so \(V-C\) and \(V'-C'\) are periodic. Almost every line has \(q\in H^1\) and \(V\in C^1\). At contact points the supporting slope of \(C\) equals \(V'\). On each interval of \(C<V\) the minorant is affine.

The line projection is \(C'-\alpha t\). Indeed for any admissible competitor with monotone lift \(g(t)=s(t)+\alpha t\), periodic integration by parts gives \[
              \int_0^L(V'-C')(g-C')=-\int_0^L(V-C)\,dg\le0.
\] Here \(dC'\) is supported on contact. This is the metric projection inequality on each line and after integration over the transverse torus. The minorant is measurable, for instance from its chord formula, and the projected periodic parts are in \(L^2\) by comparison to the admissible zero function on each line. Hence \(r=V'-C'\). On contact \(r=0\) and its weak line derivative is zero almost everywhere there. On each open gap \(\partial_e r=V''=e^THe\). On contact \(e^THe=\alpha+\partial_e P_{{\cal K}_e}q\ge0\) almost everywhere. Consequently \[
              (e^THe)_-\le|\partial_e r|\le |e^THe|\quad\text{a.e.}
\] Together with (3) this proves (2). Lattice directions are dense, and both sides of (2) are continuous in \(e\) for each \(L^2\) Hessian. All constants are independent of the period lattice.

For the closure assertion, the directional right side of (2) is at most \(\sqrt f\,\|\chi H\|_2\), uniformly in \(e\). A fixed finite \(\varepsilon\)-net of directions approximates the negative part of the least eigenvalue with pointwise error at most \(2\varepsilon|H|\). The global \(L^2\) bound follows from Parseval. Let the core error tend to zero first, then \(\varepsilon\) tend to zero. All negative eigenvalues vanish in \(L^2\). In space one also obtains a rate. At a point write \(a=(\lambda_{\min}(H))_-\), \(b=|H|\). For \(a>0\), directions having squared sine of the angle to a least eigenvector at most \(a/[2(a+b)]\) satisfy \(e^THe\le-a/2\). The two caps occupy at least \(a/(8b)\) of the sphere. Thus \[
    \int a^2\le\left(32\int_{\{\chi=1\}}|H|\right)^{2/3}
                     \left(\int |H|^2\right)^{1/3}.                      \tag{4}
\] This follows by averaging (2) and applying Hölder to \(a^3/b\), taken as zero where \(b=0\).

The angular conclusions also apply directly to stationary ergodic indicators. Periodize a typical sample on growing cubes. For the zero-mean indicator, first insert a smooth Fourier cutoff supported in an annulus. Both the cutoff and its composition with any of the angular multipliers have integrable Schwartz kernels. Inside the cube, their periodic convolutions and the whole-space convolutions differ by a vanishing normalized \(L^2\) error, by a boundary-layer estimate using the bounded indicator. The error from discarding the remaining frequencies is bounded by a constant times the normalized \(L^2\) error for the indicator itself, since the angular multipliers have bounded symbols. Take the cube limit first, using ergodic averages for the input remainder and its stationary transform. Then enlarge the annulus to all nonzero frequencies. The centered indicator has no zero-frequency atom, so its spectral remainder tends to zero. Changing a vanishing relative volume to restore the prescribed fraction preserves these \(L^2\) limits. This gives convergence of the two phase-weighted angular energies and transfers the cell inequalities.

\paragraph*{A polynomial inequality on positive trace-one Hessians}

In space put \(X=\chi-f\), \(\Phi=\Delta^{-1}X\), \(\Psi=\Delta^{-2}X\), and \(T=D^4\Psi\). Denote the trace-free second and fourth harmonics by \(Q,U\), in the unnormalized tensor convention. Thus \(H=BI/3-Q\), \(\int|Q|^2=2fm/3\) and \(\int|U|^2=8fm/35\). On the matrix, \[
          T=T_0+U,\qquad T_0=\tfrac67\operatorname{sym}(I\otimes Q)
                       -\tfrac f5\operatorname{sym}(I\otimes I),\qquad
          Q=I/3-H.
\] Symmetrization averages over all index permutations. Let \(e_2,e_3\) denote the second and third elementary symmetric functions of \(H\) and \(T_{kl}\) the matrix with entries \(T_{ijkl}\). For every \(H\succeq0\) of trace one and every symmetric trace-free fourth tensor \(U\), \[
 |U|^2+2\sum_{kl} H:\operatorname{cof}(T_{kl})
 \ \ge\ a_0+a_2e_2+\tfrac35e_3,\qquad
 a_0=\frac{14f^2-28f+23}{315},\quad a_2=\frac{28f-30}{105}.              \tag{5}
\]

Diagonalize \(H=\operatorname{diag}(h_1,h_2,h_3)\). The homogeneous quadratic form on \(U\) in (5) is nonnegative. At \(H=\operatorname{diag}(0,0,1)\) it is the sum over \(k,l\) of \[
     (U_{11kl}+U_{22kl})^2+U_{33kl}^2+2U_{13kl}^2+2U_{23kl}^2 .
\] Linearity in trace-one \(H\) and permutation of the axes prove this on the entire eigenvalue simplex. The full expression is invariant under each coordinate reflection. The reflection-invariant block contains all the linear terms and its orthogonal complement contributes a nonnegative quadratic form.

The even block has the three entries \(x_{ij}=U_{iijj}\) for \(i<j\), with \(U_{iiii}=-\sum_{j\ne i}x_{ij}\). On it the difference of the sides of (5) equals \(x^TAx+2b^Tx+c_0\). The diagonal \(ij,ij\) entry of \(A\) is \(8(h_i+h_j)\). Its off-diagonal entry for two pairs sharing the index \(k\) is \(2h_k\). The remaining entries are \[
               b_{ij}=-(h_i-h_j)^2/7,\qquad
               c_0=(3-10e_2+9e_3)/735 .
\] Direct expansion from \(T_{iiii}=U_{iiii}+6(1/3-h_i)/7-f/5\) and \(T_{iijj}=x_{ij}+(2/3-h_i-h_j)/7-f/15\) verifies this identity.

In the interior, \(A\) is positive definite by strict diagonal dominance. Order the eigenvalues and put \(h=(r,r+d,r+d+t)/S\), \(S=3r+2d+t\) for \(r,d,t\ge0\). The determinant of the bordered matrix with upper block \(A\), last column \((b,c_0)\) and last row \((b^T,c_0)\) is \(4P/(3S^6)\), where \[
\begin{aligned}
P={}&28r^4(d^2+dt+t^2)
 +r^3(72d^3+108d^2t+116dt^2+40t^3)\\
 &+r^2(60d^4+120d^3t+144d^2t^2+84dt^3+12t^4)\\
 &+r(16d^5+40d^4t+56d^3t^2+44d^2t^3+12dt^4)
   +3d^2t^2(d+t)^2 .
\end{aligned}
\] The Schur complement is nonnegative. Continuity extends (5) to the faces and vertices.

\paragraph*{Integration and equality conditions}

The elementary minors of a constant plus a periodic Hessian integrate to the minors of its mean, by integration by parts. The same holds after polarization. Each \(T_{kl}\) is a periodic Hessian of mean zero. Therefore \(\int\sum_{kl}H:\operatorname{cof}(T_{kl})=0\). At an exact edge with \(H=0\) on the core, the matrix integrals of \(e_2,e_3\) are \(m^2/3,m^3/27\). Integrating (5) gives \[
        \int_{\{\chi=0\}} |U|^2
             \ \ge\ ma_0+\frac{m^2a_2}{3}+\frac{m^3}{45}
             =\frac{fm(16-7f)}{315}.
\] Divide by the total energy \(8fm/35\) to obtain (1).

For an arbitrary sequence on which the second angular fraction tends to zero, \(\|\chi H\|_2=\|\chi Q\|_2\to0\). By (2) the distance of the matrix Hessians to the positive trace-one simplex tends to zero in \(L^2\). One can project their eigenvalues onto this simplex, clipping the negative part and rescaling the positive part to sum to one. The distance is at most a fixed multiple of the norm of the negative part. Apply (5) to these projected Hessians, using the same \(U\) and recomputing \(T_0\). All integrated errors tend to zero. Indeed each differs by at most the \(L^2\) distance times a uniform \(L^2\) bound for a polynomial of degree at most two in \(H,U\) and the projected Hessian. The core errors in the minor identities tend to zero by the same estimate. The required uniform fourth moments for the Fourier fields were proved for arbitrary geometries. This also justifies the integrations by smooth approximation.

More quantitatively \(\|\chi H\|_2=(2fm z_2/3)^{1/2}\), so (4) bounds the \(L^2\) projection error by \(O_f(z_2^{1/6})\). The same calculation gives \(z_4\le7(8+f)/72+C_f z_2^{1/6}\) also off the edge, with a finite constant independent of the geometry.

For equality in (1), conditional on the matrix, \(H\) must be isotropic, planar with spectrum \((0,1/2,1/2)\), or uniaxial with spectrum \((0,0,1)\). These are the only zeros of \(P\) on the ordered simplex. The two minor averages force their conditional probabilities to be \(m^2,\ 4mf/3,\ f(2+f)/3\), respectively. These conditions also hold for an extremizing limiting law.

Write \(Y_4(n)\) for the trace-free part of \(n^{\otimes4}\). At the isotropic state necessarily \(U=0\). At the planar state with null axis \(n\), necessarily \(U=-Y_4(n)/12\). At the uniaxial state \(nn^T\), \[
                        U=-Y_4(n)/8+W_\perp ,
\] where \(W_\perp\) is a symmetric trace-free fourth tensor supported in the perpendicular plane. In fact for both nonisotropic states the displayed axial terms solve \(Ax+b=0\) in their even blocks. At a uniaxial state the square decomposition above shows that the whole quadratic kernel consists of precisely these two real transverse harmonics. At rank two the kernels for two different axes intersect only in zero. The axial and transverse contributions are orthogonal. Subtracting the energies in the two required axial terms from the total matrix energy at equality, using \(|Y_4(n)|^2=8/35\), forces \[
       \mathbb E\bigl[|W_\perp|^2\mid H\text{ uniaxial on the matrix}\bigr]
                         =\frac{50-23f}{360(2+f)} .
\] At half volume this conditional energy would be \(77/1800\). Whether a geometry can simultaneously realize these field values, their proportions and the differential constraints remains open.

A single periodic cell with finitely many smooth phase interfaces cannot attain (1) at the exact bulk edge. Indeed, on each open matrix component the entries of \(H\) are harmonic since \(\Delta\psi=1\). Equality confines their squared Frobenius norm to the three distinct values \(1/3,1/2,1\), so continuity makes this norm constant on each component. Its Laplacian is \(2|\nabla H|^2\), giving constant \(H\) there. The gradient of \(\psi\) is continuous across a smooth core boundary and is constant on the adjacent core component. The affine matrix gradient therefore has zero tangential derivative along the common boundary, forcing \(H\) to have rank at most one. Every matrix component has such a boundary patch when \(0<f<1\) and the interfaces are finitely many disjoint smooth surfaces. This contradicts the positive isotropic share \(m^2\) required for equality. The argument leaves open the possibility of saturation in sequences with progressively finer structure.
\phantomsection\label{endfile:edge-ceiling-proof}

\subsection{Rational C15 geometry and the continuum error estimate}\label{detail:edge-mesh-proof}\label{file:edge-mesh-proof}
\paragraph*{Geometry and angular energy}

In eighths of the conventional cubic period start from \((0,0,0),(2,2,2),(5,5,5),(5,7,7),(7,5,7),(7,7,5)\), add the four face-centering translations, and subtract \((1,1,1)\) to put the inversion center at zero. The 24 Voronoi cells comprise 16 with 20 vertices and eight with 28. Rational plane intersections verify every neighbor inequality, squared circumradius at most \(97/24\) in these units, and interior circumcenters at each flag. The 3,264 flags have successive orthogonal edge vectors. The verifier checks all coordinate permutations and signed coordinate reflections on the sites. A single reflected axis, or its complementary pair, requires a glide by a quarter period in the other two coordinates. The space group leaves 17 inequivalent flags. Their determinant weights sum to the conventional cubic volume.

Let \(u_1\ge u_2\ge u_3\) be coordinates along the three successive edges in a flag, with \(0\le u_i\le1\). At core fraction \(f\), the rank \(s\) is the number of coordinates exceeding \(f\). Write \(P_s\) for the projection onto those \(s\) axes, \(p=1-f\), and \(\theta=1-s/3\). The rank volumes are \(\binom3s p^s f^{3-s}\). The convex Voronoi potential has Hessian \(P_s/3\). Consequently \[
 H=(pI-P_s)/3,\qquad
 \widehat H_{ij}(k)=n_i n_j\widehat\theta(k)\quad(k\ne0),\qquad n=k/|k|.
\] The mean of \(H\) is zero. In rank one use simple plies along the first flag axis. In rank two extrude an independently rotated planar bulk-edge texture in the plane perpendicular to the last axis. It has core fraction \(1/3\) and certified planar statistic \(z_{\rm pl}>0.973\). Rank zero is pure core and rank three pure matrix. The second harmonic cancels on every final core, so the construction has the prescribed bulk for every choice of distinct positive phase bulks and shears. Isotropizing the grains preserves the angular energies.

Let \(R(n)=\operatorname{STF}(n^{\otimes4})\), where STF denotes orthogonal projection onto symmetric trace-free tensors, and use the ordinary Frobenius inner product. Then \(|R(n)|^2=8/35\). Let \(A\) have multiplier \(R(n)\) on nonzero modes of \(\theta\) and zero mean. The local fine-field mean \(B\) is \(R(n_1)/3\) in rank one, \(R(n_3)/4\) in rank two and zero in the pure ranks. The actual fourth statistic is \[
 z_4=\frac{35}{8fp}\int\theta|A+B|^2+\frac{35p}{96}z_{\rm pl}.       \tag{1}
\] Independently rotating the planar texture gives its axial mean and its spin-four variance in (1). All volume integrals here use the normalized periodic cell.

\paragraph*{A continuum lower bound from gradients and curls}

Use the orthogonal projection \(P\) onto periodic mean-zero gradients, acting on the index \(j\) of a vector with three symmetric labels \(abc\). Define \[
 G_{j;abc}=\tfrac13(\delta_{ja}H_{bc}+\delta_{jb}H_{ac}+\delta_{jc}H_{ab}),\qquad T=PG .
\] In Fourier space \(T_{jabc}\) has multiplier \(n_jn_an_bn_c\) on nonzero modes. Thus \(A=\operatorname{STF}T\) and \[
                      \|T\|^2=\operatorname{Var}\theta=fp/3.
\] The last equality follows from the binomial rank volumes. This known norm makes it possible to bound an error in the gradient field without estimating the regularity of its jumps.

Choose any periodic gradients \(U=\nabla u\), \(D=\nabla v\), any componentwise solenoidal \(C=\operatorname{curl}c+C_0\) and any symmetric test tensor \(V\). Put \(W=\operatorname{STF}(\theta V)\). Orthogonality gives the exact identity \[
 \langle T,W\rangle=\langle U,W\rangle+\langle G-U,D\rangle
                           +\langle T-U,W-D-C\rangle .
\] Compute \(E=fp/3-2\langle G,U\rangle+\|U\|^2\) and \(F=\|W-D-C\|^2\). They bound the last term in absolute value by \(\sqrt{EF}\). Completing the square with \(V\) in the first integral of (1) now gives \[
 \int\theta|A+B|^2\ \ge 2\langle U,W\rangle+2\langle G-U,D\rangle
  +\int\theta(2V:B-|V|^2)-2\sqrt{EF}.                         \tag{2}
\] This inequality holds for arbitrary conforming potentials, with no stationarity or numerical solver tolerance assumption.

The implementation subdivides each reference flag by the Cartesian integer grid in its ordered coordinates and triangulates using successively increasing coordinates in each permutation. The affine tetrahedra match on faces since they induce the same ordered refinement on the common flag vertices. An optional piecewise-affine change in each coordinate maps the reference rank knot to \(f\). It preserves the triangulation on common faces and makes the rank interfaces element faces. On each element the gradient of a nodal potential is the sum of its three successive nodal differences times the reciprocal edge vectors. Nodes are identified under the cubic space group and face-centering translations with the signs for polar and axial tensor components. A node fixed by symmetry with opposite signs has value zero. These extensions give periodic piecewise-linear potentials, including across fundamental flags.

The primary solve proposes \(U\) close to \(PG\). The test \(V\) is the componentwise dyadic rounding of \(\operatorname{STF}\operatorname{sym}U+B\). The two adjoint solves propose the gradient and curl parts of \(W\). The constant \(C_0\) is a dyadic cubically invariant tensor. It is allowed as a solenoidal field.

\paragraph*{Checking the arithmetic}

All flags, volumes, projections \(P_s\), reciprocal edge vectors, \(G\) and \(B\) are rational. The nodal potential denominator is \(2^{30}\) and the test denominator is \(2^{25}\). The sparse matrix and solves only propose these integers. The final evaluation of (2) uses integer arithmetic and rational volume sums. Set \(S=2^{36}\). Each gradient component is the sum of three rational directional derivatives. Flooring each to the \(1/S\) grid produces less than three grid units of component error. The integer formula for \(35\operatorname{STF}\) uses the six single traces with coefficient \(-5\) and the three double traces with coefficient one. It gives \(W\) rounded down to one unit. The exact components of \(G,B\) cost at most one unit when floored.

For the primary error use \(E=fp/3+\int(|U-G|^2-|G|^2)\), with four units added to each absolute component of the first difference and an exact rational subtraction for the second term. In the adjoint defect the curl difference costs at most six units, the gradient three and \(W\) one. The code adds ten units to each absolute residual component before squaring. For the unsigned approximation to the other terms of (2), the checked stored components of \(U,D,V,W,G,B\) have absolute value at most two in physical units. The component errors in \(U:W\) and \((G-U):D\) are bounded by \(9/S\) and \(21/S\) per product. There are 81 entries counting multiplicity. Including \(V:B\), subtracting \(5600/S\) bounds the total linear-term error. Each dot product is split into base-\(2^{18}\) digits with an integer overflow check before summation. An integer square root rounds \(\sqrt{EF}\) upwards.
\paragraph*{Power cells and rigorous upper and lower energies}

Use squared power distance \(|x-c_j|^2-2\omega_j\) in a cube of side one. Give the eight diamond sites weight \(\omega=1/125\), the sixteen other sites weight zero. It changes the dual power centres while preserving the Delaunay incidences and all cubic symmetries.

In a nested Delaunay simplex flag denote its ordinary circumcentres by \(c_i\) and power centres by \(b_i\) for \(0\le i\le3\), with \(b_0=c_0\) at the site. Both differences \(c_i-c_{i-1}\) and \(b_i-b_{i-1}\) lie in the same one-dimensional orthogonal complement between consecutive affine spans. Write \(b_i-b_{i-1}=s_i(c_i-c_{i-1})\). For the new incident site at vector \(v\) from the base, solve \[
 (b_i-b_0)\cdot v=|v|^2/2+\omega_0-\omega_v .
\] The verifier checks every incident equality, strict inequality for every other candidate and positivity of every \(s_i\). At the submitted weight positivity is tested for every incident flag. The above tests over all flags keep each dual centre in the relative interior of its primal face: in any facet flag the signed distance into the next face is a positive multiple of the original distance. The covering-radius check excludes more distant sites.

In a rank-\(r\) region take \((1-f)\) times the primal flag simplex on vertices \(c_0,\ldots,c_r\) plus \(f\) times the orthogonal dual flag simplex on \(b_r,\ldots,b_3\). This is the power-diagram form of the convex-potential construction. It gives the same local projector and the same phase proportion as before. An ordered reference flag with coordinate knot \(t_0\) maps affinely in each piece. Coordinate stretches are \((1-f)/(1-t_0)\) above the knot and \(fs_i/t_0\) below. On the common rank boundary the two prescriptions agree because the rank-\(r\) corner is \((1-f)c_r+fb_r\). Faces shared by neighboring flags use the same circumcentres and power centres. The original reference nodes therefore identify a conforming periodic mesh. Exact volume sums give \(\binom3r(1-f)^r f^{3-r}\) at each rank.

For the two-sided calculation use the projection \(T=PG\), the exact coarse fourth field \(A=\operatorname{STF}T\), the compatible approximation \(U\), the dyadic test \(V\), the local mean fine field \(B\), and the adjoint defect of the original continuum certificate. Write \(E=\|T-U\|^2\). Completing the square shows that the terms omitted by its lower bound consist of a bilinear remainder of absolute value at most \(2\sqrt{EF_{\rm adj}}\) and the nonnegative squared residual \(\int\theta|A+B-V|^2\). The latter is at most \[
                       \big(\sqrt{4E/7}+10^{-6}\big)^2 .
\] Indeed the verifier checks \(|\operatorname{STF}\operatorname{sym}U+B-V|<10^{-6}\) in every element. At each Fourier frequency \(T-U\) is a gradient in its distinguished index. On such tensors the squared norm of the symmetric trace-free projection is at most \(4/7\) times the squared norm of the input. Its adjoint contracts a trace-free fourth tensor with the wave normal. Taking that normal as \(e_3\), the squared singular values are \((16-j^2)/28\), \(0\le j\le4\), with multiplicities \(1,2,2,2,2\). The verifier forms the two Gram matrices on the nine independent harmonic components and checks their generalized eigenvalues in rational arithmetic.

For a prism with transverse axes \(u,v\), let \(F=(u+iv)^{\otimes4}\). Start with a rotated planar texture of core fraction \(g\ge h=1/3\), angular core statistic \(z(g)\), and zero conditional mean planar fourth field. Coat it to fraction \(h\) with two perpendicular normals in this plane. On the final core all original core has survived. Its planar fourth field gains a constant of magnitude \(g-h\) and chosen orientation; the variance inherited from the seed is \((1-g)z(g)\). Thus for phase angle \(\phi\) the spatial mean \(B\) gains \((g-h)\operatorname{Re}(e^{i\phi}F)/8\), and the fine core variance in tensor units becomes \((1-g)z(g)/8\). This identity also accounts for the cross term with the coarse field in each element. The bare choice has \(g=h\) and \(z(h)=0.973\). The submitted witness also uses the pure core coating $g=1$, finite splits and stopped fillings. The next section constructs their means and variances.

The mesh selects one candidate in each prism element by its proposed core energy. Each offered spin tensor has exact rational entries and squared norm \(1/8\) including the factor \(1/8\) in its definition. The verifier re-evaluates the chosen geometry, the weighted linear pairing and the sum of the conditional variances exactly, with outward rounding of fixed-point gradients. The selection requires no optimality claim. Finer scales in each element realize the chosen planar coating and preserve the exact matrix-bulk constraint.

For clarity, let $M_0$ be the sum of the two linear pairings and the
local test energy before the rounding allowance, and let $E_{\rm fine}$
be the sum of conditional variances weighted by local core volume. With
$S=2^{36}$, put $M_-=M_0-5600/S$. The certified endpoints before
normalization are
\[
 M_-+E_{\rm fine}-2\sqrt{EF_{\rm adj}},\qquad
 M_-+E_{\rm fine}+11200/S+2\sqrt{EF_{\rm adj}}
                      +(\sqrt{4E/7}+10^{-6})^2 .
\]
Each square root is rounded upward in a term increasing the error
allowance. Multiplying by $35/[8f(1-f)]$ gives the interval printed
above. Integer fields, filling indices and the mesh orbit list together
specify the released witness.
\phantomsection\label{endfile:edge-mesh-proof}

\subsection{Realization of the fine C15 prism fields}\label{detail:mesh-fillings}\label{file:mesh-fillings}
The following calculations close the material part of the mesh
certificate. For rank two the core share in a planar cross-section
is $h=1/3$. Only menu entries actually selected by the saved
index arrays enter the bound. We first describe the two phase
planar fields and the conversion of their lower bounds to exact
attainable values in the closure.
\paragraph*{Geometry and two-scale fields}

The core has fraction \(f\) and the matrix has fraction \(p=1-f\). On a unit-pitch triangular lattice center a regular hexagon of apothem \(f/2\) at every vertex. Put pure core in the hexagons, pure matrix in the small triangles, and fine equal-width plies of both phases in the connecting rectangles. The ply normal points along the segment joining neighboring hexagon centers. The three areas are \(f^2,p^2,2fp\). The traceless second-harmonic field vanishes on all the core. The divergence-free hydrostatic stress certificate for this tiling proves the exact matrix bulk for every positive phase contrast. Its sixfold rotational symmetry makes the effective elasticity isotropic.

Use \(T=\partial_z/\partial_{\bar z}\) at nonzero Fourier frequencies and zero on the constant. Let \(\theta\) equal \(1,0,1/2\) on the hexagons, triangles and rectangles. In a wall let \(q=e^{-2i\alpha}\) where \(\alpha\) is the angle of the normal. The wall directions are multiples of \(\pi/3\), so \(q^2=\bar q\). The exact coarse field \(U_c=T\theta\) is \(-q/2\) in a wall and zero elsewhere. To see this directly, take an affine-periodic complex potential whose full gradient maps each hexagon by reflection and translation, each triangle to a point, and interpolates across each rectangle with the rank-one average of the neighboring gradients. Equivalently, if \(\gamma\) is a hexagon center, set \(F(z)=\overline{z-\gamma}+f\bar\gamma\) there. In a wall interpolate between its facing hexagon traces, which agree after projecting away the center-line direction. At the tip in the triangle of lattice vertices \(v_1,v_2,v_3\) the value from each adjacent hexagon is \(f\overline{(v_1+v_2+v_3)/3}\). These formulas agree across every coarse edge, and \(F_{\bar z}=\theta\), \(F_z=U_c\) with mean zero for \(F_z\). On the fine plies the core increment \(q/2\) cancels \(U_c\).

Write \(W_c=T U_c\). In a core ply the full fourth field is \(W_c+q^2/2\), and in a matrix ply it is \(W_c-q^2/2\). Pure regions have just \(W_c\). Separated Fourier scales give these increments, with the coarse cell and its material fractions kept fixed. Put \(w=1-\theta\), set \(d=q^2/2\) on the walls and zero elsewhere, and normalize integrals to unit cell area. The core's fraction \(z\) of the fourth-field energy is exactly \[
           z=1-\frac{D}{fp},\qquad D=\int w|W_c-d|^2 .                 \tag{4}
\] Parseval gives total energy \(fp\) for the fine two-phase indicator. Thus (4) requires an upper bound on a relatively small integral. It avoids estimating the high-frequency tail of a voxelized interface.

\paragraph*{Two conforming trial potentials}

Use any continuous periodic piecewise-affine complex fields \(G,J\). Define \(P=G_z\), \(E=G_{\bar z}-U_c\), \(D_0=\int w|P-d|^2\) and \(R=w(P-d)\). Let \(e=\int|E|^2\), \(a=\int|R-\overline{J_{\bar z}}|^2\) and \(t=\operatorname{Re}\int E\,J_z\). Then \[
         D_0-2t-2\sqrt{ea}\ \le\ D\ \le\ D_0-2t+2\sqrt{ea}+e .          \tag{5}
\] In fact \(W_c=P-TE\) and \(T^*\overline{J_{\bar z}}=\overline{J_z}\). Expand (4), replace \(\operatorname{Re}\langle TE,R\rangle\) by \(t\) with error at most \(\sqrt{ea}\), and bound \(0\le\int w|TE|^2\le e\). With \(G\) alone the triangle inequality also gives \(\sqrt D\le\sqrt{D_0}+\sqrt e\). These estimates hold for each trial without any convergence assumption about a differential-equation solver.

To produce trials, cut each polygon into triangles from its center, then split each coarse edge into \(N\) equal parts. The verifier identifies the periodic nodes and checks the two pure areas and matching traces on every refined side. Numerical least squares proposes nodal coordinates for \(G\) and \(J\) written as \(A+i\sqrt3\,B\). Round \(A,B\) to multiples of \(10^{-7}\). In lattice coordinates \((x,y/\sqrt3)\) all vertices and gradients used in (5) are rational apart from the explicit \(i\sqrt3\) factor. The verifier sums squared integer derivatives times exact triangle areas and bounds square roots upward by integer square roots. For \(f=1/2\), \(N=48\) gives \(z>9371/10000\). In comparison all hierarchical laminates at the same bulk have \(z\le\sqrt2-1/2<0.914214\).

Mix with the isotropic matrix laminate at the same fraction and exact bulk to attain \(z=9371/10000\). Use the two conditional Gaussian bounds for the spectral measure of the effective shear. For (1) the minimum-bulk construction has shear at least \(6221322/5870693\), and the maximum-bulk construction has shear at most \(117077433/110479064\). The isotropic matrix laminates supply the other ends of both edge intervals. For all hierarchical laminates at the two edges, substituting \(0.914214\) in the opposite Gaussian bounds instead gives \(M<1.059724950\) at \(83/81\) and \(M>1.059725468\) at \(705/688\). Each inequality strictly separates its endpoint from the whole shear strip in (1). The same substitutions with core and matrix data \((2,2),(1,1)\) give (3).

There is also an open weak-contrast family of this distinction. For half \((k,m),(k+\epsilon h,m+\epsilon s)\) with \(s>0\), put \(R=2m^2h/(k^2s)\). Whenever \(2\sqrt2-2<R<0.8742\), all sufficiently small \(\epsilon>0\) admit a positive interval of shears with the entire exact bulk range, while hierarchical laminates admit no such shear. Indeed the two core choices attain fourth statistic \(0.9371\) on their respective matrix-bulk lines. Their cubic-coordinate intervals overlap when \(R<2(0.9371)-1\). The two laminate ceilings prevent overlap when \(R>2(\sqrt2-1/2)-1\). Strictness absorbs the uniform fourth-order remainder.
\paragraph*{Optimal local cancellation for a prescribed coarse field}

Use the second- and fourth-harmonic Fourier symbols \(e^{-2i\alpha},e^{-4i\alpha}\) at nonzero frequencies, with zero constant modes. Write \(U,W\) for the corresponding fields of an indicator \(\chi\). Conjugating both symbols gives the same energy conventions as the spectral bounds. Consider a fine texture of core fraction \(t\), put \(p=1-t\), and add a constant coarse fourth field \(C\) of magnitude \(r\). With integrals normalized to unit area, every texture satisfies \[
  \int(1-\chi)|C+W|^2\ \ge\ p(r-t)_+^2.                              \tag{1}
\] Indeed \(\int W=0\) and \(|\int \chi W|=|\langle\chi-t,W\rangle|\le tp\) by Parseval. Jensen's inequality in the matrix proves (1).

For \(r\ge t\) a rank-two sequential laminate attains equality, with the additional constraint \(U=0\) on its core. Choose perpendicular normals with equal normal weights and fourth symbol \(s=C/r\). Their second symbols cancel on the core. Their fourth symbols agree, so every fine increment of \(W\) equals \(s\) times the change in local core fraction. This gives \(W=ps\) in the core and \(W=-ts\) in every matrix ply. In particular both the core second field and the total matrix fourth field cancel when \(r=t\).

For \(r<t\), begin a planar staircase at fraction \(t\), taking a first step down to \(b\in[0,t]\) with axes of fourth symbol \(s\). Continue the surviving fraction-\(b\) branch with a finely rotated staircase. The survivor has volume \(p^2/(1-b)^2\), fourth increment \(-(t-b)s\), and zero second increment. Its finer staircase has zero conditional mean fourth field in its matrix and core statistic \[
             Z_{\rm pl}(b)=(1-b)[-b-2\log(1-b)]/b,\qquad Z_{\rm pl}(0)=1 .
\] The stopped matrix plies all have \(W=-ts\). Their mass is \(m_b=p(t-b)/(1-b)\). Summing the two matrix contributions gives the attained energy \[
 \begin{aligned}
 D_b(C,s)&=m_b|C-ts|^2
    +\frac{p^2}{1-b}\big\{|C-(t-b)s|^2+b[1-Z_{\rm pl}(b)]\big\}\\
 &=p|C|^2-2m_b\operatorname{Re}(\bar s C)+c_b,\\
 c_b&=m_bt^2+\frac{p^2}{1-b}\big\{(t-b)^2+b[1-Z_{\rm pl}(b)]\big\}.
 \end{aligned}                                                        \tag{2}
\] Here every core point still has zero second field. In the second line \(s\) may be any point of the unit disk, by mixing two orientations at identical local fraction. The constant \(c_b\) is unchanged in that mixture.

In fact \(b=(t-r)_+\) and \(s=C/r\) minimize the shifted matrix energy among all these planar staircase paths (any orientation at \(r=0\)). To see this, parameterize a path by \(u(x)=\int_t^x s(y)\,dy\), with \(|s|=1\) on each step. The stopping core sees \(u+(1-x)u'\) and has mass density \(F'(x)=p^2(1+x)/(1-x)^3\). Integrating this field against \(F'\) by parts gives its integral over the core as \(-p^2[u(0)+2\int_0^t u(x)/(1-x)^3\,dx]\). The matrix integral has the opposite sign. Combining this with the unshifted staircase energy identity gives \[
 D(C)=tp[1-Z_{\rm pl}(t)]-t|C|^2
   +p^2\left[|C+u(0)|^2+2\int_0^t\frac{|C+u(x)|^2}{(1-x)^3}\,dx\right].
\] Since \(|u(x)|\le t-x\), the smallest possible integrand has \(|C+u(x)|=(r-t+x)_+\). Take \(u(x)=-s\min(t-x,r)\), using rapidly alternating frames on the constant segment. This achieves the lower envelope everywhere in the scale-separated limit. For \(r<t\) optimality here concerns this specified staircase class.
For the unshifted prism take honeycomb fraction $1/3$, tile parameter
$1/15$ and 48 subdivisions. The periodic gradients, potentials and the
complementary fourth-field trials just described are checked by rational
arithmetic. Their certified statistic exceeds $0.973$.
An isotropic matrix laminate at this fraction has zero conditional
fourth-field energy. Mix the two materials at the same fraction and
isotropize their rotations. Their second fields vanish on the core,
and equal fractions contribute no additional coarse Fourier fluctuation.
The fourth energies therefore interpolate. This realizes the value
$0.973$ exactly in the closure and restores zero orientational mean
by rotation.

For the stopped prisms the possible stopping fractions are
$g=j/48$, $1\le j\le15$, and $g=3/10$. Four scalar lower bounds
for a density at $a=3/10$ are
\[
 (Z,B,C,D)=(0.978610,0.781940,0.999410,0.872900).
\]
Here $Z$ is the normalized core energy, $B$ and $C$ are the
conditional pairing and squared norm of the interpolated disk
direction, and $D$ is the whole-cell fourth energy. Use the
Jackson interpolation of Lemma~\ref{p14p:transfer} and the grid
calculation in the proof of Proposition~\ref{p14p:strip} with
input degrees $(56,99)$ and mean $3/10$. The additional integer
test checks all four outward comparisons directly at the same
$1024\times2048$ nodes, including positivity of every normal
measure. The zero-frequency replacement uses the rational prescribed
mean. One fixed-point unit covers its truncation in the FFT.
For $g\le a$ set $c=g/a$. For $g>a$ set $c=(1-g)/(1-a)$.
Linear transport of the density and second field gives the following
attained lower bounds on the conditional fourth energy:
\[
 \begin{cases}
 c^2(1-a)Z+2c(1-c)B+(1-c)^2C,& g\le a,\\
 c^2[(1-c)D+ca(1-a)Z]/g,&g>a .
 \end{cases}
\]
In the first case scale the original core density by $c$ and add the
remaining mass in the disk direction. Expanding the squared norm gives
the three nonnegative coefficients displayed. In the second case scale
the complementary density by $c$ and the fourth field by $c$, which
gives the stated conditional average. Divide the bound by $1-g$ to get
the normalized statistic. At each of the sixteen fractions take the larger
of this bound and the ordinary staircase value proved above, then round
down to a multiple of $10^{-8}$. These two constructions suffice for
every chosen stopping seed. A rational series lower bound for
$-\log(1-g)$ suffices in the staircase expression because all its terms
are positive.
\paragraph*{A split that keeps the matrix fourth field stationary to first order}

Read a laminate from coarse to fine. At a split of angular normal \(q=e^{2i\alpha}\), a change \(s\) in conditional core fraction adds \(qs\) to \(U\) and \(q^2s\) to \(W\). A sequential fine coating from core fraction \(g\) to a pure core adds \((1-g)\langle q\rangle\) on that core. If all its layer normals have angular coordinates \(v,-v\), the fine fourth field is \(v^2(\chi-g)\) in both phases. This holds at finite rank in the scale-separation limit.

\begin{enumerate}
\def\labelenumi{\arabic{enumi}.}
\tightlist
\item
  Split along \(q=1\) into fractions \(y=(1+h)/2\) and \(x=h/2\), of volumes \(h\) and \(p\). Complete the \(y\)-branch by layering along \(q=-1\). On its core \(U=0\), and on its matrix \(W=-h\).
\item
  On the \(x\)-branch the arriving fields are \(U=W=-x\). Take \(0<\ell<x\), let \(b=1/2\), and split along \(q=-a\) with \(a=e^{i\theta}\), \(0\le\theta<\pi/2\). The child fractions \(g=\ell,b\) occur with probabilities \(\nu_\ell=(b-x)/(b-\ell)\), \(\nu_b=1-\nu_\ell\). With \(s=g-x\), the arriving fields are \(U_g=-x-as\), \(W_g=-x+a^2s\).
\item
  In each child use the two angular normals \(v_g,-v_g\), where \(v_g=(x+as)/|x+as|\). Choose their weights with mean \((x+as)/(1-g)\). For the higher child \(|x+as|\le g=1-g\), and for the lower child \(|x+as|\le 2x-g=h-g<1-g\). Hence these weights are nonnegative. The final core has \(U=0\) exactly. Average this tree with its mirror at the same \(h\) to make \(c\) real. This operation introduces zero coarse increment and preserves \(z\).
\end{enumerate}

On a final matrix point of either inner branch the displacement of \(W\) from its square value is exactly \[
 D_g=W_g-gv_g^2+h
   =\frac{xs(a-1)^2(a+1)}{a(x+\bar a s)}.                              \tag{2}
\] The outer branch has \(D=0\). In the conditional matrix distribution the three masses are \(\nu_\ell(1-\ell),\nu_b(1-b),h/2\), with conjugation on half the cells. Write \(A=\mathbb E_B\operatorname{Re}D\) and \(B=\mathbb E_B|D|^2\). Parseval and the mean of \(W\) give the exact formulas \[
 c=1-A/h,\qquad z=p+2A-B/h,\qquad
 z+tc-(p+t)=[(2h-t)A-B]/h.                                             \tag{3}
\] For each fixed \(\ell\), \(D_g=d_g\theta^2+O(\theta^3)\) with \(d_g=-2x(g-x)/g\). Its conditional mean coefficient is \[
 m_\ell=\mathbb E_B d=\frac{h(h-2\ell)(1-h)}{2\ell}>0.
\] Thus \(A=m_\ell\theta^2+O(\theta^4)\) and \(B=O(\theta^4)\). For every \(t<2h\), (3) is strictly positive at sufficiently small nonzero tilt. For \(t\ge2h\), the prior unrestricted inequality \(z\le1-hr^2\) has maximum \(p+t\) after adding \(tr\), attained by the untilted square. This proves both directions of the transition.
\paragraph*{A second realization by stopping an orthogonal staircase}

The earlier shifted staircase gives a second proof of the strict transition and the lower curvature, with stronger quantitative gains at finite tilt. Write \(d=h-t/2\) for \(0\le t\le2h\). Take its first orthogonal step from \(h\) to \(b\le h\), with real fourth spin \(+1\). The stopped matrix mass is \(p(h-b)/(1-b)\), with \(W=-h\). The surviving branch has volume \(p^2/(1-b)^2\), receives fourth increment \(-(h-b)\) and receives zero second increment. Insert a planar exact-bulk seed there. Rotate it before insertion so its core mean parameter \(c_b\) is real in this frame. The matrix energy in a real background of size \(h-d\) is \[
 \frac{p(h-b)}{1-b}d^2+
 \frac{p^2}{1-b}\big[(b-d)^2+b(1-z_b-2(b-d)c_b)\big].
\] Indeed the seed sees residual background \(b-d\), its matrix mean is \(-bc_b\), and its integrated matrix energy is \(b(1-b)(1-z_b)\). Any completed planar tree of statistics \(z,c\) has matrix energy \(p(h-d)^2+hp[1-z-t\operatorname{Re}c]\) in this background. Equating the expressions gives the attained support \[
 p+t+\frac{p}{h}\frac{b}{1-b}
       [\,2d-b-1+z_b+2(b-d)c_b\,].                                    \tag{4}
\] In particular a rotation-averaged seed has \(c_b=0\) and the completed grain has \(c=(h-b)/(h(1-b))\). This relates the oriented edge at \(h\) to attainable seeds over a range of other fractions.

Set \(b=d\) and choose the earlier rotated staircase seed \(z_b=(1-b)[-b-2\log(1-b)]/b\), with its limiting value at zero. Then (4) gives \[
 p+t+\frac{2p}{h}[-\log(1-d)-d].
\] Its gain over the square is strictly positive for \(d>0\). Its expansion has the same quadratic coefficient as the all-laminate upper bound, which independently gives the sharp support curvature above. The finite tree in (2) attains the tangent at fixed rank and attains the same limiting curvature by varying its lower child fraction.

\paragraph*{Which fillings enter the submitted mesh.}
There are 3,264,000 prism elements in the certificate. Of these,
1,890,462 use the unshifted honeycomb, 10,416 the pure-core seed,
924 the finite split, and 1,362,198 a stopped seed. The masks
therefore use only the planar constructions verified above. The
independent checker validates the menu index in each element before
taking its rational energy and variance.
\phantomsection\label{endfile:mesh-fillings}

\subsection{Two further complete finite-phase sections}\label{route:two-further-complete-finite-phase-sections}
The following sections retain the two narrow strips for which explicit
Fourier seeds and transfer estimates give both bulk endpoints.
\subsection{A planar section at positive moduli}\label{route:a-planar-section-at-positive-moduli}

\begin{theorem}\label{p14p:strip}
For equal fractions of the planar phases $(\kappa,\mu)=(1,1)$ and
$(36/35,16/15)$, at every prescribed shear in
$[1.0325242911,1.0325242947]$ the complete closed isotropic bulk
section is $[143/141,2216/2185]$.
\end{theorem}

\begin{proof}\phantomsection\label{proof:p14p:strip}
Call the core phase $(l,n)$ and the matrix $(k,m)$, and write $f$ for the
core fraction, $p=1-f$. We use both assignments. On the torus let $T_j$
have symbol $e^{2ij\alpha}$ at each nonzero wavevector of angle $\alpha$
and zero at frequency zero. For the core indicator $\chi$ put
$U=T_1\chi$, $W=T_2\chi$, $z=\langle\chi|W|^2\rangle/(fp)$.
Limits of the angular energies below refer to actual indicators along
separated-scale sequences. Lemma~\ref{p14:hydro} gives the matrix-side bulk
$K_s=[f/(l+m)+p/(k+m)]^{-1}-m$ when
$\|\chi U\|_2\to0$, since the two components of the traceless Hessian in that
formula are proportional to $\operatorname{Re}U,\operatorname{Im}U$.

\smallskip\noindent\emph{How the fourth angular field controls shear.}
Use Frobenius-orthonormal hydrostatic and two deviatoric coordinates
$(I,\operatorname{diag}(1,-1),\left(\begin{smallmatrix}0&1\\1&0\end{smallmatrix}\right))/\sqrt2$
and divide stiffnesses by two. In the reference medium with the two
given bulks and common shear $m$, constant deviatoric strain is an
exact solution. Let $\Gamma$ send a polarization $F$ to the zero-mean
compatible strain $Y$ solving
$\langle Y:C_{\rm ref}V\rangle=\langle F:V\rangle$ for all such $V$.
Restrict its input and output to deviators in the core with the orthogonal
projection $P_D$. The compression $B=P_D\Gamma P_D$ is self-adjoint
and $0\le B\le 1/m$, since the strain energy controls $m$ times the squared
deviatoric norm. For a unit deviatoric mean $e$, the projected field after
the shear change $\delta=n-m$ solves
$(I+\delta B)P_D E=P_D e$. Averaging the corresponding two energies,
and using the spectral theorem with the two vectors $P_De$ of squared
norm $f$, gives a probability $\nu$ on $[0,b]$, $b=1/m$, with
\begin{equation}\label{p14p:measure}
 M(\delta)=m+f\delta\int(1+\delta v)^{-1}\,d\nu(v).
\end{equation}
Here $M$ is the normalized trace on the shear space, equal to physical
shear for an isotropic tensor.

Set $a=(k+m)^{-1}$, $g=(a+b)/2$, $c=(b-a)/2$.
At constant reference bulk $k$ the symbol of $\Gamma$ is
$C_0^{-1}-C_0^{-1}vv^t C_0^{-1}/(v^t C_0^{-1}v)$, where
$C_0=\operatorname{diag}(k,m,m)$ and
$v=(1,-\cos2\alpha,-\sin2\alpha)^t$ spans the orthogonal complement
of compatible strains. Applied to the two core deviatoric loads, its
hydrostatic row is $a(\operatorname{Re}U,\operatorname{Im}U)$, and its
shear block is
\[
 g(\chi-f)I-c
 \begin{pmatrix}\operatorname{Re} W&\operatorname{Im} W\\
                 \operatorname{Im} W&-\operatorname{Re} W\end{pmatrix}.
\]
Changing just the core bulk changes the solution by $O(\|\chi U\|_2)$
in $L^2$, as one sees by subtracting the two weak equations and bounding the right hand
side by $|l-k|$ times the core hydrostatic strain norm. The displayed
traceless matrix has square $|W|^2 I$. Taking the mean and mean squared
norm of the two compressed fields therefore gives, along a subsequence
with weakly convergent spectral probabilities,
\[
 h:=\int v\,d\nu=pg,\qquad
 \int (v-h)^2\,d\nu=pc^2 z .
\]
The same limiting measure gives \eqref{p14p:measure} for every
$\delta>-m$ by continuity of the integrands.

We need two inverse moments. Positivity of $M(\delta)$ as
$\delta\downarrow-m$ gives $\int (b-v)^{-1}d\nu\le m/f=:C_b$.
For the other end let $H$ denote the symmetric field of mean $f I/2$
whose nonconstant Fourier coefficients are
$(\xi\otimes\xi/|\xi|^2)\widehat\chi(\xi)$, and set $Q=H-\chi I/2$.
For any unit traceless load $F$ the trial
$E=F+s[2\operatorname{sym}(FH)-fF]$ is compatible, since each
oscillatory mode is a symmetric product with $\xi$. Its deviator is
$(1+s(\chi-f))F$ and its trace is $2s F:Q$. At $s=-1/p$ the core
deviator vanishes. On the edge $\|\chi Q\|_2\to0$, and
$\langle|Q|^2\rangle=fp/2$ by Parseval. Averaging the energies on
the two loads gives $M(\delta)\le m/p+kf/(2p)$ for every
$\delta>-m$. Sending $\delta\to\infty$ in \eqref{p14p:measure}
proves $\int v^{-1}d\nu\le (m+k/2)/p=:C_0$.

Put
\[
 P(v)=v^2-(1+p)gv+pab,\qquad
 q=1+(1+p)g\delta+pab\delta^2 .
\]
$P(0)>0$, $P(h)=-pc^2<0$ and $P(b)=fbc>0$. Let $\nu_*$ put mass
on the two roots so that its mean is $h$. Its variance is $pc^2$. Evaluating the two-node transform gives
$G_*=(1+g\delta)/q$, in particular $q>0$.
At $\delta\to\infty$ and $\delta\downarrow-m$ the same formula gives its
inverse moments $C_0,C_b$. Write $A=P(v)^2/[v(b-v)]$,
$d_0=P(0)^2/b$, $d_b=P(b)^2/b$. Its partial fractions consist of
$-v^2$, an affine function, $d_0/v$ and $d_b/(b-v)$.
Thus $pc^2(1-z)$ equals the sum of the three nonnegative quantities
\[
 \int A\,d\nu,\quad d_0(C_0-\textstyle\int v^{-1}d\nu),\quad
 d_b(C_b-\textstyle\int(b-v)^{-1}d\nu).
\]
Set $\lambda=\delta^2(1+b\delta)/q^2$. The partial fractions of
$\lambda A/(1+\delta v)$ have the extra term $-(1+\delta v)^{-1}$,
since $\lambda A(-1/\delta)=-1$, and the polynomial part is affine.
Subtracting integrals against $\nu_*$ and $\nu$ shows that
$G_*-\int(1+\delta v)^{-1}d\nu$ is $\lambda$ times the same three
quantities with respective weights $(1+\delta v)^{-1}$ inside the first
integral, $1$, and $(1+b\delta)^{-1}$. In particular
\begin{equation}\label{p14p:error}
 |M-(m+f\delta G_*)|\le
 \frac{fp c^2|\delta|^3}{q^2}\max(1,n/m)(1-z).
\end{equation}

\smallskip\noindent\emph{A core with zero second field and large fourth field.}
On a rectangular torus of side ratio $15:26$, the finite certificate uses
a real trigonometric polynomial $\theta$ of rectangular degree $(114,199)$,
mean $1/2$, and coefficients in $2^{-42}(\mathbb Z+i\mathbb Z)$.
On the $N_x\times N_y=1024\times2048$ grid, put
$r=1-\theta$, $u=T_1\theta$, $w=T_2\theta$.
The certificate verifies $\theta>0$ and $r>|u|$ at every node.
For a probability on unit complex normals $q'=e^{2i\alpha}$ with
mean $m_1=-u/r$ one can prescribe its second moment as
$m_1^2+(1-|m_1|^2)s$, $|s|\le1$. Indeed, for $|s|=1$ take the two
intersections of the line through $m_1$ parallel to $\sqrt{s}$ with
the circle and choose the weights with mean $m_1$. The centered
second moment follows from $\int|q'-m_1|^2=1-|m_1|^2$.
Mixture gives the interior of the disk. With one admissible $s$ at each
node define
\[
 V=w+u^2/r+(r-|u|^2/r)s .
\]

In each coordinate of grid size $N$ and original degree $d$, set
$m_0=\lfloor(N-d-2)/4\rfloor$, $D=2m_0-2$. Normalize the kernel
$|\sum_{j=0}^{m_0-1}e^{ijx}|^4$ to have integral mean one. Its Fourier
coefficients are the autocorrelation of the integer tent $m_0-|j|$,
$|j|<m_0$, divided by the coefficient at zero. Let $J$ interpolate
nodal data with the translated product kernels divided by $N_xN_y$.
The translates form a nonnegative partition of unity since $D<N$. Given $d+D<N$, interpolation of a polynomial of degree $d$ only
multiplies its original frequencies by the kernel coefficients.
It therefore commutes with $T_1,T_2$ on $\theta$.
Use coarse core density $\Theta=J\theta$, and at each point use the
normal probability $J(r\mu)/(1-\Theta)$, where $\mu$ denotes the nodal
probabilities just constructed. Its unnormalized first moment is $-Ju$,
and the sum of its unnormalized second moment and $Jw$ is $JV$.

To realize these moments, on a fine matrix-coated sequential laminate
with core fraction $x$, the $j$th angular field on the surviving core
receives $(1-x)\int(q')^j d\mu$ beyond the coarse field. In fact, along
the nested core branch choose successive conditional fractions from $x$
to $1$ increasing by $(1-x)$ times the masses of its normal atoms.
Write these fractions as $x=b_0<\cdots<b_R=1$. At stage $i$ the
indicator $\lambda_i$ of the retained layers has mean
$\rho_i=b_{i-1}/b_i$ and normal angle $\alpha_i$. Telescoping gives
\[
 \chi-x=\prod_{i=1}^R\lambda_i-\prod_{i=1}^R\rho_i
 =\sum_{i=1}^R\left(\prod_{k<i}\lambda_k\right)
 (\lambda_i-\rho_i)\left(\prod_{k>i}\rho_k\right).
\]
On the surviving core the $i$th term equals
$(1-\rho_i)b_i=b_i-b_{i-1}$. In the separated-scale limit $T_j$
multiplies that term by $e^{2ij\alpha_i}$ up to an error tending to
zero in $L^2$. Indeed, for a slow Fourier mode $k$ multiplied
by a fast nonzero mode $Nh$, the multiplier at $k+Nh$ tends to that at
$h$. Finite Fourier sums followed by $L^2$ approximation prove the
rule for layers. More explicitly, bounded Fejér sums approximate all
slow indicator factors and all single-layer indicators in $L^2$, and their
bounded products converge in the iterated scale limit by integration
over independent fast periods. Since the multipliers have $L^2$ norm
one, the field errors tend to zero. For convergence of the restricted quadratic energies,
replace the fields by the bounded polynomial
approximants and then use Cauchy-Schwarz on the field difference.
Approximate $\Theta$ and the atom masses by their averages on fine
boxes (for the latter average the unnormalized measures). Insert the local
laminates at these box fractions by Lemma~\ref{p13:locality}. The
volume-weighted sum of the box fractions is exactly the original mean.
The multipliers of the averaged densities
converge in $L^2$, and the local first and second moments converge
uniformly since $1-\Theta$ is bounded below. General ply directions
follow by rational approximation. These constructions give a sequence
with vanishing core second field and limiting fourth statistic
\[
 z=4\langle\Theta|JV|^2\rangle>0.947060 .
\]
To enforce isotropy at every test modulus, fill one sector of a regular
hexagon with whole fine periods of one texture. Fill the remaining strip
of that sector to the same core fraction. Rotate this entire construction
into the five other sectors by successive multiples of $60^\circ$ and
repeat the hexagon by translations. Each resulting pattern has the prescribed
fraction and sixfold symmetry. Its elastic tensor is isotropic, since the
rotation acts by $120^\circ$ on the deviatoric plane. The equal sector
fractions make the coarse angular increments zero. The fine angular energies
are preserved by rotation and the layer localization above. For a Fourier
mode of a rotated fine pattern, take the nearest vector $k_L$ of the hexagonal
reciprocal lattice to its frequency $Lh$. The modulation
$1_{\rm sector}(x)e^{i(Lh-k_L)\cdot x}$ ranges over a compact set in $L^2$
because the offsets are bounded. Its Fourier truncations thus converge
uniformly on that set. The multiplier on each retained shift of $k_L$
tends to that at $h$. The filled strips have vanishing measure, so after
bounded Fourier approximation they do not affect these limits.

We put the hexagonal response in the square-periodic closure at the exact
fraction. Fix one hexagonal pattern, scale it by $1/L$, and retain the whole
scaled hexagons inside the unit square. Fill the uncovered part to the same
fraction and periodize the square. This changes the pattern only within
$O(1/L)$ of the boundary. For an imposed mean strain $E$ and mean stress $T$,
let $w$ be the hexagon-periodic displacement corrector and $\psi$ the periodic
Airy potential for the mean-zero part of the minimizing stress. The latter
belongs to $H^2$ by Fourier inversion, since every symmetric stress mode
at $k\ne0$ annihilates $k$ and is a multiple of $k^\perp\otimes k^\perp$.
Choose a smooth cutoff $\zeta_L$, zero near the boundary of the square and
one beyond a strip of width $c/L$ containing all incomplete cells, with
$|D^j\zeta_L|\le C_j L^j$ for $j=1,2$. The two square-periodic trials are
\[
 E+e\big(\zeta_L w(Lx)/L\big),\qquad
 T+\operatorname{cof}D^2\big(\zeta_L\psi(Lx)/L^2\big).
\]
They have the imposed means and the second is equilibrated. In the strip
the squared norms integrate to $O(1/L)$, because it meets $O(L)$ cells
and the product rule bounds the integral on each by $C L^{-2}$ times
the squared $H^1$ and $H^2$ corrector norms on one fixed period, including
the constant loads. In the remaining whole cells the trials give the original
cell energies. The strain and complementary upper bounds therefore squeeze
the effective tensors to the hexagonal tensor. The nearest-frequency argument above, now on the square reciprocal lattice, also preserves the limiting core angular energies because the replaced strip has vanishing area. For scalar conductivity
the same cutoff applies to the potential and to the stream function of its
mean-zero current, using their gradient and rotated gradient respectively.
One $L$ suffices for any finite list of tests on the same geometry.
Diagonal selection gives the joint limit at all positive rational test
shears. On each compact positive shear interval uniform ellipticity bounds
the squared norms of unit-load minimizing strains, and comparison at two
nearby shears gives a common Lipschitz bound. The limit consequently has
the hexagonal response at every positive test shear.

The certificate in \path{certificates/baseline/p14planar/check.py} uses the dyadic coefficients in \path{certificates/baseline/p14planar/density.npz}. It asserts the last strict inequality and every nodal feasibility test. Since $d+2D<N$, the integral is the
exact nodal average of $4\Theta|JV|^2$. The verifier encloses this average
by integer radix-two transforms. In each butterfly it stores a complex
dyadic midpoint in units of $2^{-55}$ and carries a common Euclidean
error in those units. Roots have error at most four units by half-angle
square-root bounds and interval multiplication. If the two input errors
are at most $E$ and the right input midpoint before rotation has $\ell^1$ norm at most $B$
in physical units, the new radius is $2E+4\lceil B\rceil+2$.
For a normalized forward step divide this radius by two, round up, and
add two units for the final truncation. This follows by using the
unit modulus of the exact root. Multiplication by the initial rational
Fourier symbols with frequencies $(26k,15l)$ and by the rational kernel
weights is enclosed by integer division. The proposed dyadic disk
directions have squared norm checked at most one. The feasibility
test puts both the midpoint $(r,u)$ and its exact value in the convex
cone $|u|\le r$, with $r>0$ throughout their joining segment.
Along that segment both $u^2/r$ and $|u|^2/r$ have differential
bounded by $2|du|+|dr|$. Thus if $e_w,e_u,e_r$ enclose the input
errors, the radius for $V$ is $e_w+4e_u+3e_r$ plus six units for
integer truncations. The third $e_r$ accounts for the linear term $rs$.
Positivity of $J$ bounds the propagated input errors before adding
the midpoint-transform errors. If the final nodal errors in $\Theta,JV$
are $e_0,e_1$ in physical units, the verifier checks that the
$\ell^1$ norm of the latter midpoint is at most $3$ and $e_1<1$.
It subtracts
$4(9e_0+7e_1)$ from the dyadic average and compares the exact integers.
Only the degree-$(114,199)$ coefficients enter this certificate.

Finally, use also the isotropic matrix laminate with angular moments
$\int q'\,d\mu=\int(q')^2\,d\mu=0$ (three equally spaced values on the
unit circle suffice). At fraction $f$ it has $U=W=0$ on the core by the
layer rule. Thus the spectral variance is zero and its shear is
$M_s=m+f\delta/(1+pg\delta)$. For $(k,m)=(1,1)$ and
$(l,n)=(36/35,16/15)$, substitution at $f=1/2$ gives $K_s=143/141$,
$M_s<1.0325242911$ and the lower bound from \eqref{p14p:error}
exceeds $1.0325242947$. Interchanging the phases gives
$K_s=2216/2185$, $M_s>1.0325242947$ and the upper bound from
\eqref{p14p:error} is below $1.0325242911$. These are rational comparisons
using $z>947060/10^6$. Lemma~\ref{p14:mix} fills the common shear interval
at each of the two extreme bulks, then fills the interval between these
bulks at each such shear. Its bulk bounds exclude every exterior value.
\phantomsection\label{proofend:p14p:strip}\end{proof}

\begin{lemma}\label{p14p:transfer}
For each $a\in\{1/2,3/5\}$ there is a real trigonometric polynomial
$\theta$ of mean $a$ on the $15:26$ rectangular torus, of rectangular
degree $(114,199)$ for $a=1/2$ and $(56,99)$ for $a=3/5$.
Let $T_j$ have symbol $e^{2ij\alpha}$ on wavevectors of angle $\alpha$
and vanish on constants. At every node of the $1024\times2048$ grid,
$\theta>0$ and $r=1-\theta>|u|$, where $u=T_1\theta$.
There exist nodal choices $s,d$ in the closed unit disk such that
\[
 V=T_2\theta+u^2/r+(r-|u|^2/r)s,\qquad
 \Theta=J\theta,\quad V_c=JV,\quad S=Jd.
\]
Here $J$ takes the grid average against translated products of the
unit-mean kernels proportional to $|\sum_{j=0}^{m_0-1}e^{ijx}|^4$,
with $m_0=\lfloor(N-b-2)/4\rfloor$ on an axis of grid size $N$
and density degree $b$. The nodal ply probabilities $\mu$ can be
chosen with first moment $-u/r$ and second moment $(V-T_2\theta)/r$.
The sequential laminates with unnormalized ply measures $J(r\mu)$
have vanishing core second field in the limit and core fourth field
$V_c$ at coarse density $\Theta$. The same interpolation gives $|S|\le1$.
For these fields define
\[
 z=\frac{\langle\Theta|V_c|^2\rangle}{a(1-a)},\quad
 B=\frac{\langle\Theta\operatorname{Re}(V_c\overline S)\rangle}{a},
 \quad C=\frac{\langle\Theta|S|^2\rangle}{a},\quad
 D_0=\langle|V_c|^2\rangle.
\]
The following entries are simultaneous strict lower bounds.
\[
\begin{array}{c|rrrr}
 a&z&B&C&D_0\\ \hline
 1/2&0.947060&0.61840&0.99962&0.68090\\
 3/5&0.91570&0.53150&0.99970&0.54780
\end{array}
\]
\end{lemma}
\begin{proof}\phantomsection\label{proof:p14p:transfer}
Use the same interpolation, with the preceding density at half fraction
and a second polynomial of degrees $(56,99)$ and exact mean $3/5$ from \path{certificates/baseline/p14planar/density6.npz} in the second row. The nonconstant coefficients have the same dyadic denominator.
For $a=3/5$ the midpoint of the constant coefficient in units of
$2^{-55}$ is $\lfloor(3/5)2^{55}\rfloor$. The transform starts with
radius one unit, enclosing the exact mean. The integer tests establish
$\theta>0$ and $r>|u|$ at every node for both polynomials.
At each node choose another dyadic direction in the unit disk near the
direction of the proposed value of $V$, giving the nodal field $d$. Let $S$ be its positive interpolant,
and set $V_c=JV$. The three added integrands have degrees at most $d+2D$
in each coordinate, so the same quadrature is exact. For node errors
$e_0,e_1,e_s$ in $\Theta,V_c,S$, use $|V_{c,\rm mid}|\le3$,
$|S_{\rm mid}|\le2$, $|S|\le1$, $e_1<1$. The allowances on the
unnormalized averages for $B,C,D_0$ are then respectively
$6e_0+e_1+3e_s$, $4e_0+3e_s$ and $7e_1$. The verifier checks these
midpoint bounds and compares all four lower sums against the rational
thresholds at each fraction. The grid verifier and the spatial strip
calculation both read \path{certificates/baseline/p14planar/seed-bounds.txt}, which
lists the fractions in tenths and the four bounds in millionths.
\phantomsection\label{proofend:p14p:transfer}\end{proof}

\begin{proposition}\label{p14:last}
For equal amounts of the spatial phases $(500,1)$ and $(625,5/4)$,
at each shear in the interval $[1.11941936,1.11942636]$, the closed interval of attainable bulk moduli is
\[ [1879500/3383,\;376125/677]. \]
\end{proposition}
\begin{proof}\phantomsection\label{proof:p14:last}
Use the spatial layers of Proposition~\ref{p14:sections}, replacing the
planar material inside each stopping layer. Lemma~\ref{p13:locality} preserves the volume-weighted fraction at each substitution. Let $a$ be one of the seed
fractions $1/2,3/5$ in Lemma~\ref{p14p:transfer}. Write $\Theta,V,S$
for its core density, core fourth field and interpolated disk direction
($|S|\le1$), and abbreviate the four seed integrals by
\[
 A=\langle\Theta|V|^2\rangle/a,\quad
 B=\langle\Theta\operatorname{Re}(V\overline S)\rangle/a,\quad
 C=\langle\Theta|S|^2\rangle/a,\quad D_0=\langle|V|^2\rangle.
\]
For any desired transverse fraction $h$ the conditional fourth-field energy
$(1-h)z$ is at least
\[
 \begin{cases}
 c^2 A+2c(1-c)B+(1-c)^2 C, &h\le a,\quad c=h/a,\\
 c^2[(1-c)D_0+caA]/h,     &h\ge a,\quad c=(1-h)/(1-a).
 \end{cases}
\]
Indeed in the first case use density $c\Theta$, scale the original
unnormalized ply weights by $c$ and add mass $1-c$ with zero first
moment and second moment $S$. Such a distribution exists by the unit
moment disk at first moment zero. The second field still cancels on the
core and the fourth field becomes $cV+(1-c)S$. In the second case set
the density to $1-c(1-\Theta)$ and scale all original unnormalized ply
weights by $c$. Both angular multipliers kill constants, so the core
second field again vanishes and the fourth becomes $cV$. Integration
gives both formulas. The realization and sixfold rotation in the proof
of Theorem~\ref{p14p:strip} apply to the modified densities.

These planar materials have the scalar conductivity \eqref{p14:scalar}
required for axial shear. To see this directly let $\Pi$ denote the
orthogonal projection onto periodic mean-zero gradients. At core fraction
$h$ its Fourier symbol gives
$\Pi(\chi a)=[(\chi-h)I/2+Q]a$ for constant $a$, with $Q$ as in
Lemma~\ref{p14:hydro}. Under mean electric field $v$ use field
$v-\Pi(\chi a)/m$ and flux $mv+(I-\Pi)(\chi a)$, where
$a=(n-m)v/[1+(1-h)(n-m)/(2m)]$. The constitutive residual is zero
in the matrix and bounded in $L^2$ on the core by $|n-m||a|\,\|\chi Q\|_2/m$.
It tends to zero by cancellation of the second field. Testing the weak
conductivity equation on the field difference bounds its $L^2$ norm by
the residual divided by $\min(m,n)$. The mean flux tends to $mv+ha$,
as required. The same cancellation gives the transverse dilation field
used in Proposition~\ref{p14:sections}, where the constant axial strain in
the extrusion changes only its planar trace drive, as the calculation
of the block $A,B,C$ shows.

At each stopping fraction $h$, take the greater of the two witnesses
from the seed integrals and divide by $1-h$. Estimate \eqref{p14p:error},
with planar phases $(k+m/3,m)$ and $(l+n/3,n)$, then bounds the texture's
physical shear $P$ from below for $n>m$ and from above for $n<m$ by using
the corresponding side of its error interval. Insert that bound in the
extruded normal block and take the same 256 spatial steps with the same
outward rounding. The rational check \path{certificates/baseline/p14sections/last_strip.py}
asserts every hydrostatic row sum and the final one-sided shears
$1119441455156/10^{12}$ from below with matrix $(500,1)$ and
$1119418144104/10^{12}$ from above with matrix $(625,5/4)$.
They enclose the stated strip together with their respective isotropic
matrix-laminate shears, as checked by the formula for $g$ in
Proposition~\ref{p14:sections}. Apply Lemma~\ref{p14:mix}.
\phantomsection\label{proofend:p14:last}\end{proof}

\end{document}